\documentclass[twoside]{book}
\usepackage{akp1} 
\usepackage{biber}
\usepackage{hyper} 
\usepackage{letters}

\usepackage{notation}

\usepackage{hyphenation}
\usepackage{envmnt}
\usepackage{notation-vitya}
\usepackage{tikz-cd}
\usepackage[toctitles]{titlesec}

\makeindex

\begin{document}

\title{Alexandrov geometry: foundations}
\date{}
\author{Stephanie Alexander, Vitali Kapovitch,\\ and Anton Petrunin}
\maketitle
\thispagestyle{empty}

\chapter*{Preface}

Alexandrov spaces are defined via axioms similar to those given by Euclid.
The Alexandrov axioms replace certain  equalities with inequalities. 
Depending on the signs of the inequalities, we obtain Alexandrov spaces with {}\emph{curvature bounded above} (CBA) and {}\emph{curvature bounded below} (CBB).
The definitions of the two classes of spaces are similar, but their properties and known applications are quite different.

The goal of this book is to give a comprehensive exposition of the structure theory of Alexandrov spaces 
with curvature bounded above and below.
It includes all the basic material as well as selected topics inspired by considering the two contexts simultaneously.
We only consider the intrinsic theory, leaving applications aside. 
Our presentation is linear,
with a few exceptions where topics are deferred to later chapters to streamline the exposition.
This book includes material \emph{up to the definition of dimension}.
Another volume still in preparation will cover further topics.

\section*{Brief history}

The first synthetic description of curvature is due to Abraham Wald \cite{wald}; 
it was given in a lone publication on a ``coordinateless description of Gauss surfaces'' published in 1936.
In 1941, similar definitions were rediscovered by Alexandr Alexandrov \cite{alexandrov:def}.

In Alexandrov's work the first fruitful applications of this approach were given.
Mainly: {}\emph{Alexandrov's embedding theorem} \cite{alexandrov-1941,alexandrov-1941convex}, which describes closed convex surfaces in Euclidean 3-space,
and the {}\emph{gluing theorem} \cite{alexandrov-1946}, which gave a flexible tool to modify non-negatively curved metrics on a sphere.
These two results together gave an intuitive geometric tool to study embeddings and bending of surfaces in Euclidean space and changed the subject dramatically.
They formed the foundation of the branch of geometry now called {}\emph{Alexandrov geometry}.

\parbf{Curvature bounded below.}
The theory grew out of studying intrinsic and extrinsic geometry of convex surfaces without the smoothness condition.
It was developed by Alexandr Alexandrov
and his school.
Here is a very incomplete list of contributors to the subject:
Yuriy  Borisov,
Yuriy  Burago,
Boris Dekster,
Iosif  Liberman,
Sergey  Olovyanishnikov,
Aleksey  Pogorelov,
Yuriy  Reshetnyak,
Yuriy  Volkov,
Viktor  Zalgaller.

The first result in higher dimensional Alexandrov spaces was the splitting theorem.
It was proved by Anatoliy Milka \cite{milka-line} and appeared in 1967.
Milka used a global definition similar to the one used in this book. 

{\sloppy 

In the 80's the interest in convergence of Riemannian manifolds spurred by \emph{Gromov's compactness theorem} \cite{gromov-MS} turned attention toward the singular spaces that can occur as limits of Riemannian manifolds.
Immediately it was recognized that if the manifolds have a uniform lower sectional curvature bound, then the limit spaces have a lower curvature bound in the sense of Alexandrov.
There followed during the 90's an explosion of work on intrinsic theory of Alexandrov spaces  starting with papers of Yuriy Burago, Grigory Perelman, and Michael Gromov  \cite{burago-gromov-perelman,perelman:spaces2}.
Similar ideas were developed independently by Karsten Grove
and Peter Petersen, whose work was not converted into a publication, and also by
Conrad Plaut~\cite{plaut-preprint}.

}

Around the same time an implicit application of higher-dimensional Alexandrov geometry was given by Michael Gromov in his bound on Betti numbers \cite{gromov:betti}.
Another implicit application, which essentially used Alexandrov geometry before it was was actually introduced, given later by Wu-Yi Hsiang and Bruce Kleiner in their paper on non-negatively curved 4-manifolds with infinite symmetry groups \cite{hsiang-kleiner}.
The work of Hsiang and Kleiner and its extension by Karsten Grove and Burkhard Wilking \cite{grove-wilking} are some of the most beautiful applications of this branch of Alexandrov geometry.

The above activity was very much related to so-called {}\emph{comparison geometry},
a branch of differential geometry that compares Riemannian manifolds  to  spaces of constant curvature.
In addition to the already-mentioned {}\emph{Gromov's compactness theorem},
the following results had a big influence on the development of Alexandrov geometry:
{}\emph{Toponogov comparison theorem} \cite{toponogov-globalization+splitting}, which is a generalization of the theorem of Alexandrov \cite{alexandrov-comparison};
{}\emph{Toponogov splitting theorem} \cite{toponogov-globalization+splitting}, which is a generalization of Cohn-Vossen's theorem \cite{cohn-vossen_line};
{}\emph{Finiteness theorems} of
Cheeger
and
Grove--Petersen \cite{cheeger-finiteness,grove-petersen:finiteness};
Gromov's bound on the number of generators of the fundamental group 
\cite[1.5]{gromov:almost-flat};
and 
{}\emph{Yamaguchi fibration theorem} \cite{yamaguchi-fibration}.

Let us give a list of available introductory texts on Alexandrov spaces with curvature bounded below: 
\begin{itemize}
\item The first introduction to Alexandrov geometry is given in the original paper of Yuriy Burago, Michael Gromov, and Grigory Perelman \cite{burago-gromov-perelman} 
and its extension \cite{perelman:spaces2} written by Perelman.
\item A brief and reader-friendly introduction was written by Katsuhiro Shiohama \cite[Sections 1--8]{shiohama}.
\item \cite[Chapter 10]{burago-burago-ivanov} gives another reader-friendly introduction, written by Dmiti Burago, Yuriy Burago, and Sergei Ivanov.
\end{itemize}
In addition, let us mention two surveys, one by Conrad Plaut \cite{plaut:survey} and the other by the third author \cite{petrunin:survey}.

{\sloppy

\parbf{Curvature bounded above.}
The study of  spaces with curvature bounded above started later,
inspired by analogy with the theory of curvature bounded below.
The first paper on the subject was written by Alexandrov \cite{alexandrov:strong-angle}, appearing in 1951.
An analogous weaker definition was considered earlier by Herbert Busemann \cite{busemann-CBA}.

}

Contributions to the subject were made by
Valerii Berestovskii, 
Arne Beurling, 
Igor Nikolaev,
Dmitry Sokolov,
Yuriy Reshetnyak,
Samuel Shefel; this list is not complete as well.
The most fundamental results were obtained by Yuriy Reshetnyak.
They include his {}\emph{majorization theorem} and {}\emph{gluing theorem}.
The gluing theorem states that if two non-positively curved spaces have isometric convex sets, then the space obtained by gluing these sets along an isometry is also non-positively curved.

The development of Alexandrov geometry was greatly influenced by the {}\emph{Hadamard--Cartan theorem}.
Its original formulation states that the exponential map at any point of a complete Riemannian manifold with nonpositive sectional curvature is a covering.
In particular, it implies that the universal cover is diffeomorphic to Euclidean space of the same dimension. 
See further discussion below (\ref{thm:hadamard-cartan}).

An influential implicit application of Alexandrov spaces with curvature bounded above can be seen in {}\emph{Euclidean buildings}, introduced by Jacques Tits as a means to study algebraic groups.

Here is a list of available texts covering the basics of Alexandrov spaces with curvature bounded above: 
\begin{itemize}
\item The book of Martin Bridson and Andr\'e Haefliger \cite{bridson-haefliger} gives the most comprehensive introduction available today. 
\item The lecture notes of Werner Ballmann \cite{ballmann:lectures, ballmann:notes} include a brief 
and clear
introduction.
\item \cite[Chapter 9]{burago-burago-ivanov} gives another reader-friendly introduction, by Yuriy Burago, Dmitry Burago, and Sergei Ivanov.
\item A book  by the three authors of the present volume  \cite{alexander-kapovitch-petrunin-CAT} gives an introduction aiming at reaching interesting applications and theorems with a minimum of preparation.
\item The book of Jürgen Jost \cite{jost:book} gives a more analytic viewpoint to the subject.
\end{itemize}

One of the most striking applications of $\CAT0$ spaces was given by Dmitry Burago, Sergei Ferleger, and Alexey Kononenko \cite{burago-ferleger-kononenko1998-1},
who used them to study {}\emph{billiards}; this idea was developed further in \cite{burago-ferleger-kononenko1998-2,burago-ferleger-kononenko1998-3,burago-ferleger-kononenko1998-4,burago-ferleger-kononenko2000,burago-ferleger-kononenko2001}. 
Another beautiful application is the construction of {}\emph{exotic aspherical manifolds} by Michael Davis \cite{davis:aspherical}; related results are surveyed in \cite{davis:exotic,charney-davis-1995}.
Both of these topics are discussed in \cite{alexander-kapovitch-petrunin-CAT}.
The study of group actions on $\CAT 0$ spaces and $\CAT 0$ cube complexes played a key role in the proof of the {}\emph{virtually fibered conjecture} that a finite cover of  every closed hyperbolic 3-manifold fibers over the circle.

\parbf{Satellites and successors.}
Surfaces with {}\emph{bounded integral curvature} were studied by Alexandrov's school.
An excellent book on the subject was written by Alexandr Alexandrov and Viktor Zalgaller \cite{aleksandrov-zalgaller}; see also a more up-to-date survey by Yuriy Reshetnyak \cite{reshetnyak:2D}.

Spaces with {}\emph{two-sided bounded curvature} is another subject already studied  by Alexandrov's school;
a good survey is written by Valerij Berestovskij and Igor Nikolaev \cite{berestovskii-nikolaev}.

A spin-off of the idea of synthetically defining upper curvature bounds 
was given by Michael Gromov \cite{gromov:hyp-groups}. 
He  defined so-called  {}\emph{$\delta$-hyperbolic spaces}, which satisfy   a coarse version of the  negative curvature condition, applying  in particular to discrete metric spaces.
This notion and its various generalizations such as semi-hyperbolicity (a coarse version of non-positive curvature) and relative hyperbolicity have  led to the emergence of the subject of {}\emph{geometric group theory}, which relates geometric properties of groups to their algebraic ones.
This is a well-developed subject with a large number of subfields and applications, such as the theory of small cancellation groups, automatic groups,  mapping class groups, automorphisms of free groups, isoperimetric inequalities on groups, actions on $\R$-trees, Gromov's boundaries of groups.

{\sloppy

The so-called {}\emph{curvature dimension condition} introduced by John Lott, C\'edric Villani, and Karl-Theodor Sturm gives a synthetic description of Ricci curvature bounded below; see the book of Villani \cite{villani} and references therein.

Alexandrov geometry influenced the development of {}\emph{analysis on metric spaces}. 
An excellent book on the subject was written by Juha Heinonen, Pekka Koskela, Nageswari Shanmugalingam, and Jeremy Tyson~\cite{heinonen-koskela-shanmugalingam-tyson}.

}

\section*{Acknowledgment}
We thank 
Semyon Alesker,
Valerii Berestovskii,
I. David Berg,
Richard Bishop, 
Yuriy Burago, 
Alexander Christie,
Nicola Gigli,
Sergei Ivanov,
Bernd Kirchheim, 
Bruce Kleiner, 
Rostislav Matveyev,
John Lott, 
Greg Kuperberg, 
Nikolai Kosovsky, 
Nina Lebedeva,
Wilderich Tuschmann, 
and
Sergio Zamora Barrera.
A special thanks to Alexander Lytchak, who contributed deeply to this book during the long process of writing, but refused to be one of us.

We thank the mathematical institutions where we worked on this book: 
BIRS, 
MFO,
Henri Poincar\'{e} Institute,
IHES,
University of Cologne, 
Max Planck Institute for Mathematics,
EIMI.


During the long writing of this book, we were partially supported by the following grants:
Stephanie Alexander --- 
Simons Foundation grant 209053;
Vitali Kapovitch ---  NSF grant DMS-0204187, NSERC Discovery grants, and Simons Foundation grant 390117;
Anton Petrunin --- 
NSF grants
DMS-0406482,
DMS-0905138,
DMS-1309340,
DMS-2005279,
Simons Foundation grants 
245094 and 584781,
and
Minobrnauki, grant 075-15-2022-289.

\tableofcontents

\part{Preliminaries}
\chapter{Model plane}

\section{Trigonometry}\label{model}

Given a real number $\kappa$, the \index{plane!$\kappa$-plane}\emph{model $\kappa$-plane} will be a complete simply connected 2-dimensional Riemannian manifold of constant curvature~$\kappa$.

The model $\kappa$-plane will be denoted by $\Lob2\kappa$\index{$\Lob2\kappa$}.
\begin{itemize}
\item If $\kappa>0$, $\Lob2\kappa$ is isometric to a sphere of radius $\tfrac{1}{\sqrt{\kappa}}$; the unit sphere $\Lob21$ will be also denoted by $\mathbb{S}^2$.
\item If $\kappa=0$, $\Lob2\kappa$ is the Euclidean plane, which is also denoted by $\EE^2$. 
\item If $\kappa<0$, $\Lob2\kappa$ is the Lobachevsky plane with curvature $\kappa$.
\end{itemize}

Let \index{$\varpi\kappa$}$\varpi\kappa=\diam\Lob2\kappa$, so 
$\varpi\kappa=\infty$ if $\kappa\le0$ and $\varpi\kappa=\pi/\sqrt{\kappa}$ if $\kappa>0$;
$\varpi{}$ is just a cursive form of $\pi$.

The distance between points $x,y\in \Lob2\kappa$ will be denoted by $\dist{x}{y}{}$\index{$\dist{p}{q}{}$}, and $[x y]$\index{$[{p}{q}]$} 
will denote the geodesic segment connecting $x$ and $y$. 
The segment $[x y]$ is uniquely defined for $\kappa\le 0$ and for $\kappa>0$ it is defined uniquely if $\dist{x}{y}{}<\varpi\kappa=\pi/\sqrt{\kappa}$.

A triangle in $\Lob2\kappa$ with vertices $x,y,z$ will be denoted by $\trig x y z$\index{$\trig{p}{q}{r}$}.
Formally, a triangle is an ordered set of its sides, so $\trig x y z$ is just a short notation for the triple $([y z],[z x],[x y])$.

The angle of $\trig x y z$ at $x$ will be denoted by $\mangle\hinge xyz$\index{$\mangle$}.

By $\modtrig\kappa\{a,b,c\}$\index{$\modtrig\kappa$!$\modtrig\kappa\{{a},{b},{c}\}$} we denote a triangle in 
$\Lob2\kappa$ with side lengths $a,b,c$, so 
$\trig x y z=\modtrig\kappa\{a,b,c\}$ means that $x,y,z\in \Lob2\kappa$ are such that 
\[\dist{x}{y}{}=c,\quad \dist{y}{z}{}=a,\quad \dist{z}{x}{}=b.\]
For $\modtrig\kappa\{a,b,c\}$ to be defined, the sides $a,b,c$ must satisfy the triangle inequality. If $\kappa>0$, we 
require 
in addition that $a+b+c<2\cdot\varpi\kappa$; 
otherwise $\modtrig\kappa\{a,b,c\}$ is considered to be undefined.

\parbf{Trigonometric functions.}
We will need three \textit{trigonometric functions} in $\Lob2\kappa$: {}\emph{cosine}, {}\emph{sine}, and {}\emph{modified distance}, denoted by $\cs\kappa$, $\sn\kappa$, and $\md\kappa$  respectively. 

They are defined as the solutions of the following initial value problems respectively:
\[
\begin{cases}
 x''+\kappa\cdot x=0,\\
 x(0)=1,\\
 x'(0)=0.
 \end{cases} 
 \quad 
 \begin{cases}
 y''+\kappa\cdot y=0,\\
 y(0)=0,\\
 y'(0)=1.
 \end{cases} 
\quad
 \begin{cases}
 z''+\kappa\cdot z=1,\\
 z(0)=0,\\
 z'(0)=0.
 \end{cases} 
\]

Namely, we set $\cs\kappa(t)=x(t)$, $\sn\kappa(t)=y(t)$, and 
\[
\md\kappa(t)=
\begin{cases}
z(t)& \textrm{ if } 0\le t\le \varpi\kappa,
\\
\tfrac{2}{\kappa}& \textrm{ if } t> \varpi\kappa.
\end{cases}
\]

Here are the tables which relate our trigonometric functions to the standard ones, where 
we take $\kappa>0$:\index{$\md\kappa$}\index{$\sn\kappa$}\index{$\cs\kappa$}
{\small
\begin{align*}
&\sn{\pm\kappa}=\tfrac{1}{\sqrt{\kappa}}\cdot \sn{\pm1}({x}\cdot {\sqrt{\kappa}});
&
&\cs{\pm\kappa}=\cs{\pm1}({x}\cdot {\sqrt{\kappa}});
&
&\md{\pm\kappa}
=
\tfrac{1}{\kappa}\cdot \md{\pm1}({x}\cdot {\sqrt{\kappa}});
\\
&\sn{-1} x=\sinh x;
&
&\cs{-1} x=\cosh x;
&
&\md{-1} x=\cosh x-1;
\\
&\sn{0} x=x;
&
&\cs{0} x=1;
&
&\md{0} x=\tfrac{1}{2}\cdot x^2; 
\\
&\sn{1} x=\sin x;
&
&\cs{1} x= \cos x;
&
&\md{1} x
=
\left\{\begin{smallmatrix}
1-\cos x&\text{for}\ x\le \pi,
\\
2&\text{for}\ x >\pi.
\end{smallmatrix}
\right.
\end{align*}
}

{

\begin{wrapfigure}{r}{44mm}
\centering
\includegraphics{mppics/pic-105}
\end{wrapfigure}

Note that
\[
\mdk(x)=\int\limits_0^x
\snk(t)\cdot\dd t \quad \text{ for }\quad x\le \varpi\kappa.
\]

Let $\phi$ be the angle of $\modtrig\kappa\{a,b,c\}$ 
opposite to $a$.
In this case, we will write \label{page:model-side}%
\index{$\side\kappa$!$\side\kappa \{{a};{b},{c}\}$}%
\index{$\tangle\mc\kappa$!$\tangle\mc\kappa\{a;b,c\}$}
\[a
=
\side\kappa\{\phi;b,c\}
\quad \text{or}\quad 
\phi
=
\tangle\mc\kappa\{a;b,c\}.\]

}

The functions $\side\kappa$ and $\tangle\mc\kappa$ will be called respectively the \index{model side}\emph{model side} and the \index{model angle}\emph{model angle}.
Let 
\[
\side\kappa\{\phi;b,-c\}
=\side\kappa\{\phi;-b,c\}
\df\side\kappa\{\pi-\phi;b,c\};\]
in this way we define $\side\kappa\{\phi;b,c\}$ when one of the numbers $b$ and $c$ is negative. 

\pagebreak

\begin{thm}{Properties of standard functions}\label{md-equalities}

\begin{subthm}{md-diff-eq}
For fixed $a$ and $\phi$, the function $y(t)=\md\kappa\left(\side\kappa\{\phi;a,t\}\right)$
 satisfies the following differential equation:
\[y''+\kappa\cdot y=1.\]
\end{subthm}

\begin{subthm}{sn-diff-eq}
Let $\alpha\:[a,b]\to\Lob2\kappa$ be a unit-speed geodesic, and $A$ be the image of a complete geodesic. If $f(t)$ is the distance from $\alpha(t)$ to $A$, the function 
$y(t)=\sn\kappa (f(t))$
 satisfies the following differential equation:
\[y''+\kappa\cdot y=0\]
for $y\ne 0$.
\end{subthm}

\begin{subthm}{increase}
For fixed $\kappa$, $b$, and $c$, the function 
\[a\mapsto\tangle\mc\kappa\{a;b,c\}\]
is increasing and defined on a real interval.
Equivalently, the function
\[\phi\mapsto\side\kappa \{\phi;b,c\}\]
is increasing and defined if $b,c<\varpi\kappa$, and $\phi\in[0,\pi]$.
(Formally speaking, if $\kappa>0$ and $b+c\ge \varpi\kappa$, it is defined only for $\phi\in[0,\pi)$, but $\side\kappa \{\phi;b,c\}$ can be extended to $[0,\pi]$ as a continuous function.)
\end{subthm}

\begin{subthm}{k-decrease}
For fixed $\phi,a,b,c$, the function
\[\kappa\mapsto \tangle\mc\kappa\{a;b,c\}\quad \text{and}\quad \kappa\mapsto \side\kappa \{\phi;b,c\}\]
are respectively nondereasing (in fact, increasing, if $|b\z-c|\z<a\z< b+c$)
and nonincreasing (in fact, increasing, if $0\z<\phi\z<\pi$).
\end{subthm}

\begin{subthm}{lem:alex-0}(Alexandrov's lemma)
Assume that for real numbers $a$, $b$, $a'$, $b'$, $x$, and $\kappa$, the following two expressions are defined:
\[\tangle\mc\kappa\{a;b,x\}+\tangle\mc\kappa\{a';b',x\}-\pi,\quad
\tangle\mc\kappa\{a';b+b',a\}-\tangle\mc\kappa\{x;a,b\},\]
Then they have the same sign.
\end{subthm}
\end{thm}

All the properties except Alexandrov's lemma \ref{SHORT.lem:alex-0} can be shown by direct calculation. Alexandrov's lemma is reformulated in \ref{lem:alex} and is proved there.

\parbf{Cosine law.}
The formulas $a
=
\side\kappa\{\phi;b,c\}$ and
$\phi
=
\tangle\mc\kappa\{a;b,c\}$
can be rewritten using the cosine law in $\Lob2\kappa$:
\[\cos\phi
=\left\{\begin{aligned}
&\frac{b^2+c^2-a^2}{2\cdot b\cdot c}
&\text{if}&\quad \kappa=0,
\\
&\frac{\cs\kappa a-\cs\kappa b\cdot\cs\kappa c}{\kappa\cdot\sn\kappa b\cdot\sn\kappa c}
&\text{if}&\quad \kappa\ne0.
\end{aligned}\right.\]

However, rather than using these explicit formulas, we mainly will use
the properties of $\tangle\mc\kappa$ and $\side\kappa$ listed in \ref{md-equalities}.

\section{Hemisphere lemma}\label{curves-in-model}

\begin{thm}{Hemisphere lemma}
\label{lem:hemisphere}
For $\kappa>0$, any closed path of length less than $2\cdot \varpi\kappa$ (respectively, at most $2\cdot \varpi\kappa$) in $\Lob2\kappa$ lies in an open (respectively, closed) hemisphere. 
\end{thm}

\parit{Proof.}
Applying rescaling, we may assume that $\kappa=1$, and thus $\varpi\kappa=\pi$ and $\Lob2{\kappa}=\mathbb{S}^2$.
Let $\alpha$ be a closed curve in $\mathbb{S}^2$ of length $2\cdot\ell$.

\begin{wrapfigure}{o}{35 mm}
\vskip-0mm
\centering
\includegraphics{mppics/pic-110}
\end{wrapfigure}

Assume $\ell<\pi$.
Let $\check\alpha$ be a subarc of $\alpha$ of length $\ell$, with endpoints $p$ and $q$. 
Since $\dist{p}{q}{}\le\ell<\pi$, there is a unique geodesic $[pq]$ in $\mathbb{S}^2$. 
Let $z$ be the midpoint of $[pq]$. 
We claim that $\alpha$ lies in the open hemisphere centered at $z$. 
If not, $\alpha$ intersects the boundary great circle of this hemisphere; let $r$ be a point in the intersection.
Without loss of generality, we may assume that $r\in\check\alpha$. 
The arc $\check\alpha$ together with its reflection in $z$ form a closed curve of length $2\cdot \ell$ that contains $r$ and its antipodal point~$r'$.
Thus 
\[\ell=\length \check\alpha\ge \dist{r}{r'}{}=\pi,\] 
a contradiction.

If $\ell=\pi$, then either $\alpha$ is a local geodesic, and hence a great circle, 
or $\alpha$ may be strictly shortened by substituting a geodesic arc for a subarc of $\alpha$ 
whose endpoints $p^1,p^2$ are arbitrarily close to a point $p$ on $\alpha$.
In both cases $\alpha$ lies in a closed hemisphere;
the former case is trivial, and in the latter case, $\alpha$ lies in a closed hemisphere obtained as a limit of closures of open hemispheres containing the shortened curves as $p^1,p^2$ approach~$p$.
\qeds


\begin{thm}{Exercise}\label{exr-crofton}
Give a proof of the hemisphere lemma (\ref{lem:hemisphere}) based on Crofton's formula.
\end{thm}

\chapter{Metric spaces}
In this chapter we fix conventions and notations. We are assuming that the reader is familiar with basic notions in metric geometry.

\section{Metrics and their relatives}
\label{sec:metric spaces}

\parbf{Definitions.}
Let $\II$ be a subinterval of $[0,\infty]$.
A function $\rho$ defined on $\spc{X}\times\spc{X}$ is called an $\II$-valued metric if the following conditions hold:
\begin{itemize}
 \item $\rho(x,x)=0$ for any $x$;
 \item $\rho(x,y)\in \II$ for any pair $x\ne y$;
 \item $\rho(x,y)+\rho(x,z)\ge \rho(y,z)$ for any triple of points $x,y, z$.
\end{itemize}
The value $\rho(x,y)$ is also called the \index{distance}\emph{distance} between $x$ and $y$.

The above definition will be used for four choices of interval $\II$: $(0,\infty)$, $(0,\infty]$, $[0,\infty)$, and $[0,\infty]$.
Any $\II$-valued metric can be referred to briefly as a metric;
the interval should be apparent from context but by default, a metric is $(0,\infty)$-valued. 
If we need to be more specific we may also use the following names:
\begin{itemize}
\item a $(0,\infty)$-valued metric may be called a \index{genuine metric}\emph{genuine metric}.
\item a $(0,\infty]$-valued metric may be called an \index{$\infty$-metric}\emph{$\infty$-metric}.
\item a $[0,\infty)$-valued metric may be called  a \index{genuine pseudometric}\emph{genuine pseudometric}.
\item A $[0,\infty]$-valued metric may be called a \index{pseudometric}\emph{pseudometric} or \index{$\infty$-pseudometric}\emph{$\infty$-pseudometric}.
\end{itemize}

A metric space is a set equipped with a metric.
The distance between points $x$ and $y$ in a metric space $\spc{X}$ will  usually be denoted by \[\dist{x}{y}{}\qquad\text{or}\quad \dist{x}{y}{\spc{X}};\]
the latter will be used if we need to emphasize that we are working in the space $\spc{X}$.

The function $\distfun{x}{}{}\:\spc{X}\to\RR$ defined as 
\[\distfun{x}{}{}\:y\mapsto\dist{x}{y}{}\]
will be called the \index{distance function}\emph{distance function} from $x$. 

Any  subset $A$ in a metric space $\spc{X}$ will be also considered as a \index{subspace}\emph{subspace};
that is, a metric space with the metric defined by restricting the metric of $\spc{X}$ to $A\times A\subset \spc{X}\times\spc{X}$.

The \index{direct product}\emph{direct product} $\spc{X}\times \spc{Y}$ of two metric spaces $\spc{X}$ and $\spc{Y}$ is defined as the metric space carrying  the metric
\[
\dist{(p,\phi)}{(q,\psi)}{} =\sqrt{\dist[2]{p}{q}{} + \dist[2]{\phi}{\psi}{}}
\]
for $p,q\in \spc{X}$ and $\phi,\psi\in \spc{Y}$. 

A map between two metric spaces is called an \index{isometry}\emph{isometry} if it is a bijection and preserves distances between points.

\parbf{Zero and infinity.}
Genuine metric spaces are the main objects of study in this book.
However, the  generalizations above are useful
in  various definitions and constructions.
For example, the construction of length metric (see Section~\ref{sec:intrinsic}) uses infinite distances.
The following definition gives another example.

\begin{thm}{Definition}\label{def:disjoint-union}
Assume $\{\spc{X}_\alpha\}_{\alpha\in\IndexSet}$ is a collection of $\infty$-metric spaces.
The disjoint union 
$$\bm{X}=\bigsqcup_{\alpha\in\IndexSet}\spc{X}_\alpha$$ 
has a natural $\infty$-metric on it defined as follows:
given two points $x\in\spc{X}_\alpha$ and $y\in\spc{X}_\beta$,
let 
\[
\begin{matrix}
\dist{x}{y}{\bm{X}}=\infty&\text{if}\quad\alpha\ne\beta,
\\
\dist{x}{y}{\bm{X}}=\dist{x}{y}{\spc{X}_\alpha}&\text{if}\quad\alpha=\beta.
\end{matrix}
\]
The resulting $\infty$-metric space $\bm{X}$ will be called the \index{disjoint union of metric spaces}\emph{disjoint union} of $\{\spc{X}_\alpha\}_{\alpha\in\IndexSet}$, denoted by \[\bigsqcup_{\alpha\in\IndexSet}\spc{X}_\alpha.\]
\end{thm}

Now let us give examples showing that vanishing and infinite distance between distinct points is useful.

Suppose a set ${\spc{X}}$ comes with a set of metrics $\dist{}{}{\alpha}$ for $\alpha\in\IndexSet$.
Then 
\[\dist{x}{y}{}=\sup\set{\dist{x}{y}{\alpha}}{\alpha\in\IndexSet}\]
is in general only an $\infty$-metric;
that is, even if the metrics $\dist{}{}{\alpha}$ are genuine, then $\dist{}{}{}$ might be $(0,\infty]$-valued.

Let $\spc{X}$ be a set,
$\spc{Y}$ be a metric space, and  
 $\map\:\spc{X}\to\spc{Y}$ be a map.
If $\map$ is not injective,
then the {}\emph{pullback}
\[\dist{x}{y}{\spc{X}}=\dist{\map(x)}{\map(y)}{\spc{Y}}\]
defines only a pseudometric on $\spc{X}$.

\parbf{Corresponding metric space and metric component.}
The following two observations show that
nearly any question about metric spaces can be reduced to a question about genuine metric spaces.

Assume $\spc{X}$ is a pseudometric space.
Set
$x\sim y$ if $\dist{x}{y}{}=0$. 
Note that if $x\sim x'$, then $\dist{y}{x}{}=\dist{y}{x'}{}$ for any $y\in\spc{X}$.
Thus, $\dist{}{}{}$ defines a metric on the
quotient set $\spc{X}/{\sim}$.
This way we obtain a metric space $\spc{X}'$.
The space $\spc{X}'$ is called the 
\emph{corresponding metric space} for the pseudometric space $\spc{X}$.
Often we do not distinguish between $\spc{X}'$ and~$\spc{X}$.

Set $x\approx y$ if and only if $\dist{x}{y}{}<\infty$;
this is another equivalence relation on $\spc{X}$.
The equivalence class of a point $x\in\spc{X}$ will be called the \index{metric component}\emph{metric component} 
 of $x$; it will be denoted by $\spc{X}_x$\index{ $\spc{X}_x$}.
One could think of $\spc{X}_x$ as  $\oBall(x,\infty)_{\spc{X}}$, the open ball centered at $x$ and radius $\infty$ in $\spc{X}$; see definition below.

It follows that any $\infty$-metric space is a {}\emph{disjoint union} of genuine metric spaces, the metric components of the original $\infty$-metric space; see Definition~\ref{def:disjoint-union}

To summarize this discussion: Given a $[0,\infty]$-valued metric space $\spc{X}$, we may pass to the corresponding $(0,\infty]$-valued metric space $\spc{X}'$ and break the latter  into a disjoint union of metric components, each of which is  a genuine metric space.

\section{Notations}
\label{sec:notations}

\parbf{Balls.}
Given $R\in[0,\infty]$ and a point $x$ in a metric space $\spc{X}$, the sets
\begin{align*}
\oBall(x,R)&=\set{y\in \spc{X}}{\dist{x}{y}{}<R},
\\
\cBall[x,R]&=\set{y\in \spc{X}}{ \dist{x}{y}{}\le R}
\end{align*}
are called respectively the \index{open ball}\emph{open} and the \index{closed ball}\emph{closed  balls} of radius $R$ with center $x$.

If we need to emphasize that these balls are taken in the space $\spc{X}$,
we write\index{$\oBall(x,R)_{\spc{X}}$}\index{$\cBall[x,R]_{\spc{X}}$}
\[\oBall(x,R)_{\spc{X}}\quad\text{and}\quad \cBall[x,R]_{\spc{X}}\]
respectively.

Since in the model space $\Lob{m}{\kappa}$ all balls of the same radius are isometric, 
often we will not need to specify the center of the ball,
and may write 
\[\oBall(R)_{\Lob{m}{\kappa}}\quad\text{and}\quad\cBall[R]_{\Lob{m}{\kappa}}\] respectively.

A set $A\subset\spc{X}$ is called \index{bounded set}\emph{bounded} if $A\subset\oBall(x,R)$ for some $x\in\spc{X}$ and $R<\infty$.

\parbf{Distances to sets.}
For subset $A\subset \spc{X}$, 
let us denote the distance from $A$ to a point $x$ in $\spc{X}$ by $\distfun{A}{x}$\index{$\distfun{A}{x}$};
that is,
\[\distfun{A}{x}
\df
\inf\set{\dist{a}{x}{}}{a\in A}.\] 

For any subset $A\subset\spc{X}$, the sets \index{$\oBall(A,R)$}\index{$\cBall[A,R]$}
\begin{align*}
\oBall(A,R)&=\set{y\in \spc{X}}{ \distfun{A}{y}{}<R},
\\
\cBall[A,R]&=\set{y\in \spc{X}}{\distfun{A}{y}{}\le R}
\end{align*}
are called respectively the  \index{open tubular neighborhood}\emph{open} and \index{closed tubular neighborhood}\emph{closed $R$-neighborhoods} of $A$.

\parbf{Diameter, radius, and packing.}
Let $\spc{X}$ be a metric space.
Then the \index{diameter}\emph{diameter} of $\spc{X}$ is  defined as
\[\diam \spc{X}=\sup\set{\dist{x}{y}{}}{x,y\in \spc{X}}.\]

The \index{radius}\emph{radius} of $\spc{X}$ is  defined as
\[\rad \spc{X}=\inf\set{R>0}{\oBall(x,R)= \spc{X}\ \text{for some}\ x\in\spc{X}}.\]

The \index{$\eps$-pack}\emph{$\eps$-pack} of $\spc{X}$ (or \index{packing number}\emph{packing number}) is the maximal number  (possibly infinite) of points in $\spc{X}$ at distance $>\eps$ from each other;  it is denoted by $\pack_\eps\spc{X}$\index{$\pack_\eps\spc{X}$}.
If $m=\pack_\eps\spc{X}<\infty$, then a set $\{x^1,x^2,\dots,x^m\}$ in $\spc{X}$ 
such that $\dist{x^i}{x^j}{}>\eps$ is called a \index{maximal packing}\emph{maximal $\eps$-packing} in $\spc{X}$.

\parbf{G-delta sets.}
Recall that an arbitrary union of open balls in a metric space is called an \index{open set}\emph{open set}.
A subset of a metric space is called a \index{G-delta set}\emph{G-delta set} if it can be presented as an intersection of  a countable number of open subsets. 

\begin{thm}{Baire's theorem}
Let $\spc{X}$ be a complete metric space 
and $\{\Omega_n\}$, $n\in \NN$, be a collection of open dense subsets of $\spc{X}$.
Then $\bigcap_{n\in \NN}\Omega_n$ is dense in $\spc{X}$.
\end{thm}

\parbf{Proper spaces.}
A metric space $\spc{X}$ is called \index{proper space}\emph{proper} if all closed bounded sets in $\spc{X}$ are compact. 
This condition is equivalent to each of the following statements:
\begin{enumerate}
\item For some (and therefore any) point $p\in \spc{X}$ and any $R<\infty$, 
the closed ball $\cBall[p,R]\subset\spc{X}$ is compact. 
\item The function $\distfun{p}{}{}\:\spc{X}\to\RR$ is proper for some (and therefore any) point $p\in \spc{X}$.
\end{enumerate}

We will also often use the following two classical statements:

\begin{thm}{Proposition}
Proper metric spaces are separable and second countable.
\end{thm}

\begin{thm}{Proposition}\label{compact=seq-compact}
Let $\spc{X}$ be a metric space.
Then the following are equivalent

\begin{subthm}{}
$\spc{X}$ is compact;
\end{subthm}

\begin{subthm}{}
$\spc{X}$ is \index{sequentially compact}\emph{sequentially compact}; that is, any sequence of points in $\spc{X}$ contains a convergent subsequence;
\end{subthm}

\begin{subthm}{}
$\spc{X}$ is complete and for any $\eps>0$ there is a \index{$\eps$-net}\emph{finite $\eps$-net} in $\spc{X}$; that is, there is a finite collection of points $p_1,\ldots,p_{N}$ such that $\bigcup_i\oBall(p_i,\eps)=\spc{X}$.
\end{subthm}

\begin{subthm}{}
$\spc{X}$ is complete and for any $\eps>0$ there is a {}\emph{compact $\eps$-net} in $\spc{X}$; that is, $\oBall(K,\eps)=\spc{X}$ for a compact set $K\subset \spc{X}$.
\end{subthm}

\end{thm}

\section{Length spaces}\label{sec:intrinsic}

A \index{curve}\emph{curve} in a metric space $\spc{X}$ is a continuous map $\alpha\:\II\to \spc{X}$, where $\II$ is a {}\emph{real interval} (that is, an arbitrary convex subset of $\RR$).

\begin{thm}{Definition}\label{def:length}
Let $\spc{X}$ be a metric space.
Given a curve $\alpha\: \II\to \spc{X}$, we define its \index{length}\emph{length} as 
\[
\length \alpha \df \sup\set{\sum_{i\ge 1} \dist{\alpha(t_i)}{\alpha (t_{i-1})}{}}{t_0,\dots,t_n\in \II,t_0\le\ldots\le t_n}.
\]
\end{thm}

The following lemma is an easy exercise.

\begin{thm}{Lower semicontinuity of length}\label{thm:semicont-of-length}
Assume $\alpha_n\:\II\to \spc{X}$ is a sequence of curves that converges pointwise to a curve $\alpha_\infty\:\II\to \spc{X}$.
Then 
\[\length\alpha_\infty\le \liminf_{n\to\infty}\length\alpha_n.\]

\end{thm}

Given two points $x$ and $y$ in a metric space $\spc{X}$,
consider the value
\[\yetdist{x}{y}{}=\inf_{\alpha}\{\length\alpha\},\]
where infimum is taken for all paths $\alpha$ from $x$ to $y$.

It is easy to see that $\yetdist{}{}{}$ defines a $(0,\infty]$-valued metric on  $\spc{X}$;
it will be called the \index{induced length metric}\emph{induced length metric} on $\spc{X}$.
Clearly 
\[\yetdist{x}{y}{}\ge \dist{x}{y}{}\]
for any $x,y\in \spc{X}$.

It  easily follows from the definition that the length of a curve $\alpha$ with respect to $\yetdist{}{}{}$ is equal to the length of $\alpha$ with respect to $\dist{}{}{}$.
In particular, iterating the construction produces the same metric $\yetdist{}{}{}$.

\begin{thm}{Definition}\label{def:length-space}
If $\yetdist{x}{y}{}=\dist{x}{y}{}$ for any pair of points $x,y$ in a metric space $\spc{X}$, then $\dist{}{}{}$ is called \index{length metric}\emph{length metric}, and $\spc{X}$ is called a \index{length space}\emph{length space}.
\end{thm}

In other words, a metric space $\spc{X}$ is a
\emph{length space}
if for any $\eps>0$ and any two points $x,y\in \spc{X}$ with $\dist{x}{y}{}<\infty$ there is a path $\alpha\:[0,1]\to\spc{X}$ connecting $x$ to $y$ such that 
\[\length\alpha<\dist{x}{y}{}+\eps.\]

In this book, most of the time we consider length spaces.
If $\spc{X}$ is a length space, 
and $A\subset \spc{X}$, then the set $A$ comes with the inherited metric from $\spc{X}$, 
which might be not a length metric.
The corresponding length metric on $A$ will be denoted by $\dist{}{}{A}$.

\parbf{Variations of the definition.}
We will need the following variations of Definition~\ref{def:length-space}:
\index{R-geodesic space}
\index{R-length space}

\begin{itemize}
\item Assume $R>0$.
If $\yetdist{x}{y}{}=\dist{x}{y}{}$ for any pair $\dist{x}{y}{}<R$, then $\spc{X}$ is called an \index{$R$-length space}\emph{$R$-length space}.
\item If any point in $\spc{X}$ admits a neighborhood  $\Omega$ such that $\yetdist{x}{y}{}=\dist{x}{y}{}$ for any pair of points $x,y\in \Omega$
then  $\spc{X}$ is called a \index{length space!locally length space}\emph{locally length space}.
\item A metric space is called 
\index{geodesic space}\emph{geodesic}
if for any two points $x,y$ with $\dist{x}{y}{}<\infty$ there is a geodesic $[x y]$ in $\spc X$
\item Assume $R>0$. A metric space is called 
\index{$R$-geodesic space}\emph{$R$-geodesic}
if for any two points $x,y$ such that $\dist{x}{y}{}<R$ there is a geodesic $[x y]$ in $\spc X$.
\end{itemize}

Note that the notions of $\infty$-length spaces and length spaces are the same.
Clearly, any geodesic space is a length space 
and any $R$-geodesic space is $R$-length.

\begin{thm}{Example} 
Consider a metric space $\spc{X}$ obtained by gluing a countable collection of disjoint intervals $\II_n$ of length $1+\tfrac1n$ where for each $\II_n$ one end is glued to $p$ and the other to $q$.
Then $\spc{X}$ carries a natural complete length metric such that $\dist{p}{q}{}=1$, but there is no geodesic connecting $p$ to~$q$.
\end{thm}

\begin{thm}{Exercise}\label{ex:complete=>complete}
Let $\spc{X}$ be a metric space
and $\yetdist{}{}{}$ be the length metric on it.
Show the following:
\begin{subthm}{ex:complete=>complete:complete}
If $\spc{X}$  is complete, then  $(\spc{X},\yetdist{}{}{})$ is complete.
\end{subthm}

\begin{subthm}{ex:complete=>complete:compact}
If $\spc{X}$ is compact, then $(\spc{X},\yetdist{}{}{})$ is geodesic.
\end{subthm}
\end{thm}

\begin{thm}{Exercise}\label{ex:no-geod}
Give an example of a complete length space such that no pair of distinct points can be joined by a geodesic.
\end{thm}

\begin{thm}{Exercise}\label{ex:compact+connceted}
Let $\spc{X}$ be a complete length space.
Show that for any compact subset $K$ in $\spc{X}$
there is a compact path-connected subset $K'$ that contains $K$.  
\end{thm}

\begin{thm}{Definition}
Consider two points $x$ and $y$ in a metric space $\spc{X}$.

\begin{enumerate}[(i)]
\item A point $z\in \spc{X}$ is called a \index{midpoint}\emph{midpoint} between $x$ and $y$ if 
\[\dist{x}{z}{}=\dist{y}{z}{}=\tfrac12\cdot\dist[{{}}]{x}{y}{}.\]
\item Assume $\eps\ge 0$.
A point $z\in \spc{X}$ is called  an \index{$\eps$-midpoint}\emph{$\eps$-midpoint} between $x$ and $y$
if 
\[\dist{x}{z}{}\le\tfrac12\cdot\dist[{{}}]{x}{y}{}+\eps,
\quad\text{and}\quad
\dist{y}{z}{}\le\tfrac12\cdot\dist[{{}}]{x}{y}{}+\eps.\]
\end{enumerate}

\end{thm}

Note that a $0$-midpoint is the same as a midpoint.

The following lemma was essentially proved by Karl Menger \cite[Section 6]{menger}.

\begin{thm}{Lemma}\label{lem:mid>geod}
Let $\spc{X}$ be a complete metric space.

\begin{subthm}{lem:mid>length}
Assume that for any pair of points $x,y\in \spc{X}$, and any $\eps>0$
there is an $\eps$-midpoint $z$.
Then  $\spc{X}$ is a length space.
\end{subthm}

\begin{subthm}{lem:mid>geod:geod}
Assume that for any pair of points $x,y\in \spc{X}$ 
there is a midpoint $z$.
Then  $\spc{X}$ is a geodesic space.
\end{subthm}

\begin{subthm}{lem:mid>geod:R}
If for some $R>0$, the assumptions \ref{SHORT.lem:mid>length} or \ref{SHORT.lem:mid>geod:geod} hold only for pairs of points $x,y\in \spc{X}$ such that $\dist{x}{y}{}<R$, 
then  $\spc{X}$ is an $R$-length or an $R$-geodesic space respectively.

\end{subthm}

\end{thm}

\parit{Proof.}
Fix a pair of points $x,y\in \spc{X}$.
Let $\eps_n=\frac\eps{2^{2\cdot n}}$.

Set $\alpha(0)=x$,
and $\alpha(1)=y$.
Let $\alpha(\tfrac12)$ be an $\eps_1$-midpoint between $\alpha(0)$ and $\alpha(1)$.
Further, let $\alpha(\frac14)$ 
and $\alpha(\frac34)$  be $\eps_2$-midpoints 
between the pairs $(\alpha(0),\alpha(\tfrac12)$ 
and $(\alpha(\tfrac12),\alpha(1)$ respectively.
Applying the above procedure recursively,
on the $n$-th step we define $\alpha(\tfrac{\kay}{2^n})$,
for every odd integer $\kay$ such that $0<\tfrac\kay{2^n}<1$, to be an $\eps_{n}$-midpoint between the already defined
$\alpha(\tfrac{\kay-1}{2^n})$ and $\alpha(\tfrac{\kay+1}{2^n})$.

This way we define $\alpha(t)$ for all dyadic rationals $t$ in $[0,1]$.
If $t\in[0,1]$ is not a dyadic rational, consider a sequence of dyadic rationals $t_n\to t$ as $n\to\infty$.
By completeness of $\spc{X}$, the sequence $\alpha(t_n)$ converges;
let $\alpha(t)$ be its limit.
It is easy to see that $\alpha(t)$
does not depend on the choice of the sequence $t_n$,
and $\alpha\:[0,1]\to\spc{X}$ is a path from $x$ to $y$.
Moreover,
\[\begin{aligned}
\length\alpha&\le \dist{x}{y}{}+\sum_{n=1}^\infty 2^{n-1}\cdot\eps_n\le
\\
&\le \dist{x}{y}{}+\tfrac\eps2.
\end{aligned}
\eqlbl{eq:eps-midpoint}
\]
Since $\eps>0$ is arbitrary, we have \ref{SHORT.lem:mid>length}.

To prove \ref{SHORT.lem:mid>geod:geod}, 
one should repeat the same argument 
taking midpoints instead of $\eps_n$-midpoints.
In this case, \ref{eq:eps-midpoint} holds for $\eps_n=\eps=0$.

The proof of \ref{SHORT.lem:mid>geod:R} is obtained by a straightforward modification of the proofs above.
\qeds

Since in a compact set a sequence of $\tfrac1n$-midpoints $(z_n)$ contains a convergent subsequence,
Lemma \ref{lem:mid>geod} implies the following.

\begin{thm}{Corollary}
A proper length space is geodesic.
\end{thm}

{\sloppy

\begin{thm}{Hopf--Rinow theorem}\label{thm:Hopf-Rinow}
Any complete, locally compact length space is proper.
\end{thm}

}

\parit{Proof.}
Let $\spc{X}$ be a locally compact length space.
Given $x\in \spc{X}$, denote by $\rho(x)$ the supremum of all $R>0$ such that
the closed ball $\cBall[x,R]$ is compact.
Since $\spc{X}$ is locally compact,
$$\rho(x)>0\quad\text{for any}\quad x\in \spc{X}.\eqlbl{eq:rho>0}$$
It is sufficient to show that $\rho(x)=\infty$ for some (and therefore any) point $x\in \spc{X}$.

Assume the contrary; that is, $\rho(x)<\infty$.

\begin{clm}{} $B=\cBall[x,\rho(x)]$ is compact for any $x$.
\end{clm}

Indeed, $\spc{X}$ is a length space;
therefore for any $\eps>0$, 
the set $\cBall[x,\rho(x)\z-\eps]$ is a compact $\eps$-net in $B$.
Since $B$ is closed and hence complete, it is compact by Proposition~\ref{compact=seq-compact}.
\claimqeds

\begin{clm}{} $|\rho(x)-\rho(y)|\le \dist{x}{y}{\spc{X}}$ for any $x,y\in\spc{X}$;
in particular $\rho\:\spc{X}\to\RR$ is a continuous function.
\end{clm}

Indeed, assume the contrary; that is, $\rho(x)+\dist{x}{y}{}<\rho(y)$ for $x,y\in \spc{X}$. 
Then 
$\cBall[x,\rho(x)+\eps]$ is a closed subset of $\cBall[y,\rho(y)]$ for $\eps>0$.
Then  compactness of $\cBall[y,\rho(y)]$ implies compactness of $\cBall[x,\rho(x)+\eps]$, a contradiction.\claimqeds

Set $\eps=\min_{y\in B}\{\rho(y)\}$; 
the minimum is defined since $B$ is compact.
From \ref{eq:rho>0}, we have $\eps>0$.

Choose a finite $\tfrac\eps{10}$-net $\{a_1,a_2,\dots,a_n\}$ in $B$.
The union $W$ of the closed balls $\cBall[a_i,\eps]$ is compact.
Clearly 
$\cBall[x,\rho(x)+\frac\eps{10}]\subset W$.
Therefore $\cBall[x,\rho(x)+\frac\eps{10}]$ is compact,
a contradiction.
\qeds

\begin{thm}{Exercise}\label{exercise from BH}
Construct a geodesic space that is locally compact,
but whose completion is neither geodesic nor locally compact.
\end{thm}

\section{Convex sets}

\begin{thm}{Definition} 
\label{def:convex-set}
Let $\spc{X}$ be a geodesic space and $A\subset\spc{X}$.

$A$ is 
\index{convex set}\emph{convex}
if for every two points $p,q\in A$ any geodesic $[pq]$ lies in~$A$.

$A$ is 
\index{convex set!weakly convex set}\emph{weakly convex}
if for every two points $p,q\in A$
there is a geodesic $[pq]$ 
that lies in $A$.

We say  that $A$ is 
\index{convex set!totally convex set}\emph{totally convex}
if for every two points $p,q\in A$, every local geodesic from $p$ to $q$ lies in $A$.

If for some $R\in (0,\infty]$ these definitions are applied only for pairs of points such that $\dist{p}{q}{}<R$ and only for the geodesics of length $<R$,
then $A$ is called respectively 
\index{$R$-convex set}\emph{$R$-convex},
\emph{weakly $R$-convex},
or \emph{totally $R$-convex}.

A set $A\subset\spc{X}$ is called 
\index{convex set!locally convex set}\emph{locally convex}
if every point $a\in A$ admits an open neighborhood $\Omega\ni a$
such that for every two points $p,q\in A\cap\Omega$ every geodesic $[pq]\subset \Omega$ lies in $A$.
Similarly one defines locally weakly convex and locally totally convex sets. 
\end{thm}

\parbf{Remarks.} Let us state a few observations that easily follow  from the definition.
\begin{itemize}
\item The notion of ({}\emph{weakly}) {}\emph{convex set} is the same as ({}\emph{weakly}) {}\emph{$\infty$-convex set}.
\item The inherited metric on a weakly convex set coincides with its length metric.
\item Any open set is locally convex by definition.
\end{itemize}

The following proposition states that weak convexity survives under ultralimits.
An analogous statement about convexity does not hold;
for example, there is a sequence of convex discs in $\mathbb{S}^2$ that converges to a hemisphere, which is not convex.

\begin{thm}{Proposition}\label{prop:weak-convex-stable}
Let $\spc{X}_n,$ be a sequence of geodesic spaces.
Let $\o$ be an ultrafilter on $\mathbb N$ (see Definition~\ref{def:ultrafilter}).
Assume that  $A_n\subset \spc{X}_n$ is a sequence of weakly convex sets, 
$\spc{X}_n\to \spc{X}_\o$,
and $A_n\to A_\o\subset \spc{X}_\o$ as $n\to\o$.
Then $A_\o$ is a weakly convex set of $\spc{X}_\o$.
\end{thm}

\parit{Proof.}
Fix $x_\o,y_\o\in A_\o$.
Consider sequences $x_n,y_n\in A_n$ such that $x_n\to x_\o$ and $y_n\to y_\o$ as $n\to\o$.

Denote by $\alpha_n$ a geodesic path from $x_n$ to $y_n$ that lies in $A_n$.
Let
\[\alpha_\o(t)=\lim_{n\to \o}\alpha_n(t).\]

It remains to observe that $\alpha_\o$ is a geodesic path that lies in $A_\o$.
\qeds

\section{Quotient spaces}\label{sec:quotient}

\parbf{Quotient spaces.} Assume $\spc{X}$ is a metric space
with an equivalence relation $\sim$.
Note that given a family of pseudometrics $\rho_\alpha$ on $\spc{X}/{\sim}$,
their least upper bound
\[\rho(x,y)=\sup_\alpha\{\,\rho_\alpha(x,y)\,\}\]
is also a pseudometric.
If for the projections $\spc{X}\to(\spc{X}/{\sim},\rho_\alpha)$ are \index{short map}\emph{short} (that is, \textit{distance non-increasing}), then so is $\spc{X}\to(\spc{X}/{\sim},\rho)$.

It follows that 
the quotient space $\spc{X}/{\sim}$ admits a natural quotient pseudometric;
this is the maximal pseudometric on  $\spc{X}/{\sim}$ 
that makes the quotient map 
$\spc{X}\to\spc{X}/{\sim}$ short.
The corresponding metric space will be also denoted as $\spc{X}/{\sim}$
and will be called the \index{quotient space}\emph{quotient space} of $\spc{X}$ by the equivalence relation $\sim$.

In general, the points of the metric space $\spc{X}/{\sim}$
are formed by equivalence classes in $\spc{X}$
for a wider equivalence relation.
However, in most of the cases we will consider, 
the set of equivalence classes will coincide with the set of points in the metric space $\spc{X}/{\sim}$.

\begin{thm}{Proposition}\label{prop:length-X}
Let $\spc{X}$ be a length space and 
$\sim$ be an equivalence relation on $\spc{X}$. Then $\spc{X}/{\sim}$
is a length space.
\end{thm}

\parit{Proof.}
Let $\spc{Y}$ be an arbitrary metric space.
Since $\spc{X}$ is a length space,
a map $f\:\spc{X}\to \spc{Y}$ is short if and only if 
\[\length(f\circ\alpha)\le \length\alpha\]
for any curve
$\alpha\:\II\to \spc{X}$.
Denote by $\yetdist{}{}{}$ the length metric on $\spc{Y}$.
It follows that if $f\:\spc{X}\to \spc{Y}$ is short
then so is 
$f\:\spc{X}\to (\spc{Y},\yetdist{}{}{})$.

Consider the quotient map 
$f\:\spc{X}\to\spc{X}/{\sim}$.
Recall that the space $\spc{X}/{\sim}$ is defined by the maximal pseudometric that makes $f$ short.

Denoting by $\yetdist{}{}{}$ the length metric on $\spc{X}/{\sim}$,
it follows that
\[f\:\spc{X}\to(\spc{X}/{\sim},\yetdist{}{}{})\]
is also short.

Note that 
\[\yetdist{x}{y}{}\ge\dist{x}{y}{\spc{X}/{\sim}}\]
for any $x,y\in \spc{X}/{\sim}$.
From maximality of $\dist{}{}{\spc{X}/{\sim}}$, we get
\[\yetdist{x}{y}{}=\dist{x}{y}{\spc{X}/{\sim}}\]
for any $x,y\in \spc{X}/{\sim}$;
that is, $\spc{X}/{\sim}$ is a length space.
\qeds  

\parbf{Group actions.}
Assume a group $G$ acts on a metric space $\spc{X}$.
Consider a relation $\sim$ on $\spc{X}$
defined by $x\sim y$ if there is $g\in G$ such that $x=g\cdot y$.
Note that $\sim$ is an equivalence relation.

In this case, the quotient space $\spc{X}/{\sim}$ will also be denoted by $\spc{X}/G$, and can be regarded as the space of $G$-orbits in $\spc{X}$.

Assume that a group $G$ acts on $\spc{X}$ by isometries.
Then the distance between orbits $G\cdot x$ and $G\cdot y$ in $\spc{X}/G$
can be defined directly: 
\[\dist{G\cdot x}{G\cdot y}{\spc{X}/G}=\inf\set{\dist{x}{g\cdot y}{\spc{X}}
=
\dist{g^{-1}\cdot x}{y}{\spc{X}}}{g\in G}.\]

If the $G$-orbits are closed, then $\dist{G\cdot x}{G\cdot y}{\spc{X}/G}=0$ if and only if $G\cdot x=G\cdot y$.
In this case, the quotient space $\spc{X}/G$ is a genuine metric space.

The following proposition follows from the definition of a quotient space:

\begin{thm}{Proposition}\label{prop:submetry-X/G}
Assume $\spc{X}$ is a metric space and a  group $G$ acts on $\spc{X}$ by isometries.
Then the projection $\pi\co \spc{X}\to\spc{X}/G$ is a \emph{submetry};
that is, $\pi(\oBall(p,r))=\oBall(\pi(p),r)$ for any $p\in \spc{X}, r>0$ (see Definition~\ref{def:submetry}).
\end{thm}

\section{Gluing and doubling}\label{sec:doubling}

\parbf{Gluing.}
Recall that the disjoint union of metric spaces can be also considered as a metric space; see Definition~\ref{def:disjoint-union}.
Therefore the quotient space construction works as well for an equivalence relation on the disjoint union of metric spaces.

Consider two metric spaces $\spc{X}_1$ and $\spc{X}_2$
with subsets $A_1\subset\spc{X}_1$ and $A_2\subset\spc{X}_2$,
and a bijection $\phi\:A_1\to A_2$.
Consider the minimal equivalence relation on $\spc{X}_1\sqcup\spc{X}_2$
such that $a\sim \phi(a)$ for any $a\in A_1$.
In this case, the corresponding quotient space 
$(\spc{X}_1\sqcup\spc{X}_2)/{\sim}$ will be called the \index{gluing}\emph{gluing of $\spc{X}$ and $\spc{Y}$ along $\phi$} and denoted by
\[\spc{X}_1\sqcup_{\phi}\spc{X}_2.\]

Note that if the map $\phi\:A_1\to A_2$ is distance-preserving,
then the inclusions $\iota_i\:\spc{X}_i\to \spc{X}_1\sqcup_{\phi}\spc{X}_2$ are also distance-preserving, and 
\[\dist{\iota_1(x_1)}{\iota_2(x_2)}{\spc{X}_1\sqcup_{\phi}\spc{X}_2}
=
\inf_{a_2=\phi(a_1)}
\{\,\dist{x_1}{a_1}{\spc{X}_1}+\dist{x_2}{a_2}{\spc{X}_2}\,\}\]
for any $x_1\in \spc{X}_1$ and $x_2\in \spc{X}_2$.

\parbf{Doubling.}
Let $\spc{V}$ be a metric space 
and $A\subset \spc{V}$ be a closed subset.
A metric space $\spc{W}$ glued from two copies of $\spc{V}$ along $A$ is called the \index{doubling}\emph{doubling of $\spc{V}$ in $A$}.

The space $\spc{W}$ is completely described by the following properties:
\begin{itemize}
\item The space $\spc{W}$ contains $\spc{V}$ as a subspace; 
in particular the set $A$ can be treated as a subset of $\spc{W}$.
\item There is an isometric involution of $\spc{W}$ which is called \index{reflection}\emph{reflection in $A$};
it will be denoted by $x\mapsto x'$.
\item For any $x\in \spc{W}$ we have $x\in \spc{V}$ or $x'\in \spc{V}$ and 
\[
\dist{x'}{y}{\spc{W}}
=\dist{x}{y'}{\spc{W}}
=\inf_{a\in A}\{\dist{x}{a}{\spc{V}}+\dist{a}{y}{\spc{V}}\}
\]
for any $x,y\in \spc{V}$.
\end{itemize}

The image of $\spc{V}$ under reflection in $A$ will be denoted by $\spc{V}'$.
The subspace $\spc{V}'$ is an isometric copy of $\spc{V}$.
Clearly $\spc{V}\cup\spc{V}'=\spc{W}$ and $\spc{V}\cap\spc{V}'=A$.
Moreover, $a=a'$ $\iff$ $a\in A$.

The following proposition follows directly from the definitions.

\begin{thm}{Proposition}\label{prop:doubling}
Assume $\spc{W}$ is the doubling of the metric space $\spc{V}$ in its closed subset $A$.
Then: 

\begin{subthm}{}
If $\spc{V}$ is a complete length space, then so is $\spc{W}$.
\end{subthm}
 
\begin{subthm}{}
If $\spc{V}$ is proper, then so is $\spc{W}$.
In this case, for any $x,y\in\spc{V}$ there is $a\in A$ such that 
\[\dist{x}{a}{\spc{V}}+\dist{a}{y}{\spc{V}}=\dist{x}{y'}{\spc{W}}.\]
\end{subthm}

\begin{subthm}{prop:doubling:projection}
Given $x\in \spc{W}$, let $\bar x=x$ if $x\in \spc{V}$,
and $\bar x=x'$ otherwise. The map $\spc{W}\to\spc{V}$ defined by $x\mapsto \bar x$ is short and length-preserving.
In particular, if $\gamma$ is a geodesic in $\spc{W}$ with ends in $\spc{V}$, then $\bar\gamma$ is a geodesic in $\spc{V}$ with the same ends.
\end{subthm}
\end{thm}

\section{Kuratowsky embedding}\label{Kuratowsky embedding}

Given a metric space $\spc{X}$, 
let us denote by $\Bnd(\spc{X},\RR)$\index{ $\Bnd(\spc{X},\RR)$}  the space of all bounded functions on $\spc{X}$ equipped with the sup-norm
\[\|f\|=\sup_{x\in\spc{X}}\{|f(x)|\}.\]

\parbf{Kuratowski embedding.}
Given a point $p\in\spc{X}$, consider the 
map $\kur_p\:\spc{X}\to\Bnd(\spc{X},\RR)$ 
defined by $\kur_p x=\distfun{x}{}{}-\distfun{p}{}{}$.
The map $\kur_p$ will be called the \index{Kuratowski map}\emph{Kuratowski map at $p$}.

From the triangle inequality, we have
\[\|\kur_p x-\kur_p y\|
=
\sup_{z\in\spc{X}}\{|\dist{x}{z}{}-\dist{y}{z}{}|\}
=
\dist{x}{y}{}.\]
Therefore, for any $p\in\spc{X}$, the Kuratowski map gives a distance-preserving map $\kur_p\:\spc{X}\hookrightarrow\Bnd(\spc{X},\RR)$.
Thus we can (and often will) consider the space $\spc{X}$ as a subset of  $\Bnd(\spc{X},\RR)$.

\begin{thm}{Exercise}\label{ex:compact-in-lenght}
Show that any compact metric space is isometric to a subspace in a compact length space.
\end{thm}


\chapter{Maps and functions}

Here we introduce several classes of maps between metric spaces and develop a language to describe various notions of convexity/concavity of real-valued functions on general metric spaces.

\section{Submaps}\label{sec:submaps}

We will often need maps and functions defined on subsets of a metric space.
We call them \index{submap}\emph{submaps} and \index{subfunction}\emph{subfunctions}.
Thus, given  metric spaces $\spc{X}$ and $\spc{Y}$, 
a submap \index{$\subto$}$\map\:\spc{X}\subto \spc{Y}$ is a map defined on a subset $\Dom\map\subset \spc{X}$.

A submap $\map\:\spc{X}\subto \spc{Y}$ is \index{continuous submap}\emph{continuous} if the inverse image of any open set is open.
Since $\Dom\map=\map^{-1}(\spc{Y})$, the domain $\Dom\map$ \index{$\Dom\map$ } of a continuous submap is open.
The same holds for upper and lower semicontinuous functions $f\:\spc{X}\subto \RR$ since they are  continuous functions for a special topology on $\RR$.

(Continuous partially defined maps could be defined via closed sets; namely, one could require that inverse images of closed sets are closed.
While this condition is equivalent to continuity for functions defined on the whole space,
it is different for partially defined functions. 
In particular, with this definition the domain of a continuous submap would have to be closed.)

\section{Lipschitz conditions}

\begin{thm}{Lipschitz maps}
Suppose $\spc{X}$ and $\spc{Y}$ are metric spaces, 
$\map\:\spc{X}\subto\spc{Y}$ is a continuous submap,  
and $\Lip\in\RR$.

\begin{subthm}{}
The submap $\map$ is called \index{Lipschitz map}\emph{$\Lip$-Lipschitz} if
\[\dist{\map(x)}{\map(y)}{\spc{Y}}
\le
\Lip\cdot
\dist{x}{y}{\spc{X}}\]  
for any two points $x,y\in\Dom \map$.

\begin{itemize}
 \item $1$-Lipschitz maps will be also called \index{short map}\emph{short}.
\end{itemize}

\end{subthm}

\begin{subthm}{}
We say that $\map$ is \index{Lipschitz map}\emph{Lipschitz} if it is $\Lip$-Lipschitz for a constant~$\Lip$.
The minimal such constant is denoted by $\lip\map$.
\end{subthm}

\begin{subthm}{}
We say that $\map$ is \emph{locally Lipschitz} 
if any point $x\in\Dom \map$ admits a neighborhood 
$\Omega\subset \Dom\map$ such that the restriction $\map|_\Omega$ is Lipschitz.
\end{subthm}

\begin{subthm}{}
Given $p\in\Dom \map$, we denote by $\lip_p\map$ the infimum of the real values $\Lip$ such that
$p$ admits  a neighborhood 
$\Omega\subset \Dom\map$ such that the restriction $\map|_\Omega$ is $\Lip$-Lipschitz.
\end{subthm}
\end{thm}

Note that $\map\:\spc{X}\to\spc{Y}$ is $\Lip$-Lipschitz if and only if
\[\map(\oBall(x, R)_{\spc{X}})\subset\oBall(\map(x),\Lip\cdot R)_{\spc{Y}}\]
for any $R\ge 0$ and $x\in \spc{X}$.
A dual version of this property is considered in the following definition.

\begin{thm}{Definitions}
Let 
$\map\:\spc{X}\to\spc{Y}$ be a map between metric spaces,  
and $\Lip\in\RR$.
\begin{subthm}{}
The map $\map$ is called \index{$\Lip$-co-Lipschitz map}\emph{$\Lip$-co-Lipshitz} if 
\[\map(\oBall(x,\Lip\cdot R)_{\spc{X}})\supset\oBall(\map(x),R)_{\spc{Y}}\]
for any $x\in \spc{X}$ and $R>0$.
\end{subthm}

\begin{subthm}{}
The map $\map$ is called \index{co-Lipschitz map}\emph{co-Lipschitz} if it is $\Lip$-co-Lipschitz
for some constant $\Lip$.
The minimal such constant is denoted by $\colip\map$.

\end{subthm}
\end{thm}

From the definition of co-Lipschitz maps we get the following:

\begin{thm}{Proposition}
Any co-Lipschitz map is open and surjective.
\end{thm}

In other words,  $\Lip$-co-Lipschitz maps 
can be considered as a quantitative version of open maps.
For that reason they are also called $\Lip$-open \cite{burago-gromov-perelman}.
Also, be aware that some authors
refer to our $\Lip$-co-Lipschitz maps
as $\tfrac1\Lip$-co-Lipschitz.

\begin{thm}{Proposition}\label{prop:colip=>complete}
Let $\spc{X}$ and $\spc{Y}$ be metric spaces such that $\spc{X}$ is complete, and let
$\map\: \spc{X}\to\spc{Y}$ be a continuous co-Lipschitz map. 
Then $\spc{Y}$ is complete.
\end{thm}

\parit{Proof.}
Choose a Cauchy sequence $y_n$ in $\spc{Y}$.
Passing to a subsequence if necessary, we may assume that $\dist{y_n}{y_{n+1}}{\spc{Y}}< \tfrac1{2^n}$ for each~$n$.
It is sufficient to show that $y_n$ converges in $\spc{Y}$.

Denote by $\Lip$ a co-Lipschitz constant for $\map$.
Note that  there is a sequence $x_n$ in $\spc{X}$
such that
\[\map(x_n)=y_n\quad \text{and}\quad \dist{x_n}{x_{n+1}}{\spc{X}}< \tfrac{\Lip}{2^n}\eqlbl{eq:colip+1/2n}\]
for each $n$. 
Indeed, such a sequence can be constructed recursively. 
Assuming that the points $x_1,\dots,x_{n-1}$ are already constructed, 
the existence of a sequence $x_n$ satisfying \ref{eq:colip+1/2n}
follows since $\map$ is $\Lip$-co-Lipschitz.

Notice that the sequence $x_n$ is Cauchy.
Since $\spc{X}$ is complete, $x_n$ converges in $\spc{X}$; denote its limit by $x_\infty$ 
and set $y_\infty= \map(x_\infty)$.
Since $\map$ is continuous,
$y_n\to y_\infty$ as $n\to\infty$.
Hence the result.
\qeds

\begin{thm}{Lemma}\label{lem:lip-approx}
Let $\spc{X}$ be a metric space and $f\:\spc{X}\to\RR$ be a continuous function.
Then for any $\eps>0$ there is a locally Lipschitz function $f_\eps\:\spc{X}\to\RR$
such that $|f(x)-f_\eps(x)|<\eps$ for any $x\in \spc{X}$.
\end{thm}

\parit{Proof.}
Assume that $f\ge 1$.
Construct a continuous positive function $\rho\:\spc{X}\to \RR_{>0}$ such that 
\[\dist{x}{y}{}<\rho(x)\quad \Rightarrow\quad |f(x)-f(y)|<\eps.\]
Consider the function
\[
f_\eps(y)
=
\sup\set{f(x)
\cdot
\left(1-\tfrac{\dist{x}{y}{}}{\rho(x)}\right)}%
{x\in\spc{X}}.
\]
It is straightforward to check that each $f_\eps$ is locally Lipschitz and $0\le f_\eps-f<\eps$.

Since any continuous function can be presented as the difference of two continuous functions bounded below by $1$, the result follows.
\qeds

\section{Isometries and submetries}\label{sec:quotient-CBB}

\begin{thm}{Isometry}\label{def:isometry}
Let $\spc{X}$ and $\spc{Y}$ be metric spaces
and $\map\:\spc{X}\to \spc{Y}$ be a map
\begin{subthm}{}
The map $\map$ is \index{distance-preserving map}\emph{distance-preserving} if
$$\dist{\map(x)}{\map(x')}{\spc{Y}}=\dist{x}{x'}{\spc{X}}$$
for any $x,x'\in \spc{X}$.
\end{subthm}

\begin{subthm}{}
A distance-preserving bijection $\map$ is called an \index{isometry}\emph{isometry}.
\end{subthm}

\begin{subthm}{}
The spaces $\spc{X}$ and $\spc{Y}$ are called \emph{isometric} (briefly $\spc{X}\iso \spc{Y}$)
 if there is an isometry  $\map\:\spc{X}\to \spc{Y}$.
\end{subthm}

\end{thm}

\begin{thm}{Submetry}\label{def:submetry}
A map $\sigma\:\spc{X}\to\spc{Y}$ between the metric spaces $\spc{X}$ and $\spc{Y}$
is called a \index{submetry}\emph{submetry} if 
\[\sigma(\oBall(p,r)_\spc{X})=\oBall(\sigma(p),r)_{\spc{Y}}\]
for any $p\in \spc{X}$ and $r\ge 0$.
\end{thm}

Note $\sigma\:\spc{X}\to\spc{Y}$ is a submetry if it is 1-Lipschitz and 1-co-Lipschitz.

Note also that any submetry is an onto map.

The main source of examples of submetries comes from isometric group actions.
Namely, assume $\spc{X}$ is a metric space and $G$ is a subgroup of isometries of $\spc{X}$.
Denote by $[x]=G\cdot x$ the $G$-orbit of $x\in\spc{X} $ and $\spc{X}/G$ the set of all $G$-orbits;
let us equip it with the pseudometric defined by
\[\dist{[x]}{[y]}{\spc{X}/G}=\inf\set{\dist{g\cdot x}{h\cdot y}{\spc{X}}}{g,h\in G}.\]
Note that if all the $G$-orbits form closed sets in $\spc{X}$,
then $\spc{X}/G$ is a genuine metric space.

\begin{thm}{Proposition}\label{prop:submet/G}
Let $\spc{X}$ be a metric space.
Assume that a group $G$  acts on $\spc{X}$ by isometries  
and in such a way that every $G$-orbit is closed.
Then the projection map $\spc{X}\to \spc{X}/G$ is a submetry.
\end{thm}

\parit{Proof.}
We need to show that the map $x\mapsto[x]=G\cdot x$ is $1$-Lipschitz and $1$-co-Lipschitz.
The co-Lipschitz part follows directly from the definitions of Hausdorff distance and co-Lipschitz maps.

Assume $\dist{x}{y}{\spc{X}}< r$; equivalently $\oBall(x,r)_{\spc{X}}\ni y$.
Since the action $G\acts \spc{X}$ is isometric, 
$\oBall(g\cdot x,r)_{\spc{X}}\ni g\cdot y$ for any $g\in G$.

In particular, the orbit $G\cdot y$ lies in the open $r$-neighborhood of the orbit $G\cdot x$.
In the same way we can prove that the orbit $G\cdot x$ lies in the open $r$-neighborhood of the orbit $G\cdot y$. 
That is, the Hausdorff distance between the orbits $G\cdot x$ and $G\cdot y$ is less than $r$
or, equivalently, $\dist{[x]}{[y]}{\spc{X}/G}< r$.
Since $x$ and $y$ are arbitrary, the map $x\mapsto[x]$ is $1$-Lipschitz.
\qeds

\begin{thm}{Proposition}
\label{prop:submet-length}
Let  $\spc{X}$ be a length space 
and $\sigma\:\spc{X}\to \spc{Y}$ be a submetry.
Then $\spc{Y}$ is a length space.
\end{thm}

\parit{Proof.}
Fix $\eps>0$ and a pair of points $x,y\in \spc{Y}$.

Since $\sigma$ is $1$-co-Lipschitz, there are points $\hat x,\hat y\in \spc{X}$
such that 
$\sigma(\hat x)\z=x$,
$\sigma(\hat y)\z=y$, 
and $\dist{\hat x}{\hat y}{\spc{X}}<\dist{x}{y}{\spc{Y}}+\eps$.

Since $\spc{X}$ is a length space, 
there is a curve $\gamma$ 
joining $\hat x$ to $\hat y$ in ${\spc{X}}$
such that
\[\length\gamma\le \dist{x}{y}{\spc{Y}}+\eps.\]

The curve $\sigma\circ\gamma$ joins $x$ to $y$.
Since $\sigma$ is $1$-Lipschitz,
and by the above,
\begin{align*}
\length\sigma\circ\gamma&\le \length\gamma\le
\\
&\le\dist{x}{y}{\spc{Y}}+\eps.
\end{align*}
Since $\eps>0$ is arbitrary,
$\spc{Y}$ is a length space.
\qeds

\section{Speed of curves}
\label{sec: speed}

Let $\spc{X}$ be a metric space.
Recall that a \emph{curve}  
in $\spc{X}$ is a continuous map $\alpha\:\II\to \spc{X}$, where $\II$ is a real interval. 
A curve is called \index{Lipschitz map!Lipschitz curve}\emph{Lipschitz} or \emph{locally Lipschitz} if $\alpha$ is Lipschitz or locally Lipschitz respectively. 
 Length of curves is defined in \ref{def:length}.

The following theorem follows from \cite[2.7]{burago-burago-ivanov}.

\begin{thm}{Theorem}\label{thm:speed}
Let $\spc{X}$ be a metric space  
and $\alpha\:\II \to \spc{X}$ be a locally Lipschitz
curve. 
Then the speed function
\[\speed_{t_0}\alpha
=
\lim_{\substack{t\to t_0+\\s\to t_0-}}\frac{\dist{\alpha(t)}{\alpha(s)}{}}{|t-s|}\] 
is defined for almost all $t_0 \in \II$, and 
\[\length\alpha=\int\limits_\II \speed_{t}\alpha\cdot\dd t,\]
where $\int$ denotes the Lebesgue integral.
\end{thm}

A curve $\alpha\:\II\to\spc{X}$ is \index{unit-speed curve}\emph{unit-speed}
if 
\[b-a=\length(\alpha|_{[a,b]})\]
for any subinterval $[a,b]\subset\II$.
If $\alpha$ is Lipschitz, then, according to the above theorem, this is equivalent to 
\[\speed\alpha\ae 1.\]

The following generalization of the standard Rademacher theorem 
on differentiability almost everywhere of Lipschitz maps between smooth manifolds \cite[5.5.2]{burago-burago-ivanov} was proved by Bernd Kirchheim \cite{kirchheim}. 

The conclusion of the standard Rademacher theorem does not make sense for maps to a metric space since the target might have no linear structure.
But the theorem does not hold even if we assume that the target is a Banach space.
For example, consider the map $[0,1]\to L^1[0,1]$, defined by $x\mapsto \chi_{[0,x]}$, where $\chi_A$ denotes the characteristic function of $A$.
This map is distance-preserving and in particular Lipschitz,
but its differential is undefined at any point.

\begin{thm}{Theorem}\label{thm:Rademacher-md}
Let $\spc{X}$ be a metric space 
and $f\:\RR^n \subto \spc{X}$ be $1$-Lipschitz. 
Then for almost all $x\in\Dom f$ there is a pseudonorm 
$\lVert*\rVert_x$ on $\RR^n$ such that
\[\dist{f(y)}{f(z)}{\spc{X}}=\lVert z-y\rVert_x+o(|y-x|+|z-x|).\]
\end{thm}

Given $f$, the (pseudo)norm $\lVert*\rVert_x$ in the above theorem 
will be called its \index{differential of a metric}\emph{differential of the induced metric} at $x$, or \index{metric differential}\emph{metric differential} at $x$.

\section{Convex real-to-real functions}\label{sec:conv-real}

In this section we will discuss generalized solutions
of the following differential inequalities
\[y''+\kappa\cdot  y\ge \lambda
\quad \text{and respectively}
\quad y''+\kappa\cdot  y\le \lambda
\eqlbl{eq:sec:conv-real*}\]
for fixed $\kappa,\lambda\in\RR$.
The solutions $y\:\RR\subto\RR$ are only assumed to be upper (respectively lower) semicontinuous subfunctions.

The inequalities  \ref{eq:sec:conv-real*} are understood in the sense of distributions.
That is, for any smooth function $\phi$ with compact support $\supp\phi\subset\Dom y$ the following inequality should be satisfied:
\[\begin{aligned}
\int\limits_{\Dom y}\left[y(t)\cdot\phi''(t)+\kappa\cdot  y(t)\cdot\phi(t)-\lambda\right]\cdot\dd t
&\ge 0,
\\
\text{respectively}\quad &\le0.
\end{aligned}
\eqlbl{eq:distr-conc}\]
The integral is understood in the sense of Lebesgue;
in particular the inequality \ref{eq:distr-conc}
makes sense for any Borel-measurable subfunction $y$.
The proofs of the following propositions are straightforward.

\begin{thm}{Proposition}
Let $\II\subset\RR$ be an open interval and $y_n\:\II\to\RR$ be a sequence of solutions of one of the inequalities in \ref{eq:sec:conv-real*}.
Assume $y_n(t)\to y_\infty(t)$ as $n\to\infty$ for any $t\in \II$.
Then $y_\infty$ is a solution of the same inequality in \ref{eq:sec:conv-real*}.
\end{thm}

Assume $y$ is a solution of one of the inequalities in \ref{eq:sec:conv-real*}.
For $t_0\in \Dom y$, let us define the \emph{right (left) derivative } $y^+(t_0)$ ($y^-(t_0)$) at $t_0$ by
\[y^\pm(t_0)=\lim_{t\to t_0\pm} \frac{y(t)-y(t_0)}{|t-t_0|}.\]
Note that our sign convention for $y^-$ is not standard --- for $y(t)=t$ we have
$y^+(t)=1$ and $y^-(t)=-1$.

\begin{thm}{Proposition}\label{prop:derivative-of-convex-function}
Let $\II\subset\RR$ be an open interval and $y\:\II\to\RR$ be a solution of an inequality in \ref{eq:sec:conv-real*}.
Then  $y$ is locally Lipschitz; its right and left derivatives $y^+(t_0)$ and $y^-(t_0)$ are defined
for any $t_0\in\II$.
Moreover 
\[y^+(t_0)+y^-(t_0)\ge 0
\quad \text{or respectively}
\quad y^+(t_0)+y^-(t_0)\le 0.\]
\end{thm}

The next theorem gives a  number of equivalent ways to define such  generalized solutions.

\begin{thm}{Theorem}\label{y''=<1-ky}
Let $\II$ be an open real interval and $y\:\II\to\RR$ be a locally Lipschitz function.
Then the following conditions are equivalent:
\begin{subthm}{}$y''\ge \lambda-\kappa\cdot  y$ (respectively $y''\le \lambda-\kappa\cdot  y).$
\end{subthm}

\begin{subthm}{barrier}(barrier inequality) For any $t_0\in \II$, 
there is a solution $\bar y$ 
of the ordinary differential equation $\bar y''=\lambda-\kappa\cdot  \bar y$ 
with $\bar y(t_0)\z= y(t_0)$ such that $\bar y\ge y$ (respectively $\bar y\le y$) for all $t\in [t_0-\varpi\kappa,t_0\z+\varpi\kappa]\cap \II$.

The function $\bar y$ is called a {}\emph{lower} (respectively {}\emph{upper}) \index{barrier}\emph{barrier} of $y$ at $t_0$.
\end{subthm}

\begin{subthm}{y''-mono} (Jensen's inequality)\index{Jensen's inequality}
For any pair of values $t_1<t_2$ in $\II$ such that $|t_2-t_1|<\varpi\kappa$,  the unique solution $z(t)$ of \[z''\z=\lambda-\kappa\cdot  z\] such that
\[z(t_1)=y(t_1),\quad z(t_2)=y(t_2)\] 
satisfies $y(t)\le z(t)$ (respectively $y(t)\ge z(t)$) for all $t\in[t_1,t_2]$.
\end{subthm}

Further, the following property holds:

\begin{subthm}{barrier'} 
Suppose $y''\le \lambda-\kappa\cdot  y$. Let $t_0\in\II$, and   $\bar y$ be a solution of  the
 ordinary differential equation $\bar y''=\lambda-\kappa\cdot  \bar y$ 
such that  $\bar y(t_0)\z= y(t_0)$ and 
$y^+(t_0)\le \bar y^+(t_0)\le -y^-(t_0)$.
(Note that such a $\bar{y}$ is unique if $y$ is differentiable at $t_0$.) 

Then $\bar y\ge y$  for all $t\in [t_0-\varpi\kappa,t_0+\varpi\kappa]\cap \II$; that is, $\bar{y}$ is a barrier of $y$ at $t_0$. (Similarly, by reversing inequalities, for $y''\ge \lambda-\kappa\cdot  y$.) 
\end{subthm}
\end{thm}

Note that Theorem ~\ref{y''=<1-ky} implies that $y$ satisfies $y''\ge \lambda$ ($y''\le \lambda$)  on an interval $\II\subset\RR$  if and only if $y(t)-\frac{\lambda}{2}\cdot t^2$ is convex (concave) on $\II$.

We will often need the following fact about convergence of derivatives of convex functions:

{\sloppy 

\begin{thm}{Two-shoulder lemma}\label{lem:der-conv-lim}
Let $\II$ be an open interval 
and $f_n\:\II\to\RR$ be a sequence of 
convex functions. 
Assume the functions $f_n$ pointwise converge to a function $f_\infty\:\II\to\RR$.
Then for any $t_0\in \II$,
\[f_\infty^\pm(t_0)\le \liminf_{n\to\infty}f^\pm_n(t_0).\]
\end{thm}

}

\parit{Proof.}
Since the $f_n$ are convex, we have $f^+_n(t_0)+f^-_n(t_0)\ge0$, and for any~$t$,
\[f_n(t)\ge f_n(t_0)\pm f^\pm(t_0)\cdot (t-t_0).\]
Passing to the limit, we get
\[f_\infty(t)\ge f_\infty(t_0)+\left[\limsup_{n\to\infty}f^+_n(t_0)\right]\cdot (t-t_0)\]
for $t\ge t_0$, and 
\[f_\infty(t)\ge f_\infty(t_0)-\left[\limsup_{n\to\infty}f^-_n(t_0)\right]\cdot (t-t_0)\]
for $t\le t_0$.
Hence the result.
\qeds

\begin{thm}{Corollary}
\label{cor:der-conv-lim}
Let $\II$ be an open interval 
and $f_n\:\II\to\RR$ be a sequence of functions such that $f_n''\le \lambda$ that converge pointwise to a function $f_\infty\:\II\to\RR$.
Then: 
\begin{subthm}{} If $f_\infty$ is differentiable at $t_0\in \II$, then
\[f_\infty'(t_0)=\pm\lim_{n\to\infty} f^\pm_n(t_0).\]
\end{subthm}

\begin{subthm}{} If all $f_n$ and $f_\infty$ are differentiable at $t_0\in \II$, then
\[f_\infty'(t_0)=\lim_{n\to\infty} f'_n(t_0).\]
\end{subthm}
\end{thm}

\parit{Proof.} Set $\hat f_n(t)=f_n(t)-\tfrac{\lambda}{2}\cdot t^2$ and $\hat f_\infty(t)=f_\infty(t)-\tfrac\lambda2\cdot t^2$.
Note that the $\hat f_n$ are concave and $\hat f_n\to \hat f_\infty$ pointwise.
Thus the theorem follows from the two-shoulder lemma (\ref{lem:der-conv-lim}).\qeds

\section{Convex functions on a metric space}\label{sec:conv-fun}

In this section we define different types of convexity/concavity
in the context of metric spaces; it will be mostly used for geodesic spaces.
The notation refers to the corresponding second-order ordinary differential inequality. 

\begin{thm}{Definition}\label{def:lam-convex}
Let $\spc{X}$ be a metric space.
We say that an upper semicontinuous subfunction $f\:\spc{X}\subto(-\infty,\infty]$ 
satisfies the inequality
\[f''+\kappa\cdot  f\ge \lambda\]
if for \emph{any} unit-speed geodesic $\gamma$ in $\Dom f$, 
the real-to-real function $y(t)\z= f\circ\gamma(t)$
satisfies 
\[y''+\kappa\cdot  y\ge \lambda\]
in the domain $\set{t}{y(t)<\infty}$;
see the definition in Section~\ref{sec:conv-real}.

We say that a lower semicontinuous subfunction $f\:\spc{X}\subto[-\infty,\infty)$ 
satisfies the inequality
\[f''+\kappa\cdot  f\le \lambda\]
if the subfunction $h=-f$ 
satisfies 
\[h''-\kappa\cdot  h\ge -\lambda.\]

Functions satisfying the inequalities
\[f''\ge \lambda\quad\text{and}\quad f''\le \lambda\]
are called 
\index{convex function}\index{concave function}\index{$\lambda$-convex function}\index{$\lambda$-concave function}\emph{$\lambda$-convex} and \emph{$\lambda$-concave} respectively.

$0$-convex and $0$-concave subfunctions will also be called \emph{convex} and \emph{concave} respectively.

If $f$ is $\lambda$-convex for $\lambda>0$, then $f$ will be called \index{strongly convex function}\emph{strongly convex};
correspondingly, if $f$ is $\lambda$-concave for $\lambda<0$, then $f$ will be called \index{strongly concave function}\emph{strongly concave}.

If for any point $p\in\Dom f$ 
there is a neighborhood $\Omega\ni p$ and a real number $\lambda$
such that the restriction $f|_\Omega$ is $\lambda$-convex (or $\lambda$-concave),
then $f$ is called \index{semiconvex function}\index{semiconcave function}\emph{semiconvex} (respectively \emph{semiconcave}).\end{thm}

Various authors define the class of $\lambda$-convex ($\lambda$-concave) functions differently. 
Their definitions may correspond to $\pm\lambda$-convex ($\pm\lambda$-concave) or $\pm\tfrac\lambda2$-convex ($\pm\tfrac\lambda2$-concave) functions in our definitions.

\begin{thm}{Proposition}\label{prop:conv-comp}
Let $\spc{X}$ be a metric space.
Assume that $f\:\spc{X}\subto \RR$ is a semiconvex subfunction
and $\phi\:\RR\to\RR$ is a nondecreasing semiconvex function.
Then the composition $\phi\circ f$ is a semiconvex subfunction.
\end{thm}

The proof is straightforward.

\chapter{Ultralimits}

Here we introduce ultralimits of sequences of points, metric spaces, and functions.
Our presentation is based on \cite{kleiner-leeb}.

Ultralimits are closely related to Gromov--Hausdorff limits.
We use them only as a canonical way to pass to  convergent subsequences.
We could avoid using them at the cost of saying ``pass to a convergent subsequence'' too many times.
Doing this might be cumbersome and it obscure ideas of the proof;
see for example the proof of the globalization theorem for general $\Alex{}$ spaces.
Also the use of ultralimits is convenient when dealing with $\CAT{}$ spaces due to the lack of compactness results.

\section{Ultrafilters}

We will need the existence of a selective ultrafilter $\o$ that will be fixed once and  for all.
The existence follows from the axiom of choice and the continuum hypothesis.

\parbf{Measure-theoretic definition.}
Recall that $\NN$ denotes the set of natural numbers, $\NN=\{1,2,\dots\}$.

\begin{thm}{Definition}\label{def:ultrafilter}
A finitely additive measure $\o$ 
on  $\NN$ 
is called an \index{ultrafilter}\emph{ultrafilter} if it satisfies 
\begin{subthm}{}
$\o(S)=0$ or $1$ for any subset $S\subset \NN$.
\end{subthm}
An ultrafilter $\o$ is called 
\index{ultrafilter!nonprincipal ultrafilter}\index{nonprincipal ultrafilter}\emph{nonprincipal} if in addition 
\begin{subthm}{}
$\o(F)=0$ for any finite subset $F\subset \NN$.
\end{subthm}
A nonprincipal ultrafilter $\o$ is called 
\emph{selective}\index{ultrafilter!selective ultrafilter}\index{selective ultrafilter} if in addition 
\begin{subthm}{}
for any partition of $\NN$ into sets $\{C_\alpha\}_{\alpha\in\IndexSet}$ such that $\o(C_\alpha)\z=0$ for each $\alpha$, 
there is a set $S\subset \NN$ such that $\o(S)=1$ and $S\cap C_\alpha$ is a one-point set for each $\alpha\in\IndexSet$.
\end{subthm}
\end{thm}

If $\o(S)=0$ for some subset $S\subset \NN$,
we say that $S$ is \index{$\o$-small}\emph{$\o$-small}. 
If $\o(S)=1$, we say that $S$ contains \index{$\o$-almost all}\emph{$\o$-almost all} elements of $\NN$.

\begin{thm}{Advanced exercise}\label{ex:ultrakatetov}
Let $\o$ be an ultrafilter and $f\:\NN\z\to \NN$.
Suppose that $\o(S)\le \o(f^{-1}(S))$ for any set $S\subset \NN$.
Show that $f(n)=n$ for $\o$-almost all $n\in\NN$.
\end{thm}

\parbf{Classical definition.}
More commonly, a nonprincipal ultrafilter is defined as a collection, say $\mathfrak{F}$, of sets in $\NN$ such that
\begin{enumerate}
\item\label{filter:supset} if $P\in \mathfrak{F}$ and $Q\supset P$, then $Q\in \mathfrak{F}$,
\item\label{filter:cap} if $P, Q\in \mathfrak{F}$, then $P\cap Q\in \mathfrak{F}$,
\item\label{filter:ultra} for any subset $P\subset\NN$, either $P$ or its complement is an element of $\mathfrak{F}$,
\item\label{filter:non-prin} if $F\subset \NN $ is finite, then $F\notin \mathfrak{F}$.
\end{enumerate}

Setting 
\[P\in\mathfrak{F}\quad\iff\quad\o(P)=1\] 
makes these two definitions equivalent.

A nonempty collection of sets $\mathfrak{F}$ that does not include the empty set and satisfies only conditions \ref{filter:supset} and \ref{filter:cap} is called a \index{filter}\emph{filter}; 
if in addition $\mathfrak{F}$ satisfies Condition~\ref{filter:ultra} it is called an \index{ultrafilter}\emph{ultrafilter}.
From Zorn's lemma, it follows that every filter is contained in an ultrafilter.
Thus there is an ultrafilter $\mathfrak{F}$ contained in the filter of all complements of finite sets; clearly this $\mathfrak{F}$ is nonprincipal.

The existence of a selective ultrafilter follows from the continuum hypothesis;
this was proved by Walter Rudin \cite{rudin}.

\parbf{Stone--\v{C}ech compactification.}
Given a set $S\subset \NN$, consider the subset $\Omega_S$ of all ultrafilters $\o$ such that $\o(S)=1$.
It is straightforward to check that the sets $\Omega_S$ for all $S\subset \NN$ form a topology on the set of ultrafilters on~$\NN$. 
The resulting space is called the \index{Stone--\v{C}ech compactification}\emph{Stone--\v{C}ech compactification} of~$\NN$; it is usually denoted by $\beta\NN$\index{ $\beta\NN$} .

There is a natural embedding $\NN\hookrightarrow\beta\NN$ defined by 
$n\mapsto\o_n$, where $\o_n$ is the principal ultrafilter such that $\o_n(S)=1$ if and only if $n\in S$. 
Using this embedding, we can (and will) consider $\NN$ as a subset of $\beta\NN$.

The space $\beta\NN$ is the maximal compact Hausdorff space that contains $\NN$  as an everywhere dense subset.
More precisely, for any compact Hausdorff space $\spc{X}$ 
and a  map $f\:\NN\to \spc{X}$, there is a unique continuous map $\bar f\:\beta\NN\to X$ such that the restriction $\bar f|_\NN$ coincides with $f$. 

\section{Ultralimits of points}
\label{ultralimits}

Fix an ultrafilter $\o$.
Assume $x_n$ is a sequence of points in a metric space $\spc{X}$. 
Define an  \index{ultralimit of points}\emph{$\o$-limit} of $x_n$ to be a point $x_\o$ 
such that for any $\eps>0$, $\o$-almost all elements of $x_n$ lie in $\oBall(x_\o,\eps)$; 
that is,
\[\o\set{n\in\NN}{\dist{x_\o}{x_n}{}<\eps}=1.\]
In this case, we write 
\[x_\o=\lim_{n\to\o} x_n
\quad \text{or}\quad 
x_n\to x_\o\quad \text{as}\quad n\to\o.\]

Also, if $\spc{X}=\RR$ we write $\lim_{n\to\o} x_n=\pm\infty$ if 
\[\o\set{n\in\NN}{\pm x_n>L}=1\] for any $L\in\RR$.

It easily follows from the definition that  $\o$-limits are unique if they exist. 
For example, if $\o$ is the principal ultrafilter such that $\o(\{n\})=1$ for some $n\in\NN$, then
$x_\o=x_n$.

Note that $\o$-limits of a sequence and its subsequences may differ.
For example, in general
\[\lim_{n\to\o}x_n
\ne
\lim_{n\to\o}x_{2\cdot n}.\]

The sequence $x_n$ can be regarded as a map $\NN\to\spc{X}$.
If $\spc{X}$ is compact, then this map can be uniquely extended to a continuous map to the Stone--\v{C}ech compactification $\beta\NN$ of $\NN$.
Then $x_\o$ is the image of~$\o$. 

\begin{thm}{Proposition}\label{prop:ultra/partial}
Let $\o$ be a nonprincipal ultrafilter.
Assume $x_n$ is a sequence of points in a metric space $\spc{X}$
and $x_n\to  x_\o$ as $n\to\o$.
Then there is a subsequence of $x_n$ that converges to $x_\o$ in the usual sense.

Moreover, if $\o$ is selective,
then the subsequence $(x_n)_{n\in S}$ can be chosen so that $\o(S)=1$.
\end{thm}

\parit{Proof.}
Given $\eps>0$, 
let $S_\eps=\set{n\in\NN}{\dist{x_n}{x_\o}{}<\eps}$.

Note that $\o(S_\eps)=1$ for any $\eps>0$.
Since $\o$ is nonprincipal, the set $S_\eps$ is infinite.
Therefore we can choose an increasing sequence $n_\kay$
such that $n_\kay\in S_{\frac1\kay}$ for each $\kay\in \NN$.
Clearly $x_{n_\kay}\to x_\o$ as $\kay\to\infty$.

Now assume that $\o$ is selective.
Consider the sets
\begin{align*}
C_\kay&=\set{n\in\NN}{\tfrac1{\kay}<\dist{x_n}{x_\o}{}\le \tfrac1{\kay-1}},
\intertext{where we assume $\tfrac10=\infty$, and the set }
C_\infty&=\set{n\in\NN}{x_n=x_\o}.
\end{align*}

Note that $\o(C_\kay)=0$ for any $\kay\in\NN$.

If $\o(C_\infty)=1$, we can take the subsequence consisting of the $x_n$, $n\in C_\infty$.

Otherwise, discarding all empty sets among $C_\kay$ and $C_\infty$ gives a partition of $\NN$ into a countable collection of $\o$-small sets.
Since $\o$ is selective, we can choose a set $S\subset\NN$ such that
$S$ meets each set of the partition at one point and $\o(S)=1$.
Clearly the subsequence consisting of the $x_n$, $n\in S$
converges to $x_\o$ in the usual sense.
\qeds

The following proposition 
is analogous to the statement that any sequence in a compact metric space 
has a convergent subsequence;
it can be proved in the same way.

\begin{thm}{Proposition}\label{prop:ultra/compact}
Let $\spc{X}$ be a compact metric space.
Then
any sequence of points $x_n$ in $\spc{X}$ has a unique $\o$-limit $x_\o$.

In particular, a bounded sequence of real numbers has a unique $\o$-limit. 
\end{thm}

The following lemma is an ultralimit analog of the Cauchy convergence test.

\begin{thm}{Lemma}\label{lem:X-X^w}
Let $x_n$ be a sequence of points in a complete metric space~$\spc{X}$. 
If for each subsequence $y_n$ of $x_n$, 
the $\o$-limit 
\[y_\o=\lim_{n\to\o}y_{n}\in \spc{X}\]
is defined and does not depend on the choice of a subsequence, 
then the sequence $x_n$ converges in the usual sense.
\end{thm}

\parit{Proof.} Assume the contrary. 
Then for some $\eps>0$, there is a subsequence $y_n$ of $x_n$ such that $\dist{x_n}{y_n}{}\ge\eps$ for all $n$.

It follows that $\dist{x_\o}{y_\o}{}\ge \eps$, a contradiction.
\qeds

\begin{thm}{Exercise}\label{ex:linear}
Recall that $\ell^\infty$ denotes the space of bounded sequences of real numbers.
Show that there is a linear functional $L\:\ell^\infty\to\RR$ such that
for any sequence $\bm{s}=(s_1,s_2,\dots)\in S$ the image $L(\bm{s})$ is a partial limit of $s_1,s_2,\dots$
\end{thm}

\begin{thm}{Exercise}\label{ex:ultrakatetov+}
Suppose that $f\:\NN\to\NN$ is a map such that 
\[\lim_{n\to\o}x_n=\lim_{n\to\o}x_{f(n)}\]
for any bounded sequence $x_n$ of real numbers.
Show that $f(n)=n$ for $\o$-almost all $n\in\NN$.
\end{thm}

\section{Ultralimits of spaces}\label{sec:Ultralimit of spaces}

Fix a selective ultrafilter $\o$ on the set of natural numbers.

Let $\spc{X}_n$ be a sequence of metric spaces.
Consider all sequences
$x_n\in \spc{X}_n$.
On the set of all such sequences,
define a pseudometric  by the formula
\[\dist{(x_n)}{(y_n)}{}
=
\lim_{n\to\o} \dist{x_n}{y_n}{}.
\eqlbl{eq:olim-dist}\]
Note that the $\o$-limit on the right-hand side is always defined 
and takes  value in $[0,\infty]$. 

Let $\spc{X}_\o$ be the corresponding metric space; 
that is, the underlying set of $\spc{X}_\o$ is formed by equivalence  classes of sequences of points $x_n\in\spc{X}_n$ 
defined by the relation
\[(x_n)\sim(y_n)
\quad \iff\quad 
\lim_{n\to\o} \dist{x_n}{y_n}{}=0,\]
and the distance is defined as in \ref{eq:olim-dist}.

The space $\spc{X}_\o$ is called the \index{ultralimit of spaces}\emph{$\o$-limit} of $\spc{X}_n$.
Typically  $\spc{X}_\o$ will denote the  
$\o$-limit of a sequence $\spc{X}_n$;
we may also write  
\[\spc{X}_n\to\spc{X}_\o\quad \text{as}\quad  n\to\o\quad \text{or}\quad \spc{X}_\o=\lim_{n\to\o}\spc{X}_n.\]

Given a sequence  $x_n\in \spc{X}_n$,
we will denote by $x_\o$ its equivalence class, which is a point in $\spc{X}_\o$;
in this case, we may write
\[x_n\to x_\o \quad \text{as}\quad  n\to\o\quad \text{or}\quad x_\o=\lim_{n\to\o} x_n.\]

\begin{thm}{Observation}\label{obs:ultralimit-is-complete}
The $\o$-limit of any sequence of metric spaces is complete. 
\end{thm}

\parit{Proof.}
Let $\spc{X}_n$ be a sequence of metric spaces and $\spc{X}_n\to\spc{X}_\o$ as $n\to\o$.
Choose a Cauchy sequence $x_n$ in $\spc{X}_\o$.
Passing to a subsequence, we can assume that $\dist{x_k}{x_{m}}{\spc{X}_\o}<\tfrac1{k}$ for any $k<m$.

Let us choose points $x_{n,m}\in\spc{X}_n$ such that for any fixed $m$ we have $x_{n,m}\to x_m$ as $n\to\o$.
Note that for any $k<m$ the inequality $\dist{x_{n,k}}{x_{n,m}}{}<\tfrac1{k}$ holds for $\o$-almost all $n$.
It follows that we can choose a nested sequence of sets 
\[\NN= S_1\supset S_2\supset\dots\] 
such that 
\begin{itemize}
\item $\o(S_m)=1$ for each $m$, 
\item $\bigcap_m S_m=\emptyset$, and
\item $\dist{x_{n,k}}{x_{n,l}}{}<\tfrac1{k}$ for $k<l\le m$ and $n\in S_m$.
\end{itemize}

Consider the sequence $y_n=x_{n,m(n)}$, where $m(n)$ is the largest value such that $n\in S_{m(n)}$.
Denote by $y_\o\in \spc{X}_\o$ the $\o$-limit of $y_n$.

Observe that $|y_m-x_{n,m}|<\tfrac1{m}$ for $\o$-almost all $n$.
It follows that $|x_m-y_\o|\le \tfrac1{m}$ for any $m$.
Therefore, $x_n\to y_\o$ as $n\to \infty$.
That is, any Cauchy sequence in $\spc{X}_\o$ converges.
\qeds

\begin{thm}{Observation}\label{obs:ultralimit-is-geodesic}
The $\o$-limit of any sequence of length spaces is geodesic. 
\end{thm}

\parit{Proof.}
If $\spc{X}_n$ is a sequence of length spaces, then for any sequence of pairs $(x_n, y_n)$ in $\spc{X}_n$ there is a sequence of $\tfrac1n$-midpoints $z_n$.

Let $x_n\to x_\o$, $y_n\to y_\o$, and $z_n\to z_\o$ as $n\to \o$.
Note that $z_\o$ is a midpoint between $x_\o$ and $y_\o$ in $\spc{X}_\o$.

By Observation~\ref{obs:ultralimit-is-complete}, $\spc{X}_\o$ is complete.
Applying Lemma~\ref{lem:mid>geod} we obtain the statement.
\qeds

A geodesic space $\spc{T}$ is called a \index{metric tree}\emph{metric tree} if any pair of points in $\spc{T}$ are connected by a unique geodesic,
and the union of any two geodesics $[xy]_{\spc{T}}$, and $[yz]_{\spc{T}}$ contain the geodesic $[xz]_{\spc{T}}$.
The latter means that any triangle in $\spc{T}$ is a tripod;
that is, for any three points $x$, $y$, and $z$ there is a point $p$ such that 
\[[xy]\cup[yz]\cup[zx]=[px]\cup[py]\cup[pz].\]

\begin{thm}{Exercise}\label{ex:Asym(Lob)}
Let $\spc{T}$ be a metric component of the ultralimit of $\Lob2n$ as $n\to\o$.

\begin{subthm}{ex:Asym(Lob):metric-tree}
Show that $\spc{T}$ is a complete metric tree.
\end{subthm}

\begin{subthm}{ex:Asym(Lob):homogeneous}
Show that $\spc{T}$ is homogeneous; that is, given two points $s,t\in \spc{T}$ there is an isometry of $\spc{T}$ that maps $s$ to $t$.
\end{subthm}

\begin{subthm}{ex:Asym(Lob):continuum}
Show that $\spc{T}$ has \index{degree}\emph{continuum degree} at any point;
that is, for any point $t\in \spc{T}$ the set of connected components of the complement $\spc{T}\setminus\{t\}$ has cardinality continuum.
\end{subthm}

\end{thm}

\parbf{Ultrapower.} If all the metric spaces in a sequence are identical, $\spc{X}_n\z=\spc{X}$, 
the $\o$-limit 
$\lim_{n\to\o}\spc{X}_n$
is denoted by $\spc{X}^\o$
and called the \index{ultrapower} $\o$-power of $\spc{X}$.
 
By Theorem~\ref{thm:ultra-GH},
there is a distance-preserving map
$\iota\:\spc{X}\hookrightarrow \spc{X}^\o$, where $\iota(y)$ is the equivalence class of the constant sequence $y_n=y$. 

The image $\iota(\spc{X})$ might be a proper subset of $\spc{X}^\o$.
For example, $\RR^\o$ has pairs of points at distance $\infty$ from each other, while each metric component of $\RR^\o$ is isometric to $\RR$.

According to Theorem~\ref{thm:ultra-GH}, 
if $\spc{X}$ is compact then $\iota(\spc{X})=\spc{X}^\o$;
in particular, $\spc{X}^\o$ is isometric to $\spc{X}$.
If $\spc{X}$ is proper, then $\iota(\spc{X})$ forms a metric component of~$\spc{X}^\o$.

The embedding $\iota$ allows us to treat $\spc{X}$ as a subset of its ultrapower~$\spc{X}^\o$. 

\begin{thm}{Observation}\label{obs:ultrapower-is-geodesic}
Let $\spc{X}$ be a complete metric space. 
Then $\spc{X}^\o$ is a geodesic space if and only if $\spc{X}$ is a length space.
\end{thm}

\parit{Proof.}
Assume $\spc{X}^\o$ is a geodesic space.
Then any pair of points $x,y\in \spc{X}$ has a midpoint $z_\o\in\spc{X}^\o$.
Fix a sequence of points $z_n\in  \spc{X}$ such that $z_n\to z_\o$ as $n\to \o$.

Note that 
$\dist{x}{z_n}{\spc{X}}\to \tfrac12\cdot \dist{x}{y}{\spc{X}}$
and 
$\dist{y}{z_n}{\spc{X}}\to \tfrac12\cdot \dist{x}{y}{\spc{X}}$
as 
$n\to\o$.
In particular, for any $\eps>0$, the point $z_n$ is an $\eps$-midpoint between $x$ and $y$ for $\o$-almost all $n$.
It remains to apply Lemma~\ref{lem:mid>geod}.

The if part follows from Observation~\ref{obs:ultralimit-is-geodesic}.
\qeds

Note that the proof above together with Lemma~\ref{lem:X-X^w} imply the following:

\begin{thm}{Corollary}\label{cor:two-geodesics-in-ultrapower}
Assume $\spc{X}$ is a complete length space 
and $p,q\in\spc{X}$ cannot be joined by a geodesic in $\spc{X}$.  
Then there are at least continuum distinct geodesics between $p$ and $q$ 
in the ultrapower $\spc{X}^\o$.
\end{thm}

\begin{thm}{Exercise}\label{ex:isom-ultrapower}
Let $\spc{X}$ be a countable set with discrete metric;
that is $\dist{x}{y}{\spc{X}}=1$ if $x\ne y$.
Show that 

\begin{subthm}{ex:isom-ultrapower:no}
$\spc{X}^\o$ is not isometric to $\spc{X}$.
\end{subthm}

\begin{subthm}{ex:isom-ultrapower:yes}
$\spc{X}^\o$ is  isometric to $(\spc{X}^\o)^\o$.
\end{subthm}

\end{thm}

\begin{thm}{Exercise}\label{ex:ultrapower(ultrapower)}
Given a nonprincipal ultrafilter $\o$, construct an ultrafilter $\o_1$ such that 
\[\spc{X}^{\o_1}\iso(\spc{X}^\o)^\o\]
for any metric space~$\spc{X}$.

\end{thm}

\begin{thm}{Exercise}\label{ex:notproper-limit}
Construct a proper metric space $\spc{X}$ such that $\spc{X}^\o$ is not proper;
that is, there is a point $p\in \spc{X}^\o$ and $R<\infty$ such that the closed ball $\cBall[p,R]_{\spc{X}^\o}$ is not compact.
\end{thm}

{\sloppy

\section{Ultralimits of sets}

Let $\spc{X}_n$ be a sequence of metric spaces and $\spc{X}_n\to \spc{X}_\o$
as $n\to \o$.

For a sequence of sets $\Omega_n\subset \spc{X}_n$,
consider the maximal set $\Omega_\o\subset \spc{X}_\o$ such that 
for any $x_\o\in\Omega_\o$ and any sequence $x_n\in\spc{X}_n$ such that $x_n\to x_\o$ as $n\to \o$, we have $x_n\in\Omega_n$ for $\o$-almost all $n$.

The set $\Omega_\o$ is called the  \index{ultralimit of sets}\emph{open $\o$-limit} of $\Omega_n$;
we could also write  $\Omega_n\to \Omega_\o$ as $n\to\o$ or $\Omega_\o=\lim_{n\to\o}\Omega_n$. 

{\sloppy

Applying Observation~\ref{obs:ultralimit-is-complete} to the sequence of complements $\spc{X}_n\setminus \Omega_n$, we see that $\Omega_\o$ is open for any sequence $\Omega_n$.

This definition can be applied to arbitrary sequences of sets,
but we will apply it only for sequences of open sets.

}

\section{Ultralimits of functions}\label{sec:Ultralimits of functions}

Recall that a family of submaps (see section \ref{sec:submaps}) between metric spaces $\{f_\alpha\co \spc{X}\subto\spc{Y}\}_{\alpha\in\mathcal A}$ is called \index{equicontinuous family}\emph{equicontinuous} if for any $\eps>0$ there is $\delta>0$ such that for any $p,q\in\spc{X}$ with $\dist{p}{q}{}<\delta$ and any $\alpha\in\mathcal A$ we have $\dist{f_\alpha(p)}{f_\alpha(q)}{}<\eps$.

Let $f_n\:\spc{X}_n\subto\RR$ be a sequence of subfunctions.

Set $\Omega_n=\Dom f_n$.
Consider the open $\o$-limit set $\Omega_\o\subset \spc{X}_\o$ of $\Omega_n$.

Assume there is a subfunction $f_\o\:\spc{X}_\o\subto\RR$
that satisfies the following conditions: 
(1) $\Dom f_\o=\Omega_\o$, (2) if $x_n\to x_\o\in \Omega_\o$ for a sequence of points $x_n\in\spc{X}_n$, then $f_n(x_n)\to f_\o(x_\o)$ as $n\to\o$.
In this case, the subfunction $f_\o\:\spc{X}_\o\to\RR$ is said to be the $\o$-limit of $f_n\:\spc{X}_n\to\RR$.

The following lemma gives a mild condition on a sequence of functions $f_n$
guaranteeing the existence of its $\o$-limit.

\begin{thm}{Lemma}
Let $\spc{X}_n$ be a sequence of metric spaces
and $f_n\:\spc{X}_n\subto\RR$ be a sequence of subfunctions.

Assume that  for any positive integer $\kay$, there is an open set $\Omega_n(\kay)\subset \Dom f_n$
such that the restrictions $f_n|_{\Omega_n(\kay)}$ are uniformly bounded and equicontinuous
and the open $\o$-limit of $\Omega_n(n)$ coincides with the open $\o$-limit of $\Dom f_n$.
Then the $\o$-limit $f_\o$ of $f_n$ is defined;
moreover $f_\o$ is a continuous subfunction.

In particular, if the functions $f_n$ are uniformly bounded and equicontinuous, then its $\o$-limit $f_\o$ is defined, bounded and uniformly continuous.
\end{thm}

The proof is straightforward.

{\sloppy

\begin{thm}{Exercise}\label{ex:nonconvex-limit}
Construct a sequence of compact length spaces 
$\spc{X}_n$  
with a converging sequence of $\Lip$-Lipschitz concave (see Definition \ref{def:lam-convex}) functions $f_n\:\spc{X}_n\to\RR$ such that
the $\o$-limit $\spc{X}_\o$ of $\spc{X}_n$ is compact
and the $\o$-limit $f_\o\:\spc{X}_\o\to\RR$ of $f_n$ is not concave.
\end{thm}

}

If $f\:\spc{X}\subto\RR$ is a subfunction, 
the $\o$-limit of the constant sequence $f_n=f$ is called the $\o$-power of $f$ and is denoted by $f^\o$.
So
\[f^\o\:\spc{X}\subto\RR,\quad f^\o(x_\o)=\lim_{n\to\o} f(x_n).\]

Evidently, if $f^\o$ is defined, then $f$ is continuous.

Recall that we treat $\spc{X}$ as a subset of its $\o$-power $\spc{X}^\o$.
Note that $\Dom f=\spc{X}\cap \Dom f^\o$.
Moreover, 
if $\oBall(x,\eps)_{\spc{X}}\subset \Dom f$
then $\oBall(x,\eps)_{\spc{X}^\o}\subset \Dom f^\o$.

\chapter{Space of spaces}

In this chapter we discuss the
Gromov--Hausdorff convergence of metric spaces.

To the best of our knowledge, Hausdorff convergence of subsets of a fixed metric space was first introduced by Felix Hausdorff \cite{hausdorff}, 
and a couple of years later an equivalent definition was given by Wilhelm Blaschke \cite{blaschke}.
A further refinement of this definition was introduced by Zdeněk Frolík \cite{frolik}
and then rediscovered by Robert Wijsman \cite{wijsman}.
However this refinement was a step in the direction of the so-called {}\emph{closed convergence} introduced by Hausdorff in the original book. 
For that reason we call it Hausdorff convergence
instead of
\emph{Hausdorff--Blashcke--Frol\'{\i}k--Wijsman convergence}.

Gromov--Hausdorff convergence was first introduced by David Edwards \cite{edwards}
and rediscovered later by Michael Gromov \cite{gromov-polynomial-growth}.
It was an essential tool in Gromov's proof that any group of polynomial growth has  a nilpotent subgroup of finite index.
Other versions of convergence of metric spaces
were considered earlier, but each time
the definition was limited to very specific types of problems.

The definition of Gromov--Hausdorff convergence of metric spaces uses 
the notion of Hausdorff convergence.
Gromov--Hausdorff convergence means that a sequence of metric spaces admits a sequence of distance-preserving embeddings into a common ambient metric space so that their images converge in the Hausdorff sense.
Our definition of  Gromov--Hausdorff convergence and  Gromov--Hausdorff distance differ somewhat from the standard definition.

\section{Convergence of subsets}

Let $\spc{X}$ be a metric space and $A\subset \spc{X}$.
Recall that the distance from $A$ to a point $x$ in $\spc{X}$
is given by
$$\distfun Ax\df\inf\set{\dist{a}{x}{}}{a\in A}.$$
By this definition, we have $\distfun{\emptyset}x=\infty$ for any $x$.

\begin{thm}{Definition of Hausdorff convergence}\label{def:hausdorff-coverge}
Given a sequence of closed sets $A_n$ in a metric space $\spc{X}$, 
a closed set $A_\infty\subset \spc{X}$ is called the Hausdorff limit of $A_n$,
briefly $A_n\Hto A_\infty$, if 
$$\distfun{A_n}x\to\distfun{A_\infty}x\quad \text{as}\quad n\to\infty$$
for any fixed $x\in \spc{X}$.

In this case the, sequence of closed sets $A_n$ is said to be {}\emph{converging} or \index{Hausdorff convergence}\emph{converging in the sense of Hausdorff}.
\end{thm}

\begin{thm}{Selection theorem}
Let $\spc{X}$ be a proper metric  space
and $A_n$ be a sequence of closed sets in $\spc{X}$.
Then  $A_n$ has a converging subsequence in the sense of Hausdorff.
\end{thm}

\parit{Proof.}
Since $\spc{X}$ is proper,
we can choose a countable dense set $\{x_1,x_2, \dots\}$ in $\spc{X}$.

If the sequence $a_n=\distfun{A_n}x_\kay$ is unbounded for some $\kay$,
then we can pass to a subsequence of $A_n$ such that 
$\distfun{A_n}x_\kay\to \infty$ as $n\to\infty$ for \emph{any}~$\kay$.
The obtained sequence converges to the empty set.

Now suppose that $a_n=\distfun{A_n}x_\kay$ is bounded for each $\kay$. 
In this case, passing to a subsequence of $A_n$,
we can assume that $\distfun{A_n}x_\kay$ converges as $n\to\infty$ for any fixed $\kay$.

Note that for each $n$, the function $\distfun{A_n}\:\spc{X}\to\RR$ is 1-Lipschitz and nonnegative.
Therefore the sequence $\distfun{A_n}$ converges pointwise to a 1-Lipschitz nonnegative function $f\:\spc{X}\to\RR$.

Set $A_\infty=f^{-1}(0)$.
Since $f$ is 1-Lipschitz, 
$\distfun{A_\infty}y\ge f(y)$ for any $y\in \spc{X}$.
It remains to show that $\distfun{A_\infty}y\le f(y)$ for any $y$.

Assume the  contrary,
that is, $f(z)<R<\distfun{A_\infty}z$ for $z\in \spc{X}$ and $R>0$.
Then for any sufficiently large $n$ there is a point $z_n\in A_n$ such that
$\dist{x}{z_n}{}\le R$.
Since $\spc{X}$ is proper, we can pass to a partial limit $z_\infty$ of $z_n$ as $n\to\infty$.

Clearly  $f(z_\infty)=0$, that is, $z_\infty\in A_\infty$.
On the other hand, 
\[\distfun{A_\infty}y\le\dist{z_\infty}{y}{}\le R<\distfun{A_\infty}y,\] 
a contradiction.
\qeds

\section{Convergence of spaces}

\begin{thm}{Definition}\label{def:comp-metr}
Let $\set{\spc{X}_\alpha}{\alpha\in\IndexSet}$ be a set of metric spaces.
A metric space $\bm{X}$
is called a \index{common space}\emph{common space} of $\set{\spc{X}_\alpha}{\alpha\in\IndexSet}$ if its underlying set is formed by the disjoint union
$$\bigsqcup_{\alpha\in\IndexSet} \spc{X}_\alpha$$ 
and each inclusion $\iota_\alpha\:\spc{X}_\alpha\hookrightarrow\bm{X}$
is distance-preserving.
\end{thm}

\begin{thm}{Definition}\label{def:GH}
Let $\bm{X}$ be a common space for proper metric spaces
$\spc{X}_1,\spc{X}_2,\dots$, and $\spc{X}_\infty$.
Assume that $\spc{X}_n$ forms an open set in $\bm{X}$ for each $n<\infty$ and 
$\spc{X}_n\Hto \spc{X}_\infty$ in $\bm{X}$ as $n\to\infty$.

Then the topology $\GH$ of $\bm{X}$ is called a \index{Gromov--Hausdorff convergence}\emph{Gromov--Hausdorff convergence}
and we write $\spc{X}_n\GHto \spc{X}_\infty$ or $\spc{X}_n\xto{\GH} \spc{X}_\infty$;
the latter notation is used if we need to consider  the specific Gromov--Hausdorff convergence $\GH$.
The space $\spc{X}_\infty$ is called the {}\emph{limit space} of the sequence $\spc{X}_n$ along $\GH$.
\end{thm}

When we write $\spc{X}_n\GHto \spc{X}_\infty$ we mean that we made a choice of a Gromov--Hausdorff convergence.

Note that for a fixed sequence $\spc{X}_n$ of metric spaces, one may construct different Gromov--Hausdorff convergences, say $\spc{X}_n\xto{\GH} \spc{X}_\infty$ and $\spc{X}_n\xto{\GH'} \spc{X}_\infty'$,  and their limit spaces $\spc{X}_\infty$ and $\spc{X}_\infty'$ need not be isometric to each other. 
For example, for the constant sequence $\spc{X}_n\iso\RR_{\ge0}$, 
one may take $\spc{X}_\infty\iso\RR_{\ge0}$.
In this case, a point in the disjoint space $\bm{X}$ can be regarded as a pair $(x,n)\in \RR_{\ge0}\times (\ZZ_>\cup \{\infty\})$ 
and the metric on $\bm{X}$ can be defined by
$$\dist{(x,n)}{(y,m)}{\bm{X}}\df|\tfrac1n+\tfrac1m|+|x-y|,$$
where we assume that $0=\tfrac1\infty$.
On the other hand, one can take $\spc{X}_\infty'\iso\RR$,
and consider the metric
\begin{align*}
\dist{(x,n)}{(y,m)}{\bm{X}'}
&=|\tfrac1n-\tfrac1m|+|(x-n)-(y-m)|,
\\
\dist{(x,n)}{(y,\infty)}{\bm{X}'}
&=\tfrac1n+|(x-n)-y|,
\\
\dist{(x,\infty)}{(y,\infty)}{\bm{X}'}
&=|x-y|,
\end{align*}
where $n, m<\infty$.

\begin{thm}{Induced convergences}
Suppose $\spc{X}_n\xto{\GH}\spc{X}_\infty$
as in Definition \ref{def:GH},
and $\iota_n\:\spc{X}_n\hookrightarrow\bm{X}$, $\iota_\infty\:\spc{X}_\infty\hookrightarrow\bm{X}$ are the corresponding inclusions.

\begin{subthm}{}
A sequence of points $x_n\in\spc{X}_n$ converges to $x_\infty\in\spc{X}_\infty$ (briefly, $x_n\to x_\infty$ or $x_n\xto{\GH} x_\infty$) 
if $\dist{x_n}{x_\infty}{\bm{X}}\to 0$.
\end{subthm}

\begin{subthm}{}
A sequence of closed sets 
$\mathfrak{C}_n\subset \spc{X}_n$ 
converges to a closed  set 
$\mathfrak{C}_\infty\subset \spc{X}_\infty$ (briefly, $\mathfrak{C}_n\to \mathfrak{C}_\infty$ or $\mathfrak{C}_n\xto{\GH} \mathfrak{C}_\infty$)
if $\mathfrak{C}_n\Hto\mathfrak{C}_\infty$ as subsets of $\bm{X}$.
\end{subthm}

\begin{subthm}{}
A sequence of open sets $\Omega_n\subset \spc{X}_n$ 
converges to an open set $\Omega_\infty\subset \spc{X}_\infty$
(briefly, $\Omega_n\to \Omega_\infty$ 
or $\Omega_n\xto{\GH} \Omega_\infty$)
if the complements $\spc{X}_n\setminus \Omega_n$ converge to the complement $\spc{X}_\infty\setminus \Omega_\infty$ as closed sets.
\end{subthm}

\begin{subthm}{} Let $\spc{X}_n\xto{\GH} \spc{X}_\infty$ and $\spc{Y}_n\xto{\theta} \spc{Y}_\infty$. 
A sequence of submaps (where a submap is a map defined on a subset; see Section~\ref{sec:submaps}) $\map_n\:\spc{X}_n\subto \spc{Y}_n$ converges to a submap $\map_\infty\:\spc{X}_\infty\subto \spc{Y}_\infty$ if the following conditions holds
\begin{itemize}
\item $\Dom\map_n\to \Dom\map_\infty$ as a sequence of open sets.

\item for any $x_\infty\in \Dom \map_\infty$ and any sequence $x_n\in \spc{X}_n$ such that $x_n\to x_\infty$, we have
\[\spc{Y}_n\ni \map _n(x_n)\xto\theta \map_\infty(x_\infty)\in\spc{Y}_\infty\] 
as $n\to\infty$.
\end{itemize}
\end{subthm}

\begin{subthm}{} Given a sequence of measures $\mu_n$ on $\spc{X}_n$,
we say that $\mu_n$ weakly converges to a measure $\mu_\infty$ on $\spc{X}_\infty$ 
(briefly, $\mu_n\xto{}
\mu_\infty$ or $\mu_n\xto{\GH}
\mu_\infty$) 
if the pushforward measures of $\mu_n$
weakly converge to the pushforward measure of $\mu_\infty$.

In other words, 
if for any continuous function $\phi\:\bm{X}\to\RR$ with a compact support, we have 
\[\int\limits_{\spc{X}_n} 
\phi\circ\iota_n
\cdot
\mu_n
\to 
\int\limits_{\spc{X}_\infty}
\phi\circ\iota_\infty
\cdot\mu_\infty\]
as $n\to\infty$.
\end{subthm}
\end{thm}

\parbf{Liftings.}
Given a Gromov--Hausdorff convergence 
$\spc{X}_n\GHto \spc{X}_\infty$
and a point $p_\infty\in\spc{X}_\infty$, any sequence of points $p_n\in\spc{X}_n$ such that $p_n\GHto p_\infty$  will be called a \index{Gromov--Hausdorff convergence!lifting of a point}\emph{lifting} of $p_\infty$.
The point $p_n\in \spc{X}_n$ will be called a {}\emph{lifting} of $p_\infty$ to $\spc{X}_n$.
We will also say that $\distfun{p_n}{}{}\:\spc{X}_n\to \RR$ 
is a \index{Gromov--Hausdorff convergence!lifting of a distance function}\emph{lifting} 
of the distance function $\distfun{p_\infty}{}{}\:\spc{X_\infty}\to \RR$.
Clearly $\distfun{p_n}{}{}\GHto\distfun{p_\infty}{}{}$.

Note that liftings are not uniquely defined.

Similarly, we may refer to liftings of a point array
$\bm{p}_\infty\z =(p_\infty^1,p_\infty^2,\dots,p_\infty^\kay)$
and of the corresponding distance map 
$\distfun{\bm{p}_\infty}{}{}\:\spc{X}_\infty\to\RR^\kay$,
$$\distfun{\bm{p}_\infty}{}{}\:x\mapsto(\dist{p_\infty^1}{x}{},\dist{p_\infty^2}{x}{},\dots,\dist{p_\infty^k}{x}{}).$$

\section{Gromov's selection theorem}

\begin{thm}{Gromov's selection theorem}\label{thm:gromov-selection}
Let $\spc{X}_n$ be a sequence of proper metric spaces 
with marked points $x_n\in \spc{X}_n$.
Assume that for any $R>0$, $ \eps>0$, there is $N=N(R,\eps)\in\ZZ_{>0}$ 
such that for each $n$
the ball $\cBall[x_n,R]\subset \spc{X}_n$ admits a finite $\eps$-net with at most $N$ points.
Then a subsequence of $\spc{X}_n$ admits a Gromov--Hausdorff convergence 
such that the sequence of marked points $x_n\in\spc{X}_n$ 
converges.
\end{thm}

\parit{Proof.}
By the main assumption, there is a sequence of integers $M_1\z<M_2<\dots$
such that in each space $\spc{X}_n$
there is a sequence of points $z_{i,n}\in\spc{X}_n$ for which  
\[\dist{z_{i,n}}{x_n}{\spc{X}_n}\le \kay+1\quad \text{if}\quad i\le M_\kay\]
and $\{z_{1,n},\dots,z_{M_\kay,n}\}$ is a $\tfrac1\kay$-net in $\cBall[x_n,\kay]_{\spc{X}_n}$.

Passing to a subsequence, we may assume that the sequence \[\ell_n=\dist{z_{i,n}}{z_{j,n}}{\spc{X}_n}\] 
converges for any $i$ and $j$.

Let us consider a countable set of points $\spc{W}=\{w_1,w_2,\dots\}$
equipped with the pseudometric defined by 
\[\dist{w_i}{w_j}{\spc{W}}
=
\lim_{n\to\infty}\dist{z_{i,n}}{z_{j,n}}{\spc{X}_n}.\]
Let $\hat{\spc{W}}$ be the metric space corresponding to $\spc{W}$.
Denote by
$\spc{X}_\infty$ the completion of $\hat{\spc{W}}$.

It remains to construct a metric on the disjoint union of \[\bm{X}=\spc{X}_\infty\sqcup\spc{X}_1\sqcup\spc{X}_2\sqcup\dots\] 
satisfying definitions \ref{def:comp-metr} and \ref{def:GH}.

Such a  metric can be defined as follows.
Fix a sequence $\eps_\kay\to0+$
and let $N_\kay$ be the minimal integer such that
\[\dist{w_i}{w_j}{\spc{W}}
\lege
\dist{z_{i,n}}{z_{j,n}}{\spc{X}_n}\pm\eps_\kay
\]
if $i,j\le N_\kay$ and $n\ge N_\kay$. 
Let us equip $\bm{X}$ with the maximal metric such that all the inclusions $\iota_n\:\spc{X}_n\to\bm X$  and $\iota_\infty\:\spc{X}_\infty\to\bm X$ are isometric and 
$
\dist{z_{i,n}}{w_i}{}\le \eps_\kay
$
for $i\le N_\kay$ and $n\ge N_\kay$.
It is easy to verify  that such a metric on $\bm X$ satisfies  \ref{def:comp-metr} and \ref{def:GH}.
\qeds

\section{Convergence of compact spaces.}

\begin{thm}{Definition}
Let $\spc{X}$ and $\spc{Y}$ be metric spaces. A map $f\:\spc{X}\to\spc{Y}$
is called an \index{isometry!$\eps$-isometry}\emph{$\eps$-isometry}
if the following two conditions hold:
\begin{subthm}{}
$\Im f$ is an $\eps$-net in $\spc{Y}$.
\end{subthm}

\begin{subthm}{}
$\bigl|\dist{f(x)}{f(x')}{\spc{Y}}-\dist{x}{x'}{\spc{X}}\bigr|\le \eps$ for any $x,x'\in\spc{X}$.
\end{subthm}

\end{thm}

\begin{thm}{Lemma}\label{lem:almost-isom}
Let $\spc{X}_1,\spc{X}_2,\dots$, and  $\spc{X}_\infty$ be metric spaces and $\eps_n\to\0+$ as $n\to\infty$.
Suppose that either 

\begin{subthm}{}\label{lem:almost-isom-a}
for each $n$ there is an $\eps_n$-isometry $f_n\:\spc{X}_n\to\spc{X}_\infty$, or
\end{subthm}

\begin{subthm}{}\label{lem:almost-isom-b}
for each $n$ there is an $\eps_n$-isometry $h_n\:\spc{X}_\infty\to\spc{X}_n$.
\end{subthm}

Then there is a Gromov--Hausdorff convergence $\spc{X}_n\GHto \spc{X}_\infty$.
\end{thm}

\parit{Proof.}
To prove part \eqref{lem:almost-isom-a} let us construct a common space $\bm{X}$ for the spaces $\spc{X}_1,\spc{X}_2,\dots$, and $\spc{X}_\infty$
by taking the metric $\rho$ on the disjoint union $\spc{X}_\infty\sqcup\spc{X}_1\sqcup\spc{X}_2\sqcup\dots$ that is defined the following way:
\begin{align*}
\rho(x_n,y_n)&=\dist{x_n}{y_n}{\spc{X}_n},\quad \rho(x_\infty,y_\infty)=\dist{x_\infty}{y_\infty}{\spc{X}_\infty},
\\
\rho(x_n,x_\infty)&=\inf\set{\dist{x_n}{y_n}{\spc{X}_n}+\eps_n+\dist{x_\infty}{f(y_n)}{\spc{X}_\infty}}{{y_n}\in \spc{X}_n},
\\
\rho(x_n,x_m)&=\inf\set{\rho(x_n,y_\infty)+\rho(x_m,y_\infty)}{y_\infty\in\spc{X}_\infty},
\end{align*}
where we assume that $x_m\in \spc{X}_m$, $x_n\in \spc{X}_n$, and $x_\infty\in \spc{X}_\infty$. 

It remains to observe that $\rho$ is indeed a metric and 
$\spc{X}_n\Hto \spc{X}_\infty$ in~$\bm{X}$.

The proof of the second part is analogous; one only needs to change one line in the definition of $\rho$ to the following:
\[\rho(x_n,x_\infty)=\inf\set{\dist{x_n}{h(y_\infty)}{\spc{X}_n}+\eps_n+\dist{x_\infty}{y_\infty}{\spc{X}_\infty}}{{y_\infty}\in \spc{X}_\infty}.\]
\qedsf

\begin{thm}{Definition}\label{def: inequality-of-spaces}
 Given two compact spaces $\spc{X}$ and $\spc{Y}$, we will write 
\begin{itemize}
\item $\spc{X}\le \spc{Y}$ if there is a noncontracting map $\map\:\spc{X}\to \spc{Y}$.
\item $\spc{X}\le \spc{Y}+\eps$ if there is a map $\map\:\spc{X}\to \spc{Y}$ such that for any $x,x'\in \spc{X}$ we have
\[\dist{x}{x'}{}\le \dist{\map(x)}{\map(x')}{}+\eps.\]
\end{itemize}

\end{thm}

\begin{thm}{Lemma}\label{lem:>=-isometry}
Let $\spc{X}$ and $\spc{Y}$ be two metric spaces and $\spc{X}$ be compact. Then
\[
\spc{X}\ge\spc{Y}\ge\spc{X}
\quad \iff\quad 
\spc{X}\iso\spc{Y}.
\]

\end{thm}

The following proof was suggested by Travis Morrison.

\parit{Proof.}
Let $f\: \spc{X} \to \spc{Y}$ 
and $g\: \spc{Y} \to \spc{X}$ be noncontracting mappings.
It is sufficient to prove that $h  = g\circ f\:\spc{X}\to \spc{X}$ is an isometry. 

Given any pair of points $x,y\in \spc{X}$, 
let $x_n\z=h^{\circ n}(x)$ and $y_n\z=h^{\circ n}(y)$.
Since $\spc{X}$ is compact, one can choose an increasing sequence of integers $n_\kay$
such that both sequences $x_{n_i}$ and $y_{n_i}$
converge.
In particular, both of these sequences  are 
Cauchy, 
that is,
\[
\dist{x_{n_i}}{x_{n_j}}{},\dist{y_{n_i}}{y_{n_j}}{}\to 0
\]
as $\min\{i,j\}\to\infty$.
Since $h$ is noncontracting, we have
\[
\dist{x}{x_{|n_i-n_j|}}{}\le \dist{x_{n_i}}{x_{n_j}}{}.
\]
It follows that  
there is a sequence $m_i\to\infty$ such that
\[
x_{m_i}\to x\quad \text{and}\quad y_{m_i}\to y\quad \text{as}\quad i\to\infty.
\eqlbl{eq:x_l->x}
\]

Let $\ell_n=\dist{x_n}{y_n}{}$.
Since $h$ is noncontracting, the sequence $\ell_n$ is nondecreasing.
On the other hand, 
from \ref{eq:x_l->x} it follows that $\ell_{m_i}\to\dist{x}{y}{}=\ell_0$ as $m_i\to\infty$;
that is, $\ell_n$ is a constant sequence.
In particular, $\ell_0=\ell_1$ for any $x$ and $y$ in $\spc{X}$,
so $h$ is a distance-preserving map.

Thus $h(\spc{X})$ is isometric to $\spc{X}$.
From \ref{eq:x_l->x}, $h(\spc{X})$ is everywhere dense.
Since $\spc{X}$ is compact, $h(\spc{X})=\spc{X}$.
\qeds

The \index{Gromov--Hausdorff distance}\emph{Gromov--Hausdorff distance} between isometry classes of compact metric spaces $\spc{X}$ and $\spc{Y}$, is defined by
\[\GHdist(\spc{X},\spc{Y})
\df
\inf\set{\eps>0}{\spc{X}\le \spc{Y}+\eps\ \text{and}\ \spc{Y}\le \spc{X}+\eps}.
\]
The Gromov--Hausdorff distance turns the set of all isometry classes of compact metric spaces into a metric space.
The following theorem shows that convergence in this space coincides with the Gromov--Hausdorff convergence defined above.

\begin{thm}{Theorem} Let $\spc{X}_1,\spc{X}_2,\dots$, and $\spc{X}_\infty$ be compact metric spaces.
Then there is a convergence $\spc{X}_n\GHto \spc{X}_\infty$ if and only if
$\GHdist(\spc{X}_n,\spc{X}_\infty)\to 0$ as $n\to\infty$.

\end{thm}

\parit{Proof; if part.}
Suppose $a_n\:\spc{X}_\infty\to \spc{X}_n$
and $b_n\:\spc{X}_n\to \spc{X}_\infty$ are sequences of maps such that
\[\dist{a_n(x)}{a_n(y)}{\spc{X}_\infty}
\ge
\dist{x}{y}{\spc{X}_n}-\delta_n,\]
\[\dist{b_n(v)}{b_n(w)}{\spc{X}_n}
\ge
\dist{v}{w}{\spc{X}_\infty}-\delta_n\]
for any $x,y\in \spc{X}_n$, $v,w\in\spc{X}_\infty$, and a sequence $\delta_n\to0+$.

Fix $\eps>0$ and choose a maximal $\eps$-packing $\{x^1,x^2,\dots,x^\kay\}$ in $\spc{X}_\infty$ such that 
$\sum_{i<j}\dist{x^i}{x^j}{}$ is maximal.
Note that 
\[\dist{a_n\circ b_n(x^i)}{a_n\circ b_n(x^j)}{}\ge\dist{x^i}{x^j}{}-2\cdot\delta_n.\]
Since $\sum_{i<j}\dist{x^i}{x^j}{}$ is maximal, 
\[\dist{a_n\circ b_n(x^i)}{a_n\circ b_n(x^j)}{}\to\dist{x^i}{x^j}{}\]
for all $i$ and $j$ as $n\to\infty$.
For all large $n$,
we have $2\cdot\delta_n<\dist{x^i}{x^j}{}-\eps$,
and so 
\[\dist{b_n(x^i)}{b_n(x^j)}{\spc{X}_n}>\eps
\quad\text{and}\quad
\dist{a_n\circ b_n(x^i)}{a_n\circ b_n(x^j)}{\spc{X}_\infty}>\eps\] 
for all $i\ne j$.
Therefore for each large $n$, 
the set $\{a_n\circ b_n(x^i)\}$ is a maximal $\eps$-packing and hence an $\eps$-net in $\spc{X}_\infty$.

Since $\{a_n\circ b_n(x^i)\}$ is an $\eps$-net in $\spc{X}_\infty$, we have that 
for any $y_n\in\spc{X}_n$ there is $x^i$ such that $\dist{a_n\circ b_n(x^i)}{a_n(y_n)}{}<\eps$.
Thus $\dist{b_n(x^i)}{y_n}{}<\eps+\delta_n$, 
that is, $\{b_n(x^i)\}$ is a $(\eps+\delta_n)$-net in $\spc{X}_n$.

Given $y\in \spc{X}_n$, choose $x^i$ so that $\dist{b_n(x^i)}{y_n}{}<\eps+\delta_n$ and define $h_n(y)=a_n\circ b_n(x^i)$.
Observe that $h_n$ is a $3\cdot\eps$-isometry for all large $n$.
Since $\eps>0$ is arbitrary, there is a sequence of $\eps_n$-isometries $\spc{X}_n\to\spc{X}_\infty$ such that $\eps_n\to\0+$ as $n\to\infty$.
It remains to apply \ref{lem:almost-isom}.

\parit{Only-if part.}
Assume $\spc{X}_n\xto{\GH}\spc{X}_\infty$.
Fix $\eps>0$, and choose a maximal $\eps$-packing $\{x^1,x^2,\dots,x^\kay\}$ in $\spc{X}_\infty$.
For each $x^i$, 
choose a sequence $x^i_n\in\spc{X}_n$ 
such that $x^i_n\to x^i$.
Define a map $a_n\: \spc{X}_n\to\spc{X}_\infty$
such that $a_n(x^i_n)=x_n$.
Note that for all large $n$, we have $\dist{x^i_n}{x^j_n}{}>\eps$.
For each point $z\in \spc{X}_\infty$, choose $x^i$ so that $\dist{z}{x^i}{}<\eps$. 
Define a map $b_n\:\spc{X}_\infty\to\spc{X}_n$ by setting 
$b_n(z)=x^i_n$.
Observe that 
\[\dist{b_n(y)}{b_n(z)}{\spc{X}_n}+3\cdot\eps>\dist{y}{z}{\spc{X}_\infty}\]
for all large $n$.

In the same way we can construct a map $a_n\:\spc{X}_n\to\spc{X}_\infty$ such that 
\[\dist{a_n(y)}{a_n(z)}{\spc{X}_\infty}+3\cdot\eps>\dist{y}{z}{\spc{X}_n}.\]
Hence $\GHdist(\spc{X}_n,\spc{X}_\infty)\to 0$ as $n\to \infty$.
\qeds

The following theorem states that the isometry class of a Gromov--Hausdorff limit is uniquely defined if it is compact. 

\begin{thm}{Theorem}\label{thm:GH-compact}
Let $\spc{X}_1,\spc{X}_2,\dots$, and $\spc{X}_\infty$ and $\bar{\spc{X}}_\infty$ be metric spaces
such that $\spc{X}_n\xto{\GH}\spc{X}_\infty$, 
$\spc{X}_n\xto{\bar\GH}\bar{\spc{X}}_\infty$.

Assume that $\bar{\spc{X}}_\infty$ is compact.
Then $\spc{X}_\infty\iso \bar{\spc{X}}_\infty$.
\end{thm}

\parit{Proof.}
For each point $x_\infty\in\spc{X}_\infty$,
choose  liftings $x_n\in\spc{X}_n$.

Choose a nonprincipal ultrafilter $\o$ on $\mathbb N$.
Define $\bar x_\infty\in \bar{\spc{X}}_\infty$ as the $\o$-limit of $x_n$ with respect to $\bar \tau$.
We claim that the map $x_\infty\to \bar x_\infty$ is an isometry.

Indeed, by the definition of Gromov--Hausdorff convergence, 
\[\dist{\bar x_\infty}{\bar y_\infty}{\bar{\spc{X}}_\infty}
=
\lim_{n\to\o}\dist{x_n}{y_n}{\spc{X}_n}
=
\dist{x_\infty}{y_\infty}{\spc{X}_\infty}.
\]
Thus the map $x_\infty\to\bar x_\infty$ gives a distance-preserving map
$\map\:\spc{X}_\infty\hookrightarrow\bar{\spc{X}}_\infty$.
In particular,  
$\spc{X}_\infty$ is compact.
Switching $\spc{X}_\infty$ and $\bar{\spc{X}}_\infty$ and applying the same argument, 
we get an isometric embedding 
$\bar{\spc{X}}_\infty\hookrightarrow\spc{X}_\infty$.
Now the result follows from Lemma~\ref{lem:>=-isometry}.
\qeds

\begin{thm}{Exercise}\label{ex:GH-SC}
\begin{subthm}{ex:GH-SC:circle}
Show that a sequence of compact simply connected length spaces cannot converge to a circle.
\end{subthm}

\begin{subthm}{ex:GH-SC:nonsc-limit}
Construct a sequence of compact simply connected length spaces that converges to a compact non-simply connected space.
\end{subthm}
\end{thm}

\begin{thm}{Exercise}\label{ex:sphere-to-ball}
\begin{subthm}{ex:sphere-to-ball:2}
Show that a sequence of length metrics on the 2-sphere cannot converge to the unit disk.
\end{subthm}

\begin{subthm}{ex:sphere-to-ball:3}
Construct a sequence of length metrics on the 3-sphere that converges to a unit 3-ball.
\end{subthm}

\end{thm}

\begin{thm}{Exercise}\label{ex:GH-proper-marked}
Let $\spc{X}_n$ be a sequence of metric spaces that admits 
two Gromov--Hausdorff convergences
$\GH$ and $\GH'$.
Assume 
$\spc{X}_n\xto{\GH}\spc{X}_\infty$ and $\spc{X}_n\xto{\GH'}\spc{X}_\infty'$.
Show that if $\spc{X}_\infty$ is proper and there is a sequence of points $x_n\in \spc{X}_n$ 
that converges in both
$\GH$ and $\GH'$, then $\spc{X}_\infty\iso\spc{X}_\infty'$.
\end{thm}

\sectionmark{Ultralimits revisited}
\section{Ultralimits revisited}

Recall that $\o$ denotes an ultrafilter of the set of natural numbers.

\begin{thm}{Theorem}\label{thm:ultra-GH}
Assume $\spc{X}_n$ is a sequence of complete metric spaces. 
Let $\spc{X}_n\to \spc{X}_\o$ as $n\to\o$,
and let $\spc{Y}_n$ 
be a sequence of subspaces of the $\spc{X}_n$ such that $\spc{Y}_n\GHto\spc{Y}_\infty$. 
Then there is a distance-preserving map 
$\iota:\spc{Y}_\infty\to \spc{X}_\o$.

Moreover:

\begin{subthm}{thm:ultra-GH:a}
If $\spc{X}_n\GHto \spc{X}_\infty$ 
and $\spc{X}_\infty$ is compact, then 
$\spc{X}_\infty$ is isometric to $\spc{X}_\o$.
\end{subthm}

\begin{subthm}{thm:ultra-GH:b}
If $\spc{X}_n\GHto \spc{X}_\infty$ 
and $\spc{X}_\infty$ is proper, then 
$\spc{X}_\infty$ is isometric to a metric component of $\spc{X}_\o$.
\end{subthm}

\end{thm}

\parit{Proof.} 
For each point $y_\infty\in \spc{Y}_\infty$ 
choose a lifting $y_n\in \spc{Y}_n$.
Pass to the $\o$-limit $y_\o\in \spc{X}_\o$ of $y_n$.
Clearly for any $y_\infty,z_\infty\in \spc{Y}_\infty$, 
we have 
\[\dist{y_\infty}{z_\infty}{\spc{Y}_\infty}=\dist{y_\o}{z_\o}{\spc{X}_\o};\] 
that is, the map $y_\infty\mapsto y_\o$ gives a distance-preserving map $\iota:\spc{Y}_\infty\to \spc{X}_\o$.

\parit{\ref{SHORT.thm:ultra-GH:a}$+$\ref{SHORT.thm:ultra-GH:b}.}
Fix $x_\o\in \spc{X}_\o$.
Choose a sequence $x_n$ of points in $\spc{X}_n$,  
such that $x_n\to x_\o$ as $n\to\o$. 

Denote by $\bm{X}=\spc{X}_\infty\sqcup\spc{X}_1\sqcup\spc{X}_2\sqcup\dots$ the common space for the convergence $\spc{X}_n\GHto \spc{X}_\infty$,
as in the definition of Gromov--Hausdorff convergence.
Note that $x_n$ is a sequence of points in~$\bm{X}$.

If the $\o$-limit $x_\infty$ of $x_n$ in $\bm{X}$ exists, 
it must lie in $\spc{X}_\infty$. 

The point $x_\infty$, if defined, does not depend on the choice of $x_n$.
Indeed, if $y_n\in\spc{X}_n$ is another sequence such that $y_n\to x_\o$ as $n\to\o$, then 
\[
\dist{y_\infty}{x_\infty}{}=\lim_{n\to\o}\dist{y_n}{x_n}{}=0;
\]
therefore, $x_\infty=y_\infty$.

This way we obtain a map $\nu\:x_\o\to x_\infty$, defined on  $\Dom\nu \subset\spc{X}_\o$.
By construction of $\iota$, 
we have $\iota\circ\nu(x_\o)=x_\o$ for any $x_\o\in \Dom\nu$.

Finally note that if $\spc{X}_\infty$ is compact, then $\nu$ is defined on all of $\spc{X}_\o$;
this proves \ref{SHORT.thm:ultra-GH:a}.

If $\spc{X}_\infty$ is proper, choose any point $z_\infty\in \spc{X}_\infty$
and set $z_\o=\iota(z_\infty)$.
For any point $x_\o\in \spc{X}_\o$ at finite distance from $z_\o$,
for the sequence $x_n$ 
as above we have that $\dist{z_n}{x_n}{}$ is bounded for $\o$-almost all $n$.
Since $\spc{X}_\infty$ is proper, $\nu(x_\o)$ is defined;
in other words, $\nu$ is defined on the metric component of~$z_\o$.
Hence \ref{SHORT.thm:ultra-GH:b} follows.
\qeds


\chapter{The ghost of Euclid}

\section{Geodesics, triangles and hinges}
\label{sec:geods}

\parbf{Geodesics and their relatives.}
Let $\spc{X}$ be a metric space and $\II\subset \R$\index{$\II$} be an interval. 
A globally distance-preserving map $\gamma\:\II\to \spc{X}$ is called a \index{unit-speed geodesic}\emph{unit-speed geodesic}.
(Various authors call it differently: \index{shortest path}\emph{shortest path}, \index{minimizing geodesic}\emph{minimizing geodesic}.)
In other words, $\gamma\:\II\to \spc{X}$ is a unit-speed geodesic if the equality
\[\dist{\gamma(s)}{\gamma(t)}{\spc{X}}=|s-t|\]
holds for any pair $s,t\in \II$.

A unit-speed geodesic between $p$ and $q$ in $\spc{X}$ will be denoted by $\geod_{[p q]}$\index{$\geod_{[{p}{q}]}$}.
We will always assume $\geod_{[p q]}$ is parametrized starting at $p$; 
that is, $\geod_{[p q]}(0)=p$ and $\geod_{[p q]}(\dist{p}{q}{})=q$.
The image of $\geod_{[p q]}$ will be denoted by $[p q]$\index{$[{p}{q}]$} and called a \index{geodesic}\emph{geodesic}.
The term \index{geodesic}\emph{geodesic} will also be used for a linear reparametrization of a unit-speed geodesic.
With a slight abuse of notation, we will use the notation $[p q]$ also for the class of all linear reparametrizations of $\geod_{[p q]}$.

A unit-speed geodesic $\gamma\:\RR_{\ge0}\to \spc{X}$ is called a \index{half-line}\emph{half-line}.

A unit-speed geodesic  $\gamma\:\RR\to \spc{X}$ is called a \index{line}\emph{line}.

A piecewise geodesic curve is called \index{polygonal line}\emph{polygonal line};
we may say \textit{polygonal line $p_1,\dots,p_n$} meaning polygonal line with edges $[p_1p_2],\dots,[p_{n-1}p_n]$.
A closed polygonal line will be also called a polygon.

We may write $[p q]_{\spc{X}}$ 
to emphasize that the geodesic $[p q]$ is in the space~${\spc{X}}$.
Also, we use the following short-cut notation:
\begin{align*}
\mathopen{]} p q \mathclose{[}&=[pq]\setminus\{p,q\},
&
\mathopen{]} p q ]&=[pq]\setminus\{p\},
&
[ p q \mathclose{[}&=[pq]\setminus\{q\}.
\end{align*}

In general, a geodesic between $p$ and $q$ need not exist and if it exists, it need not to be unique.
However,  once we write $\geod_{[p q]}$ or $[p q]$ we mean that we have fixed a choice of a geodesic.

A constant-speed geodesic $\gamma\:[0,1]\to\spc{X}$ is called a \index{geodesic!geodesic path}\emph{geodesic path}.
Given a geodesic $[p q]$,
we denote by $\geodpath_{[pq]}$ the corresponding geodesic path;
that is, 
$$\geodpath_{[pq]}(t)\z\equiv\geod_{[pq]}(t\cdot\dist[{{}}]{p}{q}{}).$$

A curve $\gamma\:\II\to \spc{X}$ is called a \index{geodesic!local geodesic}\emph{local geodesic} if for any $t\in\II$ there is a neighborhood $U\ni t$ in $\II$ such that the restriction $\gamma|_U$ is a constant-speed geodesic.
If $\II=[0,1]$, then $\gamma$ is called a \emph{local geodesic path}.

\begin{thm}{Proposition}\label{prop:busemann}
Suppose $\spc{X}$ is a metric space and $\gamma\:[0,\infty)\to \spc{X}$ is a half-line. 
Then the \index{Busemann function}\emph{Busemann function} $\bus_\gamma\:\spc{X}\to \RR$ 
\[\bus_\gamma(x)=\lim_{t\to\infty}\dist{\gamma(t)}{x}{}- t\eqlbl{eq:def:busemann*}\]
is defined
and $1$-Lipschitz.
\end{thm}

\parit{Proof.}
By the triangle inequality, the function $t\mapsto\dist{\gamma(t)}{x}{}- t$ is nonincreasing.  
Clearly $\dist{\gamma(t)}{x}{}- t\ge-\dist{\gamma(0)}{x}{}$.
Thus the limit in \ref{eq:def:busemann*} is defined,
and it is 1-Lipschitz as a limit of 1-Lipschitz functions.
\qeds

\begin{thm}{Example} 
If $\spc{X}$ is a Euclidean space and $\gamma(t)=p+t\cdot v$ where $v$ is a unit vector,
then
\[\bus_\gamma(x)=\langle x-p, v\rangle.\]
\end{thm}

\parbf{Triangles.}
For a triple of points $p,q,r\in \spc{X}$, a choice of a triple of geodesics $([q r], [r p], [p q])$ will be called a \index{triangle}\emph{triangle}, and we will use the short notation 
$\trig p q r=([q r], [r p], [p q])$\index{$\trig{p}{q}{r}$}.
Again, given a triple $p,q,r\in \spc{X}$, there may be no triangle 
$\trig p q r$, simply because one of the pairs of these points cannot be joined by a geodesic.  Or there may be many different triangles, any of which can be denoted by $\trig p q r$.
Once we write $\trig p q r$, it means we have chosen such a triangle; 
that is, made a choice of each $[q r]$, $[r p]$, and $[p q]$.

The value 
$\dist{p}{q}{}+\dist{q}{r}{}+\dist{r}{p}{}$ 
will be called the \index{perimeter!perimeter of a triangle}\emph{perimeter of triangle} $\trig p q r$;
it obviously coincides with perimeter of the triple $p$, $q$, $r$ as defined below.

\parbf{Hinges.}
Let $p,x,y\in \spc{X}$ be a triple of points such that $p$ is distinct from $x$ and $y$.
A pair of geodesics $([p x],[p y])$ will be called a  \index{hinge}\emph{hinge}, and will be denoted by 
$\hinge p x y=([p x],[p y])$\index{$\hinge{{p}}{{q}}{{r}}$}.


\section{Model angles and triangles}\label{sec:mod-tri/angles}

Let $\spc{X}$ be a metric space, 
$p,q,r\in \spc{X}$, 
and $\kappa\in\RR$. 
Let us define the \index{model triangle}\emph{model triangle} $\trig{\tilde p}{\tilde q}{\tilde r}$ 
(briefly, 
$\trig{\tilde p}{\tilde q}{\tilde r}=\modtrig\kappa(p q r)$%
\index{$\modtrig\kappa$!$\modtrig\kappa({p}{q}{r})$}) to be a triangle in the model plane $\Lob2\kappa$ such that
\[\dist{\tilde p}{\tilde q}{}=\dist{p}{q}{},
\quad \dist{\tilde q}{\tilde r}{}=\dist{q}{r}{},
\quad \dist{\tilde r}{\tilde p}{}=\dist{r}{p}{}.\]

In the notation of Section~\ref{model}, 
$\modtrig\kappa(p q r)=\modtrig\kappa\{\dist{q}{r}{},\dist{r}{p}{},\dist{p}{q}{}\}$.

If $\kappa\le 0$, the  model triangle is  always defined, that is, it exists and is unique up to an isometry of $\Lob2\kappa$.
If $\kappa>0$, the model triangle is said to be defined if in addition
\[\dist{p}{q}{}+\dist{q}{r}{}+\dist{r}{p}{}< 2\cdot\varpi\kappa;\]
here $\varpi\kappa$ denotes the diameter of the model space $\Lob2\kappa$.
In this case, the model triangle also exists and is unique up to an isometry of $\Lob2\kappa$.
The value $\dist{p}{q}{}+\dist{q}{r}{}+\dist{r}{p}{}$ will be called the \index{perimeter!perimeter of a triple}\emph{perimeter of the triple} $p$, $q$, $r$.

If for  $p,q,r\in \spc{X}$,
$\trig{\tilde p}{\tilde q}{\tilde r}=\modtrig\kappa(p q r)$ is defined 
and $\dist{p}{q}{},\dist{p}{r}{}>0$, the angle measure of 
$\trig{\tilde p}{\tilde q}{\tilde r}$ at $\tilde  p$ will be called the \index{model angle}\emph{model angle} of the triple $p$, $q$, $r$, and will be denoted by
$\angk\kappa p q r$%
\index{$\tangle\mc\kappa$!$\angk\kappa{{p}}{{q}}{{r}}$}\index{$\tangle\mc\kappa$!$\angk\kappa{{p}}{{q}}{{r}}$}.

In the notation of Section~\ref{model}, 
\[\angk\kappa p q r=\tangle\mc\kappa\{\dist{q}{r}{};\dist{p}{q}{},\dist{p}{r}{}\}.\]

\begin{wrapfigure}{r}{25mm}
\vskip-0mm
\centering
\includegraphics{mppics/pic-605}
\end{wrapfigure}

\begin{thm}{Alexandrov's lemma}
\index{Alexandrov's lemma}
\label{lem:alex}  
Let $p,q,r,z$ be distinct points in a metric space such that $z\in \mathopen{]}p r\mathclose{[}$ and 
\[\dist{p}{q}{}+\dist{q}{r}{}+\dist{r}{p}{}< 2\cdot\varpi\kappa.\]
Then 
the following expressions have the same sign:
\begin{subthm}{lem-alex-difference}
$
\angk\kappa p q r-\angk\kappa p q z$,
\end{subthm} 

\begin{subthm}{lem-alex-angle}
$\angk\kappa z q p
+\angk\kappa z q r -\pi$.
\end{subthm}

Moreover,
\[\angk\kappa q p r \ge \angk\kappa q p z +  \angk\kappa q z r,\]
with equality if and only if the expressions in \ref{SHORT.lem-alex-difference} and \ref{SHORT.lem-alex-angle} vanish.
\end{thm}

\parit{Proof.}
By the triangle inequality, 
\[
\dist{p}{q}{}+\dist{q}{z}{}+\dist{z}{p}{}\le \dist{p}{q}{}+\dist{q}{r}{}+\dist{r}{p}{}< 2\cdot\varpi\kappa.
\]
Therefore the model triangle $\trig{\tilde p}{\tilde q}{\tilde z}=\modtrig\kappa p q z$ is defined.
Take 
a point $\tilde r$ on the extension of 
$[\tilde p \tilde z]$ beyond $\tilde z$ so that $\dist{\tilde p}{\tilde r}{}=\dist{p}{r}{}$ (and therefore $\dist{\tilde p}{\tilde z}{}=\dist{p}{z}{}$). 
 
From monotonicity of the function $a\mapsto\tangle\mc\kappa\{a;b,c\}$ (\ref{increase}), 
the following expressions have the same sign:
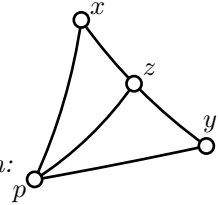
\begin{enumerate}[(i)]
\item $\mangle\hinge{\tilde p}{\tilde q}{\tilde r}-\angk\kappa{p}{q}{r}$;
\item $\dist{\tilde p}{\tilde r}{}-\dist{p}{r}{}$;
\item $\mangle\hinge{\tilde z}{\tilde q}{\tilde r}-\angk\kappa{z}{q}{r}$.
\end{enumerate}

\begin{wrapfigure}{r}{25mm}
\vskip-15mm
\centering
\includegraphics{mppics/pic-610}
\bigskip
\includegraphics{mppics/pic-611}
\vskip-0mm
\end{wrapfigure}

Since 
\[\mangle\hinge{\tilde p}{\tilde q}{\tilde r}=\mangle\hinge{\tilde p}{\tilde q}{\tilde z}=\angk\kappa{p}{q}{z}\]
and
\[ \mangle\hinge{\tilde z}{\tilde q}{\tilde r}
=\pi-\mangle\hinge{\tilde z}{\tilde p}{\tilde q}
=\pi-\angk\kappa{z}{p}{q},\]
the first statement follows.

For the second statement, let us redifine $\tilde r$;
construct $\trig{\tilde q}{\tilde z}{\tilde r}=\modtrig\kappa q z r$ on the opposite side of $[\tilde q\tilde z]$ from $\trig{\tilde p}{\tilde q}{\tilde z}$.  
Since
\begin{align*}
\dist{\tilde p}{\tilde r}{}&\le \dist{\tilde p}{\tilde z}{} + \dist{\tilde z}{\tilde r}{}=
\\
&=\dist{p}{z}{}+\dist{z}{r}{}=
\\
&=\dist{p}{r}{},
\end{align*}
we have
\begin{align*}
\angk\kappa{q}{p}{z} + \angk\kappa{q}{z}{r} 
&
= 
\mangle\hinge{\tilde q}{\tilde p}{\tilde z}+ \mangle\hinge{\tilde q}{\tilde z}{\tilde r} 
=
\\
&
= 
\mangle\hinge{\tilde q}{\tilde p}{\tilde r}
\le
\\
&\le  \angk\kappa q p r.
\end{align*}
Equality holds if and only  if $\dist{\tilde p}{\tilde r}{}=\dist{p}{r}{}$, 
as required.\qeds


\section{Angles and the first variation}\label{sec:angles}

Given a hinge $\hinge p x y$, we define its \index{angle}\emph{angle} to be \index{$\mangle$!$\mangle\hinge{{p}}{{q}}{{r}}$}
\[\mangle\hinge p x y
\df
\lim_{\bar x,\bar y\to p} \angk\kappa p{\bar x}{\bar y},\eqlbl{eq:angle-def}\]
for $\bar x\in\mathopen{]}p x]$ and $\bar y\in\mathopen{]}p y]$, if this limit exists.

Similarly to $\angk\kappa p{x}{y}$, 
we will use the short notation\index{$\side\kappa$!$\side\kappa \hinge{{p}}{{q}}{{r}}$}
\[\side\kappa \hinge p x y=
\side\kappa \left\{\mangle\hinge p x y;\dist{p}{x}{},\dist{p}{y}{}\right\},\]
where the right-hand side is defined in Section~\ref{model}.  
The value $\side\kappa \hinge p x y$ will be called the \index{model side}\emph{model side}
 of the hinge $\hinge p x y$.

\begin{thm}{Lemma}\label{lem:k-K-angle}
Let $p,x,y$ be a triple of points in a metric space with perimeter $\ell$.
Then for any $\kappa,\Kappa\in\RR$,
\[|\angk\Kappa p{x}{y}-\angk\kappa p{x}{y}|
\le 
100(|\Kappa|+|\kappa|)\cdot\ell^2,
\eqlbl{eq:k-K},\]
 whenever the left-hand side is defined.
\end{thm}

Lemma~\ref{lem:k-K-angle} implies that 
the definition of angle is independent of $\kappa$.
In particular, one can take $\kappa=0$ in \ref{eq:angle-def};
thus the angle can be calculated from the  cosine law:
\[\cos\angk{0}{p}{x}{y}
=
\frac{\dist[2]{p}{x}{}+\dist[2]{p}{y}{}-\dist[2]{x}{y}{}}{2\cdot \dist[{{}}]{p}{x}{}\cdot\dist[{{}}]{p}{y}{}}.\]

\parit{Proof.}
The function $\kappa\mapsto \angk\kappa p{x}{y}$ is nondecreasing (\ref{k-decrease}).
Thus, for $\Kappa>\kappa$, we have
\begin{align*}
0\le \angk\Kappa p{x}{y}-\angk{\kappa}p{x}{y}
&\le \angk\Kappa p{x}{y}+\angk\Kappa {x}p{y}+\angk\Kappa {y}p{x}-
\\
&\quad-\angk\kappa p{x}{y}-\angk\kappa {x}p{y}-\angk\kappa {y}p{x}
= 
\\
&=\Kappa\cdot\area\modtrig\Kappa(pxy)-\kappa\cdot\area\modtrig\kappa(pxy).
\end{align*}
Note that for $ \kappa\ge 0$ a triangle of perimeter $\ell$ in  $\Lob2\kappa$ lies in a ball of radius $2\cdot \ell$,  which easily implies that  $\area\modtrig\kappa(pxy)\le 100\cdot\ell^2$.
For $\kappa<0$ one gets the same estimate by a direct computation in the hyperbolic plane.

Therefore
\begin{align*}
\area\modtrig\kappa(pxy)&\le 100\cdot\ell^2, 
&
\area\modtrig\Kappa(pxy)&\le 100\cdot\ell^2.
\end{align*}
Thus \ref{eq:k-K} follows.
\qeds

\begin{thm}{Triangle inequality for angles}
\label{claim:angle-3angle-inq}
Let  $[px^1]$, $[px^2]$, and $[px^3]$ be three geodesics in a metric space.
If all of the angles $\alpha^{i j}=\mangle\hinge p {x^i}{x^j}$ are defined then they satisfy the triangle inequality:
\[\alpha^{13}\le \alpha^{12}+\alpha^{23}.\]

\end{thm}

\parit{Proof.}
Since $\alpha^{13}\le\pi$, we can assume that $\alpha^{12}+\alpha^{23}< \pi$.
Set $\gamma^i\z=\geod_{[px^i]}$.
Given any $\eps>0$, for all sufficiently small $t,\tau,s\in\RR_{\ge0}$ we have
\begin{align*}
\dist{\gamma^1(t)}{\gamma^3(\tau)}{}
&\le 
\dist{\gamma^1(t)}{\gamma^2(s)}{}+\dist{\gamma^2(s)}{\gamma^3(\tau)}{}<\\
&<
\sqrt{t^2+s^2-2\cdot t\cdot  s\cdot \cos(\alpha^{12}+\eps)} +
\\
&\quad+\sqrt{s^2+\tau^2-2\cdot s\cdot \tau\cdot \cos(\alpha^{23}+\eps)}\le
\end{align*}

\begin{wrapfigure}{o}{30 mm}
\vskip-0mm
\centering
\includegraphics{mppics/pic-615}
\vskip3mm
\end{wrapfigure}

Below we define 
$s(t,\tau)$ so that for 
$s=s(t,\tau)$, this chain of inequalities can be continued as follows:
\[\le
\sqrt{t^2+\tau^2-2\cdot t\cdot \tau\cdot \cos(\alpha^{12}+\alpha^{23}+2\cdot \eps)}.
\]

Thus for any $\eps>0$, 
\[\alpha^{13}\le \alpha^{12}+\alpha^{23}+2\cdot \eps.\]
Hence the result follows.

To define $s(t,\tau)$, consider three half-lines $\tilde \gamma^1$, $\tilde \gamma^2$, $\tilde \gamma^3$ on a Euclidean plane starting at one point, such that
$\mangle(\tilde \gamma^1,\tilde \gamma^2)\z=\alpha^{12}+\eps$,
$\mangle(\tilde \gamma^2,\tilde \gamma^3)\z=\alpha^{23}+\eps$,
and $\mangle(\tilde \gamma^1,\tilde \gamma^3)\z=\alpha^{12}\z+\alpha^{23}\z+2\cdot \eps$.
We parametrize each half-line by the distance from the starting point.
Given two positive numbers $t,\tau\in\RR_{\ge0}$, let $s=s(t,\tau)$ be 
the number such that 
$\tilde \gamma^2(s)\in[\tilde \gamma^1(t)\ \tilde \gamma^3(\tau)]$. 
Clearly $s\le\max\{t,\tau\}$, so $t,\tau,s$ may be taken sufficiently small.
\qeds 

\begin{thm}{Exercise}\label{ex:adjacent-angles}
Prove that the sum of adjacent angles is at least $\pi$.

More precisely,
suppose that the hinges $\hinge pxz$ and $\hinge pyz$ are \index{adjacent hinges}\emph{adjacent};
that is, the side $[pz]$ is shared and the union of two sides $[px]$ and $[py]$ is a geodesic $[xy]$.
Then 
\[\mangle\hinge pxz+\mangle\hinge pyz\ge \pi\]
whenever each angle on the left-hand side is defined.
\end{thm}

The above inequality can be strict.
For example in a metric tree angles between any two different edges coming out of the same vertex are all equal to $\pi$.

\begin{thm}{First variation inequality}\label{lem:first-var}
Assume that for a hinge $\hinge q p x$, 
the angle $\alpha=\mangle\hinge q p x$ is defined. Then
\[\dist{p}{\geod_{[qx]}(t)}{}
\le
\dist{q}{p}{}-t\cdot \cos\alpha+o(t).\]

\end{thm}

\parit{Proof.}
Take a sufficiently small $\eps>0$.
For all sufficiently small $t>0$, we have 
\begin{align*}
 \dist{\geod_{[qp]}(t/\eps)}{\geod_{[qx]}(t)}{}
&\le 
\tfrac{t}{\eps}\cdot \sqrt{1+\eps^2 -2\cdot \eps\cdot \cos\alpha}+o(t)\le
\\
&\le \tfrac{t}{\eps} -t\cdot \cos\alpha + t\cdot \eps.
\end{align*}
Applying the triangle inequality, we get 
\begin{align*}
\dist{p}{\geod_{[qx]}(t)}{}
&\le \dist{p}{\geod_{[qp]}(t/\eps)}{}+\dist{\geod_{[qp]}(t/\eps)}{\geod_{[qx]}(t)}{}
\le 
\\
&\le
\dist{p}{q}{} -t\cdot \cos\alpha + t\cdot \eps
\end{align*}
for any $\eps>0$ and all sufficiently small $t$.
Hence the result.
\qeds

\section{Space of directions} 
\label{sec:tangent-space+directions}

Let $\spc{X}$ be a metric space.
If the angle $\mangle\hinge pxy$ is defined for any hinge $\hinge pxy$ in $\spc{X}$,
then we will say that the space $\spc{X}$ has \index{angle!defined angles}\emph{defined angles}.

Let us note that this is a strong condition. For example a Banach space which has defined angles must be Hilbert.

 
Let $\spc{X}$ be a space with defined angles. For $p\in \spc{X}$,
consider the set $\mathfrak{S}_p$ 
of all nontrivial unit-speed geodesics starting at $p$.
By \ref{claim:angle-3angle-inq}, the triangle inequality holds for $\mangle$ on $\mathfrak{S}_p$,
that is, $(\mathfrak{S}_p,\mangle)$ 
forms a pseudometric space.

The metric space corresponding to  $(\mathfrak{S}_p,\mangle)$ is called the \index{direction!space of geodesic directions}\emph{space of geodesic directions} at $p$, denoted by $\Sigma'_p$ or $\Sigma'_p\spc{X}$.
The elements of $\Sigma'_p$ are called \emph{geodesic directions} at $p$.
Each geodesic direction is formed by an equivalence class of geodesics starting from $p$ 
for the equivalence relation 
\[[px]\sim[py]\quad \iff\quad \mangle\hinge pxy=0;\]
the direction of $[px]$ is denoted by $\dir px $.\index{$\dir{p}{q}$}

The completion of $\Sigma'_p$ is called the \index{direction!space of directions}\emph{space of directions} at $p$ and is denoted by $\Sigma_p$ or $\Sigma_p\spc{X}$.
The elements of $\Sigma_p$ are called \index{direction}\emph{directions} at $p$.

\section{Tangent space}\label{sec: tangent space}

The \index{cone}\emph{Euclidean cone} $\spc{Y}=\Cone\spc{X}$ 
over a metric space $\spc{X}$
is defined as the metric space whose underlying set consists of
equivalence classes in
$[0,\infty)\times \spc{X}$ with the equivalence relation ``$\sim$'' given by $(0,p)\sim (0,q)$ for any points $p,q\in\spc{X}$,
and whose metric is given by the cosine rule
\[
\dist{(s,p)}{(t,q)}{\spc{Y}} 
=
\sqrt{s^2+t^2-2\cdot s\cdot t\cdot \cos\theta},
\]
where $\theta= \min\{\pi, \dist{p}{q}{\spc{X}}\}$.
The point in $\spc{Y}$ that corresponds $(t,x)\in[0,\infty)\times \spc{X}$ will be denoted by $t\cdot x$.

The point in  $\Cone{\spc{X}}$ formed by the equivalence class of $\{\0\}\times\spc{X}$ is called the \index{tip of a  cone}\emph{tip of the cone} and is denoted by $\0$ or $\0_{\spc{Y}}$.
For $v\in\spc{Y}$ the distance $\dist{\0}{v}{\spc{Y}}$ is called the norm of $v$ and is denoted by $|v|$ or $|v|_{\spc{Y}}$.

The \index{scalar product}\emph{scalar product} $\<v,w\>$
of two vectors $v=s\cdot p$ and $w=t\cdot q$
is defined by 
\[\<v,w\>
\df |v|\cdot|w|\cdot\cos\theta;
\]
we set $\<v,w\>\df0$ if $v=\0$ or $w=\0$.

\begin{thm}{Example}
$\Cone\SS^n$ is isometric to $\R^{n+1}$.
If $G< \O(n+1)$ is a closed subgroup, then $\Cone(\SS^n/G)$ is isometric to $\R^{n+1}/G$.
\end{thm}

The Euclidean cone $\Cone\Sigma_p$ over the space of directions $\Sigma_p$ is called the \index{tangent space}\emph{tangent space} at $p$ and denoted by $\T_p$ or $\T_p\spc{X}$.
The elements of $\T_p\spc{X}$ will be called \index{tangent vector}\emph{tangent vectors} at $p$
(despite the fact that $\T_p$ is only a cone --- not a vector space).

The tangent space $\T_p$ could be also defined directly, without introducing the space of directions.
To do so, consider the set $\mathfrak{T}_p$ of all geodesics starting at $p$, with arbitrary speed.
Given $\alpha,\beta\in \mathfrak{T}_p$,
set 
\[\dist{\alpha}{\beta}{\mathfrak{T}_p}
=
\lim_{\eps\to0} 
\frac{\dist{\alpha(\eps)}{\beta(\eps)}{\spc{X}}}\eps.
\eqlbl{eq:dist-in-T_p}.\]
If the angles in $\spc{X}$ are defined, then so is
the limit in \ref{eq:dist-in-T_p}, and we obtain a pseudometric on $\mathfrak{T}_p$.

The corresponding metric space admits a natural isometric identification with the cone $\T'_p=\Cone\Sigma'_p$.
The vectors of $\T'_p$ are the equivalence classes for the relation 
\[\alpha\sim\beta\quad \iff\quad \dist{\alpha(t)}{\beta(t)}{\spc{X}}=o(t).\]
The completion of $\T'_p$ is therefore naturally isometric to $\T_p$.
A vector in $\T'_p$ that corresponds to the geodesic path $\geod_{[pq]}$ is called \index{logarithm}\emph{logarithm of $[pq]$} and denoted by $\ddir pq$.\index{$\ddir {p}{q}$}

\section{Velocity of curves}

\begin{thm}{Definition}\label{def:right-derivative}
Let $\spc{X}$ be a metric space,
$a>0$,
and $\alpha\:[0,a)\to \spc{X}$ be a function, not necessarily continuous, such that $\alpha(0)=p$.
We say that $v\in\T_p$ is the \index{right derivative}\emph{right derivative} of $\alpha$ at $0$,
briefly $\alpha^+(0)\z=v$\index{$\alpha^+$}, if for some (and therefore any) sequence of vectors $v_n\in\T'_p$ such that $v_n\to v$ as $n\to\infty$,
and corresponding geodesics $\gamma_n$, 
we have 
\[\limsup_{\eps\to0+}\frac{\dist{\alpha(\eps)}{\gamma_n(\eps)}{\spc{X}}}{\eps}\to 0\quad \text{as}\quad n\to\infty.\]

We define right and left derivatives $\alpha^+(t_0)$ and $\alpha^-(t_0)$
of $\alpha$ at $t_0\in\II$ by 
\[\alpha^\pm(t_0)=\check\alpha^+(0),\] where $\check\alpha(t)=\alpha(t_0\pm t)$.
\end{thm}

The sign convention is not quite standard; if $\alpha$ is a smooth curve in a Riemannian manifold then we have
$\alpha^+(t)=-\alpha^-(t)$.

Note that if $\gamma$ is a geodesic starting at $p$ 
and the tangent vector $v\in\T_p'$ corresponds to $\gamma$, 
then $\gamma^+(0)=v$.

\begin{thm}{Exercise}\label{ex:tangent-vect=o(t)}
Assume $\spc{X}$ is a metric space with defined angles,
and let $\alpha,\beta\:[0,a)\to\spc{X}$ 
be two maps such that the right derivatives $\alpha^+(0)$, $\beta^+(0)$ are defined and $\alpha^+(0)=\beta^+(0)$.
Show that
\[\dist{\alpha(t)}{\beta(t)}{\spc{X}}=o(t).\]
\end{thm}

\begin{thm}{Proposition}
Let $\spc{X}$ be a metric space with defined angles and $p\in \spc{X}$.
Then for any tangent vector $v\in\T_p\spc{X}$ there is a map $\alpha\:[0,\eps)\to \spc{X}$ such that $\alpha^+(0)=v$.
\end{thm}

\parit{Proof.}
If $v\in \T_p'$, then for the corresponding geodesic $\alpha$ we have $\alpha^+(0)\z=v$.

Given $v\in \T_p$, construct a sequence $v_n\in\T'_p$ 
such that $v_n\to v$, and let $\gamma_n$ be a sequence of corresponding geodesics.

The needed map $\alpha$ can be found among the maps such that $\alpha(0)\z=p$ and
\[\alpha(t)=\gamma_n(t)\quad \text{if}\quad \eps_{n+1}\le t<\eps_n,\]
where $\eps_n$
is a decreasing sequence converging to $0$ as $n\to\infty$.
In order to satisfy the conclusion of the proposition, one has to choose the sequence $\eps_n$ converging to $0$ very fast.
Note that in this construction $\alpha$ is not continuous.
\qeds

\begin{thm}{Definition}\label{def:diff-curv}
Let 
$\spc{X}$ be a metric space 
and $\alpha\:\II\to \spc{X}$ be a curve.

For $t_0\in\II$, 
if $\alpha^+(t_0)$ or $\alpha^-(t_0)$ or both are defined,
we say respectively that $\alpha$ is 
\index{curve!differentiable curve}
\index{curve!left/right differentiable}
\emph{right} or \emph{left} or \emph{both-sided differentiable} at $t_0$.
In the exceptional cases where $t_0$ is the left (respectively right) end of $\II$, $\alpha$ is by definition left (respectively right) differentiable at $t_0$.

If $\alpha$ is both-sided differentiable at $t$, and 
\[|\alpha^+(t)|=|\alpha^-(t)|=\tfrac12\cdot\dist{\alpha^+(t)}{\alpha^-(t)}{\T_{\alpha(t)}},\] then we say that $\alpha$ is \emph{differentiable} at $t$.
\end{thm}

\begin{thm}{Exercise}\label{ex:both-sided-diff}
Assume $\spc{X}$ is a metric space with defined angles.
Show that any geodesic $\gamma\:\II\to\spc{X}$ is differentiable everywhere.
\end{thm}

Recall that the speed of a curve is defined in \ref{thm:speed}.

\begin{thm}{Exercise}\label{ex:diff}
Let $\alpha$ be a curve in a metric space with defined angles.
Suppose that $\speed_t\alpha$, $\alpha^+(t)$, and $\alpha^-(t)$ are defined.

Show that $\alpha$ is differentiable at $t$.
\end{thm}

\section{Differential}\index{differential of a function}

\begin{thm}{Definition}\label{def:differential}
Let $\spc{X}$ be a metric space with defined angles, and
$f\:\spc{X}\subto\RR$ be a subfunction.
For $p\in\Dom f$, a function $\phi\:\T_p\to\RR$ is called the \index{differential}\emph{differential} of $f$ at $p$
(briefly $\phi=\dd_pf$) if for any map $\alpha\:\II\to \spc{X}$ such that $\II$ is a real interval, $\alpha(0)=p$,  and $\alpha^+(0)$ is defined, we have \[(f\circ\alpha)^+(0)=\phi(\alpha^+(0)).\]
\end{thm}

\begin{thm}{Proposition}\label{prop:differential}
Let $f\:\spc{X}\subto\RR$ be a locally Lipschitz semiconcave subfunction
on a metric space $\spc{X}$ with defined angles.
Then the differential $\dd_pf$ is uniquely defined for any $p\in\Dom f$. Moreover, 

\begin{subthm}{prop:differential:lip}
The differential $\dd_pf\:\T_p\to\RR$ is Lipschitz and 
\[\lip\dd_pf\le \lip_pf;\]
that is, the Lipschitz constant of $\dd_pf$ does not exceed the Lipschitz constant of $f$ in any neighborhood of $p$. 
\end{subthm}

\begin{subthm}{prop:differential:homo}
$\dd_pf\:\T_p\to\RR$ is a positive homogeneous function;
that is, for any $r\ge 0$ and $v\in\T_p$ we have 
\[r\cdot\dd_pf(v)=\dd_pf(r\cdot v).\]
\end{subthm}

\begin{subthm}{prop:differential:ultra}
The differential $\dd_pf\:\T_p\to\RR$ is the restriction of ultradifferential defined in Section~\ref{sec:ultradifferential};
that is,
\[\dd_pf=\dd^\o_pf|_{\T_p}.\]
\end{subthm}

\end{thm}

\parit{Proof.}
Passing to a subdomain of $f$ if necessary,
we can assume that $f$ is $\Lip$-Lipschitz and $\lambda$-concave for some $\Lip,\lambda\in\RR$.

Take a geodesic $\gamma$  in $\Dom f$ starting at $p$.
Since $f\circ\gamma$ is $\lambda$-concave,
the right derivative $(f\circ\gamma)^+(0)$ is defined.
Since $f$ is  $\Lip$-Lipschitz, we have
\[|(f\circ\gamma)^+(0)-(f\circ\gamma_1)^+(0)|
\le
\Lip\cdot\dist[{{}}]{\gamma^+(0)}{\gamma_1^+(0)}{}\eqlbl{gam-bargam}\]
for any other geodesic $\gamma_1$ starting at $p$.

Define $\phi\:\T'_p\to\RR\:\gamma^+(0)\mapsto(f\circ\gamma)^+(0)$.
From \ref{gam-bargam}, $\phi$ is an $\Lip$-Lipschtz function defined on $\T_p'$.
Thus we can extend $\phi$ to all of  $\T_p$ as an $\Lip$-Lipschitz function. 

{\sloppy 

It remains to check that $\phi$ is the differential of $f$ at $p$.
Assume $\alpha\:[0,a)\to\spc{X}$ is a map such that $\alpha(0)=p$ and $\alpha^+(0)=v\in \T_p$.
Let $\gamma_n\in\Gamma_p$ be a sequence of geodesics as in the definition \ref{def:right-derivative};
that is, if 
\[v_n=\gamma^+_n(0)\quad \text{and}\quad a_n= \limsup_{t\to0+}{\dist{\alpha(t)}{\gamma_n(t)}{}}/{t}\] 
then $a_n\to 0$ and $v_n\to v$ as $n\to\infty$.
Then 
\[\phi(v)=\lim_{n\to\infty}\phi(v_n),\] \[f\circ\gamma_n(t)=f(p)+\phi(v_n)\cdot t+o(t),\] 
\[|f\circ\alpha(t)-f\circ\gamma_n(t)|
\le
\Lip\cdot\dist[{{}}]{\alpha(t)}{\gamma_n(t)}{}.\]
Hence 
\[f\circ\alpha(t)=f(p)+\phi(v)\cdot t+o(t).\]

The last part follows from the definitions of differential and ultradifferential; see Section~\ref{sec:Ultralimits of functions}.
\qeds

}

\section{Ultratangent space}
\label{sec: ultradiff}

Fix a selective ultrafilter $\o$ on the set of natural numbers.

For a metric space $\spc{X}$ and $r>0$,
we will denote by $r\cdot\spc{X}$ its \index{blowup}\emph{$r$-blowup},
which is a metric space with the same underlying set as $\spc{X}$ and the metric multiplied by $r$.
The tautological bijection $\spc{X}\to r\cdot\spc{X}$ will be denoted by $x\mapsto x^r$, 
so 
\[\dist{x^r}{y^r}{}
=
r\cdot\dist[{{}}]{x}{y}{}\] 
for any $x,y\in \spc{X}$.

The \emph{$\o$-blowup} $\o\cdot\spc{X}$ of $\spc{X}$ is defined to be the $\o$-limit
of the $n$-blowups $n\cdot\spc{X}$; that is,
\[\o\cdot\spc{X}
\df
\lim_{n\to\o} n\cdot\spc{X}.\]

Given a point $x\in \spc{X}$, we can consider the sequence $x^n$, where $x^n\in n\cdot\spc{X}$ is the image of $x$ under $n$-blowup.
Note that if $x\ne y$, then 
\[\dist{x^\o}{y^\o}{\o\cdot\spc{X}}=\infty;\]
that is, 
$x^\o$ and $y^\o$ 
belong to different metric components of $\o\cdot\spc{X}$.

The metric component of $x^\o$ in $\o\cdot\spc{X}$ is called the \index{ultratangent space}\emph{ultratangent space} of $\spc{X}$ at $x$ and  is denoted by $\T^\o_x\spc{X}$ or $\T^\o_x$.

Equivalently, the ultratangent space $\T^\o_x\spc{X}$ can be defined as follows.
Consider all the sequences of points $x_n\in \spc{X}$ such that the sequence 
 $(n\cdot\dist{x}{x_n}{\spc{X}})$ is bounded.
Define the pseudodistance between two such sequences as 
\[\dist{(x_n)}{(y_n)}{}
=
\lim_{n\to\o}n\cdot\dist{x_n}{y_n}{\spc{X}}.\]
Then $\T^\o_x\spc{X}$ is the corresponding metric space.

Tangent spaces (see section \ref{sec: tangent space}) as well as ultratangent spaces 
generalize the notion of tangent spaces on Riemannian manifolds.
In  the simplest cases these two notions define the same space.
However in general they are different and are both useful ---
often a lack of a property in one is compensated by the other. 

It is clear from the definition that a tangent space has a cone structure.
On the other hand, in general an ultratangent space does not have a cone structure.
Hilbert's cube $\prod_{n=1}^\infty[0,2^{-n}]$ is an example.
We remark that Hilbert's cube is a $\Alex{0}$ as well as a $\CAT{0}$ Alexandrov space.

The next theorem shows that the tangent space $\T_p$ can be (and often will be) considered as a subset of  $\T^\o_p$.

\begin{thm}{Theorem}\label{thm:tangent-ultratangent}
\label{thm:T-in-T^w} 
Let $\spc{X}$ be a metric space with defined angles.
Then for any $p\in \spc{L}$, there is a distance-preserving map 
\[\iota:\T_p\hookrightarrow \T^\o_p\] 
such that for any geodesic $\gamma$ starting at $p$
we have 
\[\gamma^+(0)\stackrel{\iota}{\mapsto} \lim_{n\to\o}[\gamma(\tfrac1n)]^n.\]

\end{thm}

\parit{Proof.}
Given $v\in \T'_p$, 
choose a geodesic $\gamma$ that starts at $p$ and  such that $\gamma^+(0)\z=v$.
Set $v^n=[\gamma(\tfrac1n)]^n\in n\cdot \spc{X}$ and 
\[v^\o=\lim_{n\to\o}v^n.\]

Note that the value $v^\o\in\T^\o_p$ does not depend on the choice of $\gamma$;
that is, if $\gamma_1$ is another geodesic starting at $p$ such that $\gamma_1^+(0)=v$,
then 
\[\lim_{n\to\o}v^n=\lim_{n\to\o}v_1^n,\]
where $v_1^n=[\gamma_1(\tfrac1n)]^n\in n\cdot \spc{X}$.
The latter follows since
\[\dist{\gamma(t)}{\gamma_1(t)}{\spc{X}}=o(t),\]
and therefore $\dist{v^n}{v_1^n}{n\cdot \spc{X}}\to 0$ as $n\to\infty$.

Set $\iota(v)=v^\o$.
Since angles between geodesics in $\spc{X}$ are defined, for any $v,w\in \T_p'$ we have
$n\cdot\dist[{{}}]{v_n}{w_n}{}\to\dist{v}{w}{}$.
Thus $\dist{v_\o}{w_\o}{}=\dist{v}{w}{}$; that is, $\iota\:\T_p'\to\T_p$ is a distance-preserving map.

Since $\T_p'$ is dense in $\T_p$,
we can extend $\iota$ to a distance-preserving map $\T_p\to \T^\o_p$.
\qeds

\section{Ultradifferential}\label{sec:ultradifferential}

Given a function $f\:\spc{L}\to\RR$, consider the sequence of functions $f_n\:n\cdot\spc{L}\to\RR$ defined by 
\[f_n(x^n)=n\cdot(f(x)-f(p)),\]
where $x\mapsto x^n$ denotes the natural map $\spc{L}\to n\cdot\spc{L}$.
While $n\cdot(\spc{L},p)\to(\T^\o,\0)$ as $n\to\o$, 
the functions $f_n$ converge to the \emph{$\o$-differential} of $f$ at $p$.
It will be denoted by $\dd_p^\o f$:
\[\dd_p^\o f\:\T_p^\o\to\RR,\quad \dd_p^\o f=\lim_{n\to\o} f_n.\] 

Clearly, the \index{ultradifferential}$\o$-differential $\dd_p^\o f$ of a locally Lipschitz subfunction $f$ is defined and Lipschitz at each point $p\in \Dom f$.

\section{Remarks}
\label{page:upper-angle}

Spaces with defined angles include $\CAT{}{}$ and $\Alex{}{}$ spaces;
see \ref{angle} and \ref{cor:monoton-cba:angle=inf}.

For general metric spaces, angles may not exist, \index{$\mangle$!$\mangle^\text{up}$}
and given a hinge $\hinge p x y$  it is more natural to consider the \index{angle!upper angle}\emph{upper angle}  defined by
\[\mangle^\text{up}\hinge p x y
\df
\limsup_{\bar x,\bar y\to p} \angk\kappa p{\bar x}{\bar y},\]
where $\bar x\in\mathopen{]}p x]$ and $\bar y\in\mathopen{]}p y]$.
The triangle inequality (\ref{claim:angle-3angle-inq}) holds for upper angles as well.
\chapter{Dimension theory}\label{ch:dim}

\section{Definitions}\label{sec:prelim:dim}

In this section, we give definitions of different types of dimension-like invariants of metric spaces and state general relations between them.
The proofs of most of the statements in this section can be found in the book of Witold Hurewicz and Henry Wallman \cite{hurewicz-wallman}; 
the rest follow directly from the definitions.

\begin{thm}{Hausdorff dimension}
\label{def:HausDim}\index{dimension!Hausdorff dimension}\index{Hausdorff dimension}
Let $\spc{X}$ be a metric space. 
Its Hausdorff dimension is defined as
\[\HausDim\spc{X}=\sup\set{\alpha\in\RR}{\HausMes_\alpha(\spc{X})>0},\]
where $\HausMes_\alpha$ denotes the $\alpha$-dimensional Hausdorff measure.
\end{thm}

Let $\spc{X}$ be a metric space and $\{V_\beta\}_{\beta\in\IndexSet[2]}$
 be an open cover of $\spc{X}$.
Let us recall two notions in general topology:
\begin{itemize}

\item The \index{order of a cover}\emph{order} of $\{V_\beta\}$ is the supremum of all integers $n$ such that there is a collection of $n+1$ elements of $\{V_\beta\}$ with nonempy intersection.

\item An open cover $\{W_\alpha\}_{\alpha\in\IndexSet}$ of $\spc{X}$ is called a \index{refinement of a cover}\emph{refinement} of  $\{V_\beta\}_{\beta\in\IndexSet[2]}$ if for any $\alpha\in\IndexSet$ there is $\beta\in\IndexSet[2]$ such that $W_\alpha\subset V_\beta$.

\end{itemize}

\begin{thm}{Topological dimension}\label{def:TopDim}\index{dimension!topological dimension}\index{topological dimension}
Let $\spc{X}$ be a metric space. 
The topological dimension of $\spc{X}$ is defined to be the minimum of nonnegative integers $n$ 
such that for any open cover of $\spc{X}$ there is a finite open refinement with order~$n$.

If no such $n$ exists, the topological dimension of $\spc{X}$ is infinite.

The topological dimension of $\spc{X}$ will be denoted by $\TopDim\spc{X}$.
\end{thm}

The invariants satisfying the following two statements \ref{dim-axiom-norm} and \ref{dim-axiom-sigma} are commonly called ``dimension'';
for that reason we call these statements axioms.

\begin{thm}{Normalization axiom}
\label{dim-axiom-norm} For any $m\in\ZZ_{\ge0}$,
\[\TopDim\EE^m=\HausDim\EE^m=m.\]

\end{thm}

\begin{thm}{Cover axiom}\label{dim-axiom-sigma} 
If $\{A_n\}_{n=1}^\infty$ is a countable closed cover of $\spc{X}$, then
\begin{align*}
\TopDim \spc{X}&=\sup\nolimits_n\{\TopDim A_n\},
\\
\HausDim \spc{X}&=\sup\nolimits_n\{\HausDim A_n\}.
\end{align*}

\end{thm}

\parbf{On product spaces.} 
Recall that the direct product $\spc{X}\times\spc{Y}$ of metric spaces $\spc{X}$ and $\spc{Y}$ is defined in Section \ref{sec:metric spaces}. 
Direct product satisfies the following two inequalities:
\begin{align*}
\TopDim  (\spc{X}\times\spc{Y})
&\le 
\TopDim \spc{X}+ \TopDim\spc{Y}
\intertext{and}
\HausDim (\spc{X}\times\spc{Y})
&\ge 
\HausDim \spc{X}+ \HausDim\spc{Y}.
\end{align*}

These inequalities might be strict.
For  topological dimension,  strict inequality holds for a pair of Pontryagin surfaces \cite{pontyagin-surface}.
For Hausdorff dimension, an example was constructed by Abram Besicovitch and Pat Moran \cite{besicovitch-moran}.

\medskip
 
The following theorem follows from \cite[theorems V 8 and VII 2]{hurewicz-wallman}.

\begin{thm}{Szpilrajn's theorem}\label{thm:szpilrajn} 
Let $\spc{X}$ be a separable metric space.
Assume $\TopDim\spc{X}\ge m$. Then $\HausMes_m \spc{X}>0$.

In particular, 
$\TopDim\spc{X}\le\HausDim\spc{X}$.
\end{thm}

In fact it is true that for any separable metric space $\spc{X}$ we have
\[\TopDim\spc{X}=\inf\{\HausDim\spc{Y}\},\]
where the infimum is taken over all metric spaces $\spc{Y}$  homeomorphic to $\spc{X}$.

\begin{thm}{Definition}
Let $\spc{X}$ be a metric space
and $F\:\spc{X}\to\RR^m$ be  a continuous map.
A point $\bm{z}\in \Im F$ is called a \emph{stable value} of $F$
if there is $\eps>0$ such that $\bm{z}\in\Im F'$ 
for any \emph{$\eps$-close} to $F$ continuous map $F'\:\spc{X}\to\RR^m$,
that is, $|F'(x)-F(x)|<\eps$ for all $x\in \spc{X}$.
\end{thm}

The next theorem follows from \cite[theorems VI 1$\&$2]{hurewicz-wallman}.
(This theorem also holds for non-separable metric spaces \cite{nagata}, \cite[3.2.10]{engelking}). 

\begin{thm}{Stable value theorem}\label{thm:stable-value}
Let $\spc{X}$ be a separable metric space.
Then $\TopDim\spc{X}\ge m$ if and only if there is a map $F\:\spc{X}\to\RR^{m}$ with a stable value.
\end{thm}

\begin{thm}{Proposition}\label{thm:HausDim+Lip}
Suppose $\spc{X}$ and $\spc{Y}$ are metric spaces 
and $\map \:\spc{X}\to \spc{Y}$ satisfies
\[\dist{\map (x)}{\map (x')}{}\ge \eps\cdot\dist[{{}}]{x}{x'}{}\]
for fixed $\eps>0$ and any pair $x,x'\in \spc{X}$.
Then
\[\HausDim \spc{X}\le \HausDim \spc{Y}.\]

In particular, if there is a Lipschitz onto map $\spc{Y}\to \spc{X}$, then  
\[\HausDim \spc{X}\le \HausDim \spc{Y}.\]

\end{thm}

\section{Linear dimension} 

In addition to $\HausDim$ and $\TopDim$, 
we will use the so-called linear dimension}.
It will be applied only to  Alexandrov spaces and to their open subsets (in cases both of curvature bounded below and curvature bounded above).
As we shall see, in all these cases $\LinDim$  behaves nicely and  is easy to work with.

Recall that a \index{cone map}\emph{cone map} is a map between cones respecting the cone multiplication.

\begin{thm}{Definition of linear dimension}\label{def:lin-dim}\index{dimension!linear dimension}
Let $\spc{X}$ be a metric space with defined angles. 
The \index{linear dimension}\emph{linear dimension} of $\spc{X}$ (denoted by $\LinDim\spc{X}$\index{$\LinDim$}) is defined as the exact upper bound on $m\in\ZZ_{\ge0}$
such that there is a distance-preserving cone embedding $\EE^m\hookrightarrow \T_p\spc{X}$
for some $p\in \spc{X}$; here $\EE^m$ denotes the $m$-dimensional Euclidean space 
and $\T_p\spc{X}$ denotes the tangent space of $\spc{X}$ at $p$ (defined in Section~\ref{sec:tangent-space+directions}).
\end{thm}

Note that $\LinDim$ takes values in $\ZZ_{\ge0}\cup\{\infty\}$.
 
The linear dimension $\LinDim$ has no immediate relations to $\HausDim$ and $\TopDim$.
Also, $\LinDim$ does not satisfy the cover axiom (\ref{dim-axiom-sigma}).
Note that
\[\LinDim(\spc{X}\times \spc{Y})
=
\LinDim\spc{X}+ \LinDim\spc{Y}
\eqlbl{eq:inverse-product-axiom}\] 
for any two metric spaces $\spc{X}$ and $\spc{Y}$ with defined angles. 

The following exercise is based on a construction of Thomas Foertsch and Viktor Schroeder \cite{schroeder-foetch};
it shows that the condition on existence of  angles in \ref{eq:inverse-product-axiom} cannot be removed.

\begin{thm}{Exercise}\label{ex:schroeder-foetch}
Construct metrics $\rho_1$ and $\rho_2$ on $\RR^{10}$ defined by norms, such that $(\RR^{10},\rho_i)$ do \emph{not} contain an isometric copy of $\EE^2$ but
$(\RR^{10},\rho_1)\times (\RR^{10},\rho_2)$ has an isometric copy of $\EE^{10}$.
\end{thm}

\parbf{Remarks.}
Linear dimension was first introduced by Conrad Plaut \cite{plaut:survey}
under the name \index{local dimension}\index{dimension!local dimension}\emph{local dimension}. 
\index{geometric dimension}\index{dimension!geometric dimension}\emph{Geometric dimension}, introduced  by Bruce Kleiner \cite{kleiner} is closely related; 
it coincides 
 with the linear dimension for $\Alex{}$ and $\CAT{}$ spaces.

One can extend the definition to arbitrary metric spaces.
To do this one should modify  the definition of tangent space
and take an arbitrary $n$-dimensional Banach space instead of the Euclidean $n$-space.
For Alexandrov spaces (either $\Alex{}$ or $\CAT{}$) this modification is equivalent to our definition.

\part{Fundamentals}

\chapter{Fundamentals of curvature bounded below}
\chaptermark{Fundamentals of CBB}

\section{Four-point comparison} \label{sec:angle}

Recall (Section~\ref{sec:mod-tri/angles}) that the model angle $\angk\kappa p{x}{y}$ is defined if 
\[\dist{p}{x}{}+\dist{p}{y}{}+\dist{x}{y}{}<\varpi\kappa;\]
here $\varpi\kappa$ denotes the diameter of model space $\Lob2\kappa$.

\index{$\Alex{}$}
\begin{thm}{Four-point comparison}
\label{df:1+3}
A quadruple of points $p,x^1,x^2,x^3$ in a metric space satisfies \index{$\Alex{}$!$\Alex\kappa$ comparison}\emph{$\Alex\kappa$ comparison} 
if 
\[\angk\kappa p{x^1}{x^2}
+\angk\kappa p{x^2}{x^3}
+\angk\kappa p{x^3}{x^1}\le 2\cdot\pi.\eqlbl{Yup-kappa}\]
or at least one of the model angles $\angk\kappa p{x^i}{x^j}$ is not defined.
\end{thm}

\begin{thm}{Definition}
\label{df:cbb1+3}
Let $\spc{L}$ be a metric space.

\begin{subthm}{}
$\spc{L}$ is 
\index{$\Alex{}$!$\Alex{\kappa}$ space}
$\Alex{\kappa}$
if any quadruple $p,x^1,x^2,x^3\in \spc{L}$ satisfies  $\Alex\kappa$ comparison.
\end{subthm}

\begin{subthm}{}
$\spc{L}$ is 
\index{$\Alex{}$!locally $\Alex{\kappa}$ space}
\emph{locally $\Alex{\kappa}$} 
if any point $q\in \spc{L}$ admits a neighborhood $\Omega\ni q$ such that any quadruple $p,x^1,x^2,x^3\in \Omega$ satisfies  $\Alex\kappa$ comparison.
\end{subthm}

\begin{subthm}{}
$\spc{L}$  is a  
\index{$\Alex{}$!$\Alex{}$ space}
$\Alex{}$ space if  $\spc{L}$  is $\Alex{\kappa}$ for some $\kappa\in\RR$.
\end{subthm}
\end{thm}

\parbf{Remarks} 
\begin{itemize}

\item $\Alex\kappa$ length spaces are often called \emph{spaces with curvature $\ge\kappa$ in the sense of Alexandrov}.  These spaces will usually be denoted by $\spc{L}$, for $\spc{L}$ower curvature bound.

\item In the definition of $\Alex{\kappa}$, when $\kappa>0$ most authors assume in addition that the diameter is at most the model diameter~$\varpi\kappa$.
For a complete length space, the latter means that it is not isometric to one of the exceptional spaces, see \ref{diam-k>0}. 
We do not make this assumption.
In particular, we consider the real line to have curvature $\ge 1$.

\item If $\kappa<\Kappa$, then any complete length $\Alex{\Kappa}$ space is $\Alex{\kappa}$.
Moreover directly from the definition it follows that if $\Kappa\le 0$, then any $\Alex{\Kappa}$ space is $\Alex{\kappa}$.
However, in the case $\Kappa>0$ the latter statement does not hold and the former statement is not trivial; it will be proved in \ref{cor:CAT>k-sence}.
\end{itemize}

\begin{thm}{Exercise}\label{ex:(3+1)-expanding}
Let $\spc{L}$ be a metric space and $\kappa\le 0$.
Show that $\spc{L}$ is $\Alex\kappa$
if for any quadruple of points $p,x^1,x^2,x^3\in \spc{L}$ 
there is a quadruple of points $q,y^1,y^2,y^3\in\Lob{2}\kappa$
such that 
\[\dist{p}{x^i}{}=\dist{q}{y^i}{} 
\quad \text{and}\quad \dist{x^i}{x^j}{}\le\dist{y^i}{y^j}{}\] 
for all $i$ and $j$.
\end{thm}

The exercise above is a special case of (1+\textit{n})-point comparison (\ref{thm:pos-config}).

Recall that $\o$ denotes a selective ultrafilter on $\NN$, which is fixed once and  for all.
The following proposition follows directly from the definition of $\Alex\kappa$ comparison and the definitions of $\o$-limit and $\o$-power given in Section~\ref{ultralimits}.

\begin{thm}{Proposition}\label{prp:A^omega}
Let $\spc{L}_n$ be a $\Alex{\kappa_n}$ space for each $n$.
Assume $\spc{L}_n\to \spc{L}_\o$ 
and $\kappa_n\to\kappa_\o$ as $n\to\o$.
Then $\spc{L}_\o$ is $\Alex{\kappa_\o}$.

Moreover, a metric space $\spc{L}$ is $\Alex\kappa$ if and only if so is
its ultrapower $\spc{L}^\o$.
\end{thm}

\begin{thm}{Theorem}\label{thm:submetry-CBB}
Let $\spc{L}$ be a $\Alex{\kappa}$ space, $\spc{M}$ be a metric space, and $\sigma\:\spc{L}\to\spc{M}$ be a submetry.
Assume $p,x^1,x^2,x^3$ is a quadruple of points in $\spc{M}$ such that 
$\dist{p}{x^i}{}<\tfrac{\varpi\kappa}2$ for any $i$.
Then the  quadruple satisfies $\Alex{\kappa}$ comparison.

In particular, 
\begin{subthm}{}
The space $\spc{M}$ is locally $\Alex{\kappa}$.
Moreover, any open ball of radius $\tfrac{\varpi\kappa}4$ in $\spc{M}$ is $\Alex{\kappa}$.
\end{subthm}

\begin{subthm}{}
If $\kappa\le 0$, then  $\spc{M}$ is $\Alex{\kappa}$.
\end{subthm}
\end{thm}

Corollary~\ref{cor:submetry-cbb} gives a stronger statement; it states that if $\spc{L}$ is a complete length space, then $\spc{M}$ is always $\Alex{\kappa}$.
The theorem above together with Proposition~\ref{prop:submet/G}
imply the following:

\begin{thm}{Corollary}\label{thm:CBB/G}
Assume that  the group $G$ acts isometrically on a $\Alex{\kappa}$ space $\spc{L}$ and has closed orbits.
Then the quotient space $\spc{L}/G$ is locally $\Alex\kappa$. 
\end{thm}

\begin{thm}{Example}
If $G< \O(n+1)$ is a closed subgroup, then $\SS^n/G$ is $\Alex{1}$ and $\R^{n+1}/G$ is $\Alex{0}$.
\end{thm}

\parit{Proof of \ref{thm:submetry-CBB}.}
Fix a quadruple of points $p,x^1,x^2,x^3\in \spc{M}$ such that 
$\dist{p}{x^i}{}<\tfrac{\varpi\kappa}2$ for any $i$.
Choose an arbitrary $\hat p\in \spc{L}$ such that $\sigma(\hat{p})=p$.

Since $\sigma$ is submetry, we can choose the points $\hat{x}^1,\hat{x}^2,\hat{x}^3\in \spc{L}$ such that $\sigma(\hat x_i)=x_i$ and
\[\dist{p}{x^i}{\spc{M}}
\lege
\dist{\hat{p}}{\hat{x}^i}{\spc{L}}
\pm\delta\]
for all $i$ and any fixed $\delta>0$.

Note that 
\[\dist{x^i}{x^j}{\spc{M}}
\le
\dist{\hat{x}^i}{\hat{x}^j}{\spc{L}}\le \dist{p}{x^i}{\spc{M}}+\dist{p}{x^j}{\spc{M}}+2\cdot\delta\]
for all $i$ and $j$.

Since $\dist{p}{x^i}{}<\tfrac{\varpi\kappa}2$,
we can choose $\delta>0$ above so that the angles $\angk\kappa {\hat{p}}{\hat{x}^i}{\hat{x}^j}$ are defined.
Moreover, given $\eps>0$, the value $\delta$ can be chosen in such a way that the inequality
\[\angk\kappa p{x^i}{x^j}
<
\angk\kappa {\hat{p}}{\hat{x}^i}{\hat{x}^j}+\eps
\eqlbl{eq:angles-M-L}\]
holds for all $i$ and $j$.

By $\Alex\kappa$ comparison in $\spc{L}$,
we have
\[\angk\kappa {\hat{p}}{\hat{x}^1}{\hat{x}^2}
+\angk\kappa {\hat{p}}{\hat{x}^2}{\hat{x}^3}
+\angk\kappa {\hat{p}}{\hat{x}^3}{\hat{x}^1}
\le 
2\cdot\pi.\]
Applying  \ref{eq:angles-M-L}, 
we get 
\[\angk\kappa p{x^1}{x^2}
+\angk\kappa p{x^2}{x^3}
+\angk\kappa p{x^3}{x^1}< 2\cdot\pi+3\cdot\eps.\]
Since $\eps>0$ is arbitrary we have 
\[\angk\kappa p{x^1}{x^2}
+\angk\kappa p{x^2}{x^3}
+\angk\kappa p{x^3}{x^1}\le 2\cdot\pi;\]
that is,
the $\Alex\kappa$ comparison holds for this quadruple in  $\spc{M}$.
\qeds

\section{Geodesics}

Recall that general complete length spaces might have no geodesics;
see Exercise~\ref{ex:no-geod}.

\begin{thm}{Exercise}\label{ex:nongeod-cbb}
Construct a complete length $\Alex0$ space that is not geodesic.
\end{thm}

We are going to show that all complete length $\Alex{}$ spaces have plenty of geodesics in the following sense. Recall that a subset of a topological space is called G-delta if it is a countable intersection of open sets.

\begin{thm}{Definition}\label{def:alm-geod}
A metric space $\spc{X}$ is called \index{G-delta geodesic space}\emph{G-delta geodesic} 
if for any point $p\in \spc{X}$ there is a dense G-delta set $W_p\subset\spc{X}$ such that for any $q\in W_p$ there is a geodesic $[p q]$.

A metric space $\spc{X}$ is called {}\emph{locally G-delta geodesic} 
if for any point $p\in \spc{X}$ there is a G-delta set $W_p\subset\spc{X}$ such that
$W_p$ is dense in a neighborhood of $p$ 
and for any $q\in W_p$ there is a geodesic $[p q]$.
\end{thm}

\begin{thm}{Definition}\label{def:straight}
Let $\spc{X}$ be a metric space 
and $p\in \spc{X}$.
A point $q\in \spc{X}$ is called \index{straight point}\emph{$p$-straight} (briefly, $q\in \Str(p)$\index{$\Str(p)$}) if
\[\limsup_{r\to q}\frac{\dist{p}{r}{}-\dist{p}{q}{}}{\dist{q}{r}{}}=1.\]

For an array of points $x^1,x^2,\dots,x^\kay$, 
we use the notation
\[\Str(x^1,x^2,\dots,x^\kay)=\bigcap_i\Str(x^i).\]
\end{thm}

\begin{thm}{Theorem}\label{thm:almost.geod}
Let $\spc{L}$ be a complete length $\Alex{}$ space and $p\in \spc{L}$.
Then the set  $\Str(p)$ is a dense G-delta set.
Moreover, for any $q\in \Str(p)$ there is a unique geodesic $[p q]$.

In particular, $\spc{L}$ is G-delta geodesic.
\end{thm}

This theorem was proved by by Conrad Plaut \cite[Th. 27]{plaut:survey}.

\parit{Proof.}
Given a positive integer $n$, 
consider the set $\Omega_n$ of all points $q\in \spc{L}$ such that
\[(1-\tfrac{1}{n})\cdot\dist[{{}}]{q}{r}{}<
\dist{p}{r}{}-\dist{p}{q}{}
<\tfrac{1}{n}\]
for some $r\in \spc{L}$.
Clearly $\Omega_n$ is open; 
let us show that $\Omega_n$ is dense in $\spc{L}$.

Assuming the contrary, there is a point $x\in \spc{L}$ such that 
\[\oBall(x,\eps)\cap \Omega_n=\emptyset\]
for $\eps>0$.
Since $\spc{L}$ is a length space, 
for any $\delta>0$, there exists a point $y\in \spc{L}$ such that 
\[\dist{x}{y}{}<\tfrac\eps2+\delta
\quad\text{and}\quad
\dist{p}{y}{}<\dist{p}{x}{}-\tfrac\eps2+\delta.
\]
If $\eps$ and $\delta$ are sufficiently small, then
\[(1-\tfrac{1}{n})\cdot\dist[{{}}]{y}{x}{}
<
\dist{p}{x}{}-\dist{p}{y}{}<\tfrac{1}{n};\] 
that is, $y\in\Omega_n$, 
a contradiction.

Note that $\Str(p)=\bigcap_{n}\Omega_n$;
therefore, $\Str(p)$ is a dense G-delta set.

Assuming $q\in \Str(p)$,
let us show that there is a unique geodesic connecting $p$ and $q$.
Note that it is sufficient to show that for all sufficiently small
$t>0$ there is a unique point $z$ such that 
\[t
=
\dist{q}{z}{}
=
\dist{p}{q}{}-\dist{p}{z}{}.
\eqlbl{eq:thm:connect-1}\]

First let us show uniqueness. 
Assume $z$ and $z'$ both satisfy \ref{eq:thm:connect-1}.
Take a sequence $r_n\to q$ such that 
\[\frac{\dist{p}{r_n}{}-\dist{p}{q}{}}{\dist[{{}}]{q}{r_n}{}}
\to 1.\] 
By the triangle inequality, 
\[\dist{z}{r}{}-\dist{z}{q}{},
\quad
\dist{z'}{r}{}-\dist{z'}{q}{}
\ge 
\dist{p}{r}{}-\dist{p}{q}{};\] 
thus, as $n\to\infty$,
\[\frac{\dist{z}{r_n}{}-\dist{z}{q}{}}{\dist{q}{r_n}{}},
\quad 
\frac{\dist{z'}{r_n}{}-\dist{z'}{q}{}}{\dist{q}{r_n}{}}
\to 1.\]
Therefore $\angk\kappa q z{r_n}\to\pi$ and $\angk\kappa q{z'}{r_n}\to\pi$.
(Here we use that $t$ is small, otherwise if $\kappa>0$ the angles might be undefined.)
 
From $\Alex\kappa$ comparison (\ref{df:cbb1+3}), $\angk\kappa q z {z'}=0$ and thus $z=z'$.

The proof of existence is similar.
Choose a sequence $r_n$ as above.
Since $\spc{L}$ is a complete length space, 
there is a sequence $z_\kay\in \spc{L}$ such that $\dist{q}{z_\kay}{}\to t$ and $\dist{p}{q}{}-\dist{p}{z_\kay}{}\to t$ as $\kay\to\infty$.
Then 
\[
\lim_{n\to\infty}
\lim_{\kay\to\infty}
\angk\kappa q{z_\kay}{r_n}
=\pi.\] 
Thus, for any $\eps>0$ and sufficiently large $n,\kay$, we have $\angk\kappa q{z_\kay}{r_n}\z>\pi-\eps$.
From $\Alex\kappa$ comparison (\ref{df:cbb1+3}), for all large $\kay$ and $j$, we have $\angk\kappa q{z_\kay}{z_j}\z<2\cdot\eps$ and thus 
\[\dist{z_\kay}{z_j}{}<\eps\cdot\Const(\kappa,t);\]
that is, $z_n$ is a Cauchy sequence, and its limit $z$ satisfies \ref{eq:thm:connect-1}.
\qeds

\begin{thm}{Exercise}\label{ex:almost.geod}
Let $\spc{L}$ be a complete length $\Alex{}$ space and $A\subset\spc{L}$ be a closed subset.
Show that there is a dense G-delta set $W\subset\spc{L}$ such that
for any $q\in W$, there is a unique geodesic $[pq]$ with
$p\in A$ that realizes the distance from $q$ to $A$; that is, $\dist{p}{q}{}=\distfun{A}{q}{}$.
\end{thm}

\begin{thm}{Exercise}\label{ex:G-delta-not-thru}
Construct a complete length $\Alex{}$ space $\spc{L}$
with an everywhere dense G-delta set $A$
such that 
$A\cap \mathopen{]}xy\mathclose{[}=\emptyset$
for any geodesic $[xy]$ in $\spc{L}$. 
\end{thm}


\section{More comparisons}\label{sec:more-angles}

The following theorem makes it easier to use Euclidean intuition in the Alexandrov 
setting.
\begin{thm}{Theorem}
\label{thm:defs_of_alex} 
If $\spc{L}$ is a $\Alex\kappa$ space, 
then the following conditions hold for all $p,x,y\in \spc{L}$, provided the model triangle $\modtrig\kappa(p x y)$ is defined.

\begin{subthm}{2-sum} 
(adjacent angle comparison\index{comparison!adjacent angle comparison}) for any geodesic $[x y]$ and $z\in \mathopen{]}x y\mathclose{[}$, $z\ne p$ we have
\[\angk\kappa z p x
+\angk\kappa z p y\le \pi.\]
\end{subthm}

\begin{subthm}{point-on-side}
(\index{comparison!point-on-side comparison}point-on-side comparison)
for any geodesic $[x y]$ and $z\in \mathopen{]}x y\mathclose{[}$, we have
\[\angk\kappa x p y\le\angk\kappa x p z;\]
or, equivalently, 
\[\dist{\tilde p}{\tilde z}{}\le \dist{p}{z}{},\]
where $\trig{\tilde p}{\tilde x}{\tilde y}=\modtrig\kappa(p x y)$, $\tilde z\in\mathopen{]} \tilde x\tilde y\mathclose{[}$, $\dist{\tilde x}{\tilde z}{}=\dist{x}{z}{}$.
\end{subthm}

\begin{subthm}{angle}(hinge comparison\index{comparison!hinge comparison})
\index{hinge comparison}
for any hinge $\hinge x p y$, the angle 
$\mangle\hinge x p y$ is defined and 
\[\mangle\hinge x p y\ge\angk\kappa x p y,\]
or equivalently
\[\side\kappa \hinge x p y\ge\dist{p}{y}{}.\]
Moreover, 
\[\mangle\hinge z p y + \mangle\hinge z p x\le\pi\]
for any two adjacent hinges $\hinge z p y$ and $\hinge z p x$.
\end{subthm}

Moreover, in each case, the converse holds if $\spc{L}$ is G-delta geodesic.
That is, if one of the conditions \ref{SHORT.2-sum}, \ref{SHORT.point-on-side}), or  \ref{SHORT.angle} holds in a  G-delta geodesic space $\spc{L}$, then $\spc{L}$ is $\Alex\kappa$.
\end{thm}

A slightly stronger form of \ref{SHORT.angle} is given in \ref{lem:devel-glob}.
See also Problem~\ref{open:hinge-}.

\parit{Proof; \ref{SHORT.2-sum}).}
Since $z\in \mathopen{]}x y\mathclose{[}$, we have $\angk\kappa z x y=\pi$. 
Thus, $\Alex\kappa$ comparison
\[\angk\kappa z x y
+\angk\kappa z p x
+\angk\kappa z p y\le2\cdot\pi\]
implies
\[\angk\kappa z p x
+\angk\kappa z p y
\le\pi.\]

\parit{\ref{SHORT.2-sum}$\Leftrightarrow$\ref{SHORT.point-on-side}.} 
Follows from Alexandrov's lemma (\ref{lem:alex}).

\parit{\ref{SHORT.2-sum}$+$\ref{SHORT.point-on-side}$\Rightarrow$\ref{SHORT.angle}.} 
From \ref{SHORT.point-on-side} we get that for $\bar p\in\mathopen{]}xp]$ and $\bar y\in\mathopen{]}xy]$, the function $(\dist{x}{\bar p}{},\dist{x}{\bar y}{})\mapsto\angk\kappa x{\bar p}{\bar y}$ is nonincreasing in each argument.
In particular, 
$\mangle\hinge x p y\z=\sup\{\angk\kappa x{\bar p}{\bar y}\}$
 is defined and is
at least $\angk\kappa x p y$.

\begin{wrapfigure}{o}{30 mm}
\vskip-0mm
\centering
\includegraphics{mppics/pic-805}
\end{wrapfigure}

From above and \ref{SHORT.2-sum}, it follows that 
\[\mangle\hinge z p y + \mangle\hinge z p x\le\pi.\]

\parit{Converse.}
Assume first that $\spc{L}$ is geodesic.
Consider a point  $w\in \mathopen{]} p z \mathclose{[}$ close to $p$.
From \ref{SHORT.angle}, it follows that 
\[\mangle\hinge w x z+ \mangle\hinge w x{p}\le\pi\quad \text{and}\quad \mangle\hinge w y z + \mangle\hinge w y{p}\le\pi.\]
Since $\mangle\hinge w x y\le \mangle\hinge w x p +\mangle\hinge w y{p}$ (see \ref{claim:angle-3angle-inq}), we get 
\[\mangle\hinge w x z+ \mangle\hinge w y z +\mangle\hinge w x y
\le
2\cdot\pi.\]
Applying the first inequality in \ref{SHORT.angle}, 
\[\angk\kappa w x z
+ \angk\kappa w y z 
+\angk\kappa w x y
\le
2\cdot\pi.\]
Passing to the limits  $w\to p$, we have
\[\angk\kappa p x z 
+ \angk\kappa p y z 
+\angk\kappa p x y
\le
2\cdot\pi.\]

If $\spc{L}$ is only G-delta geodesic, we can apply the  above arguments to sequences of points $p_n,w_n\to p$, $x_n\to x$, $y_n\to y$ such that $[p_nz]$ exists, $w_n\in\mathopen{]}z p_n\mathclose{[}$ and  $[x_nw_n]$, $[y_n w_n]$ exist,  and then pass to the limit as $n\to\infty$.
\qeds

\begin{thm}{Exercise}\label{mink+alex=euclid} 
Let $\spc{L}$ be  $\RR^m$ with a metric defined by a norm.
Show that $\spc{L}$ is a complete length $\Alex{}$ space if and only if $\spc{L}\iso\EE^m$.
\end{thm}

\begin{thm}{Exercise}\label{ex:cbb-geod-overlap}
Assume $\spc{L}$ is a complete length $\Alex{}$ space, and $[px]$, $[py]$ be two geodesics  in the same geodesic direction $\xi\in \Sigma'_p$.
Show that 
\[[px]\subset [py]\quad \text{or}\quad [px]\supset [py].\]

\end{thm}

\begin{thm}{Angle-sidelength  monotonicity}\label{cor:monoton} 
Let $p$, $x$ , $y$ be points in a complete length $\Alex{\kappa}$ space $\spc{L}$.
Suppose that the model triangle
$\modtrig\kappa(p x y)$ is defined
and there is a geodesic $[xy]$.
Then for $\bar y\in\mathopen{]}xy]$ the function 
\[\dist{x}{\bar y}{}\mapsto \angk\kappa x p{\bar y}\] 
is nonincreasing.

In particular, if a geodesic $[x p]$ exists and $\bar p\in \mathopen{]}x p]$, then
\begin{subthm}{cor:monoton:2-sides}
the function 
\[(\dist{x}{\bar y}{},\dist{x}{\bar p}{})\mapsto \angk\kappa x {\bar p}{\bar y}\] is nonincresing in each argument
\end{subthm}
 
\begin{subthm}{cor:monoton:sup}
The angle $\mangle\hinge{x}{p}{y}$ is defined and 
\[\mangle\hinge{x}{p}{y}
=
\sup\set{\angk\kappa x {\bar p}{\bar y}}{
\bar p\in\mathopen{]}xp],\ 
\bar y\in\mathopen{]}xy]}.\]

\end{subthm}
\end{thm}

The proof is contained in the first part of \ref{SHORT.2-sum}$+$\ref{SHORT.point-on-side}$\Rightarrow$\ref{SHORT.angle} of the proof above.

\begin{thm}{Exercise}\label{ex:equality-alexlemma} 
Let $\spc{L}$ be a $\Alex{\kappa}$ space,
$p,x,y\in \spc{L}$
and $v,w\in \mathopen{]}xy\mathclose{[}$.
Prove that  
\[
\angk\kappa xyp=\angk\kappa xvp
\quad \iff\quad 
\angk\kappa xyp=\angk\kappa xwp.
\]

\end{thm}

\begin{thm}{Advanced exercise} \label{ex:urysohn}
Construct a geodesic space $\spc{X}$ that is not $\Alex0$, 
but meets the following condition: for any 3 points $p,x,y\in \spc{X}$ there is a geodesic $[x y]$ such that for any $z\in \mathopen{]}x y\mathclose{[}$
\[\angk{0}{z}{p}{x}+\angk{0}{z}{p}{y}
\le
\pi.\]
\end{thm}

\begin{thm}{Advanced exercise}\label{ex:lebedeva-petrunin}
Let $\spc{L}$ be a complete length space such that for any quadruple $p,x,y,z\in\spc{L}$ 
the following inequality holds
\[\dist[2]{p}{x}{}+\dist[2]{p}{y}{}+\dist[2]{p}{z}{}
\ge
\tfrac13\cdot
\left[
\dist[2]{x}{y}{}+\dist[2]{y}{z}{}+\dist[2]{z}{x}{}
\right].\eqlbl{eq:berg-nikolaev}\]
Prove that $\spc{L}$ is $\Alex0$.

Construct a 4-point metric space $\spc{X}$ that satisfies inequality \ref{eq:berg-nikolaev} for any relabeling of its points by $p,x,y,z$, such that $\spc{X}$ is not $\Alex{0}$.
\end{thm}

Assume that for a given triangle $\trig{x^1}{x^2}{x^3}$ in a metric space its $\kappa$-model triangle
$\trig{\tilde x^1}{\tilde x^2}{\tilde x^3}=\modtrig\kappa({x^1}{x^2}{x^3})$ is defined.
We say the triangle $\trig{x^1}{x^2}{x^3}$ is \emph{$\kappa$-thick} \index{$\kappa$-thick} if the natural map (see definition \ref{def:k-thin}) 
$\trig{\tilde x^1}{\tilde x^2}{\tilde x^3}\to \trig{x^1}{x^2}{x^3}$ is distance non contracting.

\begin{thm}{Exercise}\label{ex:fat-triangle}
Prove that any triangle with perimeter $<\varpi\kappa$ 
in a $\Alex{\kappa}$ space is \index{$\kappa$-thick}\emph{$\kappa$-thick}.
\end{thm} 

\begin{thm}{Exercise}\label{ex:busemann}

\begin{subthm}{}
Show that any $\Alex0$ space $\spc{L}$ satisfies the following condition:
for any three points $p,q,r\in \spc{L}$, if $\bar q$ and $\bar r$ are midpoints of geodesics $[p q]$ and $[p r]$ respectively, then $2\cdot\dist[{{}}]{\bar q}{\bar r}{}\ge\dist{q}{r}{}$.
\end{subthm}

\begin{subthm}{} Show that there is a metric on $\RR^2$ defined by a norm that satisfies the above condition, but is not $\Alex0$.
\end{subthm}

\end{thm}

\parbf{Remarks.} 
Monotonicity of the model angle with respect to adjacent sidelengths (\ref{cor:monoton}) was named the \index{convexity property}\emph{convexity property} by Alexandrov.


\section{Function comparison} \label{sec:func-comp-CBB}

In this section we will translate the angle comparison definitions (Theorem~\ref{thm:defs_of_alex}) 
to a concavity-like property of the distance functions as defined in Section~\ref{sec:conv-fun}.
This is a conceptual step ---
we reformulate a global geometric condition into an infinitesimal condition on distance functions.

\begin{thm}{Theorem}\label{thm:conc} 
Let $\spc{L}$ be a complete length space. 
Then the following 
statements are equivalent:

\begin{subthm}{main-def'} $\spc{L}$ is $\Alex\kappa$.
\end{subthm}

\begin{subthm}{comp-kappa}(function comparison\index{comparison!function comparison}) $\spc{L}$ is  G-delta geodesic and for any $p\in \spc{L}$, the function $f=\md\kappa\circ\distfun{p}{}{}$ satisfies the differential inequality
\[f''\le 1-\kappa\cdot f.\]
in $\oBall(p,\varpi\kappa)$.
\end{subthm}
\end{thm}

\begin{thm}{Corollary}
A complete G-delta geodesic space $\spc{L}$ is $\Alex{0}$ if and only if for any $p\in \spc{L}$, the function $\distfun[2]{p}{}{}\:\spc{L}\to\RR$ 
is $2$-concave.
\end{thm} 

\parit{Proof of \ref{thm:conc}.} 
Let $[x y]$ be a geodesic in $\oBall(p,\varpi\kappa)$ and $\ell=\dist{x}{y}{}$.
Consider the model triangle $\trig{\tilde p}{\tilde x}{\tilde y}=\modtrig\kappa(p x y)$.
Set \begin{align*} 
\tilde r(t)&=\dist{\tilde p}{\geod_{[\tilde x\tilde y]}(t)}{},
& 
r(t)&=\dist{p}{\geod_{[xy]}(t)}{}.                           \end{align*}
Clearly $\tilde r(0)=r(0)$ and $\tilde r(\ell)=r(\ell)$. 
Set $\tilde f=\md\kappa\circ\tilde r$ and $f=\md\kappa\circ r$.
From \ref{md-diff-eq} we get that $\tilde f''=1-\kappa\cdot  \tilde f$.

Note that the point-on-side comparison (\ref{point-on-side}) for point $p$ and geodesic $[x y]$ is equivalent to $\tilde r\le r$.
Since $\md\kappa$ is increasing on $[0,\varpi\kappa)$, 
$\tilde r\le r$ is equivalent to $\tilde f\le f$.
The latter is Jensen's inequality (\ref{y''-mono}) for the function
$t\mapsto\md\kappa\dist[{{}}]{p}{\geod_{[x y]}(t)}{}$ on the interval $[0,\ell]$. 
Hence the result.
\qeds

Recall that Busemann functions are defined in Proposition \ref{prop:busemann}.

\begin{thm}{Exercise}\label{ex:busemann-CBB}
{\sloppy 
Let $\spc{L}$ be a complete length $\Alex{\kappa}$ space
and $\bus_\gamma\:\spc{L}\to \RR$ be the Busemann function for a half-line $\gamma\:[0,\infty)\to \spc{L}$.

}

\begin{subthm}{ex:busemann-CBB:a}
If $\kappa=0$, then the Busemann function $\bus_\gamma$ is  concave.
\end{subthm}

\begin{subthm}{ex:busemann-CBB:b}
If $\kappa=-1$, then the function 
$f=\exp\circ\bus_\gamma$ 
satisfies 
\[f''- f\le 0.\]
\end{subthm}

\end{thm}

Exercise~\ref{ex:busemann-CBA} is an analogous statement for upper curvature bound.

\section{Development}

In this section we reformulate the function comparison using a more geometric language based on the definition of development given below.

This definition appears in \cite{alexandrov:devel}
and an earlier form of it can be found in \cite{liberman}.
The definition is somewhat lengthy, but it defines a useful comparison object for a curve. 
Often it is easier to write proofs in terms of function comparison but think in terms of developments.

\begin{thm}{Lemma-definition}\label{lem:devel}\label{def:devel}
Let $\kappa\in \RR$, 
$\spc{X}$ be a metric space, 
$\gamma\:\II\to \spc{X}$ be a $1$-Lipschitz curve,
$p\in \spc{X}$,
and $\tilde p\in\Lob2{\kappa}$.
Assume $0<\dist{p}{\gamma(t)}{} < \varpi\kappa$ for
all $t\in \II$.
Then there exists a unique up to rotation curve
$\tilde \gamma\: \II\to \Lob2\kappa$, parametrized by arc-length, 
such that
$\dist{\tilde p}{\tilde \gamma(t)}{}\z=\dist{p}{\gamma(t)}{}$ for all $t$
and the direction of
$[\tilde p\tilde \gamma(t)]$ monotonically turns around $\tilde p$ counterclockwise as $t$ increases.

\smallskip

If $p$, $\tilde p$, $\gamma$, and $\tilde \gamma$ are as above,
then $\tilde \gamma$ is called the \index{development}\emph{$\kappa$-development} of $\gamma$ with respect to $p$; 
the point $\tilde p$ is called the \index{development!basepoint of a development}\emph{basepoint} of the development.
When we say that the $\kappa$-development of $\gamma$ with respect to $p$ is defined, we always assume that $0<\dist{p}{\gamma(t)}{} < \varpi\kappa$ for
all 
$t\in \II$.
\end{thm}

\parit{Proof.}
Consider the functions $\rho$, $\theta\:\II\to\RR$ defined as 
\begin{align*}
\rho(t)
&=\dist{p}{\gamma(t)}{},
&
\theta(t)
&=
\int\limits_{t_0}^{t}\frac{\sqrt{1-(\rho')^2}}{\sn\kappa\rho},
\end{align*}
where $t_0\in\II$ is a fixed number and $\int$ denotes Lebesgue integral.
Since $\gamma$ is $1$-Lipshitz, so is $\rho(t)$, and thus the function $\theta$ is defined and nondecreasing.

It is straightforward to check that $(\rho,\theta)$ uniquely describe $\tilde \gamma$ in polar coordinates on $\Lob2{\kappa}$ with center at $\tilde p$.
\qeds
                                                                                                                                                                                                                                                                                                                                                                        
We need the following analogs of sub- and super-graphs 
 and convex/concave functions, adapted to polar coordinates in $\Lob2{\kappa}$.

\begin{wrapfigure}{o}{38 mm}
\vskip-0mm
\centering
\includegraphics{mppics/pic-810}
\end{wrapfigure}

\begin{thm}{Definition}\label{def:convex-devel}
Let $\tilde \gamma\:\II\to\Lob2\kappa$ be a curve and $\tilde p\in\Lob2\kappa$ be such that there is a unique geodesic $[\tilde p\,\tilde \gamma(t)]$ for any $t\in\II$ and the direction of $[\tilde p\,\tilde \gamma(t)]$ turns monotonically as $t$ grows.

The set formed by all geodesics from  $\tilde p$ to the points on $\tilde \gamma$ is called the \index{development!subgraph/supergraph} \emph{subgraph} of $\tilde \gamma$ with respect to $\tilde p$.

The set of all points $\tilde x\in\Lob2{\kappa}$ such that a geodesic $[\tilde p\tilde x]$ intersects $\tilde \gamma$ is called the \emph{supergraph} of $\tilde \gamma$ with respect to $\tilde p$.

The curve $\tilde \gamma$ is called \index{convex/concave curve with respect to a point}\emph{convex (concave) with respect to} $\tilde p$ if the subgraph (supergraph) of $\tilde \gamma$ with respect to $\tilde p$ is convex.

The curve $\tilde \gamma$ is called 
\emph{locally convex (concave) with respect to $\tilde p$} 
if for any interior value $t_0$ in $\II$ there is a subsegment $(a,b)\subset\II$, $(a,b)\z\ni t_0$, such that the restriction $\tilde \gamma|_{(a,b)}$ is convex (concave) with respect to~$\tilde p$.
\end{thm}

Note that if $\kappa>0$, then the supergraph of a curve is the subgraph with respect to the opposite point. 

For developments,
all the notions above will be considered with respect to their basepoints.
In particular, if $\tilde \gamma$ is a development, we will say it is \emph{(locally) convex} if it is (locally) convex with respect to its basepoint.

{\sloppy 

\begin{thm}{Development comparison\index{comparison!development comparison}}\label{thm:devel} 
\index{development!comparison}
A complete G-delta geodesic space
$\spc{L}$ is $\Alex\kappa$ if and only if for any point $p\in \spc{L}$ and any geodesic $\gamma$ in $\oBall(p,\varpi\kappa)\setminus\{p\}$, 
its $\kappa$-development with respect to $p$ is convex.
\end{thm}

}

A simpler proof of the only-if part can be built on the adjacent angle comparison (\ref{2-sum}).
We use a longer proof since it also implies the short hinge lemma (\ref{lem:devel-glob}).

\parit{Proof; only-if part.}  
Let $\gamma\:[0,T]\to \oBall(p,\varpi\kappa)\setminus \{p\}$ be a unit-speed geodesic in $\spc{L}$.

Consider a fine partition 
\[0=t_0<t_1<\dots<t_n=T.\]
Set $x_i=\gamma(t_i)$ and choose a point 
\[p'\in \Str(x_0,x_1,\dots,x_n)\] 
sufficiently close to $p$;
recall that geodesics $[p'x_i]$ exist for all $i$
(see \ref{thm:almost.geod}).

{\sloppy 

Let us construct a chain of model triangles 
$\trig{\tilde p'}{\tilde x_{i-1}}{\tilde x_i}\z=\modtrig\kappa(p' x_{i-1}x_i)$ in such a way that direction $[\tilde p'\tilde x_i]$ turns counterclockwise as $i$ grows.
By the hinge comparison (\ref{angle}), we have
\[\begin{aligned}
\mangle\hinge {\tilde x_i}{\tilde x_{i-1}}{\tilde p'}
+\mangle\hinge{\tilde x_i}{\tilde x_{i+1}}{\tilde p'}
&=
\angk\kappa{x_i}{x_{i-1}}{p'}
+\angk\kappa{x_i}{x_{i+1}}{p'}
\le
\\
&\le
\mangle\hinge{x_i}{x_{i-1}}{p'}
+\mangle\hinge{x_i}{x_{i+1}}{p'}
\le
\\
&\le
\pi.
\end{aligned}
\eqlbl{eq:thm:devel*}
\]

}

\begin{wrapfigure}[4]{r}{31 mm}
\vskip-10mm
\centering
\includegraphics{mppics/pic-815}
\end{wrapfigure}

Further, since $\gamma$ is a unit-speed geodesic, we have 
\[\sum_{i=1}^n\dist{x_{i-1}}{x_i}{}\le \dist{p'}{x_0}{}+\dist{p'}{x_n}{}.\eqlbl{eq:thm:detr*}\]
Since $p'\notin\gamma$, the development comparison implies that $\tilde p'$ does not lie on the polygonal line $\tilde x_0\dots\tilde x_n$.

If $\kappa\le 0$, then \ref{eq:thm:detr*} implies that
\[\theta\df
\sum_{i=1}^n\mangle\hinge{\tilde p'}{\tilde x_i}{\tilde x_{i-1}}\le\pi.
\eqlbl{eq:thm:devel**}\]

In the case $\kappa>0$, the proof of \ref{eq:thm:devel**} requires more work.
Applying rescaling, we can assume that $\kappa=1$.
Since $\gamma$ lies in $B_\pi(p')$, the point-on-side comparison implies that the antipodal point of $\tilde p'$ does not lie on the polygonal line $\tilde x_0\dots\tilde x_n$.

Consider the space $L$ glued from the solid model triangles $\trig{\tilde p'}{\tilde x_{0}}{\tilde x_1},\z\dots,\trig{\tilde p'}{\tilde x_{n-1}}{\tilde x_n}$ along the corresponding sides.
Note that $\theta$ is the total angle of $L$ at $\tilde p'$.
We can assume that $L$ has nonempty interior.
Otherwise all the triangles are degenerate, and therefore $\theta=0$;
the latter holds since $\tilde p'$ is not on the polygonal line $\tilde x_0\dots\tilde x_n$.

Consider a minimizing geodesic $[\tilde x_0\tilde x_n]_L$.
By \ref{eq:thm:detr*} we may assume that $\tilde p'\notin[\tilde x_0\tilde x_n]_L$.
Further if the geodesic $[\tilde x_0\tilde x_n]_L$ contains one of the points
$\tilde x_1,\dots,\tilde x_{n-1}$, then it coincides with the polygonal line $\tilde x_0\dots\tilde x_n$.
(In particular, we have equality in \ref{eq:thm:devel*} for each~$i$.)
In this case, the sum in the left-hand side of \ref{eq:thm:detr*} must be at most $\pi$; otherwise $[\tilde x_0\tilde x_n]_L$ is not minimizing.
Therefore \ref{eq:thm:devel**} follows.
In the remaining case $[\tilde x_0\tilde x_n]_L$ meets the boundary of $L$ only at its ends.
In this case, $\dist{\tilde x_0}{\tilde x_n}{L}\le \pi$; otherwise $[\tilde x_0\tilde x_n]_L$ is not minimizing.
Whence \ref{eq:thm:devel**} follows.

Inequalities \ref{eq:thm:devel*} and \ref{eq:thm:devel**} imply that the polygon $[\tilde p'\tilde x_0\tilde x_1\dots \tilde x_n]$ is convex.

Let us take finer and finer partitions and pass to the limit of the polygon $\tilde p'\tilde x_0\tilde x_1\dots \tilde x_n$ as $p'\to p$.
We obtain a convex curvilinear triangle formed by a curve $\tilde \gamma\:[0,T]\to\Lob2{\kappa}$ --- the limit of polygonal line $\tilde x_0\tilde x_1\dots \tilde x_n$ 
and two geodesics $[\tilde p'\,\tilde \gamma(0)]$,
$[\tilde p'\,\tilde \gamma(T)]$.
Since $[\tilde p'\tilde x_0\tilde x_1\dots \tilde x_n]$ is convex,
 the natural parametrization of $\tilde x_0\tilde x_1\dots \tilde x_n$ 
converges to the natural parametrization of $\tilde \gamma$. 
Thus $\tilde \gamma$ is the $\kappa$-development of $\gamma$ with respect to $p$.
This proves the only-if part of \eqref{thm:devel}.


\parit{If part.}  
Assuming convexity of the development, we will prove the point-on-side comparison (\ref{point-on-side}). 
We can assume that $p\notin [x y]$; otherwise the statement is trivial.

Set $T=\dist{x}{y}{}$ and $\gamma(t)=\geod_{[x y]}(t)$; note that $\gamma$ is a geodesic in $\oBall(p,\varpi\kappa)\setminus \{p\}$.
Let $\tilde \gamma\:[0,T]\to\Lob2\kappa$ be the $\kappa$-development with base $\tilde p$ of $\gamma$ with respect to $p$.
Take a partition $0=t_0<t_1<\dots<t_n=T$, and set 
\[\tilde y_i=\tilde \gamma(t_i)\quad \text{and}\quad \tau_i=\dist{\tilde y_0}{\tilde y_1}{}+\dist{\tilde y_1}{\tilde y_2}{}+\dots+\dist{\tilde y_{i-1}}{\tilde y_i}{}.\] 
Since $\tilde \gamma$ is convex, for a fine partition we have that the polygonal line $\tilde y_0\tilde y_1\dots\tilde y_n$ is also convex.
Applying Alexandrov's lemma (\ref{lem:alex}) inductively to pairs of model triangles 
\[\modtrig\kappa
\{\tau_{i-1},\dist{\tilde p}{\tilde y_0}{},\dist{\tilde p}{\tilde y_{i-1}}{}\}
\]
and 
\[\modtrig\kappa
\{\dist{\tilde y_{i-1}}{\tilde y_i}{},\dist{\tilde p}{\tilde y_{i-1}}{},\dist{\tilde p}{\tilde y_{i}}{}\}\]
we obtain that the sequence 
$\tilde \mangle\mc\kappa\{\dist{\tilde p}{\tilde y_{i}}{};\dist{\tilde p}{\tilde y_0}{},\tau_i\}$ is non increasing.

For finer and finer partitions we have 
\[\max\nolimits_i\{|\tau_i-t_i|\}\to0.\] 
Thus, the point-on-side comparison (\ref{point-on-side}) follows.
\qeds

Note that in the proof of if part we could use a slightly weaker version of the hinge comparison (\ref{angle}).
Namely, we proved the following lemma, which will be needed later in the proof of the globalization theorem (\ref{thm:glob}).

\begin{thm}{Short hinge lemma}\label{lem:devel-glob}
Let $\spc{L}$ be a complete G-delta geodesic space such that for any hinge $\hinge x p y$ in $\spc{L}$ the angle $\mangle\hinge x p y$ is defined, 
and 
\[\mangle\hinge x p y+\mangle\hinge x p z\le \pi\] 
for any two adjacent hinges.

Assume that  for any hinge $\hinge x p y$ in $\spc{L}$ we have
\[\dist{p}{ x}{}+\dist{x}{ y}{}
<\varpi\kappa
\quad\Rightarrow\quad 
\mangle\hinge x p y
\ge\angk\kappa x p y.\]
Then $\spc{L}$ is $\Alex\kappa$.
\end{thm}


\section{Local definitions and globalization}\label{sec:loc}

In this section we discuss locally $\Alex{\kappa}$ spaces.
In particular, we prove the {}\emph{globalization theorem}: equivalence of local and global definitions for complete length spaces.

The following theorem summarizes equivalent definitions of locally $\Alex{\kappa}$ spaces

\begin{thm}{Theorem}\label{thm:=def-loc}
Let $\spc{X}$ be a complete length space and $p\in \spc{X}$.
Then the following conditions are equivalent:
\begin{subthmN}
{curv>=k}(local $\Alex\kappa$ comparison) there is $R_{\hbox{\rm\scriptsize\ref{SHORT.curv>=k}}}>0$ such that the comparison 
\[\angk\kappa q{x^1}{x^2}
+\angk\kappa q{x^2}{x^3}
+\angk\kappa q{x^3}{x^1}
\le 2\cdot\pi\]
holds for any $q,x^1,x^2,x^3\in \oBall(p,R_{\hbox{\rm\scriptsize\ref{SHORT.curv>=k}}})$.
\end{subthmN}

\begin{subthmN}{def:kirszbraun-loc}(local Kirszbraun property) 
there is $R_{\hbox{\rm\scriptsize\ref{SHORT.def:kirszbraun-loc}}}>0$ 
such that for any 3-point subset $F_3$ and any 4-point subset $F_4\supset F_3$ in $\oBall(p,R_{\hbox{\rm\scriptsize\ref{SHORT.def:kirszbraun-loc}}})$, any short map $f\:F_3\to\Lob2\kappa$ can be extended to a short map $\bar f\:F_4\to\Lob2\kappa$ (so $f=\bar f|_{F_3}$).
\end{subthmN}

\begin{subthmN}{conc-loc} (local function comparison) there is 
$R_{\hbox{\rm\scriptsize\ref{SHORT.conc-loc}}}>0$ such that 
$\oBall(p,R_{\hbox{\rm\scriptsize\ref{SHORT.conc-loc}}})$ 
is G-delta geodesic and for any 
$q\in \oBall(p,R_{\hbox{\rm\scriptsize\ref{SHORT.conc-loc}}})$, 
the function $f=\md\kappa\circ\distfun{q}{}{}$ satisfies 
$f''\le 1-\kappa\cdot  f$ in
$\oBall(p,R_{\hbox{\rm\scriptsize\ref{SHORT.conc-loc}}})$.
\end{subthmN}

\begin{subthmN}{2-sum-loc} (local adjacent angle comparison) there is
$R_{\hbox{\rm\scriptsize\ref{SHORT.2-sum-loc}}}>0$ 
such that $\oBall(p,R_{\hbox{\rm\scriptsize\ref{SHORT.2-sum-loc}}})$ 
is G-delta geodesic, and if $q$ and a geodesic $[x y]$ lie in
$\oBall(p,R_{\hbox{\rm\scriptsize\ref{SHORT.2-sum-loc}}})$ 
and $z\in \mathopen{]}x y\mathclose{[}$, then
\[\angk\kappa z q x
+\angk\kappa z q y
\le \pi.\]
\end{subthmN}

\begin{subthmN}{monoton-loc} (local point-on-side comparison) 
there is $R_{\hbox{\rm\scriptsize\ref{SHORT.monoton-loc}}}>0$ 
such that $\oBall(p,R_{\hbox{\rm\scriptsize\ref{SHORT.monoton-loc}}})$ is G-delta geodesic and if $q$ and a geodesic $[x y]$ lie in $\oBall(p,R_{\hbox{\rm\scriptsize\ref{SHORT.monoton-loc}}})$ 
and $z\in \mathopen{]}x y\mathclose{[}$, then
\[\angk\kappa x q y
\le\angk\kappa x q z;\]
or, equivalently, 
\[\dist{\tilde p}{\tilde z}{}\le \dist{p}{z}{},\]
where $\trig{\tilde p}{\tilde x}{\tilde y}=\modtrig\kappa(p x y)$, $\tilde z\in\mathopen{]} \tilde x\tilde y\mathclose{[}$, $\dist{\tilde x}{\tilde z}{}=\dist{x}{z}{}$.
\end{subthmN}

\begin{subthmN}{angle-loc} (local hinge comparison) 
there is $R_{\hbox{\rm\scriptsize\ref{SHORT.angle-loc}}}>0$ such that $\oBall(p,R_{\hbox{\rm\scriptsize\ref{SHORT.angle-loc}}})$ is G-delta geodesic and if $x\in \oBall(p,R_{\hbox{\rm\scriptsize\ref{SHORT.angle-loc}}})$, then 
for any hinge $\hinge x q y$, the angle 
$\mangle\hinge x q y$ is defined, and
\[\mangle\hinge x q y+ \mangle\hinge x q z\le\pi\]
for any two adjacent hinges.
Moreover, if a hinge $\hinge x q y$ lies in $\oBall(p,R_{\hbox{\rm\scriptsize\ref{SHORT.angle-loc}}})$, then 
\[\mangle\hinge x q y
\ge\angk\kappa x q y,\]
or, equivalently,
\[\side\kappa \hinge x q y\ge\dist{q}{ y}{}.\]
\end{subthmN}

\begin{subthmN}{def:devel-alex-loc}(local development comparison) 
there is $R_{\hbox{\rm\scriptsize\ref{SHORT.def:devel-alex-loc}}}>0$ 
such that $\oBall(p,R_{\hbox{\rm\scriptsize\ref{SHORT.def:devel-alex-loc}}})$ 
is G-delta geodesic, and if a geodesic $\gamma$ lies in $\oBall(p,R_{\hbox{\rm\scriptsize\ref{SHORT.def:devel-alex-loc}}})$ and $q\in 
\oBall(p,R_{\hbox{\rm\scriptsize\ref{SHORT.def:devel-alex-loc}}})
\setminus \gamma$, then the $\kappa$-development $\tilde \gamma$ with respect to $q$ is convex.
\end{subthmN}
Moreover, for each pair $i,j\in \{1,2,\dots,7\}$ we can assume that 
\[R_i>\tfrac{1}{9}\cdot R_j.\]
\end{thm}

The proofs of each of these equivalences repeat the proofs of the corresponding global equivalences in localized form; see the proofs of Theorems \ref{thm:defs_of_alex}, \ref{thm:conc}, \ref{def:devel}, \ref{thm:kirsz-def}.

\begin{thm}{Globalization theorem}\label{thm:glob} 
Any complete length locally $\Alex\kappa$ space is $\Alex\kappa$.
\end{thm}

{\sloppy

In the two-dimensional case this theorem was proved by Paolo Pizzetti \cite{pizzetti};
later it was reproved independently by Alexandr Alexandrov \cite{alexandrov:devel}. 
Victor Toponogov \cite{toponogov-globalization+splitting} proved it for Riemannian manifolds of all dimensions.
In the above generality, the theorem first appears in the paper of Michael Gromov, Yuriy Burago, and Grigory Perelman \cite{burago-gromov-perelman}; 
simplifications and modifications were given by Conrad Plaut \cite{plaut:dimension}, Katsuhiro Shiohama \cite{shiohama}, and in the book of Dmitry Burago, Yuriy Burago, and Sergei Ivanov \cite{burago-burago-ivanov}.
A generalization for non-complete but geodesic spaces was obtained by the third author \cite{petrunin:globalization}; namely it solves the following exercise:

}

\begin{thm}{Advanced exercise}\label{ex:noncomplete-globalization}
Any locally $\Alex\kappa$ geodesic space is $\Alex\kappa$.
\end{thm}

The following corollary of the globalization  theorem says that the expression ``space with curvature $\ge \kappa$'' makes sense.

\begin{thm}{Corollary}\label{cor:CAT>k-sence}
Let $\spc{L}$ be a complete length space.
Then $\spc{L}$ is $\Alex\Kappa$ if and only if $\spc{L}$ is $\Alex\kappa$ for any $\kappa<\Kappa$.
\end{thm}

\parit{Proof.}
Note that if $\Kappa\le 0$, this statement follows directly from the  definition of an Alexandrov space (\ref{df:cbb1+3}) and monotonicity of the function $\kappa\mapsto\angk\kappa x y z$ (\ref{k-decrease}).

The if part also follows directly from the definition.

For $\Kappa>0$, the angle $\angk{\Kappa} x y z$ might be undefined while $\angk{\kappa} x y z$ is defined.
However, $\angk{\Kappa} x y z$ is defined if $x$, $y$, and $z$ are sufficiently close to each other.
Thus, if $\Kappa>\kappa$, then any $\Alex\Kappa$ space is locally $\Alex\kappa$.
It remains to apply the  globalization theorem.
\qeds

\begin{thm}{Corollary}\label{cor:submetry-cbb}
Let $\spc{L}$ be a complete length $\Alex\kappa$ space.
Assume that a space  $\spc{M}$ is the target space of a submetry from $\spc{L}$.
Then $\spc{M}$ is a complete length space $\Alex\kappa$ space.

In particular, if $G\acts \spc{L}$ is an isometric group action with closed orbits, then the quotient space $\spc{L}/G$ is a complete length $\Alex\kappa$ space.
\end{thm}

\parit{Proof.}
This follows from the globalization theorem and Theorem \ref{thm:submetry-CBB}.
\qeds

Our proof of the globalization theorem (\ref{thm:glob}) is based on presentations in \cite{plaut:dimension} and \cite{burago-burago-ivanov};
this proof was rediscovered independently by Urs Lang and Viktor Schroeder \cite{lang-schroeder:globalization}.
We will need the short hinge lemma (\ref{lem:devel-glob}) together with the following two lemmas.
The following lemma says that if comparison holds for all small hinges, then it holds for slightly bigger hinges near the given point.

\begin{thm}{Key lemma}\label{key-lem:globalization} 
Let $\kappa\in \RR$, 
$0<\ell\le\varpi\kappa$, 
$\spc{X}$ be a complete geodesic space 
and $p\in \spc{X}$ be a point 
such that $\oBall(p,2\cdot\ell)$ is locally $\Alex\kappa$. 

Assume that for any point 
$q\in \oBall(p,{\ell})$ the comparison
\[\mangle\hinge x y q
\ge\angk\kappa x y q\]
holds for any hinge $\hinge x y q$ with 
$\dist{x}{y}{}+\dist{x}{q}{}
<
\frac{2}{3}\cdot\ell$.
Then the comparison
\[\mangle\hinge x p q
\ge\angk\kappa x p q\] 
holds for any hinge $\hinge x p q$ with $\dist{x}{ p}{}+\dist{x}{q}{}<\ell$.
\end{thm}

\parit{Proof.} 
It is sufficient to prove the inequality
\[\side\kappa \hinge x p q
\ge\dist{p}{q}{}\eqlbl{eq:thm:=def-loc*}\] 
for any hinge $\hinge x p q$ with $\dist{x}{p}{}+\dist{x}{q}{}<\ell$.

Fix $q$.
Consider a hinge $\hinge x p q$ such that 
\[\tfrac{2}{3}\cdot\ell \le\dist{p}{x}{}\z+\dist{x}{q}{}< \ell.\]
First we  construct a new smaller hinge $\hinge{x'}p q$ with
\[
\dist{p}{x}{}+\dist{x}{q}{}\ge\dist{p}{x'}{}+\dist{x'}{q}{},
\eqlbl{eq:thm:=def-loc-fourstar}\]
such that 
\[\side\kappa \hinge x p q
\ge\side\kappa \hinge{x'}p q.
\eqlbl{eq:thm:=def-loc-fivestar}\]

Assume $\dist{x}{q}{}\ge\dist{x}{p}{}$; otherwise switch the roles of $p$ and $q$ in the following construction.
Take $x'\in [x q]$ such that 
\[\dist{p}{x}{}+3\cdot\dist[{{}}]{x}{x'}{}
=\tfrac{2}{3}\cdot\ell. \eqlbl{3|xx'|}\]
Choose a geodesic $[x' p]$ and consider the  hinge $\hinge{x'}p q$ formed by $[x'p]$ and $[x' q]\subset [x q]$.
(In fact by \ref{cor:unique-geod-cbb} the condition $[x' q]\subset [x q]$ always holds.)
Then \ref{eq:thm:=def-loc-fourstar} follows from the triangle inequality.

Further, note that we have $x,x'\in \oBall(p,\ell)\cap \oBall(q,\ell)$ and moreover
\begin{align*}
\dist{p}{x}{}\z+\dist{x}{x'}{}&<\tfrac{2}{3}\cdot\ell,
&
\dist{p}{x'}{}\z+\dist{x'}{x}{}&<\tfrac{2}{3}\cdot\ell.
\end{align*}
In particular, 
\[\mangle\hinge x p{x'}
\ge\angk\kappa x p{x'}
\quad \text{and}\quad 
\mangle\hinge {x'}p x
\ge\angk\kappa {x'}p x.
\eqlbl{eq:thm:=def-loc-threestar}\]

{

\begin{wrapfigure}{r}{30 mm}
\vskip-0mm
\centering
\includegraphics{mppics/pic-820}
\vskip-4mm
\end{wrapfigure}

Now, let 
$\trig{\tilde x}{\tilde x'}{\tilde p}=\modtrig\kappa(x x' p)$.
Take $\tilde  q$ on the extension of $[\tilde  x\tilde  x']$ beyond $x'$ such that $\dist{\tilde x}{\tilde q}{}=\dist{x}{q}{}$ (and therefore $\dist{\tilde x'}{\tilde q}{}=\dist{x'}{q}{}$).
From \ref{eq:thm:=def-loc-threestar},
\[\mangle\hinge x p q
=\mangle\hinge  x p{x'}\ge\angk\kappa x p{x'}\quad \Rightarrow\quad 
\side\kappa \hinge x q p\ge\dist{\tilde p}{\tilde q}{}.\]
Hence
\begin{align*}
\mangle\hinge{\tilde x'}{\tilde p}{\tilde q}&= 
\pi
-\angk\kappa{x'}p x
\ge
\\
&\ge
\pi-\mangle\hinge{x'}p x
=
\\
&=
\mangle\hinge{x'}p q,
\end{align*}
and \ref{eq:thm:=def-loc-fivestar} follows.

}

Let us continue the proof.
Set $x_0=x$.
Let us apply inductively the above construction to get a sequence of hinges  $\hinge{x_n}p q$ with $x_{n+1}=x_n'$.
From \ref{eq:thm:=def-loc-fivestar}, we have that the sequence  $s_n\z=\side\kappa \hinge{x_n}p q$ is nonincreasing.
\begin{figure}[ht!]
\centering
\includegraphics{mppics/pic-825}
\end{figure}

The sequence might terminate at $x_n$ only if $\dist{p}{x_n}{}+\dist{x_n}{q}{}\z< \tfrac{2}{3}\cdot\ell $.
In this case, by the assumptions of the lemma, $\side\kappa \hinge{x_n}p q\ge\dist{p}{q}{}$.
Since the sequence $s_n$ is nonincreasing, inequality \ref{eq:thm:=def-loc*} follows.

Otherwise, the sequence $r_n=\dist{p}{x_n}{}+\dist{x_n}{q}{}$ is nonincreasing and $r_n\ge\tfrac{2}{3}\cdot\ell$ for all $n$.
Note that by construction, the distances
$\dist{x_n}{x_{n+1}}{}$, $\dist{x_n}{p}{}$, and $\dist{x_n}{q}{}$ are bounded away from zero for all large $n$.
Indeed, since on each step we move $x_n$ toward to the point $p$ or $q$ that is further away, the distances $\dist{x_n}{p}{}$ and $\dist{x_n}{q}{}$ become about the same.
Namely, by \ref{3|xx'|}, we have that $\dist{p}{x_n}{}-\dist{x_n}{q}{}\le \tfrac29\cdot\ell$ for all large $n$.
Since $\dist{p}{x_n}{}+\dist{x_n}{q}{}\ge \tfrac23\cdot\ell$, we have $\dist{x_n}{p}{}\ge \tfrac\ell{100}$ and $\dist{x_n}{q}{}\ge \tfrac\ell{100}$.
Further, since $r_n\ge\tfrac{2}{3}\cdot\ell$, \ref{3|xx'|} implies that $\dist{x_n}{x_{n+1}}{}>\tfrac\ell{100}$.

Since the sequence $r_n$ is nonincreasing, it converges.
In particular, $r_n-r_{n+1}\to 0$ as $n\to\infty$.
It follows that $\angk\kappa{x_n}{p_n}{x_{n+1}}\to \pi$,
where $p_n=p$ if $x_{n+1}\in [x_nq]$, and otherwise $p_n=q$.
Since $\mangle\hinge{x_n}{p_n}{x_{n+1}}\ge\angk\kappa{x_n}{p_n}{x_{n+1}}$, we have
$\mangle\hinge{x_n}{p_n}{x_{n+1}}\to \pi$  as $n\to\infty$.

It follows that
\[r_n-s_n=\dist{p}{x_n}{}+\dist{x_n}{q}{}-\side\kappa \hinge{x_n}p q\to 0.\] 
(Here we used that $\ell\le\varpi\kappa$.) 
Together with the triangle inequality
\[
\dist{p}{x_n}{}+\dist{x_n}{q}{}\ge\dist{p}{q}{}
\]
this yields
\[\lim_{n\to\infty}\side\kappa \hinge{x_n}p q\ge \dist{p}{q}{}.\]
Applying monotonicity of the sequence  $s_n=\side\kappa \hinge{x_n}p q$, we obtain \ref{eq:thm:=def-loc*}.
\qeds

The final part of the proof above resembles the \emph{cat's cradle construction} introduced by the first author and Richard Bishop \cite{alexander-bishop:h-c}.

The following lemma works in all complete spaces; it will be used as a substitute for the  existence of a minimum point of a continuous function on a compact space.

\begin{thm}{Lemma on almost minimum}\label{lem:alm-min}
Let $\spc{X}$ be a complete metric space.
Suppose $r\:\spc{X}\to \RR$ is a function, $p\in \spc{X}$,   and $\eps>0$.
Assume that the function $r$ is strictly positive in $\cBall[p,\tfrac{1}{\eps^2}\cdot r(p)]$ and
$\varliminf_{n}r(x_n)>0$ for any convergent sequence 
$x_n\to x\in \cBall[p,\tfrac{1}{\eps^2}\cdot r(p)]$. 

Then there is a point $p^*\in \cBall[p,\tfrac{1}{\eps^2}\cdot r(p)]$ such that 

\begin{subthm}{}$r(p^*)\le r(p)$ and
\end{subthm}

\begin{subthm}{}$r(x)> (1-\eps)\cdot r(p^*)$ 
for any $x\in \cBall[p^*,\tfrac{1}{\eps}\cdot r(p^*)]$.
\end{subthm}
\end{thm}

\parit{Proof.} 
Assume the statement is wrong. 
Then for any $x\in \oBall(p,\tfrac{1}{\eps^2}\cdot r(p))$ with $r(x)\le r(p)$, there is a point $x'\in \spc{X}$ such that 
\[\dist{x}{x'}{}<\tfrac{1}\eps\cdot r(x)
\quad \text{and}\quad 
r(x')\le (1-\eps)\cdot r(x).\]
Take $x_0=p$ and consider a sequence of points $x_n$ such that $x_{n+1}\z=x_n'$.
Clearly 
\[\dist{x_{n+1}}{x_n}{}
\le
\tfrac{r(p)}{\eps}\cdot(1-\eps)^n
\quad \hbox{and}\quad 
r(x_n)\le r(p)\cdot(1-\eps)^n.\] 
In particular, $\dist{p}{x_n}{}<\tfrac{1}{\eps^2}\cdot r(p)$.
Therefore the sequence $x_n$ is Cauchy,
$x_n\to x\in \cBall[p,\tfrac{1}{\eps^2}\cdot r(p)]$
and
$\lim_{n}r(x_n)=0$, a contradiction.
\qeds

\parit{Proof of the globalization theorem (\ref{thm:glob}).} 
Exactly the same argument as in the proof of Theorem~\ref{thm:almost.geod} 
shows that $\spc{L}$ is G-delta geodesic.
By Theorem~\ref{angle-loc}, 
for any hinge $\hinge x p y$ in $\spc{L}$ the angle $\mangle\hinge x p y$ is defined 
and 
\[\mangle\hinge x p y+\mangle\hinge x p z\le \pi\] 
for any two adjacent hinges.

Let us denote by $\ComRad(p,\spc{L})$ 
(which stands for \index{comparison radius}\emph{comparison radius} of $\spc{L}$ at $p$) 
the maximal value (possibly $\infty$) such that the comparison 
\[\mangle\hinge x p y
\ge\angk\kappa x p y\]
holds for any hinge $\hinge x p y$ with $\dist{p}{x}{}+\dist{x}{y}{}< \ComRad(p,\spc{L})$.

As follows from \ref{conc-loc}, $\ComRad(p,\spc{L})>0$ for any $p\in\spc{L}$ and 
$$\liminf_{n\to\infty}\ComRad(p_n,\spc{L})>0$$ 
for any converging sequence of points $p_n\to p$.
That makes it possible to apply the lemma on almost minimum (\ref{lem:alm-min}) to the function $p\mapsto \ComRad(p,\spc{L})$.

According to the short hinge lemma (\ref{lem:devel-glob}), it is sufficient  to show that 
\[s_0=\inf_{p\in\spc{L}}\ComRad(p,\spc{L})
\ge \varpi\kappa\quad 
\text{for any}\quad 
p\in \spc{L}.
\eqlbl{eq:thm:=def-loc-star-star}\]
We argue by contradiction, assuming that  \ref{eq:thm:=def-loc-star-star} does not hold.

The rest of the proof is easier for geodesic spaces 
and easier still for compact spaces.
Thus we give three different arguments  for each of these cases.

\parit{Compact case.}
Assume $\spc{L}$ is compact.

By Theorem~\ref{conc-loc},  $s_0>0$.
Take a point $p^*\in \spc{L}$ such that $r^*\z=\ComRad(p^*,\spc{L})$ is sufficiently close to $s_0$
($p^*$ such that  $s_0\z\le r^*\z<\min\{\varpi\kappa,\tfrac32\cdot s_0\}$ will do).
Then the key lemma (\ref{key-lem:globalization}) applied to $p^*$ and $\ell$ slightly bigger than $r^*$ (say, such that $r^*<\ell<\min\{\varpi\kappa,\tfrac32\cdot s_0\}$) implies that
\[\mangle\hinge x{p^*}q
\ge\angk\kappa x{p^*}q\]
for any hinge $\hinge x{p^*}q$ such that $\dist{p^*}{x}{}+\dist{x}{q}{}<\ell$.
Thus $r^*\ge\ell$, a contradiction.

\parit{Geodesic case.}
Assume $\spc{L}$ is geodesic.

Fix a small $\eps>0$ ($\eps=0.0001$ will do). 
Apply the lemma on almost minimum (\ref{lem:alm-min}) to find a point $p^*\in \spc{L}$ such that 
\[r^*\z=\ComRad(p^*,\spc{L})\z<\varpi\kappa\] 
and 
\[\ComRad(q,\spc{L})\z> (1-\eps)\cdot r^*\eqlbl{comrad(q)}\] 
for any $q\in\cBall[p^*,\tfrac{1}{\eps}\cdot r^*]$. 

Applying the key lemma (\ref{key-lem:globalization}) for $p^*$ and  $\ell$ slightly bigger than $r^*$ leads to a contradiction.

\parit{General case.}
Let us construct $p^*\in\spc{L}$ as in the previous case.
Since $\spc{L}$ is not geodesic, we cannot apply the key lemma directly.
Instead, let us pass to the ultrapower $\spc{L}^\o$, which
 is a geodesic space (see \ref{obs:ultralimit-is-geodesic}).

In Theorem~\ref{thm:=def-loc},
inequality \ref{comrad(q)} implies that 
condition  \ref{curv>=k} holds for some fixed $R_{\hbox{\rm\scriptsize\ref{SHORT.curv>=k}}}=\tfrac{r^*}{100}
>0$ at any point $q\in\cBall[p^*,\tfrac1{2\cdot \eps}\cdot r^*]\subset\spc{L}$.
Therefore a similar statement is true in the ultrapower $\spc{L}^\o$;
that is,
for any point
$q_\o\in\cBall[p^*,\tfrac1{2\cdot \eps}\cdot r^*]\subset\spc{L}^\o$, 
condition~\ref{curv>=k} holds for, say, $R_{\hbox{\rm\scriptsize\ref{SHORT.curv>=k}}}=\tfrac{r^*}{101}$.

Note that $r^*\ge\ComRad(p^*,\spc{L}^\o)$.
Therefore we can apply the lemma on almost minimum  
at the point $p^*$ to the function $x\z\mapsto\ComRad(x,\spc{L}^\o)$
and $\eps'=\sqrt\eps=0.01$.

For the resulting point $p^{{*}{*}}\in \spc{L}^\o$, we have $r^{{*}{*}}\z=\ComRad(p^{{*}{*}},\spc{L})\z<\varpi\kappa$, 
and 
$\ComRad(q_\o,\spc{L}^\o)\z> (1-\eps')\cdot r^{{*}{*}}$ for any $q_\o\in\cBall[p^{{*}{*}},\tfrac{1}{\eps'}\cdot r^{{*}{*}}]$. 
Thus applying the key lemma (\ref{key-lem:globalization}) for $p^{{*}{*}}$ and for $\ell$ slightly bigger than $r^{{*}{*}}$ leads to a contradiction.
\qeds


\section{Properties of geodesics and angles}\label{sec:prop.geod}

\parbf{Remark.} 
All proofs in this section can be easily modified to use only the local definition of $\Alex{}$ spaces without use of the globalization theorem (\ref{thm:glob}).
 
\begin{thm}{Geodesics do not split}\label{thm:g-split}
In a $\Alex{}$ space, geodesics do not bifurcate.

More precisely, let $\spc{L}$ be a $\Alex{}$ space and $[p x]$, $[p y]$ be two geodesics. Then:
\begin{subthm}{nonsplit} If  there is $\eps>0$ such that $\geod_{[p x]}(t)=\geod_{[p y]}(t)$ 
for all $t\in[0,\eps)$, 
then $[p x]\subset [p y]$ or $[p y]\subset [p x]$.
\end{subthm}

\begin{subthm}{angle=0}
If $\mangle\hinge p x y=0$, then $[p x]\subset [p y]$ or $[p y]\subset [p x]$.
\end{subthm}
\end{thm}

\begin{thm}{Corollary}\label{cor:unique-geod-cbb}
Let $\spc{L}$ be a $\Alex{}$ space.
Then the restriction of any geodesic in $\spc{L}$ to a proper segment is the unique minimal geodesic joining its endpoints.
\end{thm}

In case $\kappa\le 0$, the proof is easier, since the model triangles are always defined.
To deal with $\kappa>0$ we have to argue locally.

\begin{wrapfigure}{r}{25 mm}
\vskip-0mm
\centering
\includegraphics{mppics/pic-830}
\end{wrapfigure}

\parit{Proof of \ref{thm:g-split}; \ref{SHORT.nonsplit}.}
Let $t_{\max}$ be the maximal value 
such that $\geod_{[p x]}(t)=\geod_{[p y]}(t)$ for all $t\in [0,t_{\max})$.
Since geodesics are continuous, $\geod_{[p x]}(t_{\max})\z=\geod_{[p y]}(t_{\max})$.
Let
\[q=\geod_{[p x]}(t_{\max})\z=\geod_{[p y]}(t_{\max}).\]
We must show that $t_{\max}=\min\{\dist{p}{x}{},\dist{p}{y}{}\}$.

If that is not true, choose a sufficiently small $\eps>0$ such that the  points
\[x_\eps=\geod_{[p x]}(t_{\max}+\eps)\quad 
\text{and}\quad  
  y_\eps=\geod_{[p y]}(t_{\max}+\eps)\] 
are distinct.
Let
\[z_\eps=\geod_{[p x]}(t_{\max}-\eps)
=\geod_{[p y]}(t_{\max}-\eps).\]
Clearly, $\angk\kappa q{z_\eps}{x_\eps}=\angk\kappa q{z_\eps}{y_\eps}=\pi$.
Thus from the $\Alex\kappa$ comparison (\ref{df:cbb1+3}), $\angk\kappa q{x_\eps}{y_\eps}=0$ and thus $x_\eps=y_\eps$, a contradiction.

\parit{\ref{SHORT.angle=0}.}
From hinge comparison \ref{angle}, 
\[\mangle\hinge p x y=0\quad \Rightarrow\quad \angkk\kappa p{\geod_{[p x]}(t)}{\geod_{[p y]}(t)}=0\] 
and thus ${\geod_{[p x]}(t)}={\geod_{[p y]}(t)}$ for all small $t$. 
Therefore we can apply \ref{SHORT.nonsplit}.
\qeds

\begin{thm}{Adjacent angle lemma}\label{lem:sum=pi}
Let $\spc{L}$ be a $\Alex{}$ space.
Assume that two hinges $\hinge z x p$ and $\hinge z y p$ in $\spc{L}$ are adjacent.
Then 
\[\mangle\hinge z p y + \mangle\hinge z p x=\pi. \]

\end{thm}

\parit{Proof.}
From the hinge comparison (\ref{angle}) we have that both angles 
$\mangle\hinge z p y$ and $\mangle\hinge z p x$ are defined and 
\[\mangle\hinge z p y + \mangle\hinge z p x\le\pi.\]
Clearly $\mangle\hinge z x y=\pi$.
Thus the result follows from the triangle inequality for angles (\ref{claim:angle-3angle-inq}).
\qeds

\begin{thm}{Angle semicontinuity}\label{lem:ang.semicont-cbb}
Suppose  $\spc{L}_n$ is a sequence of $\Alex\kappa$ spaces and $\spc{L}_n\to \spc{L}_\o$ as $n\to\o$.
Assume that a sequence of hinges $\hinge{p_n}{x_n}{y_n}$ in $\spc{L}_n$ converges to a hinge $\hinge{p_\o}{x_\o}{y_\o}$ in  $\spc{L}_\o$.
Then 
\[\mangle\hinge{p_\o}{x_\o}{y_\o}
\le 
\lim_{n\to\o} \mangle\hinge{p_n}{x_n}{y_n}.\]

\end{thm}

\begin{wrapfigure}{r}{23 mm}
\vskip-0mm
\centering
\includegraphics{mppics/pic-835}
\end{wrapfigure}

\parit{Proof.}
From \ref{cor:monoton},
\[\mangle\hinge{p_\o}{x_\o}{y_\o}
=
\sup\set{\angk\kappa{p_\o}{\bar x_\o}{\bar y_\o}}{\bar x_\o \in \mathopen{]}p_\o x_\o],\ \bar y_\o\in \mathopen{]}p_\o x_\o]}.\]

For fixed $\bar x_\o \in \mathopen{]}p_\o x_\o]$ 
and $\bar y_\o\in \mathopen{]}p_\o x_\o]$,
choose $\bar x_n\in \mathopen{]} p x_n ]$ and $\bar y_n\in \mathopen{]} p y_n ]$ so that $\bar x_n\to \bar x_\o$ 
and $\bar y_n\to \bar y_\o$ as $n\to\o$.
Clearly 
\[\angk\kappa{p_n}{\bar x_n}{\bar y_n}
\to 
\angk\kappa{p_\o}{\bar x_\o}{\bar y_\o}\] 
as $n\to\o$.

From the hinge comparison (\ref{angle}), $\mangle\hinge{p_n}{x_n}{y_n}\ge \angk\kappa{p_n}{\bar x_n}{\bar y_n}$.
Hence the result.
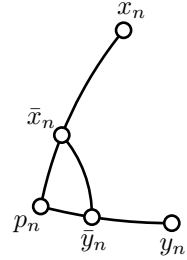
\qeds

\begin{thm}{Angle continuity}\label{cor:ang.cont-cbb}
Let $\spc{L}_n$  be a sequence of complete length $\Alex\kappa$ spaces,
and $\spc{L}_n\to \spc{L}_\o$ as $n\to\o$.
Assume that sequences of points $p_n, x_n, y_n$ in $\spc{L}_n$ 
converge to points $p_\o, x_\o, y_\o$ in  $\spc{L}_\o$ as $n\to\o$,
and the following two conditions hold:
\begin{subthm}{}
$p_\o\in \Str(x_\o)$,
\end{subthm}
\begin{subthm}{}
$p_\o\in \Str(y_\o)$ or $y_\o\in \Str(p_\o)$.
\end{subthm}

Then 
\[\mangle\hinge{p_\o}{x_\o}{y_\o}
=
\lim_{n\to\o} \mangle\hinge{p_n}{x_n}{y_n}.\]

\end{thm}

\parit{Proof.}
By Corollary~\ref{cor:CAT>k-sence},
we may assume that $\kappa\le 0$.

By Plaut's theorem (\ref{thm:almost.geod}),
the hinge 
$\hinge{p_\o}{x_\o}{y_\o}$
is uniquely defined.
Therefore the hinges 
$\hinge{p_n}{x_n}{y_n}$
converge to  
$\hinge{p_\o}{x_\o}{y_\o}$
as $n\to\o$.
Hence by the angle semicontinuity (\ref{lem:ang.semicont-cbb}), 
we have
\[
\mangle\hinge{p_\o}{x_\o}{y_\o}
\le
\lim_{n\to\o} \mangle\hinge{p_n}{x_n}{y_n}.
\]

It remains to show that 
\[
\mangle\hinge{p_\o}{x_\o}{y_\o}
\ge
\lim_{n\to\o} \mangle\hinge{p_n}{x_n}{y_n}.
\eqlbl{eq:ang-semicon->}
\]

Fix $\eps>0$.
Since $p_\o\in\Str(x_\o)$,
 there is a point $q_\o\in\spc{L}_\o$
such that 
\[\angk{\kappa}{p_\o}{x_\o}{q_\o}>\pi-\eps.\]
The hinge comparison  (\ref{angle}) implies that
\[\mangle\hinge{p_\o}{x_\o}{q_\o}>\pi-\eps.
\eqlbl{eq:>pi-eps}\]
By the triangle inequality for angles (\ref{claim:angle-3angle-inq}),
\[
\begin{aligned}
\mangle\hinge{p_\o}{x_\o}{y_\o}
&\ge \mangle\hinge{p_\o}{x_\o}{q_\o}-
\mangle\hinge{p_\o}{y_\o}{q_\o}>
\\
&>\pi-\eps-
\mangle\hinge{p_\o}{y_\o}{q_\o}.
\end{aligned}
\eqlbl{eq:ang-trig}
\]

Note that we can assume in addition that $q_\o\in\Str(p_\o)$.
Choose $q_n\in\spc{L}_n$
such that $q_n\to q_\o$ as $n\to\o$.
Note that by angle semicontinuity
 we again have
\[
\begin{aligned}
\mangle\hinge{p_\o}{x_\o}{q_\o}
&\le
\lim_{n\to\o} \mangle\hinge{p_n}{x_n}{q_n},
\\
\mangle\hinge{p_\o}{y_\o}{q_\o}
&\le
\lim_{n\to\o} \mangle\hinge{p_n}{y_n}{q_n}.
\end{aligned}
\eqlbl{eq:semicont}
\]

By the $\Alex\kappa$ comparison (\ref{df:cbb1+3}) and \ref{cor:monoton:sup}, 
\[\mangle\hinge{p_n}{x_n}{y_n}
+\mangle\hinge{p_n}{y_n}{q_n}
+\mangle\hinge{p_n}{x_n}{q_n}
\le 2\cdot\pi\]
for all $n$.
Together with \ref{eq:semicont}, \ref{eq:>pi-eps} and \ref{eq:ang-trig}, 
this implies
\begin{align*}
\lim_{n\to\o} \mangle\hinge{p_n}{x_n}{y_n}
&\le
2\cdot\pi-\lim_{n\to\o} \mangle\hinge{p_n}{x_n}{q_n}-\lim_{n\to\o} \mangle\hinge{p_n}{y_n}{q_n}
\le 
\\
&\le 
2\cdot\pi
-\mangle\hinge{p_\o}{x_\o}{q_\o}
-\mangle\hinge{p_\o}{y_\o}{q_\o}
<
\\
&<\mangle\hinge{p_\o}{x_\o}{y_\o}+2\cdot\eps.
\end{align*}
Since $\eps>0$ is arbitrary, \ref{eq:ang-semicon->} follows.
\qeds

\begin{thm}{First variation formula}\label{1st-var+}
Let $\spc{L}$ be a complete length $\Alex{}$ space.
For any point $q$ and any geodesic $[px]$ in $\spc{L}$ with $p\ne q$, we have 
\[\dist{q}{\geod_{[p x]}(t)}{}
=
\dist{q}{p}{}-t\cdot\cos\phi+o(t),
\eqlbl{eq:1st-var+*}\]
where $\phi$ is the infimum of angles between $[px]$ and all geodesics from $p$ to $q$ in the ultrapower $\spc{L}^\o$.
\end{thm}

\parbf{Remark.}
If $\spc{L}$ is a proper space, then $\spc{L}^\o=\spc{L}$; see Section~\ref{ultralimits}.
Therefore the infimum $\phi$ is achieved on a particular geodesic from $p$ to $q$.

\medskip

As a corollary we obtain the following classical  result:

\begin{thm}{Strong angle lemma}\label{lem:strong-angle}
Let $\spc{L}$ be a complete length $\Alex{}$ space and $p\ne q\in \spc{L}$ be such that there is unique geodesic from $p$ to $q$ in the ultrapower $\spc{L}^\o$.
Then for any hinge  $\hinge  p q x$ we have
\[\mangle\hinge p q x
=
\lim_{
\substack{
\bar x\to p
\\
\bar x\in\,\mathopen{]}px]}}\angk\kappa p q{\bar x}\eqlbl{eq:1st-var+***}\]
for any $\kappa\in\RR$ such that $\dist{p}{q}{}<\varpi\kappa$.

In particular, \ref{eq:1st-var+***} holds if $p\in \Str(q)$ as well as if $q\in \Str(p)$. 
\end{thm}

\parbf{Remark.}
\begin{itemize}
\item The above lemma is essentially due to Alexandrov.
The right-hand side in \ref{eq:1st-var+***} is called the \index{angle!strong angle}\index{strong angle}\emph{strong angle} of the  hinge $\hinge p q x$. 
Note that in a general metric space the angle and the strong angle of the same hinge might differ.

\item As follows from Corollary~\ref{cor:two-geodesics-in-ultrapower}, 
if there is a unique geodesic $[p q]$ in the ultrapower $\spc{L}^\o$, then $[p q]$ lies in $\spc{L}$.
\end{itemize}

\parit{Proof of \ref{lem:strong-angle}.}
The first statement follows directly from the first variation formula (\ref{1st-var+}) 
and the definition of model angle (see Section~\ref{sec:angles}).

The second statement follows from Plaut's theorem (\ref{thm:almost.geod}) applied to $\spc{L}^\o$.
(Note that according to Proposition~\ref{prp:A^omega}, $\spc{L}^\o$ is a complete length $\Alex{}$ space.)
\qeds

\parit{Proof of \ref{1st-var+}.}
By Corollary~\ref{cor:CAT>k-sence}, we can assume that $\kappa\le 0$.
The inequality 

\[\dist{q}{\geod_{[p x]}(t)}{}
\le\dist{q}{p}{}-t\cdot\cos\phi+o(t)\]
follows from the first variation inequality (\ref{lem:first-var}).
Thus, it is sufficient to show that
\[\dist{q}{\geod_{[p x]}(t)}{}
\ge\dist{q}{p}{}-t\cdot\cos\phi+o(t).\]
Assume the contrary. Then there is $\eps>0$ such that  $\phi+\eps<\pi$,
and for a sequence $t_n\to 0+$ we have
\[\dist{q}{\geod_{[p x]}(t_n)}{}
<
\dist{q}{p}{}-t_n\cdot\cos(\phi-\eps).
\eqlbl{eq:phi-eps}\]

Let $x_n=\geod_{[p x]}(t_n)$.
Clearly 
\[\angk\kappa {x_n}p q>\pi-\phi+\tfrac\eps2\]
for all large $n$.

Assume $\spc{L}$ is geodesic. 
Choose a sequence of geodesics $[x_n q]$.
Let $[x_n q]\to [pq]_{\spc{L}^\o}$  as $n\to\o$ (in general $[pq]$ might lie in $\spc{L}^\o$).
Applying both parts of hinge comparison (\ref{angle}), 
we have $\mangle\hinge {x_n}{q}{x}\z<\phi-\tfrac\eps2$ for all large~$n$.
According to \ref{lem:ang.semicont-cbb}, the angle between $[pq]$ and $[px]$ is at most $\phi-\tfrac\eps2$, a contradiction.

Finally, if $\spc{L}$ is not geodesic, choose a sequence $q_n\in\Str(x_n)$, such that $q_n\to q$ and the inequality 
\[\angk\kappa{x_n}{p}{q_n}\z>\pi-\phi+\tfrac\eps2\] still holds.
Then the same argument as above shows that $[x_n q_n]$ $\o$-converges to a geodesic  $[pq]_{\spc{L}^\o}$ from $p$ to $q$  having angle at most $\phi-\tfrac{\eps}{2}$ with $[px]$.
\qeds


\section{On positive lower bound}\label{sec:positive.bound}

In this section we consider $\Alex{\kappa}$ spaces for $\kappa>0$.
Applying rescaling we can assume that $\kappa=1$.

The following theorem states that if one ignores a few exceptional spaces, then the diameter of a space with positive lower curvature bound is bounded.
Many authors (but not us) exclude these spaces in the definition of Alexandrov space with positive lower curvature bound.

\begin{thm}{On diameter of a space}\label{diam-k>0}
Let $\spc{L}$ be a complete length $\Alex1$ space. 
Then either 
\begin{subthm}{} $\diam \spc{L}\le \pi$; 
\end{subthm}

\begin{subthm}{} or $\spc{L}$ is isometric to one of the following exceptional spaces: 
\begin{enumerate}
\item real line $\RR$,
\item a half-line $\RR_{\ge0}$,
\item a closed interval $[0,a]\in \RR$, $a>\pi$,
\item a circle $\mathbb{S}^1_a$ of length $a>2\cdot\pi$.
\end{enumerate}
\end{subthm}
\end{thm}

\parit{Proof.} 
Assume that $\spc{L}$ is a geodesic space and $\diam \spc{L}>\pi$. 
Choose $x, y\in \spc{L}$ so that $\dist{x}{y}{}=\pi+\eps$, $0<\eps<\tfrac{\pi}{4}$. 
By moving $y$ slightly, we can also assume that the  geodesic $[x y]$ is unique;
to prove this, use either Plaut's theorem (\ref{thm:almost.geod}) 
or the fact that  geodesics do not split (\ref{thm:g-split}).
Let $z$ be the midpoint of the geodesic $[x y]$.

Consider the function $f=\distfun{x}{}{}+\distfun{y}{}{}$.
As follows from Lemma~\ref{concave-pi/2}, 
$f$ is concave in $\oBall(z,\tfrac{\eps}{4})$.  
Let $p\in\oBall(z,\tfrac{\eps}{4})$.  
Choose a geodesic $[z p]$, 
and let $h(t)=f\circ\geod_{[z p]}(t)$ and $\ell=\dist{z}{p}{}$.
Clearly $h$ is concave.
From the adjacent angle lemma (\ref{lem:sum=pi}), we have $h^+(0)=0$. 
Therefore $h$ is nonincreasing which means that \[\dist{x}{p}{}+\dist{y}{p}{}
=
h(\ell)\le h(0)
=
\dist{x}{y}{}.\]  
Since the geodesic $[x y]$ is unique this means that $p\in [x y]$, and hence $\oBall(z,\tfrac{\eps}{4})$ only contains points of $[x y]$.

Since in $\Alex{}$ spaces, geodesics do not bifurcate (\ref{nonsplit}), 
it follows that all of $\spc{L}$ coincides with the maximal extension of $[x y]$ as a local geodesic $\gamma$ 
(which might not be minimizing).
In other words, $\spc{L}$ is isometric to a 1-dimensional Riemannian manifold with possibly nonempty boundary.
From this, it is easy to see that $\spc{L}$ falls into one of the exceptional spaces described in the theorem.

Lastly, if $\spc{L}$ is not geodesic and $\diam \spc{L}>\pi$, then the above argument applied to $\spc{L}^\o$ yields that each metric component of $\spc{L}^\o$ is isometric to one of the exceptional spaces. 
As all of those spaces are proper, $\spc{L}$ is a metric component in $\spc{L}^\o$.
\qeds

\begin{thm}{Lemma}\label{concave-pi/2}
Let $\spc{L}$ be a complete length $\Alex1$ space and $p\in \spc{L}$.
Then $\distfun{p}{}{}\:\spc{L}\to\RR$ is concave in $\oBall(p,\pi)\setminus \oBall(p,\tfrac\pi2)$.

In particular, if $\diam\spc{L}\le\pi$,
then the complements $\spc{L}\setminus \oBall(p,r)$ and $\spc{L}\setminus \cBall[p,r]$ are convex for any $r>\tfrac\pi2$.

\end{thm}
\parit{Proof.}
This is a consequence of \ref{comp-kappa}. 
\qeds

\begin{thm}{Exercise}\label{ex:fixed-point}
Let $\spc{L}$ be a complete length $\Alex1$ space that is not exceptional (that is, $\diam\spc{L}\le\pi$).
Assume that a group $G$ acts on  $\spc{L}$ by isometries, has closed orbits, and 
\[\diam(\spc{L}/G)>\tfrac\pi2.\]
Show that the action of $G$ has a fixed point in $\spc{L}$.
\end{thm}

\begin{thm}{Advanced exercise}\label{ex:kleiner}
Let $\spc{L}$ be a complete length $\Alex{1}$ space
Show that $\spc{L}$ contains at most 3 points with space of directions $\le\tfrac12\cdot\mathbb{S}^n$
(see Definition \ref{def: inequality-of-spaces}).
\end{thm}

{\sloppy 

\begin{thm}{On perimeter of a triple}\label{perim-k>0}
Suppose  
$\spc{L}$ is a complete length $\Alex1$ space
and $\diam \spc{L}\le \pi$.
Then the perimeter of any triple of points $p,q,r\in \spc{L}$ is at most $2\cdot\pi$.
\end{thm}

}

\parit{Proof.}  
Arguing by contradiction, suppose 
\[\dist{p}{q}{}+\dist{q}{r}{}+\dist{r}{p}{}> 2\cdot\pi\eqlbl{eq:perimeter-of-triange<2pi}\] 
for $p,q,r\in \spc{L}$. 
Rescaling the space slightly, we can assume that $\diam\spc{L}<\pi$,
but the inequality \ref{eq:perimeter-of-triange<2pi} still holds.
By Corollary \ref{cor:CAT>k-sence},
after rescaling $\spc{L}$ is still $\Alex1$.

Since $\spc{L}$ is G-delta geodesic (\ref{thm:almost.geod}), it is sufficient to consider the case when there is a geodesic $[q r]$. 

\begin{wrapfigure}{r}{35 mm}
\vskip-4mm
\centering
\includegraphics{mppics/pic-840}
\end{wrapfigure}

First note that since $\diam \spc{L}<\pi$, by \ref{comp-kappa} the function \[y(t)=\md1\dist[{{}}]{p}{\geod_{[q r]}(t)}{}\]
satisfies the differential inequality $y''\le 1- y$.

Take $z_0\in [q r]$ so that the restriction $\distfun{p}{}{}|_{[q r]}$ attains its maximum at $z_0$, 
and set $t_0=\dist{q}{z_0}{}$ so $z_0=\geod_{[q r]}(t_0)$.
Consider the following model configuration: two geodesics $[\tilde p\tilde z_0]$, $[\tilde q\tilde r]$ in $\mathbb{S}^2$ such that 
\begin{align*}
\dist{\tilde p}{\tilde z_0}{}&=\dist{p}{z_0}{},
&  
\dist{\tilde q}{\tilde r}{}&=\dist{q}{r}{},
\\ 
\dist{\tilde z_0}{\tilde q}{}&=\dist{z_0}{q}{},
&  
\dist{\tilde z_0}{\tilde r}{}&=\dist{z_0}{q}{}
\end{align*}
and 
\[\mangle\hinge{\tilde z_0}{\tilde q}{\tilde p}
=\mangle\hinge{\tilde z_0}{\tilde r}{\tilde p}
=\tfrac\pi2.\]
Clearly,
$\bar y(t)=\md1\dist[{{}}]{\tilde p}{\geod_{[\tilde q\tilde r]}(t)}{}$ 
satisfies $\bar y''=1-\bar y$ and $\bar y'(t_0)=0$,
$\bar y(t_0)=y(t_0)$.
Since $z_0$ is a maximum point, 
$y(t)\le y(t_0)+o(t-t_0)$;
thus, $\bar y(t)$ is a barrier for 
$y(t)=\md1\dist[{{}}]{p}{\geod_{[q r]}(t)}{}$ at 
$t_0$ by \ref{barrier'}.
From the barrier inequality \ref{barrier'}, we get 
\[\dist{\tilde p}{\geod_{[\tilde q\tilde r]}(t)}{}
\ge 
\dist{p}{\geod_{[q r]}(t)}{},\]
and hence $\dist{\tilde p}{\tilde q}{}\ge\dist{p}{q}{}$ and $\dist{\tilde p}{\tilde r}{}\ge\dist{p}{r}{}$.

Therefore 
$\dist{p}{q}{}+\dist{q}{r}{}+\dist{r}{p}{}$ cannot exceed the perimeter of the  spherical triangle $\trig{\tilde p}{\tilde q}{\tilde r}$. 
In particular,
\[\dist{p}{q}{}+\dist{q}{r}{}+\dist{r}{p}{}\le 2\cdot\pi\]
--- a contradiction.
\qeds

Let $\kappa>0$.
Consider the following extension $\angk{\kappa+}{{*}}{{*}}{{*}}$ 
of the model angle function $\angk\kappa{{*}}{{*}}{{*}}$.
This definition works well for $\Alex{}$ spaces; for $\CAT{}$ spaces there is a similar but definition.
Some authors define the comparison angle to be $\angk{\kappa+}{{*}}{{*}}{{*}}$.

\begin{thm}{Definition of extended angle}\label{def:extended-angle}
Suppose $p,q,r$ are points in a metric space, and $p\ne q$, $p\ne r$. 
Let
\[\angk{\kappa+} p q r=\sup\set{\angk{\Kappa} p q r}{\Kappa\le\kappa}.\]
The value $\angk{\kappa+} p q r$ is called the \index{model angle!extended model angle}\emph{extended model angle} of the triple $p$, $q$, $r$.
\end{thm}

\begin{thm}{Extended angle comparison}
Let $\kappa>0$ 
and $\spc{L}$ be a complete length 
$\Alex\kappa$ space.
Then for any hinge 
$\hinge p q r$ we have 
$\mangle\hinge p q r
\ge\angk{\kappa+} p q r$.

Moreover, the extended model angle  $\angk{\kappa+} p q r$ can be calculated using the following rule:

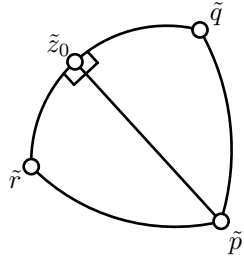
\begin{subthm}{} $\angk{\kappa+} p q r=\angk{\kappa} p q r$ if $\angk{\kappa} p q r$ is defined;
\end{subthm}

\begin{subthm}{} $\angk{\kappa+} p q r=\angk{\kappa+} p r q=0$ if $\dist{p}{q}{}+\dist{q}{r}{}=\dist{p}{r}{}$;
\end{subthm}

\begin{subthm}{} $\angk{\kappa+} p q r=\pi$ if none of the above is applicable. 
\end{subthm}
\end{thm}

\parit{Proof.}
From Corollary~\ref{cor:CAT>k-sence}, $\Kappa<\kappa$ implies that any complete length $\Alex{\Kappa}$ space is $\Alex{\kappa}$; 
thus the extended angle comparison follows from the definition.

The rule for calculating extended angle is an easy consequence of its definition.
\qeds


\section{Remarks}

The question whether the first part of \ref{angle} suffices to conclude that $\spc{L}$ is $\Alex\kappa$ is a long-standing open problem (possibly dating back to Alexandrov);
it was first stated in \cite[footnote in 4.1.5]{burago-burago-ivanov}.

\begin{thm}{Open question}\label{open:hinge-}
Let $\spc{L}$ be a complete geodesic space (you can also assume that $\spc{L}$ is homeomorphic to $\mathbb{S}^2$ or $\RR^2$) 
such that for any hinge $\hinge x p y$ in $\spc{L}$, 
the angle $\mangle\hinge x p y$ is defined and 
\[\mangle\hinge x p y\ge\angk0 x p y.\]
Is it true that $\spc{L}$ is $\Alex{0}$?
\end{thm}

\parbf{Examples and constructions.}
Let us list important sources of examples of $\Alex{}$ spaces.
We do not provide all the proofs and some proofs are deferred to later chapters.

\textit{Complete Riemannian manifolds} with sectional curvature at least $\kappa$, their \textit{Gromov--Hausdorff limits}, as well as their \textit{ultralimits}, are $\Alex\kappa$.
This statement follows from \ref{prp:A^omega}, \ref{thm:ultra-GH} and the Toponogov comparison which is a partial case of the globalization theorem (\ref{thm:glob}).
For example, 
if $M$ is a Riemannian manifold of nonnegative sectional curvature,
then the limit of its rescalings $\tfrac1n\cdot M$ as $n\to \infty$ is $\Alex0$; this is the so-called \textit{asymptotic cone} of $M$.

Most of applications to  Riemannian geometry are based on the described sources and the following corollary of Gromov's selection theorem (\ref{thm:gromov-selection}) and the Bishop--Gromov inequality.

\begin{thm}{Gromov compactness theorem}
Let $M_n$ be a sequence of Riemannian manifolds with sectional curvature at least $\kappa$. 
Then, for any choice of marked points $p_n\in M_n$,
a subsequence of $M_n$ admits a Gromov--Hausdorff convergence 
such that the corresponding subsequence of $p_n$ converges.
\end{thm}

By Corollary \ref{cor:submetry-cbb}, the \textit{target of submetry} from $\Alex\kappa$ is $\Alex\kappa$.
In particular, if $M$ is a Riemannian manifold with sectional curvature at least $\kappa$ and 
$G$ is a closed subgroup of isometries on $M$,
then the \textit{quotient space} $M/G$ is $\Alex\kappa$.

Yet another source is given by \textit{convex hypersurfaces}.
Namely, suppose $M$ is a complete Riemannian manifold of sectional curvature at least~$\kappa$.
Let $N$ be a closed (as a subset) hypersurface in $M$.
Then $N$ with the induced inner metric is $\Alex{\kappa}$.
In particular, any closed  convex hypersurface in $\R^n$ is $\Alex{0}$.

The smooth case of this statement follows from the Gauss formula.
This has been generalized by Sergei Buyalo \cite{buyalo} to the nonsmooth case and sharpened by the authors \cite{alexander-kapovitch-petrunin-buyalo}.

Let us mention that an analogous statement about convex hypersurfaces in $\Alex{\kappa}$ space is completely open.
Namely, it is unknown if \textit{the boundary of a complete finite-dimensional $\Alex{\kappa}$ length space has to be $\Alex{\kappa}$}.

Further, $\Alex{}$ spaces behave nicely with respect to several natural constructions.
For example, the product of $\Alex{0}$ spaces is a $\Alex{0}$.
Also, the Euclidean cone over $\Alex{1}$ space is a $\Alex{0}$.
These are the first examples of the so-called \textit{warped products} that are discussed in Chapter~\ref{chapter:warped products};
a general statement is given in \ref{thm:warp-CBB}.
Perelman's doubling theorem can be considered as a partial case;
it states that if $\spc{L}$ is a finite-dimensional $\Alex{\kappa}$ length space with nonempty boundary, then its \textit{doubling across the boundary} is $\Alex{\kappa}$ as well.

More conceptually, \textit{Wasserstein space} of order $2$ over $\Alex0$ space is $\Alex0$.
Also, there is a natural metric on the \textit{space of metric-measure spaces} that makes it $\Alex0$ space;
it was constructed by Karl-Theodor Sturm \cite{sturm-2012}.

Among less important examples, let us mention \textit{polyhedral spaces},
an if-and-only-if condition is given in \ref{thm:poly-CBB}.
Also, in addition to Perelman's doubling theorem, 
there are several versions of gluing theorems \cite{petrunin1997, ge-li, Kosovski};
they give conditions that guarantee that gluing of two (or more) spaces is $\Alex{\kappa}$.

\chapter{Fundamentals of curvature bounded above}
\chaptermark{Fundamentals of CBA}

\section{Four-point comparison.} \label{sec:cba-def}

\index{$\CAT{}$}
\begin{thm}{Four-point comparison}\label{def:2+2}
A quadruple of points $p^1,p^2,x^1,x^2$ in a metric space 
satisfies 
\index{$\CAT{}$!$\CAT\kappa$ comparison}
\emph{$\CAT\kappa$ comparison}
if
  
\begin{subthm}{}
$\angk{\kappa}{p^1}{x^1}{x^2} 
\le 
\angk{\kappa}{p^1}{p^2}{x^1}+\angk{\kappa}{p^1}{p^2}{x^2}$, or
\end{subthm}

\begin{subthm}{}
$\angk{\kappa} {p^2}{x^1}{x^2}\le \angk{\kappa} {p^2}{p^1}{x^1} + \angk{\kappa} {p^2}{p^1}{x^2}$, or
\end{subthm}

\begin{subthm}{}
one of the six model angles 
\begin{align*}
\angk{\kappa}{p^1}{x^1}{x^2},\quad&\angk{\kappa}{p^1}{p^2}{x^1},\quad\angk{\kappa}{p^1}{p^2}{x^2},
\\
\angk{\kappa}{p^2}{x^1}{x^2},\quad&\angk{\kappa}{p^2}{p^1}{x^1},\quad\angk{\kappa}{p^2}{p^1}{x^2}
\end{align*}
is undefined.
\end{subthm}
\end{thm}

\begin{wrapfigure}{r}{30 mm}
\vskip-0mm
\centering
\includegraphics{mppics/pic-905}
\end{wrapfigure}

Here is a more intuitive formulation.

\begin{thm}{Reformulation}\label{def:2+2-reformulated}
Let $\spc{X}$ be a metric space.
A quadruple $p^1,p^2,x^1,x^2\in \spc{X}$ satisfies 
$\CAT\kappa$ comparison if one of the following holds:
\begin{subthm}{}
One of the triples 
$(p^1,p^2,x^1)$ 
or 
$(p^1, p^2, x^2)$ 
has perimeter $>2\cdot\varpi\kappa$.
\end{subthm}

\begin{subthm}{}
If $\trig{\tilde p^1}{\tilde p^2}{\tilde x^1}
=
\modtrig\kappa(p^1 p^2 x^1)$ 
and
$\trig{\tilde p}{\tilde p^2}{\tilde x^2}
\z=
\modtrig\kappa p^1 p^2 x^2$, then
\[\dist{\tilde x^1}{\tilde z}{}+\dist{\tilde z}{\tilde x^2}{}\ge \dist{x^1}{x^2}{},\]
for any $\tilde z\in[\tilde p^1\tilde p^2]$.

\end{subthm}

\end{thm}

\begin{thm}{Definition}
\label{def:ccat}
Let $\spc{U}$ be a metric space.

\begin{subthm}{}
$\spc{U}$ is 
\index{$\CAT{}$!$\CAT{\kappa}$ space} 
$\CAT{\kappa}$ 
if any quadruple $p^1,p^2,x^1,x^2\in \spc{X}$  satisfies  $\CAT{\kappa}$ comparison.
\end{subthm}

\begin{subthm}{}
$\spc{U}$ is 
\index{$\CAT{}$!locally $\CAT{\kappa}$ space}
\emph{locally $\CAT{\kappa}$} 
if any point $q\in \spc{U}$ admits a neighborhood $\Omega\ni q$ such that any quadruple $p^1,p^2,x^1,x^2\in \spc{X}$  satisfies  $\CAT\kappa$ comparison.
\end{subthm}

\begin{subthm}{}
$\spc{U}$  is a  
\index{$\CAT{}$!$\CAT{}$ space}
$\CAT{}$ space if  $\spc{U}$  is $\CAT{\kappa}$ for some $\kappa\in\RR$.
\end{subthm}
\end{thm}



The condition $\spc{U}$ is $\CAT\kappa$ should be understood as ``$\spc{U}$ has global curvature $\le\kappa$''.
In Proposition~\ref{prop:inherit-bound}, it will be shown that this formulation makes sense; 
in particular, if $\kappa\le\Kappa$, then any $\CAT\kappa$ space is $\CAT\Kappa$.

This terminology was introduced by Michael Gromov;  
$\CAT{}$ stands for \'Elie Cartan, Alexandr Alexandrov, and Victor Toponogov.
Originally these spaces were called \index{$\mathfrak{R}_\kappa$ domain}\emph{$\mathfrak{R}_\kappa$ domains};
this is Alexandrov's terminology and is still in use.

\begin{thm}{Exercise}\label{ex:ccat-(3+1)}
Let $\spc{U}$ be a metric space.
Show that $\spc{U}$ is $\CAT\kappa$
if and only if every quadruple of points in $\spc{U}$ admits a labeling by $(p,x^1,x^2,x^3)$ such that the three angles 
$\angk\kappa p{x^1}{x^2}$,
$\angk\kappa p{x^2}{x^3}$ and
$\angk\kappa p{x^1}{x^3}$
satisfy all three triangle inequalities or one of these angles is undefined.
\end{thm}

\begin{thm}{Exercise}\label{ex:sba-2+2-short}
Show that $\spc{U}$ is $\CAT\kappa$
if and only if for any quadruple of points 
$p^1,p^2,x^1,x^2$ in $\spc{U}$ such that
$\dist{p^1}{p^2}{},\dist{x^1}{x^2}{}\le \varpi\kappa$,
there is a quadruple $q^1,q^2,y^1,y^2$ in $\Lob m\kappa$
such that 
\begin{align*}
\dist{q^1}{q^2}{}&=\dist{p^1}{p^2}{},
&
\dist{y^1}{y^2}{}&=\dist{x^1}{x^2}{},
&
\dist{q^i}{y^j}{}&\le \dist{p^i}{x^j}{}
\end{align*}
for any $i$ and $j$.
\end{thm}

\begin{thm}{Advanced exercise}\label{ex:berg-nikolaev}
Let $\spc{U}$ be a complete length space such that for any quadruple $p,q,x,y\in\spc{L}$ 
the following inequality holds
\[\dist[2]{p}{q}{}+\dist[2]{x}{y}{}\le \dist[2]{p}{x}{}
+\dist[2]{p}{y}{}+\dist[2]{q}{x}{}+\dist[2]{q}{y}{}.
\eqlbl{eq:berg-nikolaev-CAT}\]
Prove that $\spc{U}$ is $\CAT0$.

Construct a 4-point metric space $\spc{X}$ that satisfies inequality \ref{eq:berg-nikolaev-CAT} for any relabeling of its points by $p,q,x,y$, and such that $\spc{X}$ is not $\CAT{0}$.
\end{thm}

The next proposition follows directly from Definition \ref{def:ccat} and the definitions of ultralimit and ultrapower;
see Section~\ref{ultralimits} for the related definitions.
Recall that $\o$ denotes a fixed selective ultrafilter on $\NN$.

\begin{thm}{Proposition}\label{prop:CAT-olim}
\label{prop:CAT^omega}
Let $\spc{U}_n$ be a $\CAT{\kappa_n}$ space for each $n\in\NN$.
Assume $\spc{U}_n\to \spc{U}_\o$ and $\kappa_n\to\kappa_\o$ as $n\to\o$.
Then $\spc{U}_\o$ is $\CAT{\kappa_\o}$.

Moreover, a metric space $\spc{U}$ is $\CAT\kappa$ if and only if so is its ultrapower~$\spc{U}^\o$.

\end{thm} 

\section{Geodesics}

\begin{thm}{Uniqueness of geodesics}\label{thm:cat-unique}\label{thm:cat-complete} 
In a complete length $\CAT\kappa$ space, pairs of points at distance $<\varpi\kappa$ are joined by unique geodesics, and these geodesics depend continuously on their endpoint pairs.
\end{thm}

\parit{Proof.} 
Fix a complete length $\CAT\kappa$ space $\spc{U}$.
Fix two points $p^1,p^2\in \spc{U}$  such that 
\[\dist{p^1}{p^2}{\spc{U}}<\varpi\kappa.\]

Choose a sequence of approximate midpoints $z_n$ between $p^1$ and $p^2$;
that is,  
\[\dist{p^1}{z_n}{},\dist{p^2}{z_n}{}
\to\tfrac12\cdot\dist[{{}}]{p^1}{p^2}{}
\quad\text{as}\quad n\to\infty.
\eqlbl{eq:to|p1p2|/2}\]

By the law of cosines, $\angk{\kappa} {p^1}{z_n}{p^2}$ and $\angk{\kappa} {p^2}{z_n}{p^1}$ are arbitrarily small when $n$ is sufficiently large.

Let us apply $\CAT\kappa$  comparison (\ref{def:2+2}) to the quadruple $p^1$, $p^2$, $z_n$, $z_\kay$ with large $n$ and $\kay$.
We conclude that  $\angk{\kappa} {p}{z_n}{z_\kay}$ is arbitrarily small when $n,\kay$ are sufficiently large and $p$ is either $p^1$ or $p^2$.  
By \ref{eq:to|p1p2|/2} and the law of cosines, the sequence $z_n$ converges.  

Since $\spc{U}$ is complete, the sequence $z_n$ converges to a midpoint between $p^1$ and $p^2$. 
By Lemma~\ref{lem:mid>geod} we obtain  the existence of a geodesic $[p^1p^2]$.

Now suppose $p^1_n\to p^1$, $p^2_n\to p^2$ as $n\to\infty$.
Let $z_n$ be the midpoint of a geodesic $[p^1_n p^2_n]$ and $z$ be the midpoint of a geodesic $[p^1p^2]$.  

It suffices to show that 
\[\dist{z_n}{z}{}\to0
\quad \text{as}\quad 
n\to\infty.
\eqlbl{eq:z_n->z}\]

By the triangle inequality, the $z_n$ are approximate midpoints between $p^1$ and $p^2$.
Apply the $\CAT\kappa$ comparison (\ref{def:2+2}) to the quadruple $p^1$, $p^2$, $z_n$,~$z$. 
For $p=p^1$ or $p=p^2$, we see that $\angk{\kappa} {p}{z_n}{z}$ is arbitrarily small when $n$ is sufficiently large.  
By the law of cosines, \ref{eq:z_n->z} follows.
\qeds

\begin{thm}{Exercise}\label{ex:CAT-mnfld=>ext.geod}
Let $\spc{U}$ be a complete length $\CAT{}$ space.
Assume $\spc{U}$ is a topological manifold.
Show that any geodesic in $\spc{U}$ can be extended 
as a two-side infinite local geodesic.

Moreover the same holds for any locally geodesic locally $\CAT{}$ space $\spc{U}$ with nontrivial local homology groups at any point;
the latter holds in particular if $\spc{U}$ is a homological manifold.
\end{thm}

\begin{thm}{Exercise}\label{ex:complete-space-of-dir}
Assume $\spc{U}$ is a locally compact geodesic $\CAT{}$ space with extendable geodesics;
that is, any geodesic in $\spc{U}$ can be extended to a both-sided infinite local geodesic.

Show that the space of geodesic directions $\Sigma_p'$ is complete for any $p\in \spc{U}$.
\end{thm}

By the uniqueness of geodesics (\ref{thm:cat-unique}),
we have the following.

\begin{thm}{Corollary}\label{cor:cat-ccat}
Any  complete length $\CAT\kappa$ space is $\varpi\kappa$-geodesic.

\end{thm}

\begin{thm}{Proposition}\label{cor:cat-completion} 
The completion $\bar{\spc{U}}$ of any geodesic $\CAT{\kappa}$ space $\spc{U}$ is a complete length $\CAT\kappa$ space.

Moreover, $\spc{U}$ is a geodesic $\CAT\kappa$ space
if and only if there is a complete length $\CAT\kappa$ space $\bar{\spc{U}}$ that contains a $\varpi\kappa$-convex dense set isometric to $\spc{U}$.
\end{thm}

\parit{Proof.} 
By Theorem \ref{thm:cat-complete},
in order to show that  $\bar{\spc{U}}$ is $\CAT{\kappa}$,
it is sufficient to verify that the completion of a length space is a length space; 
this is straightforward.

For the second part, note that the completion $\bar{\spc{U}}$
contains the original space $\spc{U}$ as a dense $\varpi\kappa$-convex subset, and the metric on $\spc{U}$ coincides with the induced length metric from $\bar{\spc{U}}$.
\qeds

Here is a corollary from Proposition~\ref{cor:cat-completion}
and Theorem~\ref{thm:cat-unique}.

\begin{thm}{Corollary}\label{cor:cat-unique}
Let $\spc{U}$ be a  $\varpi\kappa$-geodesic $\CAT\kappa$ space.
Then pairs of points in $\spc{U}$ at distance less than $\varpi\kappa$ are joined by unique geodesics, and these geodesics depend continuously on their endpoint pairs.

Moreover for any pair of points $p,q\in \spc{U}$ and any value
\[\Lip>\sup\set{\frac{\sn\kappa r}{\sn\kappa \dist{p}{q}{}}}{0\le r\le \dist{p}{q}{}}\]
there are neighborhoods $\Omega_p\ni p$ and $\Omega_q\ni q$ such that the map
\[(x,y,t)\mapsto \geodpath_{[xy]}(t)\]
is $\Lip$-Lipschitz in $\Omega_p\times \Omega_q\times[0,1]$.
\end{thm}

\parit{Proof.}
By Proposition~\ref{cor:cat-completion}, any geodesic $\CAT{\kappa}$ space is isometric to a convex dense subset of a complete length $\CAT\kappa$ space.
It remains to apply  Theorem~\ref{thm:cat-unique}.
\qeds


\section{More comparisons}\label{sec:cat-angles}

Here we give a few reformulations of Definition~\ref{def:ccat}.

\begin{wrapfigure}{r}{25 mm}
\vskip-0mm
\centering
\includegraphics{mppics/pic-910}
\end{wrapfigure}

\begin{thm}{Theorem}
\label{thm:defs_of_cat} 
If $\spc{U}$ is a $\CAT\kappa$ space, then 
the following conditions hold for all triples $p,x,y\in \spc{U}$ of perimeter $<2\cdot\varpi\kappa$:

\begin{subthm}{cat-2-sum} (adjacent angle comparison\index{comparison!adjacent angle comparison}) for any geodesic $[x y]$ and $z\in \mathopen{]}x y\mathclose{[}$, we have
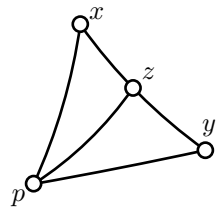
\[\angk\kappa z p x
+\angk\kappa z p y\ge \pi.\]
\end{subthm}

\begin{subthm}{cat-monoton}
(point-on-side comparison\index{comparison!point-on-side comparison}) 
for any geodesic $[x y]$ and $z\in \mathopen{]}x y\mathclose{[}$, we have
\[\angk\kappa x p y\ge\angk\kappa x p z,\]
or equivalently, 
\[\dist{\tilde p}{\tilde z}{}\ge \dist{p}{z}{},\]
where $\trig{\tilde p}{\tilde x}{\tilde y}=\modtrig\kappa(p x y)$, $\tilde z\in\mathopen{]} \tilde x\tilde y\mathclose{[}$, $\dist{\tilde x}{\tilde z}{}=\dist{x}{z}{}$.
\end{subthm}

\begin{subthm}{cat-hinge}(hinge comparison\index{comparison!hinge comparison})
for any hinge $\hinge x p y$, the angle 
$\mangle\hinge x p y$ exists and
\[\mangle\hinge x p y\le\angk\kappa x p y,\]
or equivalently,
\[\side\kappa \hinge x p y\le\dist{p}{y}{}.\]
\end{subthm}

Moreover, if  $\spc{U}$ is  $\varpi\kappa$-geodesic, then the converse holds in each case.  

\end{thm}

\parbf{Remark.}
\label{22remark}
In the following proof, the part \ref{SHORT.cat-hinge}$\Rightarrow$\ref{SHORT.cat-2-sum})
only requires that the $\CAT\kappa$ comparison (\ref{def:2+2}) hold for any quadruple, and does not require the existence of geodesics at distance $<\varpi\kappa$. 
The same is true of the parts \ref{SHORT.cat-2-sum}$\Leftrightarrow$\ref{SHORT.cat-monoton} and
\ref{SHORT.cat-monoton}$\Rightarrow$\ref{SHORT.cat-hinge}.  
Thus the conditions \ref{SHORT.cat-2-sum}, \ref{SHORT.cat-monoton}) and \ref{SHORT.cat-hinge} are valid for any metric space (not necessarily a length space) that satisfies $\CAT\kappa$ comparison (\ref{def:2+2}). 
The converse does not hold; for example, all these conditions are 
vacuously true in a 
totally disconnected space, while 
$\CAT\kappa$ comparison is not.

\parit{Proof; \ref{SHORT.cat-2-sum}}. 
Since the perimeter of $p,x,y$ is $<2\cdot \varpi\kappa$, so is the perimeter of any subtriple of $p,z,x,y$ by the triangle inequality. 
By Alexandrov's lemma (\ref{lem:alex}), 
\[\angk\kappa p z x +\angk\kappa p z y  < \angk{\kappa} p x y \quad \text{or}\quad  \angk\kappa z p x  +\angk\kappa z p y  =\pi.\]
In the former case, the $\CAT\kappa$ comparison (\ref{def:2+2}) applied to the quadruple $p, z, x, y$ implies
\[\angk\kappa z p x  +\angk\kappa z p y  \ge \angk{\kappa} z x y =\pi.\]

\parit{\ref{SHORT.cat-2-sum}$\Leftrightarrow$\ref{SHORT.cat-monoton}.}
Follows from  Alexandrov's lemma (\ref{lem:alex}).

\parit{\ref{SHORT.cat-monoton}$\Rightarrow$\ref{SHORT.cat-hinge}.} 
By \ref{SHORT.cat-monoton}, for $\bar p\in\mathopen{]}x p]$ and $\bar y\in\mathopen{]}x y]$ the function $(\dist{x}{\bar p}{},\dist{x}{\bar y}{})\mapsto\angk\kappa x{\bar p}{\bar y}$ is nondecreasing in each argument.
In particular, 
$\mangle\hinge x p y\z=\inf\angk\kappa x{\bar p}{\bar y}$.
Thus $\mangle\hinge x p y$ exists and is
at most $\angk\kappa x p y$. 

\parit{Converse.} Assume $\spc{U}$ is $\varpi\kappa$-geodesic. 
Let us first show that in this case \ref{SHORT.cat-hinge}$\Rightarrow$\ref{SHORT.cat-2-sum}.

\begin{wrapfigure}{r}{30 mm}
\vskip-4mm
\centering
\includegraphics{mppics/pic-915}
\end{wrapfigure}

Indeed, by \ref{SHORT.cat-hinge} and the triangle inequality for angles (\ref{claim:angle-3angle-inq}),
\[\angk\kappa z p x
+\angk\kappa z p y \ge \mangle\hinge z p x
+\mangle\hinge z p y \ge \pi.\]
It remains to prove the converse for \ref{SHORT.cat-monoton}.

Given a quadruple  $p^1,p^2,x^1,x^2$ whose subtriples have perimeter $<2\cdot\varpi\kappa$, we must verify the $\CAT\kappa$ comparison (\ref{def:2+2}).
In $\Lob2\kappa$, construct the model triangles  $\trig{\tilde p^1}{\tilde p^2}{\tilde x^1} = \modtrig\kappa(p^1 p^2 x^1 )$ 
and $\trig{\tilde p^1}{\tilde p^2}{\tilde x^2}\z= \modtrig\kappa(p^1 p^2 x^2)$, lying on either side of a common segment $[\tilde p^1 \tilde p^2]$.
We may suppose 
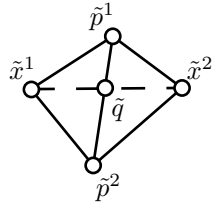
\[\angk{\kappa} {p^1}{p^2}{x^1}+\angk{\kappa} {p^1}{p^2}{x^2}
\le
\pi
\quad \text{and}\quad 
\angk{\kappa}{p^2}{p^1}{x^1}+\angk{\kappa} {p^2}{p^1}{x^2}
\le 
\pi,\] 
since otherwise $\CAT\kappa$ comparison holds trivially.  
Then $[\tilde p^1 \tilde p^2]$ and $[\tilde x^1 \tilde x^2]$ intersect, say at $\tilde q$.  

By assumption, there is a geodesic $[p^1 p^2]$.
Choose $q\in[p^1 p^2]$ corresponding to $\tilde q$; 
that is, $\dist{p^1}{q}{}=\dist{\tilde p^1}{\tilde q}{}$.
Then 
\[\dist{x^1}{x^2}{} \le \dist{x^1}{q}{} + \dist{q}{x^2}{} \le \dist{\tilde x^1}{\tilde q}{} + \dist{\tilde q}{\tilde x^2}{} = \dist{\tilde x^1}{\tilde x^2}{},\]
where the second inequality follows from \ref{SHORT.cat-monoton}. 
By monotonicity of the function $a\mapsto\tangle\mc\kappa\{a;b,c\}$ (\ref{increase}),
\begin{align*}
\angk{\kappa} {p^1}{x^1}{x^2} \le  \mangle\hinge{ \tilde p^1}{ \tilde x^1}{ \tilde x^2}
= \angk{\kappa} {p^1}{p^2}{x^1} + \angk{\kappa} {p^1}{p^2}{x^2}.
\end{align*}
\qedsf

Let us display a corollary of the proof of \ref{thm:defs_of_cat},
namely, monotonicity of the model angle with respect to adjacent sidelengths. 

\begin{thm}{Angle-sidelength  monotonicity}\label{cor:monoton-cba} 
Suppose $\spc{U}$ is a $\varpi\kappa$-geodesic $\CAT\kappa$ space, and 
$p,x,y\in \spc{U}$ have  perimeter $<2\cdot \varpi\kappa$.
Then for $\bar y\in\mathopen{]}x y]$, the function 
\[\dist{x}{\bar y}{}\mapsto \angk\kappa x p{\bar y}\] 
is nondecreasing.

In particular, if $\bar p\in \mathopen{]}x p]$, then
\begin{subthm}{two-mono-cba}the function 
\[(\dist{x}{\bar y}{},\dist{x}{\bar p}{})\mapsto \angk\kappa x {\bar p}{\bar y}\] is nondecreasing in each argument,
\end{subthm}
 
\begin{subthm}{cor:monoton-cba:angle=inf} 
$\mangle\hinge{x}{p}{y}
=
\inf\set{\angk\kappa x {\bar p}{\bar y}}{
\bar p\in\mathopen{]}x p],\ 
\bar y\in\mathopen{]}x y]}.$
\end{subthm}
\end{thm}

\begin{thm}{Exercise}\label{mink+CAT=euclid} 
Let $\spc{U}$ be  $\RR^m$ with the metric defined by a norm.
Show that $\spc{U}$ is a complete length $\CAT{}$ space if and only if $\spc{U}\iso\EE^m$.
\end{thm}

\begin{thm}{Exercise}\label{ex:convexity-CAT0}
Assume $\spc{U}$ is a geodesic $\CAT0$ space.
Show that for any two geodesic paths 
$\gamma,\sigma\:[0,1]\to \spc{U}$
the function 
\[t\mapsto\dist{\gamma(t)}{\sigma(t)}{}\] 
is convex.
\end{thm}

\begin{thm}{Proposition}
\label{prop:inherit-bound}
Assume $\kappa<\Kappa$.
Then any complete length $\CAT\kappa$ space is $\CAT\Kappa$.

Moreover a space $\spc{U}$ is $\CAT\kappa$ if  $\spc{U}$ is $\CAT\Kappa$ for all $\Kappa>\kappa$.
\end{thm}

\parit{Proof.}
The first statement follows from Corollary \ref{cor:cat-ccat}, the adjacent-angles comparison (\ref{cat-2-sum}) and the monotonicity of the function $\kappa\mapsto\angk\kappa x y z$ (\ref{k-decrease}).

The second statement follows since the function $\kappa\mapsto\angk\kappa x y z$ is continuous.
\qeds


\section{Thin triangles} \label{sec:thin-triangle}

In this section we define thin triangles
and use them to characterize $\CAT{}$ spaces.
Inheritance for thin triangles with respect to decomposition
is the main result of this section.
It will lead to two fundamental constructions:  
Alexandrov's patchwork globalization  (\ref{thm:alex-patch}) 
and Reshetnyak gluing (\ref{thm:gluing}).
 
\begin{thm}{Definition of $\bm\kappa$-thin triangles}\label{def:k-thin}
Let $\trig{x^1}{x^2}{x^3}$ be a triangle of perimeter $<2\cdot \varpi\kappa$ in a metric space
and
$\trig{\tilde x^1}{\tilde x^2}{\tilde x^3}\z=\modtrig\kappa({x^1}{x^2}{x^3})$.
Consider the \emph{natural map} $\trig{\tilde x^1}{\tilde x^2}{\tilde x^3}\to \trig{x^1}{x^2}{x^3}$ 
that sends a point $\tilde z\in[\tilde x^i\tilde x^j]$ to the corresponding point $z\in[x^ix^j]$
(that is, such that $\dist{\tilde x^i}{\tilde z}{}=\dist{x^i}{z}{}$ and therefore $\dist{\tilde x^j}{\tilde z}{}=\dist{x^j}{z}{}$).

We say the triangle $\trig{x^1}{x^2}{x^3}$ is \index{thin triangle}\emph{$\kappa$-thin} if the natural map $\trig{\tilde x^1}{\tilde x^2}{\tilde x^3}\to \trig{x^1}{x^2}{x^3}$ is short.
\end{thm}

\begin{thm}{Exercise}\label{ex:equality-for-thin}
Let $\spc{U}$ be a $\varpi\kappa$-geodesic $\CAT\kappa$ space.
Let $\trig xyz$ be a triangle in $\spc{U}$
and $\trig{\tilde x}{\tilde y}{\tilde z}$ be its model triangle in $\Lob{2}{\kappa}$.
Prove that the natural map $f\:\trig{\tilde x}{\tilde y}{\tilde z}\to \trig xyz$ 
 is distance-preserving if and only if one of the following conditions hold:

\begin{subthm}{ex:equality-for-thin:angle}
$\mangle\hinge x y z= \angk\kappa x y z$,
\end{subthm}

\begin{subthm}{ex:equality-for-thin:vertex-base}
$\dist{x}{w}{}=\dist{\tilde x}{\tilde w}{}$ for some  $\tilde w\in]\tilde y\tilde z[$ and
$w= f(\tilde w)$,   
\end{subthm}

\begin{subthm}{ex:equality-for-thin:side-side} 
$\dist{v}{w}{}=\dist{\tilde v}{\tilde w}{}$ for some  
$\tilde v\in \mathopen{]}\tilde x \tilde y\mathclose{[}$,
$\tilde w\in\mathopen{]}\tilde x \tilde z\mathclose{[}$
and $v=f(\tilde v)$, $w=f(\tilde w)$.
\end{subthm} 

\end{thm}

{\sloppy 

\begin{thm}{Proposition}\label{prop:k-thin}
Let $\spc{U}$ be a $\varpi\kappa$-geodesic space. 
Then $\spc{U}$ is  $\CAT\kappa$
if and only if every triangle of perimeter $<2\cdot \varpi\kappa$ in $\spc{U}$  is $\kappa$-thin.
\end{thm}

}

\parit{Proof.}
The if part follows from the point-on-side comparison (\ref{cat-monoton}).  
The only-if part follows from angle-sidelength  monotonicity (\ref{two-mono-cba}).
\qeds

\begin{thm}{Corollary}\label{cor:loc-geod-are-min}
Suppose $\spc{U}$ is a $\varpi\kappa$-geodesic $\CAT\kappa$ space.  
Then any local geodesic in $\spc{U}$ of length $<\varpi\kappa$ is length-minimizing.
\end{thm}

\parit{Proof.}
Suppose $\gamma\:[0,\ell]\to\spc{U}$ is a local geodesic  that is not minimizing, with $\ell<\varpi\kappa$.
Choose $a$ to be the maximal value 
such that $\gamma$ is minimizing on $[0,a]$.
Further choose $b>a$ so that $\gamma$ is minimizing on~$[a,b]$.

Since triangle $\trig{\gamma(0)}{\gamma(a)}{\gamma(b)}$ is $\kappa$-thin, we have
\[\dist{\gamma(a-\eps)}{\gamma(a+\eps)}{}<2\cdot\eps\]
for all small $\eps>0$,
a contradiction.
\qeds

Now let us formulate the main result of this section.
The inheritance lemma states that  in any metric space, a triangle is $\kappa$-thin if it decomposes into $\kappa$-thin triangles. 
In contrast, $\kappa$-thickness of triangles (\ref{ex:fat-triangle}) is not inherited in this way.

\begin{wrapfigure}{r}{25 mm}
\vskip-0mm
\centering
\includegraphics{mppics/pic-920}
\end{wrapfigure}

\begin{thm}{Inheritance lemma}
\label{lem:inherit-angle} 
In a metric space, consider a triangle $\trig p x y$ that \index{decomposed triangle}\emph{decomposes} 
into two triangles $\trig p x z$ and $\trig p y z$;
that is, $\trig p x z$ and $\trig p y z$ have common side $[p z]$, and the sides $[x z]$ and $[z y]$ together form the side $[x y]$ of $\trig p x y$.

If the triangle $\trig p x y$ has perimeter $<2\cdot\varpi\kappa$
and both triangles $\trig p x z$ and $\trig p y z$ are $\kappa$-thin, then triangle $\trig p x y$ is  $\kappa$-thin.
\end{thm} 

The following model-space lemma is  extracted from Lemma 2 in \cite{reshetnyak:major}.

\begin{thm}{Lemma}\label{lem:quadrangle}
Let $\trig{\tilde p}{\tilde x}{\tilde y}$ be a triangle in $\Lob2{\kappa}$ and $\tilde z\in[\tilde x\tilde y]$.
Consider the solid triangle $\tilde D=\Conv\trig{\tilde p}{\tilde x}{\tilde y}$.  
Construct  points $\dot p, \dot x, \dot z, \dot y\in \Lob2{\kappa}$ such that 
\begin{align*}
\dist{\dot p}{\dot x}{}&=\dist{\tilde p}{\tilde x}{},
&
\dist{\dot p}{\dot y}{}&=\dist{\tilde p}{\tilde y}{},
&
\dist{\dot p}{\dot z}{}&\le \dist{\tilde p}{\tilde z}{},
\\
\dist{\dot x}{\dot z}{}&=\dist{\tilde x}{\tilde z}{},
&
\dist{\dot y}{\dot z}{}&=\dist{\tilde y}{\tilde z}{},
\end{align*}
where points $\dot x$ and $\dot y$ lie on either side of $[\dot p\dot z]$.
Set 
\[\dot D=\Conv\trig {\dot p}{\dot x}{\dot z}\cup \Conv\trig {\dot p} {\dot y} {\dot z}.\]

Then there is a short map $F\:\tilde D\to \dot D$ that maps $\tilde p$, $\tilde x$, $\tilde y$ and $\tilde z$ to $\dot p$, $\dot x$, $\dot y$ and $\dot z$ respectively.
\end{thm}

{

\begin{wrapfigure}{r}{41 mm}
\vskip-0mm
\centering
\includegraphics{mppics/pic-925}
\end{wrapfigure}

\parit{Proof.} 
By Alexandrov's lemma (\ref{lem:alex}), 
there are nonoverlapping triangles 
$\trig{\tilde p}{\tilde x}{\tilde z_y}\iso\trig {\dot p}{\dot x}{\dot z}$ 
and 
$\trig{\tilde p}{\tilde y}{\tilde z_x}\iso\trig {\dot p}{\dot y}{\dot z}$
 inside triangle $\trig{\tilde p}{\tilde x}{\tilde y}$.

Connect points in each pair
$(\tilde z,\tilde z_x)$, 
$(\tilde z_x,\tilde z_y)$ 
and $(\tilde z_y,\tilde z)$ 
with arcs of circles centered at 
$\tilde y$, $\tilde p$, and $\tilde x$ respectively. 
Define $F$ as follows.

}

\begin{itemize}
\item Map  $\Conv\trig{\tilde p}{\tilde x}{\tilde z_y}$ isometrically onto  $\Conv\trig {\dot p}{\dot x}{\dot y}$;
similarly map $\Conv \trig{\tilde p}{\tilde y}{\tilde z_x}$ onto $\Conv \trig {\dot p}{\dot y}{\dot z}$.

\item If $w$ is in one of the three circular sectors, say at distance $r$ from the center of the circle, let $F(w)$ be the point on  
$[\dot p \dot z]$, 
$[\dot x \dot z]$,
or $[\dot y \dot z]$ whose distance from the left-hand endpoint of the segment is $r$.
\item Finally, if $w$ lies in the remaining curvilinear triangle $\tilde z \tilde z_x \tilde z_y$, 
set $F(w) = \dot z$. 
\end{itemize}
By construction, $F$ meets the conditions of the lemma. 
\qeds

\parit{Proof of \ref{lem:inherit-angle}.}
Construct model triangles $\trig{\dot p}{\dot x}{\dot z}\z=\modtrig\kappa(p x z)$ 
and $\trig {\dot p} {\dot y} {\dot z}\z=\modtrig\kappa(p y z)$ so that $\dot x$ and $\dot y$ lie on opposite sides of $[\dot p\dot z]$.

\begin{wrapfigure}{r}{30 mm}
\vskip-0mm
\centering
\includegraphics{mppics/pic-930}
\end{wrapfigure}

Suppose
\[\angk\kappa{z}{p}{x}+\angk\kappa{z}{p}{y}
<
\pi.\]
Then for some point $\dot w\in[\dot p\dot z]$, we have \[\dist{\dot x}{\dot w}{}+\dist{\dot w}{\dot y}{}
<
\dist{\dot x}{\dot z}{}+\dist{\dot z}{\dot y}{}=\dist{x}{y}{}.\]
Let $w\in[p z]$ correspond to $\dot w$; that is, $\dist{z}{w}{}=\dist{\dot z}{\dot w}{}$. 
Since $\trig p x z$ and $\trig p y z$ are $\kappa$-thin, we have 
\[\dist{x}{w}{}+\dist{w}{y}{}<\dist{x}{y}{},\]
contradicting the triangle inequality. 

Thus 
\[\angk\kappa{z}{p}{x}+\angk\kappa{z}{p}{y}
\ge
\pi.\]
By Alexandrov's lemma (\ref{lem:alex}), this is equivalent to 
\[\angk\kappa x p z\le\angk\kappa x p y.
\eqlbl{eq:for|pz|}\]

Let $\trig{\tilde  p}{\tilde  x}{\tilde  y}=\modtrig\kappa (p x y)$ 
and $\tilde  z\in[\tilde  x\tilde  y]$ correspond to $z$; that is, $\dist{x}{z}{}=\dist{\tilde  x}{\tilde  z}{}$.
Inequality~\ref{eq:for|pz|} is equivalent to $\dist{ p}{ z}{}\le \dist{\tilde  p}{\tilde  z}{}$.
Hence  Lemma~\ref{lem:quadrangle} applies.  Therefore 
there is a short map $F$ that  sends 
$\trig{\tilde  p}{\tilde  x}{\tilde  y}$ to $\dot D=\Conv\trig {\dot p}{\dot x}{\dot z}\cup \Conv\trig {\dot p} {\dot y} {\dot z}$ 
in such a way that 
$\tilde p\mapsto \dot p$,
$\tilde x\mapsto \dot x$,
$\tilde z\mapsto \dot z$
and
$\tilde y\mapsto \dot y$.

By assumption, the natural maps $\trig {\dot p} {\dot x} {\dot z}\to\trig p x z$ and $\trig {\dot p} {\dot y} {\dot z}\to\trig p y z$ are short.  
By composition,  the natural map from $\trig{\tilde  p}{\tilde  x}{\tilde  y}$ to $\trig p y z$ is short, as claimed.
\qeds


\section{Function comparison} \label{sec:func-comp}

\index{comparison!function comparison}
In this section we give analytic and geometric ways of viewing the point-on-side comparison (\ref{cat-monoton}) as a convexity condition.

First we obtain a corresponding differential inequality for the distance function in $\spc{U}$;
see Section~\ref{sec:conv-fun} for the definition.
 
\begin{thm}{Theorem}\label{thm:function-comp} 
Suppose $\spc{U}$ is a $\varpi\kappa$-geodesic space. 
Then the following are equivalent:
\begin{subthm}{function-comp-cat} 
$\spc{U}$ is $\CAT\kappa$,
\end{subthm}
\begin{subthm}{function-comp}
for any $p\in \spc{U}$, the function $f=\md\kappa\circ\distfun{p}{}{}$ satisfies 
\[f''+\kappa \cdot f\ge 1\] 
in $\oBall(p,\varpi\kappa)$.
\end{subthm}\end{thm}

\begin{thm}{Corollary}
A geodesic space $\spc{U}$ is $\CAT{0}$ if and only if for any $p\in \spc{U}$, the function $\distfun[2]{p}{}{}\:\spc{U}\to\RR$ is $2$-convex.
\end{thm}

\parit{Proof of \ref{thm:function-comp}.}
Fix a sufficiently short geodesic $[x y]$ in $\oBall(p,\varpi\kappa)$.
We can assume that the model triangle $\trig{\tilde p}{\tilde x}{\tilde y}\z=\modtrig\kappa(p x y)$ is defined. 
Let \begin{align*} 
\tilde r(t)&=\dist{\tilde p}{\geod_{[\tilde x\tilde y]}(t)}{},
& 
r(t)&=\dist{p}{\geod_{[xy]}(t)}{}.                           \end{align*}
Let $\tilde f=\md\kappa\circ\tilde r$ and $f=\md\kappa\circ r$.
By \ref{md-diff-eq}, we have $\tilde f''\z=1-\kappa\cdot  \tilde f$.
Clearly $\tilde f(t)$ and $f(t)$ agree at $t=0$ and $t=\dist{x}{y}{}$. 
The point-on-side comparison (\ref{cat-monoton}) is the condition $r(t)\le\tilde r(t)$  for all $t\in[0,\dist{x}{y}{}]$.
Since $\md\kappa$ is increasing on $[0,\varpi\kappa)$, then $r\le \tilde r$ and $f\le \tilde f$ are equivalent.
Thus the claim follows by Jensen's inequality (\ref{y''-mono}).
\qeds

\begin{thm}{Corollary}\label{cor:convex-balls}
Suppose $\spc{U}$ is a $\varpi\kappa$-geodesic $\CAT\kappa$ space.
Then any ball (closed or open) of radius $R<\tfrac{\varpi\kappa}2$ in $\spc{U}$ is convex.

Moreover, any open ball of radius $\tfrac{\varpi\kappa}2$ is convex
and any closed ball of radius $\tfrac{\varpi\kappa}2$ is $\varpi\kappa$-convex.
\end{thm}

\parit{Proof.}
Suppose $p\in\spc{U}$, $ R\le\varpi\kappa/2$,  and two points 
$x$ and $y$ lie in $\cBall[p, R]$ or $\oBall(p, R)$.
By the triangle inequality, if $\dist{x}{y}{}<\varpi\kappa$, then any
 geodesic $[x y]$ lies in $\oBall(p, \varpi\kappa)$.
 
By the function comparison (\ref{thm:function-comp}), 
the geodesic $[x y]$ lies in $\cBall[p,R]$ or $\oBall(p,R)$ respectively.

Thus any ball (closed or open) of radius $R<\tfrac{\varpi\kappa}2$ is $\varpi\kappa$-convex.
This implies convexity unless there is a pair of points in the ball at distance at least $\varpi\kappa$.
By the  triangle inequality, the latter is possible only for the closed ball of radius $\tfrac{\varpi\kappa}2$.
\qeds

Recall that Busemann functions are defined in Proposition \ref{prop:busemann}.
The following exercise is analogous to Exercise~\ref{ex:busemann-CBB}.

\begin{thm}{Exercise}\label{ex:busemann-CBA}{\sloppy 
Let $\spc{U}$ be a complete length $\CAT\kappa$ space
and $\bus_\gamma\:\spc{U}\to \RR$ be the Busemann function for a half-line $\gamma\:[0,\infty)\to \spc{L}$.

}

\begin{subthm}{}
If $\kappa=0$, then the Busemann function $\bus_\gamma$ is  convex.
\end{subthm}

\begin{subthm}{}
If $\kappa=-1$, then the function $f=\exp\circ\bus_\gamma$ satisfies
\[f''- f\ge 0.\]
\end{subthm}

\end{thm}

\section{Development}\label{sec:development-CBA}
 
Geometrically,  the development construction (\ref{def:devel}) translates distance comparison into a local convexity statement for subsets of $\Lob2\kappa$.
Recall that a curve in $\Lob2\kappa$ is \emph{(locally) concave} with respect to $p$ if (locally) its supergraph with respect to $p$ is a convex subset of $\Lob2\kappa$; see Definition~\ref{def:convex-devel}.

\begin{thm}{Development criterion\index{comparison!development comparison}}\label{thm:concave-devel} 
For a $\varpi\kappa$-geodesic space $\spc{U}$,
the following statements hold:

\begin{subthm}{locally-concave-dev}
For any $p\in \spc{U}$ and any geodesic $\gamma\:[0,T]\to\oBall(p,\varpi\kappa)$, suppose the $\kappa$-development $\tilde \gamma$ in $\Lob2\kappa$ of $\gamma$ with respect to $p$ is locally concave. 
Then $\spc{U}$ is $\CAT\kappa$.
\end{subthm}

\begin{subthm}{concave-dev} 
If $\spc{U}$ is $\CAT\kappa$, then for any geodesic $\gamma\:[0,T]\to\spc{U}$ and $p\in \spc{U}$
such that the triangle $\trig{p}{\gamma(0)}{\gamma(T)}$ has perimeter $<2\cdot\varpi\kappa$,
the $\kappa$-development $\tilde \gamma$ in $\Lob2\kappa$ of $\gamma$ with respect to $p$ is concave. 
\end{subthm}

\end{thm}

\parit{Proof; \ref{SHORT.locally-concave-dev}.}  
Let  $\gamma=\geod_{[x y]}$ and $T=\dist{x}{y}{}$. 
Let $\tilde \gamma\:[0,T]\to\Lob2\kappa$ be the concave $\kappa$-development based at $\tilde p$ of $\gamma$ with respect to $p$. 
Let us show that the function  
\[t\mapsto \angk\kappa x p{\gamma(t)}
\eqlbl{eq:ang-nondecreasing}\]   
is nondecreasing. 

For a partition $0=t^0<t^1<\dots<t^n=T$, let 
\[\tilde y^i=\tilde \gamma(t^i)\quad \text{and}\quad \tau^i=\dist{\tilde y^0}{\tilde y^1}{}+\dist{\tilde y^1}{\tilde y^2}{}+\dots+\dist{\tilde y^{i-1}}{\tilde y^i}{}.\]  
Since $\tilde \gamma$ is locally concave, 
for a sufficiently fine partition the polygonal line $\tilde y^0\tilde y^1\dots\tilde y^n$ is  locally concave with respect to $\tilde p$. 
Alexandrov's lemma (\ref{lem:alex}), applied inductively to pairs of triangles  $\modtrig\kappa \{\tau^{i-1},\dist{p}{\tilde y^0}{},\dist{p}{\tilde y^{i-1}}{}\}$ and  $\modtrig\kappa\{\dist{\tilde y^{i-1}}{\tilde y^i}{}, \dist{p}{\tilde y^{i-1}}{},\dist{p}{\tilde y^{i}}{}\}$, shows that the sequence  $\tilde \mangle\mc\kappa\{\dist{\tilde p}{\tilde y^{i}}{};\dist{\tilde p}{\tilde y^0}{},\tau^i\}$ is nondecreasing.

Taking finer partitions and passing to the limit, we get
\[\max\nolimits_i\{|\tau^i-t^i|\}\to0.\] 
Therefore \ref{eq:ang-nondecreasing} and 
the point-on-side comparison (\ref{cat-monoton}) follows.

\parit{\ref{SHORT.concave-dev}.}  
Consider a partition $0=t^0<t^1<\dots<t^n=T$, and 
let $x^i\z=\gamma(t^i)$. Construct a chain of model triangles  $\trig{\tilde p}{\tilde x^{i-1}}{\tilde x^i}=\modtrig\kappa(p x^{i-1}x^i)$ with the direction of $[\tilde p\tilde x^i]$ turning counterclockwise as $i$ grows. 
By the angle comparison (\ref{cat-hinge}),
\[\mangle\hinge{\tilde x^i}{\tilde x^{i-1}}{\tilde p}+\mangle\hinge{\tilde x^i}{\tilde x^{i+1}}{\tilde p}\ge\pi.\eqlbl{eq1:concave-devel*}
\] 
Since $\gamma$ is a geodesic, 
 \[\length \gamma = \sum_{i=1}^n\dist{x^{i-1}}{x^i}{}\le \dist{p}{x^0}{}+\dist{p}{x^n}{}.
\eqlbl{eq2:concave-devel*}
\]  
By repeated application of Alexandrov's lemma (\ref{lem:alex}), and inequality~\ref{eq2:concave-devel*}, 
\[\sum_{i=1}^n\mangle\hinge{\tilde p}{\tilde x^{i-1}}{\tilde x^i}
\le
\angk\kappa p{x^0}{x^n}\le\pi.\] 
Then by \ref{eq1:concave-devel*},  the polygonal line $\tilde p\tilde x^0\tilde x^1\dots \tilde x^n$  are concave with respect to~$\tilde p$.

Note that  under finer partitions, the polygonal line $\tilde x^0\tilde x^1\dots \tilde x^n$ approach the development of $\gamma$ with respect to $p$.
Since the polygonal lines are convex, their lengths converge to the length of $\gamma$.
Hence the result. 
\qeds


\section{Patchwork globalization}\label{sec:patchwork}

If $\spc{U}$ is a $\CAT\kappa$ space, then it is locally $\CAT\kappa$.
The converse does not hold even for complete length space.
For example, $\mathbb{S}^1$ is locally isometric to $\RR$, and so
is locally $\CAT0$, but it is easy to find a quadruple of points in $\mathbb{S}^1$ that violates $\CAT0$ comparison.  

The following theorem was essentially proved by Alexandr Alexandrov \cite[Satz 9]{alexandrov:devel}; 
it gives a global condition on geodesics that is  necessary and sufficient for a locally $\CAT\kappa$ space to be globally $\CAT\kappa$. 
The proof uses thin-triangle decompositions 
and the inheritance lemma (\ref{lem:inherit-angle}). 

\begin{thm}{Patchwork globalization theorem}\label{thm:alex-patch}
For any complete length space~$\spc{U}$, the following two statements are equivalent:

\begin{subthm}{thm:alex-patch:ccat}
$\spc{U}$ is $\CAT\kappa$.
\end{subthm}
 
\begin{subthm}{thm:alex-patch:geo-uni}
$\spc{U}$ is locally $\CAT\kappa$; moreover,  pairs of points in $\spc{U}$ at distance $<\varpi\kappa$ are joined by unique geodesics, and these geodesics depend continuously on their endpoint pairs.
\end{subthm}

\end{thm}

Note that the implication \ref{SHORT.thm:alex-patch:ccat}$\Rightarrow$\ref{SHORT.thm:alex-patch:geo-uni} follows from Theorem~\ref{thm:cat-unique}.

\begin{thm}{Corollary}\label{cor:k-for-k}
Let $\spc{U}$ be a complete length  space 
and $\Omega\subset\spc{U}$ be an open locally $\CAT\kappa$ subset. 
Then for any point $p\in \Omega$ there is $R>0$ such that $\cBall[p,R]$ is a convex subset of $\spc{U}$ 
and $\cBall[p,R]$ is $\CAT\kappa$.
\end{thm}

\parit{Proof.}
Fix $R>0$ such that $\CAT\kappa$ comparison holds in $\oBall(p,R)$.

We may assume that $\oBall(p,R)\subset\Omega$ and $R<\varpi\Kappa$.
The same argument as in the proof of the theorem on uniqueness of geodesics (\ref{thm:cat-unique}) 
shows that any two points in $\cBall[p,\tfrac R2]$ can be joined by a unique geodesic that depends continuously on the endpoints.

The same argument as in the proof of Corollary \ref{cor:convex-balls} shows that $\cBall[p,\tfrac R2]$ is a convex set.
Then \ref{SHORT.thm:alex-patch:geo-uni}$\Rightarrow$\ref{SHORT.thm:alex-patch:ccat} of the patchwork globalization theorem implies that $\cBall[p,\tfrac R2]$ is $\CAT\kappa$.
\qeds

The proof of patchwork globalization uses the following construction:

\begin{thm}{Definition (Line-of-sight map)} \label{def:sight}
Let  $p$ be a point and $\alpha$ be a curve of finite length in  a length space $\spc{U}$. 
Let $\bar\alpha:[0,1]\to\spc{U}$ be the constant-speed parametrization of $\alpha$.
If $\gamma_t\:[0,1]\to\spc{U}$ is a geodesic path from $p$ to $\bar\alpha(t)$, we say that the map $[0,1]\times[0,1]\to\spc{U}$ defined by
\[(t,s)\mapsto\gamma_t(s)\]
is a \index{line-of-sight map}\emph{line-of-sight map} for $\alpha$ with respect to $p$.
\end{thm}

Note that a line-of-sight map is closely related to geodesic homotopy (Section~\ref{sec:Hadamard--Cartan}).

\parit{Proof of \ref{thm:alex-patch}.}
It only remains to prove \ref{SHORT.thm:alex-patch:geo-uni}$\Rightarrow$\ref{SHORT.thm:alex-patch:ccat}.

Let $[p x y]$ be a triangle of perimeter $<2\cdot\varpi\kappa$  in $\spc{U}$. 
According to \ref{prop:k-thin} and \ref{prop:inherit-bound}, it is sufficient to show the triangle $\trig p x y$ is $\kappa$-thin.

Since pairs of points at distance $<\varpi\kappa$ are joined by unique geodesics and these geodesics depend continuously on their endpoint pairs, there is a unique and continuous line-of-sight map for  $[x y]$ with respect to~$p$.    

For a partition \[0\z=t^0\z<t^1\z<\z\dots\z<t^N=1,\] 
let $x^{i,j}=\gamma_{t^i}(t^j)$. 
\begin{figure}[!ht]
\vskip0mm
\centering
\includegraphics{mppics/pic-935}
\end{figure}
Since the line-of-sight map is continuous, we may assume each triangle $\trig{x^{i,j}}{x^{i,j+1}}{x^{i+1,j+1}}$ and $\trig{x^{i,j}}{x^{i+1,j}}{x^{i+1,j+1}}$ is $\kappa$-thin 
(see Proposition~\ref{prop:k-thin}).

Now we show that the $\kappa$-thin property propagates to $\trig p x y$, by repeated application of the inheritance lemma (\ref{lem:inherit-angle}):
\begin{itemize}
\item 
First, for fixed $i$, 
sequentially applying the lemma shows  that the triangles 
$\trig{x}{x^{i,1}}{x^{i+1,2}}$, 
$\trig{x}{x^{i,2}}{x^{i+1,2}}$, 
$\trig{x}{x^{i,2}}{x^{i+1,3}}$,
and so on are $\kappa$-thin. 
\end{itemize}
In particular, for each $i$, the long triangle $\trig{x}{x^{i,N}}{x^{i+1,N}}$ is $\kappa$-thin.
\begin{itemize} 
\item 
Applying the lemma again shows that the  triangles $\trig{x}{x^{0,N}}{x^{2,N}}$, $\trig{x}{x^{0,N}}{x^{3,N}}$, and so on are $\kappa$-thin. 
\end{itemize}
In particular, $\trig p x y=\trig{p}{x^{0,N}}{x^{N,N}}$ is $\kappa$-thin.
\qeds

The following exercise implies that if the space is proper, then one can drop the condition on continuous dependence of geodesics in the formulation of patchwork globalization.

\begin{thm}{Exercise}\label{ex:patchwork}
\begin{subthm}{ex:patchwork:proper}
Suppose pairs of points in a geodesic space $\spc{U}$ are joined by unique geodesics.
Show that if $\spc{U}$ is proper, then 
these geodesics depend continuously on their endpoint pairs.
\end{subthm}

\begin{subthm}{ex:patchwork:complete}
Construct an example of a complete geodesic space $\spc{U}$ such that 
pairs of points in $\spc{U}$ are joined by unique geodesics, but
these geodesics do not depend continuously on their endpoint pairs.
\end{subthm}
\end{thm}


\section{Angles}
\label{sec:angles-cba}

Recall that $\o$ denotes a selective nonprincipal ultrafilter on $\NN$, see Section~\ref{ultralimits}. 

\begin{thm}{Angle semicontinuity}\label{lem:ang.semicont}
Suppose $\spc{U}_1,\spc{U}_2,\dots$ is a sequence of $\varpi\kappa$-geodesic $\CAT\kappa$ spaces
and $\spc{U}_n\to \spc{U}_\o$ as $n\to\o$.
Assume that a sequence of hinges $\hinge{p_n}{x_n}{y_n}$ in $\spc{U}_n$ converges to a hinge $\hinge{p_\o}{x_\o}{y_\o}$ in $\spc{U}_\o$ as $n\to\o$.
Then 
\[\mangle\hinge{p_\o}{x_\o}{y_\o}
\ge 
\lim_{n\to\o} \mangle\hinge{p_n}{x_n}{y_n}.\]

\end{thm}


\parit{Proof.}
By the angle-sidelength monotonicity (\ref{cor:monoton-cba}),
\[\mangle\hinge{p_\o}{x_\o}{y_\o}
=
\inf\set{\angk\kappa{p_\o}{\bar x_\o}{\bar y_\o}}{\bar x_\o \in \mathopen{]}p_\o x_\o],\ \bar y_\o\in \mathopen{]}p_\o y_\o]}.\]

For fixed $\bar x_\o \in \mathopen{]}p_\o x_\o]$ 
and $\bar y_\o\in \mathopen{]}p_\o x_\o]$,
choose $\bar x_n\in \mathopen{]} p x_n ]$ and $\bar y_n\in \mathopen{]} p y_n ]$ so that $\bar x_n\to \bar x_\o$ 
and $\bar y_n\to \bar y_\o$ as $n\to\o$.
Clearly 
\[\angk\kappa{p_n}{\bar x_n}{\bar y_n}
\to 
\angk\kappa{p_\o}{\bar x_\o}{\bar y_\o}\] 
as $n\to\o$.

By the angle comparison (\ref{cat-hinge}), $\mangle\hinge{p_n}{x_n}{y_n}\le \angk\kappa{p_n}{\bar x_n}{\bar y_n}$.
Hence the result.
\qeds

Now we verify that the first variation formula 
holds in the $\CAT{}$ setting. 
Compare it to the first variation inequality (\ref{lem:first-var}) which holds for general metric spaces and to the strong angle lemma (\ref{1st-var+}) for $\Alex{}$ spaces. 

\begin{thm}{Strong angle lemma}
\label{lem:strong-angle-cba}
Let $\spc{U}$ be a $\varpi\kappa$-geodesic $\CAT\kappa$  space.
Then for any hinge  $\hinge  p q y$ in $\spc{U}$, 
we have
\[\mangle\hinge p q y
=
\lim_{\bar y\to p}
\set{\angk\kappa p q{\bar y}}{\bar y\in\mathopen{]}py]}
\eqlbl{eq:cba-1st-var+***}\]
for any $\kappa\in\RR$ such that $\dist{p}{q}{}<\varpi\kappa$.
\end{thm}

\parit{Proof.} 
By angle-sidelength  monotonicity  (\ref{cor:monoton-cba}), the right-hand side is defined and bigger than or  equal to the left-hand side. 

By Lemma~\ref{lem:k-K-angle}, we may take $\kappa = 0$ in \ref{eq:cba-1st-var+***}.  
By the cosine law and the first variation inequality (\ref{lem:first-var}),  
the right-hand side is less than or equal to the left-hand side.
\qeds

\begin{thm}{First variation}\label{thm:1st-var-cba}
Let $\spc{U}$ be a $\varpi\kappa$-geodesic $\CAT\kappa$  space.
For any nontrivial geodesic $[py]$ in $\spc{U}$ and point $q\ne p$ such that  $\dist{p}{q}{}<\varpi\kappa$, we have 
\[\dist{q}{\geod_{[p y]}(t)}{}
=
\dist{q}{p}{}-t\cdot\cos\mangle\hinge p q y+o(t).
\]
\end{thm}

\parit{Proof.}
The first variation equation is equivalent to the strong angle lemma (\ref{lem:strong-angle-cba}), as follows from the cosine law.
\qeds

\begin{thm} {First variation (both-endpoints version)}\label{cor:both-end-first-var-cba}
Assume that $\spc{U}$ is a $\varpi\kappa$-geodesic $\CAT\kappa$ space.
Then for any nontrivial geodesics $[py]$ and $[qz]$ in $\spc{U}$  such that $p\ne q$ and $\dist{p}{q}{}<\varpi\kappa$, we have 
\[
\dist{\geod_{[p y]}(t)}{\geod_{[q z]}(\tau)}{}
=
\dist{q}{p}{} - t\cdot\cos\mangle\hinge p q y - \tau\cdot\cos\mangle\hinge q p z+o(t+\tau).
\]
\end{thm}

\parit{Proof.}
By \ref{cat-hinge},
\[\begin{aligned}
&\dist{\geod_{[p y]}(t)}{\geod_{[q z]}(\tau)}{} \ge
\\
&\ge
\dist{q}{\geod_{[p y]}(t)}{} - \tau\cdot\cos\mangle\hinge q  {\geod_{[p y]}(t)} z +o(\tau)\ge\\
&\ge\dist{q}{p}{} - t\cdot\cos\mangle\hinge p q y + o(t) -  \tau\cdot\cos\mangle\hinge q  {\geod_{[p y]}(t)} z +o(\tau)=\\
&= \dist{q}{p}{} - t\cdot\cos\mangle\hinge p q y -  \tau\cdot\cos\mangle\hinge q  p z +o(t+\tau).
\end{aligned}
\]
Here the final equality follows from   
\[
\lim_{t\to 0}\mangle\hinge q  {\geod_{[p y]}(t)} z = \mangle\hinge q  p z.
\eqlbl{eq:2-side-variation}
\]
The angle semicontinuity (\ref{lem:ang.semicont}) implies ``$\le$'' in \ref{eq:2-side-variation}, and ``$\ge$'' holds by the triangle inequality for angles, since angle comparison (\ref{cat-hinge}) gives 
\[
\lim_{t\to 0}\mangle\hinge q p  {\geod_{[p y]}(t)} = 0.
\]

The opposite inequality follows from \ref{thm:1st-var-cba} and the triangle inequality
\[\dist{\geod_{[p y]}(t)}{\geod_{[q z]}(\tau)}{}
\le
\dist{\geod_{[p y]}(t)}{m}{}+\dist{m}{\geod_{[q z]}(\tau)}{},\]
where $m$ is the midpoint of $[pq]$.
\qeds

We have given elementary proofs of the first-variation statements \ref{lem:strong-angle-cba}, \ref{thm:1st-var-cba} and \ref{cor:both-end-first-var-cba}.
Note however that the no-conjugate-point theorem \ref{thm:no-conj-pt} not only provides proofs of these statements but also extends the statements from geodesics in $\CAT\kappa$ spaces to local geodesics in locally $\CAT\kappa$ spaces as follows:
 
\begin{thm}{First variation for local geodesics}\label{cor:1st-var++cba}
Let $\gamma_t\:[0,1]\to \spc{U}$ be a continuous family of local geodesics in a locally $\CAT\kappa$.
Set
$\alpha(t)=\gamma_t(0)$ and $\beta(t)=\gamma_t(1)$.
Suppose that $\gamma_0$ is unit-speed and $\alpha^+(0)$ and $\beta^+(0)$ are defined.
Then 
\[\length\gamma_t
=
\length\gamma_0
-
(\langle\alpha^+(0),\gamma_0^+(0)\rangle+\langle\beta^+(0),\gamma_0^-(1)\rangle)\cdot t
+
o(t).\]

\end{thm}

\section{Reshetnyak gluing theorem}\label{sec:cba-gluing}

The following theorem was proved by Yuriy Reshetnyak \cite{reshetnyak:major}, assuming $\spc{U}^1$, $\spc{U}^2$ are proper and complete. 
In the following form, the theorem appears in the book of Martin Bridson and Andr\'e Haefliger \cite{bridson-haefliger}.

\begin{thm}{Reshetnyak gluing theorem}\label{thm:gluing}
Suppose 
$\spc{U}^1$, $\spc{U}^2$ are 
$\varpi\kappa$-geodesic spaces 
with isometric complete $\varpi\kappa$-convex sets $A^i\subset\spc{U}^i$.  Let $\iota\:A^1\to A^2$ be an isometry.
Let $\spc{W}=\spc{U}^1\sqcup_{\iota}\spc{U}^2$;
that is, $\spc{W}$ is the gluing of $\spc{U}^1$ and  $\spc{U}^2$ along $\iota$ (see Section~\ref{sec:quotient}).

Then: 
\begin{subthm}{gluing0}
Both canonical mappings $\jmath_i\:\spc{U}^i\to\spc{W}$ are distance-preserving 
and the images $\jmath_i(\spc{U}^i)$ are $\varpi\kappa$-convex subsets in $\spc{W}$.
\end{subthm}

\begin{subthm}{gluing2}
If $\spc{U}^1, \spc{U}^2$ are $\CAT\kappa$,
then so is $\spc{W}$.
\end{subthm} 
\end{thm}

\parit{Proof.} 
Part \ref{SHORT.gluing0}
follows directly from $\varpi\kappa$-convexity of the $A^i$.

\parit{\ref{SHORT.gluing2}.} 
According to \ref{SHORT.gluing0},
we can identify $\spc{U}^i$ with its image $\jmath_i(\spc{U}^i)$ in $\spc{W}$;
in this way, the subsets $A^i\subset \spc{U}^i$ will be identified and denoted further by $A$.
Thus   $A=\spc{U}^1\cap \spc{U}^2\subset \spc{W}$,
and $A$ is $\varpi\kappa$-convex in $\spc{W}$.

Part \ref{SHORT.gluing2} can be reformulated as follows:

\begin{thm}{Reformulation of \ref{gluing2}}\label{thm:gluing2-reformulated}
Let $\spc{W}$ be a 
length space having two 
$\varpi\kappa$-convex subsets $\spc{U}^1,\spc{U}^2\subset\spc{W}$ such that
$\spc{W}=\spc{U}^1\cup\spc{U}^2$.
Assume the subset $A=\spc{U}^1\cap \spc{U}^2$ is complete and $\varpi\kappa$-convex in $\spc{W}$, and $\spc{U}^1$, $\spc{U}^2$ are $\CAT\kappa$ spaces.
Then $\spc{W}$ is a $\CAT\kappa$ space.
\end{thm}

\begin{clm}{}\label{clm:geod-gluing}
If $\spc{W}$ is $\varpi\kappa$-geodesic, then $\spc{W}$ is $\CAT\kappa$.
\end{clm}

Indeed, 
according to \ref{prop:k-thin},
it is sufficient to show that any triangle $\trig {x^0}{x^1}{x^2}$ of perimeter $<2\cdot \varpi\kappa$ 
in $\spc{W}$ is $\kappa$-thin.
This is obviously true if all three points $x^0$, $x^1$, $x^2$ lie in a single $\spc{U}^i$.
Thus, without loss of generality, we may assume that $x^0\in\spc{U}^1$ and $x^1,x^2\in\spc{U}^2$.
\begin{figure}[!ht]
\vskip-0mm
\centering
\includegraphics{mppics/pic-940}
\end{figure}

Choose points $z^1,z^2\in A=\spc{U}^1\cap\spc{U}^2$ 
lying respectively on the sides $[x^0x^1], [x^0x^2]$.
Note that all distances between any pair of points from $x^0$, $x^1$, $x^2$, $z^1$, $z^2$ are less than $\varpi\kappa$.
Therefore
\begin{itemize}
\item triangle $\trig{x^0}{z^1}{z^2}$ lies in $\spc{U}^1$,
\item both triangles $\trig{x^1}{z^1}{z^2}$ and $\trig{x^1}{z^2}{x^2}$ lie in $\spc{U}^2$.
\end{itemize}
In particular, each triangle $\trig{x^0}{z^1}{z^2}$,
$\trig{x^1}{z^1}{z^2}$, $\trig{x^1}{z^2}{x^2}$ is $\kappa$-thin.

Applying the inheritance lemma for thin triangles (\ref{lem:inherit-angle}) twice, 
we get that $\trig {x^0}{x^1}{z^2}$ 
and consequently $\trig {x^0}{x^1}{x^2}$ is $\kappa$-thin.
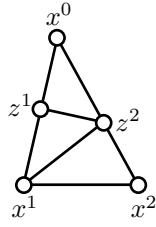
\claimqeds

\begin{clm}{}\label{clm:geod-gluing0 }
$\spc{W}$ is $\CAT\kappa$ if $\kappa\le0$.
\end{clm}
By \ref{clm:geod-gluing} it suffices to prove that $\spc{W}$ is geodesic.

For $p^1\in \spc{U}^1$, $p^2\in \spc{U}^2$, we may choose a sequence $z_n\in A$ such that $\dist{p^1}{z_n}{}+\dist{p^2}{z_n}{}$
 converges to $\dist{p^1}{p^2}{}$, and $\dist{p^1}{z_n}{}$ and $\dist{p^2}{z_n}{}$ converge.  
 Since $A$ is complete, it suffices to show $z_n$ is a Cauchy sequence.  
 In that case, the limit point $z$ of $z_n$ satisfies $\dist{p^1}{z}{}+\dist{p^2}{z}{}=\dist{p^1}{p^2}{}$, so the geodesics $[p^1z]$ in $\spc{U}^1$ and $[p^2z]$ in $\spc{U}^2$ together give a geodesic $[p^1p^2]$ in $\spc{U}$.  
 
 Suppose $z_n$ is not a Cauchy sequence.
 Then there are subsequences  $x_n$ and $y_n$ of $z_n$ satisfying  $\lim\dist{x_n}{y_n}{}>0$.
 Let $m_n$ be the midpoint of $[x_ny_n]$.
 Since $\dist{p^1}{m_n}{}+\dist{p^2}{m_n}{} \ge \dist{p^1}{p^2}{}$, and  $\dist{p^1}{x_n}{}+\dist{p^2}{x_n}{}$ and  $\dist{p^1}{y_n}{}+\dist{p^2}{y_n}{}$
 converge to $\dist{p^1}{p^2}{}$, then for  any $\epsilon >0$, we may assume (taking subsequences and possibly relabeling $p^1$ and $p^2$)
 \[
 \dist{p^1}{m_n}{}
 \ge
 \dist{p^1}{x_n}{}-\epsilon,
 \qquad
 \dist{p^1}{m_n}{}
 \ge
 \dist{p^1}{y_n}{}-\epsilon.
 \]
 
 Since triangle $\trig{p^1}{x_n}{y_n}$ is thin, the analogous inequalities hold for the Euclidean model triangle  $\trig{\tilde p^1}{\tilde x_n}{\tilde y_n}$.  
 Then there is a nondegenerate limit triangle $\trig{p}{x}{y}$ in the Euclidean plane satisfying $\dist{p}{x}{}=\dist{p}{y}{}\le\dist{p}{m}{}$ where $m$ is the midpoint of $[xy]$.  This  contradiction proves the claim.
\claimqeds

Finally suppose $\kappa>0$; by rescaling, take $\kappa=1$. Consider the Euclidean cones $\Cone\spc{U}^i$ (see Section \ref{sec: tangent space}).
By Theorem \ref{thm:warp-curv-bound:cbb:S}, $\Cone\spc{U}^i$ is a $\CAT0$ space for $i=1,2$.

Geodesics contained in the complement of the tip of $\Cone\spc{U}^i$ project to geodesics of length $<\pi$ in $\spc{U}^i$. 
It follows that $\Cone A$ is convex in $\Cone\spc{U}^1$ and $\Cone\spc{U}^2$.
By the cone distance formula, $\Cone A$ is complete since $A$ is complete.
  
Gluing along $\Cone A$ and applying \ref{clm:geod-gluing} and \ref{clm:geod-gluing0 } for $\kappa=0$, we find that 
$\Cone\spc{W}$ is a $\CAT0$ space.  By Theorem \ref{thm:warp-curv-bound:cbb:S}, $\spc{W}$ is a $\CAT1$ space.
\qeds

\begin{thm}{Exercise}\label{ex:two-rays}
Let $Q$ be the nonconvex subset of the plane bounded by two half-lines $\gamma_1$ and $\gamma_2$ with a common starting point and angle $\alpha$ between them.
Assume $\spc{U}$ is a complete length $\CAT0$ space and $\gamma_1',\gamma_2'$ are two half-lines in $\spc{U}$ with a common
starting point and angle $\alpha$ between them.
Show that the space glued from $Q$ and $\spc{U}$ along the corresponding half-lines is a $\CAT{0}$ space.
\end{thm}

\begin{thm}{Exercise}\label{ex:reshetnyak-doubling}
Suppose $\spc{U}$ is a complete length $\CAT0$ space and $A\subset \spc{U}$ is a closed subset.
Assume that the doubling of $\spc{U}$ in $A$ is $\CAT0$. 
Show that $A$ is a convex set of $\spc{U}$.
\end{thm}

\begin{thm}{Exercise}\label{ex:glue-spherical-suspension}
Let  $\spc{U}$ be a complete length $\CAT1$ space and $K\subset \spc{U}$ be a closed $\pi$-convex set.
Assume $K\subset \cBall[p,\tfrac\pi2]$ for $p\in K$.
Show that there is a decreasing continuous one-parameter family of closed convex sets $K_t$ for $t\in[0,1]$ such that $K_0=\cBall[p,\tfrac\pi2]$ and $K_1=K$.

(Decreasing means with respect to inclusion; that is $K_{t_0}\supset K_{t_1}$ if $t_0\le t_1$.
Continuous means with respect to Hausdorff distance; that is $K_t\Hto K_{t_0}$ as $t\to t_0$.)
\end{thm}

\begin{thm}{Exercise}\label{ex:AUB}
Let $A$ and $B$ be two closed convex sets in a complete length $\CAT0$ space.
Assume $A\cap B\ne\emptyset$.
Show that the union $A\cup B$ equipped with induced length metric is $\CAT0$. 
\end{thm}


\section{Space of geodesics}\label{sec:geod-space}

In this section we prove a no-conjugate-point theorem for spaces with upper curvature bounds and derive from it a number of statements about local geodesics.
These statements will be used to prove the Hadamard--Cartan theorem (\ref{thm:hadamard-cartan}) and the lifting globalization theorem (\ref{thm:globalization-lift}), in much the same way as  the exponential map is used in Riemannian geometry.

\begin{thm}{Proposition}\label{prop:geo-complete}
{\sloppy 
Let $\spc{U}$ be a locally $\CAT\kappa$ space.
 Let $\gamma_n\:[0,1]\to\spc{U}$ be a sequence of local geodesic paths converging to a path $\gamma_\infty\:[0,1]\to\spc{U}$.
Then $\gamma_\infty$ is a local geodesic path.
Moreover 
\[\length\gamma_n\to\length\gamma_\infty\]
as $n\to\infty$.

}
\end{thm}

\parit{Proof.} 
Fix $t\in[0,1]$.
By Corollary~\ref{cor:k-for-k}, we may choose $R$ satisfying $0<R<\varpi\Kappa$,
and such that
the ball $\spc{B}=\oBall(\gamma_\infty(t),R)$ is a convex subset of $\spc{U}$ and forms a $\CAT\kappa$ space.

A local geodesic segment  with length less than $R/2$ that intersects $\oBall(\gamma_\infty(t),R/2)$ cannot leave $\spc{B}$, and hence  is  minimizing by Corollary~\ref{cor:loc-geod-are-min}.
In particular, for all sufficiently large $n$, 
if subsegment of $\gamma_n$ has length less than $R/2$ and contains $\gamma_n(t)$, then it is a geodesic.

Since $\spc{B}$ is $\CAT\kappa$, geodesic segments in $\spc{B}$ depend uniquely and continuously on their endpoint pairs by Theorem~\ref{thm:cat-unique}.  
Thus there is a subinterval $\II$ of $[0,1]$
that contains a neighborhood of $t$ in $[0,1]$
and such that $\gamma_n|_\II$ is minimizing for all large $n$.
It follows that the restriction $\gamma_\infty|_\II$ is a geodesic,
and therefore $\gamma_\infty$ is a local geodesic.
\qeds

The following theorem was proved by the first author and Richard Bishop \cite{alexander-bishop:h-c}.
In analogy with Riemannian geometry, the main statement of the following theorem could be restated as: 
\textit{In a space of curvature $\le\kappa$, two points cannot be conjugate along a local geodesic of length $<\varpi\kappa$.}

\begin{thm}{No-conjugate-point theorem}
\label{thm:no-conj-pt}{\sloppy 
Suppose $\spc{U}$ is a locally complete, length, locally $\CAT\kappa$ space.
Let $\gamma\:[0,1]\to\spc{U}$ be a local geodesic path with length $<\varpi\kappa$.
Then for some neighborhoods $\Omega^0\ni \gamma(0)$ and $\Omega^1\ni\gamma(1)$, 
there is a unique continuous map from the direct product $\Omega^0\times \Omega^1\times[0,1]$ to $\spc{U}$, 

\[(x,y,t)\mapsto\gamma_{x y}(t),\]  
such that 
$\gamma_{x y}\:[0,1]\to\spc{U}$ is a local geodesic path with 
$\gamma_{x y}(0)=x$ and 
$\gamma_{x y}(1)=y$ for each $(x,y)\in\Omega^0\times\Omega^1$,
and the family $\gamma_{x y}$ contains $\gamma$.
Moreover, we can assume that the map 
\[(x,y,t)\mapsto\gamma_{x y}(t)\:\Omega^0\times\Omega^1\times[0,1]\to\spc{U}\] 
is $\Lip$-Lipschitz
for any
$\Lip>\max\set{\tfrac{\sn\kappa r}{\sn\kappa \ell}}{0\le r\le \ell}$.

}
\end{thm}

{\sloppy

\begin{thm}{Patchwork along a geodesic}
\label{lem:patch}
Let $\spc{U}$ be a locally complete, length, locally $\CAT\kappa$ space, 
and $\alpha\:[a,b]\to\spc{U}$ be a local geodesic.

Then there is a complete length $\CAT\kappa$ space $\spc{N}$
with an open set $\hat\Omega\subset \spc{N}$,
a local geodesic $\hat\alpha\:[a,b]\to\hat\Omega$,
and an open locally distance-preserving map 
$\map\:\hat\Omega\looparrowright\spc{U}$ such that
$\map\circ\hat\alpha=\alpha$.

Moreover if $\alpha$ is simple, then one can assume in addition that $\map$ is an open embedding;
thus $\hat\Omega$ is locally isometric to a neighborhood of $\Omega=\map(\hat\Omega)$ of $\alpha$.
\end{thm}

}

This lemma and its proof were suggested by Alexander Lytchak.
The proof proceeds by piecing together $\CAT{\kappa}$  neighborhoods of points on a curve to construct a new $\CAT{\kappa}$ space.  
Exercise~\ref{ex:cats-cradle} is inspired by the original idea of the proof of the no-conjugate-point theorem (\ref{thm:no-conj-pt}) given in \cite{alexander-bishop:h-c}.

\parit{Proof.} 
According to Corollary~\ref{cor:k-for-k},
we can choose $r>0$ such that 
for any $t\in[a,b]$ the closed ball
$\cBall[\alpha(t),r]$ is a convex set that forms a complete length $\CAT\kappa$ space.

Choose balls $\spc{B}_i=\cBall[\alpha(t_i),r]$
for some partition $a\z=t_0<t_1\z<\z\dots\z<t_n\z=b$
in such a way that 
\[\Int\spc{B}_i\supset \alpha([t_{i-1},t_i])\]
for all $i>0$.
We can assume in addition that $\spc{B}_{i-1}\cap \spc{B}_{i+1}\subset \spc{B}_{i}$ if $0<i<n$.

Consider the disjoint union $\bigsqcup_i\spc{B}_i=\set{(i,x)}{x\in\spc{B}_i}$ with the minimal equivalence relation $\sim$ such that $(i,x)\sim(i-1,x)$ for all $i>0$.
Let $\spc{N}$ be the space obtained by gluing the $\spc{B}_i$ along $\sim$.
Note that $A_i=\spc{B}_i\cap\spc{B}_{i-1}$ is convex in $\spc{B}_i$ and in $\spc{B}_{i-1}$.
Applying the Reshetnyak gluing theorem (\ref{thm:gluing}) several times, 
we conclude that $\spc{N}$ is a complete length $\CAT\kappa$ space.

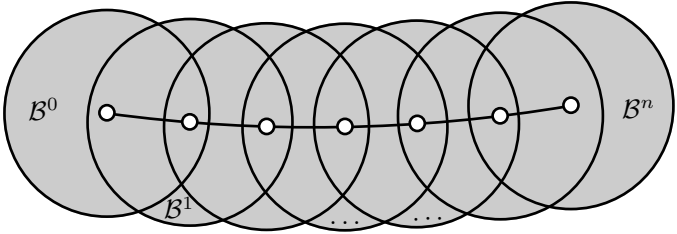
\begin{figure}[!ht]
\vskip-0mm
\centering
\includegraphics{mppics/pic-945}
\end{figure}

For $t\in[t_{i-1},t_i]$, let $\hat\alpha(t)$  be the equivalence class of $(i,\alpha(t))$ in $\spc{N}$.
Let $\hat\Omega$ be the $\eps$-neighborhood of $\hat\alpha$ in $\spc{N}$, where $\eps>0$ is chosen so that $\oBall(\alpha(t),\eps)\subset\spc{B}_i$ for all $t\in[t_{i-1},t_i]$.

Define $\map\:\hat\Omega\to\spc{U}$
by sending the equivalence class of $(i,x)$ to $x$.
It is straightforward to check that $\map\:\spc{N}\to\spc{U}$, $\hat\alpha\:[a,b]\to\spc{N}$ and $\hat\Omega\subset\spc{N}$ satisfy the conclusion of the main part of the lemma.

To prove the final statement in the lemma,
we only have to choose $\eps>0$ so that in addition, $\dist{\alpha(\tau)}{\alpha(\tau')}{}>2\cdot\eps$ if $\tau\le t_{i-1}$ and $t_i\le\tau'$ for some $i$.
\qeds

\parit{Proof of \ref{thm:no-conj-pt}.}
Apply patchwork along $\gamma$ (\ref{lem:patch}). 
\qeds

The No-conjugate-point theorem (\ref{thm:no-conj-pt}) allows us to move a local geodesic  
so that its endpoints follow given trajectories.
The following corollary describes how this process might terminate. 

\begin{thm}{Corollary}\label{cor:geo-hom}{\sloppy 
Let $\spc{U}$ be a locally complete, length, locally $\CAT\kappa$ space.
Suppose $\gamma\:[0,1]\to\spc{U}$ is a local geodesic with length $< \varpi\kappa$.  Let $\alpha^i\:[0,1]\to \spc{U}$, for $i=0,1$, be curves starting at $\gamma(0)$ and $\gamma(1)$ respectively.  

}

Then there is a uniquely determined pair consisting of an interval $\II $ satisfying $0\in \II\subset[0,1]$, and a continuous family of local geodesics $\gamma_t\:[0,1]\to \spc{U}$ for  $t\in \II$, such that  

\begin{subthm}{cor:geo-hom-length}
$\gamma_0=\gamma$, $\gamma_t(0)=\alpha^0(t)$, $\gamma_t(1)=\alpha^1(t)$, and $\gamma_t$ has length $< \varpi\kappa$,
\end{subthm} 

\begin{subthm}{cor:geo-hom-cauchy}
if $\II\ne [0,1]$, then $\II=[0, a)$, where either $\gamma_t$ converges uniformly to a local geodesic $\gamma_a$ of length $\varpi\kappa$, or 
for some fixed $s\in [0,1]$ the curve $\gamma_t(s):[0,a)\to\spc{U}$ is a Lipschitz curve with no limit 
as $t\to a-$.
\end{subthm}

\end{thm}

\parit{Proof.} Uniqueness follows from  Theorem \ref{thm:no-conj-pt}.

Let $\II$ be the maximal interval for which there is a family $\gamma_t$ satisfying condition \ref{SHORT.cor:geo-hom-length}. 
By Theorem~\ref{thm:no-conj-pt}, such an interval exists and is open in $[0,1]$.  Suppose $\II\ne[0,1]$.
Then  $\II=[0,a)$ for $0<a\le 1$.

For each fixed $s\in [0,1]$, define the curve $\alpha_s:[0,a)\to\spc{U}$ by $\alpha_s(t)\z=\gamma_t(s)$. 
By Theorem~\ref{thm:no-conj-pt}, 
each $\alpha_s$ is locally Lipschitz.  

If $\alpha_s$ for some value of $s$ does not converge as $t\to a-$, then condition \ref{SHORT.cor:geo-hom-cauchy} holds.
If each $\alpha_s$  converges as $t\to a-$, then $ \gamma_t$ converges as $t\to a-$, say to $\gamma_a$.
By  Proposition~\ref{prop:geo-complete}, $\gamma_a$ is a local geodesic and\[\length\gamma_t\to\length\gamma_a\le \varpi\kappa.\]
By maximality of $\II$, $\length\gamma_a=\varpi\kappa$ and so condition \ref{SHORT.cor:geo-hom-cauchy} again holds.
\qeds

\begin{thm}{Corollary}\label{cor:homotopy-from-p}
Let $\spc{U}$ be a complete locally $\CAT\kappa$ length space, and 
$\alpha\:[0,1]\to \spc{U}$ be a path of length $< \varpi\kappa$ that starts at $p$ and ends at $q$.
Then:  

\begin{subthm}{cor:homotopy-from-p-exist}
There is a unique homotopy of local geodesic paths $\gamma_t\:[0,1]\to \spc{U}$
such that $\gamma_0(t)=\gamma_t(0)=p$ and $\gamma_t(1)=\alpha(t)$ for any~$t$.
\end{subthm}

\begin{subthm}{cor:homotopy-from-p-length}
For any $t\in[0,1]$, 
\[\length\gamma_t\le\length(\alpha|_{[0,t]}),\]
and equality holds for given $t$ if and only if the restriction $\alpha|_{[0,t]}$ is a reparametrization of $\gamma_t$.
\end{subthm}

Moreover, instead of completeness of $\spc{U}$, one can assume that the subspace 
\[W=\set{x\in \spc{U}}{\dist{x}{p}{}+\dist{x}{q}{}\le \ell}\] 
is complete.

\end{thm}

\parit{Proof.}
By Corollary \ref{cor:geo-hom}, taking  $\alpha^0(t)=p$ and  $\alpha^1(t)=\alpha(t)$ for all $t\in [0,1]$, there is an interval $\II$ such that \ref{SHORT.cor:homotopy-from-p-exist} holds for all $t\in\II$, and either $\II=[0,1]$ or $\II=[0,a)$ for $a\le 1$.

By patchwork along a curve (\ref{lem:patch}), the values of $t$ for which condition \ref{SHORT.cor:homotopy-from-p-length} holds form an open subset of $\II$ containing $0$; clearly this subset is also closed in $\II$.
Therefore \ref{SHORT.cor:homotopy-from-p-length} holds on all of $\II$. 
 
Corollary \ref{cor:geo-hom} implies that
$\II=[0,1]$.
Indeed if $\II=[0,a)$, then either $\length\gamma_t\to\varpi\kappa$ as $t\to a-$,
or for some fixed $s\in [0,1]$ the Lipschitz curve $\gamma_t(s)$ has no limit as $t\to a-$.
Since $\length\alpha<\varpi\kappa$, \ref{cor:geo-hom} implies that neither of these is possible.
\qeds


\section{Lifting globalization}\label{sec:cat-globalize}

The Hadamard--Cartan theorem (\ref{thm:hadamard-cartan}) states that 
the universal metric cover of a complete locally $\CAT0$ space is $\CAT0$.
The lifting globalization theorem gives an appropriate generalization of the above statement to arbitrary curvature bounds;
it could be also described as a global version of Gauss's lemma.

\begin{thm}{Lifting globalization theorem}
\label{thm:globalization-lift}
Suppose $\spc{U}$ is a complete length locally $\CAT\kappa$ space and  $p\in\spc{U}$.
Then there is a complete $\CAT\kappa$ length space $\spc{B}$, 
with a point $\hat p$ such that 
there is a locally distance-preserving map $\map\:\spc{B}\to\spc{U}$
such that $\map(\hat p)=p$ and the following lifting property holds: 
for any path $\alpha\:[0,1]\to\spc{U}$ with $\alpha(0)=p$ and $\length\alpha<\varpi\kappa/2$, 
there is a unique path $\hat\alpha \:[0,1]\to \spc{B}$ such that $\hat\alpha(0)=\hat p$ 
and $\map\circ\hat\alpha\equiv\alpha$.
\end{thm}

Note that the lifting property implies that $\map(\spc{B})\supset\oBall(p,\varpi\kappa/2)$ and by completeness $\map(\spc{B})\supset\cBall[p,\varpi\kappa/2]$.
Also, since $\spc{B}$ is $\CAT\kappa$, the closed ball $\cBall[\hat p,\tfrac{\varpi\kappa}2]_{\spc{B}}$ is a weakly convex set in $\spc{B}$ (see \ref{cor:convex-balls});
in particular $\cBall[\hat p,\tfrac{\varpi\kappa}2]_{\spc{B}}$ is a complete length $\CAT\kappa$ space.
Therefore we can assume in addition that $\dist{\hat p}{\hat x}{}\le \varpi\kappa/2$ for any $\hat x\in\spc{B}$;
or equivalently
\[\cBall[\hat p,\tfrac{\varpi\kappa}2]_{\spc{B}}=\spc{B}.\]

Before proving the theorem we state and prove its corollary.

\begin{thm}{Corollary}\label{cor:loc-CAT(k)}
Suppose $\spc{U}$ is a complete length locally $\CAT\kappa$ space.
Then for any $p\in\spc{U}$ there is $\rho_p>0$
such that $\cBall[p,\rho_p]$ is a complete length $\CAT\kappa$ space.

Moreover, we can assume that $\rho_p<\tfrac{\varpi\kappa}2$
for any $p$ and the function $p\mapsto\rho_p$ is 1-Lipschitz.
\end{thm}

\parit{Proof.} 
Assume $\map\:\spc{B}\to \spc{U}$ 
and $\hat p\in \spc{B}$
are provided by the lifting globalization theorem
(\ref{thm:globalization-lift}).

Since $\map$ is local isometry,
we can choose $r>0$ so that the restriction of $\map$ to $\cBall[\hat p,r]$ is distance-preserving.
By the lifting globalization, the image  $\Phi(\cBall[\hat p,r])$ coincides with the ball
$\cBall[p,r]$.
This proves the first part of the theorem.

To prove the second part, let us choose $\rho_p$ to be the maximal value $\le\tfrac{\varpi\kappa}2$ such that $\cBall[p,\rho_p]$ is a complete length $\CAT\kappa$ space.
By Corollary~\ref{cor:convex-balls}, the ball
\[\cBall[q,\rho_p-\dist{p}{q}{}]\] 
is weakly convex in $\cBall[p,\rho_p]$.
Therefore  
\[\cBall[q,\rho_p-\dist{p}{q}{}]\] is a complete length $\CAT\kappa$ space
for any $q\in \oBall(p,\rho_p)$.
In particular, $\rho_q\ge \rho_p-\dist{p}{q}{}$ for any $p,q\in\spc{U}$.
Hence the second statement follows.
\qeds

The proof of the lifting globalization theorem relies heavily on the properties of the space of local geodesic paths discussed in Section~\ref{sec:geod-space}.
The following lemma  is a key step in the proof;
it was proved by the first author and Richard Bishop \cite{alexander-bishop:cbc}. 

\begin{thm}{Radial lemma}\label{lem:radial-glob}
Let $\spc{U}$ be a length locally $\CAT\kappa$ space,
and suppose $p\in\spc{U}$, $R\le\varpi\kappa$.
Assume the ball  $\cBall[p,\bar{R}]$ is complete for any $\bar{R}<R$, and  there is a unique geodesic path, $\geodpath_{[p x]}$, from $p$ to any point $x\in\oBall(p,R)$ 
that depends continuously on $x$.
Then $\oBall(p,\tfrac R2)$ is a $\varpi\kappa$-geodesic $\CAT\kappa$ space.
\end{thm}
 
\parit{Proof.}
Without loss of generality, we may assume  $\spc{U}=\oBall(p,R)$.

Set $f=\md\kappa\circ\distfun{p}{}{}$.  Let us show that
\[f''+\kappa\cdot f\ge 1.
\eqlbl{eq:rad-conv}\]

Fix $z\in \spc{U}$.
We will apply the no-conjugate-point theorem (\ref{thm:no-conj-pt}) for the unique geodesic path $\gamma$
from $p$ to $z$.  
The  notations $\Omega^0$, 
$\Omega^1$,
$\gamma_{x y}$, $\spc{N}$, $\hat{x}$, $\hat{y}$ will be as in the formulation of the lifting globalization theorem (\ref{thm:globalization-lift});
in particular, $z\in\Omega^1$.

By assumption,
$\gamma_{p y}=\geodpath_{[p y]}$ for any $y\in\Omega^1$. 
Consequently,
 $f(y)\z=\md\kappa\dist[{{}}]{\hat{p}}{\hat{y}}{\spc{N}}$.
Applying the function comparison (\ref{thm:function-comp}) in $\spc{N}$,
we have that $f''+\kappa\cdot f\ge 1$ in $\Omega^1$;
whence \ref{eq:rad-conv} follows.
\claimqeds

Fix $r<\tfrac R2$. Proving the following claim takes most of the remaining proof:

\begin{clm}{}\label{clm:B-is-convex}
$\cBall[ p,r]$ is a convex set in $\spc{U}$.
\end{clm}

Choose arbitrary $x,z\in \cBall[ p,r]$.
First note that \ref{eq:rad-conv} implies the following claim.

\begin{clm}{}\label{clm:B-is-almost-convex}
If $ \gamma\:[0,1]\to\spc{U}$ 
is a local geodesic path from $x$ to $z$ and  
$\length \gamma\z<\varpi\kappa$,  
then $\length \gamma \le 2\cdot r$ 
and $ \gamma$ lies completely in $\cBall[ p,r]$.
\end{clm}

Note that  $\dist{x}{z}{}<\varpi\kappa$.
Thus to prove Claim~\ref{clm:B-is-convex}, it is sufficient to show that there is a geodesic path from $x$ to $z$.
Note that by assumption $\cBall[p,2\cdot r]$ is complete.
Therefore Corollary~\ref{cor:homotopy-from-p} implies the following:

\begin{clm}{}\label{clm:loc-geod<path}
Given a path $\alpha\:[0,1]\to\spc{U}$ from $x$ to $z$ with $\length\alpha<2\cdot r$,
there is a local geodesic path $\gamma$ from $x$ to $z$ such that
\[\length\gamma\le\length\alpha.\]

\end{clm}

Further, let us prove the following:

\begin{clm}{}\label{clm:unique-loc-geod}
There is a unique local geodesic path $\gamma_{x z}$ in $\cBall[ p,r]$ from $x$ to $z$.
\end{clm}

Denote by $\Delta_{x z}$ the set of all local geodesic paths in $\cBall[ p,r]$ from $x$ to $z$.
By Corollary \ref{cor:geo-hom}, there is a  bijection $\Delta_{x z}\to\Delta_{p p}$.
According to \ref{eq:rad-conv}, 
$\Delta_{p p}$ contains only the constant path.
Claim~\ref{clm:unique-loc-geod} follows.

Note that 
claims~\ref{clm:B-is-almost-convex}, 
\ref{clm:loc-geod<path} 
and \ref{clm:unique-loc-geod}
imply that $\gamma_{x z}$ is minimizing; hence Claim~\ref{clm:B-is-convex}.

Further, Claim~\ref{clm:B-is-almost-convex} and the no-conjugate-point theorem (\ref{thm:no-conj-pt}) together 
imply that the map $(x,z)\mapsto\gamma_{x z}$ is continuous.

Therefore by the patchwork globalization theorem (\ref{thm:alex-patch}), 
$\cBall[ p,r]$ is a $\varpi\kappa$-geodesic $\CAT\kappa$ space.

Since
\[\oBall( p,R)
=
\bigcup_{r < R}\cBall[ p,r],\] 
then $\oBall( p,R)$ is convex in $\spc{U}$ and 
$\CAT\kappa$ comparison holds  for any quadruple in $\oBall( p,R)$.
Therefore $\oBall( p,\varpi\kappa/2)$ is $\CAT\kappa$.
\qeds

In the following proof, we construct a space $\mathfrak{G}_p$ of  local geodesic paths that start at $p$.
The space $\mathfrak{G}_p$ comes with 
a marked point $\hat p$ 
and the endpoint map $\map\:\mathfrak{G}_p\to\spc{U}$.
One can think of
the map $\map$ as an analog of the exponential map $\exp_p$ in the Riemannian geometry;
in this case,
the space $\mathfrak{G}_p$ corresponds to the ball of radius $\varpi\kappa$ in the tangent space at $p$, equipped with the metric pulled back by $\exp_p$.

We are going to set $\spc{B}=\oBall(\hat p,\varpi\kappa/2)\subset \mathfrak{G}_p$,
and use the radial lemma (\ref{lem:radial-glob}) to prove that $\spc{B}$ is a $\varpi\kappa$-geodesic $\CAT\kappa$ space.

\parit{Proof of \ref{thm:globalization-lift}.}
Suppose $\hat\gamma$ is a homotopy of local geodesic paths that start at $p$.  Thus the map 
\[\hat\gamma\:(t,\tau)\mapsto\hat\gamma_t(\tau)\:[0,1]\times[0,1]\to\spc{U}\] 
is continuous,
and the following holds for each $t$:
\begin{itemize}
\item $\hat\gamma_t(0)=p$,
\item $\hat\gamma_t\:[0,1]\to\spc{U}$ is a local geodesic path in $\spc{U}$.
\end{itemize}

Denote by $\theta(\hat\gamma)$ the length traced by the ends of $\hat\gamma_t$;
that is, $\theta(\hat\gamma)$ is the length of the path $t\mapsto\hat\gamma_t(1)$.

Let $\mathfrak{G}_p$ be the set of all local geodesic paths 
with length $<\varpi\kappa$ in $\spc{U}$ that start at $p$.
Denote by $\hat p\in \mathfrak{G}_p$ the constant path $\hat p(t)\equiv p$.
Given $\alpha,\beta\in \mathfrak{G}_p$, define
\[
\dist{\alpha}{\beta}{\mathfrak{G}_p}
=
\inf_{\hat\gamma} \{\theta(\hat\gamma)\},\]
with the exact lower bound taken along all homotopies 
$\hat\gamma\:[0,1]\times[0,1]\to\spc{U}$ 
such that 
$\hat\gamma_0=\alpha$, 
$\hat\gamma_1=\beta$ 
and $\hat\gamma_t\in \mathfrak{G}_p$ for all $t\in[0,1]$.

By the no-conjugate-point theorem (\ref{thm:no-conj-pt}), we have $\dist{\alpha}{\beta}{\mathfrak{G}_p}>0$ for distinct $\alpha$ and $\beta$;
that is,

\begin{clm}{}
$\dist{{*}}{{*}}{\mathfrak{G}_p}$ is a metric on $\mathfrak{G}_p$.
\end{clm}

Further, again from the no-conjugate-point theorem (\ref{thm:no-conj-pt}), we have

\begin{clm}{}\label{clm:loc-iso}
The map
\[\map\:\xi\mapsto\xi(1)\:\mathfrak{G}_p\to\spc{U}\]
is a local isometry.
In particular, $\mathfrak{G}_p$ is locally $\CAT\kappa$.
\end{clm}

Let $\alpha\:[0,1]\to\spc{U}$ be a path, $\length\alpha<\varpi\kappa$, and $\alpha(0)=p$.
The homotopy constructed in Corollary~\ref{cor:homotopy-from-p} can be regarded as a path in $\mathfrak{G}_p$, say $\hat\alpha\:[0,1]\to \mathfrak{G}_p$,
such that $\hat\alpha(0)=\hat p$ and $\map\circ\hat\alpha=\alpha$;
in particular $\hat\alpha_t(1)\equiv\alpha(t)$ for any $t$. 
By \ref{clm:loc-iso}, 
\[\length(\hat\alpha)_{\mathfrak{G}_p}=\length(\alpha)_{\spc{U}}.\]
Moreover, it follows that $\alpha$ is a local geodesic path of $\spc{U}$  if and only if $\hat\alpha$ is a local geodesic path of $\mathfrak{G}_p$.

Further, from Corollary~\ref{cor:homotopy-from-p},
for any $\xi\in \mathfrak{G}_p$ and path $\hat\alpha\:[0,1]\to\mathfrak{G}_p$ from $\hat p$ to $\xi$,
we have 
\begin{align*}
\length\hat\alpha
&=\length\map\circ\hat\alpha
\ge
\\
&\ge
\length\xi
=
\\
&=\length\hat\xi
\end{align*}
where equality holds only if $\hat\alpha$ is a reparametrization of $\hat\xi$.
In particular, 
\[\dist{\hat p}{\xi}{\mathfrak{G}_p}=\length\xi
\eqlbl{eq:dist=length}\] 
and
$\hat\xi\:[0,1]\to \mathfrak{G}_p$ is the unique geodesic path from $\hat p$ to $\xi$.
Clearly, the map $\xi\mapsto\hat\xi$ is continuous.

By \ref{eq:dist=length} and Proposition~\ref{prop:geo-complete},

\begin{clm}{}\label{clm:complete-B} 
For any $\bar R<\varpi\kappa$, the closed ball
$\cBall[\hat p,\bar R]$ in $\mathfrak{G}_p$ is complete.
\end{clm}

Take $\oBall(\hat p,\varpi\kappa/2)$ and $\map$ constructed above.
According to the radial lemma (\ref{lem:radial-glob}), $\oBall(\hat p,\varpi\kappa/2)$ is a $\varpi\kappa$-geodesic $\CAT\kappa$ space.
The map $\map$ extends to its completion $\spc{B}=\cBall[\hat p,\varpi\kappa/2]$. 
All the remaining statements are already proved.
\qeds

\section{Reshetnyak majorization}\label{sec:resh-kirz}

\begin{thm}{Definition}\label{def:majorize}
Let $\spc{X}$ be a metric space,
$\tilde \alpha$ be a simple closed curve of finite length  in $\Lob2{\kappa}$,
and $D\subset\Lob2{\kappa}$ be a closed region bounded by $\tilde \alpha$.
A length-nonincreasing map $F\:D\to\spc{X}$ is called \index{majorizing map}\emph{majorizing} if it is length-preserving on $\tilde \alpha$.

In this case, we say that $D$ \emph{majorizes} the curve $\alpha=F\circ\tilde \alpha$ under the map $F$.
\end{thm}

The following proposition is a consequence of the definition.

\begin{thm}{Proposition}
\label{prop:majorize-geodesic} 
Let  $\alpha$  be a closed curve in a metric space $\spc{X}$.
Suppose $D\subset\Lob2{\kappa}$ majorizes $\alpha$ under $F\: D \to \spc{X}$.  
Then any geodesic subarc of $\alpha$ is the image under $F$ of a subarc of $\partial_{\Lob2{\kappa}} D$ that is geodesic in the length metric of $D$.

In particular, if $D$ is convex, then the corresponding subarc is a geodesic in $\Lob2{\kappa}$.
\end{thm}

\parit{Proof.} For a geodesic subarc $\gamma\:[a,b]\to\spc{X}$ of $\alpha=F\circ\tilde \alpha$, set
\begin{align*}
\tilde r&=\dist{\tilde \gamma(a)}{\tilde \gamma(b)}{D},
&
\tilde \gamma &= (F|_{\Fr D})^{-1}\circ\gamma,
\\
s&=\length \gamma,
&
\tilde s&= \length \tilde \gamma.
\end{align*}
Then
\[\tilde r\ge r = s =\tilde s\ge\tilde r.\]
Therefore $\tilde s=\tilde r$.
\qeds

\begin{thm}{Corollary}\label{cor:maj-triangle}
Let $\trig p x y$ be a triangle of perimeter $<2\cdot\varpi\kappa$ in a metric space $\spc{X}$. Assume a convex region $D\subset \Lob2\kappa$ majorizes $\trig p x y$.
Then $D=\Conv\trig{\tilde p}{\tilde x}{\tilde y}$ for a model triangle $\trig{\tilde p}{\tilde x}{\tilde y}=\modtrig\kappa(p x y)$, and the majorizing map sends  $\tilde p$, $\tilde x$ and $\tilde y$ respectively to $p$, $x$ and $y$.
\end{thm}

Now we come to the main theorem of this section.

\begin{thm}{Majorization theorem}
\label{thm:major}
Any closed curve $\alpha$ with length smaller than $2\cdot \varpi\kappa$ in  a $\varpi\kappa$-geodesic $\CAT\kappa$ space is majorized by a convex region in $\Lob2\kappa$. \end{thm}

This theorem was proved by Yuriy Reshetnyak \cite{reshetnyak:major};
our proof uses a trick that we learned from the lectures of Werner Ballmann \cite{ballmann:lectures}.
Another proof can be built on Kirszbraun's theorem (\ref{thm:kirsz+}), but it works only for complete spaces.

The case when $\alpha$ is a triangle, say $\trig p x y$, is the base  and is nontrivial.
In this case, by Corollary~\ref{cor:maj-triangle}, the majorizing convex region has to be isometric to $\Conv\trig{\tilde p}{\tilde x}{\tilde y}$, where $\trig{\tilde p}{\tilde x}{\tilde y}=\modtrig\kappa(p x y)$.  
There is a majorizing map for $\trig p x y$ whose image $W$ is the image of the line-of-sight map (definition \ref{def:sight}) for $[x y]$ from  $p$,
but as one can see from the following example, the line-of-sight map is not majorizing in general.

\begin{wrapfigure}{r}{30 mm}
\vskip-0mm
\centering
\includegraphics{mppics/pic-950}
\end{wrapfigure}

\parbf{Example.} Let $\spc{Q}$ be a solid quadrangle $[p x z y]$ in $\EE^2$ formed by two congruent triangles, which is non-convex at $z$ (as in the picture).  
Equip $\spc{Q}$ with the length metric. 
Then $\spc{Q}$ is $\CAT0$
by Reshetnyak gluing  (\ref{thm:gluing}). 
For triangle ${\trig p x y}_\spc{Q}$ in $\spc{Q}$ and its model triangle $\trig{\tilde p}{\tilde x}{\tilde y}$ in $\EE^2$,  
we have 
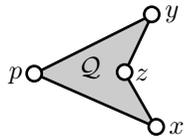
\[\dist{\tilde x}{\tilde y}{}=\dist{x}{y}{\spc{Q}}=\dist{x}{z}{}+\dist{z}{y}{}.\]
Then the map $F$ defined by matching line-of-sight parameters satisfies $F(\tilde x)=x$ and $\dist{x}{F(\tilde w)}{}>\dist{\tilde x}{\tilde w}{}$ if $\tilde w$ is near the midpoint $\tilde z$ of $[\tilde x\tilde y]$ and lies on $[\tilde p\tilde z]$. 
Indeed, by the first variation formula (\ref{1st-var+}), for $\eps=1-s$ we have
\[\dist{\tilde x}{\tilde w}{}
=\dist{\tilde x}{\tilde \gamma_\frac12(s)}{}
=\dist{x}{z}{}+o(\eps)\] and 
\[\dist{x}{F(\tilde w)}{}
=\dist{x}{\gamma_\frac12(s)}{}
=\dist{x}{z}{}-\eps\cdot\cos\mangle\hinge z p x+o(\eps).\]  
Thus $F$ is not majorizing.

\medskip

In  the following proofs, $x^1 \dots x^n$ ($n\ge 3$) denotes a polygonal line $x^1,\dots,x^n$, and $[x^1\dots x^n ]$ denotes the corresponding (closed) polygon.
For a subset $R$ of the ambient metric space,
we denote by $[x^1\dots x^n ]_R$ a polygon in the length metric of $R$.

Our first lemma gives a model space construction based on repeated application of Lemma~\ref{lem:quadrangle} from the proof of the inheritance.
Recall that convex and concave curves with respect to a point are defined in~\ref{def:convex-devel}.

\begin{thm}{Lemma}\label{lem:majorize-subgraph}
In $\Lob2{\kappa}$, let  
$\beta$ be a curve from $x$ to $y$ 
that is concave with respect  to $p$.
Let $D$  be the subgraph of $\beta$ with respect to $p$.
Assume 
\[\length\beta\z+\dist{p}{x}{}+\dist{p}{y}{}<2\cdot\varpi\kappa.\]
\begin{subthm}{curvilinear} 
Then $\beta$ forms a geodesic $[x y]_D$ in $D$ and therefore $\beta$, $[p x]$ and $[p y]$ form a triangle 
${\trig p x y}_D$ in the length metric of $D$.
\end{subthm}
\begin{subthm} {short-to-subgraph}
Let $\trig{\tilde p}{\tilde x}{\tilde y}$ be the model triangle for 
${\trig p x y}_D$.
Then there is a short map $G\:\Conv\trig{\tilde p}{\tilde x}{\tilde y}\to D$ such that $\tilde p\mapsto p$, $\tilde x\mapsto x$, $\tilde y\mapsto y$, and $G$ is length-preserving on each side of $\trig{\tilde p}{\tilde x}{\tilde y}$.
In particular, $\Conv\trig{\tilde p}{\tilde x}{\tilde y}$ majorizes triangle $[p x y]_D$ in $D$ under~$G$.
\end{subthm}
\end{thm}

\parit{Proof.}
We prove the lemma for a polygonal line $\beta$;
the general case then follows by approximation.
Namely, since $\beta$ is concave 
it can be approximated by polygonal lines that are concave with respect to $p$, 
with their lengths converging to $\length \beta$. 
Passing to a partial limit we will obtain the needed map $G$.  

Suppose $\beta=x^0x^1\dots x^n$ is a polygonal line with $x^0=x$ and $x^n=y$.
Consider a sequence of polygonal lines $\beta_i=x^0x^1\dots x^{i-1}y_i$ such that $\dist{p}{y_i}{}=\dist{p}{y}{}$ and 
$\beta_i$ has same length as $\beta$; 
that is, 
\[\dist{x^{i-1}}{y_i}{}=\dist{x^{i-1}}{x^{i}}{}+\dist{x^{i}}{x^{i+1}}{}+\dots+\dist{x^{n-1}}{x^n}{}.\]

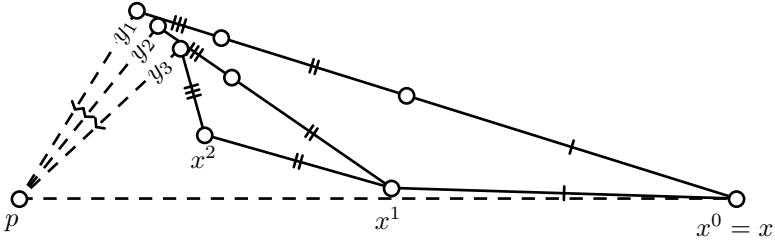
\begin{figure}[!ht]
\vskip-0mm
\centering
\includegraphics{mppics/pic-955}
\end{figure}

Clearly $\beta_n=\beta$.
Sequentially applying Alexandrov's lemma (\ref{lem:alex}) shows that each of the polygonal lines $\beta_{n-1}, \beta_{n-2},\dots,\beta_1$ is concave with respect to $p$.
Let $D_i$ be the subgraph of $\beta_i$ with respect to $p$.
Applying Lemma~\ref{lem:quadrangle} gives a short map $G_i\:D_{i}\to D_{i+1}$ that maps $y_{i}\mapsto y_{i+1}$ and does not move $p$ and $x$ (in fact,  $G_i$ is the identity everywhere except on $\Conv\trig{p}{x^{i-1}}{y_i}$).
Thus the composition 
\[G_{n-1}\circ\dots\circ G_1\: D_1\to D_n\] 
is short.
The result follows since $D_1\iso\Conv\trig{\tilde p}{\tilde x}{\tilde y}$.\qeds

\begin{thm}{Lemma}\label{lem:majorize-triangle}
Let $\trig{p}{x}{y}$ be a triangle of perimeter $<2\cdot\varpi\kappa$ in a $\varpi\kappa$-geodesic $\CAT\kappa$ space $\spc{U}$.
In $\Lob2{\kappa}$, let $\tilde \gamma$ be the $\kappa$-development of $[x y]$ with respect to $p$, where $\tilde \gamma$ has basepoint $\tilde p$ and subgraph $D$.
Consider the map $H\:D\to\spc{U}$ that sends the point with parameter $(t,s)$ under the line-of-sight map for $\tilde \gamma$ with respect to $\tilde p$, to the point with the same parameter under the line-of-sight map $f$ for $[x y]$ with respect to $p$.
Then $H$ is  length-nonincreasing.
In particular, $D$ majorizes triangle $\trig p x y$.
\end{thm}

\parit{Proof.}
Let $\gamma=\geod_{[x y]}$ and $T=\dist{x}{y}{}$. 
As in the proof of the development criterion (\ref{thm:concave-devel}), take a partition 
\[0=t^0<t^1<\dots<t^n=T,\]
and set $x^i=\gamma(t^i)$. 
Construct a chain of model triangles  $\trig{\tilde p}{\tilde x^{i-1}}{\tilde x^i}\z=
\modtrig\kappa(p x^{i-1} {x^i})$, with $\tilde x^0=\tilde x$ and the direction of $[\tilde p\tilde x^i]$ turning counterclockwise as $i$ grows.  
Let $D_n$ be the subgraph with respect to $\tilde p$ of the polygonal line $\tilde x^0\dots \tilde x^n$.

Let  $\delta_n$ be the maximum radius of a circle inscribed in any of the triangles $\trig{\tilde p}{\tilde x^{i-1}}{\tilde x^i}$.  

Now we construct a map $H_n \: D_n\to\spc{U}$  that increases distances by at most  $2\cdot\delta_n$.
Suppose $w\in D_n$.
Then $w$ lies on or inside some triangle $\trig{\tilde p}{\tilde x^{i-1}}{\tilde x^i}$.  
Define $H_n(w)$ by first mapping $w$ to a nearest point on $\trig{\tilde p}{\tilde x^{i-1}}{\tilde x^i}$ (choosing one if there are several), followed by the natural map to the triangle  $\trig {p}{x^{i-1}}{ x^i}$. 

Since triangles in $\spc{U}$ are $\kappa$-thin (\ref{prop:k-thin}), the restriction of $H_n$ to each triangle $\trig{\tilde p}{\tilde x^{i-1}}{\tilde x^i}$ is short.   
Then the triangle inequality implies that the restriction of $H_n$ to 
\[U_n=\bigcup_{1\le i\le n}\trig{\tilde p}{\tilde x^{i-1}}{\tilde x^i}\]
is short with respect to the length metric on $D_n$. 
Since nearest-point projection from $D_n$ to $U_n$ increases the $D_n$-distance between two points by at most $2\cdot\delta_n$, the map $H_n$ also increases the $D_n$-distance by at most $2\cdot\delta_n$. 

Consider converging sequences $v_n\to v$ and $w_n\to w$ such that $v_n,w_n\in D_n$ and therefore $v,w\in D$.
Note that 
\[\dist{H_n(v_n)}{H_n(w_n)}{} \le \dist{v_n}{w_n}{D_n} + 2\cdot\delta_n,\eqlbl{eq:|H(v)-H(w)|}\]
for each $n$.
Since $\delta_n\to 0$ and geodesics in $\spc{U}$ vary continuously with their endpoints (\ref{thm:alex-patch}), we have $H_n(v_n)\to 
H(v)$ and $H_n(w_n)\to H(w)$.
Therefore the left-hand side in \ref{eq:|H(v)-H(w)|} converges to $\dist{H(v)}{H(w)}{}$ and the right-hand side converges to $\dist{v}{w}{D}$, it follows that $H$ is short.
\qeds

\parit{Proof of \ref{thm:major}.}
We begin by proving the theorem in case $\alpha$ is polygonal.

First suppose $\alpha$ is a triangle, say $\trig p x y$.
By assumption, the perimeter of $\trig p x y$ is less than
$2\cdot\varpi\kappa$.
This is the base case for the induction.

 Let $\tilde \gamma$ be the $\kappa$-development of $[x y]$ with respect to $p$, where $\tilde \gamma$ has basepoint $\tilde p$ and subgraph $D$.
By the development criterion (\ref{thm:concave-devel}),  $\tilde \gamma$ is concave.
By Lemma~\ref{lem:majorize-subgraph},  there is a short map $G\:\Conv\modtrig\kappa(p x y)\to D$.
Further, by Lemma~\ref{lem:majorize-triangle},  $D$ majorizes $\trig p x y$ under a majorizing map $H\:D\to\spc{U}$. Clearly $H\circ G$ is a majorizing map for $\trig p x y$.

\begin{wrapfigure}{r}{40 mm}
\vskip-1mm
\centering
\includegraphics{mppics/pic-960}
\vskip0mm
\end{wrapfigure}

Now we claim that any closed $n$-gon $[x^1x^2 \dots x^n ]$ of perimeter less than $2\cdot \varpi\kappa$ in a $\CAT{\kappa}$ space  is majorized by a convex polygonal region \[R_n=\Conv[\tilde x^1\tilde x^2\dots\tilde x^n]\]
under a map $F_n$ such that $F_n\:\tilde x^i\mapsto x^i$ for each $i$.

Assume the statement is true for $(n-1)$-gons, $n\ge 4$.  
Then  $[x^1 x^2 \dots x^{n-1}]$  is majorized by a convex polygonal region 
\[R_{n-1}=\Conv[\tilde x^1 \tilde x^2,\dots \tilde x^{n-1}],\] 
in $\Lob2\kappa$ under a map $F_{n-1}$ satisfying $F_{n-1}(\tilde x^i)=x^i$ for all $i$. 
Take $\dot x^n\in\Lob2{\kappa}$ such that $\trig{\tilde x^1}{\tilde x^{n-1}}{\dot x^n}=\modtrig\kappa(x^1 x^{n-1} x^n)$ 
and this triangle lies on the other side of $[\tilde x^1\tilde x^{n-1}]$ from $R_{n-1}$.  
Let $\dot R\z=\Conv\trig{\tilde x^1}{\tilde x^{n-1}}{\dot x^n}$, 
and $\dot F\:\dot R\to \spc{U}$ be a majorizing map for $\trig { x^1}{x^{n-1}}{ x^n}$ as provided above.

Set 
$R= R_{n-1}\cup \dot R$, where $R$ carries its length metric.
Since $F_n$ and $F$ agree on $[\tilde x^1 \tilde x^{n-1}]$, we may define $F\:R\to\spc{U}$ by 
\[
F(x)=
\begin{cases}
F_{n-1}(x),\quad & x\in R_{n-1},\\
\dot F(x),\quad & x\in \dot R.\\
\end{cases}
\]
Then $F$ is length-nonincreasing and is a majorizing map for $[x^1 x^2 \dots x^n ]$ (as in Definition~\ref{def:majorize}).

If $R$ is a convex subset of $\Lob2\kappa$, we are done. 

If $R$ is not convex,  the total internal angle of $R$ at $\tilde x^1$ or $ \tilde x^{n-1} $ or both is $>\pi$.  
By relabeling we may suppose this holds for $\tilde x^{n-1}$.  

The region $R$ is obtained by gluing $R_{n-1}$ to $\dot R$ by $[x^1x^{n-1}]$.
Thus, by Reshetnyak gluing (\ref{thm:gluing}), $R$ carrying its length metric is a $\CAT{\kappa}$-space.  
Moreover $[\tilde x^{n-2}\tilde x^{n-1}]\cup[\tilde x^{n-1} \dot x^n]$ is a geodesic of $R$.
Thus $[\tilde x^1 \tilde x^2 \dots \tilde x^{n-2} \dot x^n]_R$ is a closed $(n-1)$-gon in $R$, to which the induction hypothesis applies. The resulting short map from a convex region in $ \Lob2\kappa$ to~$R$, followed by $F$,  is the desired majorizing map.

\medskip

Note that in fact we have proved the following:

\begin{clm}{}
Let $F_{n-1}$ be a majorizing map for the polygon $[x^1x^2\dots x^{n-1}]$,
and $\dot F$ be a majorizing map for the triangle $[x^1x^{n-1}x^{n}]$.
Then there is a majorizing map $F_n$ for the polygon $[x^1x^2\dots x^n]$
such that \[\Im F_{n+1}= \Im F_n\cup\Im \dot F.\]

\end{clm}

We now use this claim to prove the theorem for general curves.

Assume $\alpha\:[0,\ell]\to\spc{U}$ is an  arbitrary closed curve with natural parameter.
Choose a sequence of partitions $0=t^0_n<t^2_n<\dots<t^n_n=\ell$
so that:
\begin{itemize}
\item The set $\{t_{n+1}^i\}_{i=0}^{n+1}$ 
is obtained from the set  $\{t_n^i\}_{i=0}^n$ by adding one element.
\item For a sequence $\eps_n\to0+$,
we have $t^i_n-t^{i-1}_n<\eps_n$ for all $i$.
\end{itemize}

Inscribe in $\alpha$ a sequence of polygons $P_n$ with vertexes $\alpha(t^i_n)$.
Apply the claim above, to get a sequence of majorizing maps $F_n\:R_n\to\spc{U}$.
Note that for all $m>n$ we have
\begin{itemize}
\item $\Im F_m$ lies in an  $\eps_n$-neighborhood of $\Im F_n$,
\item $\Im F_m\setminus \Im F_n$ lies in an  $\eps_n$-neighborhood of $\alpha$.
\end{itemize}
It follows that the set
\[K=\alpha\cup\left(\bigcup_n\Im F_n\right)\]
is compact.
Therefore the sequence $(F_n)$
has a partial limit as $n\to\infty$; 
say $F$.
Clearly $F$ is a majorizing map for $\alpha$.
\qeds

If $p_1\dots p_n$ is a polygon, then values $\theta_i=\pi-\mangle\hinge{p_i}{p_{i-1}}{p_{i+1}}$ for all $i\pmod n$ are called \emph{external angles} of the polygon.
The following exercise is a generalization of Fenchel's theorem.

\begin{thm}{Exercise}\label{ex:fenchel}
Show that the sum of external angles of any polygon in a complete length $\CAT0$ space cannot be smaller than $2\cdot\pi$. 
\end{thm}

\begin{thm}{Very advanced exercise}\label{ex:FM}
Suppose that a simple polygon $\beta$ in a complete length $\CAT0$ space does not bound an embedded disc.
Show that the sum of external angles of $\beta$ cannot be smaller than $4\cdot\pi$.

Give an example of such a polygon $\beta$ with the sum of external angles exactly $4\cdot\pi$.
\end{thm}

The following exercise is the rigidity case 
of the majorization theorem.

{\sloppy 

\begin{thm}{Exercise}\label{ex:isometric-majorization}
Let $\spc{U}$ be a $\varpi\kappa$-geodesic $\CAT\kappa$ space
and $\alpha\:[0,\ell]\to\spc{U}$ be a closed curve with arclength parametrization.
Assume that $\ell\z<2\cdot \varpi\kappa$
and there is a closed convex curve $\tilde \alpha\:[0,\ell]\to\Lob{2}{\kappa}$ such that 
\[\dist{\alpha(t_0)}{\alpha(t_1)}{\spc{U}}=\dist{\tilde \alpha(t_0)}{\tilde \alpha(t_1)}{\Lob{2}{\kappa}}\]
for any $t_0$ and $t_1$.
Then there is a distance-preserving map $F\:\Conv\tilde \alpha\to \spc{U}$
such that $F\:\tilde \alpha(t)\mapsto \alpha(t)$ for any $t$.
\end{thm}

}

\begin{thm}{Exercise}\label{ex:bishop}
Two majorizations $F\:D\to \spc{U}$ and $F'\:D'\to \spc{U}$ will be called \index{majorizing map!equivalent majorizations}\emph{equivalent} if $F'=F\circ\iota$ for an isometry $\iota\:D\to D'$.

Show that a closed rectifiable curve in a $\CAT0$ space has an isometric majorization map if and only if the majorization map is unique up to equivalence.
\end{thm}

The following lemma states, in particular, that in a $\CAT\kappa$ space, 
a sharp triangle comparison implies the
presence of an isometric copy of the convex hull of the model triangle.
The latter statement was proved by Alexandr Alexandrov \cite{alexandrov:devel}.
  
\begin{thm}{Arm lemma}\label{lem:arm}
Let $\spc{U}$ be a $\varpi\kappa$-geodesic $\CAT\kappa$ space, 
and $P=[x^0x^1\dots x^{n+1}]$ be a polygon of length $<2\cdot \varpi\kappa$ in $\spc{U}$.
Suppose $\tilde P=[\tilde x^0\tilde x^1\dots \tilde x^{n+1}]$ is a convex  polygon in $\Lob{2}{\kappa}$
such that 
\[
\dist{\tilde x^i}{\tilde x^{i-1}}{\Lob{2}{\kappa}}
=
\dist{x^i}{x^{i-1}}{\spc{U}}
\quad \text{and}\quad 
\mangle\hinge{x^i}{x^{i-1}}{x^{i+1}}\ge\mangle\hinge{\tilde x^i}{\tilde x^{i-1}}{\tilde x^{i+1}}
\eqlbl{eq:arm}
\]
for all $i$.
Then 

\begin{subthm}{subthm:arm-ineq}
$\dist{\tilde x^0}{\tilde x^{n+1}}{\Lob{2}{\kappa}}
\le
\dist{x^0}{x^{n+1}}{\spc{U}}$.
\end{subthm}

\begin{subthm}{subthm:arm-eq}
Equality holds in \ref{SHORT.subthm:arm-ineq} if and only if the map $\tilde x^i\mapsto x^i$ can be extended 
to a distance-preserving map of $\Conv(\tilde x^0,\tilde x^1\dots \tilde x^{n+1})$ onto $\Conv(x^0,x^1\dots x^{n+1})$.
\end{subthm}
\end{thm}

\parit{Proof; \ref{SHORT.subthm:arm-ineq}.}
By majorization (\ref{thm:major}), $P$ is majorized by a convex region $\tilde D$ in $\Lob{2}{\kappa}$.
By Proposition \ref{prop:majorize-geodesic} and the definition of angle,
$\tilde D$ is bounded by a convex polygon $\tilde P_R=[\tilde y^0\tilde y^1\dots \tilde y^{n+1}]$ that satisfies
\begin{align*}
\dist{\tilde y^i}{\tilde y^{i\pm1}}{\Lob{2}{\kappa}}
&=
\dist{x^i}{x^{i\pm1}}{\spc{U}}, \qquad \dist{\tilde y^0}{\tilde y^{n+1}}{\Lob{2}{\kappa}}
=
\dist{x^0}{x^{n+1}}{\spc{U}},
\\
& \mangle\hinge{\tilde y^i}{\tilde y^{i-1}}{\tilde y^{i+1}}\ge\mangle\hinge{x^i}{x^{i-1}}{x^{i+1}}\ge\mangle\hinge{\tilde x^i}{\tilde x^{i-1}}{\tilde x^{i+1}}
\end{align*}
for $1\le i\le n$; the last inequality follows from \ref{eq:arm}.

The classical arm lemma \cite{sabitov} gives $\dist{\tilde x^0}{\tilde x^{n+1}}{}\le \dist{\tilde y^0}{\tilde y^{n+1}}{}$.
Since $ \dist{\tilde y^0}{\tilde y^{n+1}}{}=\dist{x^0}{x^{n+1}}{}$, part \ref{SHORT.subthm:arm-ineq} follows.

\parit{\ref{SHORT.subthm:arm-eq}.}
Suppose equality holds in \ref{SHORT.subthm:arm-ineq}.
Then angles at the $j$-th vertex of~$\tilde P$, $P$, and $\tilde P_R$ are equal for $1\le j\le n$, and we may take $\tilde P=\tilde P_R$.  

Let $F\:\tilde D\to\spc{U}$ be the majorizing map for $P$, where $\tilde D$ is the convex region bounded by $\tilde P$, and $F|_{\tilde P}$ is length-preserving.  

\begin{clm}{}\label{clm:arm-triangle}
Let $\tilde x,\tilde y,\tilde z$ be three vertexes of $\tilde P$, and $x,y,z$ be the corresponding vertexes of $P$.  If $\dist{\tilde x}{\tilde y}{}=
\dist{x}{y}{}$, $\dist{\tilde y}{\tilde z}{}=
\dist{x}{z}{}$ and $\mangle\hinge{\tilde y}{\tilde x}{\tilde z} = \mangle\hinge{y}{x}{z}$, then $F|_{\Conv(\tilde x, \tilde y, \tilde z)}$ is distance-preserving.
\end{clm} 

Indeed, since $F$ is majorizing, $F$ restricts to   distance-preserving maps from $[\tilde x\tilde y]$ to $[xy]$ and $[\tilde y\tilde z]$ to $[yz]$.
Suppose $\tilde p\in [\tilde x \tilde y]$ and $\tilde q\in[\tilde y\tilde z]$.  Then 
\[
\dist{\tilde p}{\tilde q}{\Lob{2}{\kappa}}
=
\dist{F( \tilde p )}{F(\tilde q)}{\spc{U}}.
 \eqlbl{eq:arm-eq}
\]
This inequality holds in one direction by majorization, and in the other direction by the angle comparison (\ref{cat-hinge}).
By the first variation formula (\ref{cor:both-end-first-var-cba}), it follows that each pair of corresponding angles of triangles $[\tilde x \tilde y \tilde z]$ and $[x y z]$ are equal.
But then \ref{eq:arm-eq} holds for $p,q$ on any two sides of these triangles, so $F$ is distance-preserving on every geodesic of $\Conv(\tilde p, \tilde x, \tilde y)$.
Hence the claim.
\claimqeds

\begin{clm}{}\label{clm:arm-induction}
Suppose $F|_{\Conv(\tilde x^0,\tilde x^1,\dots,\tilde x^{k})}$ is distance-preserving for $2\le k\le n-1$.
Then $F|_{\Conv(\tilde x^0,\tilde x^1,\dots, \tilde x^{k+1})}$ is distance-preserving.
\end{clm}

To verify this claim, let 
\[
\tilde p=[\tilde x^{k-1}\tilde x^{k+1}] \cap [\tilde x^{k}\tilde x^{0}]
\quad\text{and}\quad
p=F(p).
\]

Note that the following maps are distance-preserving:
\begin{enumerate}
\item[(i)]
$F|_{\Conv(\tilde x^{k-1},\tilde x^k,\tilde x^{k+1})}$,

\item[(ii)]
 $F|_{\Conv(\tilde x^{k+1},\tilde x^{k-1},\tilde x^{0})}$,

\item[(iii)]
$F|_{\Conv(\tilde x^{0},\tilde x^{k},\tilde x^{k+1})}$.
\end{enumerate}
Indeed, (i) follows from \ref{clm:arm-triangle}.  
Therefore $\dist{\tilde x^{k-1}}{\tilde x^{k+1}}{}=\dist{x^{k-1} }{x^{k+1}}{}$, and so $F$ restricts to a distance-preserving map from $[\tilde x^{k-1}\tilde x^{k+1}]$ onto $[x^{k-1} x^{k+1}]$.  With the induction hypothesis in
\ref{clm:arm-induction},
 it follows that $p=[x^{k-1}x^{k+1}] \cap [x^{k}x^{0}]$, hence 
\[
\mangle\hinge{\tilde x^{k-1}}{\tilde x^{k+1}}{\tilde x^{0}} 
= \mangle\hinge{x^{k-1}}{x^{k+1}}{x^{0}}.
 \eqlbl{eq:angle-arm-eq}
\] 
Then (ii) follows from \ref{eq:angle-arm-eq} and \ref{clm:arm-triangle}.  Since $\dist{\tilde x^{k}}{\tilde x^{0}}{}=\dist{x^{k} }{x^{0}}{}$, (iii) follows from \ref{eq:angle-arm-eq} and (i). 

Let $\tilde \gamma$ be a geodesic of $\Conv(\tilde x^{0},\tilde x^0,\tilde x^1\dots \tilde x^{k+1})$.
Then $\length \tilde \gamma \z< \varpi\kappa$.
If $\tilde \gamma$ does not contain the point $\tilde p$, then by the induction hypothesis in
\ref{clm:arm-induction} and (i)+(ii)+(iii), we get that  $\gamma = F\circ\tilde \gamma$ is a local geodesic of length $<\varpi\kappa$.
By \ref{cor:loc-geod-are-min}, $\gamma$ is a geodesic.

By continuity, $F\circ\tilde \gamma$ is a geodesic for all $\tilde \gamma$;
so \ref{clm:arm-induction} follows.

The base of the induction is provided by \ref{clm:arm-triangle}.
It finishes the proof of part \ref{SHORT.subthm:arm-eq}.
\qeds
 
\begin{thm}{Exercise}\label{ex:square}
Let $\spc{U}$ be a complete length $\CAT0$ space and 
suppose 
for 4 points $x^1,x^2,x^3,x^4\in \spc{U}$
there is a convex quadrangle
$[\tilde x^1\tilde x^2\tilde x^3\tilde x^4]$
in $\EE^2$
such that 
\[\dist{x^i}{x^j}{\spc{U}}=\dist{\tilde x^i}{\tilde x^j}{\EE^2}\]
for all $i$ and $j$.
Show that $\spc{U}$ contains an isometric copy of the 
\emph{solid quadrangle}
$[\tilde x^1\tilde x^2\tilde x^3\tilde x^4]$; that is, the convex hull of $\tilde x^1,\tilde x^2,\tilde x^3,\tilde x^4$ in $\EE^2$.
\end{thm}


\section{Hadamard--Cartan theorem}\label{sec:Hadamard--Cartan}

The development of Alexandrov geometry was greatly influenced by the Hadamard--Cartan theorem.
Its original formulation states that if $M$ is a complete Riemannian manifold with nonpositive sectional curvature, 
then the exponential map at any point $p\in M$ is a covering;
in particular it implies that the universal cover of $M$ is diffeomorphic to the Euclidean space of the same dimension.

In this generality, the theorem appeared in the lectures of \'Elie Cartan \cite{cartan}.
For surfaces in the Euclidean space, 
the theorem was proved by
Hans von Mangoldt \cite{mangoldt},  
and a few years later independently by Jacques Hadamard \cite{hadamard}.

Formulations for metric spaces of different generality were proved by 
Herbert Busemann \cite{busemann-CBA},
Willi Rinow \cite{rinow}, and 
Michael Gromov \cite[p.~119]{gromov:hyp-groups}. 
A detailed proof of Gromov's statement when $\spc{U}$ is proper  was given by Werner Ballmann \cite{ballmann:cartan-hadamard}, using Birkhoff's curve-shortening.  
A proof in the non-proper 
geodesic case 
was given by the first author and Richard Bishop~\cite{alexander-bishop:h-c}.  
This proof applies more generally, to {}\emph{convex spaces} (see Exercise \ref{ex:cats-cradle}).
It was pointed out by Bruce Kleiner \cite{ballmann:lectures} 
and independently by Martin Bridson and Andr\'{e} Haefliger \cite{bridson-haefliger} that 
this proof extends to length spaces, as well as geodesic spaces, giving the following statement:

\begin{thm}{Hadamard--Cartan theorem}
\label{thm:hadamard-cartan}
Let $\kappa\le 0$, and $\spc{U}$ be a complete, simply connected length locally $\CAT\kappa$ space.
Then $\spc{U}$ is $\CAT\kappa$.
\end{thm}

\parit{Proof.} Since $\varpi\kappa=\infty$, Theorem~\ref{thm:globalization-lift} implies that there is a $\CAT\kappa$ space $\spc{B}$ and a \index{metric covering}\emph{metric covering} $\map\:\spc{B}\to\spc{U}$; that is, $\map$ is a length-preserving covering map. 

Since $\spc{U}$ is simply connected, $\map\:\spc{B}\to\spc{U}$ is an isometry --- hence the result.
\qeds

To formulate the generalized Hadamard--Cartan theorem,
we need the following definition.

\begin{thm}{Definition}\label{def:l-s.c.}
Given $\ell\in (0,\infty]$,
a metric space $\spc{X}$ is called 
$\ell$-simply connected 
if it is connected and 
any closed curve of length $<\ell$ 
is null-homotopic in the class of curves of length $<\ell$ in $\spc{X}$.
\end{thm}

Note that there is a subtle difference between 
simply connected and $\infty$-simply connected spaces;
the first states that any closed curve is null-homotopic while the second means that any rectifiable curve is null-homotopic in the class of rectifiable curves.
However, as follows from Proposition~\ref{prop:sc}, for locally $\CAT\kappa$ spaces these two definitions are equivalent.
This fact makes it possible to deduce the Hadamard--Cartan theorem directly from the generalized Hadamard--Cartan theorem.

\begin{thm}{Generalized Hadamard--Cartan theorem}\label{thm:hadamard-cartan-gen}
A complete length space
$\spc{U}$ is $\CAT\kappa$ 
if and only if $\spc{U}$ is $2\cdot\varpi\kappa$-simply connected
and $\spc{U}$ is locally $\CAT\kappa$.
\end{thm}

For proper spaces, the generalized Hadamard--Cartan theorem was proved by Brian Bowditch \cite{bowditch}.
In the proof we need the following lemma.

\begin{thm}{Lemma}
Assume $\spc{U}$ is a complete length  locally $\CAT\kappa$ space,
$\eps>0$,
and $\gamma_1,\gamma_2\:\mathbb{S}^1\to\spc{U}$ are two closed curves.
Assume 
\begin{subthm}{}
$\length\gamma_i<2\cdot\varpi\kappa-4\cdot\eps$ for $i=1,2$;
\end{subthm}
 
\begin{subthm}{} $\dist{\gamma_1(x)}{\gamma_2(x)}{}<\eps$ for any $x\in\mathbb{S}^1$, and the geodesic $[\gamma_1(x)\gamma_2(x)]$ is uniquely defined and depends continuously on $x$;
\end{subthm}

\begin{subthm}{}  $\gamma_1$ is majorized by a convex region in $\Lob2\kappa$.
\end{subthm}

Then  $\gamma_2$ is majorized by a convex region in $\Lob2\kappa$.
\end{thm}

\parit{Proof.} Let $D$ be a convex region in $\Lob2\kappa$ that majorizes $\gamma_1$ under the map $F\:D\to\spc{U}$ 
(see Definition~\ref{def:majorize}).
Denote by $\tilde \gamma_1$ 
the curve bounding $D$ 
such that $F\circ\tilde \gamma_1=\gamma_1$.
Since  
\begin{align*}
\length\tilde \gamma_1
&=
\length\gamma_1
<
\\
&<
2\cdot\varpi\kappa-4\cdot\eps,
\end{align*}
there is a point $\tilde p\in D$ such that 
$\dist{\tilde p}{\tilde \gamma(x)}{\Lob2\kappa}<\tfrac{\varpi\kappa}2-\eps$
for any $x\in\mathbb{S}^1$.
Denote by $\alpha_x$ the concatenation of the paths $F\circ\geodpath_{[p\tilde \gamma_1(x)]_{\Lob2\kappa}}$ 
and  $\geodpath_{[\gamma_1(x)\gamma_2(x)]}$ in $\spc{U}$.
Note that $\alpha_x$ depends continuously on $x$, and
$$\length\alpha_x<\tfrac{\varpi\kappa}{2}\quad \text{and}\quad \alpha_x(1)=\gamma_2(x)$$ 
hold for any $x$.

Let us apply the lifting globalization theorem (\ref{thm:globalization-lift}) for $p\z=F(\tilde p)$.
We obtain a $\varpi\kappa$-geodesic $\CAT\kappa$ space $\spc{B}$
and a locally distance-preserving map $\map\:\spc{B}\to\spc{U}$
with $\map(\hat p)=p$ for some $\hat p \in \spc{B}$, and with the lifting property for the curves starting at $p$ with length $<\varpi\kappa/2$.
Applying the lifting property for $\alpha_x$, 
we get existence of a curve $\hat\gamma_2\:\mathbb{S}^1\to \spc{B}$ such that
$$\gamma_2=\map\circ\hat\gamma_2.$$

Since $\spc{B}$ is a geodesic $\CAT\kappa$ space, we can apply the majorization theorem (\ref{thm:major}) for $\hat\gamma_2$.
The composition of the obtained majorization with $\map$ is a majorization of $\gamma_2$.
\qeds

\parit{Proof of Theorem \ref{thm:hadamard-cartan-gen}.}
The only-if part follows from the Reshetnyak majorization theorem (\ref{thm:major}).

Let  $\gamma_t$, $t\in[0,1]$ 
be a null-homotopy of curves in $\spc{U}$;
that is, $\gamma_0(x)\z=p$ for some $p\in \spc{U}$
and any $x\in\mathbb{S}^1$.
Assume further that $\length \gamma_t\z<2\cdot\varpi\kappa$ for any $t$.
To prove the if part, it is sufficient to show that $\gamma_1$ is majorized by a convex region in $\Lob2\kappa$ if $\spc{U}$ is locally $\CAT\kappa$. 

By semicontinuity of length (\ref{thm:semicont-of-length}),
we can choose  $\eps>0$ sufficiently small that
$$\length \gamma_t<2\cdot\varpi\kappa-4\cdot\eps$$
for all $t$.

By Corollary~\ref{cor:loc-CAT(k)},
we may assume in addition that
$\oBall(\gamma_t(x),\eps)$ is $\CAT\kappa$ 
for any $t$ and $x$.

Choose a partition $0=t_0<t_1<\dots<t_n=1$
so that $\dist{\gamma_{t_i}(x)}{\gamma_{t_{i-1}}(x)}{}<\eps$
for any $i$ and $x$.
According to \ref{thm:cat-unique},
for any $i$,
the geodesic $[\gamma_{t_i}(x)\gamma_{t_{i-1}}(x)]$ depends continuously on $x$.

Note that $\gamma_0=\gamma_{t_0}$ is majorized by a convex region in $\Lob2\kappa$.
Applying the lemma $n$ times, we see that the same holds for $\gamma_1\z=\gamma_{t_n}$.
\qeds

\begin{thm}{Proposition}\label{prop:sc}
Let $\spc{U}$ be a complete length locally $\CAT\kappa$ space.
Then $\spc{U}$ is simply connected if and only if it is $\infty$-simply connected.
\end{thm}

\parit{Proof; if part.}
It is sufficient to show that any closed curve in $\spc{U}$ is homotopic to a polygon.

Let $\gamma_0$ be a closed curve in $\spc{U}$.
According to Corollary~\ref{cor:loc-CAT(k)},
there is $\eps>0$ such that 
$\oBall(\gamma(x),\eps)$ is $\CAT\kappa$
for any $x$.

Choose a polygon $\gamma_1$ such that $\dist{\gamma_0(x)}{\gamma_1(x)}{}<\eps$ for any $x$.
By \ref{thm:cat-unique}, 
$\geodpath_{[\gamma_0(x)\gamma_1(x)]}$ 
is uniquely defined 
and depends continuously on~$x$.

Hence $\gamma_t(x)=\geodpath_{[\gamma_0(x)\gamma_1(x)]}(t)$ gives a homotopy from $\gamma_0$ to $\gamma_1$.

\parit{Only-if part.} The proof is similar.

Assume $\gamma_t$ is a homotopy between two rectifiable curves $\gamma_0$ and $\gamma_1$.
Fix $\eps>0$ so that the ball $\oBall(\gamma_t(x),\eps)$ is $\CAT\kappa$
for any $t$ and $x$.
Choose a partition $0=t_0<t_1<\dots<t_n=1$ 
so that 
$$\dist{\gamma_{t_{i-1}}(x)}{\gamma_{t_i}(x)}{}<\tfrac\eps{10}$$
for any $i$ and $x$.
Set $\hat\gamma_{t_0}=\gamma_0$, $\hat\gamma_{t_n}=\gamma_{t_n}$.
For each $0<i<n$, approximate $\gamma_{t_i}$ by a polygon $\hat\gamma_{i}$.

Construct the \index{geodesic homotopy}\emph{geodesic homotopy} 
from $\hat\gamma_{t_{i-1}}$ 
to $\hat\gamma_{t_i}$;  
that is,
set 
$$\hat\gamma_t
=
\geodpath_{[\hat\gamma_{t_{i-1}}(x)\hat\gamma_{t_i}(x)]}(t)$$
for $t\in [t_{i-1},t_i]$.
Since $\eps$ is sufficiently small, 
by \ref{cor:cat-unique}, we get that
$$\length\hat\gamma_t
<
10\cdot(\length\hat\gamma_{t_{i-1}}+\length\hat\gamma_{t_i})$$
for any $t\in [t_{i-1},t_i]$.
In particular, $\hat\gamma_t$ is rectifiable for all $t$.

Joining the obtained homotopies for all $i$, we obtain a homotopy from $\gamma_0$ to $\gamma_1$ in the class of rectifiable curves.
\qeds

\begin{thm}{Exercise}\label{ex:cover-branching-along-2-lines}
Let $\spc{X}$ be a double cover of $\EE^3$ that branches along two distinct lines $\ell$ and $m$.
Show that  $\spc{X}$ is $\CAT0$ if and only if $\ell$ intersects $m$ at a right angle.
\end{thm}

\begin{thm}{Exercise}\label{ex:branching}
Let $\spc{U}$ be a complete length $\CAT0$ space.
Assume $\tilde{\spc{U}}\to \spc{U}$ is a metric covering branching along a geodesic.
Show that $\tilde{\spc{U}}$ is $\CAT0$.

More generally, assume $A\subset \spc{U}$ is a closed convex subset and $f\:\spc{X}\to \spc{U}\setminus A$ is a metric cover.
Denote by $\bar{\spc{X}}$ the completion of $\spc{X}$, and 
$\bar f\:\bar{\spc{X}}\to \spc{U}$ the continuous extension of $f$.
Let $\tilde{\spc{U}}$ be the space glued from $\bar{\spc{X}}$ and $A$ by identifying $x$ and $\bar f(x)$ if $\bar f(x)\in A$.
Show that $\tilde{\spc{U}}$ is $\CAT0$.
\end{thm}

\section{Convex sets}
\label{sec:convex-CBA}

Recall that according to Corollary~\ref{cor:convex-balls}, any ball (closed or open) of radius $R<\tfrac{\varpi\kappa}2$ in a $\varpi\kappa$-geodesic $\CAT\kappa$ space is convex.
From the uniqueness of geodesics in $\CAT\kappa$ spaces (\ref{thm:cat-unique}) we get the following:

\begin{thm}{Observation}
Any weakly $\varpi\kappa$-convex set 
in a complete length $\CAT\kappa$ space is $\varpi\kappa$-convex.
\end{thm}

\begin{thm}{Closest-point projection lemma}\label{lem:closest point}{\sloppy 
Let $p$ be a point in a complete length $\CAT\kappa$ space $\spc{U}$, and $K\subset \spc{U}$ be a closed $\varpi\kappa$-convex set. 
If $\distfun{K}{p}{}<\tfrac{\varpi\kappa}2$,
then there is a unique point $p^*\in K$ that minimizes the distance to $p$;
that is, $\dist{p^*}{p}{}=\distfun{K}{p}{}$. 

}

\end{thm}

\parit{Proof.} 
Fix $r$ properly between $\distfun{K}{p}{}$ and $\tfrac{\varpi\kappa}2$.
By the function comparison (\ref{thm:function-comp}),
the function $f=\md\kappa\circ\distfun{p}{}{}$ is strongly convex in $\cBall[p,r]$.

The lemma follows from Lemma~\ref{lem:argmin(convex)} applied to the subspace $K'\z=K\cap\cBall[p,r]$ 
and the restriction $f|_{K'}$. 
\qeds

\begin{thm}{Exercise}\label{ex:closest-point-projection}
Let  $\spc{U}$ be a complete length $\CAT0$ space and $K\subset \spc{U}$ be a closed convex set.
Show that the closest-point projection $\spc{U}\to K$ is short. 
\end{thm}

\begin{thm}{Advanced exercise}\label{ex:short-retraction-CBA(1)}
Let  $\spc{U}$ be a complete length $\CAT1$ space and $K\subset \spc{U}$ be a closed $\pi$-convex set.
Assume $K\subset \cBall[p,\tfrac\pi2]$ for $p\in K$.
Show that there is a short retraction of $\spc{U}$ to~$K$. 
\end{thm}

\begin{thm}{Proposition}\label{lem:dist-to-convex}
Let $\spc{U}$  be a $\varpi\kappa$-geodesic $\CAT\kappa$ space
and $K\subset \spc{U}$ be a closed $\varpi\kappa$-convex set.
Let
\[f=\sn\kappa\circ\distfun{K}{}{}.\]
Then
\[f''+\kappa \cdot f\ge 0\]
holds in $\oBall(K,\tfrac{\varpi\kappa}2)$.
\end{thm}

\parit{Proof.}
It is sufficient to show that Jensen's inequality (\ref{y''-mono})
holds on a sufficiently short 
geodesic $[pq]$ in $\oBall(K,\tfrac{\varpi\kappa}2)$.
We may assume that 
\[\dist{p}{q}{}+\distfun{K}{p}{}+\distfun{K}{q}{}<\varpi\kappa.\eqlbl{eq:sum=<varpi}\]

For each $x\in[pq]$,
we need to find a value $h(x)\in \RR$
such that 
\[
h(p)=f(p),\qquad 
h(q)=f(q),\qquad
h(x)\le f(x)\qquad
\]
for any $x\in [pq]$,
and
\[h''+\kappa\cdot h\ge 0\eqlbl{h''+kh=<0}\]
along $[pq]$.

Denote by $p^{*}$ and $q^{*}$ the closest-point projections of $p$ and $q$ on $K$; 
they are provided by \ref{lem:closest point}.
From \ref{eq:sum=<varpi} and the triangle inequality,
we have
\[\dist{p^*}{q^*}{}<\varpi\kappa.\]
Since $K$ is $\varpi\kappa$-convex, $K\supset[p^*q^*]$;
in particular
\[\distfun{K}{x}{}\le \distfun{[p^*q^*]}{x}{}\]
for any $x\in \spc{U}$.

There is a majorizing map $F: D\to \spc U$ for quadrangle $[pp^{*}q^{*}q]$, as in Definition  \ref{def:majorize} and the Reshetnyak majorization theorem (\ref{thm:major}).
By Proposition \ref{prop:majorize-geodesic}, 
the figure $D$ is a solid convex quadrangle $[\tilde p\tilde p^*\tilde q^*\tilde q]$ in $\Lob2\kappa$ such that 
\begin{align*}
\dist{\tilde p}{\tilde p^*}{\Lob2\kappa}&=\dist{p}{p^*}{\spc{U}}
&
\dist{\tilde p}{\tilde q}{\Lob2\kappa}&=\dist{p}{q}{\spc{U}}
\\
\dist{\tilde q}{\tilde q^*}{\Lob2\kappa}&=\dist{q}{q^*}{\spc{U}}
&
\dist{\tilde p^*}{\tilde q^*}{\Lob2\kappa}&=\dist{p^*}{q^*}{\spc{U}}.
\end{align*}
Given $x\in [pq]$, denote by $\tilde x$ the corresponding point on $[\tilde p\tilde q]$.
Then
\[\distfun{[pq]}{x}{}
\le
\distfun{[\tilde p\tilde q]}{\tilde x}{}.\]
Set 
\[h(x)
=
\sn\kappa\circ
\distfun{[\tilde p\tilde q]}{\tilde x}{}.\]
By straightforward calculations, \ref{h''+kh=<0} holds
and hence the statement follows.
\qeds

\begin{thm}{Corollary}\label{cor::dist-to-convex}
Let $\spc{U}$  be a complete length $\CAT\kappa$ space
and $K\subset \spc{U}$ be a closed  locally convex set.
Then there is an open set $\Omega\supset K$
such that the function 
$f=\sn\kappa\circ\distfun{K}{}{}$
satisfies 
\[f''+\kappa\cdot f\ge 0\]
in $\Omega$.
\end{thm}

\parit{Proof.}
Fix $p\in K$.
By Corollary~\ref{cor:convex-balls},
$\cBall[p,r]$ is convex for all small $r>0$.

Since $K$ is locally convex, there is $r_p>0$ such that 
the intersection
$K'=K\cap \oBall(p,r_p)$ is convex. 

Note that 
\[\distfun{K}{x}{}=\distfun{K'}{x}{}\]
for any $x\in \oBall(p,\tfrac{r_p}2)$.
Therefore the statement holds for 
\[\Omega=\bigcup_{p\in K}\oBall(p,\tfrac{r_p}2).\]
\qedsf

\begin{thm}{Theorem}\label{thm:local-global-convexity}
Assume $\spc{U}$ is a complete length $\CAT\kappa$ space and $K\subset \spc{U}$ is a closed connected locally convex set.
Assume $\dist{x}{y}{}<\varpi\kappa$ for any $x,y\in K$.
Then $K$ is convex.

In particular, if $\kappa\le 0$, then any closed connected locally convex set in $\spc{U}$ is convex.
\end{thm}

The following proof is due to Sergei Ivanov \cite{ivanov:local-global-convexity}.

\begin{wrapfigure}{o}{50 mm}
\vskip-1mm
\centering
\includegraphics{mppics/pic-970}
\vskip0mm
\end{wrapfigure}

\parit{Proof.}
Since $K$ is locally convex,
it is locally path-connected.
Since $K$ is connected and locally path connected it is path-connected.

Fix two points $x,y\in K$. 
Let us connect $x$ to $y$ by a path $\alpha\:[0,1]\to K$.
Since $\dist{x}{\alpha(s)}{}<\varpi\kappa$ for any $s$,
Theorem~\ref{thm:cat-unique} implies that the geodesic $[x\alpha(s)]$ 
is uniquely defined and depends continuously on $s$.

Let $\Omega\supset K$ be the open set provided by Corollary~\ref{cor::dist-to-convex}.
If $[xy]\z=[x\alpha(1)]$ does not completely lie in $K$, then 
there is a value $s\in [0,1]$ such that $[x\alpha(s)]$ 
lies in $\Omega$ but does not completely lie in $K$.
By Corollary~\ref{cor::dist-to-convex},
the function $f=\sn\kappa\circ\distfun{K}{}{\spc{U}}$ 
satisfies the differential inequality
\[f''+\kappa\cdot f\ge 0\eqlbl{f''+kappa f=<0}\]
along $[x\alpha(s)]$.

Since 
\begin{align*}
\dist{x}{\alpha(s)}{}&<\varpi\kappa,
&
f(x)&=f(\alpha(s))=0,
\end{align*}
then the barrier inequality (\ref{barrier}) 
implies that $f(z)\le 0$ for $z\in [x\alpha(s)]$;
that is $[x\alpha(s)]\subset K$, a contradiction.
\qeds

\section{Remarks}

The following question was known in folklore in the 80's,
but it seems that in print
it was first mentioned by Michael Gromov \cite[6.B$_1\mathrm{(f)}$]{gromov:asymt-inv}.
We do not see any reason why it should be true, 
but we also cannot construct a counterexample.

\begin{thm}{Open question}
Let $\spc{U}$ be a complete length $\CAT0$ space and $K\subset \spc{U}$ be a compact set.
Is it true that $K$ lies in a convex compact set $\bar K\subset\spc{U}$?
\end{thm}

The question can  easily be reduced to the case when $K$ is finite;
so far it is not even known if any three points in a complete length $\CAT0$ space lie in a compact convex set.

We expect that \ref{lem:patch} can be extended to all curves, not necessarily local geodesics,
but the proof does not admit a straightforward generalization.

\parbf{About convex spaces.}
A \index{convex space}\emph{convex space} $\spc{X}$ is a geodesic space such that the function
$t\mapsto\dist{\gamma(t)}{\sigma(t)}{}$ is convex 
for any two  geodesic paths $\gamma,\sigma:[0,1]\to \spc{X}$.  
A \index{convex space!locally convex space}\emph{locally convex space} is a length space in which every point has a neighborhood that is a convex space in the restricted metric.

\begin{thm}{Exercise}\label{ex:cats-cradle}
Assume $\spc{X}$ is a convex space 
such that the angle of any hinge is defined.
Show that $\spc{X}$ is $\CAT{0}$.
\end{thm}

The following exercise gives an analog of the Hadamard--Cartan theorem for locally convex spaces;
see also \cite{alexander-bishop:h-c}.

\begin{thm}{Exercise}\label{ex:Hadamard--Cartan}
Show that a complete, simply connected, locally convex space is a convex space.
\end{thm}

\parbf{Examples and constructions.}
Let us list important sources of examples of $\CAT{}$ spaces.
This should be beneficial to the reader despite that we do not provide all the proofs and some proofs are deferred to later chapters.

\textit{Riemannian manifolds} with sectional curvature at most $\kappa$ are locally $\CAT{\kappa}$.
This statement follows from the Rauch comparison, and it is one of the main motivations for $\CAT{\kappa}$ comparison.
\textit{Hilbert spaces} are another motivating example of $\CAT{0}$ spaces.

The question of when a Riemannian \textit{manifold-with-boundary} is locally $\CAT{\kappa}$ has been completely answered by the first author, David Berg, and Richard Bishop \cite{alexander-berg-bishop1993}.
If the sectional curvatures of the interior and the outward sectional curvatures of the boundary do not exceed $\kappa$ then it is locally $\CAT\kappa$ (where an outward sectional curvature is one that corresponds to a tangent $2$-plane all of whose normal curvature vectors point outward).
In particular, if a Riemannian manifold and its boundary have sectional curvature at most $\kappa$, then it is locally $\CAT{\kappa}$.

Subsets of \textit{positive reach} in  Riemannian manifolds were studied in this context by 
the first author, Richard Bishop, and
Alexander Lytchak \cite{alexander-bishop,lytchak2004}.
In particular, any compact subset of positive reach in a Riemannian manifold is $\CAT{}$;
as usual, we assume that it is equipped with the induced length metric.

It was shown by Alexander Lytchak and Stephan Stadler \cite{lytchak-stadler2021} that any \textit{simply-connected subset} of a contractible two-dimensional $\CAT{\kappa}$ length space (equipped with induced length metric) is $\CAT{\kappa}$.
In higher dimensions things are more complicated, even for a Euclidean ambient space;
see \cite[Chapter 4]{alexander-kapovitch-petrunin-CAT} for related questions (open and solved).

By the Gauss formula smooth saddle surfaces in manifolds with sectional curvature at most $\kappa$ are locally $\CAT\kappa$.
The nonsmooth analog of this statement is wide open; this is the so-called Shefel conjecture \cite[Chapter 4]{alexander-kapovitch-petrunin-CAT}.
However, the following weaker statement was proved by the third author and Stephan Stadler \cite{petrunin-stadler}: \textit{metric-minimizing surfaces} in $\CAT\kappa$ space are locally $\CAT\kappa$;
metric minimizing means that it is impossible to decrease its length metric by a small deformation.

Note that any \textit{metric tree} (see Section~\ref{sec:Ultralimit of spaces}) is $\CAT{-\infty}$;
that is, a metric tree is $\CAT{\kappa}$ for \textit{any} $\kappa\in\RR$.
In particular, any tree (in the graph-theoretical sense) with a length metric such that every edge is isometric to a line segment is $\CAT{-\infty}$.

Suppose $A^1, \dots, A^n$ is an array of convex closed sets in the Euclidean space $\EE^m$.
Let us prepare $n+1$ copy of $\EE^m$ and glue successive pairs of spaces along $A^1, \dots, A^n$.
The resulting space is called \textit{puff pastry};
by the Reshetnyak gluing theorem it is $\CAT0$.
(This observation was used in one of the most beautiful applications of the Reshetnyak gluing theorem given by Dmitri Burago,  Serge Ferleger,
and Alexey Kononenko \cite{burago-ferleger-kononenko1998-1,burago-ferleger-kononenko1998-2,burago-ferleger-kononenko1998-3,burago-ferleger-kononenko1998-4}.
They use it to study billiards; a short survey on the subject was written by Dmitri Burago \cite{burago-1998};
see also our book \cite{alexander-kapovitch-petrunin-CAT}.)

An if-and-only-if condition on \textit{polyhedral spaces} is given in \ref{thm:PL-CAT}.
It implies the so-called Gromov's flag condition \ref{thm:flag} which provides in particular a flexible way to construct $\CAT0$ \textit{cube complexes}.
(Several applications are mentioned in Section~\ref{sec:poly-reamarks}.)

By \ref{prop:CAT-olim} and \ref{thm:ultra-GH}, the \textit{ultralimits}, as well as \textit{Gromov--Hausdorff limits} of $\CAT\kappa$ spaces are $\CAT\kappa$.
In particular, if  $\spc{U}$ is a $\CAT{0}$ space then its \textit{asymptotic cone}\index{asymptotyc cone}   defined as  the ultralimit of its rescalings $\tfrac1n\cdot \spc{U}$ as $n\to \o$ is again  $\CAT{0}$.
Unlike in the case of $\Alex{0}$ spaces, the use of ultralimits is necessary even if $\spc{U}$ is a manifold due to the lack of compactness theorem.
Asymptotic cones of $\CAT{0}$ spaces and their generalizations play an important role in geometric group theory.

\textit{Conformal deformations} of $\CAT{}$ spaces were studied by Alexander Lytchak and Stephan Stadler \cite{lytchak-stadler}.
In particular, if $\spc{U}$ is a $\CAT0$ space and $f\:\spc{U}\to\RR$ is  continuous, convex
and bounded below,
then the conformally equivalent space with conformal factor $e^f$ is $\CAT0$.

Further, $\CAT{}$ spaces behave nicely with respect to some natural constructions.
For example, the product of $\CAT{0}$ spaces is again $\CAT{0}$.
Also, the Euclidean cone over a $\CAT{1}$ space is  $\CAT{0}$.
These are the first examples of the so-called \textit{warped products} that are discussed in Chapter~\ref{chapter:warped products};
a general statement is given in \ref{thm:warp-CAT}.
Also, as it was observed by Karl-Theodor Sturm \cite[Prop. 3.10]{sturm2003},
the \textit{space of $L^2$-maps} from a measure space to a complete $\CAT0$ length space is $\CAT0$.

Among more conceptual examples, let us mention a result of Brian Clarke \cite{clarke}: the \textit{space of Riemannian metrics} on a compact, orientable smooth manifold
with respect to the $L^2$-distance is $\CAT0$;
a shorter proof of this statement was given by Nicola Cavallucci \cite{cavallucci}.
By a result of Tam\'{a}s Darvas \cite{darvas}, the \textit{space of K\"{a}hler potentials} on a compact K\"{a}hler manifold is $\CAT0$.
The \textit{Teichm\"{u}ller space} with the Weil--Petersson metric is $\CAT0$; 
the latter was shown by Sumio Yamada \cite{yamada}.
\chapter{Kirszbraun revisited}

This chapter is based on our paper \cite{alexander-kapovitch-kirszbraun}
and an earlier paper of Urs Lang and Viktor Schroeder
\cite{lang-schroeder}.

\section{Short map extension definitions}\label{sec:4pt}

\begin{thm}{Theorem}\label{thm:kirsz-def} 
A complete length space
$\spc{L}$ is $\Alex\kappa$ if and only if for any 3-point set $V_3$ and any 4-point set $V_4\supset V_3$ in $\spc{L}$, 
any short map $f\:V_3\to\Lob2\kappa$ can be extended to a short map $F\:V_4\to\Lob2\kappa$ (so $f=F|_{V_3}$).
\end{thm}

The only-if part of Theorem \ref{thm:kirsz-def} can be obtained as a corollary of Kirszbraun's theorem (\ref{thm:kirsz+}).
We present another, more elementary proof; using the following analog of Alexandrov's lemma (\ref{lem:alex}).

We say that two triangles with a common vertex \emph{do not overlap} if their convex hulls intersect only at the common vertex.

\begin{thm}{Overlap lemma}\label{lem:extend-overlap}
Let $\trig{\tilde x^1}{\tilde x^2}{\tilde x^3}$ be a triangle in $\Lob2{\kappa}$.
Let $\tilde p^1$, $\tilde p^2$, $\tilde p^3$ be points in $\Lob2{\kappa}$ such that, for any permutation $\{i,j,k\}$ of $\{1,2,3\}$, we have
\begin{enumerate}[(i)]

\item 
\label{no-overlap:px=px}
$\dist{\tilde p^i}{\tilde x^\kay}{}=\dist{\tilde p^j}{\tilde x^\kay}{}$,

\item
\label{no-overlap:orient-1}
$\tilde p^i$ and $\tilde x^i$ lie in the same closed halfspace determined by $[\tilde x^j\tilde x^\kay]$, 
\end{enumerate}

If no pair of triangles $\trig{\tilde p^i}{\tilde x^j}{\tilde x^\kay}$ overlap,
then 
\[\mangle{\tilde p^1} +\mangle {\tilde p^2}+\mangle{\tilde p^3}> 2\cdot\pi,\]
where $\mangle\tilde p^i\df\mangle\hinge{\tilde p^i}{\tilde x^\kay}{\tilde x^j}$
for a permutation $\{i,j,k\}$ of $\{1,2,3\}$.
\end{thm}

{

\begin{wrapfigure}{r}{28 mm}
\vskip-0mm
\centering
\includegraphics{mppics/pic-1005}
\vskip0mm
\end{wrapfigure}

\parbf{Remarks.}
If $\kappa\le 0$, then the overlap lemma can be proved without using condition (\ref{no-overlap:px=px}).
This follows since the sum of external angles for the hexagon
$[\tilde p^1\tilde x^2\tilde p^3\tilde x^1\tilde p^2\tilde x^3]$ and its area is $2\cdot\pi-\kappa\cdot a$, where $a$ denotes the area of the hexagon.

The diagram shows that condition (\ref{no-overlap:px=px}) is essential
in case $\kappa>0$.

}

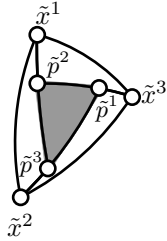
\parit{Proof.} Rotate the triangle $\trig{\tilde p^3}{\tilde x^1}{\tilde x^2}$ around $\tilde x^1$ to make $[\tilde x^1\tilde p^3]$ coincide with $[\tilde x^1\tilde p^2]$.
Let $\dot x^2$ denote the image of $\tilde x^2$ after rotation. 
Note that 
$$\mangle\hinge{\tilde x^1}{\tilde x^3}{\dot x^2}
=
\min\{\,\mangle\hinge {\tilde x^1}{\tilde x^2}{\tilde p^3}
+
\mangle\hinge {\tilde x^1}{\tilde p^2}{\tilde x^3},\,
2\cdot\pi -(\mangle\hinge {\tilde x^1}{\tilde x^2}{\tilde p^3}
+
\mangle\hinge {\tilde x^1}{\tilde p^2}{\tilde x^3})\,\}.
$$
\begin{figure}[!ht]
\vskip-3mm
\centering
\includegraphics{mppics/pic-1010}
\vskip0mm
\end{figure}
By (\ref{no-overlap:orient-1}), 
the triangles 
$\trig{\tilde p^3}{\tilde x^1}{\tilde x^2}$ 
and $\trig{\tilde p^2}{\tilde x^3}{\tilde x^1}$ do not overlap if and only if 
\[\mangle\hinge{\tilde x^1}{\tilde x^3}{\tilde x^2}
> 
\mangle\hinge{\tilde x^1}{\tilde x^3}{\dot x^2}.\eqlbl{eq:main-overlap}\]
and
\[
2\cdot\pi>\mangle\hinge {\tilde x^1}{\tilde x^2}{\tilde p^3}+\mangle\hinge {\tilde x^1}{\tilde p^2}{\tilde x^3}+\mangle\hinge {\tilde x^1}{\tilde x^2}{\tilde x^3}\eqlbl{eq:old-iii}\]
The condition \ref{eq:main-overlap} holds if and only if 
$\dist{\tilde x^2}{\tilde x^3}{}>\dist{\dot x^2}{\tilde x^3}{}$,
which in turn holds if and only if 
\[
\begin{aligned}
\mangle\tilde p^1
&> \mangle\hinge{\tilde p^2}{\tilde x^3}{\dot x^2}
\\
&=
\min\{\mangle\tilde p^3+\mangle\tilde p^2,2\cdot\pi -(\mangle\tilde p^3+\mangle\tilde p^2)\}.
\end{aligned}
\eqlbl{eq:no-overlap}\]
The inequality follows since the corresponding hinges have the same pairs of sidelengths.
(The two pictures show that both possibilities for the minimum can occur.)

Now assume $\mangle\tilde p^1 + \mangle\tilde p^2+\mangle\tilde p^3 \le 2\cdot\pi$.
Then \ref{eq:no-overlap} implies 
\[\mangle\tilde p^i>\mangle\tilde p^j + \mangle\tilde p^\kay.\]
Since no pair of triangles overlap, the same holds 
for any permutation $(i,j,\kay)$ of $(1,2,3)$.
Therefore
\[\mangle\tilde p^1+\mangle\tilde p^2+\mangle\tilde p^3>2\cdot(\mangle\tilde p^1+\mangle\tilde p^2+\mangle\tilde p^3),\]
a contradiction. 
\qeds

\parit{Proof of \ref{thm:kirsz-def}; if part.} 
Assume $\spc{L}$ is geodesic.
Consider $x^1,x^2,x^3\in \spc{L}$ such that the model triangle 
$\trig{\tilde x^1}{\tilde x^2}{\tilde x^3}=\modtrig\kappa(x^1 x^2 x^3)$ is defined.
Choose $p\in \,{]}x^1x^2{[}\,$.
Applying the short map extension property with $V_3=\{x^1,x^2,x^3\}$, $V_4=\{x^1,x^2,x^3,p\}$ and the map $f\:x^i\mapsto\tilde x^i$, we obtain the point-on-side comparison (\ref{point-on-side}).

In case $\spc{L}$ is not geodesic, pass to its ultrapower $\spc{L}^\o$.
Note that the short map extension property survives
for $\spc{L}^\o$ and recall that $\spc{L}^\o$ is geodesic (see \ref{obs:ultralimit-is-complete}).
Thus, from above, $\spc{L}^\o$ is a complete length $\Alex\kappa$ space. 
By Proposition~\ref{prp:A^omega}, $\spc{L}$ is a complete length $\Alex\kappa$ space.

\parit{Only-if part.}
Assume the contrary: 
$\spc{L}$ is complete and $\Alex\kappa$, and 
$x^1,x^2,x^3,p\in \spc{L}$ and 
$\tilde x^1,\tilde x^2,\tilde x^3\in\Lob2\kappa$ are such that
$\dist{\tilde x^i}{\tilde x^j}{}\le\dist{x^i}{x^j}{}$ for all $i,j$ but there is no point $\tilde p\in \Lob2\kappa$ such that $\dist{\tilde p}{\tilde x^i}{}\le \dist{p}{x^i}{}$ for all $i$.

Note that in this case all comparison triangles $\modtrig\kappa(p x^ix^j)$ are defined.
This is always true if $\kappa\le0$.
If $\kappa>0$, and say $\modtrig\kappa(p x^1x^2)$ is undefined, then 
\begin{align*}
\dist{p}{x^1}{}+\dist{p}{x^2}{}
&\ge 2\cdot\varpi\kappa-\dist{x^1}{x^2}{}
\ge
\\
&\ge
2\cdot\varpi\kappa-\dist{\tilde x^1}{\tilde x^2}{}\ge 
\\
&\ge 
\dist{\tilde x^1}{\tilde x^3}{}+\dist{\tilde x^2}{\tilde x^3}{}.
\end{align*}
Then the last inequality must be an equality. 
 Thus we may extend by taking $\tilde p$ on $[\tilde x^1\tilde x^3]$ or $[\tilde x^2\tilde x^3]$.

For each $i\in \{1,2,3\}$, consider a point $\tilde p^i\in\Lob2\kappa$ such that $\dist{\tilde p^i}{\tilde x^i}{}$ is minimal among points satisfying $\dist{\tilde p^i}{\tilde x^j}{}\le\dist{p}{ x^j}{}$ for all $j\ne i$. 
Clearly, every $\tilde p^i$ is inside the triangle $\trig{\tilde x^1}{\tilde x^2}{\tilde x^3}$ (that is, in $\Conv(\tilde x^1,\tilde x^2,\tilde x^3)$), and $\dist{\tilde p^i}{\tilde x^i}{}>\dist{p}{ x^i}{}$ for each $i$.
Since the function $x\mapsto\tangle\mc\kappa\{x;a,b\}$
is increasing, it follows that
\begin{enumerate}[(i)]
\item $\dist{\tilde p^i}{\tilde x^j}{}=\dist{p}{ x^j}{}$ for $i\ne j$;
\item no pair of triangles from $\trig{\tilde p^1}{\tilde x^2}{\tilde x^3}$, $\trig{\tilde p^2}{\tilde x^3}{\tilde x^1}$, $\trig{\tilde p^3}{\tilde x^1}{\tilde x^2}$ overlap in $\trig{\tilde x^1}{\tilde x^2}{\tilde x^3}$.
\end{enumerate}

As follows from the overlap lemma (\ref{lem:extend-overlap}), 
in this case 
\[\mangle\hinge {\tilde p^1}{\tilde x^2}{\tilde x^3} 
+\mangle\hinge {\tilde p^2}{\tilde x^3}{\tilde x^1}
+\mangle\hinge {\tilde p^3}{\tilde x^1}{\tilde x^2}
>2\cdot\pi.
\]
Since $\dist{\tilde x^i}{\tilde x^j}{}\le\dist{x^i}{x^j}{}$ we have
\[\mangle\hinge {\tilde p^\kay}{\tilde x^i}{\tilde x^j}
\le
\angk\kappa p{x^i}{x^j}\]
if $(i,j,k)$ is a permutation of $(1,2,3)$.
Therefore 
\[\angk\kappa p{x^1}{x^2}+\angk\kappa p{x^2}{x^3}+\angk\kappa p{x^3}{x^1}>2\cdot\pi,\]
contradicting the $\Alex\kappa$ comparison (\ref{df:cbb1+3}).
\qeds

\begin{thm}{Theorem}\label{thm:cba-kirsz-def} 
Assume any pair of points at distance $<\varpi\kappa$ in the metric space $\spc{U}$ are joined by a unique geodesic. 
Then $\spc{U}$ is $\CAT\kappa$ if and only if 
for any $3$-point set $V_3$ with perimeter $<2\cdot\varpi\kappa$
and any $4$-point set $V_4\supset V_3$ in $\Lob2\kappa$,
any short map $f\:V_3\to\spc{U}$ can be extended to a short map $F\:V_4\to\spc{U}$.
\end{thm}

Note that the only-if part of Theorem \ref{thm:cba-kirsz-def} does not follow directly from Kirszbraun's theorem, since the desired extension is in $\spc{U}$ --- not its completion.

\begin{thm}{Lemma}\label{lem:smaller-trig}
Let $x^1,x^2,x^3,y^1,y^2,y^3\in\Lob{}{\kappa}$
be points such that $\dist{x^i}{x^j}{}\ge\dist{y^i}{y^j}{}$ for all $i,j$.
Then there is a short map $\map\:\Lob{}{\kappa}\to\Lob{}{\kappa}$ such that $\map(x^i)=y^i$ for all $i$;
moreover, one can choose $\map$ so that 
\[\Im \map\subset\Conv(y^1,y^2,y^3).\]

\end{thm}

We only give an idea of the proof of this lemma;
alternatively, it can be obtained as a corollary of Kirszbraun's theorem (\ref{thm:kirsz+}) 

\parit{Idea of the proof.}
The map $\map$ can be constructed as a composition of an isometry of $\Lob{}{\kappa}$ and the following folding map:
given a halfspace $H$ in $\Lob{}{\kappa}$, consider the map $\Lob{}{\kappa}\to H$ 
that is the identity on $H$ and reflects all points outside of $H$ into $H$.
This map is a path isometry; in particular, it is short. 

The last part of the lemma can be proved by composing this map with folding maps along the sides of triangle $\trig{y^1}{y^2}{y^3}$, and passing to a partial limit.
\qeds

\parit{Proof of \ref{thm:cba-kirsz-def}; if part.}
The point-on-side comparison (\ref{cat-monoton}) follows by
taking $V_3=\{\tilde x,\tilde y,\tilde p\}$ and $V_4=\{\tilde x,\tilde y,\tilde p,\tilde z\}$ where $z\in \mathopen{]}x y\mathclose{[}$. 
It is only necessary to observe that $F(\tilde z)=z$ by uniqueness of $[x y]$.

\parit{Only-if part.}
Let $V_3=\{\tilde x^1,\tilde x^2,\tilde x^3\}$ and $V_4=\{\tilde x^1,\tilde x^2,\tilde x^3,\tilde p\}$.

Set $y^i\z=f(\tilde x^i)$ for all $i$.
We need to find a point $q\in\spc{U}$ such that $\dist{y^i}{q}{}\le\dist{\tilde x^i}{\tilde p}{}$ for all $i$.

Let $D$ be the convex set in $\Lob2\kappa$ bounded by the model triangle 
$\trig{\tilde y^1}{\tilde y^2}{\tilde y^3}\z=\modtrig\kappa{y^1}{y^2}{y^3}$;
that is, $D\z=\Conv({\tilde y^1},{\tilde y^2},{\tilde y^3})$.

Note that $\dist{\tilde y^i}{\tilde y^j}{}=\dist{y^i}{y^j}{}\le\dist{\tilde x^i}{\tilde x^j}{}$ for all $i,j$.
Applying Lemma \ref{lem:smaller-trig},
we get a short map 
$\map\:\Lob{}{\kappa}\to D$ such that 
$\map\:\tilde x^i\mapsto\tilde y^i$.

Further, by the majorization theorem (\ref{thm:major}), 
there is a short map $F\:D\to \spc{U}$ such that $\tilde y^i\mapsto y^i$ for all $i$.

Thus one can take $q=F\circ\map(\tilde p)$.
\qeds

\begin{thm}{Exercise}\label{ex:CBB+CBA}
Assume $\spc{X}$ is a complete length space that satisfies the following condition:
any 4-point subset admits a distance-preserving map to the Euclidean 3-space.

Prove that $\spc{X}$ is isometric to a closed convex subset of a Hilbert space.
\end{thm}

\begin{thm}{Exercise}\label{ex:5-point-CBA=>CBB}
Let $\spc{F}_s$ be the metric
on the 5-point set $\{p,q,x,y,z\}$ for which $\dist{p}{q}{}=s$
and all the remaining distances are equal $1$.
For which values $s$ does the space $\spc{F}_s$ admit a distance-preserving map into 
\begin{subthm}{}
a complete length $\CAT{0}$ space?
\end{subthm}
\begin{subthm}{}
a complete length $\Alex{0}$ space?
\end{subthm}
\end{thm}

The following exercise describes the first known definition of spaces with curvature bounded below;
it was given by Abraham Wald \cite{wald}. 

\begin{thm}{Exercise}\label{ex:cbb-wald}
Let $\spc{L}$ be a metric space and $\kappa\le 0$.
Prove that $\spc{L}$ is $\Alex\kappa$ if and only if any quadruple of points $p,q,r,s\in \spc{L}$ admits a distance-preserving embedding into $\Lob2{\Kappa}$ for some $\Kappa\ge\kappa$.

Is the same true for $\kappa>0$; what is the difference?
\end{thm}

\section{(1+\textit{n})-point comparison}\label{sec:1+n}

The following theorem gives a more sensitive analog of the $\Alex\kappa$ comparison (\ref{df:cbb1+3}).

\begin{thm}{(1+\textit{n})-point comparison}
\label{thm:pos-config} {\sloppy 
Let $\spc{L}$ be a complete length $\Alex{\kappa}$ space.
Then for any array $(p,x^1,\dots,x^n)$ of points in $\spc{L}$ 
there is a model array $(\tilde p,\tilde x^1,\dots,\tilde x^n)$ in $\Lob{n}\kappa$ such that

}

\begin{subthm}{}
$\dist{\tilde p}{\tilde x^i}{}=\dist{p}{x^i}{}$ for all $i$.
\end{subthm}

\begin{subthm}{}$\dist{\tilde x^i}{\tilde x^j}{}\ge\dist{x^i}{x^j}{}$ for all $i,j$.
\end{subthm}
\end{thm}

\parit{Proof.} 
It is sufficient to show that given $\eps>0$ there is an array $(\tilde p,\tilde x^1,\dots,\tilde x^n)$ in $\Lob{n}\kappa$ such that 
\[\dist{\tilde x^i}{\tilde x^j}{}\ge\dist{x^i}{x^j}{}
\quad
\text{and}
\quad
\dist{\tilde p}{\tilde x^i}{}\lege\dist{p}{x^i}{}\bigr|\pm\eps.\]
Then one can pass to a limit array for $\eps\to 0+$.

According to \ref{thm:almost.geod}, the set $\Str(x^1,\dots,x^n)$ 
 is dense in $\spc{L}$.
Thus there is a point $p'\in \Str(\tilde x^1,\dots,\tilde x^n)$ such that $\dist{p'}{p}{}\le\eps$.
According to Corollary~\ref{cor:euclid-subcone}, 
$\T_{p'}$ contains a subcone $E$ isometric to a Euclidean space 
and containing all vectors $\ddir{p'}{x^i}$.
Passing to a subspace if necessary, we may assume that $\dim E\le n$.

Mark a point $\tilde p\in \Lob{n}\kappa$ and choose a distance-preserving map
$\iota\: E\to \T_{\tilde p}\Lob{n}\kappa$.
Let 
\[\tilde x^i
=
\exp_{\tilde p}\circ\iota(\ddir{p'}{x^i}).\]
Thus $\dist{\tilde p}{\tilde x^i}{}=\dist{p'}{x^i}{}$.
Since $\dist{p}{p'}{} \le\eps$, we get
\[\dist{\tilde p}{\tilde x^i}{}\lege\dist{p}{x^i}{}
\pm
\eps.\]
From the hinge comparison (\ref{angle}) 
we have 
\[\angk\kappa{\tilde p}{\tilde x^i}{\tilde x^j}
=\mangle\hinge{\tilde p}{\tilde x^i}{\tilde x^j}
=\mangle\hinge{p'}{x^i}{x^j}\ge \angk\kappa{p'}{x^i}{x^j},\]
and thus 
\[\dist{\tilde x^i}{\tilde x^j}{}\ge \dist{x^i}{x^j}{}.\]
\qedsf

\begin{thm}{Exercise}\label{ex:sturm}
Let $(p,x_1,\dots,x_n)$ be a point array in a $\Alex0$ space.
Consider the $n{\times}n$-matrix $M$ with components 
\[m_{i,j}=\tfrac12\cdot(\dist[2]{x_i}{p}{}+\dist[2]{x_j}{p}{}-\dist[2]{x_i}{x_j}{}).\]
Show that
\[\bm{s}\cdot M\cdot \bm{s}^\top\ge 0\eqlbl{eq:sMs}\]
for any vector $\bm{s}=(s_1,\dots,s_n)$ with nonnegative components.
\end{thm}

\begin{wrapfigure}{r}{23 mm}
\vskip-3mm
\centering
\includegraphics{mppics/pic-1015}
\vskip0mm
\end{wrapfigure}

The above exercise describes the so-called \index{Lang--Schroeder--Sturm inequality}\emph{Lang--Schroeder--Sturm inequality}; it was discovered by Urs Lang and Viktor Schroeder \cite{lang-schroeder} and rediscovered by Karl-Theodor Sturm \cite{sturm}.
It turns out to be weaker than (1+\textit{n})-point comparison.
An example can be constructed by perturbing the 6-point metric isometric to a regular pentagon with its center, making its sides slightly longer and diagonals slightly shorter \cite{lebedeva-petrunin-zolotov}.
In particular, this inequality in general metric spaces (not necessarily length spaces) does not imply the inequality in the following exercise.

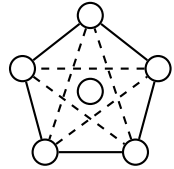
\begin{thm}{Exercise}\label{6-point-comparison}
Let $\spc{L}$ be a complete length $\Alex{\kappa}$ space.
Show that for any points $p, x_1,x_2,x_3,x_4,x_5 $ in $\spc{L}$ we have 
\[
\angk\kappa p{x_1}{x_5}
+
\angk\kappa p{x_2}{x_1}
+
\angk\kappa p{x_3}{x_2}
+
\angk\kappa p{x_4}{x_3}
+
\angk\kappa p{x_5}{x_4}\le 4\cdot \pi,
\]
assuming that the left-hand side is defined.
\end{thm}

\begin{thm}{Exercise}\label{ex:(3+1)-nonsufficient}
Give an example of a metric on a finite set that satisfies the comparison inequality 
\[\angk{0}{p}{x_1}{x_2}+\angk{0}{p}{x_2}{x_3}+\angk{0}{p}{x_3}{x_1}
\le
2\cdot\pi\]
for any quadruple of points $(p,x_1,x_2,x_3)$, 
but is not isometric to a subset of an Alexandrov space with curvature $\ge0$.
\end{thm}

\begin{thm}{Exercise}\label{ex:strut+embedding}
Let $\spc{L}$ be a complete length $\Alex{\kappa}$ space. 
Assume that a point array $(a^0,a^1,\dots,a^\kay)$ in $\spc{L}$
 is \emph{$\kappa$-strutting} (Definition \ref{def:strut-I})
 for a point $p\in\spc{L}$.
Show that there are points
$\tilde p,\tilde a^0,\dots,\tilde a^m$ in $\Lob{m+1}\kappa$ such that
\[\dist{\tilde p}{\tilde a^i}{}=\dist{p}{a^i}{}\quad \text{and}\quad \dist{\tilde a^i}{\tilde a^j}{}=\dist{a^i}{a^j}{}\]
for all $i$ and $j$.
\end{thm}

\section{Helly's theorem}\label{sec:helly}

\begin{thm}{Helly's theorem}\label{thm:helly}
Let $\spc{U}$ be a complete length $\CAT0$ space
and $\{K_\alpha\}_{\alpha\in \IndexSet}$ be an arbitrary collection of closed bounded convex subsets in~$\spc{U}$.

If 
\[\bigcap_{\alpha\in \IndexSet}K_\alpha=\emptyset,\]
then there is a finite index array $(\alpha_1,\alpha_2,\dots,\alpha_n)$ in $\IndexSet$ such that
\[\bigcap_iK_{\alpha_i}=\emptyset.\]
\end{thm}

\parbf{Remarks.}
\begin{enumerate}[(i)]
\item In general, none of the $K_\alpha$ may be compact; 
otherwise the statement is trivial.
\item If $\spc{U}$ is a Hilbert space (not necessarily separable), 
then Helly's theorem is equivalent to the following statement: if a convex bounded set is closed in the ordinary topology then it is compact in the weak topology.
One can define {}\emph{weak topology} in an arbitrary metric space by taking exteriors of closed ball as prebase.
Then 
Helly's theorem
 implies the analogous statement for complete length $\CAT0$ spaces
(compare to \cite{monod}).
\end{enumerate}

\medskip

We present the proof of Urs Lang and Viktor Schroeder \cite{lang-schroeder}.

\parit{Proof of \ref{thm:helly}.} 
Assume the contrary. Then for any finite set $F\subset \IndexSet$, 
\[K_{F}\df \bigcap_{\alpha\in F}K_{\alpha}\ne \emptyset.\]
We will construct a point $z$ such that $z\in K_\alpha$ for each $\alpha$.
Thus we will arrive at a contradiction since
\[\bigcap_{\alpha\in \IndexSet}K_\alpha=\emptyset.\]

Choose a point $p\in \spc{U}$, and let $r=\sup\{\distfun{K_{F}}{p}{}\}$ where $F$ runs over all finite subsets of $\IndexSet$.
Let $p^*_F$ be the closest point on $K_{F}$ to $p$; 
according to the closest-point projection lemma (\ref{lem:closest point}), $p^*_F$ 
exists and is unique.

Take a nested sequence of finite subsets 
$F_1\subset F_2\subset \dots$ of $\IndexSet$, such that $\distfun{K_{F_n}}{p}{}\to r$.

Let us show that the sequence $p^*_{F_n}$ is Cauchy. 
If not, then for fixed $\eps>0$, 
we can choose two subsequences $y'_n$ and $y''_n$ of $p^*_{F_n}$ 
such that $\dist{y'_n}{y''_n}{}\ge\eps$.
Let $z_n$ be the midpoint of $[y'_ny''_n]$. 
From the point-on-side comparison (\ref{point-on-side}), 
there is $\delta>0$ such that 
\[\dist{p}{z_n}{}\le \max\{\dist{p}{y'_n}{},\dist{p}{y''_n}{}\}-\delta.\]
Thus 
\[\limsup_{n\to\infty}\dist{p}{z_n}{}<r.\]
On the other hand, from convexity, each $K_{F_n}$ 
contains all $z_\kay$ with sufficiently large $\kay$, a contradiction.

Thus, $p^*_{F_n}$ converges and we can take $z=\lim_n p^*_{F_n}$.
Clearly 
\[\dist{p}{z}{}=r.\]

Repeat the above arguments for the sequence $F_n'=F_n\cup \{\alpha\}$.
As a result, we get another point $z'$ such that $\dist{p}{z}{}=\dist{p}{z'}{}=r$ and 
$z,z'\in K_{F_n}$ for all $n$.
Thus, if $z\ne z'$ the midpoint $\hat z$ of $[zz']$ would belong to all 
$K_{F_n}$, and from comparison, we would have $\dist{p}{\hat z}{}<r$, a contradiction.

Thus, $z'=z$; in particular 
$z\in K_\alpha$ for each $\alpha\in\IndexSet$.
\qeds

\section{Kirszbraun's theorem}\label{sec:kirszbraun}

A slightly weaker version of the following theorem was proved by Urs Lang and Viktor Schroeder \cite{lang-schroeder}.

\begin{thm}{Kirszbraun's theorem}
\label{thm:kirsz+}
Let
$\spc{L}$ be a complete length $\Alex{\kappa}$ space, 
$\spc{U}$ be a complete length $\CAT\kappa$ space, 
$Q\subset \spc{L}$ be arbitrary subset
and $f\: Q\to\spc{U}$ be a short map.
Assume that there is $z\in\spc{U}$ such that 
$f(Q)\subset \oBall[z,\tfrac{\varpi\kappa}{2}]_{\spc{U}}$.
Then $f\:Q\to\spc{U}$ can be extended to a short map 
$F\:\spc{L}\to \spc{U}$
(that is, there is a short map $F\:\spc{L}\to \spc{U}$ such that $F|_Q=f$).
\end{thm}
 
The condition $f(Q)\subset \oBall[z,\tfrac{\varpi\kappa}{2}]$ trivially holds for any $\kappa\le 0$ since in this case $\varpi\kappa=\infty$. 
The following example shows that this condition is needed for $\kappa>0$.

Conjecture~\ref{conj:kirsz} (if true) gives an equivalent condition for the existence of a short extension;
it states that the following example is the only obstacle.

\begin{thm}{Example}\label{example:SS_+}
Let $\mathbb{S}^m_+$ be a closed $m$-dimensional unit hemisphere. Denote its boundary, which is isometric to $\mathbb{S}^{m-1}$, by $\partial\mathbb{S}^m_+$.
Clearly, $\mathbb{S}^m_+$ is $\Alex{1}$ and $\partial\mathbb{S}^m_+$ is $\CAT1$, but the identity map ${\partial\mathbb{S}^m_+}\to \partial\mathbb{S}^m_+$ cannot be extended to a short map $\mathbb{S}^m_+\to \partial\mathbb{S}^m_+$ (there is no place for the pole).

There is also a direct generalization of this example to a hemisphere in a Hilbert space of arbitrary cardinal dimension.
\end{thm}

First we prove this theorem in the case $\kappa\le 0$ (\ref{thm:kirsz}).
In the proof of the more complicated case $\kappa>0$, we use the case $\kappa=0$.
The following lemma is the main ingredient in the proof. 

\begin{thm}{Finite$\bm{+}$one lemma}\label{lem:kirsz-neg:new}
Let $\kappa\le 0$,
$\spc{L}$ be a complete length $\Alex{\kappa}$ space, and 
$\spc{U}$ be a complete length $\CAT\kappa$ space.
Suppose
$x^1,x^2,\dots,x^n$ in $\spc{L}$ 
and $y^1,y^2,\dots,y^n$ in $\spc{U}$
are
such that $\dist{x^i}{x^j}{}\ge\dist{y^i}{y^j}{}$ for all $i,j$.

Then for any $p\in\spc{L}$, there is $q\in\spc{U}$ such that $\dist{y^i}{q}{}\le\dist{x^i}{p}{}$ for each $i$.
\end{thm}

\parit{Proof.}
It is sufficient to prove the lemma only for $\kappa=0$ and $-1$.
The proofs of these two cases are identical, only the formulas differ.
In the proof, we assume $\kappa=0$ and provide the formulas for $\kappa=-1$ in the footnotes.

From the (1+\textit{n})-point comparison (\ref{thm:pos-config}), 
there is a model configuration 
$\tilde p,\tilde x^1,\tilde x^2,\dots,\tilde x^n$ in $\Lob{n}{\kappa}$ such that
$\dist{\tilde p}{\tilde x^i}{}=\dist{p}{x^i}{}$
and $\dist{\tilde x^i}{\tilde x^j}{}\ge\dist{x^i}{x^j}{}$ 
for all $i$, $j$.
It follows that we can assume that $\spc{L}=\Lob{n}{\kappa}$.

For each $i$, consider functions 
$f^i\:\spc{U}\to\RR$ and $\tilde f^i\:\Lob{n}{\kappa}\to\RR$ 
defined as follows:%
\footnote{In case $\kappa=-1$,
\[
\begin{aligned}
&f^i=\cosh\circ\distfun{y^i}{}{},
&
&\tilde f^i=\cosh\circ\distfun{\tilde x^i}{}{}.
\end{aligned}
\leqno{(A)\mc-}\]}
\[
\begin{aligned}
&f^i=\tfrac{1}{2}\cdot\distfun[2]{y^i}{}{},
&
&\tilde f^i=\tfrac{1}{2}\cdot\distfun[2]{\tilde x^i}{}{}.
\end{aligned}
\leqno{(A)\mc0}
\]
Consider the function arrays
$\bm{f}=(f^1,f^2,\dots,f^n)\:\spc{U}\to\RR^n$ and $\bm{\tilde f}\z=(\tilde f^1,\tilde f^2,\dots,\tilde f^n)\:\Lob{n}{\kappa}\to\RR^n$.

Define
\begin{align*}
\Up \bm{f}(\spc{U})
&=
\set{\bm{v}\in\RR^{\kay+1}}{\exists\, \bm{w}\in \bm{f}(\spc{U})\ \text{such that}\ \bm{v}\succcurlyeq\bm{w}},
\\
\Min \bm{f}(\spc{U}) 
&=
\set{\bm{v}\in \bm{f}(\spc{U})}{\text{if}\ \bm{v}\succcurlyeq\bm{w}\in \bm{f}(\spc{U})\ \text{then}\ \bm{w}=\bm{v}}.
\end{align*}
(See Definition \ref{def:ordung}.)
Note it is sufficient to prove that $\bm{\tilde f}(\tilde p)\in\Up \bm{f}(\spc{U})$.

Clearly,
$(f^i)''\ge 1$.
Thus by Theorem \ref{thm:web:Up-convex}, 
the set $\Up\bm{f}(\spc{U})\subset\RR^{n}$ is convex.

Arguing by contradiction, let us assume that $\bm{\tilde f}(\tilde p)\notin\Up \bm{f}(\spc{U})$.

Then there exists a supporting hyperplane $\alpha_1\cdot x_1+\ldots \alpha_n\cdot x_n=c$ to $\Up\bm{f}(\spc{U})$, separating it from $\bm{\tilde f}(\tilde p)$.
According to Lemma~\ref{lem:Up-convex:subnormal}, 
$\alpha_i\ge 0$ for each $i$. 
So we may assume that 
$(\alpha_1,\alpha_2,\dots,\alpha_n)\in\Delta^{n-1}$
(that is, 
$\alpha_i\ge 0$ for each $i$ and $\sum\alpha_i=1$
and 
\[\sum_i\alpha_i\cdot\tilde f^i(\tilde p))
< 
\inf
\set{\sum_i\alpha_i\cdot f^i(q)}{q\in\spc{U}}.\]
The latter contradicts the following claim.

\begin{clm}{}
Given $\bm{\alpha}=(\alpha_1,\alpha_2,\dots,\alpha_n)\in\Delta^{n-1}$,
let
\begin{align*}
&h=\sum_i\alpha_i\cdot f^i
&
&h\:\spc{U}\to\RR
&
&z=\argmin h\in \spc{U}
\\
&\tilde h=\sum_i\alpha_i\cdot \tilde f^i
&
&\tilde h\:\Lob{n}{\kappa}\to\RR
&
&\tilde z=\argmin \tilde h\in \Lob{n}{\kappa}.
\end{align*}
Then 
$h(z)\le \tilde h(\tilde z)$.
\end{clm}

\parit{Proof of the claim.}
Note that $\dd_z h\ge 0$.
Thus, for each $i$, we have%
\footnote{In case $\kappa=-1$, the same calculations give
\[
\begin{aligned}
0
&\le\dots \le
-\tfrac{1}{\sinh\dist[{{}}]{z}{y^i}{}}
\cdot 
\sum_j
\alpha_j\cdot\left[\cosh\dist[{{}}]{z}{y^i}{}\cdot\cosh\dist[{{}}]{z}{y^j}{}-\cosh\dist[{{}}]{y^i}{y^j}{}\right].
\end{aligned}
\leqno{(B)\mc-}
\]

}

\[
\begin{aligned}
0
&\le (\dd_z h)(\dir{z}{y^i})
=
\\
&=
-\sum_j\alpha_j\cdot\dist[{{}}]{z}{y^j}{}\cdot\cos\mangle\hinge{z}{y^i}{y^j}
\le
\\
&\le
-\sum_j\alpha_j\cdot\dist[{{}}]{z}{y^j}{}\cdot\cos\angk0{z}{y^i}{y^j}
=
\\
&=
-\tfrac{1}{2\cdot\dist[{{}}]{z}{y^i}{}}
\cdot 
\sum_j
\alpha_j\cdot\left[\dist[2]{z}{y^i}{}+\dist[2]{z}{y^j}{}-\dist[2]{y^i}{y^j}{}\right].
\end{aligned}
\leqno{(B)\mc0}\]
In particular%
\footnote{In case $\kappa=-1$, the same calculations give
\[
\begin{aligned} 
\sum_{i}\alpha_i\cdot\left[\sum_j
\alpha_j\cdot\left[\cosh\dist[{{}}]{z}{y^i}{}\cdot\cosh\dist[{{}}]{z}{y^j}{}
-\cosh\dist[{{}}]{y^i}{y^j}{}\right]
\right]\le0
\end{aligned}.
\leqno{(C)\mc-}
\]
},
\[
\begin{aligned}
\sum_{i}
\alpha_i
\cdot
\left[\sum_j
\alpha_j
\cdot
\left[\dist[2]{z}{y^i}{}+\dist[2]{z}{y^j}{}-\dist[2]{y^i}{y^j}{}\right]
\right]\le 0,
\end{aligned}
\leqno{(C)\mc0}
\]
or%
\footnote{In case $\kappa=-1$,
\[(h(z))^2\le
\sum_{i,j}
\alpha_i\cdot\alpha_j
\cdot
\cosh\dist[{{}}]{y^i}{y^j}{}. \leqno{(D)\mc-}\]
}
\[2\cdot h(z)
\le
\sum_{i,j}
\alpha_i\cdot\alpha_j
\cdot
\dist[2]{y^i}{y^j}. \leqno{(D)\mc0}\]

Note that if $\spc{U}\iso\Lob{n}{\kappa}$, 
then all inequalities in $(B,C,D)$ are sharp.
Thus the same argument as above, repeated for $\tilde x^1,\tilde x^2,\dots,\tilde x^n$ in $\Lob{n}{\kappa}$,
gives%
\footnote%
{In case $\kappa=-1$,
\[(\tilde h(\tilde z))^2
=
\sum_{i,j}
\alpha_i\cdot\alpha_j
\cdot
\cosh\dist[{{}}]{\tilde x^i}{\tilde x^j}{}.
\leqno{(E)\mc-}\]
}
\[
2\cdot \tilde h(\tilde z)
=
\sum_{i,j}
\alpha_i\cdot\alpha_j
\cdot
\dist[2]{\tilde x^i}{\tilde x^j}{}. 
\leqno{(E)\mc0}
\]
Note that 
\[\dist{\tilde x^i}{\tilde x^j}{}
\ge
\dist{x^i}{x^j}{}\ge\dist{y^i}{y^j}{}\]
for all $i$, $j$.
Thus, $(D)$ and $(E)$ imply the claim.
\qedqeds

\begin{thm}{Kirszbraun's theorem for nonpositive bound}
\label{thm:kirsz}
Let
$\kappa\le0$,
$\spc{L}$ be a complete length $\Alex{\kappa}$ space, 
$\spc{U}$ be a complete length $\CAT\kappa$ space, 
$Q\subset \spc{L}$ be arbitrary subset
and $f\: Q\to\spc{U}$ be a short map.
Then there is a short extension 
$F\:\spc{L}\to \spc{U}$ of $f$;
that is, there is a short map $F\:\spc{L}\to \spc{U}$ such that $F|_Q=f$.
\end{thm}

\parbf{Remark.}
If $\spc{U}$ is proper, then we do not need Helly's theorem (\ref{thm:helly}); compactness of closed balls in $\spc{U}$ is sufficient in this case.

\parit{Proof of \ref{thm:kirsz}.} 
By Zorn's lemma, we can assume 
that $Q\subset\spc{L}$ is a maximal set;
that is, $f\:Q\to\spc{U}$ does not admit a short extension to any larger set $Q'\supset Q$.

Let us argue by contradiction.
Assume that $Q\ne \spc{L}$; 
choose $p\in \spc{L}\setminus Q$.
Then
\[\bigcap_{x\in Q} \cBall[f(x),\dist{p}{x}{}]
=
\emptyset.\]

Since $\kappa\le 0$, the balls are convex; 
thus, by Helly's theorem (\ref{thm:helly}), 
one can choose points $x^1,x^2,\dots, x^n$ in $Q$ such that
\[\bigcap_i \cBall[y^i,\dist{x^i}{p}{}]
=
\emptyset,
\eqlbl{eq:cap=cBalls=0}\]
where $y^i=f(x^i)$.
Finally note that \ref{eq:cap=cBalls=0} contradicts the finite+one lemma (\ref{lem:kirsz-neg:new}).\qeds

\parit{Proof of Kirszbraun's theorem (\ref{thm:kirsz+}).} 
The case $\kappa\le 0$ is already proved in \ref{thm:kirsz}.
Thus it remains to prove the theorem only in case $\kappa>0$.
After rescaling we can assume that $\kappa=1$
and therefore $\varpi\kappa=\pi$.

Since $\cBall[z,\pi/2]_{\spc{U}}$ is a complete length $\CAT\kappa$ space, we can assume $\spc{U}=\cBall[z,\pi/2]_{\spc{U}}$. 
In particular, $\diam\spc{U}\le\pi$.

Further, any two points $x,y\in \spc{U}$ such that $\dist{x}{y}{}<\pi$ are joined by a unique geodesic;
if $\dist{x}{y}{}=\pi$, then the concatenation of 
$[x z]$ and $[z y]$ as a geodesic from $x$ to $y$.
Hence $\spc{U}$ is geodesic.

We may also assume that $\diam\spc{L}\le\pi$.
Otherwise $\spc{L}$ is one-dimensional (see \ref{diam-k>0});
in this case the result follows since $\spc{U}$ is geodesic.

\medskip

Assume the theorem is false. Then 
there is a set $Q\subset \spc{L}$, 
a short map $f\: Q\to \spc{U}$, and 
$p\in \spc{L}\setminus Q$ such that 
\[\bigcap_{x\in Q}
\cBall[f(x),\dist{x}{p}{}]=\emptyset.
\eqlbl{eq:cap-of-balls}\]

We will apply \ref{thm:kirsz} for $\kappa=0$ to the Euclidean cones $\mathring{\spc{L}}=\Cone \spc{L}$ and $\mathring{\spc{U}}\z=\Cone \spc{U}$. 
Note that 
\begin{itemize}
\item ${\mathring{\spc{U}}}$ is a complete length $\CAT0$ space (see \ref{thm:warp-curv-bound:cbb:S}),
\item since $\diam \spc{L}\le \pi$ we have ${\mathring{\spc{L}}}$ is $\Alex{0}$ (see \ref{thm:warp-curv-bound:cbb:a}).
\end{itemize}
Further, we will view the spaces $\spc{L}$ and $\spc{U}$ as unit spheres in $\mathring{\spc{L}}$ and $\mathring{\spc{U}}$ respectively.
In the cones $\mathring{\spc{L}}$ and $\mathring{\spc{U}}$ we will use 
``$|{*}|$'' for distance to the tip, denoted by $\0$, 
``$\cdot$'' for cone multiplication,
``$\mangle(x,y)$'' for $\mangle\hinge{\0}{x}{y}$, 
and ``$\<x,y\>$'' for $|x|\cdot|y|\cdot\cos\mangle\hinge{\0}{x}{y}$.
In particular,
\begin{itemize}
\item $\dist{x}{y}{\spc{L}}=\mangle(x,y)$ for any $x,y\in\spc{L}$,
\item $\dist{x}{y}{\spc{U}}=\mangle(x,y)$ for any $x,y\in\spc{U}$,
\item for any $y\in \spc{U}$, we have
\[\mangle(z,y)\le\tfrac\pi2.
\eqlbl{eq:=<pi/2}\]

\end{itemize}
Let $\mathring{Q}=\Cone Q\subset \mathring{\spc{L}}$ and let $\mathring f\:\mathring{Q}\to \mathring{\spc{U}}$ be the natural cone extension of~$f$; 
that is, 
$y=f(x)$ $\Rightarrow$ $t\cdot y=\mathring f(t\cdot x)$ 
for $t\ge0$.
Clearly $\mathring f$ is short.

\begin{wrapfigure}{r}{60mm}
\vskip-0mm
\centering
\includegraphics{mppics/pic-1020}
\vskip0mm
\end{wrapfigure}

Applying \ref{thm:kirsz} for $\mathring f$, 
we get a short extension map $\mathring F\:\mathring{\spc{L}}\to\mathring{\spc{U}}$. 
Let $s=\mathring F(p)$.
Then 
\[\dist{s}{\mathring f(w)}{}
\le 
\dist{p}{w}{}
\eqlbl{eq:clm:kirszbraun-curv=1-rad-star}\]
for any $w\in \mathring Q$.
In particular, $|s|\le 1$.
Applying \ref{eq:clm:kirszbraun-curv=1-rad-star} 
for $w=t\cdot x$ and $t\to\infty$ we have
\[\<f(x),s\>\ge \cos\mangle(p,x)\eqlbl{eq:<,>=<}\]
for any $x\in Q$.

By comparison,
the geodesics $\geod_{[s\ t\cdot z]}$ converge as $t\to\infty$
and the limit is a half-line;
denote it by $\alpha\:[0,\infty)\to \mathring{\spc{U}}$.
From \ref{eq:=<pi/2}, 
the function $t\mapsto\<f(x),\alpha(t)\>$ is nondecreasing. 
From \ref{eq:<,>=<}, for
the necessarily unique point $\bar s$ on the half-line $\alpha$ such that $|\bar s|=1$, we also have 
\[\<f(x),\bar s\>\ge \cos\mangle(p,x)
\quad\text{or, equivalently,}\quad
\mangle(\bar s,f(x))
\le 
\mangle(p,f(x))
\]
for any $x\in Q$,
in contradiction to \ref{eq:cap-of-balls}.
\qeds

\begin{thm}{Exercise}\label{ex:flat-in-CAT}
Let $\spc{U}$ be $\CAT0$. 
Assume there are two point arrays, $(x^0,x^1,\dots,x^\kay)$ in $\spc{U}$ and $(\tilde x^0,\tilde x^1,\dots,\tilde x^\kay)$ in $\EE^m$, such that 
$\dist{x^i}{x^j}{\spc{U}}\z=\dist{\tilde x^i}{\tilde x^j}{\EE^m}$ for each $i,j$, and 
for any point $z_0\in\spc{U}$ there is $i>0$ such that $\dist{z_0}{x_i}{}\ge\dist{x_0}{x_i}{}$.

Prove that there is a subset $Q\subset\spc{L}$ isometric to a convex set in $\EE^m$ and containing all the points $x^i$.
\end{thm}

\begin{thm}{Exercise}\label{ex:flat-in-CBB}
Let
$(x^0,x^1,\dots,x^\kay)$ in $\spc{L}$ and $(\tilde x^0,\tilde x^1,\dots,\tilde x^\kay)$ in $\EE^m$
be two point arrays in complete length $\Alex{0}$ space $\spc{L}$.
Assume that 
$\dist{x^i}{x^j}{\spc{L}}=\dist{\tilde x^i}{\tilde x^j}{\EE^m}$ for each $i,j$
and
$\tilde x^0$ lies in the interior of $\Conv(\tilde x^1,\dots,\tilde x^\kay)$.

Prove that there is a subset $Q\subset\spc{L}$ isometric to a convex set in $\EE^m$ and containing all the points $x^i$.
\end{thm}


The following statement we call \emph{(2\textit{n}+2)-point comparison}.

\begin{thm}{Exercise}\label{CBA-n-point}
Let $\spc{U}$ be a complete length $\CAT\kappa$ space.
Consider $x,y\in \spc{U}$ and an array $( (p^1,q^1)$, $(p^2,q^2),\dots,(p^n,q^n) )$ of pairs of points in $\spc{U}$, such that there is a model configuration
$\tilde x$, $\tilde y$ and array of pairs $( (\tilde p^1,\tilde q^1)$, $(\tilde p^2,\tilde q^2),\dots,(\tilde p^n,\tilde q^n) )$ in $\Lob{3}\kappa$ with the following properties:
\begin{subthm}{}
$\trig{\tilde x}{\tilde p^1}{\tilde q^1}=\modtrig\kappa x p^1q^1$
and 
$\trig{\tilde y}{\tilde p^n}{\tilde q^n}=\modtrig\kappa y p^n q^n$,
\end{subthm}

\begin{subthm}{}
the simplex $\tilde p^i\tilde p^{i+1}\tilde q^i\tilde q^{i+1}$ is a model simplex
 of $p^ip^{i+1}q^iq^{i+1}$
for all $i$;
that is,
$\dist{\tilde p^i}{\tilde q^i}{}
=\dist{p^i}{q^i}{}$,
$\dist{\tilde p^i}{\tilde p^{i+1}}{}
=\dist{p^i}{p^{i+1}}{}$,
$\dist{\tilde q^i}{\tilde q^{i+1}}{}
=\dist{q^i}{q^{i+1}}{}$,
$\dist{\tilde p^i}{\tilde q^{i+1}}{}
=
\dist{p^i}{q^{i+1}}{}$, 
and $\dist{\tilde p^{i+1}}{\tilde q^{i}}{}=\dist{p^{i+1}}{q^{i}}{}$.
\end{subthm}

\begin{figure}[!ht]
\vskip-3mm
\centering
\includegraphics{mppics/pic-1025}
\vskip0mm
\end{figure}

Then for any choice of $n$ points $\tilde z^i\in [\tilde p^i\tilde q^i]$,
we have
\[\dist{\tilde x}{\tilde z^1}{}+\dist{\tilde z^1}{\tilde z^2}{}+\dots+\dist{\tilde z^{n-1}}{\tilde z^n}{}+\dist{\tilde z^n}{\tilde y}{}
\ge 
\dist{x}{y}{}.\]
\end{thm}

\section{Remarks}\label{sec:kirszbraun:open}

The following problem and its relatives were mentioned by Michael Gromov \cite[1.19]{gromov-MS}.

\begin{thm}{Open problem}\label{open:n-point-CBB}
Find a necessary and sufficient condition for a finite metric space to admit distance-preserving embeddings into 

\begin{subthm}{}
some length $\Alex{\kappa}$ space,
\end{subthm}

\begin{subthm}{}
some length $\CAT{\kappa}$ space.
\end{subthm}

\end{thm}

A metric on a finite set $\{a^1,a^2,\dots,a^n\}$,
can be described by the matrix with components
\[s^{ij}
=
\dist[2]{a^i}{a^j}{},\]
which we will call the \index{associated matrix}\emph{associated matrix}.
The set of associated matrices of all metrics that admit a distance-preserving map into a $\Alex{0}$ or a $\CAT{0}$ space 
form a convex cone. 
The latter follows since the rescalings and products of $\Alex{0}$ (or $\CAT{0}$) spaces are $\Alex{0}$ (or $\CAT{0}$ respectively).
This convexity gives a bit of hope that the cone admits an explicit description.

For the 5-point $\CAT{0}$ case, the (2+2)-comparison is a necessary and sufficient condition.
This was proved by Tetsu Toyoda \cite{toyoda};
another proof was found by Nina Lebedeva and the third author \cite{lebedeva-petrunin:toyoda}.
For the 5-point $\Alex{0}$ case, the (1+4)-comparison is a necessary and sufficient condition;
it was proved by Nina Lebedeva and the third author \cite{lebedeva-petrunin:5CBB}.
Starting from the 6-point case, only some necessary and some sufficient conditions are known;
for more on the subject see \cite{alexander-kapovitch-kirszbraun, lebedeva-petrunin:toyoda, lebedeva-petrunin-zolotov}.

The following conjecture (if true) would give the right generality for Kirszbraun's theorem (\ref{thm:kirsz+}).
It states that the example \ref{example:SS_+} 
is the only obstacle to extending short maps.

\begin{thm}{Conjecture}\label{conj:kirsz}
Assume $\spc{L}$ is a complete length $\Alex1$ space,
$\spc{U}$ is a complete length $\CAT1$ space,
$Q\subset \spc{L}$ is a proper subset,
and $f\: Q\to\spc{U}$ is a short map that does not admit a short extension to any bigger set $Q'\supset Q$. 
Then: 

\begin{subthm}{}
$Q$ is isometric to a sphere in a Hilbert space (of finite or cardinal dimension).
Moreover, there is a point $p\in \spc{L}$ such that $\dist{p}{q}{}=\tfrac{\pi}{2}$ for any $q\in Q$.
\end{subthm}

\begin{subthm}{}
The map $f\:Q\to\spc{U}$ is a distance-preserving map and there is no point $p'\in \spc{U}$ such that $\dist{p'}{q'}{}=\tfrac{\pi}{2}$ for any $q'\in f(Q)$.
\end{subthm}
\end{thm}

\parbf{Curvature-free analogs.}
Let us present a collection of exercises on curvature-free analogs of Kirszbraun's theorem.
It is worthwhile to know these results despite they are far from Alexandrov geometry.

\begin{thm}{Exercise}\label{ex:petrunin-stadler}
Let $\spc{X}$ and $\spc{Y}$ be metric spaces, $A\subset \spc{X}$, and $f\:A\to \spc{Y}$ be a short map.
Assume $\spc{Y}$ is compact and for any finite set $F\subset \spc{X}$ there is a short map $F\to \spc{Y}$ that agrees with $f$ on $F\cap A$.
Then there is a short map $\spc{X}\to \spc{Y}$ that agrees with $f$ on $A$.
\end{thm}

\begin{thm}{Exercise}\label{ex:isbell}
We say that a metric space $\spc{X}$ is \index{injective metric space}\emph{injective} 
if for an arbitrary metric space $\spc Z$ 
and a subset $Q\subset\spc Z$, 
any short map $Q\to\spc{X}$ can be extended as a short map $\spc{Z}\to\spc{X}$.
\begin{subthm}{ex:isbell:embedded}
Prove that any metric space $\spc X$ admits a distance-preserving embedding into an injective metric space.
\end{subthm}

\begin{subthm}{ex:isbell:injective hull}
Use this to construct an analog of convex hull in the category of metric spaces; this is called the \index{injective hull}\emph{injective hull}.
\end{subthm}
\end{thm}


\chapter{Warped products}\label{chapter:warped products}

The warped product is a construction that produces 
a new metric space, denoted by $\spc{B}\warp f \spc{F}$,
from two metric spaces base $\spc{B}$ and fiber $\spc{F}$, and a function $f\:\spc{B}\to\RR_{\ge0}$. 

Many important constructions such as direct product, cone, spherical suspension, and join can be defined using warped products.

\section{Definition}
\label{sec:wp-def}
\label{sec:wp-properties}

First we define the warped product for length spaces and then we expand the definition to allow for arbitrary fiber $\spc{F}$.

Let $\spc{B}$ and $\spc{F}$ be length spaces and $f\:\spc{B}\to [0,\infty)$ be a continuous function.

For any path $\gamma\:[0,1]\to \spc{B}\times \spc{F}$, we write $\gamma=(\gamma_\spc{B},\gamma_\spc{F})$ where 
$\gamma_\spc{B}$ is the projection of $\gamma$ to $\spc{B}$, 
and $\gamma_\spc{F}$ is the projection to $\spc{F}$.
If $\gamma_\spc{B}$ and $\gamma_\spc{F}$ are Lipschitz, set
\[
\length_f \gamma \df \int\limits_0^1 \sqrt{
v_\spc{B}^2+ (f\circ\gamma_\spc{B})^2\cdot v_\spc{F}^2}\cdot dt,
\eqlbl{eq:length}
\]
where $\int$ is Lebesgue integral, and $v_\spc{B}$ and $v_\spc{F}$ are the speeds of $\gamma_\spc{B}$ and $\gamma_\spc{F}$ respectively. 
(Note that $\length_f \gamma\ge \length \gamma_\spc{B}$.)

Consider the pseudometric on $\spc{B}\warp f\spc{F}$ defined by
 \[
 \dist{x}{y}{}
 \df 
 \inf\set{\length_f\gamma}{\gamma(0)=x, \gamma(1)=y}
 \]
where the exact lower bound is taken for all Lipschitz paths $\gamma\:[0,1]\to \spc{B}\times\spc{F}$. 
The corresponding metric space is called the \index{warped product}\emph{warped product with base $\spc{B}$, fiber $\spc{F}$ and warping function $f$}; it will be denoted by $\spc{B}\warp f\spc{F}$.

The points in $\spc{B}\warp f\spc{F}$ can be described by corresponding pairs $(p,\phi)\in\spc{B}\times\spc{F}$. Note that if $f(p)=0$ for $p\in \spc{B}$, then $(p,\phi)\z=(p,\psi)$ for any $\phi,\psi\in \spc{F}$.

We do not claim that every Lipschitz curve in $\spc{B}\warp f \spc{F}$ may be reparametrized as the image of a Lipschitz curve in $\spc{B}\times \spc{F}$; in fact this is not true.

\begin{thm}{Proposition}
The warped product $\spc{\spc{B}}\warp{f}\spc{F}$ satisfies:

\begin{subthm}{horiz-leaf-proj}
The projection $(p,\phi_0)\mapsto p$ is a submetry.
Moreover, its restriction to any horizontal leaf $\spc{B}\times\{\phi_0\}$
is an isometry to~$\spc{B}$.
\end{subthm}

\begin{subthm}{vert-leaf-proj}
If $f(p_0)\ne0$, the projection $(p_0,\phi)\mapsto \phi$ of the vertical leaf $\{p_0\}\times \spc{F}$, with its length metric, is a homothety onto $\spc{F}$ with multiplier $\tfrac1{f(p_0)}$.
\end{subthm}

\begin{subthm}{Df>0}If $f$ achieves its (local) minimum at $p_0$, then the inclusion of the vertical leaf $\{p_0\} \times \spc{F}$ in $\spc{B}\warp{f}\spc{F}$ is (locally) distance-preserving.
\end{subthm}

\end{thm}

\parit{Proof.} 
Claim \ref{SHORT.vert-leaf-proj} follows from the $f$-length formula \ref{eq:length}.

Also, by \ref{eq:length}, the projection of
$\spc{B}\warp{f}\spc{F}$ onto $\spc{B}\times\{\phi_0\}$ given by $(p,\phi)\mapsto (p,\phi_0)$ is length-nonincreasing; hence \ref{SHORT.horiz-leaf-proj}.

Suppose $p_0$ is a local minimum point of $f$.
Then the projection $(p,\phi)\mapsto (p_0,\phi)$ of a neighborhood of the vertical leaf $\{p_0\} \times \spc{F}$ to $\{p_0\} \times \spc{F}$ is length-nonincreasing.

If $p_0$ is a global minimum point of $f$, then the same holds for the projection of whole space.
Hence \ref{SHORT.Df>0}.
\qeds

Note that any horizontal leaf is weakly convex, but does not have to be convex even if $\spc{B}\warp{f}\spc{F}$ is a geodesic space.
The latter follows since vanishing of the warping function~$f$ allows geodesics to bifurcate into distinct horizontal leaves.
For instance, if there is a geodesic with the ends in the zero set 
\[Z=\set{(p,\phi)\in\spc{B} \warp{f}\spc{F})}{f(p)=0},\] 
then there is a geodesic with the same endpoints in each horizontal leaf.

\begin{thm}{Proposition}
Suppose $\spc{B}$ and $\spc{F}$ are length spaces and $f\:\spc{B}\to [0,\infty)$ is a continuous function.
Then the warped product $\spc{B}\warp f\spc{F}$ is a length space.
\end{thm}

\parit{Proof.}
It is sufficient to show that for any $\alpha\:[0,1]\to \spc{B}\warp f\spc{F}$ there is a path 
$\beta\:[0,1]\to \spc{B}\times\spc{F}$ with the same endpoints such that 
\[\length \alpha\ge \length_f\beta.\]

If $f\circ\alpha_{\spc{B}}(t)>0$ for any $t$, then the vertical projection $\alpha_{\spc{F}}$ is defined.
In this case, let $\beta(t)=(\alpha_{\spc{B}}(t),\alpha_{\spc{F}}(t))\in \spc{B}\times\spc{F}$.
Clearly 
\[\length \alpha= \length_f\beta.\]

If $f\circ\alpha_{\spc{B}}(t_0) = 0$ for some $t_0$, let $\beta$ be the concatenation of three curves in $\spc{B}\times{\spc{F}}$,
namely: 
(1) the horizontal curve $(\alpha_\spc{B}(t),\phi)$ for $t\le t_0$,
(2) a vertical path in form $(s,\phi)$ to $(s,\psi)$
and 
(3) the horizontal curve $(\alpha_\spc{B}(t),\psi)$ for $t\ge t_0$.
By \ref{eq:length}, the $f$-length of the middle path in the concatenation is vanishing;
therefore the $f$-length of $\alpha$ cannot be smaller than length of $\alpha_{\spc{B}}$,
that is,
\[\length_f\alpha \ge 
\length\alpha_{\spc{B}}=\length_f\beta.
\]
The statement follows.
\qeds
The following theorem states that 
distance in a warped product is fiber-independent, in the sense that distances may be calculated by substituting for $\spc{F}$ a different length space:

\begin{thm}{Fiber-independence theorem}\label{thm:fiber-independence}
Consider length spaces $\spc{B}$, $\spc{F}$ and $\check{\spc{F}}$, and a locally Lipschitz function
$f\:\spc{B}\to\R_{\ge 0}$. 
Assume $p,q\in \spc{B}$, $\phi,\psi\in \spc{F}$ and $\check{\phi},\check{\psi}\in \check{\spc{F}}$.
Then 
\[
\begin{aligned}
\dist{\phi}{\psi}{\spc{F}}
&
\ge\dist{\check{\phi}}{\check{\psi}}{\check{\spc{F}}}
\\
&\Downarrow
\\
\dist{(p,\phi)}{(q,\psi)}{\spc{B}\warp{f}\spc{F}}
&\ge\dist{(p,\check{\phi})}{(q,\check{\psi})}{\spc{B}\warp{f}\check{\spc{F}}}.
\end{aligned}
\]
In particular,
\[
\dist{(p,\phi)}{(q,\psi)}{\spc{B}\warp f \spc{F}} =
\dist{(p,0)}{(q,\ell)}{\spc{B}\warp f\R},
\]
where $\ell=\dist{\phi}{\psi}{\spc{F}}$.
\end{thm}

\parit{Proof.} 
Let $\gamma$ be a path in $(\spc{B}\times \spc{F})$. 

Since $\dist{\phi}{\psi}{\spc{F}}
\ge\dist{\check{\phi}}{\check{\psi}}{\check{\spc{F}}}$,
there is a Lipschitz path $\gamma_{\check{\spc{F}}}$ 
from $\check\phi$ to $\check\psi$ in $\check{\spc{F}}$ such that
\[(\speed\gamma_{\spc{F}})(t)
\ge
(\speed\gamma_{\check{\spc{F}}})(t)\]
for almost all $t\in[0,1]$.
Consider the path $\check\gamma=(\gamma_{\spc{B}},\gamma_{\check{\spc{F}}})$ from $(p,\check\phi)$ to $(q,\check\psi)$ in $\spc{B}\warp{f}\check{\spc{F}}$.
Clearly
\[\length_f\gamma\ge \length_f\check\gamma.\]
\qedsf

\begin{thm}{Exercise}\label{ex:warp=<}
Let $\spc{B}$ and $\spc{F}$ be length spaces and $f,g\:\spc{B}\to \RR_{\ge0}$ be two locally Lipschitz nonnegative functions.
Assume $f(b)\le g(b)$ for any $b\in\spc{B}$.
Show that 
$\spc{B}\warp{f}\spc{F}\le \spc{B}\warp{g}\spc{F}$;
that is, there is a distance-noncontracting map $\spc{B}\warp{f}\spc{F}\to \spc{B}\warp{g}\spc{F}$.
\end{thm}

\section{Extended definitions}

The fiber-independence theorem implies that 
\[
\dist{(p,\phi)}{(q,\psi)}{\spc{B}\warp f \spc{F}} =
\dist{(p,0)}{(q,\dist{\phi}{\psi}{\spc{F}})}{\spc{B}\warp f\R}
\]
for any $(p,\phi),(q,\psi) \in \spc{B}\times \spc{F}$.
In particular, if $\iota\:A\to \check A$ is an isometry between two subsets
$A\subset \spc{F}$ and $\check A\subset \check{\spc{F}}$
in length spaces $\spc{F}$ and $\check{\spc{F}}$, and $\spc{B}$ is a length space, then for any warping function $f\:\spc{B}\to\RR_{\ge0}$,
the map $\iota$ induces an isometry between the sets 
$\spc{B}\warp{f} A \subset \spc{B}\warp{f} \spc{F}$ and $\spc{B}\warp{f}\check{A}\subset \spc{B}\warp{f} \check{\spc{F}}$.

This observation allows us to define the warped product $\spc{B}\z{\warp{f}} \spc{F}$ where the fiber $\spc{F}$ does not carry its length metric.
Indeed we can use Kuratowsky embedding to realize $\spc{F}$ as a subspace in a length space, say $\spc{F}'$.
Therefore we can take the warped product $\spc{B}\warp{f} \spc{F}'$
and identify $\spc{B}\warp{f} \spc{F}$ with its subspace consisting of all pairs $(b,\phi)$ such that $\phi\in \spc{F}$.
According to the Fiber-independence theorem \ref{thm:fiber-independence}, the resulting space does not depend on the choice of $\spc{F}'$.

\section{Examples}\label{direct-products}

\parbf{Direct product.}
The simplest example is the \index{direct product}\emph{direct product} $\spc{B}\times \spc{F}$, which can also be written as the warped product $\spc{B}\warp1 \spc{F}$. 
That is, for $p,q\in \spc{B}$ and $\phi,\psi\in \spc{F}$, the latter metric simplifies to
\[
\dist{(p,\phi)}{(q,\psi)}{} =\sqrt{\dist[2]{p}{q}{} + \dist[2]{\phi}{\psi}{}}.
\]

\parbf{Cones.}
The \index{Euclidean cone}\emph{Euclidean cone} $\Cone\spc{F}$ over a metric space $\spc{F}$
can be written as the warped product $[0,\infty)\warp{\id} \spc{F}$.
That is, for $s,t\in [0,\infty)$ and $\phi,\psi\in \spc{F}$, 
the metric is given by the cosine rule
\[
\dist{(s,\phi)}{(t,\psi)}{} 
=
\sqrt{s^2+t^2-2\cdot s\cdot t\cdot \cos\alpha},
\]
where $\alpha= \min\{\pi, \dist{\phi}{\psi}{}\}$.
(See Section \ref{sec: tangent space}.)

Instead of the Euclidean cosine rule, 
we may use the cosine rule in $\Lob2\kappa$:
\[
\dist{(s,\phi)}{(t,\psi)}{} 
=
\side\kappa\{\alpha;s,t\}.
\]
This way we get \index{$\kappa$-cones}\emph{$\kappa$-cones} over $\spc{F}$, denoted by $\Cone\mc\kappa\spc{F}=[0,\infty)\warp{\sn\kappa} \spc{F}$ for $\kappa\le 0$
and $\Cone\mc\kappa\spc{F}=[0,\varpi\kappa]\warp{\sn\kappa} \spc{F}$ for $\kappa>0$.

The $1$-cone $\Cone\mc1\spc{F}$ is also called the \index{spherical suspension}\emph{spherical suspension} over $\spc{F}$;
it is also denoted by $\Susp\spc{F}$.
That is,
\[
\Susp\spc{F}=[0,\pi]\warp{\sin}\spc{F}.
\]

\begin{thm}{Exercise}\label{ex:convexity-in-cone}
Let $\spc{F}$ be a length space and $A\subset \spc{F}$.
Show that $\Cone\mc\kappa A$ is convex in $\Cone\mc\kappa\spc{F}$ 
if and only if $A$ is $\pi$-convex in $\spc{F}$.
\end{thm}

\parbf{Doubling.}
The doubling space $\spc{W}$ of a metric space $\spc{V}$ on a closed subset $A\subset\spc{V}$
can be also defined as a special type of warped product.
Consider the fiber $\mathbb{S}^0$ consisting of two points with distance $2$ from each other.
Then 
\[\spc{W}\iso\spc{V}\warp{\distfun{A}{}{}}\mathbb{S}^0;\]
that is,
$\spc{W}$ is isometric to the warped product 
with base $\spc{V}$, 
fiber $\mathbb{S}^0$ and warping function $\distfun{A}{}{}$.

\section{1-dimensional base}

The following theorems provide conditions for the spaces and functions in a warped product with 1-dimensional base to have curvature bounds.
These theorems are originally due to Valerii Berestovskii \cite{berestovskii}.
They are baby cases of the characterization of curvature bounds in warped products given in \cite{alexander-bishop:warps,alexander-bishop:warp1}.

\begin{thm}{Theorem}\label{thm:warp-curv-bound:cbb}
\begin{subthm}{thm:warp-curv-bound:cbb:a}
If $\spc{L}$ is a complete length $\Alex{1}$ space and $\diam\spc{L}\le\pi$, 
then 
\begin{align*}
\Susp\spc{L}&=[0,\pi]\warp{\sin}\spc{L}\quad\text{is}\quad \Alex1,
\\
\Cone\spc{L}&=[0,\infty)\warp{\id}\spc{L}\quad\text{is}\quad \Alex0,
\\
\Cone\mc{-1}\spc{L}&=[0,\infty)\warp{\sinh}\spc{L}\quad\text{is}\quad \Alex{-1}.
\end{align*}
Moreover, the converse also holds in each of the three cases.
\end{subthm}

\begin{subthm}{thm:cbb-product}
If $\spc{L}$ is a complete length $\Alex0$ space, 
then 
\begin{align*}
\RR\times\spc{L}&\quad\text{is a complete length $\Alex0$ space},
\\
\RR\warp{\exp}\spc{L}&\quad\text{is a complete length $\Alex{-1}$ space.}
\end{align*}
Moreover, the converse also holds in each of the two cases.
\end{subthm}

\begin{subthm}{}
If $\spc{L}$ is a complete length $\Alex{-1}$ space,
if and only if the warped product $\RR\warp{\cosh}\spc{L}$ is a complete length $\Alex{-1}$ space.
\end{subthm}
\end{thm}


\begin{thm}{Theorem}\label{thm:warp-curv-bound:cat}
Let $\spc{L}$ be a metric space.
\begin{subthm}{thm:warp-curv-bound:cbb:S}
If $\spc{L}$ is $\CAT{1}$,
then 
\begin{align*}
\Susp\spc{L}&=[0,\pi]\warp{\sin}\spc{L}\quad\text{is}\quad \CAT{1},
\\
\Cone\spc{L}&=[0,\infty)\warp{\id}\spc{L}\quad\text{is}\quad \CAT{0},
\\
\Cone\mc{-1}\spc{L}&=[0,\infty)\warp{\sinh}\spc{L}\quad\text{is}\quad \CAT{-1}.
\end{align*}
Moreover, the converse also holds in each of the three cases.
\end{subthm}

\begin{subthm}{thm:warp-curv-bound:cbb:E}
If $\spc{L}$ is a complete length $\CAT0$ space,
then 
$\RR\times\spc{L}$ is $\CAT0$ 
and 
$\RR\warp{\exp}\spc{L}$ is $\CAT{-1}$.
Moreover, the converse also holds in each of the two cases.
\end{subthm}

\begin{subthm}{thm:warp-curv-bound:cbb:H}
If $\spc{L}$ is $\CAT{-1}$
if and only if 
$\RR\warp{\cosh}\spc{L}$
is $\CAT{-1}$.
\end{subthm}
\end{thm}


In the proof of the above two theorems 
we will use the following proposition.

\begin{thm}{Proposition}\label{prop:warp-examples}

\begin{subthm}{prop:warp-examples:S}
\begin{align*}
\Susp\mathbb S^{m-1}&=[0,\pi]\warp{\sin}\mathbb S^{m-1}\iso\mathbb{S}^m,
\\
\Cone\mathbb S^{m-1}&=[0,\infty)\warp{\id}\mathbb S^{m-1}\iso\EE^m,
\\
\Cone\mc{-1}\mathbb S^{m-1}&=[0,\infty)\warp{\sinh}\mathbb S^{m-1}\iso\Lob{m}{-1}.
\end{align*}
\end{subthm}

\begin{subthm}{}
\begin{align*}
\RR\times\EE^{m-1}&\iso\EE^{m},
\\
\RR\warp{\exp}\EE^{m-1}&\iso\Lob{m}{-1}.
\end{align*}
\end{subthm}

\begin{subthm}{}
\[\RR\warp{\cosh}\Lob{m-1}{-1}\iso\Lob{m}{-1}.\]
\end{subthm}

\end{thm}

The proof is left to the reader.

\parit{Proof of \ref{thm:warp-curv-bound:cbb} and \ref{thm:warp-curv-bound:cat}.}
Each proof is based on the fiber-independence theorem~\ref{thm:fiber-independence} 
and 
the corresponding statement in Proposition~\ref{prop:warp-examples}.

Let us prove the last statement in \ref{SHORT.thm:warp-curv-bound:cbb:S}; the remaining statements of this theorem are similar.
Choose an arbitrary quadruple of points 
\[(s,\phi),(t^1,\phi^1),(t^2,\phi^2),(t^3,\phi^3)\in[0,\infty)\warp{\sinh} \spc{L}.\]
Since $\diam\spc{L}\le\pi$,
the (1+3)-point comparison (\ref{thm:pos-config}) provides a quadruple of points $\psi,\psi^1,\psi^2,\psi^3\in\mathbb{S}^3$ such that 
\[\dist{\psi}{\psi^i}{\mathbb{S}^3}=\dist{\phi}{\phi^i}{\spc{L}}\] 
and
\[\dist{\psi^i}{\psi^j}{\mathbb{S}^3}\ge\dist{\phi^i}{\phi^j}{\spc{L}}\]
for all $i$ and $j$.

According to Proposition~\ref{prop:warp-examples:S}, 
\[\Cone\mc{-1}\mathbb{S}^3=[0,\infty)\warp{\sinh}\mathbb{S}^3\iso\Lob{4}{-1}.\]

Consider the quadruple of points
\[(s,\psi),(t^1,\psi^1),(t^2,\psi^2),(t^3,\psi^3)\in \Cone\mc{-1}\mathbb{S}^3=\Lob{4}{-1}.\]

By the fiber-independence theorem~\ref{thm:fiber-independence},
\[\dist{(s,\psi)}{(t^i,\psi^i)}{[0,\infty)\warp{\sinh}\mathbb{S}^3}=\dist{(s,\phi)}{(t^i,\phi^i)}{[0,\infty)\warp{\sinh}\spc{L}}\]
and
\[\dist{(t^i,\psi^i)}{(t^j,\psi^j)}{[0,\infty)\warp{\sinh}\mathbb{S}^3}\ge\dist{(t^i,\phi^i)}{(t^j,\phi^j)}{[0,\infty)\warp{\sinh}\spc{L}}\]
for all $i$ and $j$.
Since four points of $\Lob{4}{-1}$ lie in an isometric copy of $\Lob{3}{-1}$, it remains to apply Exercise \ref{ex:(3+1)-expanding}.

\medskip

Now let us prove the converse to \ref{SHORT.thm:warp-curv-bound:cbb:S}.
Choose a quadruple $\phi$, $\phi^1$, $\phi^2$, $\phi^3\in\spc{L}$ with all distances smaller than $\tfrac\pi2$.
Choose small $s>0$ and let $t_i$ be the hypotenuse in a hyperbolic right triangle with angle $\dist{\phi}{\phi^i}{\spc{L}}$ and the adjacent side $s$.
Observe that the $\Lob{}{-1}$ model angles of the quadruple 
$(s,\phi)$, $(t^1,\phi^1)$, $(t^2,\phi^2)$, $(t^3,\phi^3)$ in $[0,\infty)\warp{\sinh} \spc{L}$
at $(s,\phi)$ are the same as the $\Lob{}{1}$ model angles of the quadruple $\phi$, $\phi^1$, $\phi^2$, $\phi^3\in\spc{L}$ at $\phi$.
Whence $\Alex1$ comparison holds for $\phi,\phi^1,\phi^2,\phi^3$;
in particular, $\spc{L}$ is locally $\Alex1$.
It remains to apply the globalization theorem (\ref{thm:glob}). 

\medskip

The proof of \ref{thm:warp-curv-bound:cat} is nearly identical, but one has to apply (2+2)-comparison (\ref{ex:sba-2+2-short}).
\qeds

\begin{thm}{Exercise}\label{ex:spherical-join}
The 
\index{spherical join}\emph{spherical join} $\spc{U}\star\spc{V}$ of two metric spaces $\spc{U}$ and $\spc{V}$
is defined as the unit sphere equipped with the angle metric in the product of Euclidean cones $\Cone \spc{U}\times \Cone\spc{V}$.

Assume $\spc{U}$ and $\spc{V}$ are nonempty spaces.
\begin{subthm}{}Show that $\spc{U}\star\spc{V}$ is $\CAT1$ if and only if $\spc{U}$ and $\spc{V}$ are $\CAT1$.
\end{subthm}

\begin{subthm}{}
Suppose $\spc{U}$ and $\spc{V}$ have $\diam\le \pi$.
Show that $\spc{U}\star\spc{V}$ is $\Alex1$ if and only if $\spc{U}$ and $\spc{V}$ are $\Alex1$.
\end{subthm}

\end{thm}

\section{Remarks}

Let us formulate general results on curvature bounds of warped products 
proved by the first author and Richard Bishop \cite{alexander-bishop:warps}.

\begin{thm}{Theorem}\label{thm:warp-CBB}
Let $\spc{B}$ be a complete finite-dimensional $\Alex\kappa$ length space, and $f\:{\spc{B}}\to \RR_{\ge0}$ be a locally Lipschitz function.
Denote by $Z\subset \spc{B}$ the zero set of the restriction of $f$ to the boundary $\partial \spc{B}$ of $\spc{B}$.

Suppose that $\spc{W}$ is doubling of $\spc{B}$ along the closure of $\partial \spc{B}\setminus Z$, and $\bar f\: \spc{W}\to\RR_{\ge0}$ is the natural extension of $f$.
Assume that $\spc{W}$ is $\Alex\kappa$ 
and $\bar f''+\kappa\cdot \bar f\le 0$.

Suppose $\spc{F}$ is a complete finite-dimensional $\Alex{\kappa'}$ space.
Then the warped product $\spc{B}\warp f\spc{F}$ is $\Alex\kappa$ in the following two cases:
\begin{subthm}{thm:warp-CAT:Z=0}
If $Z= \emptyset$ and 
$\kappa'\ge \kappa\cdot f^2(b)$
for any $b\in \spc{B}$.
\end{subthm}

\begin{subthm}{thm:warp-CAT:Zne0}
If $Z\ne \emptyset$ and
$|d_zf|^2\le\kappa'$
for any $z\in Z$.
\end{subthm}

\end{thm}

We mention that in the setting of this theorem, $f$ necessarily vanishes only at boundary points if $f$ is not identically $0$.

\begin{thm}{Theorem}\label{thm:warp-CAT}
Let $\spc{B}$ be a 
complete 
$\CAT\kappa$ length space, and the function $f\:\spc{B}\to \RR_{\ge0}$ satisfy $f''+\kappa\cdot f\ge 0$, 
where $f$ is Lipschitz on bounded sets or $B$ is locally compact.
Denote by $Z\subset \spc{B}$ the zero set of~$f$.
Suppose $\spc{F}$ is a 
complete $\CAT{\kappa'}$ space.
Then the warped product $\spc{B}\warp f\spc{F}$ is $\CAT\kappa$ in the following two cases:
\begin{subthm}{}
If $Z= \emptyset$ and $\kappa'\le \kappa\cdot f^2(b)$ for any $b\in \spc{B}$.
\end{subthm}

\begin{subthm}{}
If $Z\ne \emptyset$ and $[(d_zf)\dir zb]^2\ge\kappa'$ for any minimizing geodesic $[zb]$ from $Z$ to a point $b\in\spc{B}$ and 
$\kappa'\le \kappa\cdot f^2(b)$ for any $b\in \spc{B}$ such that $\distfun{Z}{b}{}\ge \tfrac{\varpi\kappa}2$.
\end{subthm}

\end{thm}

\chapter{Polyhedral spaces}\label{chapter:polyhedral}

\section{Definitions}

\begin{thm}{Definition}\label{def:poly}
A length space $\spc{P}$ is called a
\index{polyhedral space}
\index{polyhedral space!piecewise $\Lob{}\kappa$ space}\emph{piecewise $\Lob{}\kappa$} 
if it admits a finite triangulation $\tau$ 
such that an arbitrary simplex $\sigma$ in $\tau$ is isometric to a simplex in the model space $\Lob{\dim \sigma}{\kappa}$.

By \index{polyhedral space!triangulation of a polyhedral space}\emph{triangulation} of a piecewise
$\Lob{}\kappa$ space 
we will understand a triangulation as in the definition.
If we do not wish to specify $\kappa$, we will say that $\spc{P}$ is a \index{polyhedral space}\emph{polyhedral space}. 

By rescaling we can assume that $\kappa=1$, $0$, or $-1$.
\begin{subthm}{}
Piecewise $\Lob{}1$ spaces will also be 
called 
\index{spherical polyhedral space}\emph{spherical polyhedral spaces};
\end{subthm}

\begin{subthm}{}
Piecewise $\Lob{}0$ spaces will also be 
called 
\index{Euclidean polyhedral space}\emph{Euclidean polyhedral spaces};
\end{subthm}

\begin{subthm}{}
Piecewise  $\Lob{}{-1}$ spaces will also be 
called 
\emph{hyperbolic polyhedral spaces}\index{hyperbolic polyhedral space}.
\end{subthm}
\end{thm}

Note that according to the above definition,
all polyhedral spaces are compact.
However, 
most of the statements below admit straightforward generalizations 
to \index{polyhedral space!locally polyhedral space}\emph{locally polyhedral space};
that is, complete length spaces,  
any point of which admits a closed neighborhood isometric to a polyhedral space.
The latter class of spaces includes in particular  infinite covers of polyhedral spaces.

The dimension of a polyhedral space $\spc{P}$
is defined as the maximal dimension of a simplex 
in one (and therefore any) triangulation of~$\spc{P}$.

\parbf{Links.}
Let $\spc{P}$ be a polyhedral space
and $\sigma$ be a simplex in its triangulation~$\tau$.

The simplexes that contain $\sigma$
form an abstract simplicial complex called the \index{link}\emph{link} of $\sigma$, 
denoted by $\Link_\sigma$.
If $m=\dim\sigma$,
then the set of vertexes of $\Link_\sigma$
is formed by the $(m+1)$-simplexes that contain $\sigma$;
the set of its edges are formed by the $(m+2)$-simplexes 
that contain $\sigma$, and so on.

The link $\Link_\sigma$
can be identified with the subcomplex of $\tau$ 
formed by all the simplexes $\sigma'$ 
such that $\sigma\cap\sigma'=\emptyset$ 
but both $\sigma$ and $\sigma'$ are faces of a simplex of $\tau$.

The points in $\Link_\sigma$ can be identified with the normal directions to $\sigma$ at a point in the interior of $\sigma$.
The angle metric between directions makes  $\Link_\sigma$ into a spherical polyhedral space.
We will always consider the link with this metric.

\parbf{Tangent space and space of directions.}
Let $\tau$ be a triangulation of a polyhedral space $\spc{P}$.
If a point $p\in \spc{P}$ 
lies in the interior of a $\kay$-simplex $\sigma$ of $\tau$ 
then the tangent space $\T_p\spc{P}$
is  naturally isometric to
\[\EE^\kay\times(\Cone\Link_\sigma).\]
Equivalently, the space of directions $\Sigma_p$
can be isometrically identified with the 
$\kay$-th spherical suspension over $\Link_\sigma$;
that is, 
\[\Sigma_p\iso\Susp^{\kay}(\Link_\sigma).\]

If $\spc{P}$ is an $m$-dimensional polyhedral space,
then for any $p\in \spc{P}$
the space of directions $\Sigma_p$ is a spherical polyhedral space
of dimension at most $m-1$. 

In particular, 
for any point $p$ in the interior of a simplex $\sigma$,
the isometry class of $\Link_\sigma$ and $\kay=\dim\sigma$
determine the isometry class of $\Sigma_p$ and the other way around.

A small neighborhood of $p$ is isometric to a neighborhood of the tip of the $\kappa$-cone over $\Sigma_p$.
In fact, if this property holds at any point of a compact length space $\spc{P}$
then  $\spc{P}$ is a piecewise
$\Lob{}\kappa$ space \cite{lebedeva-petrunin-poly}.

\section{Curvature bounds}

Recall that $\ell$-simply connected spaces are defined in \ref{def:l-s.c.}.

The following theorem provides a combinatorial description of polyhedral spaces with curvature bounded above.

\begin{thm}{Theorem}\label{thm:PL-CAT}
Let $\spc{P}$ be a piecewise
$\Lob{}\kappa$ space and $\tau$ be a triangulation of $\spc{P}$. Then 

\begin{subthm}{thm:PL-CAT:curc>=k}
$\spc{P}$  is locally $\CAT\kappa$
if and only if any connected component of the link of any simplex $\sigma$ in $\tau$
is $(2\cdot\pi)$-simply connected.
Equivalently, if and only if any closed local geodesic in $\Link_\sigma$ has length at least $2\cdot\pi$.
\end{subthm}

\begin{subthm}{thm:PL-CAT:CAT}
$\spc{P}$ is a complete length $\CAT\kappa$ space
if and only if $\spc{P}$ is $(2\cdot\varpi\kappa)$-simply connected and any connected component of the link of any simplex $\sigma$ in $\tau$
is $(2\cdot\pi)$-simply connected.
\end{subthm}

\end{thm}

\parit{Proof.}
We will prove only the if part;
the only-if part is evident by the generalized Hadamard--Cartan theorem (\ref{thm:hadamard-cartan-gen}) and Theorem~\ref{thm:warp-curv-bound:cbb}.

Let us apply induction on $\dim\spc{P}$.
The  {}\emph{base}  case $\dim\spc{P}=0$ is evident.

\parit{Induction Step.}
Assume that the theorem is proved in the case $\dim\spc{P}\z<m$. Suppose  $\dim\spc{P}=m$.

Fix a point $p\in\spc{P}$.
A neighborhood of $p$ 
is isometric to a neighborhood of the tip in the $\kappa$-cone over 
 $\Sigma_p$.
By Theorem~\ref{thm:warp-curv-bound:cbb:a}, 
it is sufficient to show that 
\begin{clm}{}\label{eq:Sigma-in-CAT(1)}
 $\Sigma_p$ is $\CAT1$.
\end{clm}

Note that $\Sigma_p$ is a spherical polyhedral space 
and its  links are isometric to  links of $\spc{P}$. 
By the  induction hypothesis, $\Sigma_p$ is locally $\CAT1$.
Applying the generalized Hadamard--Cartan theorem (\ref{thm:hadamard-cartan-gen}),
we get \ref{eq:Sigma-in-CAT(1)}.

Part \ref{SHORT.thm:PL-CAT:CAT} follows from the generalized Hadamard--Cartan theorem.
\qeds

\begin{thm}{Exercise}\label{ex:metric tree}
Show that any metric tree is $\CAT\kappa$ for any $\kappa$.
\end{thm}

\begin{thm}{Exercise}\label{ex:poly-unique-geodesic}
Show that if in a Euclidean polyhedral space $\spc{P}$
any two points can be connected by a unique geodesic,  
then $\spc{P}$ is $\CAT0$.
\end{thm}

The following theorem provides a combinatorial description of polyhedral spaces with curvature bounded below.

\begin{thm}{Theorem}\label{thm:poly-CBB} 
Let $\spc{P}$ be a piecewise
$\Lob{}\kappa$ space and $\tau$ be a triangulation of $\spc{P}$.
Then $\spc{P}$ is $\Alex\kappa$ if and only if the following conditions hold.

\begin{subthm}{} $\tau$ is \index{pure triangulation}\emph{pure}; 
that is, any simplex in $\tau$ is  a face of some simplex of dimension exactly $m$. 
\end{subthm}

\begin{subthm}{thm:poly-CBB:m-1}
The link of any simplex of dimension $m-1$ is formed by a single point or two points.
\end{subthm}

\begin{subthm}{thm:poly-CBB:2pi}
Any link of any simplex of dimension $m-2$
has diameter at most $\pi$.
\end{subthm}

\begin{subthm}{thm:poly-CBB:connected}
The link of any simplex of dimension $\le m-2$ is connected.
\end{subthm}

\end{thm}

\parbf{Remarks.}
Condition \ref{SHORT.thm:poly-CBB:connected}
can be reformulated in the following way:

\begin{itemize}
 \item[\ref{SHORT.thm:poly-CBB:connected}$\,'\!$] 
Any path $\gamma\:[0,1]\to \spc{P}$ can be approximated by paths
$\gamma_n$ 
that cross only simplexes of dimension $m$ and $m-1$.
\end{itemize}

Further, modulo the other conditions,
 condition \ref{SHORT.thm:poly-CBB:2pi}
is equivalent to the following:

\begin{itemize}
 \item[\ref{SHORT.thm:poly-CBB:2pi}$\,'\!$] 
The link of any simplex of dimension $m-2$ is 
isometric to a circle of length $\le 2\cdot\pi$
or a closed real interval of length $\le \pi$.
\end{itemize}

\parit{Proof.}
We will prove the if part.
The only-if part is similar and is left to the reader.

We apply induction on $m$.
The {}\emph{base} case $m=1$ follows from the assumption \ref{SHORT.thm:poly-CBB:m-1}.

\parit{Step.}
Assume that the theorem is proved for polyhedral spaces  of dimesnion less than $m$.
Suppose  $\dim\spc{P}=m$.

According to the globalization theorem (\ref{thm:glob}),
it is sufficient to show that 
$\spc{P}$  is locally $\Alex\kappa$. 

Fix $p\in \spc{P}$.
Note that a spherical neighborhood of $p$
is isometric
to a  spherical neighborhood of the tip of the tangent $\kappa$-cone 
$$\T_p\mc\kappa\z
=
\Cone\mc\kappa(\Sigma_p).$$
Hence it is sufficient to show that 
\begin{clm}{}
 $\T_p\mc\kappa$ is $\Alex\kappa$ for any $p\in \spc{P}$.
\end{clm}

By Theorem~\ref{thm:warp-curv-bound:cbb:a}, 
the latter is equivalent to 
\begin{clm}{}\label{clm:curv+diam}
$\diam\Sigma_p\le \pi$ and $\Sigma_p$ is $\Alex1$.
\end{clm}

If $m=2$, then \ref{clm:curv+diam} follows from \ref{SHORT.thm:poly-CBB:m-1}.

To prove the case $m\ge 3$,
note that $\Sigma_p$ is an $(m-1)$-dimensional spherical polyhedral space and all the conditions of the theorem hold for $\Sigma_p$.
It remains to apply the induction hypothesis.\qeds

\begin{thm}{Exercise}\label{ex:polyKk}
Assume $\spc{P}$ is a piecewise
$\Lob{}\kappa$ space and $\dim \spc{P}\ge 2$. 
Show that 

\begin{subthm}{} if $\spc{P}$ is $\Alex{\kappa'}$, then $\kappa'\le \kappa$ and $\spc{P}$ is $\Alex{\kappa}$, 
\end{subthm}

\begin{subthm}{}
if $\spc{P}$ is $\CAT{\kappa'}$, then $\kappa'\ge \kappa$ and $\spc{P}$ is $\CAT\kappa$.
\end{subthm}

\end{thm}

\section{Flag complexes}

\begin{thm}{Definition}
A simplicial complex $\mathcal{S}$ 
is \index{flag complex}\emph{flag} if whenever $\{v_0,\z\dots,v_\kay\}$
is a set of distinct vertexes of $\mathcal{S}$
that are pairwise joined by edges, then the vertexes $v_0,\dots,v_\kay$
span a $\kay$-simplex in $\mathcal{S}$.

If the above condition is satisfied for $\kay=2$, 
then we say $\mathcal{S}$ satisfies 
the \index{no-triangle condition}\emph{no-triangle condition}.
\end{thm}

Note that every flag complex is determined by its 1-skeleton.

\begin{thm}{Proposition}\label{prop:no-trig}
A simplicial complex $\mathcal{S}$ is flag if and only if 
$\mathcal{S}$, as well as all the links of all its simplexes,
satisfy the no-triangle condition.
\end{thm}

From the definition of flag complex 
we get the following:

\begin{thm}{Observation}\label{obs:link-of-flag}
Any link of a flag complex is flag.
\end{thm}

\parit{Proof of Proposition~\ref{prop:no-trig}.}
By Observation~\ref{obs:link-of-flag}, the no-triangle condition holds 
for any flag complex and all its links.

Now assume a complex $\spc{S}$ and all its links satisfy 
the no-triangle condition.
It follows that $\spc{S}$ includes a 2-simplex for each triangle.
Applying the same observation for each edge we get that $\spc{S}$ 
includes a 3-simplex for any complete graph with 4 vertexes.
Repeating this observation 
for triangles, 
4-simplexes,
5-simplexes
and so on we get that $\spc{S}$ is flag.
\qeds

\parbf{Right-angled triangulation.} 
A triangulation of a spherical polyhedral space 
is called \index{right-angled triangulation}\emph{right-angled} 
if each simplex of the triangulation is isometric 
to a spherical simplex all of whose angles are right.
Similarly, we say that a simplicial complex 
is equipped with a \index{right-angled spherical metric}\emph{right-angled spherical metric}
if it is a length metric and each simplex is isometric 
to a spherical simplex all of whose angles are right.

Spherical polyhedral $\CAT1$ spaces glued from right-angled simplexes
admit the following characterization 
discovered by Michael Gromov \cite[p. 122]{gromov:hyp-groups}.

\begin{thm}{Flag condition}\label{thm:flag}
Assume that a spherical polyhedral space $\spc{P}$
admits a right-angled triangulation $\tau$.
Then $\spc{P}$ is $\CAT1$
if and only if $\tau$ is flag.
\end{thm}

\parit{Proof; only-if part.} 
Assume there are three vertexes $v_1,v_2$ and $v_3$ of $\tau$
that are pairwise joined by edges 
but do not span a simplex.
Note that in this case 
$$\mangle\hinge{v_1}{v_2}{v_3}=\mangle\hinge{v_2}{v_3}{v_1}=\mangle\hinge{v_3}{v_1}{v_1}=\pi.$$
Equivalently,
\begin{clm}{}\label{clm:3pi/2}
The concatenation of the geodesics $[v_1v_2]$, $[v_2v_3]$ and $[v_3v_1]$
forms a closed local geodesic in $\spc{P}$. 
\end{clm}

Now assume that $\spc{P}$ is $\CAT1$.
Then by \ref{thm:warp-curv-bound:cbb:a},
$\Link_\sigma\spc{P}$ is a compact length $\CAT1$ space for every simplex $\sigma$ 
in $\tau$. 

Each of these links is a right-angled spherical complex
and
by Theorem \ref{thm:PL-CAT}, 
none
of these links can contain a geodesic circle of length less than $2\cdot\pi$. 

Therefore Proposition~\ref{prop:no-trig} and \ref{clm:3pi/2} 
imply the only-if part.

\parit{If part.} 
By Observation~\ref{obs:link-of-flag} and Theorem~\ref{thm:PL-CAT},
it is sufficient to show that any closed local geodesic $\gamma$ 
in a flag complex $\spc{S}$ with right-angled metric has length at least $2\cdot\pi$.

Fix a flag complex $\spc{S}$.
Recall that the  \index{star of vertex}\emph{star} of a vertex $v$ (briefly $\overline \Star_v$)
is formed by all the simplexes  containing $v$. Similarly, $\Star_v$,   the open star of $v$, is the union of all simplexes containing $v$ with faces opposite $v$ removed.

\begin{wrapfigure}{o}{45mm}
\vskip-0mm
\centering
\includegraphics{mppics/pic-1200}
\end{wrapfigure}

Choose a simplex $\sigma$ that contains a point of $\gamma$.
Let $v$ be a vertex of $\sigma$.
Set $f(t)=\cos\dist{v}{\gamma(t)}{}$.
Note that 
\[f''(t)+f(t)=0\] if $f(t)>0$.  
Since the zeroes of sine are $\pi$ apart,
$\gamma$ 
spends time $\pi$ on every visit to $\Star_v$.

After leaving $\Star_v$,
the local geodesic $\gamma$ must enter another simplex, 
say $\sigma'$, 
which has a vertex $v'$ not joined to $v$ by an edge.

Since $\tau$ is flag, the open  stars $\Star_v$ and $\Star_{v'}$
do not overlap.
The same argument as above shows that $\gamma$ spends time $\pi$ on every visit to $\Star_{v'}$.
Therefore the total length of $\gamma$ is at least $2\cdot\pi$.
\qeds

\begin{thm}{Exercise}\label{ex:barycenric-flag}
Show that the barycentric subdivision of any simplicial complex is flag.
Conclude that any finite  simplicial complex is homeomorphic to a compact length $\CAT1$ space.
\end{thm}

\begin{thm}{Exercise}\label{ex:tree-product}
Let $p$ be a point in a product of metric trees.
Show that a closed geodesic in the space of directions $\Sigma_p$ has length either $2\cdot\pi$ or at least $3\cdot\pi$.
\end{thm}

\begin{thm}{Exercise}\label{ex:obtuce-flag}
Assume that a spherical polyhedral space $\spc{P}$
admits a triangulation $\tau$ such that all edgelengths of all simplexes in $\tau$ are at least $\tfrac\pi2$
Show that $\spc{P}$ is $\CAT{1}$
if $\tau$ is flag.
\end{thm}

\begin{thm}{Exercise}\label{ex:short+commuting}
Let $\phi_1,\phi_2,\dots,\phi_k\:\spc{U}\to \spc{U}$ be commuting short retractions of
a complete length $\CAT0$ space; 
that is, 
\begin{itemize}
\item $\phi_i\circ\phi_i=\phi_i$ for each $i$;
\item $\phi_i\circ\phi_j=\phi_j\circ\phi_i$ for any $i$ and $j$;
\item $\dist{\phi_i(x)}{\phi_i(y)}{\spc{U}}\le \dist{x}{y}{\spc{U}}$ for each $i$ and any $x,y\in\spc{U}$.
\end{itemize}
Set $A_i=\Im \phi_i$ for all $i$.
Note that each $A_i$ is a weakly convex set.

Assume $\Gamma$ is a finite graph 
(without loops and multiple edges) 
with edges labeled by $1,2,\dots, n$.
Denote by $\spc{U}^\Gamma$ the space obtained by taking 
a copy of $\spc{U}$ for each vertex of $\Gamma$ and 
gluing two such copies along $A_i$ if the corresponding vertexes are joined by an edge labeled by $i$.

Show that $\spc{U}^\Gamma$ is $\CAT0$
\end{thm}

\parbf{The space of trees.}
The following construction is given by 
Louis Billera,
Susan Holmes,
and  Karen Vogtmann \cite{billera-holmes-vogtmann}.

Let $\spc{T}_n$ be the set of all metric trees with $n$ end-vertexes
labeled by $a_1,\dots,a_n$.
To describe one tree in $\spc{T}_n$ we may fix a topological tree $\tau$ with end vertexes $a_1,\dots,a_n$ and all the other vertexes of degree 3,  
and prescribe the lengths of $2\cdot n-3$ edges.
If the length of an edge is $0$, we assume that edge degenerates;
such a tree can be also described using a different topological tree $\tau'$.
The subset of $\spc{T}_n$ corresponding to the given topological tree $\tau$ can be identified with a convex closed cone in  $\mathbb{R}^{2\cdot n-3}$.
Equip each such subset with the metric induced from $\mathbb{R}^{2\cdot n-3}$ and consider the length metric on $\spc{T}_n$ induced by these metrics.

\begin{thm}{Exercise}\label{ex:space-of-trees}
Show that $\spc{T}_n$ with the described metric is $\CAT0$.
\end{thm}

\parbf{Cubical complexes.}
The definition of a cubical complex
mostly repeats the definition of a simplicial complex, 
with simplexes replaced by cubes.

Formally, a \index{cubical complex}\emph{cubical complex} is defined as a subcomplex 
of the unit cube in Euclidean space of large dimension;
that is, a collection of faces of the cube, that with each face contains all its sub-faces.
Each cube face in this collection 
will be called a {}\emph{cube} of the cubical complex.

Note that according to this definition, 
any cubical complex is finite,
that is, contains a finite number of cubes.

The union of all the cubes in a cubical complex $\spc{Q}$ will be called its \index{underlying space of a cubical complex}\emph{underlying space};
it will be denoted by $\spc{Q}$ or by $\ushort{\spc{Q}}$ 
if we need to emphasize that we are talking about a topological space, 
not a complex.
A homeomorphism from $\ushort{\spc{Q}}$ to a topological space $\spc{X}$ is called a \index{cubulation}\emph{cubulation of}~$\spc{X}$.

The underlying space of a cubical complex $\spc{Q}$ will be always considered with the length metric
induced from $\RR^N$.
In particular, with this metric, 
each cube of $\spc{Q}$ is isometric to the unit cube of the same dimension.

It is straightforward to construct a triangulation 
of $\ushort{\spc{Q}}$ 
such that each simplex is isometric to a Euclidean simplex.
In particular, $\ushort{\spc{Q}}$ is a Euclidean polyhedral space.

The link of each cube in a cubical complex admits a natural right-angled triangulation; 
each simplex corresponds to an adjusted cube.

\begin{thm}{Exercise}\label{ex:cubical-complex}
Show that a cubical complex $\spc{Q}$ is locally $\CAT0$ if and only if the link of each vertex in $\spc{Q}$ is flag.
\end{thm}

\section{Remarks}\label{sec:poly-reamarks}

The condition on polyhedral $\CAT\kappa$ spaces given in Theorem~\ref{thm:PL-CAT} might look easy to use, 
but in fact, it is hard to check even in simple cases.
For example, the description of those coverings of $\mathbb{S}^3$ that branch at three 
great circles and are $\CAT1$ requires quite a bit of work;
an answer is given by Ruth Charney and Michael Davis \cite{charney-davis-1993}.

Analogs of the flag condition for spherical Coxeter simplexes
could resolve the following problem. 

\begin{thm}{Braid space problem}
Consider $\CC^n$ with coordinates $z_1,\dots,z_n$.
Let us remove from $\CC^n$ the complex hyperplanes $z_i=z_j$ for all $i\ne j$,
pass to the universal cover, and consider the completion $\spc{B}_n$ 
of the obtained space.

Is it true that $\spc{B}_n$ is $\CAT0$ for any $n$?
\end{thm}

The above question has an affirmative answer for $n\le 3$ and is open for all $n\ge 4$ \cite{charney-davis-1993,panov-petrunin}.

Recall that by the Hadamard--Cartan theorem (\ref{thm:hadamard-cartan}), 
any complete length $\CAT{0}$ space is contractible.
Therefore any complete length locally $\CAT{0}$ space 
is \index{aspherical space}\emph{aspherical};
that is, has contractible universal cover.

This observation can be used together with Exercise \ref{ex:cubical-complex} to construct examples of exotic aspherical spaces;
for example, compact topological manifolds with universal cover not homeomorphic to a Euclidean space.
A survey on the subject is given by Michael Davis \cite{davis-2001}; a more elementary introduction to the subject is given by the authors \cite[Chapter 3]{alexander-kapovitch-petrunin-CAT}.

The flag condition also leads to the so-called {}\emph{hyperbolization} procedure, a flexible tool for constructing  aspherical spaces;
a good survey on the subject is given by Ruth Charney and Michael Davis \cite{charney-davis-1995}.

The $\CAT0$ property of a cube complex admits interesting (and useful) geometric descriptions if one replaces the $\ell^2$-metric with a natural $\ell^1$- or $\ell^\infty$-metric on each cube.
The following statement was proved by Brian Bowditch \cite{bowditch-2020}.

\begin{thm}{Theorem}
The following three conditions are equivalent.

\begin{subthm}{cube-2} A cube complex $Q$ equipped with  $\ell^2$-metric is $\CAT0$.
\end{subthm}

\begin{subthm}{cube-infty} A cube complex $Q$ equipped with $\ell^\infty$-metric is \index{injective metric space}\emph{injective}; that is, for any metric space $\spc{X}$ with a subset $A$, any short map $A\to (Q,\ell^\infty)$ can be extended to a short map $\spc{X}\to (Q,\ell^\infty)$.
\end{subthm}

\begin{subthm}{cube-1} A cube complex $Q$ equipped with $\ell^1$-metric is \emph{median}; that is, for any three points $x,y,z$ there is a unique point   $m$ (called the  \index{median of three points}\emph{median} of $x$, $y$, and $z$) that lies on some geodesics $[xy]$, $[xz]$ and $[yz]$.
\end{subthm}
\end{thm}

\part{Structure and tools}
\chapter{First order differentiation}\label{chap:tan}

\section{Ultratangent space} 

The following theorem is often used together with the 
observation that the ultralimit of any sequence of length spaces is geodesic (see \ref{obs:ultralimit-is-geodesic}).

\begin{thm}{Theorem}\label{thm:tan-is}
\begin{subthm}{thm:tan-is-CBB}
If $\spc{L}$ is a complete length  $\Alex{\kappa}$ space and $p\in \spc{L}$, then $\T^\o_p$ is $\Alex{0}$.
\end{subthm}

\begin{subthm}{thm:tan-is-CBA}
If $\spc{U}$ is a complete length $\CAT\kappa$ space and $p\in \spc{U}$, then $\T^\o_p$ is $\CAT0$.
\end{subthm}

\end{thm}

The proofs of both parts are nearly identical.

\parit{Proof; \ref{SHORT.thm:tan-is-CBB}.}
Since $\spc{L}$ is a complete length $\Alex{\kappa}$ space, then its blowup $n\cdot\spc{L}$
(see Section \ref{sec: ultradiff})
 is a complete length $\Alex{\kappa/{n^{2}}}$ space.
By Proposition~\ref{prp:A^omega}, the $\o$-blowup $\o\cdot\spc{L}$ is $\Alex0$
and so is $\T_p^\o$ as a metric component of~$\o\cdot\spc{L}$.

\parit{\ref{SHORT.thm:tan-is-CBA}.}
Since $\spc{U}$ is a complete length $\CAT\kappa$ space, then its blowup $n\cdot\spc{U}$ is $\CAT{\kappa/{n^{2}}}$.
By Proposition~\ref{prop:CAT^omega}, $\o\cdot\spc{U}$ is $\CAT0$
and so is $\T_p^\o$ as a metric component of~$\o\cdot\spc{U}$.
\qeds

Recall that the tangent space $\T_p$ can be considered as a subset of $\T^\o_p$ (see \ref{thm:tangent-ultratangent}).
Therefore we have the following:

\begin{thm}{Corollary}\label{cor:real-tan-is}
\begin{subthm}{cor:tan-is-CBB}
If $\spc{L}$ is a complete length $\Alex{\kappa}$ space and $p\in \spc{L}$, then $\T_p$ is $\Alex{0}$.
Moreover, $\T_p$ satisfies the (1+\textit{n})-point comparison (\ref{thm:pos-config}).
\end{subthm}

\begin{subthm}{cor:tan-is-CBA}
If $\spc{U}$ is a complete length $\CAT\kappa$ space and $p\in \spc{U}$, then $\T_p$ is $\CAT0$.
Moreover, $\T_p$ satisfies the (2\textit{n}+2)-comparison (\ref{CBA-n-point}).
\end{subthm}

\end{thm}

\begin{thm}{Proposition}\label{ddo-concave}
Assume $\spc{Z}$ is a complete length $\Alex{}$ or $\CAT{}$ space
and $f\:\spc{Z}\subto\RR$ is a semiconcave locally Lipscitz subfunction.
Then for any $p\in\Dom f$, the ultradifferential $\dd^\o_p\:\T^\o_p\to\RR$ is a concave function.
\end{thm}

\parit{Proof.}
Fix a geodesic $[x^\o y^\o]$ in $\T^\o_p$.

It is sufficient to show that for any subarc $[\bar x^\o \bar y^\o]$ of $[x^\o y^\o]$
that does not contain the ends
there is a sequence of geodesics $[\bar x^n\bar y^n]$ in $n\cdot \spc{Z}$ converging to $[\bar x^\o\bar y^\o]$.

Choose any sequences $\bar x^n,\bar y^n\in n\cdot \spc{Z}$ such that $\bar x^n\to \bar x^\o$, $\bar y^n\to \bar y^\o$ as $n\to\o$.
We can assume that there is a geodesic $[\bar x^n \bar y^n]_{n\cdot \spc{Z}}$ for any $n$; see \ref{thm:almost.geod} and \ref{thm:cat-complete}.
Note that $[\bar x^n \bar y^n]$ 
converges to $[\bar x^\o \bar y^\o]$
as $n\to\o$.
The latter holds trivially in the $\CAT{}$ case,
and the $\Alex{}$ case follows from \ref{cor:unique-geod-cbb}.
\qeds

\section{Length property of tangent space}\label{halbeisen}

\begin{thm}{Theorem}\label{thm:tanCAT}
Let $\spc{U}$ be a complete length $\CAT\kappa$ space and $p\in \spc{U}$.
Then $\T_p\spc{U}$ is a length space.
\end{thm}

This theorem together with \ref{cor:real-tan-is} imply the following.

\begin{thm}{Corollary}
For any point $p$ in a complete length $\CAT\kappa$ space, the tangent space $\T_p$ is a complete length $\CAT0$ space.
\end{thm}

\parit{Proof of Theorem \ref{thm:tanCAT}.}
Since $\T_p=\Cone \Sigma_p$, it is sufficient to show that for any hinge $\hinge pxy$ such that 
$\mangle \hinge pxy<\pi$ and any $\eps>0$, there is $z\in \spc{U}$ such that 
\[\mangle \hinge pxz<\tfrac12\cdot\mangle \hinge pxy+\eps
\quad\text{and}\quad
\mangle \hinge pyz<\tfrac12\cdot\mangle \hinge pxy+\eps.\eqlbl{eq:midpoint}\]

Fix a small $\delta>0$.
Let $\bar x\in \mathopen{]}px]$ and $\bar y\in \mathopen{]}py]$ denote the points such that 
$\dist{p}{\bar x}{}=\dist{p}{\bar y}{}=\delta$.
Let $z$ denote the midpoint between $\bar x$ and $\bar y$.

Since $\delta$ is small, we can assume that 
\[\angk\kappa p{\bar x}{\bar y}<\mangle \hinge pxy+\eps.\]
By  Alexandrov's lemma (\ref{lem:alex}), we have
\[\angk\kappa p{\bar x}{z}+\angk\kappa p{\bar y}{z}< \angk\kappa p{\bar x}{\bar y}.\]
By construction,
\[\angk\kappa p{\bar x}{z}=\angk\kappa p{\bar y}{z}.\]
Applying the angle comparison (\ref{cat-hinge}), we get \ref{eq:midpoint}.
\qeds

The following example was constructed by Stephanie Halbeisen \cite{halbeisen}.
It shows that an analogous statement does not hold for $\Alex{}$ spaces.
If the dimension is finite, such examples do not exist; 
for proper spaces the question is open, see \ref{open:Halb-proper}.

\begin{thm}{Example}\label{Halbeisen's example}
There is a complete length $\Alex{}$ space $\check{\spc{L}}$
with a point $p\in\check{\spc{L}}$ such that the space of directions $\Sigma_p\check{\spc{L}}$ is not a $\pi$-length space, and therefore the tangent space $\T_p\check{\spc{L}}$ is not a length space. 
\end{thm}

\parit{Construction.}
Let $\HH$ be a Hilbert space formed by infinite sequences of real numbers $\bm{x}=(x_0,x_1,\dots)$ with the $\ell^2$-norm
$|\bm{x}|^2=\sum_i(x_i)^2$. 
Fix $\eps=0.001$ and consider two functions $f,\check f:\HH\to\RR$:
\[f(\bm{x})=|\bm{x}|,\]
\[\check f(\bm{x})
=
\max\left\{|\bm{x}|,\max_{n\ge1}\{(1+\eps)\cdot x_n-\tfrac{1}{n}\}\right\}.\] 
Both of these functions are convex and Lipschitz, therefore their graphs in $\HH\times \RR$ equipped with its length metric form infinite-dimensional Alexandrov spaces, say $\spc{L}$  and $\check{\spc{L}}$ (this is proved formally in \ref{lem:hil-con}).

Let $p$ be the origin of $\HH\times \RR$.
Note that $\check{\spc{L}}\cap\spc{L}$ is a starshaped subset of $\HH$ with center at $p$.
Further, $\check{\spc{L}}\setminus\spc{L}$ consists of a countable number of disjoint sets
\[\Omega_n=\set{(\bm{x},\check f(\bm{x}))\in\check{\spc{L}}}{(1+\eps)\cdot x_n-\tfrac{1}{n}>|\bm{x}|}.\]
Note that $\dist{\Omega_n}{p}{}>\tfrac{1}{n}$ for each $n$.
It follows that for any geodesic $[p q]$ in $\check{\spc{L}}$,
a small subinterval $[p \bar q]\subset [p q]$ 
is a straight line segment in $\HH\times\RR$, 
and also a geodesic in $\spc{L}$.
Thus we can treat $\Sigma_p\spc{L}$ and $\Sigma_p\check{\spc{L}}$ as one set, with two angle metrics $\mangle$ and $\check\mangle$.
Let us denote by $\mangle_{\HH\times \RR}$ the angle in $\HH\times\RR$.

The space $\spc{L}$  is isometric to the Euclidean cone
over $\Sigma_p\spc{L}$ with vertex at~$p$; 
$\Sigma_p\spc{L}$ is isometric to a sphere in Hilbert space with radius~$\frac{1}{\sqrt{2}}$.
In particular, $\mangle$ is the length metric of $\mangle_{\HH\times\RR}$ on $\Sigma_p{\spc{L}}$.

Therefore in order to show that $\check \mangle$ does not define a length metric on $\Sigma_p{\spc{L}}$,
it is sufficient to construct a pair of directions $(\xi_+,\xi_-)$ such that
\[\check \mangle(\xi_+,\xi_-)<\mangle(\xi_+,\xi_-).\] 
Set $\bm{e}_0=(1,0,0,\dots)$, $\bm{e}_1=(0,1,0,\dots),\dots\in \HH$. 
Consider the following two half-lines in $\HH\times \RR$:
\[\gamma_+(t)
=
\tfrac{t}{\sqrt{2}}\cdot(\bm{e}_0,1)
\quad  \text{and}\quad 
\gamma_-(t)
=
\tfrac{t}{\sqrt{2}}\cdot(-\bm{e}_0,1),
\quad t\in[0,+\infty).\] 
They form unit-speed geodesics in both $\spc{L}$ and $\check{\spc{L}}$.
Let $\xi_\pm$ be the directions of $\gamma_\pm$ at $p$.
Denote by $\sigma_n$ the half-planes in $\HH$ 
spanned by $\bm{e}_0$ and $\bm{e}_n$;
that is, $\sigma_n\z=\set{x\cdot\bm{e}_0+y\cdot\bm{e}_n}{y\ge 0}$.
Consider a sequence of $2$-dimensional sectors $Q_n=\check{\spc{L}}\cap (\sigma_n\times \RR)$. 
For each $n$, the sector $Q_n$ intersects $\Omega_n$ and is bounded by two geodesic half-lines $\gamma_\pm$.
Note that $Q_n\GHto Q$, where  $Q$ is a solid Euclidean angle
in $\EE^2$ with angle measure $\beta<\mangle(\xi_+,\xi_-)=\tfrac\pi{\sqrt{2}}$.
Indeed, $Q_n$ is path-isometric to the subset of $\EE^3$ described by
\begin{align*}
 y\ge0 \quad 
\text{and}\quad  
&
z=\max\left\{\sqrt{x^2+y^2},
(1+\eps)\cdot y-\tfrac{1}{n} \right\}
\intertext{with length metric.
Thus its limit $Q$ is path-isometric to the subset of $\EE^3$ described by}
y\ge0
\quad \text{and}\quad  
&
z=\max\left\{\sqrt{x^2+y^2},(1+\eps)\cdot y\right\}
\end{align*}
with length metric.
In particular, for any $t,\tau\ge0$, 
\begin{align*}
\dist{\gamma_+(t)}{\gamma_-(\tau)}{\check{\spc{L}}} 
&\le 
\lim_{n\to\infty}\dist{\gamma_+(t)}{\gamma_-(\tau)}{Q_n}
=
\\ 
&=\side0 \{\beta;t,\tau\}.
\end{align*}
That is, $\check\mangle(\xi_+,\xi_-) \le \beta<\mangle(\xi_+,\xi_-)$.\qeds

\begin{thm}{Lemma}\label{lem:hil-con}
Let $\HH$ be a Hilbert space,
$f\:\HH\to \RR$ be a convex Lipschitz function 
and $S\subset \HH\times \RR$ be the graph of $f$ 
equipped with  the length metric.
Then $S$ is $\Alex{0}$.
\end{thm}

\parit{Proof.}
Recall that for a subset $X\subset \HH\times \RR$, 
we will denote by $\dist{*}{*}{X}$ the
length metric on $X$.

By the theorem of Sergei Buyalo \cite{buyalo}, sharpened by the authors in \cite{alexander-kapovitch-petrunin-buyalo},
any convex hypersurface in a Euclidean space, equipped with the length metric, is non-negatively curved.
Thus it is sufficient to show that for any 4-point set $\{x_0,x_1,x_2,x_3\}\subset S$, 
there is a finite-dimensional subspace $E\subset \HH\times \RR$ 
such that $\{x_i\}\in E$ and $\dist{x_i}{x_j}{S\cap E}$ is arbitrary close to $\dist{x_i}{x_j}{S}$.

Clearly $\dist{x_i}{x_j}{S\cap E}\ge \dist{x_i}{x_j}{S}$; 
thus it is sufficient to show that for given $\eps>0$ one can choose $E$ so that 
\[\dist{x_i}{x_j}{S\cap E}
<
\dist{x_i}{x_j}{S}+\eps.
\eqlbl{eq:claim:hil-con*}.\]

For each pair $(x_i,x_j)$, choose a polygonal line $\beta_{i j}$ connecting $x_i,x_j$ that lies under $S$ (that is, outside of $\Conv S$) in $\HH\times \RR$ 
and has length at most $\dist{x_i}{x_j}{S}+\eps$.
Let $E$ be the affine hull of all the vertexes in all $\beta_{i j}$.
Thus
\[\dist{x_i}{x_j}{S\cap E}\le \length \beta_{i j}\] 
and \ref{eq:claim:hil-con*} follows.\qeds

{\sloppy 

\begin{thm}{Exercise}\label{ex:norays}
Construct a non-compact complete geodesic $\Alex{0}$ space that contains no half-lines.
\end{thm}

}

\section{Rademacher theorem}

At the end of this section we give an extension of the Rademacher theorem (see 
Section \ref{sec: speed}) to $\Alex{}$ and $\CAT{}$ spaces (\ref{thm:Rademacher-CBB+CBA}); it was proved by Alexander Lytchak \cite{lytchak:diff}. 
The following proposition is the 1-dimensional case of the extended Rademacher theorem.

Recall that differentiable curves are defined in \ref{def:diff-curv}.

\begin{thm}{Proposition}\label{prop:Rademacher-dim=1}
Let $\alpha\:\II\to \spc{Z}$ be a locally Lipschitz curve in a complete length space.
Suppose that $\spc{Z}$ is either $\Alex{}$ or $\CAT{}$.
Then $\alpha$ is differentiable almost everywhere.
\end{thm}

The following two lemmas provide sufficient conditions for existence of the one-sided derivative of a curve in $\Alex{}$ and $\CAT{}$ spaces.
The proofs of both lemmas are similar.

\begin{thm}{Lemma}\label{lem:CBB-diff-curve}
Let $\alpha\:\II\to \spc{L}$ be a $1$-Lipschitz curve in a $\Alex{}$ space.
Suppose that for some $t_0\in \II$ and any $\eps>0$, there is a point $p$ such that $\dist{\alpha(t_0)}{p}{}<\eps$ and
\[\liminf_{t\to t_0+} \frac{\distfun{p}\circ\alpha(t)-\distfun{p}\circ\alpha(t_0)}{t-t_0}>1-\eps.\]
Then the right derivative $\alpha^+(t_0)$ is defined and $|\alpha^+(t_0)|=1$.
\end{thm}

\parit{Proof.}
Without loss of generality, we may assume that $t_0=0$.
Set $x\z=\alpha(0)$.
Fix a sequence of points $p_n\to x$ such that 
\[\liminf_{t\to 0+} \frac{\dist{p_n}{\alpha(t)}{}-\dist{p_n}{x}{}}{t}\to 1\]
as $n\to\infty$.

Observe that there are sequences $\delta_n\to 0+$ and $t_n\to 0+$ such that 
\[\angk\kappa x{\alpha(s)}{p_n}>\pi-\delta_n
\quad\text{and}\quad
(1-\delta_n)\cdot s<\dist{\alpha(s)}{x}{}\le s\eqlbl{eq:ang+dist}\]
for any $s\in(0,t_n]$.

For each $n$, choose $q_n\in \Str(x)$ sufficiently close $\alpha(t_n)$ that the inequality 
\[\angk\kappa x{q_n}{p_n}>\pi-\delta_n\]
still holds (see Definition \ref{def:straight}).

Set $\gamma_n=\geod_{[xq_n]}$.
By comparison,
\begin{align*}
\angk\kappa x{\alpha(s)}{\gamma_n(s)}&\le2\cdot\pi-\angk\kappa x{p_n}{\gamma_n(s)} - \angk\kappa x{\alpha(s)}{p_n}\le
\\
&\le 2\cdot\pi-\angk\kappa x{q_n}{p_n} - \angk\kappa x{\alpha(s)}{p_n}<
\\
&<2\cdot\delta_n. 
\end{align*}
Therefore \ref{eq:ang+dist} implies that
\[\dist{\gamma_n(s)}{\alpha(s)}{}<10\cdot\delta_n\cdot (s)\]
if $s$ is a sufficiently small and positive.
That is, $\alpha^+(0)$ is defined (see Definition \ref{def:right-derivative}).
\qeds

\begin{thm}{Lemma}\label{lem:CBA-diff-curve}
Let $\alpha\:\II\to \spc{U}$ be a $1$-Lipschitz curve in a $\CAT{}$ space.
Suppose that for some $t_0\in \II$ and any $\eps>0$ there is a point $q$ such that $\dist{\alpha(t_0)}{q}{}<\eps$ and
\[\limsup_{t\to t_0+} \frac{\distfun{q}\circ\alpha(t)-\distfun{q}\circ\alpha(t_0)}{t-t_0}<-1+\eps.\]
Then the right derivative $\alpha^+(t_0)$ is defined and $|\alpha^+(t_0)|=1$.
\end{thm}

\parit{Proof.}
Without loss of generatiy we may assume that $t_0=0$.
Set $x\z=\alpha(0)$.
Fix a sequence of points $q_n\to x$ such that 
\[\liminf_{t\to 0+} \frac{\dist{q_n}{\alpha(t)}{}-\dist{q_n}{x}{}}{t}\to -1\]
as $n\to\infty$.

Observe that there are sequences $\delta_n\to 0+$ and $t_n\to 0+$ such that 
\[\angk\kappa x{\alpha(s)}{q_n}<\delta_n
\quad\text{and}\quad
(1-\delta_n)\cdot s<\dist{\alpha(s)}{x}{}\le s\eqlbl{eq:ang++dist}\]
for any $s\in(0,t_n]$.

Without loss of generality, we may assume that $\dist{x}{q_n}{}<\varpi\kappa$ for any~$n$;
in particular, the geodesic $\gamma_n=\geod_{[xq_n]}$ is uniquely defined.

By comparison,
\begin{align*}
\angk\kappa x{\alpha(s)}{\gamma_n(s)}&\le\angk\kappa x{\alpha(s)}{q_n}<
\\
&<\delta_n. 
\end{align*}
Therefore \ref{eq:ang++dist} implies that
\[\dist{\gamma_n(s)}{\alpha(s)}{}<10\cdot\delta_n\cdot s\]
if $s$ is a sufficiently small and positive.
That is, $\alpha^+(0)$ is defined (see Definition  \ref{def:right-derivative}).
\qeds

\parit{Proof of \ref{prop:Rademacher-dim=1}.}
By the standard Rademacher theorem, we may assume that $\alpha$ has an arc-length parametrization.
In particular, $\alpha$ is 1-Lipschitz.

Recall that by Theorem~\ref{thm:speed},
\[\speed_s\alpha\ae1.\eqlbl{eq:speed=1}\]

Fix a countable dense set $T\subset\II$;
given $t\in T$, let
\[h_t(s)=\dist{\alpha(t)}{\alpha(s)}{}.\]
Note that $h_t$ is $1$-Lipschitz for each $t\in T$.
Therefore, by the standard Rademacher theorem and countability of $T$ for almost all $s\in\II$,  $h_t'(s)$ is defined for all $t\in T$.

Let
\[w^+(s)\df\limsup_{\substack{t\in T\\t\to s-}} \{h'_t(s)\}.\]
Let us show that
\[w^+(s)\ae1.\eqlbl{eq:w+=1}\]
Note that once this is proved, Lemma \ref{lem:CBB-diff-curve} implies the proposition in the $\Alex{}$ case.

For a small $\eps>0$, denote by $N_\eps^+$ the set of all points $s\in\II$ such that $w^+(s)<1-\eps$.
Note that the sets $N_\eps^+$ are measurable.

Suppose $N_\eps^+$ has positive measure. 
Let $s_0\in N_\eps^+$ be a  Lebesgue point of $\alpha$.
We may assume that $\speed_{s_0}\alpha=1$ and $h_t'(s_0)$ is defined for any $t\in T$.
Suppose $t\in T$ is sufficiently close to $s_0$ and $t<s_0$.
Since $\speed_{s_0}\alpha=1$, we have 
\[h_t(s_0) \ge (s_0-t)\cdot(1 - \eps^2).
\eqlbl{hn>=}\]
Further, there is a set $A\subset [t,s_0]$ with measure at least $(1-\eps)\cdot|s_0-t|$ such that
\[h_t'(s) < 1-\eps\]
for any $s\in A_n$.
Since $h_t$ is $1$-Lipschitz, we have
\begin{align*}
h_t(s_0)&=\int_{[t,s_0]\setminus A}h_t'(s)\cdot\dd s +\int_{A}h_t'(s)\cdot\dd s\le
\\
&\le (s_0 - t)\cdot [\eps+(1 - \eps)^2].
\end{align*}
The latter contradicts \ref{hn>=}.
Thus $w^+(s)\ge1-\eps$ almost everywhere.
Since $\eps>0$ is arbitrary, \ref{eq:w+=1} follows.

In the same way we can show that 
\[w^-(s)\ae-1,
\eqlbl{eq:w-=1}\]
where 
\[w^-(s)\df\liminf_{\substack{t\in T\\t\to s+}} \{h'_t(s)\}.\]
Then Lemma \ref{lem:CBA-diff-curve} implies the proposition in the $\CAT{}$ case.
\qeds

\begin{thm}{Extended Rademacher theorem}\label{thm:Rademacher-CBB+CBA}
Let $f\:\EE^m\subto \spc{Z}$ be a locally Lipschitz submap from a Euclidean space to a complete length space $\spc{Z}$.
Suppose that $\spc{Z}$ is either $\Alex{}$ or $\CAT{}$.
Then the differential $\dd_x f$ is defined at almost all points $x\in\Dom f$.

Moreover the differential $\dd_x f$ is \index{linear differential}\emph{linear} at almost all $x$ in the following sense: 
the image $\Im f$ is a convex subcone of $\T_{f(x)}\spc{Z}$, and
there is an isometry $\iota$ from $\Im f$ to a Euclidean space such that the composition $\iota\circ\dd_x f$ is linear.
\end{thm}

The proof is a reduction to the 1-dimensional case (\ref{prop:Rademacher-dim=1}) by standard arguments \cite{kirchheim,margulis-mostow}.

\parit{Proof.}
Without loss of generality, we may assume that $\Dom f$ is bounded and $f$ is Lipschitz.

Fix a countable dense set of vectors $\{v_i\}$ in $\EE^m$.
Fix $v_i$ and a point $p\in \Dom f$.
By Proposition~\ref{prop:Rademacher-dim=1}, the value $\dd_xf(v_i)$ is defined at $x=p+t\cdot v_i$ for almost all $t$ such that $x\in \Dom f$.
It follows that $\dd_xf(v_i)$
is defined for every $i$ on a set $A$ of full measure in $\Dom f$.
Since the metric differential of $f$ is defined almost everywhere (\ref{thm:Rademacher-md}), we have that $\dd_xf(v)$ is defined for any $v$ on a set $B$ of full measure in $\Dom f$.

Applying the definitions of metric differential and differential (\ref{thm:Rademacher-md} and \ref{def:differential}), we obtain that the image of $\dd_xf$ is a weakly convex set in $\T_{f(x)}$.
It follows that $\Im\dd_xf$ is $\Alex{0}$ or $\CAT0$ if the space $\spc{Z}$ is $\Alex{}$ or $\CAT{}$ respectively.
It remains to apply Exercise \ref{mink+alex=euclid} or \ref{mink+CAT=euclid} if the space $\spc{Z}$ is $\Alex{}$ or $\CAT{}$ respectively.
\qeds

\section{Differential}

\begin{thm}{Exercise}\label{ex:d_q dist_p(v)=-<dri p q, v>-CAT}
Let $\spc{U}$ be a complete length $\CAT\kappa$ space and $p,q\in \spc{U}$.
Assume $\dist{p}{q}{}<\varpi\kappa$.
Show that 
\[(\dd_q\distfun{p}{}{})(v)=-\<\dir q p,v\>.\]

\end{thm}

\begin{thm}{Exercise}\label{ex:d_q dist_p(v)=-<dri p q, v>}
Let $\spc{L}$ be a length $\Alex{\kappa}$ space and $p,q\in \spc{L}$ be distinct points. 
Assume  $q\in \Str(p)$ or $p\in \Str(q)$
Show that 
\[(\dd_q\distfun{p}{}{})(v)=-\<\dir q p,v\>.\]

\end{thm}

\begin{thm}{Lemma}\label{lem:d(CAT)} 
Let $\spc{U}$ be a complete length $\CAT{}$ space, 
$f\:\spc{U}\subto\RR$ be a locally Lipschitz semiconcave subfunction,  
and $p\in \Dom f$.
Then $\dd_p f$ is a Lipschitz concave function on $\T_p\spc{U}$.
\end{thm}

\parit{Proof.}
Recall that the tangent space $\T_p = \T_p\spc{U}$ can be considered as a subspace of the ultratangent space $\T_p^\o$ (\ref{thm:tangent-ultratangent}).
Since $\T_p^\o$ is $\CAT0$, \ref{thm:tanCAT} implies that $\T_p$ is a convex set in $\T_p^\o$.

By \ref{ddo-concave}, $\dd_p^\o f$ is a concave function on $\T_p^\o$.
It remains to apply that $\dd_pf=(\dd_p^\o f)|_{\T_p}$ (\ref{prop:differential:ultra}).
\qeds

As it is shown in Halbeisen's example (Section \ref{halbeisen}),  
a $\Alex{}$ space  might have tangent spaces that are not length spaces; 
thus concavity of the differential $\dd_p f$ of a semiconcave function $f$ is meaningless. 
Nevertheless, as the following lemma says, the differential $\dd_p f$ of a semiconcave function always satisfies a weaker property similar to concavity (compare to \cite[Prop. 136]{plaut:survey} and \cite[4.2]{ohta}).  
In the finite dimensional case, $\dd_p f$ is concave. 

\begin{thm}{Lemma}\label{lem:ohta} 
Let $\spc{Z}$ be a complete length space,
$f\:\spc{Z}\subto\RR$ be a locally Lipschitz semiconcave subfunction,  
and $p\in \Dom f$.
Suppose that $\spc{Z}$ is either $\Alex{}$ or $\CAT{}$.
Then for any $u,v\in \T_p$, we have
\[s\cdot \sqrt{|u|^2+2\cdot\<u,v\> +|v|^2}
\ge 
(\dd_p f)(u)+(\dd_p f)(v),\]
where
\[s=\sup\set{(\dd_p f)(\xi)}{\xi\in\Sigma_p}.\]

\end{thm}

\parit{Proof.}
If $\spc{Z}$ is $\CAT{}$, then the statement follows from \ref{lem:d(CAT)}.
Indeed, let $z$ be the midpoint of a geodesic $[uv]_{\T_p}$. 
Observe that 
\[2\cdot |z|\z=\sqrt{|u|^2+2\cdot\<u,v\> +|v|^2}.\]
Since $\dd_pf$ is concave, we have that 
\[2\cdot\dd_pf(z)\ge \dd_pf(u)+\dd_pf(v).\]
It remains to choose $\xi\in \Sigma_p$ so that $\xi\cdot |z|=z$ and observe that $s\ge \dd_p(\xi)$.

Now assume $\spc{Z}$ is $\Alex{}$.
We can assume that $\alpha=\mangle(u,v)>0$, otherwise the statement is trivial.
Moreover, we can assume that $\exp_p(t\cdot u)$
 and $\exp_p(t\cdot v)$ are defined for all small $t>0$;
the latter follows since geodesic space of directions $\Sigma'_p$ is dense in $\Sigma_p$.

{

\begin{wrapfigure}{r}{34 mm}
\vskip-7mm
\centering
\includegraphics{mppics/pic-1205}
\vskip0mm
\end{wrapfigure}

Prepare a model configuration of five points: $\tilde p$, $\tilde u$, $\tilde v$, $\tilde q$, $\tilde w\in\EE^2$ such that
\begin{itemize}
\item $\mangle\hinge{\tilde p}{\tilde u}{\tilde v}=\alpha$, 
\item $\dist{\tilde p}{\tilde u}{}=|u|$, 
\item $\dist{\tilde p}{\tilde v}{}=|v|$,
\end{itemize}
}
\begin{itemize}
\item $\tilde q$ lies on an extension of $[\tilde p\tilde v]$ so that $\tilde v$ is the midpoint of $[\tilde p\tilde q]$, 
\item $\tilde w$ is the midpoint between $\tilde u$ and ${\tilde v}$.
\end{itemize}
Note that 
\[\dist{\tilde p}{\tilde w}{}
=
\tfrac{1}{2}\cdot\sqrt{|u|^2+2\cdot\<u,v\>+|v|^2}.\]

Assume that $\spc{Z}$ is geodesic.

For all small $t>0$, construct points $u_t,v_t,q_t,w_t\in \spc{Z}$ as follows:
\begin{enumerate}[(a)]
\item $v_t=\exp_p(t\cdot v)$,\quad  $q_t=\exp_p(2\cdot t\cdot v)$
\item\label{u_t}  $u_t=\exp_p(t\cdot u)$.
\item $w_t$ is the midpoint of $[u_t v_t]$.
\end{enumerate}
Clearly $\dist{p}{u_t}{}=t\cdot |u|$, $\dist{p}{v_t}{}=t\cdot|v|$, $\dist{p}{q_t}{}=2\cdot t\cdot|v|$. 
Since $\mangle(u,v)$ is defined, 
we have $\dist{u_t}{v_t}{}=t\cdot\dist{\tilde u}{\tilde v}{}+o(t)$ 
and $\dist{u_t}{q_t}{}=t\cdot\dist{\tilde u}{\tilde q}{}+o(t)$ 
(see Theorem~\ref{angle} and Section~\ref{sec:angle}).

From the point-on-side and hinge comparisons (\ref{point-on-side}$+$\ref{angle}), we have
\[\angk\kappa{v_t}p{w_t}
\ge
\angk\kappa{v_t}p{u_t}
\ge
\mangle\hinge{\tilde v}{\tilde p}{\tilde u}+\tfrac{o(t)}t\]
and
\[\angk\kappa{v_t}{q_t}{w_t}
\ge
\angk\kappa{v_t}{q_t}{u_t}
\ge
\mangle\hinge{\tilde v}{\tilde q}{\tilde u}+\tfrac{o(t)}t.\]
Clearly, 
$\mangle\hinge{\tilde v}{\tilde p}{\tilde u}+\mangle\hinge{\tilde v}{\tilde q}{\tilde u}=\pi$. 
From the adjacent angle comparison (\ref{2-sum}), 
$\angk\kappa{v_t}p{u_t}\z+\angk\kappa{v_t}{u_t}{q_t}\le \pi$.
Hence
$\angk\kappa{v_t}p{w_t}
\to
\mangle\hinge{\tilde v}{\tilde p}{\tilde w}$ as $t\to0+$
and thus 
\[\dist{p}{w_t}{}=t\cdot\dist{\tilde p}{\tilde w}{}+o(t).\]

Since $f$ is $\lambda$-concave, we have 
\begin{align*}
2\cdot f(w_t)&\ge f(u_t)+f(v_t)+\tfrac\lambda4\cdot\dist[2]{u_t}{v_t}{}=
\\
&=2\cdot f(p)
+t\cdot [(\dd_p f)(u)+(\dd_p f)(v)]+o(t).
\end{align*}
 
Applying $\lambda$-concavity of $f$, we have
\[(\dd_p f)(\dir p{w_t})
\ge 
\frac{t\cdot[(\dd_p f)(u)+(\dd_p f)(v)]
+o(t)}{2\cdot t\cdot\dist[{{}}]{\tilde p}{\tilde w}{}+o(t)}.
\eqlbl{eq:lem:ohta*}\]
The lemma follows.

\medskip

Finally, if $\spc{Z}$ is not geodesic one needs to make two adjustments in the above construction.
Namely: 
\begin{enumerate}[(i)]
\item For the geodesic $[u_t v_t]$ to be defined, in (\ref{u_t}) one has to take  $u_t\in \Str(v_t)$, $u_t\approx\exp_p(t\cdot u)$;
more precisely, 
\[\dist{u_t}{\exp_p(t\cdot u)}{}=o(t).\] 
Thus instead of $\dist{p}{u_t}{}=t\cdot|u|$ we have 
\[\dist{p}{u_t}{}=t\cdot|u|+o(t),\] and this is sufficient for the rest of proof.
\item The direction $\dir p{w_t}$ might be undefined.
Thus in the estimate \ref{eq:lem:ohta*}, instead of $\dir p{w_t}$ one should take $\dir p{w'_t}$ for some point $w_t'\in \Str(p)$ near $w_t$ (that is, $\dist{w_t}{w_t'}{}=o(t)$).
\end{enumerate}
\qedsf


\section{Definition of gradient}\label{sec:grad-def}

\begin{thm}{Definition of gradient}\label{def:grad} 
Let $\spc{X}$ be a length space with defined angles and  
$f\:\spc{X}\subto\RR$ be a subfunction.
Suppose  for a point
$p\in\Dom f$ the differential $\dd_p f\:\T_p\to\RR$ is defined.

A tangent vector $g\in \T_p$ is called a 
\index{gradient}\emph{gradient of $f$ at $p$} 
(briefly,  $g\z=\nabla_p f$\index{$\nabla$}) if
\begin{subthm}{}
$(\dd_p f)(w)\le \<g,w\>$ for any $w\in \T_p$, and
\end{subthm}

\begin{subthm}{}
$(\dd_p f)(g) = \<g,g\> .$
\end{subthm}
\end{thm}

\begin{wrapfigure}{r}{34 mm}
\vskip-10mm
\centering
\includegraphics{mppics/pic-1215}
\vskip0mm
\end{wrapfigure}

\begin{thm}{Example}\label{l-inf-grad}
Consider the Euclidean plane with standard $(x,y)$-coordinates.
Then the function $f\:(x,y)\mapsto -|x|-|y|$ is concave;
its gradient field is sketched on the figure.

If a point does not lie on an axis, then its gradient has length $\sqrt2$ and 
takes one of four values $(\pm1,\pm1)$ depending on the quadrant of the point.
At the origin the gradient vanish, 
and on the on the remaining parts of the $x$-axis and $y$-axis it is $(\pm1,0)$ and $(0,\pm 1)$ respectively. 
\end{thm}

\begin{thm}{Exercise}\label{ex:grad(dist)}
Let $\spc{U}$ be a complete length $\CAT0$ space.
Show that
\[\nabla_p(-\distfun q)=\dir pq\]
for any pair of distinct points $p,q\in\spc{U}$.
\end{thm}

\begin{thm}{Existence and uniqueness of the gradient}\label{thm:ex-grad} 
Let $\spc{Z}$ be a complete length space
and $f\:\spc{Z}\subto\RR$ be a 
locally Lipschitz 
and 
semiconcave subfunction.
Suppose that $\spc{Z}$ is either $\Alex{}$ or $\CAT{}$.
Then for any point $p\in \Dom f$, there is a unique gradient $\nabla_p f\in \T_p$.
\end{thm}

\parit{Proof; uniqueness.} 
If $g,g'\in \T_p$ are two gradients of $f$,
then 
\begin{align*}
\<g,g\>
&=(\dd_p f)(g)\le \<g,g'\>,
&
\<g',g'\>
&=(\dd_p f)(g')\le \<g,g'\>.
\end{align*}
Therefore,
\[\dist[2]{g}{g'}{}=\<g,g\>-2\cdot\<g,g'\>+\<g',g'\>\le0.\] 
It follows that $g=g'$.

\parit{Existence.} 
Note first that if $\dd_p f\le 0$, then one can take $\nabla_p f=\0$.

Otherwise, if $s=\sup\set{(\dd_p f)(\xi)}{\xi\in\Sigma_p}>0$, 
it is sufficient to show that there is  $\overline{\xi}\in \Sigma_p$ such that 
\[
(\dd_p f)\left(\overline{\xi}\right)=s.
\eqlbl{overlinexi}
\]
Indeed, suppose $\overline{\xi}$ exists.
Then applying Lemma~\ref{lem:ohta} for $u=\overline{\xi}$, $v=\eps\cdot w$ with $\eps\to0+$, 
we get
\[(\dd_p f)(w)\le \<w,s\cdot\overline{\xi}\>\] 
for any $w\in\T_p$;
that is, $s\cdot\overline{\xi}$ is the gradient at $p$.

Take a sequence of directions $\xi_n\in \Sigma_p$, such that $(\dd_p f)(\xi_n)\to s$.
Applying Lemma~\ref{lem:ohta} for $u=\xi_n$, $v=\xi_m$, we get
\[s
\ge
\frac{(\dd_p f)(\xi_n)+(\dd_p f)(\xi_m)}{\sqrt{2+2\cdot\cos\mangle(\xi_n,\xi_m)}}.\]
Therefore $\mangle(\xi_n,\xi_m)\to0$ as $n,m\to\infty$;
that is, the sequence $\xi_n$ is Cauchy.
Clearly $\overline{\xi}=\lim_n\xi_n$ satisfies \ref{overlinexi}.
\qeds

\begin{thm}{Exercise}\label{ex:|grad|=1}
Let $p$ be a point in a complete length $\Alex{}$ space.
Show that 
\[|\nabla_x\distfun p|=1\]
for $x$ in a dense G-delta subset.
\end{thm}

\section{Calculus of gradient}\label{sec:grad-calculus}

The next lemma states that the gradient points 
in the direction of maximal slope; 
moreover, if the slope in the given direction is almost maximal, then it is almost the direction of the gradient.

\begin{thm}{Lemma}\label{lem:alm-grad}
Let $\spc{Z}$ be a complete length space,
$f\:\spc{Z}\subto\RR$ be locally Lipschitz and semiconcave, 
and $p\in \Dom f$.
Suppose that $\spc{Z}$ is either $\Alex{}$ or $\CAT{}$.

Assume $|\nabla_p f|>0$;
let $\overline{\xi}=\tfrac{1}{|\nabla_p f|}\cdot\nabla_p f$.
Then:
\begin{subthm}{near-grad} If for some $v\in\T_p$, we have 
\[|v|\le 1+\eps
\quad
\text{and}
\quad
(\dd_p f)(v) > |\nabla_p f|\cdot(1-\eps),
\]
then
\[\dist{\overline{\xi}}{v}{}<100\cdot\sqrt{\eps}.\]
\end{subthm}

\begin{subthm}{conv-to-grad} 
If $v_n\in \T_p$ is a sequence of vectors such that 
\[\limsup_{n\to\infty} |v_n|\le 1\quad  
\text{and}\quad  \liminf_{n\to\infty}(\dd_p f)(v_n)\ge |\nabla_p f|,\] 
then 
\[\lim_{n\to\infty} v_n=\overline{\xi}.\]
\end{subthm}

\begin{subthm}{alm-max} $\overline{\xi}$ is the unique maximum direction for the restriction $\dd_p f|_{\Sigma_p}$. 
In particular, 
\[|\nabla_p f|=\sup\set{\dd_p f}{\xi\in\Sigma_p f}.\]
\end{subthm}
\end{thm}

\parit{Proof.}
According to the definition of gradient,
\begin{align*}
 |\nabla_p f|\cdot(1-\eps)
&<
(\dd_p f)(v)
\le
\\
&\le\<v,\nabla_p f\>
=
\\
&=
|v|\cdot|\nabla_p f|\cdot\cos\mangle(\nabla_p f,v).
\end{align*}
Thus 
$
|v|>1-\eps$
and
$
\cos\mangle(\nabla_p f,v)>\tfrac{1-\eps}{1+\eps}.
$
Hence  \ref{SHORT.near-grad}.

Statements \ref{SHORT.conv-to-grad} and \ref{SHORT.alm-max} follow directly from \ref{SHORT.near-grad}).
\qeds

As a corollary of the above lemma and Proposition~\ref{prop:conv-comp} we obtain the following: 

\begin{thm}{Chain rule} 
Let $\spc{Z}$ be a complete length space, 
$f\:\spc{Z}\subto \RR$ be a semiconcave function,
and $\phi\:\RR\to\RR$ be a nondecreasing semiconcave function.
Suppose that $\spc{Z}$ is either $\Alex{}$ or $\CAT{}$.
Then $\phi\circ f$ is semiconcave and  $\nabla_x(\phi\circ f)=\phi^+(f(x))\cdot\nabla_x f$ for any $x\in\Dom f$.
\end{thm}

\begin{wrapfigure}{r}{33 mm}
\vskip-0mm
\centering
\includegraphics{mppics/pic-1210}
\vskip0mm
\end{wrapfigure}

The following inequalities describe an important property of the ``gradient
vector field''.

\begin{thm}{Lemma} 
\label{lem:grad-lip}
Let $\spc{Z}$ be a complete length space,   
$f\:\spc{Z}\subto\RR$ satisfy $f''+\kappa\cdot f\le \lambda$ for some $\kappa,\lambda\in\RR$.
 Let $[p q]\subset \Dom f$. 
and $\ell=\dist{p}{q}{}$.
Suppose that $\spc{Z}$ is either $\Alex{}$ or $\CAT{}$.
Then
\[\<\dir pq,\nabla_p f\>\ge
\frac
{{f(q)}-{f(p)\cdot\cs\kappa\ell}-\lambda\cdot\md\kappa\ell}
{\sn\kappa\ell}.\]

In particular, 
\begin{subthm}{lem:grad-lip:lam=0}
if $\kappa=0$, 
\[\<\dir pq,\nabla_p f\>\ge
\left({f(q)}-{f(p)}-\tfrac\lambda2\cdot\ell^2\right)/\ell;\]
\end{subthm}

\begin{subthm}{} if $\kappa=1$, $\lambda=0$ we have
\[\<\dir pq,\nabla_p f\>\ge
\left(f(q)-f(p)\cdot\cos\ell\right)/\sin\ell;\]
\end{subthm}

\begin{subthm}{} if $\kappa=-1$, $\lambda=0$ we have
\[\<\dir pq,\nabla_p f\>\ge
\left(f(q)-f(p)\cdot\cosh\ell\right)/\sinh\ell.\]
\end{subthm}
\end{thm}

\parit{Proof of \ref{lem:grad-lip}.} 
Note that 
\[\geod_{[p q]}(0)=p,\quad 
\geod_{[p q]}(\ell)=q,\quad
(\geod_{[p q]})^+(0)\z=\dir pq.\]
Thus,
\begin{align*}
\<\dir pq,\nabla_p f\>
&\ge 
d_p f(\dir pq)=
\\
&=
(f\circ\geod_{[p q]})^+(0)
\ge
\\
&\ge
\frac
{{f(q)}-{f(p)\cdot\cs\kappa\ell}-\lambda\cdot\md\kappa\ell}
{\sn\kappa\ell}.
\end{align*}
\qedsf

The following corollary states that the gradient vector field is monotonic in  a sense similar to the definition of \index{monotone operator}\emph{monotone operators} \cite{phelps}.

\begin{thm}{$\bm{\lambda}$-Monotonicity of gradient}
\label{cor:grad-lip}
Let $\spc{Z}$ be a complete length space, 
$f\:\spc{Z}\subto\RR$ be locally Lipschitz and $\lambda$-concave 
and $[p q]\subset \Dom f$.
Suppose that $\spc{Z}$ is either $\Alex{}$ or $\CAT{}$.
Then
\[
\<\dir p q,\nabla_p f\>
+
\<\dir q p,\nabla_q f\>
\ge 
-\lambda\cdot\dist[{{}}]{p}{q}{}.
\]

\end{thm}

\parit{Proof.}
Add two inequalities from \ref{lem:grad-lip:lam=0}.
\qeds

\begin{thm}{Lemma}\label{lem:close-grad}
Let $\spc{Z}$ be a complete length space, 
$f,g\:\spc{Z}\subto\RR$, 
and $p\in\Dom f\cap\Dom g$.
Suppose that $\spc{Z}$ is either $\Alex{}$ or $\CAT{}$.

Then 
\[\dist[2]{\nabla_p f}{\nabla_p g}{\T_p}
\le 
s\cdot(|\nabla_p f|+|\nabla_p g|),\]
where
\[s=\sup\set{|(\dd_p f)(\xi)-(\dd_p g)(\xi)|}{\xi\in\Sigma_p}.\]

In particular, if $f_n\:\spc{Z}\subto\RR$ is a sequence of locally Lipschitz and semiconcave subfunctions,
$p\in \Dom f_n$ for each $n$, 
and $\dd_p f_n$ converges uniformly on ${\Sigma_p}$, 
then the sequence $\nabla_p f_n\in \T_p$ converges.
\end{thm}

\parit{Proof.}
Clearly for any $v\in \T_p$, we have 
\[|(\dd_p f)(v)-(\dd_p g)(v)|\le s\cdot|v|.\]
From the definition of gradient (\ref{def:grad}) we have:
\begin{align*}
&(\dd_p f)(\nabla_p g)\le\<\nabla_p f,\nabla_p g\>,
&&(\dd_p g)(\nabla_p f)\le\<\nabla_p f,\nabla_p g\>,
\\
&(\dd_p f)(\nabla_p f)=\<\nabla_p f,\nabla_p f\>,
&&(\dd_p g)(\nabla_p g)=\<\nabla_p g,\nabla_p g\>.
\end{align*}
Therefore,
\begin{align*}
\dist[2]{\nabla_pf}{\nabla_pg}{}
&=\<\nabla_p f,\nabla_p f\>+\<\nabla_p g,\nabla_p g\>-2\cdot\<\nabla_p f,\nabla_p g\>
\le
\\
&\le (\dd_p f)(\nabla_p f)+(\dd_p g)(\nabla_p g)-
\\
&\quad -(\dd_p f)(\nabla_p g)-(\dd_p g)(\nabla_p f)
\le
\\
&\le s\cdot(|\nabla_p f|+|\nabla_p g|).
\end{align*}
\qedsf

\begin{thm}{Exercise}\label{ex:df(v)=<grad f,v>}
Let $\spc{L}$ be a complete length $\Alex\kappa$ space,
the function
$f\:\spc{L}\to\RR$ be semiconcave and locally Lipschitz,
and
$\alpha\:\II\to\spc{L}$ be a Lipschitz curve.
Show that 
\[\<\nabla_{\alpha(t)}f,\alpha^+(t)\>
=
(\dd_{\alpha(t)}f)(\alpha^+(t))\]
for almost all $t\in\II$.

\end{thm}

\section{Semicontinuity of \textbar gradient\textbar}\label{sec:grad-semicont}

In this section we collect a few consequences of the following lemma.

\begin{thm}{Ultralimit of \textbar gradient\textbar} \label{lem:gradcon}
Assume that
\begin{itemize}
\item $(\spc{Z}_n)$ is a sequence of complete length spaces and $(\spc{Z}_n,p_n) \to (\spc{Z}_\o,p_\o)$ as $n\to\o$.
Suppose that all $\spc{Z}_n$ are either $\Alex{}$ or $\CAT{}$.
\item $f_n\:\spc{Z}_n\subto \RR$ and $f_\o\:\spc{Z}_\o\subto \RR$ are locally Lipschitz and $\lambda$-concave, and $f_n\to f_\o$ as $n\to\o$.
\item $x_n\in\Dom f_n$ and $x_n\to x_\o\in \Dom f_\o$ as $n\to\o$.
\end{itemize}
Then 
\[|\nabla_{x_\o} f_\o|
\le 
\lim_{n\to \o} |\nabla_{x_n} f_n|.\]

\end{thm}

\parbf{Remarks.} The inequality might be strict.
For example, consider $\spc{Z}_n\z=\RR$, $f_n(x)=-|x|$ and $x_n\to 0+$.

 From the convergence of gradient curves (proved later in \ref{ultr-lim-g-curve}), 
one can deduce the following slightly stronger statement.
 
\begin{thm}{Proposition}\label{prop:lim|grad|=|grad|}
Assume that
\begin{itemize}
\item $\spc{Z}_n$ is a sequence of complete length spaces and $(\spc{Z}_n,p_n) \to (\spc{Z}_\o,p_\o)$ as $n\to\o$.
Suppose that all $\spc{Z}_n$ are either $\Alex{}$ or $\CAT{}$.
\item $f_n\:\spc{Z}_n\subto \RR$ and $f_\o\:\spc{Z}_\o\subto \RR$ are locally Lipschitz and $\lambda$-concave and $f_n\to f_\o$ as $n\to\o$.
\end{itemize}
Then 
\[|\nabla_{x_\o} f_\o|
=
\inf \{\lim_{n\to \o} |\nabla_{x_n} f_n|\},\]
where infimum is taken for all sequences $x_n\in\Dom f_n$ such that $x_n\to x_\o\in \Dom f_\o$ as $n\to\o$.
\end{thm}

\parit{Proof of \ref{lem:gradcon}.} 
Fix an $\eps>0$ and choose $y_\o\in \Dom f_\o$ sufficiently close to $x_\o$ that 
\[|\nabla_{x_\o} f_\o|-\eps<\frac{f_\o(y_\o)-f_\o(x_\o)}{\dist{x_\o}{y_\o}{}}.\]
Choose $y_n\in \spc{Z}_n$ such that $y_n\to y_\o$ as $n\to\o$. 
Since $\dist{x_\o}{y_\o}{}$ is sufficiently small, the $\lambda$-concavity of $f_n$ implies that
\[ |\nabla_{x_\o} f_\o|-2\cdot\eps
<
(\dd_{x_n}f_n)(\dir{x_n}{y_n})\]
for $\o$-almost all $n$.
Hence
\[
|\nabla_{x_\o} f_\o|-2\cdot\eps
\le 
\lim_{n\to \o} |\nabla_{x_n} f_n|.\]
Since $\eps>0$ is arbitrary, the proposition follows.
\qeds

Note that the distance-preserving map $\iota\:\spc{Z}\hookrightarrow \spc{Z}^\o$ induces an embedding 
\[\dd_p\iota\:\T_p \spc{Z}\hookrightarrow \T_p \spc{Z}^\o.\]
Thus, we can (and will) consider $\T_p \spc{Z}$ as a subcone of $\T_p \spc{Z}^\o$.

\begin{thm}{Corollary}\label{nablaf=mablaf^o}
Let $\spc{Z}$ be a complete length space 
and $f\:\spc{Z}\subto\RR$ be a locally Lipschitz semiconcave subfunction.
Suppose that $\spc{Z}$ is either $\Alex{}$ or $\CAT{}$.
Then 
\[\nabla_x f=\nabla_x f^\o.\]
for any point $x\in\Dom f$.
\end{thm}

\parit{Proof.} 
Note that 
\[
\begin{matrix}
\spc{Z}&\subset&\spc{Z}^\o\\
\cup&&\cup\\
\Dom f&\subset &\Dom f^\o
\end{matrix}.
\]
Applying \ref{lem:gradcon} for $\spc{Z}_n=\spc{Z}$ and $x_n=x$, we get that $|\nabla_x f|\ge|\nabla_x f^\o|$.

On the other hand, $f=f^\o|_{\spc{Z}}$, hence $\dd_p f=\dd_p f^\o|_{\T_p \spc{Z}}$.
Thus from \ref{alm-max},
$|\nabla_x f|\le|\nabla_x f^\o|$. 
Therefore
\[
|\nabla_x f|=|\nabla_x f^\o|
\eqlbl{gradfgradultraf}
\]
 for any $x\in \spc{Z}$.

Further,
\begin{align*}
|\nabla_x f|^2&=(\dd_x f)(\nabla_x f)=\\
&=\dd_xf^\o(\nabla_x f)\le\\ 
&\le\<\nabla_x f^\o,\nabla_x f\>=\\
&=|\nabla_x f^\o|\cdot|\nabla_x f|\cdot\cos\mangle(\nabla_x f^\o,\nabla_x f).
\end{align*}
Together with \ref{gradfgradultraf}, this implies $\mangle(\nabla_x f^\o,\nabla_x f)=0$ and the statement follows.
\qeds

\begin{thm}{Semicontinuity of \textbar gradient\textbar}\label{cor:gradlim} 
Let $\spc{Z}$ be a complete length space 
and $f\:\spc{Z}\subto\RR$ be a locally Lipschitz semiconcave subfunction.
Suppose that $\spc{Z}$ is either $\Alex{}$ or $\CAT{}$.
Then the function $x\mapsto|\nabla_x f|$  is lower-continuous;
that is, for any sequence $x_n\to x\in \Dom f$, we have 
\[|\nabla_x f|\le \liminf_{n\to \infty} |\nabla_{x_n} f|.\]
\end{thm}

\noi\textit{Proof.} 
According to \ref{nablaf=mablaf^o}, $|\nabla_x f|=|\nabla_x f^\o|$. 
Applying \ref{lem:gradcon} for $x_n\to x$, we obtain
\[\lim_{n\to\o}|\nabla_{x_n}f|
\ge
|\nabla_x f^\o|
=
|\nabla_x f|.\]
The same holds for an arbitrary subsequence of $x_n$ --- hence the result. \qeds

\section{Polar vectors}

Here we give a corollary of Lemma \ref{lem:close-grad}.
It will be used to prove basic properties of the
tangent space.

\begin{thm}{Anti-sum lemma}\label{lem:minus-sum} 
Let $\spc{L}$ be a complete length $\Alex{}$ space and $p\in \spc{L}$.

Given two vectors $u,v\in \T_p$, there is a unique vector $w\in \T_p$ such that
\[\<u,x\>+\<v,x\>+\<w,x\>\ge 0\]
for any $x\in \T_p$, and
\[\<u,w\>+\<v,w\>+\<w,w\>=0.\]

\end{thm}

If $\T_p$ were a length space, then the lemma would follow from the existence  of the gradient (\ref{thm:ex-grad}), applied to the function $\T_p\to \RR$ defined by $x\mapsto -(\<u,x\>+\<v,x\>)$.
However, the tangent space $\T_p$ might be not a length space; see  Halbeisen's example \ref{Halbeisen's example}.

Applying the above lemma for $u=v$, we have the following statement.

\begin{thm}{Existence of polar vector}\label{cor:polar}
Let $\spc{L}$ be a complete length $\Alex{}$ space 
and $p\in \spc{L}$. 
Given a vector $u\in \T_p$,  there is a unique vector $u^*\in\T_p$ such that $\<u^*,u^*\>+\<u,u^*\>= 0$ and
$u^*$ is \index{polar vectors}\emph{polar} to $u$;
that is,
 $\<u^*,x\>+\<u,x\>\ge 0$ for any $x\in \T_p$.

In particular, for any vector $u\in \T_p$ there is a polar vector $u^*\in\T_p$ such that
$|u^*|\le |u|$.
\end{thm}

Milka's lemma provides a refinement of this statement;
it states that in the finite-dimensional $\Alex{}$ space, any tangent vector $u$ has a polar vector $u^*$ such that $|u^*|= |u|$.

\begin{thm}{Example}
Let $\spc{L}$ be the upper half space of the Euclidean space $\EE^n$;
that is, $\spc{L}=\{(x_1,\ldots, x_n)\in \EE^n\mid x_n\ge 0\}$.
It is a complete length $\Alex{0}$ space.
For $p=0$, the tangent space $\T_p$ can be canonically identified with $\spc{L}$.
Then any vector $u=(v_1,\ldots, v_n)\in \T_p$
a unique polar vector such that $|u^*|=|u|$ which is $u^*=(-v_1, \ldots, -v_{n-1}, v_n)$.
However, if $v_n\ne 0$, then $u$ has other polar vectors, in particular $(-v_1, \ldots, -v_{n-1}, 0)$.
\end{thm}

It is instructive to solve the following exercise before reading the proof of \ref{lem:minus-sum}.

\begin{thm}{Exercise}\label{ex:d dist(grad)<0}
Let $\spc{L}$ be a complete length $\Alex\kappa$ space and $a,b,p$
be mutually distinct points in $\spc{L}$.
Prove that 
\[(\dd_p\distfun{a}{}{})(\nabla_p\distfun{b}{}{})
\le\cos\angk\kappa pab.\]
\end{thm}

\parit{Proof of \ref{lem:minus-sum}.}
Choose two sequences of points $a_n,b_n\in \Str(p)$ such that $\dir{p}{a_n}\to u/|u|$ and $\dir{p}{b_n}\to v/|v|$ as $n\to \infty$.
Consider a sequence of functions 
\[f_n=|u|\cdot\distfun{a_n}{}{}+|v|\cdot\distfun{b_n}{}{}.\]
According to Exercise~\ref{ex:d_q dist_p(v)=-<dri p q, v>}, 
\[(\dd_pf_n)(x)=-|u|\cdot\<\dir{p}{a_n},x\>-|v|\cdot\<\dir{p}{b_n},x\>.\]
Thus we have the following uniform convergence for all $x\in\Sigma_p$:
\[(\dd_pf_n)(x)\xto[n\to\infty]{}-\<u,x\>-\<v,x\>.\]
According to Lemma~\ref{lem:close-grad}, 
the sequence $\nabla_pf_n$ converges.
Let 
\[w=\lim_{n\to\infty}\nabla_pf_n.\]

By the definition of gradient,
\[\begin{aligned}
\<w,w\>&=\lim_{n\to\infty}\<\nabla_pf_n,\nabla_pf_n\>=
&&&
\<w,x\>&=\lim_{n\to\infty}\<\nabla_pf_n,x\>\ge
\\
&=\lim_{n\to\infty}(\dd_p f_n)(\nabla_p f_n)
=
&&&
&\ge
\lim_{n\to\infty}(\dd_pf_n)(x)
=
\\
&=-\<u,w\>-\<v,w\>,
&&&
&=-\<u,x\>-\<v,x\>.
\end{aligned}\]
\qedsf

\section{Linear subspace of tangent space}

\begin{thm}{Definition}\label{def:opp+Lin}
Let $\spc{L}$ be a complete length $\Alex{\kappa}$ space, $p\in \spc{L}$ and $u,v\in\T_p$.
We say that vectors $u$ and $v$ are \index{opposite vectors}\emph{opposite}\label{def:opposite:page} to each other, (briefly, $u+v=0$) if $|u|=|v|=0$ or $\mangle(u,v)=\pi$ and $|u|=|v|$.

The subcone
\[\Lin_p=\set{v\in\T_p}{\exists\ w\in\T_p\quad \text{such that}\quad w+v=0}\]
will be called the \index{linear subspace}\emph{linear subcone} of $\T_p$.
\end{thm}

The reason for the term ``linear'' will become evident in Theorem~\ref{thm:lin-subcone}.

\begin{thm}{Proposition}\label{prop:opposite}
Let $\spc{L}$ be a complete length $\Alex{}$ space and $p\in \spc{L}$.
Given two vectors $u,v\in\T_p$, the following statements are equivalent:
\begin{subthm}{opposite} $u+v=0$;
\end{subthm}
\begin{subthm}{<x,u>} $\<u,x\>+\<v,x\>=0$ for any $x\in\T_p$;
\end{subthm}
\begin{subthm}{<xi,u>} $\<u,\xi\>+\<v,\xi\>=0$ for any $\xi\in\Sigma_p$.
\end{subthm}
\end{thm}

\parit{Proof.}
The condition $u+v=0$ is equivalent to 
\[\<u,u\>=-\<u,v\>=\<v,v\>;\]
thus 
\ref{SHORT.<x,u>}$\Rightarrow$\ref{SHORT.opposite}.
Since $\T_p$ is isometric to a subset of $\T^\o_p$,
the splitting theorem (\ref{thm:splitting}) applied to $\T_p^\o$
gives \ref{SHORT.opposite}$\Rightarrow$\ref{SHORT.<x,u>}.

The equivalence  \ref{SHORT.<x,u>}$\Leftrightarrow$\ref{SHORT.<xi,u>} is trivial.
\qeds

\begin{thm}{Proposition}\label{prop:two-opp}
Let $\spc{L}$  be a complete length $\Alex{}$ space and $p\in \spc{L}$.
Then for any three vectors $u,v,w\in\T_p$, if $u+v=0$ and $u+ w=0$ then $v=w$.
\end{thm}

\parit{Proof.}
By Proposition~\ref{prop:opposite}, both $v$ and $w$ satisfy the condition in corollary~\ref{cor:polar}. 
Hence the result.\qeds

Let $u\in \Lin_p$; that is $u+v=0$ for some $v\in\T_p$.
Given $s<0$, let 
\[s\cdot u\df (-s)\cdot v.\]
This way we define multiplication of any vector in $\Lin_p$ by any real number (positive and negative).
Proposition~\ref{prop:two-opp} implies that such multiplication is uniquely defined.

\begin{thm}{Theorem}\label{thm:lin-subcone}
Let $\spc{L}$  be a complete length $\Alex{\kappa}$ space and $p\in \spc{L}$. 
Then $\Lin_p$ is a subcone of $\T_p$ isometric to a Hilbert space.
\end{thm}

Before proving the theorem, 
let us give a corollary.

\begin{thm}{Corollary}\label{cor:euclid-subcone}
Let $\spc{L}$  be a complete length $\Alex{\kappa}$ space
and $p\in \Str(x_1,\dots,x_n)$.
Then there is a subcone $E\subset \T_p$ that is isometric to a Euclidean space such that $\ddir p{x_i}\in E$ for every $i$.
\end{thm}

\parit{Proof.} 
By the definition of $\Str(x_1,\dots,x_n)$ (\ref{def:straight}), $\ddir{p}{x_i}\in \Lin_p$ for each $i$.
It remains to apply Theorem~\ref{thm:lin-subcone}.
\qeds

The main difficulty in the proof of Theorem~\ref{thm:lin-subcone} comes from the fact that in general $\T_p$ is not a length space;
see Halbeisen's example (\ref{Halbeisen's example}).
If the tangent space were a length space, the statement would follow directly from the splitting theorem (\ref{thm:splitting}).
In fact the proof of Theorem~\ref{thm:lin-subcone}
   is very circuitous --- we use the construction of the gradient, as well as the splitting theorem, namely its corollary (\ref{cor:splitting}).
Thus in order to understand our proof one needs to read most of Chapter~\ref{chap:grad}.

\parit{Proof of \ref{thm:lin-subcone}.}
First we show that $\Lin_p$ is a complete geodesic $\Alex0$ space.

Recall that $\T^\o_p$ is a complete geodesic $\Alex0$ space (see \ref{obs:ultralimit-is-geodesic} and \ref{thm:tan-is-CBB}) and $\Lin_p$ is a closed subset of $\T^\o_p$.
Thus, it is sufficient to show that the metric on $\Lin_p$ inherited from $\T^\o_p$ is a length metric.

Fix two vectors $x,y\in\Lin_p$.
Let $u$ and $v$ be such that $u+\tfrac{1}{2}\cdot x=0$ 
and $v+\tfrac{1}{2}\cdot y=0$.
Apply Lemma~\ref{lem:minus-sum} 
to the vectors $u$ and $v$;
let $w\in \T_p$ denote the obtained tangent vector.
\begin{clm}{}\label{clm:w-mid(xy)}
$w$ is a midpoint of $[x y]$.
\end{clm}

Indeed, according to Lemma~\ref{lem:minus-sum}, 
\begin{align*}
|w|^2
&=
-\<w,u\>-\<w,v\>
=
\\
&=
\tfrac{1}{2}\cdot\<w,x\>+\tfrac{1}{2}\cdot\<w,y\>.
\\
\intertext{Therefore}
\dist[2]{x}{w}{}+\dist[2]{w}{y}{}
&=2\cdot|w|^2+|x|^2+|y|^2-2\cdot\<w,x\>-2\cdot\<w,y\>=
\\
&=|x|^2+|y|^2-\<w,x\>-\<w,y\>\le
\\
&\le |x|^2+|y|^2+\<u,x\>+\<v,x\>+\<u,y\>+\<v,y\>=
\\
&=\tfrac{1}{2}\cdot|x|^2+\tfrac{1}{2}\cdot|y|^2-\<x,y\>=
\\
&=\tfrac{1}{2}\cdot\dist[2]{x}{y}{}.
\end{align*}
Thus $\dist{x}{w}{}=\dist{w}{y}{}=\tfrac{1}{2}\cdot\dist[{{}}]{x}{y}{}$ and \ref{clm:w-mid(xy)} follows.
\claimqeds

Note that for any $v\in\Lin_p$ there is a line $\ell$ that contains $v$ and $\0$.
Therefore by \ref{cor:splitting}, $\Lin_p$ is isometric to a Hilbert space.
\qeds

\section{Comments}

\begin{thm}{Open question}\label{open:Halb-proper}
Let $\spc{L}$ be a proper length $\Alex\kappa$ space.
Is it true that for any $p\in \spc{L}$, the tangent space $\T_p$ is a length space?
\end{thm}

\chapter{Dimension of CAT spaces}\label{chap:web+bary}

In this chapter we discuss constructions introduced by Bruce Kleiner \cite{kleiner}.

The material of this chapter is used mostly for $\CAT{}$ spaces, 
but the results in section~\ref{sec:web-general} find applications for finite-dimensional $\Alex{}$ spaces as well.

\section{The case of complete geodesic spaces}\label{sec:web-general}

The following construction gives a $\kay$-dimensional submanifold 
for a given ``nondegenerate'' array of $\kay+1$ strongly convex functions.

\begin{thm}{Definition}\label{def:ordung}
For two real arrays $\bm{v}$, $\bm{w}\in \RR^{\kay+1}$,
$\bm{v}=(v^0,v^1,\dots,v^\kay)$ 
and 
$\bm{w}=(w^0,w^1,\dots,w^\kay)$, 
we will write
$\bm{v}\succcurlyeq\bm{w}$ if $v^i\ge w^i$ for each $i$.
\end{thm}

Given a subset $Q\subset \RR^{\kay+1}$, 
denote by $\Up Q$ \label{PAGE.def:Up}
the smallest upper set containing $Q$, 
and by 
$\Min Q$ the set of minimal elements of $Q$ with respect to $\succcurlyeq$;
that is,
\begin{align*}
\Up Q 
&=
\set{\bm{v}\in\RR^{\kay+1}}{\exists\, \bm{w}\in Q\ \text{such that}\ \bm{v}\succcurlyeq\bm{w}},
\\
\Min Q 
&=
\set{\bm{v}\in Q}{\text{if}\ \bm{v}\succcurlyeq\bm{w}\in Q\ \text{then}\ \bm{w}=\bm{v}}.
\end{align*}

\begin{thm}{Definition}\label{def:web}
Let $\bm{f}=(f^0,f^1,\dots,f^\kay)\:\spc{X}\to \RR^{\kay+1}$ be a function array on a metric space $\spc{X}$.
The set 
\[\Web\bm{f}
\df
\bm{f}^{-1}\left[\Min\bm{f}(\spc{X})\right]
\subset 
\spc{X}\] 
will be called the \index{web}\emph{web} of $\bm{f}$.
\end{thm}

Given an array $\bm{f}=(f^0,f^1,\dots,f^\kay)$,
we denote by $\bm{f}^{\without i}$ the subarray of $\bm{f}$ with $f^i$ removed;
that is, 
\[\bm{f}^{\without i\,}\df(f^0,\dots,f^{i-1},f^{i+1},\dots,f^\kay).\]
Clearly 
$\Web\bm{f}^{\without i}\subset \Web\bm{f}$.
Define the \index{web!inner web}\emph{inner web} of $\bm{f}$ 
as 
\[\InWeb\bm{f}
=
\Web\bm{f}\setminus\left(\bigcup_{i}\Web\bm{f}^{\without i}\right).\]

We say that a function array is \emph{nondegenerate} 
if $\InWeb\bm{f}\ne \emptyset$.

\parbf{Example.} 
If $\spc{X}$ is a geodesic space, 
then $\Web(\distfun{x}{}{},\distfun{y}{}{})$ is the union of all geodesics from $x$ to $y$, and 
\[\InWeb(\distfun{x}{}{},\distfun{y}{}{})=\Web(\distfun{x}{}{},\distfun{y}{}{})\setminus\{x,y\}.\]

\parbf{Barycenters.}
Let us denote by $\Delta^\kay\subset \RR^{\kay+1}$\index{$\Delta^m$} 
the \index{standard simplex}\emph{standard $\kay$-simplex}; 
that is, $\bm{x}=(x^0,x^1,\dots,x^\kay)\in\Delta^\kay$ if $\sum_{i=0}^\kay x^i=1$ and $x^i\ge0$ for all $i$.

Let $\spc{X}$ be a metric space 
and $\bm{f}=(f^0,f^1,\dots,f^\kay)\:\spc{X}\to \RR^{\kay+1}$ be a function array.
Consider the map $\spx{\bm{f}}\:\Delta^\kay\to \spc{X}$\index{$\spx{\bm{f}}$} defined by 
\[\spx{\bm{f}}(\bm{x})=\argmin\sum_{i=0}^\kay x^i\cdot f^i,\]
where $\argmin f$\index{$\argmin$} denotes a point of minimum of $f$.
The map $\spx{\bm{f}}$ will be called a \index{barycentric simplex of function array}\emph{barycentric simplex} of $\bm{f}$.
Note that for a general function array $\bm{f}$, 
the value $\spx{\bm{f}}(\bm{x})$ might be undefined or nonuniquely defined.

It is clear from the definition that $\spx{\bm{f}^{\without i}}$ 
coincides with the restriction of $\spx{\bm{f}}$ to the corresponding facet of $\Delta^\kay$.

\begin{thm}{Theorem}\label{thm:web}
Let $\spc{X}$ be a complete geodesic space 
and $\bm{f}\z=(f^0,f^1,\dots,f^\kay)\:\spc{X}\to\RR^{\kay+1}$ 
be an array of strongly convex and locally Lipschitz functions.
Then $\bm{f}$ defines a $C^{\frac12}$-embedding 
$\Web\bm{f}\hookrightarrow\RR^{\kay+1}$.

Moreover,
\begin{subthm}{thm:web:Up-convex}
$W=\Up[\bm{f}(\spc{X})]$ is a convex closed subset of $\RR^{\kay+1}$,
and

$S=\Fr_{\RR^{\kay+1}} W$ is a convex hypersurface in $\RR^{\kay+1}$.
\end{subthm}

\begin{subthm}{thm:web:f(web)=min}
\[\bm{f}(\Web\bm{f})=\Min W \subset S\]
and
\[\bm{f}(\InWeb\bm{f})= \Int_{S}(\Min W).\]
\end{subthm}

\begin{subthm}{thm:web:bary}
The barycentric simplex 
$\spx{\bm{f}}\:\Delta^\kay\to \spc{X}$ is a uniquely defined Lipshitz map and $\Im\spx{\bm{f}}=\Web\bm{f}$.
In particular, $\Web\bm{f}$ is compact.
\end{subthm}

\begin{subthm}{thm:web:lip-const}
Let us equip $\Delta^\kay$ with the metric induced by the $\ell^1$-norm on $\RR^{\kay+1}$.
Then the Lipschitz constant of $\spx{\bm{f}}\:\Delta^\kay\to\spc{U}$ can be estimated in terms of 
positive lower bounds on $(f^i)''$ 
and Lipschitz constants of $f^i$
in a neighborhood of $\Web\bm{f}$ for all $i$.
\end{subthm}

In particular, by \ref{SHORT.thm:web:Up-convex} and \ref{SHORT.thm:web:f(web)=min}, $\InWeb\bm{f}$ is $C^{\frac12}$-homeomorphic to an open set of $\RR^\kay$.
\end{thm}

The proof is preceded by a few preliminary statements.

\begin{thm}{Lemma}\label{lem:argmin(convex)}
Suppose $\spc{X}$ is a complete geodesic space and $f\:\spc{X}\to\RR$ is a locally Lipschitz, strongly convex function. Then the minimum point 
of $f$
is uniquely defined.
\end{thm}

\parit{Proof.}
Without loss of generality, we can assume that $f$ is $1$-convex.
In particular, the following claim holds:
\begin{clm}{}\label{midpoint}
 if $z$ is a midpoint of the geodesic $[x y]$, then 
\[s\le f(z)
\le
\tfrac{1}{2}\cdot f(x)+\tfrac{1}{2}\cdot f(y)-\tfrac{1}{8}\cdot\dist[2]{x}{y}{},
\]
where $s$ is the infimum of $f$.
\end{clm}

\parit{Uniqueness.}
Assume that $x$ and $y$ are distinct minimum points of $f$. 
From \ref{midpoint} we have
\[f(z)<f(x)=f(y),\] 
a contradiction. 

\parit{Existence.}
Fix a point $p\in \spc{X}$, and
let $\Lip\in\RR$ be a Lipschitz constant of $f$ in a neighborhood of $p$.

Choose a geodesic $[px]$;
consider the function $\phi\:t\mapsto f\circ\geod_{[px]}(t)$.
Clearly $\phi$ is $1$-convex and $\phi^+(0)\ge -\Lip$.
Setting $a=\dist{p}{x}{}$, we have 
\begin{align*}
f(x)
&=
\phi(a)
\ge
\\
&\ge
f(p)-\Lip\cdot a+\tfrac{1}{2}\cdot a^2
\ge
\\
&\ge f(p)-\tfrac{1}{2}\cdot{\Lip^2}.
\end{align*}

In particular,
\begin{align*}
s
&\df
\inf\set{f(x)}{x\in \spc{X}}
\ge
\\
&\ge
f(p)-\tfrac{1}{2}\cdot{\Lip^2}.
\end{align*}

Choose a sequence of points $p_n\in \spc{X}$ such that $f(p_n)\to s$.
Applying \ref{midpoint} for $x\z=p_n$, $y\z=p_m$, we see that $p_n$ is a Cauchy sequence.
Thus the sequence $p_n$ converges to a minimum point of $f$.
\qeds

\begin{thm}{Definition}
Let $Q$ be a closed subset of $\RR^{\kay+1}$.
A vector $\bm{x}\z=(x^0,x^1,\dots,x^\kay)\in\RR^{\kay+1}$
is \index{subnormal vector}\emph{subnormal} to $Q$ at a point $\bm{v}\in Q$ 
if
\[\<\bm{x},\bm{w}-\bm{v}\>
\df
\sum_ix^i\cdot(w^i-v^i)
\ge 0\]
for any $\bm{w}\in Q$.
\end{thm}

\begin{thm}{Lemma}\label{lem:Up-convex}{\sloppy 
Let $\spc{X}$ be a complete geodesic space 
and $\bm{f}\z=(f^0,f^1,\dots,f^\kay)\:\spc{X}\to\RR^{\kay+1}$ 
be an array of strongly convex and locally Lipschitz functions.
Let $W=\Up\bm{f}(\spc{X})$.
Then: 

}
\begin{subthm}{lem:Up-convex:Up-convex}
$W$ is a closed convex set, bounded below with respect to $\succcurlyeq$.
\end{subthm}

\begin{subthm}{lem:Up-convex:subnormal}
If $\bm{x}$ is a subnormal vector to $W$, then $\bm{x}\succcurlyeq\bm{0}$.
\end{subthm}

\begin{subthm}{lem:Up-convex:surface}
 $S=\Fr_{\RR^{\kay+1}}W$ is a complete convex hypersurface in $\RR^{\kay+1}$.
\end{subthm}

\end{thm}

\parit{Proof.}
Denote by $\bar W$ the closure of $W$.

Convexity of all $f^i$ implies that
for any two points $p,q\in \spc{X}$ and $t\in[0,1]$ we have
\[(1-t)\cdot\bm{f}(p)+t\cdot \bm{f}(q)
\succcurlyeq
\bm{f}\circ\geodpath_{[p q ]}(t),
\eqlbl{n-convex}\]
where $\geodpath_{[p q]}$ denotes a geodesic path from $p$ to $q$. 
Therefore $W$, as well as $\bar W$, are convex sets in $\RR^{\kay+1}$.

Let
\[w^i=\min\set{f^i(x)}{x\in\spc{X}}.\]
By Lemma~\ref{lem:argmin(convex)}, $w^i$ is finite for each $i$.
Evidently, $\bm{w}=(w^0,w^1,\dots,w^\kay)$ is a lower bound of $\bar W$ with respect to $\succcurlyeq$.

It is clear that $W$ has nonempty interior,
and $W\ne \RR^{\kay+1}$ since $W$ is bounded below.
Therefore $S=\Fr_{\RR^{\kay+1}}W=\Fr_{\RR^{\kay+1}}\bar W$
is a complete convex hypersurface in $\RR^{\kay+1}$.

Since $\bar W$ is closed and bounded below, we also have
\[\bar W=\Up[\Min\bar W].
\eqlbl{eq:W=Up Min W}\]

Choose an arbitrary $\bm{v}\in S$.
Let $\bm{x}\in\RR^{\kay+1}$ be a subnormal vector to $\bar W$ at $\bm{v}$. 
In particular, 
$\<\bm{x},\bm{y}\>
\ge
0$ 
for any $\bm{y}\succcurlyeq\bm{0}$;
that is, $\bm{x}\succcurlyeq\bm{0}$.

Further, according to Lemma~\ref{lem:argmin(convex)}, 
the function 
\[p\mapsto\langle\bm{x},\bm{f}(p)\rangle=\sum_i x^i\cdot f^i(p)\]
has a uniquely defined minimum point, say $p$.
Clearly 
\[\bm{v}\succcurlyeq\bm{f}(p)\quad\text{and}\quad \bm{f}(p)\in \Min W.\eqlbl{eq:v>f(p)}\]

Note that for any $\bm{u}\in \bar W$ there is $\bm{v}\in S$ such that $\bm{u}\succcurlyeq\bm{v}$. 
Therefore \ref{eq:v>f(p)} implies 
\[\bar W\subset\Up[\Min W]\subset W.\]
Hence
$\bar W=W$; that is, $W$ is closed.
\qeds

\parit{Proof of \ref{thm:web}; \ref{SHORT.thm:web:Up-convex}+\ref{SHORT.thm:web:f(web)=min}.}
Without loss of generality, we may assume that all $f^i$ are $1$-convex.

Given $\bm{v}=(v^0,v^1,\dots,v^\kay)\in\RR^{\kay+1}$, consider the function 
$h_{\bm{v}}\: \spc{X}\to \RR$ defined by
\[h_{\bm{v}}(p)=\max_i\{f^i(p)-v^i\}.\]
Note that $h_{\bm{v}}$ is $1$-convex.
Let 
$$\map(\bm{v})\df\argmin h_{\bm{v}}.$$
According to Lemma~\ref{lem:argmin(convex)}, $\map(\bm{v})$ is uniquely defined.

From the definition of web (\ref{def:web}) 
we have
$\map\circ\bm{f}(p)=p$ for any $p\in \Web\bm{f}$;
that is, $\map$ is a left inverse to the restriction $\bm{f}|_{\Web\bm{f}}$.
In particular, 
\[\Web\bm{f}=\Im\map.
\eqlbl{eq:Web=Im}\]

Given $\bm{v},\bm{w}\in\RR^{\kay+1}$,
set $p=\map (\bm{v})$ and $q=\map (\bm{w})$.
Since $h_{\bm{v}}$ and $h_{\bm{w}}$ are 1-convex, we have
\begin{align*}
h_{\bm{v}}(q)
&\ge 
h_{\bm{v}}(p)+\tfrac{1}{2}\cdot\dist[2]{p}{q}{},
&
h_{\bm{w}}(p)
&\ge 
h_{\bm{w}}(q)+\tfrac{1}{2}\cdot\dist[2]{p}{q}{}.
\end{align*}
Therefore,
\begin{align*}
\dist[2]{p}{q}{}
&\le 
2\cdot\sup_{x\in\spc{X}}\{ |h_{\bm{v}}(x)-h_{\bm{w}}(x)| \}
\le
\\
&\le 
2\cdot\max_{i}\{|v^i-w^i|\}.
\end{align*}
In particular,
$\map$ is $C^{\frac{1}{2}}$-continuous.
Hence $\bm{f}|_{\Web\bm{f}}$ is a $C^{\frac{1}{2}}$-embedding.

As in Lemma~\ref{lem:Up-convex},
let $W=\Up\bm{f}(\spc{X})$ and $S=\Fr_{\RR^{\kay+1}}W$.
Then
$S$ is a convex hypersurface in $\RR^{\kay+1}$.
Clearly $\bm{f}(\Web\bm{f})\z=\Min W\subset S$.
From the definition of inner web, we have
$\bm{v}\in \bm{f}(\InWeb\bm{f})$ 
if and only if 
$\bm{v}\in S$ and
for any $i$ there is $\bm{w}=(w^0,w^1,\dots,w^\kay)\in W$ such that $w^j<v^j$ for all $j\ne i$.
Thus $\bm{f}(\InWeb\bm{f})$ is open in $S$.
That is, $\InWeb\bm{f}$ is $C^{\frac{1}{2}}$-homeomorphic to an open set in a convex hypersurface $S\subset\RR^{\kay+1}$,
and hence to an open set of $\RR^{\kay}$, as claimed.

\parit{\ref{SHORT.thm:web:bary}+\ref{SHORT.thm:web:lip-const}.}
Since $f^i$ is $1$-convex, for any $\bm{x}=(x^0,x^1,\dots,x^\kay)\in\Delta^\kay$ 
the convex combination 
\[\left(\sum_i x^i\cdot f^i\right)\:\spc{X}\to\RR\] 
is also $1$-convex.
Therefore, according to Lemma~\ref{lem:argmin(convex)}, the barycentric simplex 
$\spx{\bm{f}}$ is uniquely defined on $\Delta^\kay$.
 
For $\bm{x},\bm{y}\in\Delta^\kay$,
let 
\begin{align*}
f_{\bm{x}}
&=\sum_i x^i\cdot f^i,
&
f_{\bm{y}}
&=\sum_i y^i\cdot f^i,
\\
p
&=\spx{\bm{f}}(\bm{x}),
&
q
&=\spx{\bm{f}}(\bm{y}),
\\
s&=\dist{p}{q}{}.
\end{align*}
Note the following:
\begin{itemize}
\item The function $\phi(t)=f_{\bm{x}}\circ\geod_{[p q]}(t)$ has minimum at $0$. 
Therefore $\phi^+(0)\ge 0$.

\item The function $\psi(t)=f_{\bm{y}}\circ\geod_{[p q]}(t)$ has minimum at $s$. 
Therefore $\psi^-(s)\ge 0$.
\end{itemize}
From $1$-convexity of $f_{\bm{y}}$, we have
$\psi^+(0)+\psi^-(s)+s\le0$.

Let $\Lip$ be a Lipschitz constant for all $f^i$ in a neighborhood $\Omega\ni p$.
Then 
\[\psi^+(0)
\le 
\phi^+(0)+\Lip\cdot\|\bm{x}-\bm{y}\|_1,\] 
where $\|\bm{x}-\bm{y}\|_1=\sum_{i=0}^\kay|x^i-y^i|$.
That is, given $\bm{x}\in\Delta^\kay$, there is a constant $\Lip$ such that
\begin{align*}
\dist{\spx{\bm{f}}(\bm{x})}{\spx{\bm{f}}(\bm{y})}{}
&=
s
\le
\\
&\le 
\Lip\cdot\|\bm{x}-\bm{y}\|_1
\end{align*}
for any $\bm{y}\in\Delta^\kay$.
In particular, there is $\eps>0$ such that if $\|\bm{x}-\bm{y}\|_1<\eps,$ $\|\bm{x}-\bm{z}\|_1 <\eps$, then $\spx{\bm{f}}(\bm{y})$, $\spx{\bm{f}}(\bm{z})\in\Omega$. 
Thus the same argument as above implies 
\[\dist{\spx{\bm{f}}(\bm{y})}{\spx{\bm{f}}(\bm{z})}{}
\le \Lip\cdot\|\bm{y}-\bm{z}\|_1\]
for any $\bm{y}$ and $\bm{z}$ sufficiently close to $\bm{x}$; that is, $\spx{\bm{f}}$ is locally Lipschitz.
Since $\Delta^\kay$ is compact, $\spx{\bm{f}}$ is Lipschitz.

Clearly $\spx{\bm{f}}(\Delta^\kay)\subset \Web \bm{f}$.
It remains to show that $\spx{\bm{f}}(\Delta^\kay)\supset \Web \bm{f}$.
According to Lemma~\ref{lem:Up-convex},
$W=\Up\bm{f}(\spc{X})$ is a closed convex set in $\RR^{\kay+1}$.
Let $p\in \Web \bm{f}$. 
Clearly $\bm{f}(p)\in \Min W\subset S\z=\Fr_{\RR^{\kay+1}}W$.
Let $\bm{x}$ be a subnormal vector to $W$ at $\bm{f}(p)$.
According to Lemma~\ref{lem:Up-convex}, 
$\bm{x} \succcurlyeq\bm{0}$.
Without loss of generality, we may assume that $\sum_i x^i=1$;
that is, $\bm{x}\in \Delta^\kay$.
By Lemma~\ref{lem:argmin(convex)},
$p$ is the unique minimum point of $\sum_i x^i\cdot f^i$;
that is, $p=\spx{\bm{f}}(\bm{x})$.
\qeds

\section{The case of CAT spaces}

Let $\bm{a}=(a^0,a^1,\dots,a^\kay)$ be a point array in a metric space $\spc{U}$.
Recall that 
$\distfun{\bm{a}}{}{}$
denotes the distance map
\[(\distfun{a^0}{}{},\distfun{a^1}{}{},\dots,\distfun{a^\kay}{}{})\:\spc{U}\to\RR^{\kay+1},\]
which can be also regarded as a function array.
The \index{radius of a point array}\emph{radius} of the point array $\bm{a}$ is defined to be the radius of the set $\{a^0,a^1,\dots,a^\kay\}$;
that is,
\[\rad\bm{a}=\inf\set{r>0}{\exists z\in\spc{U}\ \text{such that}\ a^i\in \oBall(z,r)\ \text{for any}\ i}.\]

Fix $\kappa\in\RR$.
Let $\bm{a}=(a^0,a^1,\dots,a^\kay)$ be a point array of radius $<\tfrac{\varpi\kappa}2$
in a metric space $\spc{U}$.
Consider the function array $\bm{f}=(f^0,f^1,\dots,f^\kay)$ 
where 
\[f^i(x)\z=\md\kappa\dist[{{}}]{a^i}{x}{}.\]
Assuming the barycentric simplex $\spx{\bm{f}}$ is defined,
then $\spx{\bm{f}}$ is called the \index{$\kappa$-barycentric simplex}\emph{$\kappa$-barycentric simplex} for the point array $\bm{a}$;
it will be denoted by $\spx{\bm{a}}\mc\kappa$.
The points $a^0,a^1,\dots,a^\kay$ are called 
\index{vertexes of the $\kappa$-barycentric simplex}\emph{vertexes} of the $\kappa$-barycentric simplex.
Note that once we say the $\kappa$-barycentric simplex is defined, 
we automatically assume that $\rad\bm{a}<\tfrac{\varpi\kappa}2$.

\begin{thm}{Theorem}\label{thm:cat-bary-web}
Let $\spc{U}$ be a complete length $\CAT\kappa$ space
and $\bm{a}\z=(a^0,a^1,\dots a^\kay)$ be a point array with radius $<\tfrac{\varpi\kappa}{2}$.
Then: 

\begin{subthm}{thm:cat-bary-web:Lip}
The $\kappa$-barycentric simplex $\spx{\bm{a}}\mc\kappa\:\Delta^\kay\to \spc{U}$ 
is defined. 
Moreover, $\spx{\bm{a}}\mc\kappa$ is a Lipschitz map,
and if $\Delta^\kay$ is equipped with the $\ell^1$-metric, then its Lipschitz constant can be estimated in terms of $\kappa$ and the radius of $\bm{a}$ (in particular it does not depend on $\kay$).
\end{subthm}

\begin{subthm}{thm:cat-bary-web:web=Im(bary)}
$\Web(\distfun{\bm{a}}{}{})=\Im \spx{\bm{a}}\mc\kappa$.
Moreover, if a closed convex set $K\subset\spc{U}$ contains all $a^i$, then $\Web(\distfun{\bm{a}}{}{})\subset K$.
\end{subthm}

\begin{subthm}{thm:cat-bary-web:mnfld}
The restriction%
\footnote{Recall that $\distfun{\bm{a}^{\without 0}}{}{}$ denotes the array $(\distfun{a^1}{}{},\dots,\distfun{a^\kay}{}{})$.}
$\distfun{\bm{a}^{\without 0}}{}{}|_{\InWeb(\distfun{\bm{a}}{}{})}$ is an open $C^{\frac12}$-embedding in $\RR^\kay$.
Thus there is an inverse of 
$\distfun{\bm{a}^{\without 0}}{}{}|_{\InWeb(\distfun{\bm{a}}{}{})}$, say $\map\:\RR^\kay\subto\spc{U}$.

The subfunction $f=\distfun{a^0}{}{}\circ\map$ is semiconvex and locally Lipschitz.
Moreover, if $\kappa\le 0$, then $f$ is convex.
\end{subthm}

In particular, $\Web(\distfun{\bm{a}}{}{})$ is a compact set and
$\InWeb(\distfun{\bm{a}}{}{})$ is $C^{\frac12}$-homeomorphic to an open subset of $\RR^\kay$.

\end{thm}

\begin{thm}{Definition}\label{prop-def:web-embedding}
The submap $\map\:\RR^{\kay}\subto \spc{X}$ of Theorem~\ref{thm:cat-bary-web:mnfld}
will be called the \index{web embedding}\emph{$\distfun{\bm{a}}{}{}$-web embedding} 
with \index{brace}\emph{brace} $\distfun{a^0}{}{}$.
The terminology invokes Theorem~\ref{thm:cat-bary-web:mnfld}.
\end{thm}

\begin{thm}{Definition}
Let $\spc{U}$ be a complete length $\CAT\kappa$ space
and $\bm{a}\z=(a^0,a^1,\dots a^\kay)$ be a point array with radius $<\tfrac{\varpi\kappa}{2}$.
If $\InWeb(\distfun{\bm{a}}{}{})$ is nonempty, then the point array $\bm{a}$ is called \emph{nondegenerate}.
\end{thm}

Lemma~\ref{lem:nondeg-test-with-balls} will provide examples of nondegenerate point arrays,
which can be used in \ref{thm:cat-bary-web:mnfld}.

\begin{thm}{Corollary}\label{cor:LinDim>bary}
Let $\spc{U}$ be a complete length $\CAT\kappa$ space,
$\bm{a}\z=(a^0,a^1,\dots a^m)$ be a nondegenerate point array 
of radius $<\tfrac{\varpi\kappa}{2}$ in $\spc{U}$
and $\sigma=\spx{\bm{a}}\mc\kappa$ be the corresponding $\kappa$-baricentric simplex.
Then for some $\bm{x}\in \Delta^m$,
the differential $\dd_{\bm{x}}\sigma$ is linear 
and the image $\Im\dd_{\bm{x}}\sigma$
forms a subcone isometric to an $m$-dimensional Euclidean space in the tangent cone $\T_{\sigma(\bm{x})}$.
\end{thm}

\parit{Proof.}
Denote the distance map $\distfun{\bm{a}^{\without 0}}{}{}$ by $\tau\:\spc{U}\to\RR^m$.

According to Theorem~\ref{thm:cat-bary-web},
$\sigma$ is Lipschitz
and the distance map $\tau$ 
gives an open embedding of 
$\InWeb(\distfun{\bm{a}}{}{})=\sigma(\Delta^m)\setminus\sigma(\partial\Delta^m)$.
Note that $\tau$ is Lipschitz.
According to Rademacher's theorem (\ref{thm:Rademacher-CBB+CBA}), 
the differential 
$\dd_{\bm{x}}(\tau\circ\sigma)$
is linear for almost all $\bm{x}\in\Delta^m$.
Further, since $\InWeb(\distfun{\bm{a}}{}{})\z\ne\emptyset$,
the area formula \cite{karmanova} implies that $\dd_{\bm{x}}(\tau\circ\sigma)$ is surjective on a set of positive masure of points $\bm{x}\in\Delta^m$.

Note that $\dd_{\bm{x}}(\tau\circ\sigma)=(\dd_{\sigma(\bm{x})}\tau)\circ(\dd_{\bm{x}}\sigma)$.
Applying Rademacher's theorem again, we have linearity of 
$\dd_{\bm{x}}\sigma$ for almost all $\bm{x}\in\Delta^m$;
at these points $\Im\dd_{\bm{x}}\sigma$ forms a subcone isometric to a Euclidean space in $\T_{\sigma(\bm{x})}$.
Clearly the dimension of $\Im\dd_{\bm{x}}(\tau\circ\sigma)$ is at least as big as the dimension of $\Im\dd_{\bm{x}}\sigma$.
Hence the result.
\qeds

\parit{Proof of \ref{thm:cat-bary-web}.}
Fix $z\in\spc{U}$ and $r<\tfrac{\varpi\kappa}2$
such that $\dist{z}{a^i}{}<r$ for all $i$.
Note that the set $K\cap \cBall[z,r]$ is convex, closed, and contains all $a^i$.
Applying the theorem on short retract (Exercise~\ref{ex:short-retraction-CBA(1)}),
we get the second part of~\ref{SHORT.thm:cat-bary-web:web=Im(bary)}.

\medskip

The remaining statements are proved first in the case $\kappa\le 0$, 
and then the remaining case $\kappa>0$ is reduced to the case $\kappa=0$.

\parit{Case $\kappa\le 0$.}
Consider the function array $f^i=\md\kappa\circ \distfun{a^i}{}{}$.
From the definition of web (\ref{def:web}),
it is clear that $\Web(\distfun{\bm{a}}{}{})=\Web\bm{f}$.
Further, by the definition of $\kappa$-barycentric simplex,
$\spx{\bm{a}}\mc\kappa=\spx{\bm{f}}$.

All the functions $f^i$ are strongly convex (see \ref{function-comp}).
Therefore \ref{SHORT.thm:cat-bary-web:Lip}, \ref{SHORT.thm:cat-bary-web:web=Im(bary)}, and the first statements in \ref{SHORT.thm:cat-bary-web:mnfld} follow from Theorem \ref{thm:web}.

\parit{Case $\kappa>0$.}
Applying rescaling, we may assume $\kappa=1$,
so $\varpi\kappa=\varpi1\z=\pi$.

Let $\mathring{\spc{U}}=\Cone\spc{U}$.
By \ref{thm:warp-curv-bound:cbb:S}, $\mathring{\spc{U}}$ is $\CAT0$.
Let us denote by $\iota$ the natural embedding of $\spc{U}$ as the unit sphere in $\mathring{\spc{U}}$, and by 
$\proj\:\mathring{\spc{U}}\subto\spc{U}$ the submap
defined by $\proj(v)=\iota^{-1}(v/|v|)$ for all 
$v\ne \0$.
Note that there is $z\in\spc{U}$ and $\eps>0$ such that
the set 
\[K_\eps
=
\set{v\in\mathring{\spc{U}}}%
{\<\iota(z),v\>\ge\eps}\] 
contains all $\iota(a^i)$.
Then 
$\0\notin K_\eps$, 
and
the set $K_\eps$ is closed and convex.
The latter follows from Exercise~\ref{ex:busemann-CBA},
since $v\mapsto -\<\iota(z),v\>$ is a Busemann function.

Denote by $\iota(\bm{a})$ the point array $(\iota(a^0), \iota(a^1),\dots,\iota(a^\kay))$ in $\mathring{\spc{U}}$. 
From the case $\kappa=0$,
we get that $\Im \spx{\iota(\bm{a})}\mc0\subset K_\eps$.
In particular, $\Im \spx{\iota(\bm{a})}\mc0\not\ni \0$ and thus $\proj\circ\spx{\iota(\bm{a})}\mc0$ is defined.
Direct calculations show 
\[\spx{\bm{a}}\mc1
=
\proj\circ\spx{\iota(\bm{a})}\mc0
\quad\text{and}\quad
\Web(\distfun{\bm{a}}{}{})=\proj[\Web(\distfun{\iota(\bm{a})}{}{})].\]
Thus the case $\kappa=1$ of the theorem is reduced to the case $\kappa=0$,
which is proved already.
\qeds

\begin{thm}{Lemma}\label{lem:nondeg-test-with-balls}
Let $\bm{a}\z=(a^0,a^1,\dots a^\kay)$ be a point array of radius $<\tfrac{\varpi\kappa}2$
in a complete length $\CAT\kappa$ space $\spc{U}$, 
and $B^i=\cBall[a^i,r^i]$ for an array of positive reals $(r^0,r^1,\dots,r^\kay)$.
Assume that
$\bigcap_i B^i=\emptyset$,
but
$\bigcap_{i\ne j} B^i\ne \emptyset$
for any $j$.
Then $\bm{a}$ is nondegenerate. 
\end{thm}

\parit{Proof.} 
Without loss of generality, we may assume that $\spc{U}$ is geodesic and $\diam\spc{U}\z<\varpi\kappa$.
If not, choose $z\in\spc{U}$ and $r<\tfrac{\varpi\kappa}{2}$ so that
$\dist{z}{a^i}{}\le r$
for each $i$, 
and consider $\cBall[z,r]$ instead of $\spc{U}$.
The latter can be done since $\cBall[z,r]$ is convex and closed, 
so $\cBall[z,r]$ is a complete length $\CAT\kappa$ space 
and $\Web(\distfun{\bm{a}}{}{})\subset\cBall[z,r]$;
see \ref{cor:convex-balls} and \ref{thm:cat-bary-web:web=Im(bary)}.

By Theorem~\ref{thm:cat-bary-web}, $\Web(\distfun{\bm{a}}{}{})$ is a compact set;
therefore there is a point $p\in\Web(\distfun{\bm{a}}{}{})$
minimizing the function 
\[f(x)=\max_i\{\distfun{B^i}{x}{}\}=\max\{0,\dist{a^0}{x}{}-r^0,\dots,\dist{a^\kay}{x}{}-r^\kay\}.\]

By the definition of web (\ref{def:web}), 
$p$ is also the minimum point of $f$ on~$\spc{U}$.
Let us prove the following claim:

\begin{clm}{}
 $p\notin B^j$ for any $j$.
\end{clm}

Indeed, 
assume the contrary; that is, 
\[
p\in B^j
\eqlbl{eq:p-in-Bj}
\] 
for some $j$.
Then $p$ is a point of local minimum for the function 
\[h^j(x)=\max_{i\ne j}\{\distfun{B^i}{x}{}\}.\]
Hence 
\[\max_{i\ne j}\{\mangle\hinge p x {a^i}\}\ge \tfrac\pi2
\]
for any $x\in\spc{U}$.
From the angle comparison (\ref{cat-hinge}), it follows that 
$p$ is a global minimum of $h^j$ and hence
\[
p\in \bigcap_{i\ne j} B^i.
\]
The latter and \ref{eq:p-in-Bj} contradict $\bigcap_i B^i=\emptyset$. \claimqeds 
\noindent 

From the definition of web, it also follows that 
\[\Web(\distfun{\bm{a}^{\without j}}{}{})\subset \bigcup_{i\ne j}B^i.\]
Indeed, if $q\in \bigcap_{i\ne j}B^i$ and $q'\notin \bigcup_{i\ne j}B^i$,
then $\dist{a_i}{q}{}\z<\dist{a_i}{q'}{}$ for any $i\ne j$ therefore $q'\notin \Web(\distfun{\bm{a}^{\without j}}{}{})$.
Therefore the claim implies that
$p\notin\Web(\distfun{\bm{a}^{\without j}}{}{})$ for each $j$;
that is, $p\in\InWeb(\distfun{\bm{a}}{}{})$.
\qeds

\section{Dimension}\label{sec:dim-cba}

See Chapter \ref{ch:dim} for definitions of various dimension-like invariants of metric spaces.

We start with two examples.

The first example shows that the dimension of complete length $\CAT{}$ spaces is not local;
that is, such spaces might have open sets with different linear dimensions.

Such an example can be constructed by gluing at one point two Euclidean spaces of different dimensions.
According to Reshetnyak's gluing theorem (\ref{thm:gluing}), this construction gives a $\CAT{0}$ space.

The second example provides a complete length $\CAT{}$ space 
with topological dimension 1 and arbitrary large Hausdorff dimension.
Thus for complete length $\CAT{}$ spaces, one should not expect any relations between topological and Hausdorff dimensions except for the one provided by Szpilrajn's theorem (\ref{thm:szpilrajn}).

To construct the second type of example,
note that the completion of any metric tree has topological dimension 1 and is $\CAT\kappa$ for any $\kappa$.
Start with a binary tree $\Gamma$, and a sequence $\eps_n>0$ such that $\sum_n\eps_n<\infty$.
Define the metric on $\Gamma$
by prescribing the length of an edge from level $n$ to level $n+1$ to be $\eps_n$.
For an appropriately chosen sequence $\eps_n$, the completion of $\Gamma$ will contain a Cantor set of arbitrarily large Hausdorff dimension.

\medskip

The following is a version of a theorem proved by Bruce Kleiner \cite{kleiner}, with an improvement made by Alexander Lytchak \cite{lytchak:diff}.

\begin{thm}{Theorem}\label{thm:dim-infty-CBA}
For any complete length $\CAT\kappa$ space $\spc{U}$, the following statements are equivalent:

\begin{subthm}{LinDim-CBA} $\LinDim\spc{U}\ge m$.
\end{subthm}

\begin{subthm}{thm:dim-infty-CBA:bary} 
For some $z\in \spc{U}$ there is an array of $m+1$ balls $B^i\z=\oBall(a^i,r^i)$ with $a^0,a^1,\dots,a^m\in \oBall(z,\frac{\varpi\kappa}2)$ 
such that 
\[\bigcap_i B^i=\emptyset
\quad\text{and}\quad
\bigcap_{i\ne j} B^i\ne \emptyset
\quad \text{for each $j$}.\]

\end{subthm}

\begin{subthm}{thm:dim-infty-CBA:mnfld} 
There is a $C^{\frac{1}{2}}$-embedding $\map\:\cBall[1]_{\EE^m}\hookrightarrow \spc{U}$;
that is, $\map$ is bi-Hölder with exponent $\tfrac{1}{2}$.
\end{subthm}

\begin{subthm}{thm:dim-infty-CBA:TopDim}
There is a closed separable set $K\subset\spc{U}$ such that 
\[\TopDim K\ge m.\]
\end{subthm}

\end{thm}

\parbf{Remarks.}
Theorem \ref{thm:loc-lip-inverse} gives a stronger version of part \ref{SHORT.thm:dim-infty-CBA:mnfld} in the finite-dimensional case.
Namely, a complete length $\CAT{}$ space with linear dimension $m$ 
admits a bi-Lipschitz embedding $\map$ of an open set of $\RR^m$.
Moreover, the Lipschitz constants of $\map$ can be made arbitrarily close to~$1$.

\begin{thm}{Corollary}\label{cor:dim-CBA}
For any separable complete length $\CAT{}$ space $\spc{U}$, we have
\[\TopDim\spc{U}=\LinDim\spc{U}.\]

\end{thm}

Any simplicial complex can be equipped with a length metric
such that each $\kay$-simplex 
is isometric to the standard simplex
\[\Delta^\kay
=
\set{(x_0,\dots,x_\kay)\in \RR^{\kay+1}}
{x_i\ge 0,\quad x_0+\dots+x_\kay=1}\]
with the metric induced by the $\ell^1$-norm on $\RR^{\kay+1}$.
This metric will be called the \index{$\ell^1$-metric}\emph{$\ell^1$-metric} on the simplicial complex.

\begin{thm}{Lemma}\label{lem:approximation-cba}
Let $\spc{U}$ be a complete length $\CAT\kappa$ space
and $\rho\:\spc{U}\to\RR$ be a continuous positive function.
Then there is a simplicial complex $\spc{N}$ equipped with $\ell^1$-metric,
a locally Lipschitz map $\map\:\spc{U}\to \spc{N}$, 
and a Lipschitz map $\map[2]\:\spc{N}\to\spc{U}$ such that:

\begin{subthm}{lem:approximation-cba:displacement}
The displacement of the composition $\map[2]\circ\map\:\spc{U}\to\spc{U}$ is bounded by $\rho$;
that is,
\[\dist{x}{\map[2]\circ\map(x)}{}<\rho(x)\] 
for any $x\in\spc{U}$.
\end{subthm}

\begin{subthm}{lem:approximation-cba:im}
If $\LinDim\spc{U}\le m$ 
then the $\map[2]$-image of any closed simplex in $\spc{N}$ 
coincides with the image of its $m$-skeleton.
\end{subthm}

\end{thm}

\parit{Proof.}
Without loss of generality, we may assume that for any $x$ we have $\rho(x)\z<\rho_0$
for fixed $\rho_0<\tfrac{\varpi\kappa}{2}$.

By Stone's theorem, any metric space is paracompact.
Thus, we can choose a locally finite covering $\set{\Omega_\alpha}{\alpha\in\IndexSet}$ of $\spc{U}$ such that $\Omega_\alpha\subset \oBall(x,\tfrac{1}{3}\cdot\rho(x))$ for any $x\in \Omega_\alpha$. 

Denote by $\spc{N}$ the nerve of the covering $\{\Omega_\alpha\}$;
that is, $\spc{N}$ is an abstract simplicial complex with 
vertex set $\IndexSet$,
such that
$\{\alpha^0,\alpha^1,\dots,\alpha^n\}\subset\IndexSet$ 
are vertexes of a simplex if and only if
$\Omega_{\alpha^0}
\cap
\Omega_{\alpha^1}
\cap\dots\cap
\Omega_{\alpha^n}
\ne 
\emptyset$.

Fix a Lipschitz partition of unity 
$\phi_\alpha\:\spc{U}\to [0,1]$ subordinate to $\{\Omega_\alpha\}$.
Consider the map $\map\:\spc{U}\to \spc{N}$ such that the barycentric coordinate of $\map(p)$ is $\phi_\alpha(p)$.
Note that $\map$ is locally Lipschitz. 
Clearly the $\map$-preimage of any open simplex in $\spc{N}$ lies in $\Omega_\alpha$ for some $\alpha\in\IndexSet$.

For each $\alpha\in\IndexSet$, 
choose $x_\alpha\in\Omega_\alpha$.
Let us extend the map $\alpha\mapsto x_\alpha$
to a map $\map[2]\:\spc{N}\to\spc{U}$ that is $\kappa$-barycentric on each simplex.
According to Theorem~\ref{thm:cat-bary-web:Lip}, this extension exists, 
$\map[2]$ is Lipschitz, 
and its Lipschitz constant depends only on $\rho_0$ and $\kappa$.

\parit{\ref{SHORT.lem:approximation-cba:displacement}.}
Fix $x\in\spc{U}$. Denote by $\Delta$ the minimal simplex that contains $\map(x)$, 
and let $\alpha^0,\alpha^1,\dots,\alpha^n$ be the vertexes of $\Delta$.
Note that $\alpha$ is a vertex of $\Delta$ if and only if $\phi_{\alpha}(x)>0$.
Thus
\[\dist{x}{x_{\alpha^i}}{}<\tfrac{1}{3}\cdot\rho(x)\] 
for any $i$.
Therefore 
\[\diam\map[2](\Delta)
\le
\max_{i,j}\{\dist{x_{\alpha^i}}{x_{\alpha^j}}{}\}
<
\tfrac{2}{3}\cdot\rho(x).\]
In particular, 
\[\dist{x}{\map[2]\circ\map(x)}{}\le\dist{x}{x_{\alpha^0}}{}+\diam \map[2](\Delta) <\rho(x).\]

\parit{\ref{SHORT.lem:approximation-cba:im}.}
Assume the contrary;
that is, $\map[2](\spc{N})$ is not included in the $\map[2]$-image of the $m$-skeleton of $\spc{N}$.
Then for some $\kay>m$,
there is a $\kay$-simplex $\Delta^\kay$ in $\spc{N}$
such that the barycentric simplex $\sigma=\map[2]|_{\Delta^\kay}$ is nondegenerate; 
that is, 
$$W=\map[2](\Delta^\kay)\setminus\map[2](\partial\Delta^\kay)\ne \emptyset.
$$
Applying Corollary~\ref{cor:LinDim>bary}
gives $\LinDim\spc{U}\ge \kay$, a contradiction.
\qeds

\parit{Proof of \ref{thm:dim-infty-CBA}; \ref{SHORT.thm:dim-infty-CBA:bary}$\Rightarrow$\ref{SHORT.thm:dim-infty-CBA:mnfld}$\Rightarrow$\ref{SHORT.thm:dim-infty-CBA:TopDim}.}
The implication \ref{SHORT.thm:dim-infty-CBA:bary}$\Rightarrow$\ref{SHORT.thm:dim-infty-CBA:mnfld} follows from Lemma~\ref{lem:nondeg-test-with-balls}
and Theorem~\ref{thm:cat-bary-web:mnfld}, and \ref{SHORT.thm:dim-infty-CBA:mnfld}$\Rightarrow$\ref{SHORT.thm:dim-infty-CBA:TopDim} is trivial.
 
\parit{\ref{SHORT.thm:dim-infty-CBA:TopDim}$\Rightarrow$\ref{SHORT.LinDim-CBA}.}
According to Theorem~\ref{thm:stable-value}, 
there is a continuous map $f\:K\to \RR^{m}$ with a stable value.
By the Tietze extension theorem, it is possible to extend $f$ 
to a continuous map $F\:\spc{U}\to \RR^{m}$.

Fix $\eps>0$.
Since $F$ is continuous, there is a continuous positive function $\rho$ defined on $\spc{U}$ such that 
\[\dist{x}{y}{}<\rho(x)
\quad\Rightarrow\quad
|F(x)- F(y)|<\tfrac13\cdot\eps.\]
Apply Lemma~\ref{lem:approximation-cba} to $\rho$.
For the resulting simplicial complex $\spc{N}$ 
 and the maps $\map\:\spc{U}\to \spc{N}$, $\map[2]\:\spc{N}\to \spc{U}$, we have
\[|F\circ \map[2]\circ\map(x)-F(x)|<\tfrac13\cdot\eps\] 
for any $x\in \spc{U}$.

According to Lemma~\ref{lem:lip-approx},
there is a locally Lipschitz map $F_\eps\:\spc{U}\to \RR^{m+1}$ 
such that $|F_\eps(x)-F(x)|<\tfrac13\cdot\eps$ for any $x\in \spc{U}$.

Note that
$\map(K)$ is contained in a countable subcomplex of $\spc{N}$, say $\spc{N}'$.
Indeed, since $K$ is separable, there is a countable dense collection of points $\{x_n\}$ in $K$.
Denote by $\Delta_n$ the minimal simplex of $\spc{N}$ that contains $\map(x_n)$.
Then $\map(K)\subset\bigcup_i\Delta_n$.

Arguing by contradiction,
assume $\LinDim\spc{U}<m$.
By \ref{lem:approximation-cba:im},
the image $F_\eps\circ\map[2]\circ\map(K)$ lies in the $F_\eps$-image of the $(m-1)$-skeleton of $\spc{N}'$;
In particular, it can be covered by a countable collection of Lipschitz images of $(m-1)$-simplexes.
Hence
$\bm{0}\in \RR^m$ is not a stable value of the restriction $F_\eps\circ\map[2]\circ\map|_K$.
Since $\eps>0$ is arbitrary, 
then $\bm{0}\in \RR^m$ is not a stable value of $f$ --- a contradiction.

\parit{\ref{SHORT.LinDim-CBA}$\Rightarrow$\ref{SHORT.thm:dim-infty-CBA:bary}.} 
Choose $q\in \spc{U}$ such that $\T_q$ contains a subcone $E$ isometric to $m$-dimensional Euclidean space.
Note that one can choose $\eps>0$ 
and a point arrray $(\dot a^0,\dot a^1,\dots,\dot a^m)$ in $E\subset \T_q$ 
such that 
$\bigcap_i\cBall[\dot a^i,1+\eps]=\emptyset$
and $\bigcap_{i\ne j}\cBall[\dot a^i,1\z-\eps]\z{\ne }\emptyset$ for each $j$.

For each $i$ choose a geodesic $\gamma^i$ 
from $q$ that goes almost in the directions of $\dot a^i$.
Choose small $\delta>0$ and take the point $a^i$ on $\gamma^i$ at distance $\delta\cdot|\dot a^i|$ from $q$.
We get a point array 
$(a^0,a^1,\dots,a^m)$ in $\spc{U}$
such that $\bigcap_i\cBall[a^i,\delta]\z=\emptyset$
and $\bigcap_{i\ne j}\cBall[a^i,\delta]\ne \emptyset$ for each $j$.
Since $\delta>0$ can be chosen arbitrarily small, 
 \ref{SHORT.thm:dim-infty-CBA:bary} follows.
\qeds

\section{Finite-dimensional spaces}

Recall that a web embedding and its brace are defined in \ref{prop-def:web-embedding}.

{\sloppy 

\begin{thm}{Theorem}\label{thm:loc-lip-inverse}
Suppose $\spc{U}$ is a complete length $\CAT\kappa$ space such that 
$\LinDim\spc{U}=m$,
and $\bm{a}=(a^0,a^1,\dots a^m)$ is a point array in $\spc{U}$ 
with radius $<\tfrac{\varpi\kappa}{2}$.
Then 
the $\distfun{\bm{a}}{}{}$-web embedding $\map\:\RR^m\subto\spc{U}$ with brace $\distfun{a^0}{}{}$ is locally Lipschitz.
\end{thm}

}

Note that if $\bm{a}$ is degenerate,
that is, if $\InWeb(\distfun{\bm{a}}{}{})=\emptyset$, 
then
the domain of $\map$ is empty, and hence the conclusion of the theorem trivially holds.

\begin{thm}{Lemma}\label{lem:nondeg-bs-test}
Let $\spc{U}$ be a complete length $\CAT\kappa$ space,
and $\bm{a}\z=(a^0,a^1,\dots a^\kay)$ be a point array with radius $<\tfrac{\varpi\kappa}{2}$.
Then for any $p\in \InWeb(\distfun{\bm{a}}{}{})$,
there is $\eps>0$ such that 
if for some $q\in \Web(\distfun{\bm{a}}{}{})$ and $b\in\spc{U}$
we have 
\[\dist{p}{q}{}<\eps,
\quad 
\dist{p}{b}{}<\eps,
\quad\text{and}\quad 
\mangle\hinge{q}{b}{a^i}<\tfrac\pi2+\eps
\eqlbl{eq:2nd-angle}\]
for each $i$,
then the array $(b,a^0,a^1,\dots,a^m)$ is nondegenerate.
\end{thm}

\parit{Proof.}
Without loss of generality, we may assume that $\spc{U}$ is geodesic and $\diam\spc{U}<\varpi\kappa$.
If not, consider instead of $\spc{U}$,
a ball $\cBall[z,r]\subset\spc{U}$ 
for some $z\in\spc{U}$ 
and $r<\tfrac{\varpi\kappa}{2}$
such that $\dist{z}{a^i}{}\le r$ for each $i$.

From the angle comparison (\ref{cat-hinge}), it follows that 
$p\in\InWeb\bm{a}$ if and only if both of the following conditions hold:
\begin{enumerate}
\item $\max_i\{\mangle\hinge p{a^i}{u}\}\ge \tfrac\pi2$ for any $u\in\spc{U}$,
\item\label{prop:<pi/2} for each $i$ there is $u^i\in\spc{U}$ such that $\mangle\hinge p{a^j}{u^i}<\tfrac\pi2$ for all $j\ne i$.
\end{enumerate}

Due to the semicontinuity of angles (\ref{lem:ang.semicont}),
there is $\eps>0$ such that for any $x\in \oBall(p,10\cdot\eps)$ we have
\[
\mangle\hinge {x}{a^j}{u^i}
<
\tfrac\pi2-10\cdot\eps
\quad\text{for all}\quad j\ne i.
\eqlbl{eq:<pi/2-eps}\]

Now assume that for sufficiently small $\eps>0$
there are points $b\in\spc{U}$ and $q\in\Web(\distfun{\bm{a}}{}{})$ such that \ref{eq:2nd-angle} holds.
According to Theorem~\ref{thm:cat-bary-web:web=Im(bary)},
for all small $\eps>0$ we have 
\[\rad\{b,a^0,a^1,\dots,a^\kay\}<\tfrac{\varpi\kappa}2.\]

\begin{wrapfigure}{r}{33 mm}
\vskip-0mm
\centering
\includegraphics{mppics/pic-1310}
\vskip0mm
\end{wrapfigure}

Fix a sufficiently small $\delta>0$
and let 
\[v^i=\geod_{[q u^i]}(\tfrac13\cdot\delta)\quad\text{and}\quad w^i=\geod_{[v^i b]}(\tfrac23\cdot\delta).\]
Clearly
\begin{align*}
\dist{b}{w^i}{}
&=
\dist{b}{v^i}{}-\tfrac23\cdot\delta
\le
\\
&\le
\dist{b}{q}{}-\tfrac13\cdot\delta.
\\
\intertext{Further, the inequalities \ref{eq:<pi/2-eps} and \ref{eq:2nd-angle} imply}
\dist{a^j}{w^i}{}
&<
\dist{a^j}{v^i}{}+\tfrac23\cdot\eps\cdot\delta
<
\\
&<
\dist{a^i}{q}{}-\eps\cdot\delta
<
\\
&<
\dist{a^i}{q}{}
\end{align*}
for all $i\ne j$.

Set $B^i=\cBall[a^i,\dist{a^i}{q}{}]$ and $B^{m+1}=\cBall[b,\dist{a^i}{q}{}-\tfrac13\cdot\delta]$.
Clearly 
\begin{align*}
&\!\!\!\!\bigcap_{i\ne m+1} B^i=\{q\},
\\
&\bigcap_{i\ne j}B^i\ni w^j\quad\text{for}\quad j\ne m+1,
\\
&\bigcap_{i}B^i=\{q\}\cap B^{m+1}=\emptyset.
\end{align*}
Lemma~\ref{lem:nondeg-test-with-balls} finishes the proof.
\qeds

\parit{Proof of \ref{thm:loc-lip-inverse}.}
Suppose $\map$ is not locally Lipshitz.; that is, there are sequences $\bm{y}_n, \bm{z}_n\to \bm{x}\in\Dom\map$ such that
\[\frac{\dist{\map(\bm{y}_n)}{\map(\bm{z}_n}{})}{|\bm{y}_n-\bm{z}_n|}
\to\infty
\quad
\text{as}
\quad
n\to\infty.
\eqlbl{eq:nonlip}\]
Set $p=\map(\bm{x})$,
$q_n=\map(\bm{y}_n)$, 
and $b_n=\map(\bm{z}_n)$.
By \ref{prop-def:web-embedding}, $p$, $q_n$, $b_n\in\InWeb(\distfun{\bm{a}}{}{})$
and $q_n,b_n\to p$ as $n\to\infty$.
Choose $\eps>0$; note that \ref{eq:nonlip} implies
\[\mangle\hinge{q_n}{a^i}{b_n}<\tfrac\pi2+\eps
\]
for all $i>0$ and all large $n$.
Further, according to \ref{prop-def:web-embedding}, the subfunction
$(\distfun{\bm{a}^0}{}{})\circ\map$ is locally Lipschitz.
Therefore we also have 
\[\mangle\hinge{q_n}{a^0}{b_n}<\tfrac\pi2+\eps
\]
for all large $n$.
According to Lemma~\ref{lem:nondeg-bs-test}, the point array $b_n,a^0,\dots,a^\kay$ for large $n$ is nondegenerate.

Applying Corollary~\ref{cor:LinDim>bary},
we have a contradiction.
\qeds

\section{Remarks}

The following conjecture was formulated by Bruce Kleiner \cite{kleiner}, see also \cite[p.~133]{gromov:asymt-inv}.
For separable spaces, it follows from Corollary~\ref{cor:dim-CBA}.

\begin{thm}{Conjecture}
For any complete length $\CAT{}$ space $\spc{U}$, we have
\[\TopDim\spc{U}
=
\LinDim\spc{U}.\]

\end{thm}

\chapter{Dimension of CBB spaces}

As the main dimension-like invariant, we will use the linear dimension $\LinDim$; 
see Definition~\ref{def:lin-dim}. In other words, by default {}\emph{dimension} means {}\emph{linear dimension}. 

\section{Struts and rank}\label{sec:struts+rank}

Our definitions of strut 
and distance chart 
differ from the one in \cite{burago-gromov-perelman};
it is closer to Perelman's definitions \cite{perelman:spaces2,perelman:morse}.

The term ``strut'' seems to have the closest meaning to the original Russian term used by Yuriy Burago, Grigory Perelman, and Michael Gromov \cite{burago-gromov-perelman}.
In the official translation,
it appears as ``burst'', 
and in the authors' translation it was ``strainer''.
Neither seems intuitive, 
so we decided to switch to ``strut''.

\begin{thm}{Definition of struts}\label{def:strut-I}
Let $\spc{L}$ be a complete length $\Alex{}$ space.
We say that a point array $(a^0,a^1,\dots,a^\kay)$ in $\spc{L}$
 is \index{strutting point array}\emph{$\kappa$-strutting} for a point $p\in\spc{L}$ if $\angk\kappa p {a^i}{a^j}>\tfrac\pi2$ for all $i\ne j$.
\end{thm} 

Recall that the packing number is defined in \ref{sec:notations}.
The following definition is motivated by the observation that $k\z=\pack_{\pi/2}(\SS^{k-1})-1$ for any integer $k>0$. 

\begin{thm}{Definition}\label{def:rank}
Let $\spc{L}$ be a complete length $\Alex{}$ space
and $p\in \spc{L}$.
Let us define rank of $\spc{L}$ at $p$ as 
\[\rank_p=\rank_p\spc{L}
\df
\pack_{\pi/2}\Sigma_p-1.\]

\end{thm}

Thus $\rank$ takes values in $\ZZ_{\ge0}\cup\{\infty\}$.

\begin{thm}{Proposition}\label{prop:stutt}
Let $\spc{L}$ be a complete length $\Alex{\kappa}$ space 
and $p\in \spc{L}$.
Then the following conditions are equivalent:

\begin{subthm}{lem:stut<=>rank:rank}$\rank_p\ge\kay$,
\end{subthm}

\begin{subthm}{lem:stut<=>rank:strut}
there is a point array $(a^0,a^1,\dots,a^\kay)$
that is $\kappa$-strutting at $p$. 
\end{subthm}
\end{thm}

\parit{Proof of \ref{prop:stutt}, \ref{SHORT.lem:stut<=>rank:strut}$\Rightarrow$\ref{SHORT.lem:stut<=>rank:rank}.}
For each $i$,
choose a point $\acute a^i\in\Str(p)$ sufficiently close to $a^i$ (so $[p\acute a^i]$ exists for each $i$).
One can choose $\acute a^i$ so that we still have
$\angk\kappa{p}{\acute a^i}{\acute a^j}>\tfrac\pi2$ for all $i\ne j$.

From hinge comparison (\ref{angle}),
\[\mangle(\dir{p}{\acute a^j},\dir{p}{\acute a^j})
\ge
\angk\kappa{p}{\acute a^i}{\acute a^j}
>
\tfrac\pi2\]
for all $i\ne j$.
In particular, $\pack_{\pi/2}\Sigma_p\ge \kay+1$.

\parit{\ref{SHORT.lem:stut<=>rank:rank}$\Rightarrow$\ref{SHORT.lem:stut<=>rank:strut}.} 
Assume $(\xi^0,\xi^1,\dots,\xi^\kay)$ is an array of directions in $\Sigma_p$, such that $\mangle(\xi^i,\xi^j)>\tfrac\pi2$ if $i\ne j$.

Without loss of generality, 
we may assume that each direction $\xi^i$ is geodesic;
that is, for each $i$ there is a geodesic $\gamma^i$ in $\spc{L}$ such that $\gamma^i(0)\z=p$ and $\xi^i=(\gamma^i)^+(0)$.
From the definition of angle, it follows that for sufficiently small $\eps>0$ the array of points $a^i=\gamma^i(\eps)$ satisfies \ref{SHORT.lem:stut<=>rank:strut}.
\qeds

\begin{thm}{Corollary}\label{cor:rank>=k-open}
Let $\spc{L}$ be a complete length $\Alex{}$ space and $\kay\in\ZZ_{\ge0}$.
Then the set of all points in $\spc{L}$ 
with rank $\ge \kay$ is open.
\end{thm}

\parit{Proof.} Given an array of points $\bm{a}=(a^0,\dots,a^\kay)$ in $\spc{L}$, consider 
the set $\Omega_{\bm{a}}$ of all points $p\in \spc{L}$ such that array $\bm{a}$
 is $\kappa$-strutting for a point $p$.
Clearly $\Omega_{\bm{a}}$ is open.

According to Proposition~\ref{prop:stutt}, the set of points in $\spc{L}$ 
with rank $\ge \kay$ can be presented as
\[\bigcup_{\bm{a}}\Omega_{\bm{a}},\]
where the union is taken over all $\kay$-arrays $\bm{a}$ of points in $\spc{L}$.
Hence the result.
\qeds

\section{Right-inverse theorem}\label{sec:right-inverse-1}

Suppose that $\bm{a}=(a^1,\dots,a^\kay)$ is a point array in a metric space $\spc{L}$.
Recall that the map $\distfun{\bm{a}}{}{}\:\spc{L}\to\RR^n$ is defined by
\[\distfun{\bm{a}}{p}{}\df(\dist{a^1}{p}{},\dots,\dist{a^n}{p}{}).\] 

\begin{thm}{Right-inverse theorem}
\label{thm:right-inverse-function}{\sloppy 
Suppose $\spc{L}$ is a complete length $\Alex{\kappa}$ space,
$p,b\in\spc{L}$, 
and $\bm{a}=(a^1,\dots,a^\kay)$ is a point array in $\spc{L}$.

}

Assume that $(b,a^1,a^2,\dots,a^\kay)$ is $\kappa$-strutting for $p$.
Then the distance map $\distfun{\bm{a}}{}{}\:\spc{L}\to\RR^\kay$ has a right inverse defined in a neighborhood of $\distfun{\bm{a}}{p}\in\RR^\kay$;
that is, there is a submap $\map\:\RR^\kay\subto\spc{L}$ such that $\Dom \map\z\ni \distfun{\bm{a}}{p}$ and 
$\distfun{\bm{a}}{[\map(\bm{x})]}=\bm{x}$ for any $\bm{x}\in\Dom \map$.
Moreover,

\begin{subthm}{thm:right-inverse-function:Hoelder}
The map $\map$ can be chosen to be $C^{\frac{1}{2}}$-continuous (that is, Hölder continuous with exponent $\tfrac{1}{2}$) and such that 
\[\map(\distfun{\bm{a}}{p})=p.\]
\end{subthm}

\begin{subthm}{thm:right-inverse-function:open-map}
The distance map $\distfun{\bm{a}}{}{}\:\spc{L}\to\RR^\kay$ is locally co-Lipschitz (in particular, open) in a neighborhood of $p$.
\end{subthm}

\end{thm}

Part \ref{SHORT.thm:right-inverse-function:open-map} of the theorem 
is closely related to \cite[Theorem 5.4]{burago-gromov-perelman} by Yuriy Burago, Grigory Perelman, and Michael Gromov, 
but the proof presented here is different.
Yet another proof can be built on \cite[Proposition~4.3]{lytchak:open-map} by Alexander Lytchak.


\parit{Proof.} 
Fix $\eps,r,\lambda>0$ such that the following conditions hold: 
\begin{enumerate}[(i)]
\item Each distance function $\distfun{a^i}{}{}$ and $\distfun{b}{}{}$ is $\tfrac\lambda2$-concave in $\oBall(p,r)$.
\item For any $q\in \oBall(p,r)$, we have $\angk{\kappa}{q}{a^i}{a^j}>\tfrac\pi2+\eps$ for all $i\ne j$ and $\angk{\kappa}{q}{b}{a^i}\z>\tfrac\pi2+\eps$ for all $i$.
In addition, $\eps<\tfrac{1}{10}$.
\end{enumerate}

Given $\bm{x}=(x^1,x^2,\dots,x^\kay)\in \RR^\kay$, 
consider the function 
$f_{\bm{x}}\:\spc{L}\to \RR$ defined by
\[f_{\bm{x}}=\min_{i}\{h_{\bm{x}}^i\}+\eps\cdot\distfun{b}{}{},\]
where $h_{\bm{x}}^i(q)=\min\{0,\dist{a^i}{q}{}-x^i\}$.
Note that for any $\bm{x}\in\RR^\kay$, the function $f_{\bm{x}}$ is $(1+\eps)$-Lipschitz and $\lambda$-concave in $\oBall(p,r)$.
Denote by $\alpha_{\bm{x}}(t)$ the $f_{\bm{x}}$-gradient curve (see Chapter~\ref{chap:grad}) that starts at $p$.

\begin{clm}{}\label{clm:|a alpha_x|=x}
If for some $\bm{x}\in\RR^\kay$ and $t_0\le\tfrac{r}{2}$ we have
$|\distfun{\bm{a}}{p}-\bm{x}|
\le
\tfrac{\eps^2}{10}\cdot t_0$, then 
$
\distfun{\bm{a}}{[\alpha_{\bm{x}}(t_0)]}
= 
\bm{x}$.

\end{clm}

First note that \ref{clm:|a alpha_x|=x} follows if for any $q\in \oBall(p,r)$, we have
\begin{enumerate}[(i)]
\item\label{111} $(\dd_q\distfun{a^i}{}{})(\nabla_q f_{\bm{x}})<-\tfrac{1}{10}\cdot\eps^2$ if $\dist{a^i}{q}{}>x^i$ and
\item\label{222} $(\dd_q\distfun{a^i}{}{})(\nabla_q f_{\bm{x}})>\tfrac{1}{10}\cdot\eps^2$ if 
\[\dist{a^i}{q}{}-x^i=\min_j\{\dist{a^j}{q}{}\z-x^j\}<0.\]
\end{enumerate}
Indeed, since $t_0\le\tfrac{r}2$, then $\alpha_{\bm{x}}(t)\in\oBall(p,r)$ for all $t\in[0,t_0]$.
Consider the following real-to-real functions:
\[\begin{aligned}
\phi(t)
&\df
\max_{i}\{\dist{a^i}{\alpha_{\bm{x}}(t)}{}-x^i\},
\\
\psi(t)
&\df
\min_{i}\{\dist{a^i}{\alpha_{\bm{x}}(t)}{}-x^i\}.
\end{aligned}\eqlbl{eq:xy-def}\]
Then from (\ref{111}), 
we have $\phi^+<-\tfrac{1}{10}\cdot\eps^2$
if $\phi>0$ and $t\in[0,t_0]$.
Similarly, 
from (\ref{222}), 
we have $\psi^+>\tfrac{1}{10}\cdot\eps^2$
if $\psi<0$ and $t\in[0,t_0]$.
Since $|\distfun{\bm{a}}{p}-\bm{x}|
\le
\tfrac{\eps^2}{10}\cdot t_0$, it follows that $\phi(0)\le \tfrac{\eps^2}{10}\cdot t_0$ and $\psi(0)\ge -\tfrac{\eps^2}{10}\cdot t_0$.
Thus $\phi(t_0)\le 0$ and $\psi(t_0)\ge 0$.
On the other hand, from \ref{eq:xy-def} we have $\phi(t_0)\ge \psi(t_0)$.
That is, $\phi(t_0)=\psi(t_0)=0$; hence \ref{clm:|a alpha_x|=x} follows.

Thus, to prove \ref{clm:|a alpha_x|=x}, it remains to prove (\ref{111}) and (\ref{222}).
First let us prove it assuming that $\spc{L}$ is geodesic.

Note that 
\[(\dd_q\distfun{b}{}{})(\dir{q}{a^i})
\le\cos\angk{\kappa}{q}{b}{a^j}<-\tfrac\eps2
\eqlbl{inq-b}\]
for all $i$, and
\[(\dd_q\distfun{a^j}{}{})(\dir{q}{a^i})
\le
\cos\angk{\kappa}{q}{a^i}{a^j}
<
-\tfrac\eps2\eqlbl{inq-a^j}\]
for all $j\ne i$. 
Further, \ref{inq-a^j} implies
\[(\dd_q h_{\bm{x}}^j)(\dir{q}{a^i})\le 0.\eqlbl{inq-h}\]
for all $i\ne j$.
The assumption in (\ref{111}) implies
\[\dd_q f_{\bm{x}}
=
\min_{j\ne i} \{\dd_q h_{\bm{x}}^j\}+\eps\cdot(\dd_q\distfun{b}{}{}).\]
Thus
\begin{align*}
-(\dd_q\distfun{a^i}{}{})(\nabla_q f_{\bm{x}})
&\ge
\<\dir q{a^i},\nabla_q f_{\bm{x}}\>
\ge
\\
&\ge
(\dd_qf_{\bm{x}})(\dir q{a^i})
=
\\
&=
\min_{i\ne j}\{(\dd_qh_{\bm{x}}^i)(\dir q{a^i})\}+\eps\cdot(\dd_q\distfun{b}{}{})(\dir q{a^i}).
\end{align*}
Therefore (\ref{111}) follows from \ref{inq-b} and \ref{inq-h}.

The assumption in (\ref{222}) implies that $f_{\bm{x}}(q)
=
h_{\bm{x}}^i(q)+\eps\cdot\distfun{b}{}{}$ and 
\[\dd_q f_{\bm{x}}\le \dd_q \distfun{a^i}{}{}+\eps\cdot(\dd_p\distfun{b}{}{}).\] 
Therefore,
\begin{align*}
(\dd_q \distfun{a^i}{}{})(\nabla_q f_{\bm{x}})
&\ge 
\dd_qf_{\bm{x}}(\nabla_q f_{\bm{x}})
\ge 
\\
&\ge
\left[(\dd_qf_{\bm{x}})(\dir qb)\right]^2
\ge
\\
&\ge
\left[\min_i\{\cos\angk\kappa qb{a^i}\}-\eps^2\right]^2.
\end{align*}
Thus (\ref{222}) follows from \ref{inq-b}, since $\eps<\tfrac{1}{10}$. 

Therefore \ref{clm:|a alpha_x|=x} holds if $\spc{L}$ is geodesic.
If $\spc{L}$ is not geodesic,
perform the above estimate in $\spc{L}^\o$, the ultrapower of $\spc{L}$. 
(Recall that according to \ref{obs:ultralimit-is-geodesic}, $\spc{L}^\o$ is geodesic.)
This completes the proof of \ref{clm:|a alpha_x|=x}. 
\claimqeds

Set $t_0(\bm{x})=\tfrac{10}{\eps^2}\cdot|\distfun{\bm{a}}{p}-\bm{x}|$, 
giving equality in \ref{clm:|a alpha_x|=x}.
Define the submap $\map$ by
\[\map\:{\bm{x}}\mapsto \alpha_{\bm{x}}\circ t_0(\bm{x}),\quad 
\Dom\map=\oBall(\distfun{\bm{a}}{p},\tfrac{\eps^2\cdot r}{20} )\subset\RR^\kay.\]
It follows from \ref{clm:|a alpha_x|=x} that
$\distfun{\bm{a}}{[\map(\bm{x})]}=\bm{x}$ for any $\bm{x}\in\Dom\map$.

Clearly $t_0(p)=0$; thus $\map(\distfun{\bm{a}}{p})=p$.
Further, by construction of $f_{\bm{x}}$, 
\[|f_{\bm{x}}(q)-f_{\bm{y}}(q)|\le |\bm{x}-\bm{y}|,\]
for any $q\in \spc{L}$.
Therefore, according to Lemma~\ref{lem:fg-dist-est}, $\map$ is $C^{\frac{1}{2}}$-continuous.
Thus \ref{SHORT.thm:right-inverse-function:Hoelder}.

Further, note that 
\[\dist{p}{\map(\bm{x})}{}
\le (1+\eps)\cdot t_0(\bm{x})
\le\tfrac{11}{\eps^2}\cdot|\distfun{\bm{a}}{p}-\bm{x}| 
\eqlbl{co-lip}\]
holds.

The above construction may be repeated for any $p'\in \oBall(p,\tfrac{r}{4})$, $\eps'=\eps$, and $r'=\tfrac{r}{2}$.
The inequality \ref{co-lip} for the resulting map $\map'$ implies that for any $p',q \in \oBall(p,\tfrac{r}{4})$
there is $q'\in\spc{L}$ such that $\map'(q)=\map'(q')$ and 
\[\dist{p'}{q'}{}
\le 
\tfrac{11}{\eps^2}\cdot|\distfun{\bm{a}}{p'}-\bm{x}|.\]
That is, the distance map $\distfun{\bm{a}}{}{}$ is locally $\tfrac{11}{\eps^2}$-co-Lipschitz in $\oBall(p,\tfrac{r}{4})$.
\qeds

\section{Dimension theorem}\label{sec:dim>m}

The following theorem is the main result of this section.

\begin{thm}{Theorem}\label{thm:dim-infty}
Let $\spc{L}$ be a complete length $\Alex{\kappa}$ space, 
$q\in \spc{L}$, 
$R>0$ 
and $m\in \ZZ_{\ge0}$.
Then the following statements are equivalent:
\begin{subthmA}{LinDim} $\LinDim\spc{L}\ge m$.
\end{subthmA}

\begin{subthmA}{thm:dim-infty:rank}
There is a point $p\in\spc{L}$ that admits a $\kappa$-strutting array $(b,a^1,\dots,a^m)\in\spc{L}^{m+1}$.
\end{subthmA}

\begin{subthmA}{LinDim+} Let $\Euk^m$ be the set 
of all points $p\in \spc{L}$ such that there is a distance-preserving embedding $\EE^m\hookrightarrow \T_p$
that preserves the cone structure 
(see Section \ref{sec: tangent space}).
Then $\Euk^m$ contains a dense G-delta set in $\spc{L}$.
\end{subthmA}

\begin{subthmA}{TopDim}
There is a $C^{\frac{1}{2}}$-embedding
\[\cBall[1]_{\EE^m}\hookrightarrow \oBall(q,R);\]
that is, a bi-Hölder embedding with exponent $\tfrac{1}{2}$.
\end{subthmA}

\begin{subthmA}{pack} 
\[\pack_\eps \oBall(q,R)>\frac{\Const}{\eps^m}\]
for fixed $\Const>0$ and any $\eps>0$.
\end{subthmA}

\medskip

In particular:
\begin{enumerate}[(i)]
\item If $\LinDim\spc{L}=\infty$, then all the statements \ref{SHORT.LinDim+}, \ref{SHORT.TopDim}, and \ref{SHORT.pack} are satisfied for all $m\in\ZZ_{\ge0}$. 
\item 
 If the statement \ref{SHORT.TopDim} or \ref{SHORT.pack} is satisfied for some choice of $q\in \spc{L}$ and $R>0$, then it also is satisfied for any other choice of $q$ and $R$.
\end{enumerate}
\end{thm}

For finite-dimensional spaces, Theorem~\ref{thm:dim-finite} gives a stronger version 
of the theorem above.

The above theorem with the exception of statement~\ref{SHORT.TopDim} was proved by Conrad Plaut \cite{plaut:dimension}.
At that time, it was not known whether
\[\LinDim\spc{L}=\infty\quad \Rightarrow\quad \TopDim\spc{L}=\infty\]
for any complete length $\Alex\kappa$ space $\spc{L}$.
The latter implication was proved by Grigory Perelman and the third author \cite{perelman-petrunin:qg};
it was done by combining an idea of Conrad Plaut with the technique of gradient flow.
Part \ref{SHORT.TopDim} is somewhat stronger.

To prove Theorem \ref{thm:dim-infty} we will need the following three propositions.

\begin{thm}{Proposition}\label{E=T}
Let $p$ be a point in a a complete length $\Alex{\kappa}$ space
$\spc{L}$.
Assume there is a distance-preserving embedding $\iota\:\EE^{m}\hookrightarrow \T_p\spc{L}$ 
that preserves the cone structure.
Then either
\begin{subthm}{}
 $\Im\iota=\T_p\spc{L}$, or
\end{subthm}

\begin{subthm}{} there is a point $p'$ arbitrarily close to $p$ such that there is a distance-preserving embedding $\iota'\:\EE^{m+1}\hookrightarrow \T_{p'}\spc{L}$ 
that preserves the cone structure.
\end{subthm}
\end{thm}

\parit{Proof.}
Assume $\iota(\EE^{m})$ is a proper subset of $\T_p\spc{L}$.
Equivalently, there is a direction $\xi \in \Sigma_p\setminus\iota(\mathbb{S}^{m-1})$,
where $\mathbb{S}^{m-1}\subset \EE^m$ is the unit sphere. 

Fix $\eps>0$ so that $\mangle(\xi,\sigma)>\eps$ for any $\sigma\in \iota(\mathbb{S}^{m-1})$. 
Choose a maximal $\eps$-packing in $\iota(\mathbb{S}^{m-1})$;
that is, an array $(\zeta^1,\zeta^2,\dots,\zeta^n)$ of directions in $\iota(\mathbb{S}^{m-1})$ such that $n=\pack_\eps \mathbb{S}^{m-1}$ and $\mangle(\zeta^i,\zeta^j)>\eps$ for any $i\ne j$.

Choose an array $(x,z^1,z^2,\dots,z^n)$ of points in $\spc{L}$ such that
$\dir p x\approx\xi$, $\dir p{z^i}\approx\zeta^i$; 
here we write ``$\approx$'' for ``sufficiently close''.
We can choose this array so that
$\angk{\kappa}p x{z^i}>\eps$ for all $i$ 
and $\angk{\kappa}p{z^i}{z^j}>\eps$ for all $i\ne j$.
Applying Corollary \ref{cor:euclid-subcone}, there is a point $p'$ arbitrarily close to $p$ 
such that all directions $\dir{p'}x$, $\dir{p'}{z^1}$, $\dir{p'}{z^2},\dots,\dir{p'}{z^n}$
belong to an isometric copy of $\mathbb{S}^{\kay-1}$ in $\Sigma_{p'}$.
In addition, we may assume that $\angk{\kappa}{p'}x{z^i}>\eps$ and $\angk{\kappa}{p'}{z^i}{z^j}>\eps$.
From the hinge comparison (\ref{angle}),
$\mangle(\dir{p'}x,\dir{p'}{z^i})>\eps$ 
and $\mangle(\dir{p'}{z^i},\dir{p'}{z^j})>\eps$;
that is, 
\[\pack_\eps \mathbb{S}^{\kay-1}\ge n+1>\pack_\eps \mathbb{S}^{m-1}.\] 
Hence $\kay>m$.
\qeds

\begin{thm}{Proposition}\label{pack-homogeneus}
Let $\spc{L}$ be a complete length $\Alex{\kappa}$ space. Then 
 for any two points $p,\bar p\in \spc{L}$ and any $R,\bar R>0$, there is a constant $\delta\z=\delta(\kappa,R,\bar R,\dist{p}{\bar p}{})>0$ such that
\[\pack_{\delta\cdot\eps}\oBall(\bar p,\bar R)\ge \pack_{\eps}\oBall(p,R).\]

\end{thm}

\parit{Proof.} According to \ref{cor:CAT>k-sence}, we can assume that $\kappa\le 0$.

Let $n=\pack_{\eps}\oBall(p,R)$ and ${\{x^1,\dots, x^n}\}$ be a maximal $\eps$-packing in $\oBall(p,R)$;
that is, $\dist{x^i}{x^j}{}>\eps$ for all $i\ne j$.
Without loss of generality, we may assume the $x^i$ are in $\Str(\bar p)$.
Thus, for each $i$ there is a unique geodesic $[\bar p x^i]$ (see \ref{thm:almost.geod}).
Choose a factor $1>s>0$ so that $\bar R>s\cdot(\dist{p}{\bar p}{}+R)$.
For each $i$, take $\bar x^i\in[\bar p x^i]$ so that 
$\dist{\bar p}{\bar x^i}{}=s\cdot(\dist{p}{x^i}{})$.
From \ref{cor:monoton:2-sides},
\[\angkk\kappa {\bar p}{\bar x^i}{\bar x^j}{}\ge\angk\kappa {\bar p}{x^i}{x^j}.\]

\begin{wrapfigure}{r}{50 mm}
\vskip-4mm
\centering
\includegraphics{mppics/pic-1405}
\vskip0mm
\end{wrapfigure}

The cosine law gives a constant $\delta\z=\delta(\kappa,R,\bar R,\dist{p}{\bar p}{})>0$ such that 
\[\dist{\bar x^i}{\bar x^j}{}>\delta\cdot(\dist{x^i}{x^j}{})>\delta\cdot\eps\] 
for all $i\ne j$.
Hence the statement follows.
\qeds

\begin{thm}{Proposition}\label{E-comeagre} 
Let $\spc{L}$ be a complete length $\Alex{\kappa}$ space, 
$r\z<\varpi\kappa$ 
and $p\in \spc{L}$.
Assume that 
\[\pack_{\eps} \oBall(p,r)
>\pack_{\eps}\cBall[r]_{\Lob{m}{\kappa}}
\eqlbl{eq:pack>pack}\]
for $\eps>0$.
Then there is a G-delta set $A\subset \spc{L}$
that is dense in a neighborhood of $p$ and
such that $\dim\Lin_q>m$ for any $q\in A$.
\end{thm}

\parit{Proof.} 
Choose a maximal $\eps$-packing in $\oBall(p,r)$,
that is, an array $(x^1,x^2,\dots, x^n)$ of points in $\oBall(p,r)$ such that $n=\pack_\eps \oBall(p,r)$ and $\dist{x^i}{x^j}{}>\eps$ for any $i\ne j$.
Choose a neighborhood $\Omega\ni p$
such that $\dist{q}{x^i}{}<r$ for any $q\in \Omega$ and all $i$.
Let 
\[A= \Omega\cap\Str(x^1,x^2,\dots,x^n).\]
According to Theorem \ref{thm:almost.geod}, $A$ is a G-delta set that is dense in $\Omega$.

Assume $\kay=\dim\Lin_q\le m$ for $q\in A$.
Consider an array $(v^1,v^2,\dots,v^n)$ of vectors in $\Lin_q$,
where $v^i=\ddir{q}{x^i}$.
Clearly 
\[|v^i|=\dist{q}{x^i}{}<r,\] 
and from the hinge comparison (\ref{angle})
we have 
\[\side\kappa \hinge \0{v^i}{v^j}\ge \dist{x^i}{x^j}{}>\eps.\]
Note that the ball $\oBall(\0,r)_{\Lin_q}$ equipped with the metric $\rho(v,w)\z=\side\kappa \hinge \0{v}{w}$ is isometric to 
$\cBall[r]_{\Lob{\kay}{\kappa}}$.
Thus
\[
\pack_\eps\cBall[r]_{\Lob{\kay}{\kappa}}
\ge
\pack_\eps \oBall(p,r),
\]
which contradicts $\kay\le m$ and \ref{eq:pack>pack}.
\qeds

The proof of Theorem \ref{thm:dim-infty} is essentially done in \ref{E=T}, \ref{pack-homogeneus}, \ref{E-comeagre}, \ref{thm:inverse-function},
\ref{thm:right-inverse-function}; 
now we assemble the proof from these parts.

\parit{Proof of \ref{thm:dim-infty}.} 
We will prove the implications 
\[\textrm{\ref{SHORT.LinDim+} 
$\Rightarrow$ 
\ref{SHORT.LinDim} 
$\Rightarrow$ 
\ref{SHORT.thm:dim-infty:rank}
$\Rightarrow$ 
\ref{SHORT.pack}
$\Rightarrow$ 
\ref{SHORT.LinDim+} 
$\Rightarrow$ 
\ref{SHORT.TopDim} 
$\Rightarrow$ 
\ref{SHORT.pack}.}\]
The implication \ref{SHORT.LinDim+}$\Rightarrow$\ref{SHORT.LinDim} is trivial.

\parit{\ref{SHORT.LinDim}$\Rightarrow$\ref{SHORT.thm:dim-infty:rank}.}
Choose a point $p\in\spc{L}$ such that $\dim\Lin_p\ge m$.
Clearly one can choose an array $(\xi^0,\xi^1,\dots,\xi^m)$ of directions in $\Lin_p$ such that $\mangle(\xi^i,\xi^j)\z>\tfrac\pi2$ for all $i\ne j$.
Choose an array $(x^0,x^1,\dots,x^m)$ of points in $\spc{L}$ such that each $\dir{p}{x^i}$ is sufficiently close to $\xi^i$;
in particular, we have $\mangle\hinge{p}{x^i}{x^j}>\tfrac\pi2$.
Choose points $a^i\in\mathopen{]}p x^i\mathclose{]}$ sufficiently close to $p$.
This can be done so that each $\angk\kappa p{a^i}{a^j}$ is arbitrarily close to $\mangle\hinge p{a^i}{a^j}$,
in particular $\angk\kappa p{a^i}{a^j}>\tfrac{\pi}{2}$.
Finally, set $b=a^0$.

\parit{\ref{SHORT.thm:dim-infty:rank}$\Rightarrow$\ref{SHORT.pack}.} 
Let $p\in \spc{L}$ be a point that admits a $\kappa$-strutting array $(b,a^1,\dots, a^m)$ 
of points in $\spc{L}$.
The right-inverse theorem (\ref{thm:right-inverse-function:open-map})
implies that the distance map $\distfun{\bm{a}}{}{}\:\spc{L}\to\RR^m$,
\[\distfun{\bm{a}}{}{}\:x\mapsto(\dist{a^1}{x}{},\dist{a^2}{x}{},\dots,\dist{a^n}{x}{}),\]
is open in a neighborhood of $p$.
Since the distance map $\distfun{\bm{a}}{}{}$ is Lipschitz, 
for any $r>0$, there is $\Const>0$ such that
\[\pack_\eps \oBall(p,r)>\frac{\Const}{\eps^m}\]
for any $\eps>0$.
Applying \ref{pack-homogeneus}, we get a similar inequality for any other ball in $\spc{L}$;
that is, for any $q\in\spc{L}$ and $R>0$, there is $\Const'>0$ such that 
\[\pack_\eps \oBall(q,R)>\frac{\Const'}{\eps^m}.\]

\parit{\ref{SHORT.pack}$\Rightarrow$\ref{SHORT.LinDim+}.} 
Note that for any $q'\in\spc{L}$ and $R'>\dist{q}{q'}{}+R$ we have
\begin{align*}
\pack_\eps\oBall(q',R')
&\ge
\pack_\eps\oBall(q,R)
\ge
\\
&\ge
\frac{\Const}{\eps^m}
>
\\
&>
\pack_\eps\cBall[R']_{\Lob{m-1}{\kappa}}
\end{align*}
for all sufficiently small $\eps>0$.
Applying \ref{E-comeagre},
$\Euk^m$
contains a G-delta set that is dense in a neighborhood of any point $q'\in\spc{L}$.

\parit{\ref{SHORT.LinDim+}$\Rightarrow$\ref{SHORT.TopDim}.} 
Since $\Euk^m$ contains a dense G-delta set in $\spc{L}$, we can choose $p\in \oBall(q,R)$ with a distance-preserving cone embedding $\iota\:\EE^m\hookrightarrow \T_p$.

Repeating the construction in \ref{SHORT.LinDim}$\Rightarrow$\ref{SHORT.thm:dim-infty:rank}, 
we get a $\kappa$-strutting array $(p,a^1,\dots, a^m)$ for $p$.

Applying the right-inverse theorem (\ref{thm:right-inverse-function}),
we obtain a $C^{\frac{1}{2}}$-submap 
\[\map\:\RR^m\subto \oBall(q,R)\]
that is a right inverse for $\distfun{\bm{a}}{}{}\:\spc{L}\to\RR^m$ and such that $\map(\distfun{\bm{a}}{p})=p$.
In particular, $\map$ is a $C^{\frac{1}{2}}$-embedding of $\Dom\map$.

\parit{\ref{SHORT.TopDim}$\Rightarrow$\ref{SHORT.pack}.}
This proof is valid for general metric spaces;
it is based on general relations between topological dimension, Hausdorff measure and $\pack_\eps$. 

Let $W\subset\oBall(q,R)$ be the image of the embedding $\map$.
Since $\TopDim W\z=m$,
Szpilrajn's theorem (\ref{thm:szpilrajn}) implies that
\[\HausMes_m W>0.\]

Given $\eps>0$, consider a maximal $\eps$-packing of $W$, 
that is, an array $(x^1,x^2,\dots,x^n)$ of points in $W$ such that $n=\pack_\eps W$ and $\dist{x^i}{x^j}{}>\eps$ for all $i\ne j$.
Note that $W$ is covered by balls $\oBall(x^i,2\cdot\eps)$.

By the definition of Hausdorff measure,
\[\pack_\eps W
\ge
\frac{\Const}{\eps^m}\cdot\HausMes_m W\]
for a fixed constant $\Const>0$ and all small $\eps>0$.
Hence \ref{SHORT.pack} follows.
\qedsf

\section{Inverse function theorem}

\begin{thm}{Inverse function theorem}\label{thm:inverse-function}
Let $\spc{L}$ be an $m$-dimensional complete length $\Alex\kappa$ space 
and $p,b,a^1,a^2,\dots,a^m\in\spc{L}$.

Assume that the point array $\bm{a}=(b,a^1,\dots,a^m)$ is $\kappa$-strutting for~$p$.
Then
there are $R>0$ and $\eps>0$ such that:

\begin{subthm}{thm:inverse-function:strut}
For all $i\ne j$ and any $q\in\oBall(p,R)$ we have
\[\angk\kappa{q}{a^i}{a^j}>\tfrac{\pi}{2}+\eps
\quad\text{and}\quad
\angk\kappa{q}{b}{a^i}>\tfrac{\pi}{2}.\]

\end{subthm}

\begin{subthm}{thm:inverse-function:chart}
The restriction of the distance map 
\[\distfun{\bm{a}}{}{}\:x\mapsto(\dist{a^1}{x}{},\dots,\dist{a^m}{x}{})\] 
to the ball $\oBall(p,R)$
is an open $[\eps,\sqrt{m}]$-bi-Lipschitz embedding $\oBall(p,R)\hookrightarrow\RR^m$.
\end{subthm}

\begin{subthm}{thm:inverse-function:R}
The value $R$ depends only on $\kappa$, $\dist{p}{a^i}{}$, $\dist{a^i}{a^j}{}$, and $\dist{b}{a^i}{}$
for all $i$ and $j$.
\end{subthm}

\end{thm}

\begin{thm}{Definition}\label{def:dist-chart}
Suppose $\spc{L}$ is an $m$-dimensional complete length $\Alex\kappa$ space.
If a point array $(b,a^1,a^2,\dots,a^m)$ 
and the value $R$ satisfy the conditions in Theorem~\ref{thm:inverse-function}, 
then the restriction 
$\bm{x}=\distfun{\bm{a}}{}{}|_{\oBall(p,R)}$
is called a \index{distance chart}\emph{distance chart},
the restrictions $x^i=\distfun{a^i}{}{}|_{\oBall(p,R)}$ are called \index{distance chart!coordinates of a distance chart}\emph{coordinates}, 
and the restriction $y=\distfun{b}{}{}|_{\oBall(p,R)}$ is called \index{distance chart!strut of a distance chart} the \emph{strut} of the distance chart.
\end{thm}

\begin{thm}{Lemma}\label{lem:pack(S^m)+}
Let $p$ be a point in an $m$-dimensional complete length $\Alex\kappa$ space $\spc{L}$.
Assume for the directions $\xi,\zeta^1,\zeta^2,\dots,\zeta^\kay\in\Sigma_p$ the following conditions hold: 

\begin{subthm}{}
$\mangle(\xi,\zeta^i)>\tfrac\pi2-\eps$ for all $i$,
\end{subthm}

\begin{subthm}{}
$\mangle(\zeta^i,\zeta^j)>\tfrac\pi2+\eps$ for all $i\ne j$.
\end{subthm}
Then $\kay\le m$.
\end{thm}

\parit{Proof.}
Without loss of generality, we can assume that all $\xi,\zeta^1,\zeta^2,\dots,\zeta^\kay$ are geodesic directions;
let $\xi=\dir{p}{x}$ and $\zeta^i=\dir{p}{z^i}$ for all $i$.
Fix a small $r>0$, and
let $\bar x\in \mathopen{]}px]$ and $\bar z^i\in\mathopen{]}p z^i]$ be points 
such that
\[\dist{p}{\bar x}{}=\dist{p}{\bar z^1}{}=\dots=\dist{p}{\bar z^\kay}{}=r.\]
From the definition of angle,
if $r$ is sufficiently small we have
\begin{itemize}
\item $\angk\kappa{p}{\bar x}{\bar z^i}>\tfrac\pi2-\eps$ for all $i$,
and $\angk\kappa{p}{\bar z^i}{\bar z^j}>\tfrac\pi2+\eps$ for all $i\ne j$.
\end{itemize}
Choose a point $p'\in\Str(\bar x,\bar z^1,\bar z^2,\dots,\bar z^\kay)$ sufficiently close to $p$ that the above conditions still hold for $p'$; that is,
\begin{clm}{}
 $\angk\kappa{p'}{\bar x}{\bar z^i}>\tfrac\pi2-\eps$ for all $i$, and $\angk\kappa{p'}{\bar z^i}{\bar z^j}>\tfrac\pi2+\eps$ for all $i\ne j$.
\end{clm}
Set $\acute\xi=\dir{p'}{\bar x}$ and $\acute\zeta^i=\dir{p'}{\bar z^i}$ for each $i$.
By the hinge comparison (\ref{angle}), 
\begin{clm}{}
$\mangle(\acute\xi,\acute\zeta^i)>\tfrac\pi2-\eps$ for all $i$, and $\mangle(\acute\zeta^i,\acute\zeta^j)>\tfrac\pi2+\eps$ for all $i\ne j$.
\end{clm}

According to Corollary~\ref{cor:euclid-subcone}, all directions $\acute\xi,\acute\zeta^1,\acute\zeta^2,\dots,\acute\zeta^\kay$ lie in an isometric copy of the standard $n$-sphere in $\Sigma_{p'}$. Clearly $n\le m-1$.
Thus it remains to prove the following claim, which is a partial case of the lemma.

\begin{clm}{}
If $\xi,\zeta^1,\zeta^2,\dots,\zeta^\kay\in\mathbb{S}^{m-1}$,
$\dist{\xi}{\zeta^i}{}>\tfrac\pi2-\eps$ for all $i$, and
$\dist{\zeta^i}{\zeta^j}{}>\tfrac\pi2+\eps$ for all $i\ne j$,
then $\kay\le m$.
\end{clm}

For each $i$, 
let 
$\bar\zeta^i$
be the closest point to $\zeta^i$
in
$\Xi=\mathbb{S}^{m-1}\setminus \oBall(\xi,\tfrac\pi2)
\iso
\mathbb{S}^{m-1}_+$ 
(if $\zeta\in\Xi$, then $\bar\zeta^i=\zeta^i$).
By straightforward calculations, we have
\[\dist{\bar\zeta^i}{\bar\zeta^j}{}\ge \dist{\zeta^i}{\zeta^j}{}-\eps>\tfrac\pi2.\]
Thus it is sufficient to show the following claim:

\begin{clm}{}
$\pack_{\frac\pi2}\mathbb{S}^{m-1}_+= m.$
\end{clm}

Clearly, $\pack_{\frac\pi2}\mathbb{S}^{m-1}_+\ge m$.

The opposite inequality is proved by induction on $m$.
The base case $m=1$ is obvious. 
Assume $(\bar\zeta^1,\bar\zeta^2,\dots,\bar\zeta^{\kay})$ is an array of points in $\mathbb{S}^{m-1}_+$ with $\dist{\bar\zeta^i}{\bar\zeta^j}{}>\tfrac\pi2$.
Without loss of generality, we can also assume that $\bar\zeta^\kay\in\partial \mathbb{S}^{m-1}_+$.
For each $i<\kay$, 
let $\check\zeta^i
=\dir{\bar\zeta^\kay}{\bar\zeta^i}\in\Sigma_{\bar\zeta^\kay}\mathbb{S}^{m-1}_+\iso\mathbb{S}^{m-2}_+$.
By the hinge comparison (\ref{angle}), $\mangle(\check\zeta^i,\check\zeta^j)>\tfrac\pi2$ 
for all $i<j<\kay$.
Thus from the induction hypothesis we have $\kay-1\le {m-1}$.
\qeds

\parit{Proof of \ref{thm:inverse-function}; \ref{SHORT.thm:inverse-function:strut}.} 
Fix $\eps>0$ such that 
$\angk\kappa p{a^i}{a^j}\z>\tfrac\pi2+\eps$ and $\angk\kappa p{b}{a^i}\z>\tfrac\pi2+\eps$ for all $i\ne j$.
Choose $R>0$ sufficiently small that 
$\angk\kappa q{a^i}{a^j}>\tfrac\pi2+\eps$ and $\angk\kappa q{b}{a^i}>\tfrac\pi2+\eps$ for all $i\ne j$ and any $q\in\oBall(p,R)$.
Clearly, \ref{SHORT.thm:inverse-function:strut} holds for $\oBall(p,R)$.

\parit{\ref{SHORT.thm:inverse-function:chart}.}
Note that the distance map $\distfun{\bm{a}}{}{}$ is Lipschitz
and its restriction $\distfun{\bm{a}}{}{}|_{\oBall(p,R)}$ is open;
the latter follows from the right-inverse theorem (\ref{thm:right-inverse-function:open-map}).
Thus to prove \ref{SHORT.thm:inverse-function:chart}, it is sufficient to show that
\[
\max_{i}\left\{\,\bigl|\dist{a^i}{x}{}-\dist{a^i}{y}{}\bigr|\,\right\}
>
\tfrac\eps2\cdot\dist[{{}}]{x}{y}{}
\eqlbl{expend}
\]
for any $x,y\in {\oBall(p,R)}$.

According to Lemma~\ref{lem:pack(S^m)+}, 
\[
\mangle\hinge{x}y{b}
\le
\tfrac\pi2-\eps
\quad\text{or}\quad
\mangle\hinge{x}y{a^i}
\le
\tfrac\pi2-\eps
\quad \text{for some}\quad i.
\]
In the latter case,
since $\dist{x}{y}{}<2\cdot R$ and $R$ is small, 
the hinge comparison (\ref{angle}) implies 
\[
\dist{a^i}{x}{}-\dist{a^i}{y}{}>\tfrac\eps2\cdot\dist[{{}}]{x}{y}{}
\quad \text{for some}\quad i.
\eqlbl{eq:y-x}\]

If $\mangle\hinge{x}y{b}
\le
\tfrac\pi2-\eps$, then 
switching $x$ and $y$, we get
\[
\dist{a^j}{y}{}-\dist{a^j}{x}{}
>
\tfrac\eps2\cdot\dist[{{}}]{x}{y}{}
\quad \text{for some}\quad j.
\eqlbl{eq:x-y}\] 
Then \ref{eq:y-x} and \ref{eq:x-y} imply \ref{expend}.

Finally, part \ref{SHORT.thm:inverse-function:R}
follows since the angle $\angk\kappa q{a^i}{a^j}$ 
depends continuously on $\kappa$, $\dist{q}{a^i}{}$, $\dist{q}{a^j}{}$ and $\dist{a^i}{a^j}{}$.
\qeds

\section{Finite-dimensional spaces}\label{sec:dim=m}

The next theorem is a refinement of \ref{thm:dim-infty} for the finite-dimensional case;
it was essentially proved by Yuriy Burago, Grigory Perelman, and Michael Gromov \cite{burago-gromov-perelman}.

\begin{thm}{Theorem}\label{thm:dim-finite}
Suppose $\spc{L}$ is a complete length $\Alex{\kappa}$ space, 
$m$ is a nonnegative integer,
$0<R\le \varpi\kappa$, and
$q\in \spc{L}$.
Then the following statements are equivalent:
\begin{subthm}{LinDim-f} $\LinDim\spc{L}= m$.
\end{subthm}

\begin{subthm}{thm:dim-finite:rank}
$m$ is the maximal integer such that there is a point $p\in\spc{L}$ that admits a $\kappa$-strutting array $(b,a^1,\dots,a^m)$.
\end{subthm}

\begin{subthm}{LinDim+-f} $\T_p\iso \EE^m$ for any point $p$ in a dense G-delta set of $\spc{L}$.
\end{subthm}

\begin{subthm}{TopDim-f} There is an open bi-Lipschitz embedding 
\[\cBall[1]_{\EE^m}\hookrightarrow \oBall(q,R)\subset \spc{L}.\]
\end{subthm}

\begin{subthm}{pack-f} For any $\eps>0$,
\[\pack_\eps\cBall[R]_{\Lob{m}{\kappa}} \ge\pack_\eps \oBall(q,R).\]
moreover, there is $\Const=\Const(q,R)>0$ such that 
\[\pack_\eps \oBall(q,R)>\frac\Const{\eps^m}.\]
\end{subthm}

\end{thm}

Using theorems \ref{thm:dim-infty} and \ref{thm:dim-finite}, 
one can show that linear dimension is equal to many different types of dimension, such 
as {}\emph{small} and \index{inductive dimension}\emph{big inductive dimension} 
and {}\emph{upper} and \index{box counting dimensions}\emph{lower box-counting dimension}
(also known as \index{Minkowski dimension}\emph{Minkowski dimension}), 
\index{homological dimension}\emph{homological dimension} and so on.

The next two corollaries follow from \ref{SHORT.pack-f}.

\begin{thm}{Corollary}\label{cor:dim>proper}
Any finite-dimensional complete length $\Alex{}$ space is proper and geodesic.
\end{thm}

\begin{thm}{Corollary} Let $(\spc{L}_n)$ be a sequence of length $\Alex\kappa$ spaces and $\spc{L}_n\to \spc{L}_\o$ as $n\to\o$.
Assume $\LinDim L_n\le m$ for all $n$.
Then $\LinDim L_\o\le m$.
\end{thm}

\begin{thm}{Corollary}\label{dim=dim} 
Let $\spc{L}$ be a complete length $\Alex{\kappa}$ space. 
Then for any open $\Omega\subset \spc{L}$, we have
\[
\LinDim \spc{L}=
\LinDim\Omega =
\TopDim\Omega=
\HausDim\Omega,
\]
where $\TopDim$ and $\HausDim$ denote topological dimension (\ref{def:TopDim}) and Hausdorff dimension (\ref{def:HausDim}) respectively.

In particular, $\spc{L}$ is dimension-homogeneous; that is, all open sets have the same linear dimension.
\end{thm}

\parit{Proof of \ref{dim=dim}.} 
The equality
\[\LinDim \spc{L}= \LinDim\Omega\]
follows from \ref{LinDim}$\&$\ref{SHORT.LinDim+}.

If $\LinDim \spc{L}=\infty$, then
applying \ref{TopDim} for $\oBall(q,R)\subset \Omega$, we find that there is a compact subset $K\subset \Omega$ having an arbitrarily large $\TopDim K$. Therefore
\[\TopDim\Omega=\infty.\] 
By Szpilrajn's theorem (\ref{thm:szpilrajn}),
$\HausDim K\ge \TopDim K$.
Thus we also have 
\[\HausDim\Omega=\infty.\]

If $\LinDim \spc{L}=m<\infty$, then the first inequality in \ref{pack-f} 
implies that \[\HausDim \oBall(q,R)\le m.\] 
According to Corollary~\ref{cor:dim>proper}, 
$\spc{L}$ is proper and in particular has countable base. 
Thus applying Szpilrajn's theorem again, we have
\[\TopDim\Omega\le \HausDim \Omega\le m.\]
Finally, \ref{TopDim-f} implies that $m\le\TopDim\Omega$.
\qeds

\parit{Proof of \ref{thm:dim-finite}.}
The equivalence \ref{SHORT.LinDim-f}$\Leftrightarrow$\ref{SHORT.thm:dim-finite:rank} follows from \ref{thm:dim-infty}.

\parit{\ref{SHORT.LinDim-f}$\Rightarrow$\ref{SHORT.LinDim+-f}.}
If $\LinDim\spc{L}=m$, then by Theorem~\ref{thm:dim-infty}, 
$\Euk^m$ contains a dense G-delta set in $\spc{L}$.
From \ref{E=T}, it follows that $\T_p$ is isometric to $\EE^m$ for any $p\in \Euk^m$.

\parit{\ref{SHORT.LinDim+-f}$\Rightarrow$\ref{SHORT.TopDim-f}.}
This is proved in exactly the same way as implication \ref{SHORT.LinDim+}$\Rightarrow$\ref{SHORT.TopDim} of theorem \ref{thm:dim-infty}, 
but applying the existence of a distance chart (\ref{thm:inverse-function}) 
instead of the right-inverse theorem (\ref{thm:right-inverse-function}).

\parit{\ref{SHORT.TopDim-f}$\Rightarrow$\ref{SHORT.pack-f}.} 
From \ref{SHORT.TopDim-f}, it follows that there is a point $p\in\oBall(q,R)$ and $r>0$ such that
$\oBall(p,r)\subset \spc{L}$ is bi-Lipschitz homeomorphic to a bounded open set of $\EE^m$.
Thus there is $\Const>0$ such that 
\[\pack_\eps \oBall(p,r)>\frac{\Const}{\eps^m}.\eqlbl{eq:thm:dim-finite*}\]
Applying \ref{pack-homogeneus} shows that inequality \ref{eq:thm:dim-finite*}, with different constants, holds for any other ball, in particular for $\oBall(q,R)$.

Applying \ref{E-comeagre} gives the first inequality in \ref{SHORT.pack-f}.

\parit{\ref{SHORT.pack-f}$\Rightarrow$\ref{SHORT.LinDim-f}.} 
From Theorem \ref{thm:dim-infty}, we have $\LinDim\spc{L}\ge m$. 
Applying Theorem \ref{thm:dim-infty} again, if $\LinDim\spc{L}\ge m+1$ then for some $\Const>0$ and any $\eps>0$,
\[\pack_\eps \oBall(q,R)\ge \frac{\Const}{\eps^{m+1}}.\]
But
\[\frac{\Const'}{\eps^m}\ge\pack_\eps \oBall(q,R)\] 
for any $\eps>0$,
a contradiction.
\qeds

\begin{thm}{Exercise}\label{ex:compact-dimension-cbb}
Suppose $\spc{L}$ is a complete length $\Alex{}$ space and $\Sigma_p\spc{L}$ is compact for any $p\in\spc{L}$.
Prove that $\spc{L}$ is finite-dimensional.
\end{thm}

\section{One-dimensional spaces}

\begin{thm}{Theorem}\label{thm:dim=1.CBB} 
Let $\spc{L}$ be a one-dimensional complete length $\Alex\kappa$ space.
Then $\spc{L}$ is isometric to a connected complete Riemannian one-dimensional manifold with possibly non-empty boundary.
\end{thm}

\parit{Proof.}
Clearly $\spc{L}$ is connected.
It remains to show the following:
\begin{clm}{}\label{clm:1-dim-all}
For any point $p\in\spc{L}$
there is $\eps>0$ such that $\oBall(p,\eps)$ 
is isometric to either $[0,\eps)$ or $(-\eps,\eps)$.
\end{clm}

First let us show:
\begin{clm}{}\label{clm:1-dim-mid}
If $p\in\mathopen{]}x y\mathclose{[}$ for $x$, $y\in\spc{L}$ and $\eps<\min\{\dist{p}{x}{},\dist{p}{y}{}\},$
then $\oBall(p,\eps)\subset\mathopen{]}x y\mathclose{[}$.
In particular,
$\oBall(p,\eps)\iso(-\eps,\eps)$.
\end{clm}

\begin{wrapfigure}{r}{33 mm}
\vskip0mm
\centering
\includegraphics{mppics/pic-1410}
\vskip0mm
\end{wrapfigure}

Assume the contrary;
that is, there is 
$$z\in \oBall(p,\eps)\setminus\mathopen{]}x y\mathclose{[}.$$
Consider a geodesic $[p z]$, and let $q\in[p z]\cap[x y]$ be the point that maximizes the distance $\dist{p}{q}{}$.
At $q$, we have three distinct directions: 
to $x$, $y$, and $z$.
Moreover, $\mangle\hinge{q}{x}{y}=\pi$.
Thus, according to Proposition~\ref{E=T}, 
$\LinDim\spc{L}>1$, a contradiction.
\claimqeds

Now we assume no geodesic includes $p$ as a non endpoint.
Since $\LinDim\spc{L}\z=1$ there is a point $y\ne p$.

Fix a positive value $\eps<\dist{p}{y}{}$.
Let us show:

{

\begin{wrapfigure}{r}{33 mm}
\vskip0mm
\centering
\includegraphics{mppics/pic-1415}
\vskip0mm
\end{wrapfigure}

\begin{clm}{}\label{clm:1-dim-end}
$\oBall(p,\eps)\subset [p y]$;
in particular, $\oBall(p,\eps)\iso[0,\eps)$.
\end{clm}

Assume the contrary;
let $z\in \oBall(p,\eps)\setminus[p y]$.

Choose a point $w\in \mathopen{]} p y \mathclose{[}$ such that 
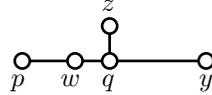
\[\dist{p}{w}{}+\dist{p}{z}{}<\eps.\]
Consider geodesic $[w z]$, and let $q\in[p y]\cap[w z]$ be the point that maximizes the distance $\dist{w}{q}{}$.
Since no geodesic includes $p$ as a non endpoint, we have $p\ne q$.
As above, $\mangle\hinge{q}{p}{y}=\pi$ 
and $\dir{q}{z}$ is distinct from $\dir{q}{p}$ and $\dir{q}{p}$.
Thus, according to Proposition~\ref{E=T}, 
$\LinDim\spc{L}>1$, a contradiction.
\claimqeds

}

Clearly $\text{\ref{clm:1-dim-mid}}+\text{\ref{clm:1-dim-end}}\Rightarrow\text{\ref{clm:1-dim-all}}$;
hence the result.
\qeds

\chapter{Gradient flow}\label{chap:grad}

Gradient flow could be considered as a nonsmooth version of first-order ordinary differential equations.
It provides a universal tool in Alexandrov geometry with most significant applications to $\Alex{}{}$ spaces. 

The theory of gradient flows of semiconvex functions on Hilbert spaces (which are of course both  $\Alex{0}{}$ and $\CAT{0}$)  is classical, see for example \cite{Brezis-book}.

 The technique of gradient flows in the context of comparison geometry takes its roots in \index{Sharafutdinov's retraction}\emph{Sharafutdinov's retraction}, 
introduced by Vladimir Sharafutdinov \cite{sharafutdinov}. 
It has been used  widely in comparison geometry
since then.
In $\Alex{}$ spaces, it was first used by Grigory Perelman and the third author \cite{perelman-petrunin:qg, petrunin:qg}.
A bit later, independently Jürgen Jost and Uwe Mayer \cite{jost,mayer} 
used the gradient flow in $\CAT{}$ spaces.
Later, Alexander Lytchak unified and generalized these two approaches
to a wide class of metric spaces \cite{lytchak:open-map}.
It was developed  yet further by Shin-ichi Ohta \cite{ohta} and by Giuseppe Sevar\'e \cite{sevare}.
It is based on the more  analytic approach suitable for the study of synthetic spaces with lower Ricci bounds was developed by Luigi Ambrosio, Nicola Gigli and  Giuseppe Sevar\'e in a general metric and metric measure setting 
\cite{ambrosio-gigli-savare}.

{\sloppy 
The following exercise is a stripped-down version of Sharfutdinov's retraction;
it gives the idea behind gradient flow.

}

\begin{thm}{Exercise}\label{ex:sharafutdinov}
Assume that a one-parameter family of convex sets $K_t\subset \EE^m$ is nested; 
that is, $K_{t_1}\supset K_{t_2}$ if $t_1\le t_2$.
Show that there is a family of short maps $\phi_t\:\EE^m\to K_t$ 
such that $\phi_t|_{K_t}=\id$ for any $t$ and $\phi_{t_2}\circ\phi_{t_1}=\phi_{t_2}$ if $t_1\le t_2$.
\end{thm}

\section{Gradient-like curves}\label{sec:gradient-like}

Gradient-like curves will be used later in the construction of gradient curves.
The latter are a special reparametrization of gradient-like curves.

\begin{thm}{Definition}\label{def:grad-like-curve}{\sloppy 
Let $\spc{Z}$ be a complete length space
and $f\:\spc{Z}\subto\RR$ be locally Lipschitz semiconcave subfunction.
Suppose that $\spc{Z}$ is either $\Alex{}$ or $\CAT{}$.

}

A Lipschitz curve $\hat\alpha\:[s_{\min},s_{\max})\to\Dom f$ will be called an  \index{gradient-like curve}\emph{$f$-gradient-like curve} if
\[\hat\alpha^+=\tfrac{1}{|\nabla_{\hat\alpha} f|}\cdot\nabla_{\hat\alpha} f;\]
that is, for any $s\in[s_{\min},s_{\max})$, the right derivative $\hat\alpha^+(s)$ is defined and
\[\hat\alpha^+(s)=\tfrac{1}{|\nabla_{\hat\alpha(s)} f|}\cdot\nabla_{\hat\alpha(s)} f.\]

\end{thm}

Note that this definition implies that $|\nabla_p f|>0$ for any point $p$ on $\hat\alpha$.

The following theorem  gives a seemingly weaker condition that is equivalent to the definition of gradient-like curve.

\begin{thm}{Theorem}\label{thm:grad-like-2nd-def}
Suppose $\spc{Z}$ is a complete length space, 
$f\:\spc{Z}\subto\RR$ is a locally Lipschitz semiconcave subfunction,
and 
$|\nabla_p f|>0$ for any $p\in\Dom f$.
Assume that $\spc{Z}$ is either $\Alex{}$ or $\CAT{}$.

A curve $\hat\alpha\:[s_{\min},s_{\max})\to\Dom f$ is an $f$-gradient-like curve if and only if it is $1$-Lipschitz and
\[\liminf_{s\to s_0+}\frac{f\circ\hat\alpha(s)-f\circ\hat\alpha(s_0)}{s-s_0}
\ge 
|\nabla_{\hat\alpha(s_0)} f|
\eqlbl{eq:thm:grad-like-2nd-def-1}\]
for almost all $s_0\in [s_{\min},s_{\max})$.
\end{thm}

\parit{Proof.} The only-if part follows directly from the definition.
To prove the if part, note that for any $s_0\in[s_{\min},s_{\max})$ we have
\begin{align*}
\liminf_{s\to s_0+}\frac{f\circ\hat\alpha(s)-f\circ\hat\alpha(s_0)}{s-s_0}
&\ge 
\liminf_{s\to s_0+}
\frac{1}{s-s_0}
\cdot
\int\limits_{s_0}^s|\nabla_{\hat\alpha(\under s)}f|\cdot\dd\under s
\ge
\\
&\ge 
|\nabla_{\hat\alpha(s_0)}f|;
\end{align*}
the first inequality follows from \ref{eq:thm:grad-like-2nd-def-1} 
and the second from lower semicontinuity of the function $x\mapsto|\nabla_x f|$, 
see \ref{cor:gradlim}.
From \ref{lem:alm-grad}, we have 
\[\hat\alpha^+(s_0)=\tfrac{1}{|\nabla_{\hat\alpha(s_0)} f|}\cdot\nabla_{\hat\alpha(s_0)} f.\]
Hence the result.
\qeds


Recall that second-order differential inequalities are understood in a barrier sense; see Section~\ref{sec:conv-real}.

\begin{thm}{Theorem} \label{thm:concave}
Let $\spc{Z}$ be a complete length space 
and
$f\:\spc{Z}\subto \RR$ be
locally Lipschitz and $\lambda$-concave. 
Suppose that $\spc{Z}$ is either $\Alex{}$ or $\CAT{}$.
Assume $\hat\alpha\:[0,s_{\max})\to\Dom f$ is an $f$-gradient-like curve.
Then 
\[(f\circ\hat\alpha)''\le\lambda\] 
everywhere on $[0,s_{\max})$.
\end{thm}

{\sloppy 
Closely related statements were proved independently by Uwe Mayer \cite[2.36]{mayer} and Shin-ichi Ohta \cite[5.7]{ohta}.

}

Before the proof, let us formulate and prove a corollary. 

\begin{thm}{Corollary}\label{cor:right-cont}
Let $\spc{Z}$ be a complete length space,
$f\:\spc{Z}\subto \RR$ be a locally Lipschitz and semiconcave function, 
and $\hat\alpha\:[0,s_{\max})\to\Dom f$ be an $f$-gradient-like curve.
Suppose that $\spc{Z}$ is either $\Alex{}$ or $\CAT{}$.
Then the function $s\mapsto |\nabla_{\hat\alpha(s)}f|$
is right-continuous; 
that is, for any $s_0\in [0,s_{\max})$ we have
\[|\nabla_{\hat\alpha(s_0)}f|
=
\lim_{s\to s_0+} |\nabla_{\hat\alpha(s)}f|.\]

\end{thm}

\parit{Proof.} Applying \ref{thm:concave} locally, we have that $f\circ\hat\alpha(s)$ is semiconcave.
The statement follows since 
\[(f\circ\hat\alpha)^+(s)
=
(\dd_p f)\left(\tfrac{1}{|\nabla_{\hat\alpha(s)}f|}\cdot\nabla_{\hat\alpha(s)}f\right)
=
|\nabla_{\hat\alpha(s)}f|.\]
\qedsf

\parit{Proof of \ref{thm:concave}.} For any $s>s_0$,
\begin{align*}
(f\circ\hat\alpha)^+(s_0)&=|\nabla_{\hat\alpha(s_0)}f|
\ge
\\
&\ge
(d_{\hat\alpha(s_0)}f)(\dir{\hat\alpha(s_0)}{\hat\alpha(s)})
\ge
\\
&\ge
\frac{f\circ\hat\alpha(s)-f\circ\hat\alpha(s_0)}{\dist{\hat\alpha(s)}{\hat\alpha(s_0)}{}}
-
\tfrac\lambda2\cdot\dist[{{}}]{\hat\alpha(s)}{\hat\alpha(s_0)}{}.
\end{align*}
Let $\lambda_+=\max\{0,\lambda\}$. 
Since $s-s_0\ge\dist{\hat\alpha(s)}{\hat\alpha(s_0)}{}$, for any $s>s_0$ we have 
\[(f\circ\hat\alpha)^+(s_0)\ge
\frac{f\circ\hat\alpha(s)-f\circ\hat\alpha(s_0)}{s-s_0}-\tfrac{\lambda_+}2\cdot(s-s_0).
\eqlbl{eq:thm:concave-1}\]
Thus $f\circ\hat\alpha$ is $\lambda_+$-concave.
That finishes the proof for $\lambda\ge 0$.
For $\lambda<0$ we get only that $f\circ\hat\alpha$ is $0$-concave.

Note that $\dist{\hat\alpha(s)}{\hat\alpha(s_0)}{}=s-s_0-o(s-s_0)$. Thus
\[(f\circ\hat\alpha)^+(s_0)\ge
\frac{f\circ\hat\alpha(s)-f\circ\hat\alpha(s_0)}{s-s_0} -\tfrac\lambda2\cdot(s-s_0)+o(s-s_0).
\eqlbl{eq:thm:concave-2}\]
Together, \ref{eq:thm:concave-1} and \ref{eq:thm:concave-2} imply that $f\circ\hat\alpha$ is $\lambda$-concave.
\qeds

\begin{thm}{Proposition}
\label{prop:grad-like-unique-past}
Let $\spc{L}$ be a complete length $\Alex{\kappa}$ space, $p,q\in\spc{L}$.
Assume $\hat\alpha\:[s_{\min},s_{\max})\to\spc{L}$ is a $\distfun{p}$-gradient-like curve such that $\hat\alpha(s)\to z\in\mathopen{]}p q\mathclose{[}$ as $s\to s_{\max}+$.
Then $\hat\alpha$ is a unit-speed geodesic that lies in $[p q]$.
\end{thm}

\parit{Proof.} 
Clearly,
\[ \tfrac{d^+}{dt}\dist[{{}}]{q}{\hat\alpha(t)}{}
\ge
-1.
\eqlbl{eq:>=-1}
\]
On the other hand,

\[\begin{aligned}
\tfrac{d^+}{dt}\dist[{{}}]{p}{\hat\alpha(t)}{}
&\ge
(\dd_{\hat\alpha(t)}\distfun{p}{}{})(\dir{\hat\alpha(t)}{q})
\ge\\
&\ge
-\cos\angk\kappa{\hat\alpha(t)}p q.
\end{aligned}
\eqlbl{eq:>=-cos}\]
Inequalities \ref{eq:>=-1} and \ref{eq:>=-cos} imply that the function $t\mapsto\angk\kappa q {\hat\alpha(t)}p $ is nondecreasing.
Hence the result.
\qeds

\section{Gradient curves}\label{sec:grad-curves:def}

In this section we define gradient curves 
and tie them tightly to gradient-like curves 
which were introduced in Section~\ref{sec:gradient-like}.

\begin{thm}{Definition}\label{def:grad-curve}{\sloppy 
Let $\spc{Z}$ be a complete length space
and $f\:\spc{Z}\subto\RR$ be a locally Lipschitz and semiconcave subfunction.
Suppose that $\spc{Z}$ is either $\Alex{}$ or $\CAT{}$.

}

A locally Lipschitz curve $\alpha\:[t_{\min},t_{\max})\to\Dom f$ will be called an \index{gradient curve}\emph{$f$-gradient curve} if
\[\alpha^+=\nabla_{\alpha} f;\]
that is, for any $t\in[t_{\min},t_{\max})$, $\alpha^+(t)$ is defined and 
$\alpha^+(t)=\nabla_{\alpha(t)} f$.
\end{thm}

The following exercise describes a global geometric property of a gradient curve without direct reference to its function.
It uses the notion of \textit{self-contracting curves} introduced by Aris Daniilidis, Olivier Ley, St\'ephane Sabourau \cite{daniilidis-ley-sabourau}.

\begin{thm}{Exercise}\label{ex:elf-contracting}
Let 
$\spc{Z}$ be a complete length space,
$f\:\spc{Z}\to\RR$  a concave locally Lipschitz function, 
and $\alpha\:\II\to\spc{Z}$  an $f$-gradient curve.
Suppose that $\spc{Z}$ is either $\Alex{}$ or $\CAT{}$.

Show that $\alpha$ is \index{self-contracting curve}\emph{self-contracting}; that is,
\[t_1\le t_2\le t_3
\quad\Longrightarrow\quad
\dist{\alpha(t_1)}{\alpha(t_3)}{\spc{Z}}\ge \dist{\alpha(t_2)}{\alpha(t_3)}{\spc{Z}}.\]
\end{thm}

The next lemma states that gradient and gradient-like curves are special reparametrizations of each other.

\begin{thm}{Lemma}\label{lem:grad--grad-like}
Let $\spc{Z}$ be a complete length space
and
$f\:\spc{Z}\subto\RR$ be a locally Lipschitz semiconcave subfunction 
such that $|\nabla_p f|>0$ for any $p\in\Dom f$.
Suppose that $\spc{Z}$ is either $\Alex{}$ or $\CAT{}$.

Assume that $\alpha\:[0,t_{\max})\to \Dom f$ is a locally Lipschitz curve 
and $\hat\alpha\:[0,s_{\max})\to \Dom f$ is its reparametrization by arc-length, 
so $\alpha\z=\hat\alpha\circ\varsigma$ for a homeomorphism $\varsigma\:[0,t_{\max})\to [0,s_{\max})$.
Then 
\begin{align*}
\alpha^+&=\nabla_\alpha f
\\
&\Updownarrow
\\
\hat\alpha^+=\frac{1}{|\nabla_{\hat\alpha} f|}\cdot\nabla_{\hat\alpha} f
\quad
&
\text{and}
\quad
\varsigma^{-1}(s)
=
\int\limits_0^s\frac{\dd\under s}{(f\circ\hat\alpha)'(\under s)
 }.
\end{align*}

\end{thm}

\parit{Proof; $(\Downarrow)$.} 
According to \ref{thm:speed},
\[
\begin{aligned}
\varsigma'(t)&\ae|\alpha^+(t)|=
\\
&=|\nabla_{\alpha(t)}f|.
\end{aligned}
\eqlbl{eq:lem:grad--grad-like-1}\]
Note that 
\begin{align*}
(f\circ\alpha)'(t)&\ae (f\circ\alpha)^+(t)=
\\
&=|\nabla_{\alpha(t)} f|^2.
\end{align*}
Setting $s=\varsigma(t)$, we have
\begin{align*}(f\circ\hat\alpha)'(s)
&\ae\frac{(f\circ\alpha)'(t)}{\varsigma'(t)}
\ae
\\
&\ae|\nabla_{\alpha(t)}f|=
\\
&=|\nabla_{\hat\alpha(s)}f|.
\end{align*}

From \ref{thm:grad-like-2nd-def}, it follows that $\hat\alpha(t)$ is an $f$-gradient-like curve; 
that is,
\[\hat\alpha^+=\frac{1}{|\nabla_{\hat\alpha} f|}\cdot\nabla_{\hat\alpha} f.\]
In particular, $(f\circ\hat\alpha)^+(s)=|\nabla_{\hat\alpha^+(s)} f|$, and by \ref{eq:lem:grad--grad-like-1},
\begin{align*}\varsigma^{-1}(s)
&=\int\limits_0^s\frac{\dd\under s}{|\nabla_{\hat\alpha(\under s)} f|}
=
\\
&=
\int\limits_0^s\frac{\dd\under s}{(f\circ\hat\alpha)'(\under s)}.
\end{align*}
\medskip

\noi{$(\Uparrow)$.}
Clearly,
\begin{align*}\varsigma(t)
&=
\int\limits_0^{t}(f\circ\hat\alpha)^+(\varsigma(\under t))\cdot\dd \under t
=
\\
&=
\int\limits_0^{t}|\nabla_{\alpha(\under t)}f|\cdot\dd \under t.
\end{align*}
According to \ref{cor:right-cont}, the function $s\mapsto|\nabla_{\hat\alpha(s)}f|$ is right-continuous.
Therefore so is the function $t\mapsto|\nabla_{\hat\alpha\circ\varsigma(t)}f|=|\nabla_{\alpha(t)}f|$.
Hence, for any $t_0\in[0,t_{\max})$ we have
\begin{align*}\varsigma^+(t_0)
&=
\lim_{t\to t_0+}
\frac1{t-t_0}\cdot\int\limits_{t_0}^t
|\nabla_{\alpha(\under t)}f|\cdot\dd\under t
=
\\
&=
|\nabla_{\alpha(t_0)}f|.
\end{align*}
Thus, we have 
\begin{align*}\alpha^+(t_0)
&=
\varsigma^+(t_0)\cdot\hat\alpha^+(\varsigma(t_0))
=
\\
&=
\nabla_{\alpha(t_0)} f.
\end{align*}
\qedsf

\begin{thm}{Exercise}\label{ex:grad-curve-condition}
Let $\spc{Z}$ be a complete length space, and 
$f\:\spc{Z}\to \RR$ be a semiconcave locally Lipschitz 
function.
Suppose that $\spc{Z}$ is either $\Alex{}$ or $\CAT{}$.
Assume $\alpha\:\II\to \spc{Z}$ is a Lipschitz curve such that 
\begin{align*}
\alpha^+(t)&\le|\nabla_{\alpha(t)}f|,
\\
(f\circ\alpha)^+(t)&\ge |\nabla_{\alpha(t)}f|^2
\end{align*}
for almost all $t$.
Show that $\alpha$ is an $f$-gradient curve.
\end{thm}

\begin{thm}{Exercise}\label{ex:grad-curve-analitic}
Let $\spc{Z}$ be a complete length space and $f\:\spc{Z}\to\RR$ be a concave locally Lipschitz function.
Suppose that $\spc{Z}$ is either $\Alex{}$ or $\CAT{}$.
Show that $\alpha\:\RR\to\spc{Z}$ is an $f$-gradient curve if and only if
\[\dist[2]{x}{\alpha(t_1)}{\spc{Z}}-\dist[2]{x}{\alpha(t_0)}{\spc{Z}}
\le 
2\cdot(t_1-t_0)\cdot  (f\circ\alpha(t_1)-f(x))\]
for any $t_1>t_0$ and $x\in\spc{Z}$.
\end{thm}

\section{Distance estimates}\label{sec:grad-curv:dist-est}

\begin{thm}{First distance estimate}\label{thm:dist-est}
Let $\spc{Z}$ be a complete length space, and 
$f\:\spc{Z}\to \RR$ be a locally Lipschitz 
 $\lambda$-concave function.
Suppose that $\spc{Z}$ is either $\Alex{}$ or $\CAT{}$.
Let $\alpha,\beta\:[0,t_{\max})\to \spc{Z}$ be two $f$-gradient curves.
Then
\[\dist{\alpha(t)}{\beta(t)}{}
\le 
e^{\lambda\cdot t}\cdot\dist[{{}}]{\alpha(0)}{\beta(0)}{}\]
for any $t$.

Moreover, the statement holds for a locally Lipschitz $\lambda$-concave subfunction $f\:\spc{Z}\subto \RR$ if  there is a geodesic $[\alpha(t)\,\beta(t)]$ in $\Dom f$ for any~$t$.
\end{thm}

\parit{Proof.} 
If $\spc{Z}$ is not geodesic, then pass to its ultrapower $\spc{Z}^\o$.

Fix a choice of geodesic $[\alpha(t)\,\beta(t)]$ for each $t$.

Setting $\ell(t)=\dist{\alpha(t)}{\beta(t)}{}$, from the first variation inequality (\ref{lem:first-var}) and the estimate in \ref{cor:grad-lip} we get
\[\ell^+(t)\le-\<\dir{\alpha(t)}{\beta(t)},\nabla_{\alpha(t)}f\>-\<\dir{\beta(t)}{\alpha(t)},\nabla_{\beta(t)}f\>\le \lambda\cdot\ell(t).\]
Here one has to apply the first variation inequality for distance to the midpoint $m$ of $[\alpha(t)\,\beta(t)]$, and apply the triangle inequality.
Hence the result. 
\qeds

\begin{thm}{Second distance estimate}\label{lem:fg-dist-est}
Let $\spc{Z}$ be a complete length space, 
$\eps>0$,  
and $f,g\:\spc{Z}\to \RR$ be two $\lambda$-concave locally Lipschitz functions such that $|f-g|<\eps$.
Suppose that $\spc{Z}$ is either $\Alex{}$ or $\CAT{}$.
Assume
$\alpha,\beta\:[0,t_{\max})\to \spc{Z}$ are respectively $f$- and $g$-gradient curves.
Let $\ell\:t\mapsto\dist{\alpha(t)}{\beta(t)}{}$.
Then 
\[\ell^+\le \lambda\cdot\ell+\tfrac{2\cdot\eps}{\ell}.\]
In particular, if  $\alpha(0)=\beta(0)$ and $t_{\max}<\infty$ then
\[\dist{\alpha(t)}{\beta(t)}{}
\le
\Const\cdot\sqrt{\eps\cdot t}\]
for a constant $\Const=\Const(t_{\max},\lambda)$.

Moreover, the same conclusion holds for locally Lipschitz  $\lambda$-concave subfunctions $f,g\:\spc{Z}\subto \RR$ if for any $t\in[0,t_{\max})$ there is a geodesic $[\alpha(t)\,\beta(t)]$ in $\Dom f\cap\Dom g$.
\end{thm}

\parit{Proof.} Set $\ell=\ell(t)=\dist{\alpha(t)}{\beta(t)}{}$.
Fix $t$, and let $p=\alpha(t)$ and $q=\beta(t)$.
From the first variation formula and \ref{lem:grad-lip},
\begin{align*}
 \ell^+
&\le -\<\dir{p}{q},\nabla_{p}f\>
-\<\dir{q}{p},\nabla_{q}g\>
\le
\\
&\le -{\left({f(q)}-{f(p)}-\lambda\cdot\tfrac{\ell^2}2\right)}/{\ell}
-{\left({g(p)}-{g(q)}-\lambda\cdot\tfrac{\ell^2}2\right)}/{\ell}\le
\\
&\le \lambda\cdot\ell+\tfrac{2\cdot\eps}{\ell}.
\end{align*}
By integrating, we get the second statement.
\qeds

\section{Existence and uniqueness}
\label{sec:grad-curv:exist}

In general, the ``past'' of gradient curves can not be determined by the ``present''.
For example, consider the concave function $f\:\RR\to\RR$, $f(x)\z=-|x|$. The 
two curves $\alpha(t)=\min\{0,t\}$ with $\beta(t)=0$
are $f$-gradient with $\alpha(t)\z=\beta(t)\z=0$ for all $t\ge0$; 
however $\alpha(t)\z\ne\beta(t)$ for all $t<0$.
Another example can be given as follows.

\begin{wrapfigure}[8]{r}{34 mm}
\vskip-0mm
\centering
\includegraphics{mppics/pic-1215}
\vskip0mm
\end{wrapfigure}

\begin{thm}{Example}
Let $f$ be as in \ref{l-inf-grad};
that is, $f\:(x,y)\mapsto-|x|-|y|$ be the concave function on the $(x,y)$-plane;
its gradient field is sketched on the figure.

Let $\alpha$ be an $f$-gradient curve that starts at $p=(x,y)$ for $x>y>0$.
Then 
\[\alpha(t)=
\begin{cases}
(x-t,y-t) &\text{for}\quad 0\le t\le  x-y,
\\
(x-t,0) &\text{for}\quad x-y\le t\le  x,
\\
(0,0) &\text{for}\quad x\le t.
\end{cases}
\]
In particular, gradient curves can merge even in the region where $|\nabla f|\z\ne 0$. 
Hence their \textit{past} cannot be uniquely determined from their \textit{present}.
\end{thm}

The next theorem shows that the future gradient curve is determined by its present.

\begin{thm}{Picard's theorem}\label{thm:picard}
Let $\spc{Z}$ be a complete length space,
$f\:\spc{Z}\subto \RR$ be a semiconcave subfunction.
Suppose that $\spc{Z}$ is either $\Alex{}$ or $\CAT{}$.
Assume $\alpha,\beta\:[0,t_{\max})\to\Dom f$ are two $f$-gradient curves 
such that $\alpha(0)=\beta(0)$.
Then $\alpha(t)=\beta(t)$ for any $t\in[0,t_{\max})$.
\end{thm}

\parit{Proof.}
Follows from the first distance estimate (\ref{thm:dist-est}).
\qeds

\begin{thm}{Local existence}\label{thm:exist-grad-curv}
Let $\spc{Z}$ be a complete length space 
and $f\:\spc{Z}\subto \RR$ be locally Lipschitz $\lambda$-concave subfunction.
Suppose that $\spc{Z}$ is either $\Alex{}$ or $\CAT{}$.
Then for any $p\in \Dom f$,
\begin{subthm}{}
if $|\nabla_pf|>0$, then for some $\eps>0$, 
there is an $f$-gradient-like curve $\hat\alpha\:[0,\eps)\to\spc{Z}$ that starts at $p$ (that is, $\hat\alpha(0)\z=p$);
\end{subthm}

\begin{subthm}{}for some $\delta>0$, there is an $f$-gradient curve $\alpha\:[0,\delta)\to \spc{Z}$ that starts at $p$ (that is $\alpha(0)=p$).
\end{subthm}
\end{thm}

This theorem was proved by Grigory Perelman and the third author \cite{perelman-petrunin:qg};
we present a simplified proof given by Alexander Lytchak \cite{lytchak:open-map}.

\parit{Proof.} 
If $|\nabla_p f|=0$, then the constant curve $\alpha(t)=p$ is $f$-gradient.

Otherwise, choose $\eps>0$ 
such that $\oBall(p,\eps)\subset\Dom f$,
the restriction $f|_{\oBall(p,\eps)}$ is Lipschitz, 
and $|\nabla_x f|>\eps$ for all $x\in \oBall(p,\eps)$;
the latter is possible due to semicontinuity of \textbar gradient\textbar\ (\ref{cor:gradlim}).

The curves $\hat\alpha$ and $\alpha$ will be constructed in the following three steps.
First we construct an $f^\o$-gradient-like curve $\hat\alpha_\o\:[0,\eps)\to\spc{Z}^\o$ as an $\o$-limit of a certain sequence of broken geodesics in $\spc{Z}$.
Second, we parametrize $\hat\alpha_\o$ as in \ref{lem:grad--grad-like}, to obtain an $f^\o$-gradient curve $\alpha_\o$ in $\spc{Z}^\o$.
Third, applying Picard's theorem (\ref{thm:picard}) together with Lemma~\ref{lem:X-X^w}, we obtain that $\alpha_\o$ lies in $\spc{Z}\subset \spc{Z}^\o$ and therefore one can take $\alpha=\alpha_\o$ and $\hat\alpha=\hat\alpha_\o$.

Note that if $\spc{Z}$ is proper, then $\spc{Z}$ is a metric component of $\spc{Z}^\o$ and $f=f^\o|_{\spc{Z}}$.
Thus, in this case, the third step is not necessary.

\parit{Step 1.}
Given $n\in \NN$, 
by an open-closed argument,
we can construct a unit-speed curve $\hat\alpha_n\:[0,\eps] \to \spc{Z}$ starting at $p$, with a partition of $[0,\eps)$ into a countable number of half-open intervals $[\varsigma_i,\bar\varsigma_i)$ 
such that for each $i$ we have 
\begin{enumerate}[(i)]
\item $\hat\alpha_n([\varsigma_i,\bar\varsigma_i])$ is a geodesic and $\bar\varsigma_i-\varsigma_i<\tfrac{1}{n}$,
\item\label{alm-grad} 
$f\circ\hat\alpha_n(\bar\varsigma_i)-f\circ\hat\alpha_n(\varsigma_i)
>
(\bar\varsigma_i-\varsigma_i)
\cdot
(|\nabla_{\hat\alpha_n(\varsigma_i)}f|-\tfrac{1}{n}).$
\end{enumerate}

Passing to a subsequence of $\hat\alpha_n$ such that $f\circ\hat\alpha_n$ uniformly converges, let 
\[h(s)=\lim_{n\to\infty}f\circ\hat\alpha_n(s).\]

Let $\hat\alpha_\o=\lim_{n\to\o}\hat\alpha_{n}$;
it is a curve in $\spc{Z}^\o$ that starts at $p\in \spc{Z}\subset \spc{Z}^\o$.

Clearly $\hat\alpha_\o$ is $1$-Lipschitz.
From (\ref{alm-grad}) and \ref{lem:gradcon}, we have
\[(f^\o\circ\hat\alpha_\o)^+(\varsigma)
\ge
|\nabla_{\hat\alpha_\o(\varsigma)}f^\o|.\]
According to \ref{thm:grad-like-2nd-def}, $\hat\alpha_\o\:[0,\eps)\to \spc{Z}^\o$  is an $f^\o$-gradient-like curve.

\parit{Step 2.}
Clearly $h(s)=f^\o\circ\alpha_\o$. 
Therefore, according to \ref{thm:concave}, $h$ is $\lambda$-concave.
Thus we can define a homeomorphism $\varsigma\:[0,\delta]\to[0,\eps]$ by 
\[{\varsigma^{-1}(s)}
=
\int\limits_0^{s}\frac{\dd\under s}{h'(\under s)},
\eqlbl{eq:thm:exist-grad-curv-1}\]

According to \ref{lem:grad--grad-like}, $\alpha(t)=\hat\alpha\circ\varsigma(t)$ is an $f^\o$-gradient curve in $\spc{Z}^\o$. 

\parit{Step 3.}
Clearly, $\nabla_p f=\nabla_p f^\o$ for any $p\in \spc{Z}\subset \spc{Z}^\o$;
more formally, if $\iota\:\spc{Z}\hookrightarrow\spc{Z}^\o$ is the natural embedding, then
$(\dd_p\iota)(\nabla_p f)=\nabla_p f^\o$.
Thus it is sufficient to show that $\alpha_\o$ lies in $\spc{Z}$.
Assume the contrary; then according to \ref{lem:X-X^w}, there is a subsequence $\hat\alpha_{n_\kay}$ such that
\[\hat\alpha_\o\not
=
\hat\alpha'_\o
\df
\lim_{\kay\to\o}\hat\alpha_{n_\kay}.\]
Clearly $h(s)=f^\o\circ\hat\alpha_\o=f^\o\circ\hat\alpha'_\o$.
Thus for $\varsigma\:[0,\delta]\to[0,\eps]$ defined by \ref{eq:thm:exist-grad-curv-1}, 
we have that both curves
$\hat\alpha_\o\circ\varsigma$ and $\hat\alpha'_\o\circ\varsigma$ are $f^\o$-gradient.
From Picard's theorem (\ref{thm:picard}), we have $\hat\alpha_\o\circ\varsigma=\hat\alpha'_\o\circ\varsigma$.
Therefore $\hat\alpha_\o=\hat\alpha'_\o$, a contradiction.
\qeds

\section{Convergence}

\begin{thm}{Ultralimit of gradient curves}\label{ultr-lim-g-curve}
Assume
\begin{itemize}
\item $\spc{Z}_n$ is a sequence of complete spaces, $\spc{Z}_n \to \spc{Z}_\o$ as $n\to\o$, and $p_n\to p_\omega$ for a sequence of points $p_n\in \spc{Z}_n$,
\item all spaces $\spc{Z}_n$ are either $\Alex\kappa$ or $\CAT\kappa$, 
\item $f_n\:\spc{Z}_n\subto \RR$ are $\Lip$-Lipschitz and $\lambda$-concave,
$f_n\to f_\o$ as $n\to\o$, and $p_\o\in\Dom f_\o$.
\end{itemize}

Then: 

\begin{subthm}{thm:convex-limit-cbb}
$f_\o$ is $\lambda$-concave.
\end{subthm}

\begin{subthm}{lim-grad-like}
If $|\nabla_{p_\o}f_\o|>0$, then there is $\eps>0$ such that, the $f_n$-gradient-like curves $\hat\alpha_n\:[0,\eps)\to\spc{Z}_n$ are defined for $\o$-almost all $n$.
Moreover, a curve $\hat\alpha_\o\:[0,\eps)\to\spc{Z}_\o$ is a gradient-like curve that starts at $p_\o$ if and only if
$\hat\alpha_n(s)\to\hat\alpha_\o(s)$ as $n\to\o$ for all $s\in[0,\eps)$.
\end{subthm}

\begin{subthm}{lim-grad}
For some $\delta>0$, the $f_n$-gradient curves $\alpha_n\:[0,\delta)\to\spc{Z}_n$ are defined for $\o$-almost all $n$.
Moreover, a curve $\alpha_\o\:[0,\delta)\to\spc{Z}_\o$ is a gradient curve that starts at $p_\o$ if and only if
$\alpha_n(t)\to\alpha_\o(t)$  as $n\to\o$ for all $t\in[0,\delta)$.
\end{subthm}
\end{thm}


Note that according to Exercise~\ref{ex:nonconvex-limit}, part \ref{SHORT.thm:convex-limit-cbb} does not hold for general metric spaces.
The idea of the proof is the same as in the proof of local existence (\ref{thm:exist-grad-curv}).

\parit{Proof of \ref{ultr-lim-g-curve}; \ref{SHORT.thm:convex-limit-cbb}.}
Fix a geodesic $\gamma_\o\:\II\to \Dom f_\o$;
we need to show that the function 
\[t\mapsto f_\o\circ\gamma_\o(t)-\tfrac\lambda 2\cdot t^2\eqlbl{eq:lambda-concave}\]
is concave.

Since the $f_n$ are $\Lip$-Lipschitz, so is $f_\o$.
Therefore it is sufficient to prove concavity in the interior of $\II$.
In particular, we can assume that $\gamma_\o$ is sufficiently short and can be extended behind its ends $p_\o$ and $q_\o$ as a minimizing geodesic.
If $\spc{Z}$ is $\Alex{}$, then by Theorem~\ref{thm:almost.geod}, $\gamma_\o$ is the unique geodesic connecting $p_\o$ to $q_\o$.
The same holds true if $\spc{Z}$ is $\CAT{}$ by the uniqueness of geodesics (\ref{thm:cat-unique}).

Construct two sequences of points $p_n,q_n\in\spc{Z}_n$ such that $p_n\to p_\o$ and $q_n\to q_\o$ as $n\to \o$.
Applying either \ref{thm:almost.geod} or \ref{thm:cat-unique},
we can assume that for each $n$ there is a geodesic $\gamma_n$ from $p_n$ to $q_n$ in $\spc{Z}_n$.

Since $f_n$ is $\lambda$-concave, the function 
\[t\mapsto f_n\circ\gamma_n(t)-\tfrac\lambda 2\cdot t^2\]
is concave.

The $\o$-limit of the sequence $\gamma_n$ is a geodesic in $\spc{Z}_\o$ from $p_\o$ to $q_\o$.
By uniqueness of such geodesics, we have that $\gamma_n\to \gamma_\o$ as $n\to \o$.
Passing to the limit, we have \ref{eq:lambda-concave}.

\parit{If part of \ref{SHORT.lim-grad-like}.}
Take $\eps>0$ so small that $\oBall(p_\o,\eps)\subset\Dom f_\o$ and $|\nabla_{x_\o}f_\o|\z>0$ for any $x_\o\in\oBall(p_\o,\eps)$ (this is possible by \ref{cor:gradlim}).

Clearly $\hat\alpha_\o$ is $1$-Lipschitz.
From \ref{lem:gradcon}, we get 
\[(f_\o\circ\hat\alpha_\o)^+(s)
\ge
|\nabla_{\hat\alpha_\o(s)}f^\o|.\]
According to \ref{thm:grad-like-2nd-def}, $\hat\alpha_\o\:[0,\eps)\to \spc{Z}^\o$  is an $f_\o$-gradient-like curve.

\parit{If part of \ref{SHORT.lim-grad}.}
Assume first that $|\nabla_{p_\o}f_\o|>0$, 
so we can apply the if part of \ref{SHORT.lim-grad-like}.
Let $h_n=f_n\circ\hat\alpha_n\:[0,\eps)\to \RR$ 
and $h_\o=f_\o\circ\hat\alpha_\o$.
From \ref{thm:concave}, the $h_n$ are $\lambda$-concave, and clearly $h_n\to h_\o$ as $n\to\o$.
Let us define reparametrizations
\begin{align*}
{\varsigma_n^{-1}(s)}
&=
\int\limits_0^{s}\frac{\dd\under s}{h_n'(\under s)},
&
{\varsigma_\o^{-1}(s)}
&=
\int\limits_0^{s}\frac{\dd\under s}{h_\o'(\under s)}.
\end{align*}
The $\lambda$-convexity of the $h_n$ implies that $\sigma_n\to\sigma_\o$ as $n\to\o$.
By \ref{lem:grad--grad-like}, 
$\alpha_n=\hat\alpha_n\circ\varsigma_n$.
Applying the if part of \ref{SHORT.lim-grad-like} together with Lemma~\ref{lem:grad--grad-like},
we get that $\alpha_\o=\hat\alpha_\o\circ\varsigma_\o$ is gradient curve.

The remaining case $|\nabla_{p_\o}f_\o|=0$ can be reduced to the one above using the following trick.
Consider the sequence of spaces $\spc{Z}_n^{\times}=\spc{Z}_n\times\RR$,
with the sequence of subfunctions $f^{\times}_n\:\spc{Z}_n^{\times}\to\RR$ defined by
\[f^{\times}_n(p,t)=f_n(p)+t.\]
Applying either \ref{thm:warp-curv-bound:cbb:E} or \ref{thm:cbb-product},  we have that
$\spc{Z}_n^{\times}$ is a $\Alex{\kappa_-}$ space for $\kappa_-=\min\{\kappa,0\}$, or $\CAT{\kappa_+}$ space for $\kappa_+=\max\{\kappa,0\}$.
Note that the $f_n^{\times}$ are $\lambda_+$-concave
for $\lambda_+=\max\{\lambda,0\}$.
Now let $\spc{Z}_\o^{\times}=\spc{Z}_\o\times\RR$,
and $f^{\times}_\o(p,t)=f_\o(p)+t$.

Clearly 
$\spc{Z}_n^{\times}\to\spc{Z}_\o^{\times}$,
$f_n^{\times}\to f_\o^{\times}$ as $n\to\o$,
and $|\nabla_xf^{\times}_\o|>0$ for any $x\in\Dom f_\o^{\times}$.
Thus for the sequence $f_n^{\times}\:\spc{Z}_n^{\times}\subto\RR$, 
we can apply the if part of \ref{SHORT.lim-grad-like}.
It remains to note that the curve $\alpha^{\times}_\o(t)=(\alpha_\o(t),t)$ is an $f^{\times}_\o$-gradient curve in $\spc{Z}^{\times}_\o$ 
if and only if $\alpha_\o(t)$ is an $f_\o$-gradient curve.

\parit{Only-if part of \ref{SHORT.lim-grad} and \ref{SHORT.lim-grad-like}.}
The only-if part of \ref{SHORT.lim-grad} follows from
the if part of \ref{SHORT.lim-grad} and Picard's theorem (\ref{thm:picard}).
Applying Lemma~\ref{lem:grad--grad-like}, we get the only-if part of \ref{SHORT.lim-grad-like}.
\qeds

Directly from  local existence (\ref{thm:exist-grad-curv}) and the distance estimates (\ref{thm:dist-est}), we obtain the following.

\begin{thm}{Global existence}\label{thm:glob-exist-grad-curv}
Let $f\:\spc{Z}\subto \RR$ be a locally Lipschitz and $\lambda$-concave subfunction on a complete length space $\spc{Z}$.
Suppose that $\spc{Z}$ is either $\Alex{}$ or $\CAT{}$.
Then for any $p\in \Dom f$, there is $t_{\max}\in(0,\infty]$ such that
there is an $f$-gradient curve $\alpha\:[0,t_{\max})\to \spc{Z}$ with $\alpha(0)=p$.
Moreover, for any sequence $t_n\to t_{\max}-$, the sequence $\alpha(t_n)$ does not have a limit point in $\Dom f$.
\end{thm}

The following theorem guarantees the existence of gradient curves for all times for the special type of semiconcave functions that play important role in the theory.
It follows from \ref{thm:glob-exist-grad-curv},
\ref{thm:concave} and \ref{lem:grad--grad-like}.

\begin{thm}{Theorem}\label{thm:comp-grad-test}
Let $\spc{Z}$ be a complete length space 
and $f\:\spc{Z}\to\RR$ satisfies 
\[f''+\kappa\cdot f\le \lambda\] 
for real constants $\kappa$ and $\lambda$.
Suppose that $\spc{Z}$ is either $\Alex{}$ or $\CAT{}$.
Then $f$ has \emph{complete gradient};
that is, for any $x\in\spc{Z}$ there is a $f$-gradient curve $\alpha\:[0,\infty)\to\spc{Z}$ that starts at~$x$.
\end{thm}

\section{Gradient flow}\label{sec:Gradient flow}

In this section we define gradient flow for semiconcave subfunctions 
and reformulate theorems obtained earlier in this chapter using this new terminology.

Let $\spc{Z}$ be a complete length space 
and $f\:\spc{Z}\subto \RR$ be a locally Lipschitz semiconcave subfunction.
Suppose that $\spc{Z}$ is either $\Alex{}$ or $\CAT{}$.
For any $t\ge 0$, we write $\GF^t_f(x)=y$ if there is an $f$-gradient curve $\alpha$ such that $\alpha(0)=x$ and $\alpha(t)=y$.
The partially defined map $\GF^t_f$ from $\spc{Z}$ to itself is called the \index{gradient flow}\emph{$f$-gradient flow} for time $t$.
 
From \ref{lem:fg-dist-est}, 
it follows that for any $t\ge 0$, the domain of definition of $\GF^t_f$ is an open subset of $\spc{Z}$; 
that is, $\GF^t_f$ is a submap.
Moreover, if $f$ is defined on all of $\spc{Z}$ and $f''+\Kappa\cdot f\le \lambda$ for constants $\Kappa,\lambda\in\RR$, 
then according to \ref{thm:comp-grad-test}, $\GF^t_f(x)$ is defined for all pairs $(x,t)\in\spc{Z}\times\RR_{\ge0}$.

Clearly $\GF^{t_1+t_2}_f=\GF_f^{t_1}\circ\GF_f^{t_2}$;
in other words, gradient flow is given by an action of the semigroup $(\RR_{\ge0},+)$.

From the first distance estimate (\ref{thm:dist-est}),
we have the following:

\begin{thm}{Proposition}\label{prop:GF-is-lip}
Let $\spc{Z}$ be a complete length $\Alex{}$ or $\CAT{}$ space 
and $f\:\spc{Z}\to \RR$ be a semiconcave function.
Then the map $x\mapsto\GF^t_f(x)$ is locally Lipschitz.

Moreover, if $f$ is $\lambda$-concave, then $\GF^t_f$ is $e^{\lambda\cdot t}$-Lipschitz.
\end{thm}

The next proposition states that gradient flow is stable under Gromov--Hausdorff convergence.
The proposition follows directly from the proposition on ultralimit of gradient curves~\ref{ultr-lim-g-curve}.

\begin{thm}{Proposition}\label{grad-curve-conv}
Supose $\spc{Z}_\infty,\spc{Z}_1,\spc{Z}_2,\dots$ are complete length $\Alex\kappa$ space, $\spc{Z}_n\xto{\GH} \spc{Z}_\infty$, and $f_n\:\spc{Z}_n\to\RR$ is a sequence of
$\lambda$-concave functions that converges to $f_\infty\:\spc{Z}_\infty\to \RR$. 
Then
$\GF_{f_n}^t\:\spc{Z}_n\to \spc{Z}_n$ converges to $\GF_{f_\infty}^t\:\spc{Z}_\infty\to \spc{Z}_\infty$.
\end{thm}


\section{Line splitting theorem}

Let $\spc{X}$ be a metric space and $A,B\subset \spc{X}$.
We will write 
\[\spc{X}=A\oplus B\]\index{$A\oplus B$}
if there are projections $\proj_A\:\spc{X}\to A$ 
and 
$\proj_B\:\spc{X}\to B$
such that 
\[\dist[2]{x}{y}{}=\dist[2]{\proj_A(x)}{\proj_A(y)}{}+\dist[2]{\proj_B(x)}{\proj_B(y)}{}\]
for any two points $x,y\in \spc{X}$.

Note that if 
\[\spc{X}=A\oplus B\]
then 
\begin{itemize}
\item $A$ intersects $B$ at a single point,
\item both sets $A$ and $B$ are convex sets in $\spc{X}$.
\end{itemize}

Recall that a line in a metric space is a both-sided infinite geodesic; thus it minimizes the length on each segment.

 {\sloppy 

\begin{thm}{Line splitting theorem}\label{thm:splitting}
Let $\spc{L}$  be a complete length $\Alex{0}$ space
and $\gamma$ be a line in $\spc{L}$. 
Then 
\[\spc{L}=\spc{L}'\oplus \gamma(\RR)\]
for a subset $\spc{L}'\subset \spc{L}$.
\end{thm}

}

For smooth $2$-dimensional surfaces, 
this theorem was proved by Stefan Cohn-Vossen \cite{cohn-vossen_line}.
For Riemannian manifolds of higher dimensions 
it was proved by Victor Toponogov \cite{toponogov-globalization+splitting}.
Then it was generalized by Anatoliy Milka \cite{milka-line}
to Alexandrov spaces; nearly the same proof is used in \cite[1.5]{burago-burago-ivanov}.

Further generalizations of the splitting theorem for Riemannian manifolds with nonnegative Ricci curvature were obtained by Jeff Cheeger and Detlef Gromoll \cite{cheeger-gromoll-split}.
This was further generalized by Jeff Cheeger and Toby Colding for limits of Riemannian manifolds with almost nonnegative Ricci curvature \cite{cheeger-colding-alm-rigidity} and to their synthetic generalizations, so-called {}\emph{RCD spaces}, by Nicola Gigli \cite{gigli2013splitting, gigli-splitting-overview}.
Jost-Hinrich Eschenburg obtained an analogous result for  Lorentzian manifolds \cite{eshenburg-split}, that is, pseudo-Riemannian manifolds of signature $(1,n)$.

We present a proof that uses gradient flow for Busemann functions. 
It is close in spirit to the proof given in \cite{cheeger-gromoll-split}.

Before going into the proof, let us state a few corollaries of the theorem.

\begin{thm}{Corollary}\label{cor:splitting}
Let $\spc{L}$ be a complete length $\Alex{0}$ space. 
Then there is an isometric splitting
\[
\spc{L}=\spc{L}'\oplus H
\]
where $H\subset \spc{L}$ is a subset isometric to a Hilbert space, and $\spc{L}'\subset \spc{L}$ is a convex subset that contains no line. 
\end{thm}

 {\sloppy 

\begin{thm}{Corollary}\label{cor:splitting-vectors}
Let $\spc{K}$ be a finite-dimensional complete length $\Alex0$ cone and $v_+,v_-\in \spc{K}$ be a pair of opposite vectors 
(that is, $v_+ + v_-=0$, see Definiton~\ref{def:opp+Lin}).
Then there is an isometry $\iota\:K\to K'\times \RR$ such that
$\iota:v_\pm\mapsto (\0',\pm|v_\pm|)$, where $K'$ is a complete length $\Alex0$ space having a cone structure with tip $\0'$.
\end{thm}

}

\begin{thm}{Corollary}\label{cor:splitting-CBB[1]}
Let $\spc{L}$ be an $m$-dimensional complete length $\Alex1$ space, $2\le m<\infty$, and $\rad\spc{L}=\pi$.
Then \[\spc{L}\iso \mathbb{S}^m.\]
 
\end{thm}

The following lemma is closely relevant to the first distance estimate (\ref{thm:dist-est}); its proof goes along the same lines.

\begin{thm}{Lemma}\label{lem:dist-estimate}
Let $\spc{L}$ be a complete length $\Alex{0}$ space.
Suppose $f\:\spc{L}\to\RR$ be a concave 1-Lipschitz function.
Consider two $f$-gradient curves $\alpha$ and $\beta$.
Then for any $t, s\ge 0$ we have
\begin{align*}
&\dist[2]{\alpha(s)}{\beta(t)}{}
\le 
\dist[2]{p}{q}{}+
2\cdot(f(p)-f(q))\cdot(s-t)+ (s-t)^2,
\end{align*}
where $p=\alpha(0)$ and $q=\beta(0)$.
\end{thm}

\parit{Proof.}
If $\spc{L}$ is not geodesic, then pass to its ultrapower $\spc{L}^\o$.

Since $f$ is 1-Lipschitz, $|\nabla f|\le1$.
Therefore 
\[f\circ\beta(t)\le f(q)+t\]
for any $t\ge0$.

Set $\ell(t)=\dist{p}{\beta(t)}{}$.
Applying \ref{lem:grad-lip:lam=0} and the first variation inequality (\ref{lem:first-var}), we get
\begin{align*}
\ell^2(t)^+
&\le 2\cdot \left(f\circ\beta(t)-f(p)\right)\le 
\\
&\le2\cdot\left(f(q)+t-f(p)\right).
\end{align*}
Therefore 
\[\ell^2(t)-\ell^2(0)\le 2\cdot\left(f(q)-f(p)\right)\cdot t + t^2.\]
It proves the needed inequality in case $s=0$.
Combining it with the first distance estimate (\ref{thm:dist-est}), we get the result in case $s\le t$.
The case $s\ge t$ follows by switching the roles of $s$ and $t$.
\qeds

\parit{Proof of \ref{thm:splitting}.} Consider two Busemann functions, $\bus_+$ and $\bus_-$, asociated with half-lines $\gamma:[0,\infty)\to \spc{L}$ and $\gamma:(-\infty,0]\to \spc{L}$ respectively; that is,
\[
\bus_\pm(x)
=
\lim_{t\to\infty}\dist{\gamma(\pm t)}{x}{}- t.
\]
According to Exercise~\ref{ex:busemann-CBB}, 
both functions $\bus_\pm$ are concave.

Fix $x\in \spc{L}$.
Note that since $\gamma$ is a line, we have 
\[\bus_+(x)+\bus_-(x)\ge0.\]

On the other hand, by \ref{comp-kappa}, 
$f(t)=\distfun[2]{x}{}{}\circ\gamma(t)$ 
is $2$-concave.
In particular, $f(t)\le t^2+at+b$ for some constants $a,b\in\RR$. 
Passing to the limit as $t\to\pm\infty$, we have \[\bus_+(x)+\bus_-(x)\le0.\]

Hence
\[
\bus_+(x)+\bus_-(x)= 0
\]
for any $x\in \spc{L}$.
In particular, the functions $\bus_\pm$ are \index{affine function}\emph{affine};
that is, they are convex and concave at the same time.

It follows that for any $x$,
\begin{align*}
|\nabla_x \bus_\pm|
&=\sup\set{\dd_x\bus_\pm(\xi)}{\xi\in\Sigma_x}=
\\
&=\sup\set{-\dd_x\bus_\mp(\xi)}{\xi\in\Sigma_x}\equiv
\\
&\equiv1.
\end{align*}
By Exercise~\ref{ex:grad-curve-condition}, a 
$1$-Lipschitz curve $\alpha$ such that $\bus_\pm(\alpha(t))=t+\Const$ is a $\bus_\pm$-gradient curve. 
In particular, $\alpha(t)$ is a $\bus_+$-gradient curve if and only if $\alpha(-t)$ is a $\bus_-$-gradient curve.
It follows that for any $t>0$, the $\bus_\pm$-gradient flows commute;
that is, 
\[\GF_{\bus_+}^t\circ\GF_{\bus_-}^t=\id_\spc{L}.\]
Setting
\[\GF^t=\left[\begin{matrix}
\GF_{\bus_+}^t&\hbox{if}\ t\ge0\\
\GF_{\bus_-}^t&\hbox{if}\ t\le0
               \end{matrix}\right.\]
defines an $\RR$-action on $\spc{L}$.

Consider the level set $\spc{L}'=\bus_+^{-1}(0)=\bus_-^{-1}(0)$;
it is a closed convex subset of $\spc{L}$, and therefore forms an Alexandrov space.
Consider the map $h\:\spc{L}'\times \RR\to \spc{L}$ defined by $h\:(x,t)\mapsto \GF^t(x)$.
Note that $h$ is onto.
Applying Lemma \ref{lem:dist-estimate} for $\GF_{\bus_+}^t$ and $\GF_{\bus_-}^t$ shows that $h$ is short and non-contracting at the same time; that is, $h$ is an isometry.
\qeds

\section{Radial curves}\label{sec:Radial curves: definition}

The radial curves are specially reparametrized gradient curves for distance functions.
This parametrization makes them behave like unit-speed geodesics in a natural comparison sense (\ref{sec:Radial comparisons}).

\begin{thm}{Definition}\label{def:rad-curv}
Assume $\spc{L}$ is a complete length $\Alex{}$ space, 
$\kappa\in\RR$, 
and $p\in \spc{L}$.
A curve 
$$\sigma\:[s_{\min},s_{\max})\to \spc{L}$$  
is called a 
\emph{$(p,\kappa)$-radial curve} 
if
$$s_{\min}
\z=
\dist{p}{\sigma(s_{\min})}{}\in(0,\tfrac{\varpi\kappa}2)$$ 
and $\sigma$ satisfies the differential equation
\[\sigma^+(s)
\z=
\frac{\tang\kappa\dist[{{}}]{p}{\sigma(s)}{}}{\tang\kappa s}
\cdot
\nabla_{\sigma(s)}\distfun{p}{}{}
\eqlbl{eq:rad}\]
for any $s\in[s_{\min},s_{\max})$, where $\tang\kappa x\df\frac{\sn\kappa x}{\cs\kappa x}$.

If $x=\sigma(s_{\min})$, we say that $\sigma$ {}\emph{starts} at  $x$.
\end{thm}

Note that according to the definition, $s_{\max}\le\tfrac{\varpi\kappa}2$.

In the the next section, we will see that  $(p,\kappa)$-radial curves 
work best for $\Alex\kappa$ spaces.

\begin{thm}{Definition}\label{def:rad-geod}
Let $\spc{L}$ be a complete length $\Alex{}$ space
and $p\in\spc{L}$.
A unit-speed geodesic  $\gamma\:\II\to \spc{L}$  is called a
\index{radial geodesic}\emph{$p$-radial geodesic} if 
$\dist{p}{\gamma(s)}{}\equiv s$.
\end{thm}

The proofs of the following two propositions follow directly from the definitions. 

\begin{thm}{Proposition}\label{prop:rad-geod}
Let $\spc{L}$ be a complete length $\Alex{}$ space
and $p\in\spc{L}$.
Assume $\tfrac{\varpi\kappa}{2}
\ge 
s_{\max}$.
Then any $p$-radial geodesic 
$\gamma\:[s_{\min},s_{\max})
\to 
\spc{L}$ 
is a $(p,\kappa)$-radial curve.
\end{thm}

\begin{thm}{Proposition}\label{prop:dist<s}
Suppose $\spc{L}$ is a complete length $\Alex{}$ space, 
$p\in\spc{L}$ 
and $\sigma\:[s_{\min},s_{\max})\to \spc{L}$ is a $(p,\kappa)$-radial curve.
Then for any $s\in [s_{\min},s_{\max})$, 
we have $\dist{p}{\sigma(s)}{}\le s$.

Moreover, 
if for some $s_0$ we have $\dist{p}{\sigma(s_0)}{}= s_0$, 
then the restriction $\sigma|_{[s_{\min},s_0]}$ is a $p$-radial geodesic.
\end{thm}

\begin{thm}{Existence and uniqueness}\label{rad-curv-exist}
Let $\spc{L}$ be a complete length $\Alex{}$ space, 
$\kappa\in\RR$, 
$p\in\spc{L}$, 
and $x\in \spc{L}$.
Assume
$0
<
\dist{p}{x}{}
<
\tfrac{\varpi\kappa}2$.
Then there is a unique $(p,\kappa)$-radial curve $\sigma\:[\dist{p}{x}{},\tfrac{\varpi\kappa}2)\to \spc{L}$ 
that starts at $x$;
that is, $\sigma(\dist{p}{x}{})=x$.
\end{thm}

\parit{Proof; existence.}
Let \index{$\itg\kappa$} 
\[\itg\kappa\:[0,\tfrac{\varpi\kappa}2)\to\RR,
\quad 
\itg\kappa (t)=\int\limits_0^t\tang\kappa\under t\cdot\dd\under t.\]
Clearly $\itg\kappa$ is smooth and increasing.
From \ref{prop:conv-comp} it follows that the composition 
\[f=\itg\kappa\circ\distfun{p}{}{}\] 
is semiconcave in $\oBall(p,\tfrac{\varpi\kappa}2)$.

According to \ref{thm:exist-grad-curv}, there is an $f$-gradient curve $\alpha\:[0,t_{\max})\to \spc{L}$ defined on the maximal interval such that $\alpha(0)=x$.

Now consider a solution $\tau(t)$ for the initial value problem $\tau'\z=(\tang\kappa\tau)^2$, $\tau(0)=r$. 
Note that $\tau(t)$ is also a gradient curve  for the function $\itg\kappa$ defined on $[0,\tfrac{\varpi\kappa}2)$.
Direct calculations show that the composition $\alpha\circ\tau^{-1}$ 
is a $(p,\kappa)$-radial curve.

\parit{Uniqueness.} Assume $\sigma^1,\sigma^2$ are two $(p,\kappa)$-radial curves that start at $x$.
Then the compositions $\sigma^i\circ\tau$ both give $f$-gradient curves.
By Picard's theorem (\ref{thm:picard}), we have
$\sigma^1\circ\tau 
\equiv 
\sigma^2\circ\tau$.
Therefore $\sigma^1(s)=\sigma^2(s)$ 
for any $s\ge r$ such that both sides are defined.
\qeds

\section{Radial comparisons}\label{sec:Radial comparisons}

In this section we show that radial curves behave in a comparison sense like unit-speed geodesics.
Recall that notation $\tangle\mc\kappa\{a;b,c\}$ is introduced in \ref{model}.

\begin{thm}{Radial monotonicity}\label{rad-mon}
Let $\spc{L}$ be a complete length $\Alex{\kappa}$ space and
$p, q$ be distinct points in $\spc{L}$.
Let $\sigma\:  [s_{\min},\tfrac{\varpi\kappa}2)\to \spc{L}$
be a $(p,\kappa)$-radial curve.
Then the function 
\[s\mapsto 
\tangle\mc\kappa\{
\dist{q}{\sigma(s)}{};
\dist{p}{q}{},
s
\}\]
is nonincreasing in its domain of definition.
\end{thm}

Radial monotonicity implies the following by straightforward calculations.

\begin{thm}{Corollary}\label{cor:rad-comp}
Let $\kappa\le0$,
$\spc{L}$ be a complete $\Alex\kappa$ space,
and $p, q\in \spc{L}$.
Let $\sigma\:[s_{\min},\infty)\to \spc{L}$ be a $(p,\kappa)$-radial curve.
Then for any $w\ge 1$, 
the function
\[
s\mapsto \tangle\mc\kappa\{\dist{q}{\sigma(s)}{};\dist{p}{q}{},w\cdot s\}
\]
is nonincreasing in its domain of definition.
\end{thm}

\begin{thm}{Radial comparison}\label{rad-comp}
Let $\spc{L}$ be a complete length $\Alex{\kappa}$ space 
and $p\in \spc{L}$.
Let $\rho\:  [r_{\min},\tfrac{\varpi\kappa}2)\to \spc{L}$
and    $\sigma\:[s_{\min},\tfrac{\varpi\kappa}2)\to \spc{L}$
be two $(p,\kappa)$-radial curves.
Let
\[\phi_{\min}=\angkk\kappa p{\rho(r_{\min})}{\sigma(s_{\min})}.
\]
Then for any $r\in[r_{\min},\tfrac{\varpi\kappa}2)$ and  $s\in[s_{\min},\tfrac{\varpi\kappa}2)$,
we have
\[
\tangle\mc\kappa\{\dist{\rho(r)}{\sigma(s)}{};r,s\}
\le \phi_{\min},
\]
or equivalently,
\[
\dist{\rho(r)}{\sigma(s)}{}
\le \side\kappa\{\phi_{\min};r,s\}.
\]

\end{thm}

We prove Theorems \ref{rad-mon} and \ref{rad-comp} simultaneously.
The proof is an application of \ref{lem:grad-lip} plus trigonometric manipulations.
We give a proof first in the simplest case $\kappa=0$,
and then in the harder case $\kappa\ne 0$.
The arguments for both cases are nearly the same, 
but the case $\kappa\ne 0$ requires an extra twist.

\parit{Proof of \ref{rad-mon} and \ref{rad-comp} in case $\kappa=0$.}
Set
\begin{align*}
R=R(r)&=\dist{p}{\rho(r)}{},
\\
S=S(s)&=\dist{p}{\sigma(s)}{},
\\
\ell=\ell(r,s)&=\dist{\rho(r)}{\sigma(s)}{},
\\
\phi=\phi(r,s)&=\tangle\mc0\{\ell(r,s);r,s\}.
\end{align*}

\begin{figure}[!ht]
\vskip-0mm
\centering
\includegraphics{mppics/pic-1505}
\vskip0mm
\end{figure}

It will be sufficient to prove the following inequalities:
\[\tfrac{\partial^+}{\partial r}\phi(s_{\min},r)\le 0,\qquad
\tfrac{\partial^+}{\partial s}\phi(s,r_{\min})\le 0\leqno(*)\mc0_\phi\]
\[
s\cdot\tfrac{\partial^+}{\partial s}\phi
+
r\cdot\tfrac{\partial^+}{\partial r}\phi\le 0.
\leqno(**)\mc0_\phi
\]

The radial monotonicity follows from $(*)\mc0_\phi$.
The radial comparison follows from  $(*)\mc0_\phi$ and $(**)\mc0_\phi$.
Indeed, one can connect $(s_{\min},r_{\min})$ and $(s_0,r_0)$ in $[s_{\min},\infty)\times[r_{\min},\infty)$ 
by a concatenation of a coordinate segment and a segment defined by $r/s=r_0/s_0$ as in the figure.
According to $(*)\mc0_\phi$ and $(**)\mc0_\phi$, $\phi$ does not increase while the  pair $(r,s)$ moves along this concatenation with nondecreasing $r$ and $s$.
Thus $\phi(r_0,s_0)\le\phi(r_{\min},s_{\min})=\phi_{\min}$.

\begin{figure}[!ht]
\vskip-0mm
\centering
\includegraphics{mppics/pic-1510}
\vskip0mm
\end{figure}

Let us rewrite the inequalities $(*)\mc0_\phi$ and $(**)\mc0_\phi$ in an equivalent form:
\[
\begin{aligned}
\tfrac{\partial^+}{\partial s}\ell(s,r_{\min})
&\le 
\cos\tangle\mc0\{r_{\min};s,\ell\},
\\
\tfrac{\partial^+}{\partial r}\ell(s_{\min},r)
&\le 
\cos\tangle\mc0\{s_{\min};r,\ell\},
\end{aligned}
\leqno(*)\mc0_\ell
\]

\[
s\cdot\tfrac{\partial^+}{\partial s}\ell
+
r\cdot\tfrac{\partial^+}{\partial r}\ell\le 
 s\cdot\cos\tangle\mc0\{r;s,\ell\}
+
r\cdot\cos\tangle\mc0\{s;r,\ell\}=\ell.
\leqno(**)\mc0_\ell
\]

Let 
\[f=\tfrac{1}{2}\cdot\distfun[2]{p}{}{}.\leqno(A)\mc0\] 
Clearly $f$ is $1$-concave, and
\[\rho^+(r)=\tfrac{1}{r}\cdot\nabla_{\rho(r)} f\quad \text{and}\quad \sigma^+(s)=\tfrac{1}{s}\cdot\nabla_{\sigma(s)} f.\leqno(B)\mc0\]
Thus from \ref{lem:grad-lip}, we have
\[\tfrac{\partial^+}{\partial r}\ell
=
-\tfrac{1}{r}\cdot\<\nabla_{\rho(r)} f,\dir{\rho(r)}{\sigma(s)}\>
\le\frac{{\ell^2}+{R^2}-{S^2}}{2\cdot\ell\cdot r}.\leqno(C)\mc0\]
Since $R(r)\le r$ and $S(s_{\min})=s_{\min}$, then
\[
\begin{aligned}
\tfrac{\partial^+}{\partial r}\ell(r,s_{\min})
&\le
\frac{{\ell^2}+r^2-s_{\min}^2}{2\cdot\ell\cdot r}
=\\
&=
\cos\tangle\mc0\{s_{\min};r,\ell\},
\end{aligned}
\leqno(D)\mc0
\]
which is the first inequality in $(*)\mc0_\ell$.
By switching $\rho$ and $\sigma$ we obtain the second inequality in $(*)\mc0_\ell$.
Further, adding $(C)\mc0$ and its mirror-inequality for $\frac{\partial^+}{\partial s}\ell$, we have
\[r\cdot\tfrac{\partial^+}{\partial r}\ell
+
s\cdot\tfrac{\partial^+}{\partial s}\ell
\le 
\frac{{\ell^2}+{R^2}-{S^2}}{2\cdot\ell }+\frac{{\ell^2}+{S^2}-{R^2}}{2\cdot\ell }
= 
\ell,
\leqno(E)\mc0\]
namely $(**)\mc0_\ell$.
\qeds

\parit{Proof of \ref{rad-mon} and \ref{rad-comp} in case $\kappa\ne 0$.} 
As before, let
\begin{align*}
R=R(r)&=\dist{p}{\rho(r)}{},&\ell&=\ell(r,s)=\dist{\rho(r)}{\sigma(s)}{},
\\
S=S(s)&=\dist{p}{\sigma(s)}{},&\phi&=\phi(r,s)=\tangle\mc\kappa\{\ell(r,s);r,s\}.
\end{align*}
It suffices to prove the following three inequalities:
\[
\begin{aligned}
&\tfrac{\partial^+}{\partial r}\phi(s_{\min},r)\le 0, 
&
&\tfrac{\partial^+}{\partial s}\phi(s,r_{\min})\le 0,
\end{aligned}
\leqno(*)\mc\pm_\phi
\]
\[
\sn\kappa s\cdot\cs\kappa S\cdot\tfrac{\partial^+}{\partial s}\phi
+
\sn\kappa r\cdot\cs\kappa R\cdot\tfrac{\partial^+}{\partial r}\phi\le 0.
\leqno(**)\mc\pm_\phi.
\]

Then radial monotonicity follows from $(*)\mc\pm_\phi$.
The radial comparison follows from $(*)\mc0_\phi$ and $(**)\mc\pm_\phi$.
Indeed, the functions $s\mapsto \sn\kappa s\cdot\cs\kappa S$ and $r\mapsto \sn\kappa r\cdot\cs\kappa R$ are Lipschitz.
Thus there is a solution for the differential equation
\[(r',s')=(\sn\kappa s\cdot\cs\kappa S,\sn\kappa r\cdot\cs\kappa R)\] 
with any initial data $(r_0,s_0)\in[r_{\min},\tfrac{\varpi\kappa}2)\times[s_{\min},\tfrac{\varpi\kappa}2)$.
(Unlike the case $\kappa=0$, the solution cannot be written explicitly.)
Since $\sn\kappa s\cdot\cs\kappa S$, $\sn\kappa r\cdot\cs\kappa R>0$, this solution $(r(t),s(t))$ must meet one of the coordinate rays
$\{r_{\min}\}\times[s_{\min},\tfrac{\varpi\kappa}2)$ or $[r_{\min},\tfrac{\varpi\kappa}2)\times\{s_{\min}\}$.
That is, one can connect the pair $(s_{\min},r_{\min})$ to $(s_0,r_0)$ by a concatenation of a coordinate segment (vertical or horizontal) and part of the solution $(r(t),s(t))$.
According to $(*)\mc\pm_\phi$ and $(**)\mc\pm_\phi$, the value of $\phi$ does not increase while the pair $(r,s)$ moves along this concatenation in direction of increasing $r$ and $s$.
Thus $\phi(r_0,s_0)\le\phi(r_{\min},s_{\min})$.

As before, we rewrite the inequalities $(*)\mc\pm_\phi$ and $(**)\mc\pm_\phi$ in terms of $\ell$:
\[
\begin{aligned}
\tfrac{\partial^+}{\partial s}\ell(s,r_{\min})
&\le 
\cos\tangle\mc\kappa\{r_{\min};s,\ell\},
\\
\tfrac{\partial^+}{\partial r}\ell(s_{\min},r)
&\le 
\cos\tangle\mc\kappa\{s_{\min};r,\ell\},
\end{aligned}
\leqno(*)\mc\pm_\ell
\]

\[
\begin{aligned}
\sn\kappa s&\cdot\cs\kappa S\cdot\tfrac{\partial^+}{\partial s}\ell
+
\sn\kappa r\cdot\cs\kappa R\cdot\tfrac{\partial^+}{\partial r}\ell\le 
\\
&\le\sn\kappa s\cdot\cs\kappa S\cdot\cos\tangle\mc\kappa\{r;s,\ell\}
+
\sn\kappa r\cdot\cs\kappa R\cdot\cos\tangle\mc\kappa\{s;r,\ell\}.
\end{aligned}
\leqno(**)\mc\pm_\ell.
\]
Let
\[f=-\tfrac{1}{\kappa}\cdot\cs\kappa\circ\distfun{p}{}{}
=
\md\kappa\circ\distfun{p}{}{}-\tfrac{1}{\kappa}.\leqno(A)\mc\pm\]
Clearly $f''+\kappa\cdot  f\le 0$ and
\[
\begin{aligned}
\rho^+(r)&=\frac{1}{\tang\kappa r\cdot\cs\kappa R}\cdot\nabla_{\rho(r)} f,
\\
\sigma^+(s)&=\frac{1}{\tang\kappa s\cdot\cs\kappa S}\cdot\nabla_{\sigma(s)} f.
\end{aligned}
\leqno(B)\mc\pm\]
Thus from \ref{lem:grad-lip}, we have
\[\begin{aligned}
\tfrac{\partial^+}{\partial r}\ell
&=
-\frac{1}{\tang\kappa r\cdot\cs\kappa R}
\cdot
\<\nabla_{\rho(r)} f,\dir{\rho(r)}{\sigma(s)}\>
\le
\\
&\le
\frac
{1}
{\tang\kappa r\cdot\cs\kappa R}
\cdot
\frac
{\cs\kappa S-\cs\kappa R\cdot\cs\kappa\ell}
{\kappa\cdot\sn\kappa\ell}
=
\\
&=
\frac
{\frac{\cs\kappa S}{\cs\kappa R}-\cs\kappa\ell}
{\kappa\cdot\tang\kappa r\cdot\sn\kappa\ell}.
\end{aligned}
\leqno(C)\mc\pm\]
Note that for all $\kappa\ne 0$,
the function $x\mapsto\frac{1}{\kappa\cdot\cs\kappa x}$ is increasing.
Thus, since $R(r)\le r$ and $S(s_{\min})=s_{\min}$, we have 
\[\begin{aligned}
\tfrac{\partial^+}{\partial r}\ell(r,s_{\min})
&\le 
\frac
{\frac{\cs\kappa s_{\min}}{\cs\kappa r}-\cs\kappa\ell}
{\kappa\cdot\tang\kappa r\cdot\sn\kappa\ell}
=
\\
&=
\frac
{{\cs\kappa s_{\min}}-\cs\kappa\ell\cdot\cs\kappa r}
{\kappa\cdot\sn\kappa r\cdot\sn\kappa\ell}=
\\
&=\cos\tangle\mc\kappa\{s_{\min};r,\ell\},
  \end{aligned}\leqno(D)\mc\pm\]
which is the first inequality in $(*)\mc\pm_\ell$ for $\kappa\ne 0$.
By switching $\rho$ and $\sigma$ we obtain the second inequality in $(*)\mc\pm_\ell$.
Further, adding $(C)\mc\pm$ and its mirror-inequality for $\tfrac{\partial^+}{\partial s}\ell$, we have
\[\begin{aligned}
\sn\kappa r\cdot\cs\kappa R\cdot\tfrac{\partial^+}{\partial r}\ell
&+
\sn\kappa s\cdot\cs\kappa S\cdot\tfrac{\partial^+}{\partial s}\ell\le
\\
&\le
\frac
{{\cs\kappa S}\cdot\cs\kappa r-\cs\kappa\ell\cdot\cs\kappa R\cdot\cs\kappa r}
{\kappa\cdot\sn\kappa\ell}
+\\
&\quad +
\frac
{{\cs\kappa R}\cdot\cs\kappa s-\cs\kappa\ell\cdot\cs\kappa S\cdot\cs\kappa s}
{\kappa\cdot\sn\kappa\ell}=
\\
&=
\sn\kappa r\cdot\cs\kappa R\cdot
\frac
{\cs\kappa s-\cs\kappa\ell\cdot\cs\kappa r}
{\kappa\cdot\sn\kappa r\cdot\sn\kappa\ell}
+\\
&\quad +
\sn\kappa s\cdot\cs\kappa S\cdot
\frac
{\cs\kappa r-\cs\kappa\ell\cdot\cs\kappa s}
{\kappa\cdot\sn\kappa s\cdot\sn\kappa\ell}
=
\\
&=\sn\kappa r\cdot\cs\kappa R\cdot\cos\tangle\mc\kappa\{r;s,\ell\}
+\\
&\quad +\sn\kappa s\cdot\cs\kappa S\cdot\cos\tangle\mc\kappa\{s;r,\ell\},
\end{aligned}
\leqno(E)\mc\pm\]
which is $(**)\mc\pm_\ell$.\qeds

\begin{thm}{Exercise}\label{ex:geodesic}
Suppose $\spc{L}$ is a complete length $\Alex\kappa$ space 
and $x,y,z\in \spc{L}$.
Assume $\angk\kappa zxy=\pi$.
Show that there is a geodesic $[xy]$
that contains $z$.
In particular, $x$ can be connected to $y$ by a minimizing geodesic.
\end{thm}

\section{Gradient exponent}\label{sec:gexp}

Let $\spc{L}$ be a complete length $\Alex{\kappa}$ space, 
$p\in \spc{L}$, 
and $\xi\in \Sigma_p$.
Consider a sequence of points $x_n\in \spc{L}$ such that $\dir p{x_n}\to \xi$.
Let $r_n=\dist{p}{x_n}{}$, and let
$\sigma_n\:[r_n,\tfrac{\varpi\kappa}2)\to \spc{L}$ be the $(p,\kappa)$-radial curve that starts at~$x_n$.

By the radial comparison (\ref{rad-comp}), 
the curves $\sigma_n\:[r_n,\tfrac{\varpi\kappa}2)\to \spc{L}$ 
converge to a curve $\sigma_\xi\:(0,\tfrac{\varpi\kappa}2)\to \spc{L}$, 
and this limit is independent of the choice of the sequence $x_n$.
Let $\sigma_\xi(0)=p$, and if $\kappa>0$ define \[\sigma_\xi(\tfrac{\varpi\kappa}2)
=
\lim_{t\to\frac{\varpi\kappa}2}\sigma_\xi(t).\]
The resulting curve $\sigma_\xi$ will be called the \index{radial curve}\emph{$(p,\kappa)$-radial curve} in direction~$\xi$.

The \index{gradient exponential map}\emph{gradient exponential map} 
$\gexp\mc\kappa_p\:\cBall[\0,\tfrac{\varpi\kappa}2]_{\T_p}\to \spc{L}$
is defined by
\[
\gexp\mc\kappa_p\: r\cdot\xi\mapsto\sigma_\xi(r).
\]
\index{$\gexp\mc\kappa_p$}

Here are properties of radial curves reformulated in terms of the gradient exponential map:

\begin{thm}{Theorem}\label{thm:prop-gexp}
Let $\spc{L}$ be a complete length $\Alex{\kappa}$ space. 
Then:
\begin{subthm}{}
If $p,q\in \spc{L}$ are points such that $\dist{p}{q}{}\le\tfrac{\varpi\kappa}2$, then for any geodesic $[pq]$ in $\spc{L}$ we have
\[\gexp\mc\kappa_p(\ddir p q)=q.\] 
\end{subthm}

\begin{subthm}{thm:prop-gexp:short} 
For any $v,w\in \cBall[\0,\tfrac{\varpi\kappa}2]_{\T_p}$,
\[\dist{\gexp\mc\kappa_p v}{\gexp\mc\kappa_p w}{}
\le
\side\kappa\hinge{\0}v w.\]
In other words, if we denote by $\mathcal{T}_{p}\mc\kappa$ the set $\cBall[\0,\tfrac{\varpi\kappa}2]_{\T_p}$ 
equipped with the metric $\dist{v}{w}{\mathcal{T}\mc\kappa_{p}}=\side\kappa\hinge{\0}v w$, 
then 
\[\gexp\mc\kappa_p:\mathcal{T}\mc\kappa_{p}\to \spc{L}\] 
is a short map.
\end{subthm}

\begin{subthm}{gexp-mono} 
Suppose
$p, q\in \spc{L}$ 
and $\dist{p}{q}{}\le \tfrac{\varpi\kappa}2$.
If $v\in\T_p$, $|v|\le 1$, and 
\[\sigma(t)=\gexp\mc\kappa_p(t\cdot v),\]
then the function
\[
s
\mapsto 
\tangle\mc\kappa(\sigma|_0^s,q)
\df
\tangle\mc\kappa\{\dist{q}{\sigma(s)}{};\dist{q}{\sigma(0)}{},s\}
\]
is nonincreasing in its domain of definition.
\end{subthm}
\end{thm}

\parit{Proof.}
Follows directly from the construction of $\gexp\mc\kappa_p$ and the radial comparison (\ref{rad-comp}).
\qeds

Applying the theorem above together with \ref{LinDim+-f},
we obtain the following.

{\sloppy 
 
\begin{thm}{Corollary}\label{cor:short-map-to-ball}
Let $p$ be a point in an $m$-dimensional complete length $\Alex\kappa$ space $\spc{L}$,
$m<\infty$, and $0<R\le\tfrac{\varpi\kappa}2$.
Then there is a short map 
$f\:\cBall[R]_{\Lob{m}{\kappa}}\to \spc{L}$
such that $\Im f= \cBall[p,R]\subset \spc{L}$.
\end{thm}

}

\begin{thm}{Exercise}\label{ex:gexp} 
Let $\spc{L}\subset\EE^2$ be the Euclidean halfplane. 
Clearly $\spc{L}$ is a two-dimensional complete length $\Alex{0}$ space.
Given a point $x\in \EE^2$, denote by $\proj(x)$ the closest point to $x$ on $\spc{L}$. 

Apply the radial comparison (\ref{rad-comp}) to show that for any interior point $p\in \spc{L}$ and any $v\in\R^2$  we have 
\[\gexp_p v=\proj(p+v).\]
\end{thm}

\begin{thm}{Exercise}\label{ex:inv-gexp}
Suppose $x,p,$ and $q$ are points in a complete length $\Alex\kappa$ space, and $x\in [pq[$.
Show that there is a unique vector $v\in\T_p$ such that $\gexp_p v=x$.
\end{thm}

%

\section{Remarks}

\subsection*{Gradient flow on Riemannian manifolds}
The gradient flow for general semiconcave functions 
on smooth Riemannian manifolds  can be introduced with much less effort.
To do this note that the distance estimates proved in the Section~\ref{sec:grad-curv:dist-est}
can be proved in the same way for gradient curves of smooth semiconcave subfunctions.
By the Greene--Wu lemma \cite{greene-wu}, 
given 
a $\lambda$-concave function $f$, 
a compact set $K\subset \Dom f$,
and $\eps>0$
there is a smooth $(\lambda-\eps)$-concave function that is 
$\eps$-close to $f$ on $K$.
Hence one can apply smoothing and pass to the limit as $\eps\to0$.
Note that by the second distance  estimate (\ref{lem:fg-dist-est}), the  limit curve obtained does not depend on the smoothing.

\subsection*{Gradient curves of a family of functions}

Gradient flow can be extended to a family of functions.
This type of flow was studied by Chanyoung Jun \cite{jun-thesis,jun:grad}, by Lucas Ferreira and Julio Valencia-Guevara \cite{ferreira-valencia}, and by Alexander Mielke, Riccarda Rossi, and Giuseppe Savar\'{e} \cite{mielke-rossi-savare}.
We will follow the simplified and generalized approach given by Alexander Lytchak and the third author \cite{lytchak-petrunin-2020}, where an application related to this type of flow is given.
The original motivation of Chanyoung Jun came from the study of pursuit-evasion problems.
Another application of this type of flow comes from the fact that
 the optimal transport plan, or equivalently geodesics in the Wasserstein metric, can be described as gradient flow for a family of semiconcave functions.
This observation was used by the third author to prove that Alexandrov spaces with nonnegative curvature have nonnegative Ricci curvature in the sense of Lott--Villani--Sturm \cite{petrunin:optimal}.

Suppose that $\spc{Z}$ is either $\Alex{}$ or $\CAT{}$.
Let $f_t$ be a family of functions defined on open subsets $\Dom f_t$ of~$\spc{Z}$.
More precisely, we assume that the parameter $t$ lies in a real interval $\II$ and 
\[\Omega=\set{(x,t)\in\spc{Z}\times \II}{x\in\Dom f_t}\]
is an open subset in $\spc{Z}\times \II$.

A family of functions $f_t$ is called \index{Lipschitz family}\emph{Lipschitz} if 
the function $(x,t)\mapsto f_t(x)$ is 
$L$-Lipschitz for some constant $L$.

A family of functions $f_t$ will be called \index{semiconcave family of functions}\emph{semiconcave} if 
the function $x\mapsto f_t(x)$ is $\lambda$-concave for each $t$.
A family $f_t$ is called \emph{locally semiconcave} if for each $(p_0,t_0)\in \Omega$ there is a neighborhood $\Omega'$ and $\lambda\in\RR$ such that the restriction of $f_t$ to $\Omega'$ is semiconcave. 

One cannot expect that a direct generalization of Definition \ref{def:grad-curve}  holds for every family of functions $f_t$; that is, gradient curves of a family $f_t$ cannot be defined as curves satisfying the equation $\alpha^+=\nabla_{\alpha} f$.

For example, consider a $1$-Lipschitz curve $\alpha$ in the real line. 
It is reasonable to assume that $\alpha$ is an $f_t$-gradient curve for the family $f_t(x)\z=-|x-\alpha(t)|$.
(Indeed $\alpha$ can be realized as a limit of  gradient curves for a family of functions obtained by smoothing $f_t$.)
On the other hand, $\alpha^+(t)$ might be undefined,
and even if it is defined, in general $\alpha^+(t)\ne0$,  while $\nabla_{\alpha(t)} f_t\equiv0$.

Instead we define an {}\emph{$f_t$-gradient curve} as a Lipschitz curve $\alpha$ that satisfies the following inequality
for any point $p$, time $t$, and
small $\eps >0$:  
\[\distfun{p}{}{}\circ\alpha(t+\eps)\le \distfun{p}{}{}\circ\alpha(t)-\eps\cdot \dd_{\alpha(t)}f_t(\dir{\alpha(t)}p)+o(\eps).\eqlbl{def:tdflow}\]
If there is no geodesic $[\alpha(t)\,p]$ then we impose no condition.

If $\alpha^+(t)=\nabla_{\alpha(t)}f_t$ for all $t$, then \ref{def:tdflow} holds by the definition of gradient (\ref{def:grad}).
On the other hand, the example above shows that the converse does not hold;
that is, \ref{def:tdflow} generalizes Definition \ref{def:grad-curve}.
The defining  inequality \ref{def:tdflow} is closely related to the so-called \index{evolution variational inequality}\emph{evolution variational inequality} \cite[Theorem 4.0.4(iii)]{ambrosio-gigli-savare}.

\begin{thm}{Distance estimate}\label{Distance estimate}
Let $f_t$ and $h_t$ be two Lipschitz families of $\lambda$-concave functions on a complete length space $\spc{Z}$, and $s\ge 0$.
Suppose that $\spc{Z}$ is either $\Alex{}$ or $\CAT{}$.
Assume $f_t$ and $h_t$ have common domain $\Omega\subset {\spc{Z}}\times \RR$, and $|f_t(x)-h_t(x)|\le s$ for any $(x,t)\in \Omega$.
Assume $t\mapsto \alpha(t)$ and $t\mapsto \beta(t)$ are $f_t$- and $h_t$-gradient curves respectively defined on a common interval $t\in [a,b)$, and let $\ell(t)\z=\dist{\alpha(t)}{\beta(t)}{\spc{Z}}$.
If for all $t$, a minimizing geodesic $[\alpha(t)\,\beta(t)]$ lies in $\set{x\in {\spc{Z}}}{(x,t)\in \Omega}$, then
\[\ell'(t)\le \lambda\cdot\ell(t)+2\cdot s/\ell(t),\]
whenever the left-hand side is defined.
Moreover,
\[\ell(t)^2+\tfrac{2\cdot s}\lambda\le(\ell(a)^2+\tfrac{2\cdot s}\lambda)\cdot e^{2\cdot\lambda\cdot (t-a)}.\]

In particular, these inequalities hold for any $t\in\II$ if $\Omega\supset B(p,2\cdot r)\times \II$ and $\alpha(t),\beta(t)\z\in B(p, r)$ for any $t\in \II$.
\end{thm}

Note that if $f_t=h_t$ then $s=0$;
in this case the second inequality can be written as
\[\ell(t)\le \ell(a)\cdot e^{\lambda\cdot (t-a)}.\eqlbl{dist-est-s=0}\]
In particular, it implies uniqueness of the future of gradient curves with given initial data.
This inequality also makes it possible to estimate the distance between two gradient curves for close functions.
In particular, it implies convergence for $f_t^n$-gradient curves if a sequence of $\Lip$-Lipschitz and $\lambda$-concave families $f^n_t$ converges uniformly as $n\to \infty$. 

\parit{Proof of \ref{Distance estimate}.}
Fix a time moment $t$ and set $f=f_t$ and $h=h_t$.
Let $m$ be the midpoint of the geodesic $[\alpha(t)\beta(t)]$.
Let $\gamma\:[0,\ell(t)]\to \spc{Z}$ be an arclength parametrization of $[\alpha(t)\beta(t)]$.
Note that $\dd_{\alpha(t)}f(\dir{\alpha(t)}{m})$ is the right derivative of $f\circ\gamma$ at $0$
and $-\dd_{\alpha(t)}h(\dir{\beta(t)}m)$ is the left derivative of $h\circ\gamma$ at $\ell(t)$.
Since $f$ and $h$ are $\lambda$-concave,
\begin{align*}
f\circ\beta(t)&\le f\circ\alpha(t)+\ell(t)\cdot \dd_{\alpha(t)}f(\dir{\alpha(t)}{m}) +\tfrac12\cdot\lambda\cdot\ell(t)^2,
\\
h\circ\alpha(t)&\le h\circ\beta(t)+\ell(t)\cdot \dd_{\alpha(t)}h(\dir{\beta(t)}m) +\tfrac12\cdot\lambda\cdot\ell(t)^2.
\end{align*}
Adding these inequalities and taking into account  $|f(x)-h(x)|<s$ for any $x$, we conclude that 
\[\dd_{\alpha(t)}f(\dir{\alpha(t)}{m})+\dd_{\alpha(t)}h(\dir{\beta(t)}m)\ge \lambda\cdot \ell(t)+2\cdot s/\ell(t).\]

Applying the triangle inequality and \ref{def:tdflow} at $m$, we obtain
\begin{align*}
\ell(t+\eps)&=\dist{\alpha(t+\eps)}{\beta(t+\eps)}{}\le
\\
&\le \dist{\alpha(t+\eps)}{m}{}+\dist{\beta(t+\eps)}{m}{}\le 
\\
&\le \dist{\alpha(t)}{m}{}-\eps\cdot \dd_{\alpha(t)}f(\dir{\alpha(t)}{m})+
\\
&\quad +\dist{\beta(t+\eps)}{m}{}-\eps\cdot \dd_{\beta(t)}h(\dir{\beta(t)}m)+o(\eps)=
\\
&=\ell(t)-\eps\cdot(\lambda\cdot \ell(t)+2\cdot s/\ell(t))+o(\eps)
\end{align*}
for $\eps>0$. The first inequality follows.

Since $\alpha$ and $\beta$ are Lipschitz, $t\mapsto \ell(t)$ is a Lipschitz function.
By Rademacher's theorem, its derivative $\ell'$ is defined almost everywhere and satisfies the fundamental theorem of calculus.
Therefore the first inequality implies the second.
\qeds

\begin{thm}{Proposition}\label{prop:def-time-dependent}
Suppose  $\spc{Z}$ is a complete length space that is either $\Alex{}$ or $\CAT{}$.
Let $f_t$ be a family of $\lambda$-concave functions for $t\in [a,b)$, where $\Dom f_t\supset B(z,2\cdot r)$ for fixed $z\in\spc{Z}$, $r>0$ and any~$t$.

Let $\alpha\:[a,b)\to B(z,r)$ be Lipschitz.
Then $\alpha$ is an $f_t$-gradient curve if and only if 
\[\begin{aligned}
&\distfun{p}{}{}\circ\alpha(t+\eps)\le 
\\
&\quad\le \distfun{p}{}{}\circ\alpha(t)-\eps\cdot \left[\frac{f_t(p)-f_t\circ\alpha(t)}{\dist{p}{\alpha(t)}{}}-\tfrac\lambda2\cdot \dist{p}{\alpha(t)}{}\right]+o(\eps)
\end{aligned}
\eqlbl{def:tdflow-plus}\]
for any $t\in [a,b)$ and $p\in B(z,r)\setminus \{\alpha (t)\}$.
\end{thm}

\parit{Proof.}
Note that the geodesics $[\alpha(t)p]$ lie in $\Dom f_t$ for any $t$.

Since $f_t$ is $\lambda$-concave, we have 
\[\dd_{\alpha(t)}f_t(\dir{\alpha(t)}p)
\ge
\frac{f(p)-f\circ\alpha(t)}{\dist{p}{\alpha(t)}{}}-\tfrac\lambda2\cdot \dist{p}{\alpha(t)}{}.\]
Hence the only-if part follows.

Given  $p\in \spc{Z}$ and $t$,
consider a point $\bar p\in [\alpha(t)p]$.
Applying \ref{def:tdflow-plus} for $\bar p$, and the triangle inequality, we have
\[\distfun{p}{}{}\circ\alpha(t+\eps)
\le
\distfun{p}{}{}\circ\alpha(t)-\eps\cdot \left[\frac{f(\bar p)-f\circ\alpha(t)}{\dist{\bar p}{\alpha(t)}{}}-\tfrac\lambda2\cdot \dist{\bar p}{\alpha(t)}{}\right]+o(\eps).\]
By taking $\bar p$ close to $\alpha(t)$,
the value $\tfrac{f(\bar p)-f\circ\alpha(t)}{\dist{\bar p}{\alpha(t)}{}}-\tfrac\lambda2\cdot \dist{\bar p}{\alpha(t)}{}$ can be made arbitrarily close to $\dd_{\alpha(t)}f_t(\dir{\alpha(t)}p)$.
Therefore, given $\delta>0$, the inequality
\[\distfun{p}{}{}\circ\alpha(t+\eps)\le \distfun{p}{}{}\circ\alpha(t)-\eps\cdot \dd_{\alpha(t)}f_t(\dir{\alpha(t)}p)+\eps\cdot\delta\]
holds for all sufficiently small $\eps>0$.
Therefore \ref{def:tdflow} holds.
\qeds

Now we are ready to formulate and prove global existence of gradient curves for time-dependent families --- an analog of \ref{thm:glob-exist-grad-curv}.

\begin{thm}{Theorem}\label{prop:time-dependent}
Suppose $\spc{Z}$ is a complete length space that is either $\Alex{}$ or $\CAT{}$.
Let
$\{f_t\}$ be a family of functions defined on an open set
\[\Omega=\set{(x,t)\in \spc{Z}\times\RR}{x\in \Dom f_t}.\]
Suppose that $f_t$ is Lipschitz and locally semiconcave.
Then for any time  $a$ and initial point $p\in \Dom f_a$, there is a unique $f_t$-gradient curve $t\mapsto\alpha(t)$ defined on a maximal semiopen interval $[a,b)$. 
Moreover, if $b<\infty$ then $(\alpha(t),t)$ escapes from any closed set $K\subset \Omega$.
\end{thm}

\parit{Proof.}
Let $L$ be a Lipschitz constant of $f_t$.
Fix $b>a$ sufficiently small  that $\Dom f_t\supset B(p,\eps\cdot L)$ for any $t\in[a,b)$.
Consider a sequence  $a=t_0<t_1\dots<t_n\z=b$, and a piecewise constant family of functions on $B(p,\eps\cdot L)$ defined by $\hat f_t=f_{t_i}$ if $t_i\le t<t_{i+1}$.

Note that $\hat f_t$ is time-independent on each interval $[t_i,t_{i+1})$.
By  \ref{thm:glob-exist-grad-curv} applied recursively on each interval $[t_i,t_{i+1})$,  the proposition holds for $\hat f_t$.
That is, there is a unique $\hat f_t$-gradient curve $\hat \alpha$ that starts at $p$ and is defined on the interval $[a,b)$.

The distance estimates (\ref{Distance estimate}) show that as the partition gets finer, the gradient curves $\hat\alpha$ form a Cauchy sequence; denote its limit by $\alpha$.
Then
\begin{align*}
\distfun{p}{}{}\circ\hat\alpha(t+\eps)
&\le 
\distfun{p}{}{}\circ\hat\alpha(t)-
\\
&\quad
-\eps\cdot \left[\frac{\hat f_t(p)-\hat f_t\circ\hat\alpha(t)}{\dist{p}{\alpha(t)}{}}-\tfrac\lambda2\cdot \dist{p}{\hat\alpha(t)}{}\right] 
+o(\eps)\le
\\
&\le 
\distfun{p}{}{}\circ\hat\alpha(t)-
\\
&\quad
-\eps\cdot \left[\frac{f_t(p)-f_t\circ\hat\alpha(t)-2\cdot s}{\dist{p}{\alpha(t)}{}}-\tfrac\lambda2\cdot \dist{p}{\hat\alpha(t)}{}\right]
+o(\eps),
\end{align*}
where 
\[s=\sup_{t,x} \{|f_t(x)-\hat f_t(x)|\}.\]
Since $s\to 0$ as $\hat\alpha\to \alpha$, then \ref{def:tdflow-plus} holds for $\alpha$;
that is, $\alpha$ is an $f_t$-gradient curve.

This proves short time existence.
Applying this argument recursively, we can find a gradient curve defined on a maximal interval $[a,b)$.
Uniqueness of this curve follows from the distance estimate \ref{dist-est-s=0}. 

Note that $\alpha$ is $L$-Lipschitz.
In particular, if $b<\infty$ then $\alpha(t)\to p'$ as $t\to b$.
If $(p',b)\in \Omega$ then we can repeat the procedure; otherwise $\alpha$ escapes from any closed set in $\Omega$. 
\qeds

\subsection*{Gradient curves for non-Lipschitz functions}\label{sec:non-lip}

In this book, we only consider gradient curves for locally Lipschitz semiconcave subfunctions;
this turns out to be sufficient for all our needs.
However, 
instead of Lipschitz semiconcave subfunctions,
it is more natural to consider upper semicontinuous semiconcave functions
with target in $[-\infty,\infty)$,
and to assume in addition that 
the functions take finite values at a dense set in the domain of definition.
Suppose that $\spc{Z}$ is a complete length space that is either $\Alex{}$ or $\CAT{}$.
The set of such subfunctions on $\spc{Z}$ will be denoted by 
$\LSCSC(\spc{Z})$ (for \textbf{l}ower semi-\textbf{c}ontinous and semi-\textbf{c}oncave).

In this section we describe the adjustments needed
to construct gradient curves for the subfunctions in $\LSCSC(\spc{Z})$.

When $\spc{Z}=H$ is a Hilbert space the theory we develop is equivalent to the classical theory of gradient flows on Hilbert space mentioned earlier \cite{Brezis-book}.

Further examples of such functions include the entropy and other closely related functionals on the Wasserstein space over a $\Alex0$ space.
Another important example is given by the Cheeger energy on metric measure spaces, its gradient flow leads to the notion of the heat flow on such spaces.
The gradient flow for these functions plays an important role in the theory of optimal transport, see \cite{villani} and references there in.

\parbf{Differential.} 
First we need to extend the definition of differential (\ref{def:differential}) to $\LSCSC$ subfunctions.

Let $\spc{Z}$ be a complete length space and $f\in\LSCSC(\spc{Z})$.
Suppose that $\spc{Z}$ is either $\Alex{}$ or $\CAT{}$.
Given a point $p\in \Dom f$ and a geodesic direction $\xi=\dir pq$, 
let 
$\hat \dd_pf(\xi)=(f\circ\geod_{[pq]})^+(0)$.
Since $f$ is semiconcave, the value $\hat \dd_pf(\xi)$ is defined if $f\circ\geod_{[pq]}(t)$ is finite at all sufficiently small values $t>0$,
but $\hat \dd_pf(\xi)$ may take value $\infty$. 
Note that $\hat \dd_pf$ is defined on a dense subset of $\Sigma_p$.

Let 
\[\dd_pf(\zeta)=\limsup_{\xi\to\zeta}\hat\dd_pf(\xi),\]
and $\dd_pf(v)=|v|\cdot \dd_pf(\xi)$ if $v=|v|\cdot\xi$ for some $\xi\in\Sigma_p$.

In other words, we define differential as the smallest 
upper semi-continuous  positive-homogeneous function $\dd_pf\:\T_p\to\RR$
such that if $\hat\dd_pf(\xi)$ is defined, then $\dd_pf(\xi)\ge \hat \dd_pf(\xi)$.

\parbf{Existence and uniqueness of the gradient.}
Note that in the proof of \ref{thm:ex-grad}, 
we used the Lipschitz condition just once,
to show that 
\begin{align*}
s&=\sup\set{(\dd_p f)(\xi)}{\xi\in\Sigma_p}=
\\
&=\limsup_{x\to p}\frac{f(x)-f(p)}{\dist{x}{p}{}}<
\\
&<\infty.
\end{align*}

The value $s$ above will be denoted by $|\nabla|_pf$.
Note that 
if the gradient $\nabla_pf$ is defined then $|\nabla|_pf=|\nabla_pf|$,
and otherwise $|\nabla|_pf=\infty$.

Summarizing the discussion above, 
we have the following.

\parbf{\ref{thm:ex-grad}$'$ Existence and uniqueness of the gradient.}
\textit{Assume $\spc{Z}$ is a complete space and $f\in \LSCSC(\spc{Z})$. 
Suppose that $\spc{Z}$ is either $\Alex{}$ or $\CAT{}$.
Then for any point $p\in \Dom f$, either there is a unique gradient $\nabla_p f\in \T_p$ 
or $|\nabla|_pf=\infty$.}

\medskip

Further, in all the results of Section~\ref{sec:grad-calculus} 
we may assume only that both $f$ and the gradient of $f$ are defined at the points under consideration. The proofs are the same.

Sections \ref{sec:grad-semicont}--\ref{sec:grad-curv:exist}
require almost no changes.
Mainly, where appropriate
one needs to change $|\nabla_p f|$ 
to $|\nabla|_pf$ 
and/or assume that the gradient is defined at the points of interest.
Also,  \ref{eq:thm:grad-like-2nd-def-1} in Theorem \ref{thm:grad-like-2nd-def}
is taken as the definition of gradient-like curve.
Then the theorem states that any  gradient-like curve $\alpha\:\II\to\spc{Z}$ satisfies Definition \ref{def:grad-like-curve} at $t\in \II$ if $\nabla_{\hat\alpha(s)} f$ is defined.
Further, Definition \ref{def:grad-curve}, should be changed to the following:

\medskip

\parbf{\ref{def:grad-curve}$'$. Definition.}
{\it Let $\spc{Z}$ be a complete length space
and $f\in\LSCSC(\spc{Z})$.
Suppose that $\spc{Z}$ is either $\Alex{}$ or $\CAT{}$.

A curve 
$\alpha\:[t_{\min},t_{\max})\to\Dom f$ will be called an  \index{gradient curve}\emph{$f$-gradient curve} if
\[\alpha^+(t)=\nabla_{\alpha(t)} f\]
when $\nabla_{\alpha(t)} f$ is defined and 
\[(f\circ\alpha)^+(t)=\infty\]
otherwise.}

\medskip

In the proof of local existence (\ref{thm:exist-grad-curv}), condition (\ref{alm-grad})
should be changed to the following condition:
\begin{itemize}

\item[{(\ref{alm-grad})}$'$]
$f\circ\hat\alpha_n(\bar\varsigma_i)-f\circ\hat\alpha_n(\varsigma_i)
>
(\bar\varsigma_i-\varsigma_i)
\cdot
\max\{n,|\nabla|_{\hat\alpha_n(\varsigma_i)}f-\tfrac{1}{n})\}.$
\end{itemize}

Any gradient curve $\alpha[0,\ell)\to\spc{Z}$
for a subfunction
$f\in \LSCSC(\spc{Z})$
satisfies the equation
\[\alpha^+(t)=\nabla_{\alpha(t)} f\]
at all values $t$, with the possible exception of $t=0$.
In particular, the gradient of $f$ is defined at all points of any 
$f$-gradient curve, with the possible exception of the initial point.

\begin{thm}{Example}
Let $\spc{X}=L^2(\R^n)$. Let $F\co \spc{X}\to \R$ be given  by $F(f)=-\int|\nabla f|^2d\mathcal L$.
Then $F$ is a concave but not locally Lipschitz functional on $\spc{X}$ and its finite precisely at functions $f\in W^{1,2}(\R^n)\subset L^2(\R^n)$.
Integration by parts shows that for any smooth $f\in W^{1,2}(\R^n)$ it holds that $\nabla_fF=\Delta f$.
The gradient flow of F is given by the heat flow $f_t$ starting at $f$ and $f_t$ is smooth for all positive $t$.
\end{thm}

\subsection*{Slower radial curves}
Let $\kappa\ge 0$. 
Assume that for some function $\psi$, the curves defined by the equation 
\[\sigma^+(s)=\psi(s,\dist{p}{\sigma(s)}{})\cdot\nabla_{\sigma(s)}\distfun{p}{}{}\]
satisfy radial comparison \ref{rad-comp}.
Then in fact the $\sigma(s)$ are radial curves; 
that is, 
\[\psi(s,\dist{p}{\sigma(s)}{})= \frac{\tang\kappa\dist[{{}}]{p}{\sigma(s)}{})}{\tang\kappa s},\]
see exercise \ref{ex:gexp}.

In case $\kappa<0$, such a function $\psi$ is not unique.
In particular, one can take curves defined by the simpler equation
\[\sigma^+(s)
=
\frac{\sn\kappa \dist[{{}}]{p}{\sigma(s)}{}}{\sn\kappa s}\cdot\nabla_{\sigma(s)}\distfun{p}{}{}
=
\frac{1}{\sn\kappa s}\cdot\nabla_{\sigma(s)}(\md\kappa\circ\distfun{p}{}{}).\]
Among all curves of that type, the radial curves for curvature $\kappa$ 
as defined in \ref{def:rad-curv} maximize the growth of $\dist{p}{\sigma(s)}{}$.


\appendix
\chapter{Semisolutions}
\parbf{\ref{exr-crofton}.}
Suppose $\alpha$ is a closed spherical curve. 
By Crofton's formula, the length of  $\alpha$  is $\pi\cdot n_\alpha$, where $n_\alpha$ denotes the average number of crossings of $\alpha$ with equators.

Since $\alpha$ is closed, almost all equators cross it at an even number of points (we assume that $\infty$ is an even number).
If $\length \alpha<2\cdot\pi$ then $n_\alpha<2$.
Therefore there is an equator that does not cross $\alpha$ --- hence the result.

\parbf{\ref{ex:complete=>complete}}; \ref{SHORT.ex:complete=>complete:complete}.
Note that any Cauchy sequence $x_n$ in $(\spc{X},\yetdist{}{}{})$ is also Cauchy in $\spc{X}$.
Since $\spc{X}$ is complete, $x_n$ converges; denote its limit by $x_\infty$.

Passing to a subsequence, we may assume that $\yetdist{x_{n-1}}{x_n}{}<\tfrac1{2^n}$.
It follows that there is a 1-Lipschitz curve $\alpha\:[0,1]\to (\spc{X},\yetdist{}{}{})$ such that $x_n=\alpha(\tfrac1{2^n})$ and $x_\infty=\alpha(0)$.
In particular, $\yetdist{x_n}{x_\infty}{}\to0$ and $n\to\infty$.

\parit{\ref{SHORT.ex:complete=>complete:compact}.}
Fix two points $x,y\in\spc{X}$ such that $\ell=\yetdist{x}{y}{}<\infty$.
Let $\alpha_n$ be a sequence of paths from $x$ to $y$ such that $\length(\alpha_n)\to\ell$ as $n\to \infty$.
Without loss of generality, we may assume that each $\alpha_n$ is $(\ell+1)$-Lipschitz.

Since $\spc{X}$ is compact, there is a partial limit $\alpha_\infty$ of $\alpha_n$ as $n\to \infty$.
By semicontinuity of length, $\length\alpha_\infty\le\ell$;
that is; $\alpha$ is a shortest path in~$\spc{X}$.

\parit{Source:} Part \ref{SHORT.ex:complete=>complete:complete} appears as a Corollary in \cite{hu-kirk}; see also \cite[Lemma 2.3]{petrunin-stadler}.

\parbf{\ref{ex:no-geod}.}
The following example was suggested by Fedor Nazarov~\cite{nazarov}.

\medskip

Consider the unit ball $(B,\rho_0)$
in the space $c_0$ of all sequences converging to zero equipped with the sup-norm.

Consider another metric $\rho_1$ which is different from $\rho_0$ by the conformal factor
\[\phi(\bm{x})=2+\tfrac{1}2\cdot x_1+\tfrac{1}4\cdot x_2+\tfrac{1}8\cdot x_3+\dots,\]
where $\bm{x}=(x_1,x_2\,\dots)\in B$.
That is, if $\bm{x}(t)$, $t\in[0,\ell]$, is a curve parametrized by $\rho_0$-length 
then its $\rho_1$-length is 
\[\length_{\rho_1}\bm{x}=\int\limits_0^\ell\phi\circ\bm{x}.\]
Note that the metric $\rho_1$ is bi-Lipschitz equivalent  to~$\rho_0$.

Assume $\bm{x}(t)$ and $\bm{x}'(t)$ are two curves parametrized by $\rho_0$-length that differ only in the $m$-th coordinate; denote them by $x_m(t)$ and $x_m'(t)$ respectively.
Note that if $x'_m(t)\le x_m(t)$ for any $t$ and 
the function $x'_m(t)$ is locally $1$-Lipschitz at all $t$ such that $x'_m(t)< x_m(t)$, then 
\[\length_{\rho_1}\bm{x}'\le \length_{\rho_1}\bm{x}.\]
Moreover this inequality is strict if $x'_m(t)< x_m(t)$ for some~$t$.

Fix a curve $\bm{x}(t)$, $t\in[0,\ell]$, parametrized by  $\rho_0$-length.
We can choose $m$ large so that $x_m(t)$ is sufficiently close to $0$ for any~$t$.
In particular, for some values $t$, we have $y_m(t)<x_m(t)$, where
\[y_m(t)=(1-\tfrac t\ell)\cdot x_m(0)
+\tfrac t\ell\cdot x_m(\ell)
-\tfrac 1{100}\cdot \min\{t,\ell-t\}.\]
Consider the curve $\bm{x}'(t)$ as above with
\[x'_m(t)=\min\{x_m(t),y_m(t)\}.\]
Note that $\bm{x}'(t)$ and $\bm{x}(t)$ have the same endpoints, and by the above
\[\length_{\rho_1}\bm{x}'<\length_{\rho_1}\bm{x}.\]
That is, for any curve $\bm{x}(t)$ in $(B,\rho_1)$, we can find a shorter curve $\bm{x}'(t)$ with the same endpoints.
In particular, $(B,\rho_1)$ has no geodesics.

\parbf{\ref{ex:compact+connceted}.}
Choose a sequence of positive numbers $\varepsilon_n\to 0$ and an $\varepsilon_n$-net $N_n$ of $K$ for each $n$.
Assume $N_0$ is a one-point set, so $\eps_0>\diam K$.
Connect each point $x\in N_{k+1}$ to a point $y\in N_{k}$ by a curve of length at most $\eps_k$.

Consider the union $K'$ of all these curves with $K$; observe that $K'$ is compact and path-connected.

\parit{Source:} This problem was suggested by Eugene Bilokopytov \cite{bilokopytov}.

\begin{wrapfigure}{r}{20 mm}
\vskip-4mm
\centering
\includegraphics{mppics/pic-5}
\end{wrapfigure}

\parbf{\ref{exercise from BH}.}
Consider the following subset of $\R^2$ equipped with the induced length metric
\[
\spc{X}
=
\bigl((0,1]\times\{0,1\}\bigr)
\cup
\bigl(\{1,\tfrac12,\tfrac13,\dots\}\times[0,1]\bigr).
\]
Note that $\spc{X}$ is locally compact and geodesic.

Its completion $\bar{\spc{X}}$ is isometric to the closure of $\spc{X}$ equipped with the induced length metric;
$\bar{\spc{X}}$ is obtained from $\spc{X}$ by adding two points $p\z=(0,0)$ and $q\z=(0,1)$.

The point $p$ admits no compact neighborhood in $\bar{\spc{X}}$ 
and there is no geodesic connecting $p$ to $q$ in~$\bar{\spc{X}}$.

\parit{Source:} This example is taken from~\cite{bridson-haefliger}.

\parbf{\ref{ex:compact-in-lenght}.}
Let $\spc{X}$ be a compact metric space.
Let us identify $\spc{X}$ with its image in $\Bnd(\spc{X},\RR)$ under the Kuratowsky embedding (Section~\ref{Kuratowsky embedding}). 
Denote by $\spc{K}$ the \textit{linear} convex hull of $\spc{X}$ in the space of bounded functions on $\spc{X}$; 
that is, $x\in \spc{K}$ if and only if $x$ cannot be separated from $\spc{X}$ by a hyperplane.

Since $\spc{X}$ is compact, so is $\spc{K}$.
It remains to observe that $\spc{K}$ is a length space since it is convex.

\parit{Remark.}
Alternatively, one can use the embedding of $\spc{X}$ into its injective hull; see \cite{isbell}.

\parbf{\ref{ex:ultrakatetov}.} 
Let $F=\set{n\in \NN}{f(n)=n}$; we need to show that $\o(F)=1$.

Consider an oriented graph $\Gamma$ with vertex set $\NN\setminus F$ such that $m$ is connected to $n$ if $f(m)=n$.
Show that each connected component of $\Gamma$ has at most one cycle.
Use it to subdivide vertices of $\Gamma$ into three sets $S_1$, $S_2$, and $S_3$ such that $f(S_i)\cap S_i=\emptyset$ for each $i$.

Conclude that $\o(S_1)=\o(S_2)=\o(S_3)=0$ and hence \[\o(F)=\o(\NN\setminus(S_1\cup S_2\cup S_3))=1.\]

\parit{Source:} 
The presented proof was given by Robert Solovay \cite{solovay}, but
the key statement is due to Miroslav Katětov \cite{katetov}.

\parbf{\ref{ex:linear}.}
Choose a nonprincipal ultrafilter $\o$ and set $L(\bm{s})=s_\o$.
It remains to observe that $L$ is linear.

\parit{Remark.}
By this exercise, $\o$ corresponds to a vector in $(\ell^\infty)^*\setminus\ell^1$. 

\parbf{\ref{ex:ultrakatetov+}.}
Use \ref{ex:ultrakatetov}.

\parbf{\ref{ex:Asym(Lob)}}; \ref{SHORT.ex:Asym(Lob):metric-tree}.
Show that there is $\delta>0$ such that sides of any geodesic triangle in $\Lob21$ intersect a disk of radius $\delta$.
Observe that $\Lob2n\z=\tfrac1{\sqrt{n}}\cdot\Lob21$, and use it to show that any geodesic triangle in $\spc{T}$ is a tripod.

\parit{\ref{SHORT.ex:Asym(Lob):homogeneous}.} Observe and use that $\Lob2n$ are homogeneous.

\parit{\ref{SHORT.ex:Asym(Lob):continuum}.} 
Choose $p_1\in \Lob21$, denote by $p_n$ the corresponding point in $\Lob2n\z=\tfrac1{\sqrt{n}}\cdot\Lob21$
Suppose $p_n\to p_\o$ as $n\to\o$; we can assume that $p_\o\in\spc{T}$.
By \ref{SHORT.ex:Asym(Lob):homogeneous}, it is sufficient to show that $p_\o$ has a continuum degree.

Choose distinct geodesics $\gamma_1,\gamma_2\:[0,\infty)\to \Lob21$ that start at a point $p_1$.
Show that the limits of $\gamma_1$ and $\gamma_2$ run in the different connected components of $\spc{T}\setminus \{p_\o\}$.
Since there is a continuum of distinct geodesics starting at $p$,
we get that the degree of $p_\o$ is at least continuum.

On the other hand, the set of sequences of points $q_n\in\Lob2n$  has cardinality continuum.
In particular, the set of points in $\spc{T}$ has cardinality at most continuum.
It follows that the degree of any vertex is at most continuum.

\parit{Remark.}
The properties \ref{SHORT.ex:Asym(Lob):homogeneous} and \ref{SHORT.ex:Asym(Lob):continuum} describe the tree $\spc{T}$ up to isometry \cite{dyubina-polterovich}.
In particular, $\spc{T}$ does not depend on the choice of the ultrafilter.

\parbf{\ref{ex:isom-ultrapower}.}
Show and use that the spaces $\spc{X}^\o$ and $(\spc{X}^\o)^\o$ have discrete metrics and both have cardinality of the continuum.

\parbf{\ref{ex:ultrapower(ultrapower)}.}
Choose a bijection $\iota\:\NN\to \NN\times \NN$.
Given a set $S\subset \NN$, consider the sequence $S_1$, $S_2,\dots$ of subsets in $\NN$ defined by $m\in S_n$ if $(m,n)\z=\iota(k)$ for some $k\in S$.
Set $\o_1(S)=1$ if and only if $\o(S_n)=1$ for $\o$-almost all $n$.
It remains to check that $\o_1$ meets the conditions of the exercise.

\parit{Comment.}
It turns out that $\o_1\ne \o$ for any $\iota$;
see the post of Andreas Blass \cite{blass}.

\begin{wrapfigure}{r}{40 mm}
\vskip-0mm
\centering
\includegraphics{mppics/pic-603}
\end{wrapfigure}

\parbf{\ref{ex:notproper-limit}.} Consider the infinite metric $\spc{T}$ tree with unit edges shown
on the diagram. Observe that $\spc{T}$ is proper.

Consider the vertex $v_\o=\lim_{n\to\o}v_n$ in the ultrapower $\spc{T}^\o$.
Observe that $\o$ has an infinite degree.
Conclude that $\spc{T}^\o$ is not locally compact.

\parbf{\ref{ex:nonconvex-limit}.}
Let $\spc{X}_n$ be the square $\{( x,y)\in \R^2, |x|\le 1, |y|\le 1\}$ with the metric induced by the $\ell^n$-norm and let $f_n(x,y)=x$ for all $n$.
Observe that $\spc{X}_\o$ is the square with the metric induced by the $\ell^\infty$-norm where the limit function $f_\o(x,y)=x$ is not concave.

\parbf{\ref{ex:GH-SC}}; \ref{SHORT.ex:GH-SC:circle}.
Suppose $\spc{X}_n\GHto \spc{X}$ and $\spc{X}_n$ are simply connected length metric spaces.
It is sufficient to show that any nontrivial covering map $f\:\tilde{\spc{X}}\to \spc{X}$ corresponds to a nontrivial covering map $f_n\:\tilde{\spc{X}}_n\to \spc{X}_n$ for large $n$.

The latter can be constructed by covering $\spc{X}_n$ by small balls that lie close to sets in $\spc{X}$ evenly covered by $f$, preparing a few copies of these sets and gluing them in the same way as the inverse images of the evenly covered sets in $\spc{X}$ are glued to obtain $\tilde{\spc{X}}$.

\begin{wrapfigure}{r}{43 mm}
\vskip-8mm
\centering
\includegraphics{mppics/pic-2}
\end{wrapfigure}

\parit{\ref{SHORT.ex:GH-SC:nonsc-limit}.}
Let $\spc{V}$ be a cone over Hawaiian earrings.
Consider the {}\emph{doubled cone} $\spc{W}$ --- two copies of $\spc{V}$ with  their base points glued (see the diagram).

The space $\spc{W}$ can be equipped with a length metric
(for example, the induced length metric from the shown embedding).

Show that $\spc{V}$ is simply connected, but $\spc{W}$ is not; the latter is a good exercise in topology.

If we delete from the earrings all small circles, then the obtained double cone becomes simply connected and it remains close to $\spc{W}$.
That is, $\spc{W}$ is a Gromov--Hausdorff limit of simply connected spaces.

\parit{Remark.}
Note that from part \ref{SHORT.ex:GH-SC:nonsc-limit}, the limit does not admit a nontrivial covering.
So, if we define the fundamental group as the inverse image of groups of deck transformations for all the coverings of the given space, then one may say that a Gromov--Hausdorff limit of simply connected length spaces is simply connected.

\parbf{\ref{ex:sphere-to-ball},}
\textit{\ref{SHORT.ex:sphere-to-ball:2}.}
Suppose that a metric on $\mathbb{S}^2$ is close to the disk $\DD^2$.
Note that $\mathbb{S}^2$ contains a circle $\gamma$ that is close to the boundary curve of $\DD^2$.
By the Jordan curve theorem, $\gamma$ divides $\mathbb{S}^2$ into two disks, say $D_1$ and $D_2$.

By \ref{ex:GH-SC:circle}, the Gromov--Hausdorff limits of $D_1$ and $D_2$ have to contain the whole $\DD^2$, otherwise the limit would admit a nontrivial covering.

Consider points $p_1\in D_1$ and $p_2\in D_2$ that are close to the center of $\DD^2$.
If $n$ is large, the distance $\dist{p_1}{p_2}{n}$ has to be very small.
On the other hand, any curve from $p_1$ to $p_2$ must cross $\gamma$, so it has length about 2 --- a contradiction.

\parit{\ref{SHORT.ex:sphere-to-ball:3}.}
Make holes in the unit 3-disc, that do not change its topology and do not change its length metric much 
and pass to its doubling in the boundary.

\parit{Source:} The exercise is taken from \cite{burago-burago-ivanov}.

\parbf{\ref{ex:GH-proper-marked}.} Modify proof of \ref{thm:GH-compact}, or apply \ref{thm:ultra-GH:b}.

\parbf{\ref{ex:adjacent-angles}.}
If $\mangle\hinge pxz+\mangle\hinge pyz< \pi$, then by the triangle inequality for angles (\ref{claim:angle-3angle-inq}) we have $\mangle\hinge pxy< \pi$.
The latter implies that $[xy]$ fails to be minimizing near $p$.

\parbf{\ref{ex:tangent-vect=o(t)}.}
By the definition of a right derivative, there is a geodesic $\gamma$ such that both limits 
\[\limsup_{\eps\to0+}\frac{\dist{\alpha(\eps)}{\gamma(\eps)}{\spc{X}}}{\eps}
\quad\text{and}\quad
\limsup_{\eps\to0+}\frac{\dist{\beta(\eps)}{\gamma(\eps)}{\spc{X}}}{\eps}\]
are arbitrarily small.
By the triangle inequality, we get
\[\limsup_{\eps\to0+}\frac{\dist{\alpha(\eps)}{\beta(\eps)}{\spc{X}}}{\eps}=0.\]

\parbf{\ref{ex:both-sided-diff}.}
Follows directly from the definition.

\parbf{\ref{ex:diff}.}
Observe that
\[\speed_t\alpha=|\alpha^+(t)|=|\alpha^-(t)|.\]
Apply Theorem~\ref{thm:speed} to show that
\[\dist{\alpha^+(t)}{\alpha^-(t)}{\T_{\alpha(t)}}=2\cdot\speed_t\alpha.\]

\parbf{\ref{ex:schroeder-foetch}.}
Choose two non-Euclidean norms $\|{*}\|_{\spc{X}}$ and $\|{*}\|_{\spc{Y}}$ on $\RR^{10}$ such that the sum $\|{*}\|_{\spc{X}}+\|{*}\|_{\spc{Y}}$ is Euclidean.
See \cite{schroeder-foetch} for more details.

\parbf{\ref{ex:(3+1)-expanding}.} 
Assume $\dist{p}{x^i}{}=\dist{q}{y^i}{}$ for each $i$.
Observe and use that
\[\dist{x^i}{x^j}{}\le\dist{y^i}{y^j}{}
\quad\iff\quad
\angk\kappa p{x^i}{x^j}\le \angk\kappa q{y^i}{y^j}.\]

\parbf{\ref{ex:(3+1)-expanding}.} Apply the four-point comparison (\ref{df:1+3}).

\parbf{\ref{ex:nongeod-cbb}.}
Modify the induced length metric on the unit sphere in an infinite-dimensional Hilbert space in small neighborhoods of a countable collection of points.
To prove that the obtained space is $\Alex0$, you may need to use the technique from Halbeisen's example (\ref{Halbeisen's example}).

\parbf{\ref{ex:almost.geod}.} Mimic the proof of Theorem~\ref{thm:almost.geod}.

\parbf{\ref{ex:G-delta-not-thru}.}
On the plane, any nonnegatively curved metric having an everywhere dense set of singular points will do the job, where 
by singular point we mean a point having total angle around it strictly smaller than $2\cdot\pi$.

Indeed, if $x_i$ is a singular point, then there is $0<\eps_i<1/20$ such that no geodesic with ends outside of $\oBall(x_i,r)$ can meet the ball $\oBall(x_i,\eps_i\cdot r)$.
The set 
\[\Omega_n=\bigcup_i \oBall(x_i,\tfrac{\eps_i}n)\]
is open and everywhere dense.
Note that $\Omega_n$ may intersect a geodesic  of length $1/n$ only within $\frac 1 {10n}$  of its endpoints.
The intersection of the $\Omega_n$ is a G-delta dense set that does not intersect the interior of any geodesic.

\parbf{\ref{mink+alex=euclid}.} 
Note that rescaling does not change the space.
Therefore if the space is $\Alex\kappa$ then it is $\Alex{\lambda\cdot\kappa}$ for any $\lambda>0$.
Passing to the limit as $\lambda\to 0$, we may assume that the space is $\Alex0$.

The point-on-side comparison (\ref{point-on-side}) for $p=v$, $x=w$, $y=-w$ and $z=0$ implies that 
\[\|v+w\|^2+\|v-w\|^2\le 2\cdot\|v\|^2+2\cdot\|w\|^2.\]
Applying the comparison for 
$p=v+w$, $x=w-v$, $y=v-w$ and $z=0$ gives the opposite inequality.
That is, the parallelogram identity
\[\|v+w\|^2+\|v-w\|^2= 2\cdot\|v\|^2+2\cdot\|w\|^2\]
holds for any vectors $v$ and $w$.
Whence the statement follows.

\parbf{\ref{ex:cbb-geod-overlap}.}
Apply the hinge comparison (\ref{angle}).

\parbf{\ref{ex:equality-alexlemma}.}
Without loss of generality, we may assume that the points $x,v,w,y$ appear on the geodesic $[xy]$ in that order.
By the point-on-side comparison (\ref{point-on-side}) we have
\begin{align*}
\angk\kappa xyp\le\angk\kappa xwp&\le \angk\kappa xvp,
\\
\angk\kappa ywp&\ge\angk\kappa yvp\ge\angk\kappa yxp.
\end{align*}
Therefore
\begin{align*}\angk\kappa xyp<\angk\kappa xwp
\quad&\Longrightarrow\quad
\angk\kappa xyp<\angk\kappa xvp,
\\
\angk\kappa yxp<\angk\kappa ywp
\quad&\Longleftarrow\quad
\angk\kappa yxp<\angk\kappa yvp.
\end{align*}

By Alexandrov's lemma (\ref{lem:alex}), we have
\begin{align*}
\angk\kappa xyp<\angk\kappa xvp
\quad&\Longleftrightarrow\quad
\angk\kappa yxp<\angk\kappa yvp,
\\
\angk\kappa xyp<\angk\kappa xwp
\quad&\Longleftrightarrow\quad
\angk\kappa yxp<\angk\kappa ywp.
\end{align*}
Hence the statement follows.

\parbf{\ref{ex:urysohn}.}
See the construction of Urysohn's space \cite[3.11$\tfrac{3}{2}_+$]{gromov-MS} or \cite{petrunin2020pure}.

\parbf{\ref{ex:lebedeva-petrunin}.}
Read \cite{lebedeva-petrunin}.

\parbf{\ref{ex:fat-triangle}.}
Apply the angle-sidelength  monotonicity (\ref{cor:monoton}) twice. 

\parbf{\ref{ex:busemann}.}
The first part follows from the angle-sidelength  monotonicity (\ref{cor:monoton}).
An example for the second part can be found among metrics on $\RR^2$ induced by a norm. (Compare to Exercise~\ref{mink+alex=euclid}.)

\parit{Remark.} This exercise is inspired by Busemann's definition \cite{busemann-CBA}.

\parbf{\ref{ex:busemann-CBB} and \ref{ex:busemann-CBA}}; \ref{SHORT.ex:busemann-CBB:a}.
By the function comparison definitions of $\Alex{0}$ space (\ref{comp-kappa}), for any $p\in \spc{L}$ and $\eps>0$ the function $\distfun{p}{}{}$ is $\eps$-concave everywhere sufficiently far from $p$.
Applying the definition of Busemann function we get the result.

The $\CAT0$ case is analogous; we have to apply (\ref{function-comp}) and use $\eps$-convexity.

\parit{\ref{SHORT.ex:busemann-CBB:b}.}
By the definition of Busemann function (see  \ref{prop:busemann}),
\begin{align*}
\exp(\bus_\gamma) 
&=  \lim_{t\to \infty} \exp(\distfun{{\gamma (t)}}{}{} - t)=
\\
&= \lim_{t\to \infty} \left[\exp (\distfun{\gamma (t)}{}{} -t)
+\exp(-\distfun{\gamma (t)}-t)\right]=
\\
&=  \lim_{t\to \infty} \left(2\cdot \exp(-t)\cdot \cosh \circ\distfun{\gamma (t)}{}{}\right).
\end{align*}

By the function comparison definitions of $\CAT\kappa$ space (\ref{function-comp}) or $\Alex{\kappa}$ space (\ref{comp-kappa}),  for any $p\in \spc{U}$ the function $f=\cosh \circ\distfun{p}{}{}$ satisfies $f''+\kappa \cdot f\ge 1$ (respectively  $f''+\kappa \cdot f\le 1$).
The result follows.

\parbf{\ref{ex:noncomplete-globalization}.}
Read \cite{petrunin:globalization}.

\parbf{\ref{ex:fixed-point}.} If $\diam(\spc{L}/G)>\tfrac\pi2$, then for some $x\in \spc{L}$ we have
\[\sup \set{\distfun{G\cdot x}(y)}{y\in \spc{L}}
>
\tfrac\pi2.\]
Use comparison to show that there is a unique point $y^{*}$ that lies at maximal distance from the orbit $G\cdot x$.
Observe that $y^{*}$ is a fixed point.

\parbf{\ref{ex:kleiner}.}
Assume there are 4 such points $x_1,x_2,x_3,x_4$.
Since the space $\spc{L}$ is $\Alex{1}$ it is also $\Alex{0}$.
By the angle comparison, the sum of the angles in any geodesic triangle in an $\Alex{0}$ space is $\ge \pi$.
Therefore the average of the $\mangle\hinge{x_i}{x_j}{x_k}$ is  larger than $\tfrac\pi3$.
On the other hand, since each $x_i$ has space of directions $\le\tfrac12\cdot\mathbb{S}^n$ and the perimeter of any triangle in $\tfrac12\cdot\mathbb{S}^n$ is at most $\pi$, the average of $\mangle\hinge{x_i}{x_j}{x_k}$ is at most $\tfrac\pi3$ --- a contradiction.

\parit{Source:} Based on the main idea in \cite{hsiang-kleiner}.

\parbf{\ref{ex:ccat-(3+1)}.}
Suppose that 
\[\angk\kappa {x^0}{x^1}{x^2}+\angk\kappa {x^0}{x^2}{x^3}<\angk\kappa {x^0}{x^1}{x^3}.\]
Show that
\[\angk\kappa {x^2}{x^0}{x^1}+\angk\kappa {x^2}{x^1}{x^3}+\angk\kappa {x^2}{x^3}{x^0}>2\cdot\pi.\]
Conclude that one can take $p=x^2$.

\parbf{\ref{ex:sba-2+2-short}.}
Modify the configuration in \ref{def:2+2-reformulated}.

\parbf{\ref{ex:berg-nikolaev}.}
Read \cite{sato};
the original proof \cite{berg-nikolaev} is harder to follow.

An example for the second part of the problem can be found among 4-point metric spaces.
It is sufficient to take four vertices of a generic convex quadrangle and increase one of its diagonals slightly;
it will still satisfy the inequality for all relabeling but will fail to meet \ref{def:2+2-reformulated}.

\parbf{\ref{ex:CAT-mnfld=>ext.geod}.}
Suppose that a geodesic $[px]$ is not extendable beyond $x$.
We may assume that $\dist{p}{x}{}<\varpi\kappa$;
otherwise move $p$ along the geodesic toward  $x$.

By the uniqueness of geodesics (\ref{thm:cat-unique}), any point $y$ in a neighborhood $\Omega\ni x$ is connected to $p$ by a unique geodesic path; denote it by $\gamma_y$.
Note that $h_t(y)=\gamma_y(t)$ defines a homotopy, called the  \index{geodesic homotopy}\emph{geodesic homotopy}, between the identity map of $\Omega$ and the constant map with value~$p$.

Since $[px]$ is not extendable, $x\notin h_t(\Omega)$ for any $t<1$.
In particular, the local homology groups vanish at $x$ --- a contradiction.

\parbf{\ref{ex:complete-space-of-dir}.}
Choose a sequence of directions $\xi_n$ at $p$;
by $\gamma_n\:\RR\to \spc{U}$ the corresponding local geodesics.
Since the space $\spc{U}$ is locally compact, we may pass to a converging subsequence of $(\gamma_n)$; it's limit is a local geodesic by Corollary~\ref{cor:loc-geod-are-min}.
Denote the limit by  $\gamma_\infty$ and its direction by $\xi_\infty$.
By comparison, $\xi_\infty$ is a limit of $(\xi_n)$.

\parbf{\ref{mink+CAT=euclid}.} Follow the solution in the \ref{mink+alex=euclid}, reversing all the inequalities.

\parbf{\ref{ex:convexity-CAT0}.} 
It is sufficient to show that if $v$ and $y$ are midpoints of geodesics $[uw]$ and $[xz]$ in $\spc{U}$, then
\[\dist{v}{y}{}\le \tfrac12\cdot(\dist{u}{x}{}+\dist{w}{z}{}).\]

\begin{wrapfigure}{r}{35 mm}
\vskip-5mm
\centering
\includegraphics{mppics/pic-10}
\end{wrapfigure}

Denote by $p$ the midpoint of $[uz]$.
Applying the angle-sidelength  monotonicity (\ref{cor:monoton-cba}) twice, we have
\[\dist{v}{p}{}\le \tfrac12\cdot\dist{w}{z}{}.\]
Similarly we have
\[\dist{y}{p}{}\le \tfrac12\cdot\dist{u}{x}{}.\]
It remains to add these two inequalities and apply the triangle inequality.

\parit{Remark.}
This inequality also follows directly from the majorization theorem (\ref{thm:major}).

\parbf{\ref{ex:equality-for-thin}.}
The only-if part is evident.
Use \ref{thm:defs_of_cat} to show that 
\ref{SHORT.ex:equality-for-thin:side-side}$\Rightarrow$\ref{SHORT.ex:equality-for-thin:vertex-base}$\Rightarrow$\ref{SHORT.ex:equality-for-thin:angle}.
By \ref{thm:defs_of_cat}, condition \ref{SHORT.ex:equality-for-thin:angle} implies that the natural map is distance-preserving on the sides $[\tilde x\tilde y]$ and $[\tilde x\tilde z]$.
Applying it again, we have that condition \ref{SHORT.ex:equality-for-thin:angle} holds for all permutations of the labels $x,y,z$.
Whence the natural map is distance-preserving on all three sides.

\parit{Remark.}
These conditions imply that the natural map can be extended to a distance-preserving map to the solid model triangle.
In fact the image of the line-of-sight map (\ref{def:sight}) is isometric to the model triangle.

\parbf{\ref{ex:busemann-CBA}.}
See the solution of Exercise~\ref{ex:busemann-CBB}.

\parbf{\ref{ex:patchwork}}; \ref{SHORT.ex:patchwork:proper}.
Suppose that $x_n\to x_\infty$, $y_n\to y_\infty$ as $n\to\infty$,
but $[x_ny_n]$ does not converge to $[x_\infty y_\infty]$.
Since the space is proper, we can pass to a subsequence such that $[x_ny_n]$ converges to another geodesic.
That is, we have at least two geodesics between $x_\infty$ and $y_\infty$.

\parit{\ref{SHORT.ex:patchwork:complete}.}
Let $\Delta_n$ be a sequence of solid spherical triangles 
with angle $\tfrac\pi4$ and adjacent sides $\pi-\tfrac1n$.
Let us glue each $\Delta_n$ to $[0,\pi]$ along an isometry of one of the longer sides.
It remains to show that the obtained space $\spc{X}$ is a needed example.

\parit{Source:}
The example \ref{SHORT.ex:patchwork:complete} is taken from \cite[Chapter I, Exercise 3.14]{bridson-haefliger}.

\parbf{\ref{ex:two-rays}.}
Subdivide $Q$ into a a half-plane $A$ bounded by the extension of $\gamma_1$ and the remaining solid angle $B$; it has angle measure measure $\pi-\alpha$.
First glue $B$ along $\gamma_2$, and then glue $A$.
Each time apply the Reshetnyak gluing theorem (\ref{thm:gluing}), to show that the obltained space is $\CAT0$.

\parbf{\ref{ex:reshetnyak-doubling}.}
Suppose that $A$ is not convex.
Then there is a geodesic $[xy]$ with ends in $A$ that does not lie in $A$ completely.
Note that $[xy]$ can be lifted to two different geodesics with the same ends  in the doubling, and apply uniqueness of geodesics (\ref{cor:cat-unique}).

\parbf{\ref{ex:glue-spherical-suspension}.}
Since $K$ is $\pi$-convex, it is $\CAT1$.
By \ref{thm:warp-curv-bound:cat}, the spherical suspension $\Susp K$ is $\CAT1$ as well.
Let us glue $\Susp K$ to $\spc{U}$  along $K$;
according to the Reshetnyak gluing theorem, the resulting space, say $\spc{U}'$, is $\CAT1$.

Consider the geodesic path $\gamma\:[0,1]$ from $p$ to a pole of the suspension in $\spc{U}'$.
Set $K_t=\spc{U}\cap\cBall[\gamma(t),\tfrac\pi2]$.
By \ref{cor:convex-balls}, $K_t$ is $\pi$-convex for any $t$;
monotonicity and continuity of the family should be evident.

\parit{Source:}
This construction was used in \cite{lytchak-petrunin-2020}.
Applying it together with Sharafutdinov retraction leads to another solution of Exercise~\ref{ex:short-retraction-CBA(1)}.

\parbf{\ref{ex:AUB}.}
Apply the Reshetnyak gluing theorem, or its reformulation \ref{thm:gluing2-reformulated}.

\parbf{\ref{ex:fenchel}.} By \ref{thm:major}, there is a majorization $F\:D\to\spc{U}$ of the polygonal line $\beta$.
Show and use that $D$ is a convex plane polygon and its external angles cannot exceed the corresponding external angle of $\beta$.

\parbf{\ref{ex:FM}.} This exercise generalizes the so-called Fáry–Milnor theorem.
An elementary proof is given in the first author and Richard Bishop \cite{alexander-bishop:fm};
another proof is given by Stephan Stadler \cite{stadler}.

\parbf{\ref{ex:isometric-majorization}.}
\textit{(Easier way.)} 
Let 
$(t,s)\mapsto \gamma_t(s)$ be the line-of-sight map 
for $\alpha$ with respect to $\alpha(0)$,
and 
$(t,s)\mapsto \tilde \gamma_t(s)$ be the line-of-sight map 
for $\tilde \alpha$ with respect to $\tilde \alpha(0)$.
Consider the map  $F\:\Conv\tilde \alpha\to \spc{U}$ such that 
$F\:\tilde \gamma_t(s)\mapsto \gamma_t(s)$.

Show that $F$ majorizes $\alpha$
and conclude that $F$ is distance-preserving.

\parit{(Harder way.)}
Prove and apply the following statement together with the Majorization theorem.
\begin{itemize}
\item Let $\alpha$ and $\beta$ be two convex curves in $\Lob2\kappa$.
Assume 
\[\length \alpha=\length\beta<2\cdot\varpi\kappa\]
and there is a short bijecction $f\:\alpha\to\beta$.
Then $f$ is an isometry.
\end{itemize}

\parbf{\ref{ex:bishop}.}
Suppose that points $p,x,q,y$ appear on the curve in that cyclic order.
Assume that the geodesics $[pq]$ and $[xy]$ do not intersect.
Use the argument in the proof of the majorization theorem (\ref{thm:major}) to show that in this case there are nonequivalent majorization maps.

Now we can assume that pairs of geodesics $[pq]$ and $[xy]$ intersect for all choices of points $p,x,q,y$ on the curve in that cyclic order.
Show that in this case the convex hull $K$ of the curve is isometric to a convex figure.

Note that the composition of a majorization map and closest point projection to $K$ is a majorization.
Show and use that the boundary of a convex figure in the plane admits a unique majorization up to equivalence.

\parit{Remark.}
A typical rectifiable closed curve in a $\CAT0$ space can be majorized by more than one convex figure.
There are two exceptions: (1) the majorization map is distance-preserving, and (2) the curve is geodesic triangle.
It is expected that there are no other exceptions;
this question was asked by Richard Bishop in a private conversation.

\parbf{\ref{ex:square}.}
Show that quadrangle $[x^1x^2x^3x^4]$ is majorized by the solid quadrangle $[\tilde x^1\tilde x^2\tilde x^3\tilde x^4]$.
Further show that the majorization is isometric;
argue as in \ref{ex:equality-alexlemma}.

\parbf{\ref{ex:cover-branching-along-2-lines}.}
If $\ell$ and $m$ do not intersect, then the double cover $\spc{X}$ is not simply connected.
In particular, by the Hadamard--Cartan theorem, $\spc{X}$ is not $\CAT0$.

If $\ell$ and $m$ intersect then $\spc{X}$ is a cone over a double cover $\Sigma$ of $\mathbb{S}^2$ branching at two pairs $(x,y)$ and $(v,w)$ of antipodal points.
Suppose $\dist{x}{v}{\mathbb{S}^2}=\ell<\tfrac\pi2$.
Note that the inverse image of $[xv]_{\mathbb{S}^2}$ is a closed geodesic of length $4\cdot\ell<2\cdot\pi$.
Therefore, by the generalized Hadamard--Cartan theorem, $\Sigma$ is not $\CAT1$. Hence $\spc{X}$ is not $\CAT0$ by Theorem \ref{thm:warp-curv-bound:cat}  on curvature of cones.

\parbf{\ref{ex:branching}.}
Let us do the second part first.
Assume $A$ has nonempty interior. 
Note that the space $\tilde{\spc{U}}$ is simply connected and locally isometric to the doubling $\spc{W}$ of $\spc{U}$ in $A$;
that is, any point in $\tilde{\spc{U}}$ has a neighborhood 
that is isometric to a neighborhood of a point in $\spc{W}$.

By the Reshetnyak gluing theorem (\ref{thm:gluing}), $\spc{W}$ is $\CAT0$.
Therefore $\tilde{\spc{U}}$ is locally $\CAT0$;
it remains to apply the Hadamard--Cartan theorem (\ref{thm:hadamard-cartan}).

Let us come back to the general case.
The above argument can be applied to a closed $\eps$-neighborhood of $A$.
After that we need to pass to a limit as $\eps\to 0$.

The first part of the problem follows since a geodesic is a convex set.

\parbf{\ref{ex:closest-point-projection}.}
Let $p\mapsto\bar p$ denote the closest-point projection to $K$.
We need to show that $\dist{\bar p}{\bar q}{}\le\dist{ p}{ q}{}$ for any $p,q\in \spc{U}$.

Assume $p\ne \bar p\ne \bar q\ne q$.
Note that in this case $\mangle\hinge {\bar p}{p}{\bar q}\ge \tfrac\pi2$ and $\mangle\hinge {\bar q}{q}{\bar p}\ge \tfrac\pi2$.
Otherwise a point on the geodesic $[\bar p\bar q]$ would be closer to $p$ or to $q$ than $\bar p$ or $\bar q$ respectively.
The latter is impossible since $K$ is convex and therefore $[\bar p\bar q]\subset K$.

Applying the arm lemma (\ref{lem:arm}), we get the statement.

The cases $p= \bar p\ne \bar q\ne q$ and $p\ne \bar p\ne \bar q= q$ can be done similarly.
The rest of the cases are trivial.

\parbf{\ref{ex:short-retraction-CBA(1)}.}
A more transparent, but less elementary solution via gradient flow is given by Alexander Lytchak and the third author~\cite{lytchak-petrunin-2020}.

\medskip

Without loss of generality, we may assume that $p\in K$.

If $\distfun{K}{x}{}\ge\pi$, then set $\map[2](x)=p$.

Otherwise, if $\distfun{K}{x}{}<\pi$, by the closest-point projection lemma (\ref{lem:closest point}), 
there is a unique point $x^*\in K$ that minimizes distance to $x$;
that is, $\dist{x^*}{x}{}=\distfun{K}{x}{}$.
Let us define $\ell_x$, $\phi_x$ and $\psi_x$ using the following identities:
\begin{align*}
\ell_x&=\dist{p}{x^*}{},
\\
\phi_x&=\tfrac\pi2-\dist[{{}}]{x^*}{x}{},
\\
\sin\psi_x&=\sin\phi_x\cdot\sin\ell_x, 
\quad 0\le \psi_x\le \tfrac\pi2.
\intertext{Let}
\map[2](x)&=\geod_{[px^*]}(\psi_x).
\end{align*}

Note that $\map[2]$ is a retraction to $K$; 
that is,
$\map[2](x)\in K$ for any $x\in \spc{U}$
and 
$\map[2](a)=a$ for any $a\in K$.

Let us show that $\map[2]$ is short.
Given $x,y\in\oBall(K,\tfrac\pi2)$, let
\begin{align*}
x'&=\map[2](x)
&
y'&=\map[2](y)
\\
r&=\dist{x}{y}{}
&
r'&=\dist{x'}{y'}{}
\\
d&=\dist{x^*}{y^*}{}
&
\alpha&=\angk1{p}{x^*}{y^*}.
\end{align*}

Note that 
\[\cos r\le 
\cos\phi_x\cdot\cos\phi_y
-
\cos d\cdot\sin\phi_x\cdot\sin\phi_y.
\eqlbl{eq:cos(r)}\]

Indeed, if $x,y\notin K$,
then 
$\mangle\hinge{x^*}{x}{y*}, 
\mangle\hinge{y^*}{y}{x*}
\ge 
\tfrac\pi2$
and
the inequality~\ref{eq:cos(r)} follows from the arm lemma (\ref{lem:arm}).
If $x\in K$ and $y\notin K$, we obtain \ref{eq:cos(r)} by the angle comparison (\ref{cat-hinge}) 
since $\mangle\hinge{y^*}{y}{x*}\ge \tfrac\pi2$.
In the same way, \ref{eq:cos(r)} is proved 
if $x\notin K$ and $y\in K$.
Finally, if $x,y\in K$, then $\phi_x=\phi_y=\tfrac\pi2$ and $r=d$;
that is, the inequality trivially holds.

Further note that
\[\cos\alpha
=
\frac{\cos d-\cos \ell_x\cdot\cos\ell_y}{\sin\ell_x\cdot\sin\ell_y}.\]
Applying the angle-sidelength  monotonicity (\ref{cor:monoton-cba}), we have
\begin{align*}
\cos r'&\ge
\cos\psi_x\cdot\cos\psi_y
-
\cos \alpha \cdot\sin\psi_x\cdot\sin\psi_y=
\\
&=
\cos\psi_x\cdot\cos\psi_y
-(\cos d-\cos \ell_x\cdot\cos\ell_y)\cdot\sin\phi_x\cdot\sin\phi_y\ge
\\
&\ge \cos\psi_x\cdot\cos\psi_y
-\cos d\cdot\sin\phi_x\cdot\sin\phi_y.
\end{align*}

Note that 
$\psi_x\le \phi_x$
and
$\psi_y\le \phi_y$;
in particular,
\[
\cos\phi_x\cdot\cos\phi_y\le \cos\psi_x\cdot\cos\psi_y.
\]
Hence 
\[\cos r'\ge \cos r;\]
that is, the restriction $\map[2]|_{\oBall(K,\tfrac\pi2)}$ is short.
Clearly $\map[2]$ is continuous. 
Since the complement of $\oBall(K,\tfrac\pi2)$ is mapped to $p$,
 $\map[2]$ is short; that is,
\[r'\le r \eqlbl{eq:cos=<cos}\]
for any $x,y\in\spc{U}$.

If we have equality in \ref{eq:cos=<cos}
then 
\[\cos\ell_x\cdot\cos\ell_y\cdot\sin\phi_x\cdot\sin\phi_y=0.\]
If $K\subset \oBall(p,\tfrac\pi2)$, then $\ell_x,\ell_y<\tfrac\pi2$, 
which implies that $x\in K$ or $y\in K$.
Without loss of generality, we may assume that $x\in K$.

It remains to show that if $y\notin K$ 
then the inequality~\ref{eq:cos=<cos}
is strict.
If $\distfun{K}{y}{}\ge\tfrac\pi2$, then \ref{eq:cos=<cos} holds since 
the left-hand side is $<\tfrac\pi2$
while the right-hand side is $\ge \tfrac\pi2$.
If $\distfun{K}{y}{}<\tfrac\pi2$, then $\phi_y>0$. Clearly $\psi_y<\phi_y$,
hence the inequality~\ref{eq:cos=<cos} is strict.
\qeds

Below you will find a geometric way to think about the given construction; 
it is close to the construction in the proof of Kirszbraun's theorem (\ref{thm:kirsz+}).

\parit{Geometric interpretation of the map $\map[2]$.}
Let $\mathring{\spc{U}}=\Cone \spc{U}$, and 
denote by $\mathring{K}$ the subcone of $\mathring{\spc{U}}$ spanned by $K$.
The space $\spc{U}$ can be naturally identified with the unit sphere in $\mathring{\spc{U}}$, 
that is, the set 
\[\set{z\in \mathring{\spc{U}}}{|z|=1}.\]

According to \ref{thm:warp-curv-bound:cat}, $\mathring{\spc{U}}$ is $\CAT0$.
Note that $\mathring{K}$ forms a convex closed subset of $\mathring{\spc{U}}$.
According to \ref{lem:closest point}, for any point $x$ there is a unique point $\hat x\in \mathring{K}$
that minimizes the distance to $x$,
that is, $\dist{\hat x}{x}{}=\distfun{K}{x}{}$.
(If $|\hat x|\ne0$, then in the notation above we have
$x^*=\tfrac1{|\hat x|}\cdot\hat x$.)

Consider the half-line $t\mapsto t\cdot p$ in  $\mathring{\spc{U}}$.
By comparison, 
for given $s\in \mathring{\spc{U}}$
the geodesics $\geod_{[s\ t\cdot p]}$ converge as $t\to\infty$ to a half-line, 
say $\alpha_s\:[0,\infty)\to \mathring{\spc{U}}$.

Note that if $|x|=1$, then $|\hat x|\le 1$.
By assumption, for any $a\in K$ the function $t\mapsto |\alpha_a(t)|$ is monotonically increasing.
Therefore there is a unique value $t_x\ge 0$ such that
$|\alpha_{\hat x}(t_x)|=1$.
Define $\map[2]\:\spc{U}\to K$
 by 
\[\map[2](x)=\alpha_{\hat x}(t_x).\]

\parbf{\ref{ex:cats-cradle}.}
Prove that the angle comparison (\ref{cat-hinge}) holds.

\parbf{\ref{ex:Hadamard--Cartan}.}
Mimic the proof of the Hadamard--Cartan theorem.

\parbf{\ref{ex:CBB+CBA}.}
Note that it is sufficient to show that any finite set of points $x^1,\dots,x^n\in\spc{X}$ lies in an isometric copy of a Euclidean polyhedron.

Observe that $\spc{X}$ is $\Alex0$ and $\CAT0$ at the same time.
Show that there is a unique point $p$ that minimizes the sum $\dist{p}{x^1}{}+\dots+\dist{p}{x^n}{}$.
Note that the vectors $v^i=\ddir{p}{x^i}$ lie in a linear subspace of $\T_p$.
Moreover if $K$ is the convex hull of $v_i$, then the origin of $\T_p$ lies in the interior of $K$ relative to its affine hull.
Finally observe that the exponential map is defined on all of $K$ and is distance-preserving.
The statement follows since the exponential map sends $v^i\mapsto x^i$ for each $i$.

\parbf{\ref{ex:5-point-CBA=>CBB}.}
The answers are $s\le \sqrt3$ and $s\le 2$ respectively.

Let us start with the $\CAT0$ case.
The upper bound $s\le \sqrt3$ follows from (2+2)-point comparison.
The Euclidean space works as an example if $s$ is smaller than the large diagonal of the double pyramid with unit side (that is, if $s\le 2\cdot\sqrt{2/3}$).
Otherwise it can be embedded into a product of the real line with a two-dimensional cone.

For the $\Alex0$ case, the needed space can be constructed by doubling  a polyhedron $K\subset\EE^3$ in its boundary.
The obtained space is $\Alex0$ by \ref{thm:poly-CBB};
the same follows from Perelman's doubling theorem \cite{perelman:spaces2}.
We assume that the points correspond to vertices of a regular tetrahedron with 3 vertices on the boundary of $K$ and one  in its interior; this point corresponds to a pair of points in the doubling at distance $s$ from each other.

\parit{Remark.}
The $\CAT0$ case also follows from \cite{toyoda,lebedeva-petrunin:toyoda}.

\parbf{\ref{ex:cbb-wald}.}
Choose a quadruple of points $p,q,r,s$. 
Suppose that it admits a distance-preserving embedding into some $\Lob2{\Kappa}$ for some $\Kappa\ge \kappa$.
Then 
\[\angk\Kappa p{q}{r}
+\angk\Kappa p{r}{s}
+\angk\Kappa p{a}{q}\le 2\cdot\pi.\]
Applying monotonicity of the function $\kappa\mapsto\angk\kappa p{q}{r}$ (\ref{k-decrease}) shows  that
\[\angk\kappa p{q}{r}
+\angk\kappa p{r}{s}
+\angk\kappa p{s}{q}\le 2\cdot\pi.\]
Since the quadruple $p,q,r,s$ is arbitrary, the if part follows.

Now let us prove the only-if part.
Denote by $\sigma$ the exact upper bound on values $\Kappa\ge \kappa$ such that all model triangles with the vertices $p,q,r,s$ are defined.

Recall that $\angk{\Kappa+} p{q}{r}$ denotes extended angle (\ref{def:extended-angle}).
Observe that if 
\[\angk{\Kappa+} p{q}{r}
+\angk{\Kappa+} p{r}{s}
+\angk{\Kappa+} p{s}{q}= 2\cdot\pi\eqlbl{eq:Kappa3}\]
for some $\sigma\ge \Kappa\ge \kappa$, then the quadruple admits a distance-preserving embedding into $\Lob2\Kappa$.

Observe that the left-hand side of \ref{eq:Kappa3} is continuous in $\Kappa$.
Since $\spc{L}$ is $\Alex\kappa$, for $\Kappa=\kappa$ the left-hand side cannot exceed $2\cdot \pi$.
Therefore it remains smaller than $2\cdot\pi$ for all $\sigma\ge \Kappa\ge \kappa$;
moreover the same holds for all permutations of the labels $p,q,r,s$.

Note that we can assume the perimeter of the triple $q,r,s$ is $2\cdot\varpi{\sigma}$, and use this and the overlap lemma (\ref{lem:extend-overlap}) to arrive at a contradiction.

According to our definition, the real line is $\Alex\kappa$ for any $\kappa\in\RR$,
but it does not satisfy the property for $\kappa>0$. 
The condition $\kappa\le 0$ was used just once to ensure that the $\kappa$-model triangles with the vertices $p,q,r,s$ are defined.
One can assume instead that perimeters of all triangles in $\spc{L}$ are at most $2\cdot\varpi\kappa$.
This condition holds for all complete length $\Alex\kappa$ spaces of dimension at least 2; see \ref{diam-k>0}.

\parbf{\ref{ex:sturm}.}
Let $\tilde p,\tilde x_1,\dots,\tilde x_n$ be the array in $\EE^n$ provided by the (1+\textit{n})-point comparison (\ref{thm:pos-config}).
We may assume that $\tilde p$ is the origin of $\EE^n$.

Consider an $n{\times}n$-matrix $\tilde M$ with components 
\[\tilde m_{i,j}=\tfrac12\cdot(\dist[2]{\tilde x_i}{\tilde p}{}+\dist[2]{\tilde x_j}{\tilde p}{}-\dist[2]{\tilde x_i}{\tilde x_j}{}).\]
Note that $\tilde m_{i,j}=\langle\tilde x_i,\tilde x_j\rangle$.
It follows that $\tilde M=A\cdot A^\top$ for an $n{\times}n$-matrix $A$ that defines a linear transformation sending the standard basis to the array $\tilde x_1,\dots,\tilde x_n$.
Therefore
\[\bm{s}\cdot \tilde M\cdot \bm{s}^\top=|A^\top\cdot \bm{s}^\top|^2 \ge 0\]
for any vector $\bm{s}$.
Further show and use that
\[\bm{s}\cdot M\cdot \bm{s}^\top\ge \bm{s}\cdot \tilde M\cdot \bm{s}^\top\]
for any vector $\bm{s}=(s_1,\dots,s_n)$ with nonnegative components.

\parbf{\ref{6-point-comparison}.} Apply the (5+1)-point comparison (\ref{thm:pos-config}).

\parbf{\ref{ex:(3+1)-nonsufficient}.}
It is sufficient to construct a metric on the set of points $\{p$, $x^1$, $x^2$, $x^3$, $x^4\}$ that does not satisfy (1+4)-point comparison but does satisfy all (1+3)-point comparisons.
To do this, set  $x^i$ to be  the vertices of a regular tetrahedron in $\EE^3$. Suppose $p$ is its center and reduce the distances $\dist{p}{x^i}{}$ slightly.

\parit{Remark.}
There are examples of 6-point metric spaces that satisfy all (1+5)-point comparisons, but do not admit embedding into a complete length $\Alex{0}$ space \cite{lebedeva-petrunin-zolotov}.

\parbf{\ref{ex:strut+embedding}.}
By the (1+\textit{n})-point comparison (\ref{thm:pos-config}), there is a point array $\tilde p,\tilde a^0,\dots,\tilde a^m\in \Lob{m+1}\kappa$ such that
\[\dist{\tilde p}{\tilde a^i}{}=\dist{p}{a^i}{}\quad \text{and}\quad \dist{\tilde a^i}{\tilde a^j}{}\ge\dist{a^i}{a^j}{}\]
for all $i$ and $j$.

For each $i$, set 
$\tilde \xi^i=\dir{\tilde p}{\tilde a^i}\in\mathbb{S}^m=\Sigma_{\tilde p}(\Lob{m+1}\kappa)$.
Note that 
\[\dist{\tilde \xi^i}{\tilde \xi^j}{\mathbb{S}^m}\ge \angk\kappa{p}{ a^i}{ a^j}>\tfrac\pi2.\]

Consider two matrices $S$ and $\tilde S$ with components
$s_{i,j}=\langle\tilde \xi^i,\xi^j\rangle$
and
$\tilde s_{i,j}=\cos[\angk\kappa{p}{a^i}{a^j}]$.
By construction, $S\ge 0$; that is $\bm{v}\cdot S\cdot \bm{v}^\top\ge 0$ for any vector $\bm{v}$.

Observe that it is sufficient to show that $\tilde S\ge 0$.
The latter follows since $s_{i,j}\le \tilde s_{i,j}\le 0$ if $i\ne j$ and
$s_{i,j}= \tilde s_{i,j}=1$ if $i=j$.

\parbf{\ref{ex:flat-in-CAT}.}
Set $\tilde Q=\Conv\{\tilde x^0,\tilde x^1,\dots,\tilde x^\kay\}$.
By Kirszbraun's theorem, the map $\tilde x^i\mapsto x^i$ can be extended to a short map $F\:\tilde Q\to\spc{L}$;
it remains to show that the map $F$ is distance-preserving.

Consider the logarithm map $G\:x\mapsto \ddir{x_0}x$; note that $G$ is short.
Observe that the composition $G\circ F$ is distance-preserving.
Therefore $F$ is distance-preserving;
in particular we can take $Q=F(\tilde Q)$.

\parbf{\ref{ex:flat-in-CBB}.}
Consider vectors $v^i=\ddir{x^0}{x^i}\in\T_{x^0}$.
Show that all the $v^i$ lie in a linear subspace of $\T_{x^0}$ and that $x^i\mapsto v^i$ is distance-preserving.
It follows that we can identify the convex hull $K$ of  the $v^i$ with the convex hull of  the $\tilde x^i$.

Note that the gradient exponential map $\gexp_{x_0}$ maps $v^i$ to $x^i$.
By assumption, 
\[\dist{v^i}{v^j}{}=\dist{x^i}{x^j}{}\eqlbl{eq:vv=xx}\]
for all $i$ and $j$.
By \ref{thm:prop-gexp}, $\gexp_{x_0}$ is a short map.
By \ref{eq:vv=xx}, $\gexp_{x_0}$ cannot be strictly short at a pair of points in $K$.
That is, $\gexp_{x_0}$ is distance-preserving on $K$.

\parbf{\ref{CBA-n-point}.}
Apply \ref{thm:kirsz} for each of the following maps
\begin{itemize}
\item $f_0\:\tilde x\mapsto x$, $\tilde p^1\mapsto p^1$, $\tilde q^1\mapsto q^1$;
\item $f_i\:\tilde p^i\mapsto p^i$, $\tilde p^{i+1}\mapsto p^{i+1}$, $\tilde q^i\mapsto q^i$, $\tilde q^{i+1}\mapsto q^{i+1}$ for $1\le i<n$;
\item $f_n\:\tilde y\mapsto y$, $\tilde p^n\mapsto p^n$, $\tilde q^n\mapsto q^n$.
\end{itemize}
Denote by $F_i$ the short extension of $f_i$.
Observe and use that $F_{i-1}(\tilde z_i)\z=F_{i}(\tilde z_i)$ for each $i$.

\parbf{\ref{ex:petrunin-stadler}.} Consider the space $\spc{Y}^{\spc{X}}$ of all maps $\spc{X}\to \spc{Y}$ equipped with the product topology.

Denote by $\mathfrak{S}_F$ the set of maps $h\in \spc{Y}^\spc{X}$ such that the restriction $h|_F$  is short and agrees with $f$ in $F\cap A$.
Note that the sets $\mathfrak{S}_F\subset \spc{Y}^\spc{X}$ are closed and any finite intersection of these sets is nonempty.

According to Tikhonov's theorem, $\spc{Y}^{\spc{X}}$ is compact.
By the finite intersection property, the intersection $\bigcap_F\mathfrak{S}_F$ for all finite sets $F\subset X$ is nonempty.
Hence the statement follows.

\parit{Source:} This statement appears in \cite{petrunin-stadler}; it is an analogous of the finite+one lemma (\ref{lem:kirsz-neg:new}).

\parbf{\ref{ex:isbell}.}
The Kuratowsky embedding is a distance-preserving map of $\spc{X}$ into the space of bounded functions $\spc{X}$ equipped with the metric induced by the sup-norm (Section~\ref{Kuratowsky embedding}).
It remains to show that the latter space is injective.

The second part of the exercise is a classical result of John Isbell \cite{isbell} which was rediscovered several times after him; for more on the subject see lecture notes of the third author \cite{petrunin2020pure}.

\parbf{\ref{ex:warp=<}.}
It is sufficient to show that the natural map $\spc{B}\warp{g}\spc{F}\to \spc{B}\warp{f}\spc{F}$ is short.
The latter follows from the fiber-independence theorem (\ref{thm:fiber-independence}).

\parbf{\ref{ex:convexity-in-cone}.}
Show and use that any geodesic path in $\Cone^\kappa\spc{F}$ projects to a reparametrized geodesic in $\spc{F}$ of length less than $\pi$.

\parbf{\ref{ex:spherical-join}.}
By \ref{thm:warp-curv-bound:cbb:a}, the space $\spc{U}$, $\spc{V}$, or $\spc{U}\star\spc{V}$ is $\Alex1$ if and only if $\Cone\spc{U}$, $\Cone\spc{V}$, or $\Cone(\spc{U}\star\spc{V})=\Cone\spc{U}\times\Cone\spc{V}$ is $\Alex0$ respectively.

By \ref{thm:warp-curv-bound:cbb:S}, the space $\spc{U}$, $\spc{V}$, or $\spc{U}\star\spc{V}$ is $\CAT1$ if and only if $\Cone\spc{U}$, $\Cone\spc{V}$, or $\Cone(\spc{U}\star\spc{V})=\Cone\spc{U}\times\Cone\spc{V}$ is $\CAT0$ respectively.

It remains to show that the product of two spaces is $\Alex0$ or $\CAT0$ if and only if each space is $\Alex0$ or $\CAT0$ respectively.

\parbf{\ref{ex:metric tree}.}
Apply Reshetnyak gluing theorem (\ref{thm:gluing}) several times.

\parbf{\ref{ex:poly-unique-geodesic}.}
Assume $\spc{P}$ is not $\CAT0$.
Then by \ref{thm:PL-CAT}, a link $\Sigma$ of some simplex contains a closed local geodesic $\alpha$ with length $4\cdot\ell<2\cdot\pi$.
We can assume that $\Sigma$ has minimal possible dimension;
then by \ref{thm:PL-CAT}, $\Sigma$ is locally $\CAT1$.

Divide $\alpha$ into two equal arcs $\alpha_1$ and $\alpha_2$.

Assume $\alpha_1$ and $\alpha_2$ are length-minimizing, and 
parametrize them by $[-\ell,\ell]$.
Fix a small $\delta>0$ and 
consider the two curves in $\Cone\Sigma$ given in polar coordinates by 
\[\gamma_i(t)=(\alpha_i(\arctan \tfrac t\delta),\sqrt{\delta^2+t^2}).\]
Show that the curves $\gamma_1$ and $\gamma_2$ are geodesics in $\Cone\Sigma$ having common endpoints.

Observe that a small neighborhood of the tip of $\Cone\Sigma$ admits a distance-preserving embedding into~$\spc{P}$.
Hence we can construct two geodesics $\gamma_1$ and $\gamma_2$ in $\spc{P}$ with common endpoints.

It remains to consider the case where $\alpha_1$ (and therefore $\alpha_2$) is not length-minimizing.

Pass to a maximal length-minimizing arc $\bar\alpha_1$ of $\alpha_1$.
Since $\Sigma$ is locally $\CAT1$, by the no-conjugate-point theorem (\ref{thm:no-conj-pt}) 
there is another geodesic $\bar\alpha_2$ in $\Sigma_p$ that shares endpoints with $\bar\alpha_1$.
It remains to repeat the above construction for the pair $\bar\alpha_1$, $\bar\alpha_2$.

\parit{Remark.}
By \ref{thm:cat-unique} the converse holds as well.
This problem was suggested by Dmitri Burago.

\parbf{\ref{ex:polyKk}.} Apply \ref{thm:tan-is}, \ref{thm:poly-CBB}, and \ref{thm:PL-CAT}.

\parbf{\ref{ex:barycenric-flag}.}
Observe and use that (1) in the barycentric subdivision every vertex corresponds to a simplex of the original triangulation,
and (2) a simplex of the subdivision corresponds to a decreasing sequence of simplexes in the original triangulation. 

\parit{Remark.}
The second statement, \textit{any finite  simplicial complex is homeomorphic to a compact length $\CAT1$ space}, is due to Valerii Berestovskii \cite{berestovskii}.

\parbf{\ref{ex:obtuce-flag}.}
Use induction on the dimension  to prove that if in a spherical simplex $\triangle$ every edge is at least $\tfrac\pi2$, then 
all dihedral angles of $\triangle$ are at least~$\tfrac\pi2$.

The rest of the proof goes along the same lines as the proof of the flag condition (\ref{thm:flag}).
The only difference is that a geodesic may spend time \textit{at least} $\pi$ on each visit to $\Star_v$.

\parit{Remark.}
It is not sufficient to assume only that all the dihedral angles of the simplexes are at least~$\tfrac\pi2$. 
Indeed, the two-dimensional sphere with the interior of a small rhombus removed is a spherical polyhedral space glued from four triangles with angles at least~$\tfrac\pi2$.
On the other hand, the boundary of the rhombus is a closed local geodesic in this space and has length less than $2\cdot\pi$.
Therefore the space cannot be $\CAT1$.

\parbf{\ref{ex:short+commuting}.}
Observe that if we glue two copies of spaces along $A_i$, then the copies of $A_j$ for some $j\ne i$ form a convex subset in the glued space.
Use this and the Reshetnyak gluing theorem (\ref{thm:gluing}) $n$ times, once for each label of the edges.

\parbf{\ref{ex:space-of-trees}.}
The space $\spc{T}_n$ has a natural cone structure whose vertex is the  completely degenerate tree --- all its edges have zero length.

Note that the space $\Sigma$
over which the cone is taken comes naturally with a triangulation 
by right-angled spherical simplexes.
Each simplex corresponds to the combinatorics of a possibly degenerate tree.

Note that the link of any simplex of this triangulation satisfies the no-triangle condition.
Indeed, fix a simplex $\triangle$ of the complex;
suppose it is described by a possibly degenerate topological tree $t$.
A triangle in the link of  $\triangle$ can be described by three ways to resolve a degeneracy of $t$ by adding one edge, where 
(1) any pair of these resolutions can be done simultaneously, but (2) all three cannot be done simultaneously.
Direct inspection shows that this is impossible.

By Proposition~\ref{prop:no-trig}, our complex is flag.
It remains to apply the flag condition (\ref{thm:flag}) and \ref{thm:warp-curv-bound:cbb:a}.

\parbf{\ref{ex:cubical-complex}.}
Apply the flag condition (\ref{thm:flag}) and Theorem~\ref{thm:warp-curv-bound:cbb:S}.

\parbf{\ref{ex:norays}.}
Consider a cube in the $\ell^2$-space defined by $|x_i|\le 1$.

\parbf{\ref{ex:d_q dist_p(v)=-<dri p q, v>-CAT} and \ref{ex:d_q dist_p(v)=-<dri p q, v>}.}
Apply the strong angle lemmas
(\ref{lem:strong-angle-cba}
and \ref{lem:strong-angle}).

\parbf{\ref{ex:grad(dist)}.} Apply \ref{ex:d_q dist_p(v)=-<dri p q, v>-CAT}.

\parbf{\ref{ex:|grad|=1}.} Apply \ref{thm:almost.geod}.

\parbf{\ref{ex:df(v)=<grad f,v>}.}
Since $\alpha$ is Lipschitz, so is $f\circ\alpha$.
By the standard Rademacher theorem, the derivative $(f\circ\alpha)'$ is defined almost everywhere.
In particular, 
\[(\dd_{\alpha(t)}f)(\alpha^+(t))+(\dd_{\alpha(t)}f)(\alpha^-(t))\ae0.\]

Further, by the extended Rademacher theorem (more precisely its 1-dimensional case; see Proposition~\ref{prop:Rademacher-dim=1}),
we have 
\[\alpha^+(t)+\alpha^-(t)\ae0.\]
In particular,
\[\<\nabla_{\alpha(t)}f,\alpha^+(t)\>+\<\nabla_{\alpha(t)}f,\alpha^-(t)\>
\ae0.\]

Finally, by the definition of gradient, we have 
\[\<\nabla_{\alpha(t)}f,\alpha^\pm(t)\>\ge (\dd_{\alpha(t)}f)(\alpha^\pm(t)).\]
Hence the result follows.

\parbf{\ref{ex:d dist(grad)<0}.}
Let us pass to the ultrapower $\spc{L}^\o\supset\spc{L}$.
Argue as in \ref{ex:d_q dist_p(v)=-<dri p q, v>} to show that there is a geodesic $[pb]_{\spc{L}^\o}$ such that 
\[\dd_p\distfun b{}(\dir pa)=-\langle\dir pb,\dir pa\rangle\]

It follows that
\begin{align*}(\dd_p\distfun{a}{}{})(\nabla_p\distfun{b}{}{})&\le-\langle\dir pa,\nabla_p\distfun{b}{}{} \rangle\le
\\
&\le -\dd_p\distfun{b}{}{}(\dir pa)=
\\
&=\langle\dir pb, \dir pa\rangle=
\\
&= \cos\mangle\hinge p a b _{\spc{L}^\o}\le 
\\
&\le \cos\angk\kappa p a b.
\end{align*}

\parbf{\ref{ex:compact-dimension-cbb}.}
Suppose that $\spc{L}$ is infinite-dimensional.
Denote by $\Omega_m\subset \spc{L}$ the set of all points $p$ with $\rank_p\ge m$.
Evidently $\Omega_1\supset \Omega_2\supset\dots$, and $\Omega_m$ is open for each $m$.

By \ref{LinDim+}, each $\Omega_m$ is dense in $\spc{L}$.
Hence there is a G-delta dense set of points $p\in\spc{L}$ such that $\rank_p=\infty$.
It follows that $\Sigma_p$ is not compact.

\parit{Source:} It was suggested by Alexander Lytchak.

\parbf{\ref{ex:sharafutdinov}.}
Choose a finite sequence $t_0<\dots<t_n$.
Denote by $\Phi_{t_i}$ the composition of the closest-point projections to $K_{t_0},\dots, K_{t_i}$.
Pass to a limit of the $\Phi_{t_i}$ as the sequence becomes denser in the parameter interval. 
Show and use that the limit $\phi_t$ does not depend on the choice of the sequences.

\parbf{\ref{ex:elf-contracting}.}
Let $\ell(t)=\dist{\alpha(t)}{\alpha(t_3)}{}$.
Note that 
\[\ell'(t)\le -\langle \nabla_{\alpha(t)}f,\dir{\alpha(t)}{\alpha(t_3)}\rangle.\]

Observe that the function $t\mapsto f\circ\alpha(t)$ is nondecreasing;
in particular, $f(\alpha(t_1))\le f(\alpha(t_2))\le f(\alpha(t_3))$.
Therefore 
\begin{align*}\langle \nabla_{\alpha(t)}f,\dir{\alpha(t)}{\alpha(t_3)}\rangle&\ge\dd_{\alpha(t)}f(\dir{\alpha(t)}{\alpha(t_3)})\ge 0
\end{align*}
for any $t\in[t_1,t_2]$.
Therefore $\ell'\le 0$ for any $t\in[t_1,t_2]$.
Hence the statement.

\parbf{\ref{ex:grad-curve-condition}.} 
Without loss of generality, we may assume that $(f\circ\alpha)'(t)>0$ for any $t$.
Let $\hat\alpha$ be the arclength reparametrization of $\alpha$.
Note that 
\[(f\circ\hat\alpha)'(s)\ge |\nabla_{\hat\alpha(s)}f|\]
almost everywhere.
Therefore, by Theorem~\ref{thm:grad-like-2nd-def}, $\hat\alpha$ is a gradient-like curve.
It remains to apply Lemma~\ref{lem:grad--grad-like}.

\parbf{\ref{ex:grad-curve-analitic}.}
Use \ref{ex:grad-curve-condition} to prove the only-if part.

To prove the if part, set $h(z)=\tfrac12\cdot\dist[2]{x}{z}{}$.
If $\alpha$ is an $f$-gradient curve, then 
\begin{align*}
(h\circ\alpha)^+&\ge \dist{\alpha(t)}{x}{}\cdot\langle\dir{\alpha(t)}{x},\nabla_{\alpha(t)}f\rangle\ge
\\
&\ge \dist{\alpha(t)}{x}{}\cdot\dd_{\alpha(t)}f(\dir{\alpha(t)}{x})\ge 
\\
&\ge f(x)-f\circ\alpha(t).
\end{align*}
It remains to integrate the inequality and observe that $f\circ\alpha$ is nondecreasing.

\parbf{\ref{ex:geodesic}.}
Consider $(x,\kappa)$- and $(z,\kappa)$-radial curves that start at $y$
and observe that they form a geodesic from $x$ to $z$.
(Compare to Exercise~\ref{ex:flat-in-CBB}.)

\parbf{\ref{ex:gexp}.} Set $q=p+v$ and $q'=\gexp_pv$. 
By radial comparison, $\dist{q'}{x}{}\le \dist{q}{x}{}$ for any $x\in \spc{L}$.
If $q\in \spc{L}$, this implies that $q=q'$.
Otherwise note that $q'$ lies on the boundary line of $\spc{L}$, and $\proj(q)$ is the only point on this line that satisfies the inequality.

\parbf{\ref{ex:inv-gexp}.}
By the angle comparison,
$|\nabla_x\distfun p|\ge-\cos \angk\kappa xpq$.

Choose a $(p,\kappa)$-radial curve $\alpha$ that starts at $p$.
Observe that 
\[(\distfun p\circ \alpha)^+(t)\ge-|\alpha^+(t)|\cdot \cos \angk\kappa {\alpha(t)}p q\]
and
\[(\distfun q\circ\alpha)^+(t)\ge -|\alpha^+(t)|.\]
Therefore $t\mapsto\angk\kappa q{\alpha(t)}p$  is nondecreasing, hence the result.

{\small\sloppy

\input{arXiv.ind}

\def\emph{\textit}

\printbibliography[heading=bibintoc]

}

\end{document}